\documentclass[11pt,reqno]{amsart}
\usepackage[T1]{fontenc}
\usepackage{lmodern,mathtools,amssymb,booktabs,microtype}
\usepackage[letterpaper,margin=1.12in]{geometry}
\usepackage{xcolor}
\usepackage[colorlinks=true,linkcolor=blue!45!black,citecolor=blue!45!black,urlcolor=blue!45!black]{hyperref}
\newtheorem{theorem}{Theorem}[section]
\newtheorem{proposition}[theorem]{Proposition}
\newtheorem{lemma}[theorem]{Lemma}
\newtheorem{stagelemma}{Lemma}[section]

\newtheorem{corollary}[theorem]{Corollary}
\theoremstyle{definition}

\theoremstyle{remark}

\numberwithin{equation}{section}
\newcommand{\T}{\mathbb T}
\newcommand{\R}{\mathbb R}
\newcommand{\C}{\mathbb C}

\newcommand{\Z}{\mathbb Z}

\newcommand{\norm}[1]{\left\|#1\right\|}
\newcommand{\bk}[1]{\langle #1\rangle}

\allowdisplaybreaks[2]
\title[Infinite cascades for cubic NLS]{Infinite cascades for the defocusing cubic nonlinear Schr\"odinger equation on the square torus}
\author{Xi Chen}
\address{Xi Chen, Department of Mathematics and Computer Science, University of Basel, Spiegelgasse 1, CH--4051 Basel, Switzerland}
\email{{\tt xi01.chen@unibas.ch}}
\date{}
\begin{document}
\begin{abstract}
For each finite \(s>1\) we construct a solution of the defocusing cubic
nonlinear Schr\"odinger equation on the square two-torus with initial
data in \(H^s\) whose \(H^s\) norm tends to infinity as \(t\to+\infty\).
The growth can be arbitrarily slow. More precisely, the squared norm
agrees, up to a bounded error at all positive times, with the number
of selected axial frequencies \(K_j\) for which \(K_j^s\le t\).
The central construction couples a finite tree of orthogonal frequency
pairs to a finite-dimensional family of one-dimensional solutions of the same equation.
Each transfer leaves a nonvanishing Fourier coefficient on the horizontal
axis. Control on an entire half-line allows successive transfers to
retain the coefficients produced earlier. At and above the index two,
we match the second-order normal frequency shift separately on each
triangle using harmonics of the axial solution.
For every fixed \(2<s<\infty\) we also obtain a solution for which
\(\log\|u(t)\|_{H^s}\asymp\sqrt{\log t/\log\log t}\) at every
sufficiently large time.
\end{abstract}
\maketitle
\section{Introduction}\label{sec:intro}

On the square torus \(\T^2\coloneqq (\R/2\pi\Z)^2\) consider
\begin{equation}\label{eq:nls}
 i\partial_tu=-\Delta u+|u|^2u,\qquad u(0)=u_0 .
\end{equation}
Here \(u(t,\mathbf x)\) is complex valued, \(\mathbf x\coloneqq (x,y)\), and
\(\Delta\coloneqq \partial_x^2+\partial_y^2\). We write \(d\mu\coloneqq d x\,d y/(2\pi)^2\)
for normalized Haar measure. For \(k=(k_x,k_y)\in\Z^2\) put
\(\langle k\rangle\coloneqq (1+|k|^2)^{1/2}\), where \(|k|\) is the Euclidean
length, and use the same bracket for integer indices in one dimension.
Our Fourier and Sobolev conventions are
\[
 \widehat u(k,t)\coloneqq \int_{\T^2}u(t,\mathbf x)e^{-ik\cdot\mathbf x}\,d\mu,
 \qquad
 \|u(t)\|_{H^s}^2\coloneqq \sum_{k\in\Z^2}\langle k\rangle^{2s}|\widehat u(k,t)|^2 .
\]
We call \((n,0)\) an axial frequency and a field independent of \(y\)
an axial field. The conserved mass \(M(u)\) and energy \(E(u)\) are the
integrals against \(d\mu\) of \(|u|^2\) and
\(\tfrac12|\nabla u|^2+\tfrac14|u|^4\), respectively. They bound the
\(H^1\) norm along a trajectory. Above \(H^1\), a sequence of smaller
and smaller amounts of mass can carry a growing Sobolev norm. The
question is whether the autonomous equation can produce this behavior
from a single initial datum.

\begin{theorem}\label{thm:main}
For every real \(1<s<\infty\) there is \(u_0\in H^s(\T^2)\) such that the
global solution of \eqref{eq:nls} satisfies
\[
                    \lim_{t\to+\infty}\|u(t)\|_{H^s}=+\infty .
\]
The datum can be chosen with arbitrarily small \(H^s\) norm.
\end{theorem}

The growth in Theorem~\ref{thm:main} comes from distinct axial Fourier
coefficients that persist once they have been formed. The next result
describes how these contributions build up in physical time. The
quantity \(\langle k\rangle^s|\widehat u(k,t)|\) is the normalized
Fourier modulus at \(k\). Its square is precisely the contribution of
this frequency to \(\|u(t)\|_{H^s}^2\). In the formula below, \(\#\)
counts the indices \(j\in\mathbb N\coloneqq \{1,2,\ldots\}\) in the indicated set.

\begin{theorem}\label{thm:clock}
Fix \(1<s<\infty\), and let \(G:[0,\infty)\to(0,\infty)\) be nondecreasing
with \(G(t)\to\infty\). One can choose a datum as in
Theorem~\ref{thm:main} and increasing positive integers \(K_j\to\infty\)
so that
\begin{equation}\label{eq:count}
 \sup_{t\ge0}\left|\|u(t)\|_{H^s}^2-
                       \#\{j:K_j^s\le t\}\right|<\infty,\qquad
 \|u(t)\|_{H^s}=o(G(t)).
\end{equation}
For each \(j\), there is a time \(t_j\) such that
\(\langle K_j\rangle^s|\widehat u((K_j,0),t)|\ge 1/2\) for every
\(t\ge t_j\). Given any nondecreasing function
\(\ell:[1,\infty)\to(0,\infty)\) tending to infinity, the choices can
be made with \(t_j\le K_j^s\ell(K_j)\) for all \(j\).
\end{theorem}

For \(s>2\) a different choice of initial datum gives a quantitative
growth rate at all sufficiently large times.

\begin{theorem}\label{thm:two-sided}
For each fixed \(2<s<\infty\) there exist \(u_0\in H^s(\T^2)\),
\(0<c_s\le C_s<\infty\) and \(t_s>e^e\) such that the solution
of \eqref{eq:nls} with initial datum \(u_0\) satisfies
\begin{equation}\label{eq:two-sided}
 \exp\!\left(c_s\sqrt{\frac{\log t}{\log\log t}}\right)
 \le \|u(t)\|_{H^s}
 \le \exp\!\left(C_s\sqrt{\frac{\log t}{\log\log t}}\right)
 \qquad (t\ge t_s).
\end{equation}
\end{theorem}

In Section~5 of his survey of Hamiltonian PDE, Bourgain asked whether
a solution of the periodic equation with data in \(H^s\) can exhibit
unbounded Sobolev growth and how fast the growth can be
\cite[p.~43]{BourgainProblem}. For the square torus, Giuliani and
Guardia \cite[Question~1]{GiulianiGuardia} formulate the existence
question as follows. Is there a solution with
\(\limsup_{t\to+\infty}\|u(t)\|_{H^s}=\infty\) for some \(s>1\)?
Theorem~\ref{thm:main} answers this question for every \(1<s<\infty\)
by giving arbitrarily small data in \(H^s\) whose solutions have
\(\|u(t)\|_{H^s}\to\infty\) as \(t\to+\infty\).

On the square torus, Colliander, Keel, Staffilani, Takaoka and Tao
\cite{CKSTT} showed that a solution can start arbitrarily small in
$H^s$, $s>1$, and reach any prescribed size. The geometry of their
construction comes from the rectangles formed by four distinct Fourier
modes satisfying the cubic resonance relation. By arranging these
rectangles in generations, they made the output of one exchange available
for the next. Guardia and Kaloshin \cite{GKcascade} obtained a transfer
time polynomial in the desired amplification factor.\footnote{For large
finite growth with quintic or analytic nonlinearities on $\T^2$, see
\cite{HPquintic,GHP}, and for a convolution potential, see
\cite{GuardiaPotential}. On irrational two-tori, rational approximations
supply nearly resonant rectangles for the cubic equation
\cite{GiulianiGuardia} and for its perturbation by a convolution
potential \cite{GiulianiPotential}.} In both results the
initial datum may depend on the target. Colliander et al.~\cite[p.~43,
after (1.6)]{CKSTT} explicitly left open the question of unbounded
growth for a single datum. The limit in Theorem~\ref{thm:main}
requires control even between successive
transfers, when energy can flow back.

The beating solutions of Grébert and Villegas-Blas \cite{GVB} and Grébert
and Thomann \cite{GrebertThomann} exhibit this return explicitly. A mode cedes energy to its partner
and recovers it later in the same orbit.
One can also follow the creation of high modes from finitely supported
data by multiphase geometric optics \cite{CDS,CarlesFaou}. A different
picture emerges when the size of the torus tends to infinity. In the
weakly nonlinear large-box limit the sum over discrete resonances gives
rise to a continuous resonant equation \cite{FGH}, with its own stationary
states and Hamiltonian dynamics \cite{GHT}. An earlier statistical approach
to high frequencies is due to Kuksin \cite{KuksinTurbulence}. For a randomly
forced equation on a cube he obtained lower bounds on spatial derivatives
averaged over realizations and time as the dispersion and diffusion tend
to zero.

For one fixed orbit, a particularly close precedent is Hani's theorem
for a family of cubic equations on $\T^2$ that includes the resonant
system \cite{Hani}. On the waveguide $\mathbb R\times\T^2$, Hani,
Pausader, Tzvetkov and Visciglia carried unbounded growth from the
resonant system to the cubic NLS by modified scattering \cite{HPTV}.
Partial harmonic trapping leads to another setting in which the
long-time dynamics of small solutions are governed by a resonant
equation \cite{HaniThomann}. On the circle, the cubic Szegő equation
has unbounded Sobolev orbits accessible through inverse spectral theory
\cite{GGszego,GGhankel}. Its flow approximates small solutions of the
cubic half-wave equation over a long finite interval \cite{GGhalf}.
For the half-wave equation itself, the author obtained exponential growth for
$1/2<s<3/2$ \cite{ChenHalf}.

Fourier restriction estimates gave Bourgain a local theory for the periodic
cubic equation \cite{BourgainRestrict,BourgainGlobal}. The same method
also yields bilinear estimates for quadratic Schrödinger equations in
two dimensions \cite{StaffQuadratic}. On the square torus, estimates
based on almost conservation extended the global cubic theory below
the energy space \cite{Iteam}, and Herr and Kwak proved that the defocusing cubic equation
is globally well-posed in $H^s(\T^2)$ for every $s>0$
\cite{HerrKwakSharp,HerrKwak}. For more regular data one can ask how
large a higher Sobolev norm can become by time $t$. Bourgain and
Staffilani obtained polynomial upper bounds \cite{BourgainGrowth,StaffGrowth},
and subsequent work on flat tori \cite{CatoireWang} and compact
manifolds \cite{PTV} showed how the bounds depend on the geometry.
For example, Deng and Germain improved the polynomial $H^2$ bound on
generic irrational three-tori for powers $3<p<5$ \cite{DengGermain}.
If the norm does tend to infinity, how slowly can it do so? For every
$1<s<\infty$ and every nondecreasing $G(t)\to\infty$,
Theorem~\ref{thm:clock} gives a solution on the square torus whose
$H^s$ norm tends to infinity and is $o(G(t))$.
Staffilani \cite[\S2]{StaffWave} discusses a growth lower bound along
a sequence of times for the resonant system and asks what growth is
possible for the full equation. For \(s>2\), Theorem~\ref{thm:two-sided}
gives the two-sided rate \eqref{eq:two-sided} for the full equation at
every sufficiently large time. Both theorems describe rates attained by
particular solutions. Bourgain's conjecture \cite[p.~44, (5.10)]{BourgainProblem}
asks for an \(O(t^\varepsilon)\) upper bound for every solution and
every \(\varepsilon>0\).

The connection with the equation on the circle is particularly useful
on $\T^2$, since any one-dimensional solution remains an exact solution
when viewed as a function of $(x,y)$. In particular, the finite-gap tori of the
one-dimensional equation sit inside the two-dimensional phase space.
Small invariant tori for NLS have also been constructed by KAM methods
\cite{KP,EK,PP}. Long-time control of motion near small tori goes back
to Bambusi \cite{BambusiStability}, and the normal forms of Bambusi
and Grébert apply to a broad class of Hamiltonian PDE \cite{BambusiGrebert}. For the cubic equation, long-time stability is known near
plane waves \cite{FGL} and near suitable small finite-gap tori \cite{MP}.
Yet one-dimensional tori can also support large transverse excursions.
For $s>1$, Guardia, Hani, Haus, Maspero and Procesi \cite{GHHMP}
obtained any prescribed amplification in $H^s$ from data arbitrarily
close in $L^2$ to such a torus. Giuliani and Scandone
\cite[Theorem~1.3 and \S1.2]{GiulianiScandone} obtained finite
\(H^s\) amplification near a plane wave, with closeness measured in
Fourier \(\ell^1\), and asked whether a growing solution can start
arbitrarily close to a plane wave in \(H^s\). Corollary~\ref{cor:near-plane}
answers this NLS question for \(2\le s<\infty\). Every \(H^s\)
neighborhood of every plane wave contains a datum whose \(H^s\) norm
tends to infinity under the flow.

Here the one-dimensional solution supplies a Fourier harmonic whose
phase must agree with that of the transverse interaction. In global Birkhoff coordinates the actions of the
one-dimensional solution are constant, and each angle advances at a frequency determined by the actions \cite{GKbook}.
Grébert and Kappeler described the symmetries of this flow in action-angle
and spectral variables \cite{GKsym}. To use it in a two-dimensional
rectangle we need both its frequencies and its ordinary Fourier
coefficients. For the latter, the nonlinear part of the Birkhoff map
gains one derivative \cite{KST}. Quantitative estimates for the map and
for the Fourier coefficients of the flow are available in
\cite{MolnarEstimates,Maspero,KSTscatt}. For the frequencies we use the
formulas for the NLS hierarchy \cite{KM} and Molnar's local inversion
of the renormalized frequency map \cite{Molnar}.
Appendix~\ref{app:molnar-four} gives the corrected spectral tail estimate for \(p=4\) and the resulting compactness statement at each fixed small action.

\subsection{Strategy of the proof}

Suppose that, for each \(j\), the normalized Fourier modulus at
\((K_j,0)\) stays above \(1/2\) from some time \(t_j\) onwards.
Once the first \(m\) such times have passed, these frequencies alone
contribute at least \(m/4\) to the squared \(H^s\) norm. We build towards
this conclusion one frequency at a time. At each stage we perturb the
initial datum by an arbitrarily small amount in \(H^s\). The new mode
eventually contributes about one to the squared norm, while the change
at each of the finitely many earlier modes stays small for all future
times. Both the permanence of the contribution and the small cost at
the initial time come from the form of a single transfer.

The starting point is a one-dimensional solution, viewed on \(\T^2\)
as a field independent of \(y\). Multiplying it by \(e^{2iMt}\), where
\(M\) is its mass, puts it in the mass gauge. Proposition~\ref{prop:bank}
chooses the solution so that in this gauge it is periodic and the
coefficient of \(e^{ix-it}\) in its space-time Fourier series is nonzero.
Why this particular harmonic? Put \(e\coloneqq (1,0)\) and \(k_v\coloneqq e+Rv\)
for integer vectors \(v\) and a large integer \(R\). If \(p,q\) are
nonzero and perpendicular, then
\[
 k_p+k_q=k_{p+q}+e,\qquad
 |k_p|^2+|k_q|^2=|k_{p+q}|^2+1 .
\]
The cubic product of the waves at \(k_p,k_q\) with the conjugate of
the chosen harmonic has exactly the spatial index and the free time
frequency of the wave at \(k_{p+q}\). The harmonic's temporal frequency
is as important here as its spatial index. We obtain it by tuning the
actions of the one-dimensional solution.

We arrange these rectangles in a binary tree with root velocity
\((P,0)\), so the root has axial Fourier index \((K,0)\) with
\(K\coloneqq 1+RP\). If the two children carry fractions
\(\zeta,1-\zeta\) of the squared length of their parent velocity,
their combined homogeneous Sobolev weight is a fraction
\(\zeta^s+(1-\zeta)^s\) of the parent's weight. We keep \(\zeta\)
away from zero and one, so for \(s>1\) this fraction is bounded by a
constant \(c_s<1\) at every split. At depth \(d\) the sum of the
primary leaf weights is therefore at most \(c_s^d\) times the root
weight. The geometry
also produces auxiliary resonances. Section~\ref{sec:tree} includes
them in the finite system and chooses a trajectory on which the
transverse amplitudes decay and the root tends to one in model units.
The auxiliary amplitudes can be chosen small after the tree is fixed.
This is how a contribution of order one at late times can have
arbitrarily small initial cost.

To realize this transfer in the full equation, fix the tree and take
the physical amplitude \(r\coloneqq (K-1)^{-s}\), so that
\(r\langle K\rangle^s\to1\) as \(R\to\infty\). Two scale comparisons
enter Proposition~\ref{prop:insertion}. The new axial coordinates
give an \(H^1\) perturbation of order \(Rr\), which tends to zero. A transverse
frequency correction of order \(R^{-2}\) accumulates a phase of order
\(R^{-2}/r\) over the exchange time \(r^{-1}\). This too tends to zero,
since \(s<2\). A surviving axial mode must remain in phase on an
entire half-line, so we prescribe its exact frequency on an enlarged
periodic axial solution. In solving for the difference, we propagate
the rapidly decaying transverse directions forwards and prescribe
terminal conditions on the complementary directions. Bounds independent
of the terminal time give an exact solution approaching the enlarged
axial solution. We then follow it back to time zero.

After the first transfer the previous solution is already
two-dimensional. If \(U\) is that solution and \(V\) its periodic axial
limit, inserting a new tree at \(V\) accounts only for the limiting
background. Proposition~\ref{prop:history} restores the decaying
difference \(U-V\) by a correction to the newly constructed solution.
Section~\ref{sec:limit} can then choose the successive changes of
initial datum summably small in \(H^s\). Uniform control in a lower
Sobolev norm protects each previously chosen Fourier coefficient at
all later stages. Continuous dependence at each fixed time transfers
those bounds to the solution from the limiting datum, giving the
divergence in Theorem~\ref{thm:main}.

For Theorem~\ref{thm:clock} we also keep track of the norm during
each exchange. The new frequencies lie where the preceding finite-stage
orbit has a uniformly small Fourier tail, so cross terms in the
squared norm can be made summable. We space the \(K_j\) far enough apart
that the time windows for successive transfers are disjoint. At any
time, all but at most one transfer have either negligible contribution
or a contribution within a summable error of one. Counting the
thresholds \(K_j^s\) therefore gives \eqref{eq:count} with bounded
error even during an exchange. The counting threshold \(K_j^s\) lies
inside that window. The time by which the root stays above the
required height has a separate bound from the finite system. At fixed
depth this bound is a constant times \((K_j-1)^s\), so we can choose
\(K_j\) large enough to absorb the constant into the prescribed
\(\ell(K_j)\). Requiring also \(G(K_j^s)\ge j\) gives the arbitrarily
slow divergence in \eqref{eq:count}. The two corollaries in Section~\ref{sec:limit} use
the same construction to identify the limiting axial motion in lower
Sobolev norms and to allow the contribution of an individual root
to increase with depth.

Sections~\ref{sec:tree}--\ref{sec:limit} give the construction for
\(1<s<2\). If \(s\ge 2\), a normal frequency shift of order
\(R^{-2}\) accumulates on the transfer time of order \(R^s\).
Section~\ref{sec:h2-extension} computes the shift and assigns a
separate harmonic to each triangle. Section~\ref{sec:super} matches
the higher-order shifts needed for \(s>2\). At each new exponent
the construction selects an initial datum in the stated Sobolev space.

For Theorem~\ref{thm:two-sided} the construction must also control
the time between successive transfers. In Section~\ref{sec:two-sided}
two families of triangles with shared children allow the root height
to grow with the depth $d_j$, while the logarithmic frequency scale
$B_j$ advances by a fixed factor. Proposition~\ref{rate:prop:joint-insert}
bounds the cost of a transfer with this choice of scales. At any time
between transfers, one earlier root already supplies the lower bound.
Propositions~\ref{rate:p:flux} and~\ref{rate:p:back} bound the
high-frequency tails of the later transfers at that same time.

\smallskip
\noindent\emph{Conventions.}
We use \(\|\,\cdot\,\|_b\) for the \(H^b\) norm, and \(dx/(2\pi)\) for
measure on \(\T\coloneqq \R/(2\pi\Z)\). On the underlying real spaces we use the real part of the
usual inner product. Bars denote complex conjugates,
\(A^*\) denotes the adjoint for the unweighted \(L^2\) or
\(\ell^2\) inner product, and \([A,B]\coloneqq AB-BA\).
In estimates \(A\lesssim B\) and \(A=O(B)\) mean \(|A|\le C B\)
for \(B\ge0\). We write \(A\asymp B\) when both inequalities hold
for nonnegative \(A,B\), and \(A=o(B)\) when \(A/B\to0\) in the
indicated limit. A subscript indicates allowed dependence of the constant.
At a single insertion the tree and the old solution are fixed before
\(R\to\infty\). Section~\ref{sec:limit} specifies the order of choices
when the tree also varies.

For a constant \(c\in\R\), the field \(e^{-ict}u(t)\) solves
\eqref{eq:nls} with the additional potential \(c\). The choice
\(c=-2M(u)\) is the \emph{mass gauge}. In that gauge the same Fourier
moduli describe the solution, while a common temporal frequency is
removed from the equations for the axial modes.

\smallskip
\noindent\emph{Letters with section-dependent meaning.}
The notation is chosen locally, and several letters are used with
different meanings in different sections. Each is defined where it is
used, and within a subsection it has the meaning introduced there.
For orientation we list the main instances; \S\,$n$ refers to
Section~$n$.

{\small
\begin{center}
\begin{tabular}{@{}lp{0.88\textwidth}@{}}
\toprule
\(p\) & a child velocity in a split \(p+q=m\) (\S\S\,\ref{sec:tree},
\ref{sec:h2-extension}, \ref{sec:super}); the regularity index of the
scaled norms: \(p=(s+2)/2\) in \S\,\ref{sec:insertion},
\(p\in\{p_0,2\}\) in \S\,\ref{sec:h2-extension}, \(p=\lceil s\rceil\)
in \S\,\ref{sec:super}, \(p>2\lceil s\rceil+10\) in
\S\,\ref{sec:two-sided}; the index of a hierarchy flow and the
positive sites of a word in the path formulas; in
\S\,\ref{rate:sec:geom} also \(k_xl_x/h_x^2\) and the numerator of
\(t=p/q\). The balance constant of Proposition~\ref{prop:tree} and the
weak regularity index of \S\,\ref{sec:h2-extension} are both called
\(p_0\).\\
\(L\) & the scale \(R^\alpha\) of the residue weight
(\S\,\ref{sec:insertion}); the propagator \(L(t,v)\)
(\S\,\ref{sec:limit}); the height of a transfer
(\S\,\ref{sec:h2-extension}); the order of the expansion
\eqref{eq:sup-empty}, with \(\mathcal L\) a power of \(1+\log N\)
(\S\,\ref{sec:super}); in \S\,\ref{sec:two-sided}, successively,
\(N^{1/\sigma}\), a relative scale and a word cutoff, a Hessian bound,
a precision constant, \(\lambda^{-C_s}\)
(\S\,\ref{rate:sec:exterior}) and the right endpoint of a strip
(\S\,\ref{rate:sec:ward-landing}). In the proof of
Lemma~\ref{lem:empty-frequency}, \(\mathcal L\) is differentiation
along the angle flow.\\
\(T\) & the mixed vertex \(T_m\) (\S\S\,\ref{sec:tree},
\ref{sec:h2-extension}); the beginning of the terminal slow-time
interval and the matrix of approximate columns
(\S\,\ref{sec:insertion}; the column matrices and frames of
\S\S\,\ref{sec:h2-extension}--\ref{sec:two-sided} are also \(T\),
\(T_h\)); the period \(T_0\) (\S\,\ref{sec:limit}); the scale of the
five weak sites (\S\S\,\ref{sec:h2-extension}, \ref{sec:super}); in
\S\,\ref{sec:two-sided} also \(T_n=\tan\arg(Pu^n)\), subsets
\(T\subset\{1,\dots,m\}\), the times \(T_j\) and the time of the
stable section.\\
\(D\) & \(D=-i\partial_x\) and \(D_x\); divisors \(D(a,b,c)\), \(D\),
\(D_F\), \(D_A\); the diagonal operators \(D_{\rm free}\), \(D_h\),
\(D_0\); the times \(D_j=K_j^s\) (\S\S\,\ref{sec:limit}--\ref{sec:super});
\(D=D\chi\) in the proof in \S\,\ref{sec:h2-extension}; defects and
component pairs \(D(p,q)\), \(D_e\) (\S\S\,\ref{sec:super},
\ref{rate:sec:geom}); \(D=C_*d\log(d+2)\) in
Lemma~\ref{rate:l:finite}, later \(D_g\); the dressing matrix
\(D(\zeta)\); \(D=\mathcal Q/2\) in the calibration of
\S\,\ref{sec:two-sided}, which is the constant \(B\) of
Lemma~\ref{lem:empty-frequency}; the derivations
\(\mathcal D=\Omega\cdot\partial_\phi\) and
\(D=\Omega(I)\cdot\partial_\phi\); the historical precision \(D_h\).\\
\(h\) & an output vertex and the amplitudes \(h_j\) of the private
axes (\S\S\,\ref{sec:tree}, \ref{sec:h2-extension},
\ref{rate:sec:geom}); the transverse error
(\S\S\,\ref{sec:insertion}, \ref{sec:h2-extension}) and the field
variation (\S\,\ref{rate:sec:ward-landing}); the five-letter weak
word (\S\S\,\ref{sec:h2-extension}, \ref{sec:super}); the cutoff
\(|n|\ge hR\) (\S\,\ref{sec:h2-extension}); the order of the
reciprocal-moment template (\S\,\ref{rate:sec:axial-mark}); a moving
head \((h_x,RY)\) (\S\,\ref{rate:sec:exterior}); the exponent \(h_w\);
the sequence \(h_j=|j|^{-1}\) (Appendix~\ref{app:molnar-four}).\\
\(b\) & the radii \(b_j\) of the auxiliary spheres or circles
(\S\S\,\ref{sec:tree}, \ref{sec:h2-extension}); the spatial index of
an old harmonic in a key \((b,\nu)\); the axial error
(\S\,\ref{sec:insertion}) and the history corrections \(b_n\)
(\S\,\ref{sec:limit}); the gauge function of
Lemma~\ref{lem:empty-frequency} and the scalar phases \(b_2,b_3\); a
Sobolev index or a decay exponent; an endpoint of the weight interval
in Lemma~\ref{lem:common-endpoint-green}; in \S\,\ref{sec:two-sided}
the exponent in \(H=e^{bd}\), the radius \(e^{-c_1D}\), the weight
\(b(k)\) and integer coefficients \(b_e\).\\
\(J\) & the companion vertex \(J(v)\); the scaled old action
displacement (\S\,\ref{sec:insertion}; \(J_{\rm old}\), \(J_o\) and
\(j\) later); the generator \(J_{\rm off}\); an interval of weights
and a bound for defects (\S\,\ref{sec:h2-extension}); the orders \(J\)
and \(J_0\); the hot actions \(J_\infty\); the matrices \(J\),
\(\mathbb J\), \(\mathcal J\); period-one actions
(\S\,\ref{sec:tree}, Appendix~\ref{app:molnar-four}).\\
\(A,B\) & the constants \(A=\mathcal Q-M^2\), \(B=\mathcal Q/2\) of
Lemma~\ref{lem:empty-frequency}
(\S\S\,\ref{sec:tree}, \ref{sec:h2-extension}, \ref{sec:super}); in
\S\,\ref{sec:insertion} the operator \(A=D_x^2+m^2\) and the Wiener
norm \(\|\cdot\|_B\); in \S\,\ref{sec:two-sided} the scale
\(B=B_*Q^j\) with \(B_o=B/Q\), the cap exponents \(A,A'\) and the axis
vertex \(A(V)\).\\
\(E\) & \(E=n^2+m^2\); the weight \(e^{-\sigma(\tau-T)}\)
(\S\,\ref{sec:insertion}); the energy \(E(u)\) and the high-frequency
energy \(E_N\); the price classes \(E\), \(E_c\), \(E_o\),
\(\mathcal E\) of \S\,\ref{sec:two-sided}.\\
\(K\) & the root indices \(K=1+RP\), \(K_j\); the position of the
moderator sites (\S\S\,\ref{sec:h2-extension}, \ref{sec:super}) and
their dilation (\S\,\ref{sec:two-sided}); the Fourier cutoffs \(K_j\)
of Proposition~\ref{rate:p:back}; divided primitives.\\
\(\lambda\) & a decay rate of the tree orbit
(\S\,\ref{sec:insertion}); the modulus of the driving harmonics
(\S\S\,\ref{sec:h2-extension}, \ref{sec:super}), equal to
\(e^{-\epsilon B}\) in \S\,\ref{sec:two-sided}; a spectral parameter
(\S\,\ref{sec:limit}, Appendix~\ref{app:molnar-four}).\\
\(\tau,\theta,u\) & the slow time is \(\tau=rt\) in
\S\,\ref{sec:insertion}, \(\theta=r\lambda t\) in
\S\S\,\ref{sec:h2-extension}, \ref{sec:super}, and \(\tau=ut\) with
\(u=r\lambda\) in \S\,\ref{sec:two-sided}. In
\S\S\,\ref{sec:h2-extension}, \ref{sec:super}, \(\tau\) is also the
common factor of the frozen frequencies \(n^2+\tau c_n\). In
\S\,\ref{sec:two-sided}, \(\theta=n^2/E_n\) and \(\theta\) is an
angle displacement, and \(u=(1+it)^{-1}\) in
\S\,\ref{rate:sec:geom}.\\
\bottomrule
\end{tabular}
\end{center}
}

\section*{Acknowledgments and use of generative artificial intelligence.}
The author acknowledges the use of ChatGPT (5.6 sol, 6 Astra, 6.1 sol and 6.2 Astra) during the research process to explore mathematical ideas, discuss possible proof strategies, and assist with the drafting and editing of the manuscript. 

\section{An axial solution and a tree of right triangles}\label{sec:tree}

Write \(x,y\) for the two spatial variables. The subspace of functions of \(x\) alone
is invariant under the cubic equation. Our construction starts with an
orbit in this subspace that is periodic in the mass gauge. We shall use one Fourier coefficient of the orbit to couple
two short sides of a right triangle of frequencies to its hypotenuse. Both the
spatial and the temporal frequency of this coefficient matter.

\subsection{A periodic axial solution}

Write \(U(t,x)\) for an axial solution before making the mass gauge
transformation. In Birkhoff coordinates \(\xi\coloneqq \Phi(U)=(\xi_n)_{n\in\mathbb Z}\),
the action \(I_n\coloneqq |\xi_n|^2\) is constant and the frequency \(\omega_n(I)\)
depends only on the action vector \(I\coloneqq (I_n)_{n\in\mathbb Z}\).
Thus \(\xi_n(t)=\xi_n(0)e^{-i\omega_n(I)t}\).
The inverse map is denoted by \(\Psi\coloneqq \Phi^{-1}\), the differential of
\(\Phi\) at zero is the Fourier transform, and the mass in normalized
Haar measure is \(M=\sum_n I_n\). We subtract the common frequency \(2M\)
by multiplying the solution by \(e^{2iMt}\). In this gauge the frequency at \(n\)
is \(n^2+\chi_n(I)\), where
\begin{equation}\label{eq:renormalized-frequency}
   \chi_n(I)\coloneqq \omega_n(I)-n^2-2M.
\end{equation}
The same canonical chart diagonalizes the Hamiltonians of the NLS hierarchy
\cite{GKbook,KM}. Write \(T_{x_0}U\coloneqq U(\cdot+x_0)\). Translation by
\(x_0\) and multiplication by \(e^{i\vartheta}\) preserve every action and act
on the Birkhoff coordinates by
\begin{equation}\label{eq:birkhoff-covariance}
 \Phi(T_{x_0}U)_n=e^{inx_0}\,\Phi(U)_n,\qquad
 \Phi(e^{i\vartheta}U)_n=e^{i\vartheta}\,\Phi(U)_n,\qquad n\in\mathbb Z .
\end{equation}
Let \(G_a\) denote either \(T_a\) or multiplication by \(e^{ia}\), and set
\(S_a=\Phi\circ G_a\circ\Psi\). Both groups preserve the spectral actions
\cite[Theorem~1.2(i)]{GKsym}, hence also the reindexed and rescaled
actions \(I_j\) used here. Since \(S_a\) is symplectic, it commutes
with the Hamiltonian flow of each \(I_j\), and thus with rotation in
each coordinate separately.

Fix a finitely supported \(\xi\) and an index \(n\) with \(\xi_n\ne0\).
The coordinate rotations act transitively on the torus with moduli
\((|\xi_j|)_j\). Thus \(c_n(a,\xi)=(S_a\xi)_n/\xi_n\) has modulus one
and is constant on this torus for each \(a\). Since \(S_b\xi\) lies on
the same torus, the group law gives
\(c_n(a+b,\xi)=c_n(a,\xi)c_n(b,\xi)\). Hence \(a\mapsto c_n(a,\xi)\)
is a continuous character of \(\R/2\pi\Z\). Along the ray \(r\xi\),
\(0<r\le1\), these characters vary continuously, so their integer
winding number is constant. A circle character is determined by its
winding number, and therefore
\[
 c_n(a,\xi)=\lim_{r\downarrow0}\frac{(S_a(r\xi))_n}{r\xi_n}
           =\frac{(D S_a(0)\xi)_n}{\xi_n}.
\]
As \(D\Phi(0)=\mathcal F\) is the Fourier transform,
\(D S_a(0)=\mathcal F G_a\mathcal F^{-1}\). The last ratio is
\(e^{ina}\) for translation and \(e^{ia}\) for phase rotation.
Coordinates with \(\xi_n=0\) remain zero by action invariance.
Density of finitely supported sequences in \(\ell^2\) and continuity
of \(S_a\) now give \eqref{eq:birkhoff-covariance} on all of \(L^2\).
In the period-one convention of \cite[p.~7214]{KST}, the first component of
the Fourier differential sends an input \(v\) to \(-\widehat v(-j)\) at
spectral index \(j\). We use \(n=-j\) and multiply both coordinate
components by \(-1\), which preserves the actions and the canonical brackets.
The chart and its inverse are real analytic on \(L^2\) and \(\ell^2\).
For \(1\le q<\infty\), \(\ell^q\) denotes the usual sequence
space on the indicated indices, with counting measure, and
\(\ell^\infty\) has the supremum norm. Use \(h^b\) for the space of sequences with norm
\((\sum_{n\in\mathbb Z}(1+n^2)^b|\xi_n|^2)^{1/2}\).
The nonlinear part of each chart means its difference from its
differential at zero. These differences gain one derivative as analytic
maps from \(H^1\) to \(h^2\) and from \(h^1\) to \(H^2\), respectively
\cite[Theorem~1.1 and Corollary~1.1]{KST}.\footnote{In the notation
of the proof of \cite[Lemma~3.1(iii)]{KST}, the estimates on the
isolating neighborhood give \(\chi_n(\tau_n+t\gamma_n/2)-1=O(|n|^{-1})\)
for \(-1\le t\le1\). Thus the estimate on p.~7225 is
\(II=O(|n|^{-1})\). The constant term in \(\xi_n^2=I+II\) comes from
\(I=1+2t_n^2\).}
This also gives the chart regularity at every real order \(1\le b\le2\)
used below. For example, the nonlinear part of \(\Psi\) is analytic
from \(h^b\) to \(H^b\) by the continuous inclusions
\(h^b\subset h^1\) and \(H^2\subset H^b\). Its linear part is the
inverse Fourier transform. The same argument applies to \(\Phi\).
We call \(F\coloneqq \{n:I_n>0\}\) the Birkhoff support. For fixed actions with
finite support, the axial torus is the image under \(\Psi\) of all
choices of \(\xi_j=\sqrt{I_j}e^{i\varphi_j}\), \(j\in F\), with all
other coordinates zero. A subscript \(F\) on a sequence denotes its
restriction to these indices. Since \(\Phi\) identifies \(H^1\) with \(h^1\),
the energy domain of the axial NLS Hamiltonian corresponds to
\(\sum_n\langle n\rangle^2 I_n<\infty\). On a finite action face of
support \(F\), its differential along that face is
\(\sum_{j\in F}\omega_j\,dI_j\).

Another property we shall need is analyticity of \(I\mapsto\chi(I)\) on an
open \(\ell^2\) neighborhood of zero, with \(d\chi(0)=-\mathrm{Id}\)
\cite[Theorem~18.1 and Corollary~18.2, with \(p=4\)]{Molnar}; for
\(p=4\) the spectral tail estimate entering these statements, and the
compactness of \(D\chi+\mathrm{Id}\) at each small action, are supplied
in Appendix~\ref{app:molnar-four}. To check the
normalization, let \(v\) solve the equation of period one with cubic coefficient
two, as in these references. Then
\(U(t,x)=\sqrt{2}\,(2\pi)^{-1}v(t/(2\pi)^2,x/(2\pi))\) solves our equation.
For finitely supported period-one actions \(J\), the part of the Hamiltonian
quadratic in \(J\) is \(2(\sum_m J_m)^2-\sum_m J_m^2\).
Its \(J_n\)-derivative is \(4\sum_m J_m-2J_n\). Subtracting the common
frequency \(4\sum_m J_m\) leaves the linear term \(-2J_n\), as in
\cite[(4.3) and the proof of Corollary~18.2]{Molnar}.
The actions of \(v\) are \((2\pi)^2I_n/2\), after indexing by the physical
Fourier mode, and its frequencies are multiplied by \((2\pi)^{-2}\) when
passing to \(U\). Thus the derivative \(-2\mathrm{Id}\) of the renormalized
frequency map in the cited normalization becomes \(-\mathrm{Id}\).

\begin{proposition}\label{prop:bank}
There are smooth axial solutions \(U\) with arbitrarily small \(H^1\) norm
and finite Birkhoff support whose mass-gauged fields \(V(t,x)\coloneqq e^{2iMt}U(t,x)\)
have period \(2\pi Q\) for some
positive integer \(Q\) and satisfy
\begin{equation}\label{eq:bank}
  g\coloneqq \frac{1}{2\pi Q}\int_0^{2\pi Q}
               \int_{\mathbb T} e^{it-ix}V(t,x)\,\frac{dx}{2\pi}\,dt\ne0 .
\end{equation}
The initial solution can be chosen with Birkhoff support
\(F=\{0,4,5,2,6\}\). The Birkhoff coordinate at the spatial frequency \(1\)
vanishes identically on this torus.
\end{proposition}

\begin{proof}
Give \(F\) the order \(S\coloneqq (0,4,5,2,6)\) and put
\(\ell\coloneqq (1,1,1,-1,-1)\). Write \({\bf1}\coloneqq (1,1,1,1,1)\) and take
powers on \(S\) componentwise. Then
\[
 \ell\cdot{\bf1}=\ell\cdot S=\ell\cdot S^2=1,\qquad
 \ell\cdot S^3=-35 .
\]
We compute the coefficient of
\(\xi_0\xi_4\xi_5\overline{\xi_2\xi_6}\) in the spatial Fourier coefficient
at \(1\) of \(\Psi(\xi)\). Its powers of the five phases are encoded
by the character \(e^{i\ell\cdot\varphi}\), where
\(\varphi\coloneqq (\varphi_j)_{j\in S}\) in the indicated order. This is the
unique monomial of lowest degree with that character. For a cubic product with unconjugated inputs at \(a,c\) and
conjugated input at \(b\), put
\[
 k\coloneqq a-b+c,\qquad D(a,b,c)\coloneqq a^2-b^2+c^2-k^2 .
\]
When \(D\ne0\), comparison of the cubic terms in the physical equation
and its diagonal Birkhoff equation gives \(1/D(a,b,c)\) per ordered triple
in the inverse map. Indeed, temporal differentiation acts on the monomial
at linear order by its signed input energy, whereas the output Laplacian
has eigenvalue \(k^2\). The resonant cubic monomials contain an action
pair, since \(D=-2(a-b)(c-b)\). Such a pair cannot occur in our quintic
monomial with five distinct indices.

At fifth order the divisor for the NLS flow vanishes. We compute the same
coefficient with the commuting flow
\(U_\tau=U_{xxx}-3|U|^2U_x\), obtained from \cite[(6)]{KM} by the
same spatial and amplitude scaling, with the third time of the
period-one equation set to \(-\tau/(2\pi)^3\). The linear term in the equation for \(i\partial_\tau U\)
acts on \(e^{inx}\) by multiplication by \(n^3\). Symmetrizing the
cubic term in its two unconjugated inputs \(A,C\) gives the
multiplier \(3(A+C)/2\) per ordered triple. In the outer cubic
\(U_A\bar U_B U_C\), where \(U_A\) denotes a spatial Fourier
coefficient, insert the cubic part of \(\Psi\) in one slot. The other
two slots use its linear part. The inner triple \((a,b,c)\) supplies
the index \(a-b+c\) for the chosen slot. Inserting it into the
conjugated slot reverses all three signs. Summing over the assignments
of the three positive indices \(0,4,5\) and the two negative indices
\(2,6\) gives
\[
\begin{array}{c|rrr}
 \text{slot of the inner triple} &1&2&3\\ \hline
 \sum 1/D(a,b,c)&-1/2&1&-1/2\\
 \sum 3(A+C)/(2D(a,b,c))&-3/4&11/2&-3/4
\end{array}
\]
For example, in the middle slot the two negative indices \(2,6\) play
the roles of \(a,c\), and \(b\) ranges over \(0,4,5\). The sum of the
outer positive indices is \(A+C=9-b\). The reciprocals \(1/D\) for
the three choices of \(b\) are \(-1/24,1/8,1/6\). Each occurs four
times, since \(a,c\) and \(A,C\) can be permuted independently. This
gives \(1\) in the first row and
\(4(-9/16+15/16+1)=11/2\) in the second.
In the first slot the inner conjugated index is \(b=2\) or \(6\).
For these choices the sums of \(1/D\) over the positive assignments
are \(1/4\) and \(-3/4\), respectively, with \(A+C=9-b\).
The third slot gives the same sums by symmetry.
The forcing for the NLS flow therefore vanishes, whereas the third-flow
forcing equals four.

\smallskip
\noindent\emph{The homological equations.} Write
\(\Psi(\xi)=\sum_n\xi_ne^{inx}+\Psi^{(3)}(\xi)+\Psi^{(5)}(\xi)
+O(\|\xi\|^7)\), with \(\Psi^{(k)}\) homogeneous of degree \(k\) in
\((\xi,\bar\xi)\). Consider a flow of the hierarchy,
\(i\partial_tU=\Lambda U+N_3(U)\) with \(\Lambda e^{inx}=\lambda_ne^{inx}\)
and \(N_3\) cubic, whose Birkhoff form is \(i\dot\xi_n=\Omega_n(I)\xi_n\)
with \(\Omega_n(I)=\lambda_n+O(I)\). Insert \(U=\Psi(\xi(t))\) and use
\(i\frac{d}{dt}(\xi_a\bar\xi_b\xi_c)
=(\Omega_a-\Omega_b+\Omega_c)\,\xi_a\bar\xi_b\xi_c\). For a monomial
without an action pair, the coefficient \(\psi_{abc}\) of
\(\xi_a\bar\xi_b\xi_c\,e^{ikx}\) in \(\Psi^{(3)}\), \(k=a-b+c\), then
satisfies
\[
 (\lambda_a-\lambda_b+\lambda_c-\lambda_k)\,\psi_{abc}=\nu_{abc},
\]
where \(\nu_{abc}\) is the coefficient of \(N_3\) on this ordered triple.
The \(O(I)\) parts of \(\Omega\) contribute only to monomials containing
an action pair. For NLS, \(\lambda_n=n^2\) and \(\nu_{abc}=1\), which is
the rule \(1/D(a,b,c)\). At degree five, the coefficient \(\psi_5\) of
\(\xi_0\xi_4\xi_5\overline{\xi_2\xi_6}\,e^{ix}\) satisfies
\[
 \Bigl(\sum_j\ell_j\lambda_{S_j}-\lambda_1\Bigr)\psi_5
 =\bigl[N_3\text{ with one slot replaced by }\Psi^{(3)}\bigr]_5 ,
\]
the right side being the coefficient of this monomial. No other term
contributes: both flows used here have no quintic part, and every term
produced by the \(O(I)\) frequency corrections contains an action pair.
For NLS the divisor is \(\ell\cdot S^2-1=0\), and the first row of the
table confirms that the right side vanishes as well. For the third flow
\(\lambda_n=n^3\), the symmetrized coefficient is \(\nu_{ABC}=3(A+C)/2\),
the right side is the sum \(-3/4+11/2-3/4=4\) of the second row, and the
divisor is \(\ell\cdot S^3-1=-36\). Since the hierarchy shares the chart
\(\Psi\), the two flows determine the same coefficient \(\psi_5\).
Corrections to the linear angle velocities carry at least one action
pair and have zero coefficient at this monomial. Its third-flow divisor
is \(\ell\cdot S^3-1^3=-36\), so its coefficient is \(-1/9\).
Writing \(\xi_j=\sqrt{I_j}e^{i\varphi_j}\), analyticity shows that the
full angular coefficient equals
\[
 \left(-\tfrac19+O\Big(\sum_{j\in F}I_j\Big)\right)
             \prod_{j\in F}\sqrt{I_j}.
\]
Every higher monomial with the same angle character is obtained from the
lowest monomial by adding action pairs.

We now impose \(\ell\cdot\chi_F=0\) so that this harmonic has time
frequency \(1\). The derivative of the left side at the origin is
\(-\ell\). Set \(I_0=I_4=I_5=I_2=\epsilon\). The derivative with respect
to \(I_6\) is \(1\), so the implicit function theorem gives
\(I_6=2\epsilon+O(\epsilon^2)\). All five actions are positive for
small \(\epsilon>0\). On a relatively open patch about this point the
angular coefficient computed above stays nonzero. The local coordinate
map \(I_F\mapsto\chi_F\) sends the patch to an open subset of
\(\ell\cdot\chi_F=0\). This hyperplane has rational coefficients, so
we may fix a point in the patch with rational \(\chi_F\).

It remains to choose an orbit on this torus. Let \(a_m\), \(m\in\Z^5\),
be the angular Fourier coefficients of the spatial Fourier mode at
\(1\) on the torus with the chosen actions. We have proved that
\(a_\ell\ne0\). From initial angles \(\varphi\), the mass-gauged
angles at time \(t\) are \(\varphi-(S^2+\chi_F)t\). Inserting their
Fourier series into \eqref{eq:bank} and averaging over the common
period \(2\pi Q\) gives
\[
 g(\varphi)=
    \sum_{m\cdot(S^2+\chi_F)=1}a_m e^{im\cdot\varphi}.
\]
Here \(Q\) is a common denominator of the frequencies, and the series
is absolutely convergent by analyticity in the angles. Several
characters can survive this time average. Since \(\ell\) is among
them, orthogonality in the initial angles gives
\(\int_{\T^5}|g(\varphi)|^2\,d\varphi/(2\pi)^5\ge |a_\ell|^2>0\).
This supplies a single periodic orbit with \(g\ne0\).
\end{proof}

At each subsequent insertion we shall add finitely many axial
coordinates and keep the angle velocities already used. This is why we
tuned the frequencies through the actions in the first construction.

\begin{lemma}\label{lem:locking}
Fix a torus of support \(F\) in a sufficiently small \(\ell^2\) action
ball, and let \(G\) be a disjoint finite set. Sufficiently small positive
actions \(I_G\) can be added while preserving \(\chi_F\) by a change of
size at most \(\|I_G\|_{\ell^2}/3\) in the old actions. By allowing the
new positive radii to vary in arbitrarily small open neighborhoods of
their prescribed values, one may in addition take all new
\(\chi_j\), \(j\in G\), rational.
\end{lemma}
\begin{proof}
Take the action ball so small that \(\|d\chi+\mathrm{Id}\|\le1/4\).
The compression of \(d\chi\) to \(F\) has an inverse of norm at most \(4/3\),
and the derivative from \(G\) to \(F\) has norm at most \(1/4\).
We need the allowed size of \(I_G\) to be independent of the locations
of the new indices. Identify both coordinate subspaces with their zero
extensions and put \(x_0\coloneqq I_F^{\rm old}\), \(y\coloneqq I_G\). For fixed \(y\),
the equation for the new old actions is the fixed point equation
\[
                 x=x+\chi_F(x+y)-\chi_F(x_0).
\]
The right-hand side has Lipschitz constant at most \(1/4\) as a
function of \(x\). At \(x=x_0\) it differs from \(x_0\) by at most
\(\|y\|_{\ell^2}/4\), since \(y\) is supported off \(F\).
It thus maps the closed ball of radius \(\|y\|_{\ell^2}/3\) about
\(x_0\) into itself. For small enough \(\|y\|_{\ell^2}\) this ball
keeps all old actions positive and \(x+y\) inside the action chart.
Both restrictions depend only on the old point and the chart
neighborhood. The contraction theorem gives the required old actions.
The full derivative on \(F\cup G\) is likewise within \(1/4\) of
minus the identity and hence invertible. In its local frequency
coordinates, fixing \(\chi_F\) therefore leaves an open set of
\(\chi_G\) values. This remains true in any open neighborhood where
the prescribed new radii are positive. The new \(\chi_G\) may therefore be taken rational,
with the new positive radii in any prescribed sufficiently small open
neighborhood of their nominal values.
\end{proof}

Keep the initial old angles fixed.
Their subsequent trajectories in the two mass gauges agree in every
old angular coordinate. By continuity of \(\Psi\) on the compact old
angle torus, the change of the axial field is uniformly small in \(H^1\)
whenever the new coordinates are sufficiently small in \(h^1\) and the
old radii change sufficiently little. In the periodic construction of
Sections~\ref{sec:insertion} and~\ref{sec:limit}, the old occupied frequencies
are rational. Choosing the new frequencies rational as well gives the two
mass-gauged fields a common period.
Averaging over this common period shows that their coefficients in \eqref{eq:bank} differ by at
most their uniform \(L^2\) distance.

We also need the unoccupied frequencies of a fixed
smooth axial torus. Here a bound at an individual mode will be more
useful than an asymptotic in a sequence norm.

\begin{lemma}\label{lem:empty-frequency}
For a smooth axial torus of finite Birkhoff support let
\[
 \mathcal{H}=\int_{\T}(|U_x|^2+\tfrac12|U|^4)\,\frac{dx}{2\pi},
 \qquad \mathcal{Q}=2\left(\sum_j I_j\omega_j-\mathcal{H}\right),
 \qquad A=\mathcal{Q}-M^2,\quad B=\mathcal{Q}/2 .
\]
For empty indices \(n\), as \(|n|\to\infty\),
\[
                 \chi_n(I)=\frac{A-B}{n^2}+O(|n|^{-3}).
\]
More generally the complexified transverse linearization at vertical
frequency \(m\) has an approximate column with carrier energy
\(E+2M+A/n^2-B/E\), where \(E=n^2+m^2\). After extracting the
carrier its two envelopes are \(1+O(|n|^{-1})\) and \(O(E^{-1})\).
Its residual in \(L^2_x\) is \(O(|n|^{-3})\), uniformly in \(m\).
The constants are locally uniform on compact families with fixed
finite support.
\end{lemma}

\begin{proof}
We work in the ungauged variables, write \(\langle\cdot\rangle\)
for spatial average, \(D=-i\partial_x\) and \(\mathcal L\) for
differentiation along the angle flow. The difference of the physical
Liouville form \(-\operatorname{Im}\langle\bar U\,dU\rangle\) and
the pullback of \(-\operatorname{Im}\sum_j\bar\xi_j\,d\xi_j\) is exact
on the Birkhoff chart about zero, by canonicity. For a general smooth
torus the same conclusion follows on the star-shaped preimage in
Birkhoff coordinates, by radial integration. Evaluating a primitive
\(G\) along the flow gives
\[
 \mathcal LG
 =\langle |U_x|^2+|U|^4\rangle-\sum_jI_j\omega_j,\qquad
 \langle|U|^4\rangle=\mathcal{Q}+2\mathcal LG.
\]
This identity holds on the whole torus, including resonant
trajectories. Put \(V_0=2(|U|^2-M)\) and \(C=U^2\). After removing
\(e^{inx+imy-i(E+2M)t}\) from both components of the doubled
linearization, prescribe a common scalar \(a\). The two residuals,
with the second equation multiplied by minus one, are
\[
 \begin{gathered}(2nD+K_+-a)f+Cg,\qquad
 (2E+2nD+K_-+a)g+\bar C f,\\
 \quad K_+=D^2+V_0-i\mathcal L,\quad K_-=D^2+4M+V_0+i\mathcal L .\end{gathered}
\]
Define \(f_1=-D^{-1}V_0/2,\ f_2=-D^{-1}K_+f_1/2\), with the
inverse zero on constants. Both arguments have mean zero since
\(\langle V_0D^{-1}V_0\rangle=0\). Set
\[
 L_2=K_+ f_2-\frac{n^2}{2E}|C|^2,\qquad
 f_3=-\tfrac12D^{-1}(L_2-\langle L_2\rangle),\qquad
 a=\langle L_2\rangle/n^2,\qquad f=1+\sum_{j=1}^3 n^{-j}f_j .
\]
For the lower row put \(g_0=-\bar C f/(2E)\),
\(S=(2nD+K_-+a)/(2E)\) and \(g=(1-S+S^2)g_0\). Its residual
is \(2E S^3g_0=O(|n|^{-3})\), since \(E\ge n^2\) and only
finitely many derivatives of the smooth torus occur. In the upper
row \(g=-\bar C/(2E)+O(|n|^{-3})\), and the definitions of
\(f_1,f_2,f_3,a\) cancel the terms through order \(n^{-2}\).
All these estimates hold with any fixed number of envelope derivatives.
Explicitly, put \(\theta=n^2/E\in(0,1]\). Since \(D1=0\), \(K_+1=V_0\),
\(a=O(n^{-2})\) and \(Cg=-\theta|C|^2/(2n^2)+O(|n|^{-3})\), the upper residual is
\[
 (2nD+K_+-a)f+Cg=(2Df_1+V_0)+n^{-1}(2Df_2+K_+f_1)
   +n^{-2}(2Df_3+L_2-n^2a)+O(|n|^{-3}) .
\]
The first bracket vanishes because \(\langle V_0\rangle=0\) (\(M\) is the spatial
mean of \(|U|^2\)) and \(f_1=-D^{-1}V_0/2\). The second vanishes because
\(\langle K_+f_1\rangle=\langle V_0f_1\rangle=-\tfrac12\langle V_0D^{-1}V_0\rangle=0\),
the Fourier terms at \(l\) and \(-l\) cancelling, so that
\(f_2=-D^{-1}K_+f_1/2\) solves \(2Df_2=-K_+f_1\). The third vanishes by the
definitions of \(f_3\) and \(a\). In the lower row,
\(2E(1+S)(1-S+S^2)g_0=2E(1+S^3)g_0\) and \(2Eg_0=-\bar Cf\), so the lower
residual is exactly \(2ES^3g_0\).

For the scalar, bilinear integration by parts gives
\[
 \langle V_0f_2\rangle
 =-\langle f_1K_+ f_1\rangle
 =\tfrac14\langle V_0^2\rangle
                  +\tfrac i2\mathcal L\langle f_1^2\rangle .
\]
Indeed \(Df_1=-V_0/2\) and \(\langle (Df_1)f_1^2\rangle=0\).
Since \(f_2\) has spatial mean zero, the mean of \(K_+f_2\) is
\(\langle V_0 f_2\rangle\). Now \(V_0=2(|U|^2-M)\) and
\(|C|^2=|U|^4\), so the formula for \(L_2\) gives
\[
 a=\left(\frac1{n^2}-\frac1{2E}\right)\langle |U|^4\rangle
       -\frac{M^2}{n^2}
       +\frac{i}{2n^2}\mathcal L\langle f_1^2\rangle .
\]
Substituting \(\langle |U|^4\rangle=\mathcal Q+2\mathcal LG\) gives
\[
 a=\frac A{n^2}-\frac B E+\mathcal L b,\qquad
 b=\left(\frac2{n^2}-\frac1E\right)G+
                    \frac{i}{2n^2}\langle f_1^2\rangle .
\]
Multiplication of both envelopes by \(e^{-ib}\) replaces \(a\)
by \(A/n^2-B/E\), with the same residual bounds;
\(b=O(|n|^{-2})\) on the fixed torus.
Here \(b\) is independent of \(x\) and may be complex. The primitive \(G\) is real,
and \(f_1=-D^{-1}V_0/2\) is purely imaginary because \(V_0\) is real,
so that \(\langle f_1^2\rangle\le0\) and
\(\operatorname{Im}b=\langle f_1^2\rangle/(2n^2)\). Accordingly the
scalar \(a\) of the preceding display has the imaginary part
\(\mathcal L\langle f_1^2\rangle/(2n^2)\), and the factor \(e^{-ib}\)
has modulus \(\exp(\langle f_1^2\rangle/(2n^2))=1+O(n^{-2})\), bounded
together with its inverse and its derivatives along the flow. Write
\(R(a;f,g)\) for the pair of residuals above. Since \(e^{-ib}\) commutes
with \(D\) and
\(\mp i\mathcal L(e^{-ib}h)=e^{-ib}(\mp i\mathcal Lh\mp(\mathcal Lb)h)\),
\[
 R(a-\mathcal Lb;e^{-ib}f,e^{-ib}g)=e^{-ib}R(a;f,g).
\]
The upper envelope is
still \(1+O(|n|^{-1})\), and the new scalar \(A/n^2-B/E\) is real.

For \(m=0\) apply the complexified differential of the Birkhoff
map to this column. At an empty index the variation of \(\omega_n\)
does not enter its equation. After factoring out the carrier time
phase, the positive \(n\)th coefficient \(c_n\) is
\(1+O(|n|^{-1})\): the input has \(H^1\) norm \(O(|n|)\), and the
one-smoothing part of the differential gives an \(O(|n|)\) bound in
\(h^2\), uniformly in the old angles. The transformed residual is
\(O(|n|^{-3})\) in \(\ell^2\) by the \(L^2\) chart bound. Thus
\[
 i\mathcal L c_n-\bigl(\chi_n-(A-B)/n^2\bigr)c_n=O(|n|^{-3}).
\]
Integrate over the full torus of old angles. The integral of
\(\mathcal L c_n\) vanishes because \(\mathcal L\) is a constant
linear combination of angle derivatives. The factor
\(\chi_n-(A-B)/n^2\) is constant on the torus and the average
of \(c_n\) is \(1+O(|n|^{-1})\). Dividing by this average proves
the empty-frequency expansion, including for a resonant angle flow.
\end{proof}

\subsection{The resonance geometry}

For a velocity \(v\in\mathbb Z^2\) the spatial frequency of a new wave
will be \(k_v\coloneqq (1,0)+Rv\), with \(R\) a positive integer. If \(p+q=h\), then
\[
          |k_p|^2+|k_q|^2-|k_h|^2=1-2R^2p\cdot q.
\]
Thus the coefficient \eqref{eq:bank} couples the three waves exactly at
the level of the free dispersion when \(p\cdot q=0\). We call
\((p,q,h)\) a triangle in this case and regard its first two entries as
unordered. A pair \(p,h\) with \(p\cdot(h-p)=0\) can also produce
\(q=h-p\). We shall keep track of both operations when choosing the
velocities.

\begin{proposition}\label{prop:tree}
For every integer \(d\ge1\) there is a finite set of integer velocities
containing a full binary tree of depth \(d\), with two children at every
internal vertex and \(d\) edges from the root to each leaf. Its root \((P,0)\) has \(P>0\), and the
remaining primary vertices (the vertices of the binary tree) are off
the horizontal axis. Each parent is the sum of its two orthogonal
children. There is an absolute \(p_0>0\) such that the two squared
child lengths lie between \(p_0\) and \(1-p_0\) times the squared parent
length. All horizontal components are nonzero and have distinct absolute
values. The following properties hold.
\begin{enumerate}
\item At a split \(p+q=m\) with \(m\) off the axis, write \(v=(v_x,v_y)\)
for either child and include the points
\[
 J(v)\coloneqq (v_y^2/v_x,-v_y),\quad A(v)\coloneqq v+J(v)\quad (v=p,q),
 \qquad T_m\coloneqq J(p)+J(q).
\]
The four pairs of legs at this split are \((p,q)\), \((p,J(p))\),
\((q,J(q))\), and \((J(p),J(q))\), with respective hypotenuses
\(m\), \(A(p)\), \(A(q)\), and \(T_m\). The root split consists of a single triangle.
These are all nondegenerate forward and reverse quadratic resonances
among pairs of the selected velocities, including those whose output
is initially sought outside the set.
\item If \(P_i,P_j\) are distinct horizontal coordinates of axial
vertices, then \(P_i(P_j-P_i)\) is never a positive rational square.
\end{enumerate}
\end{proposition}

\begin{proof}
We construct the tree over the rationals, starting with root \((1,0)\).
For an off-axis parent \(m=(X,Y)\), take a rational parameter \(t\ne0\)
and put
\[
 p\coloneqq \frac{X+tY}{1+t^2}(1,t),\qquad
 q\coloneqq \frac{Xt-Y}{1+t^2}(t,-1).
\]
These are the two orthogonal projections of \(m\). The vectors
\(p,J(q)\) lie on one of the perpendicular lines and \(q,J(p)\) on
the other. Their four possible sums are precisely the prescribed
hypotenuses. In particular, \(J(p)\) is the unique partner of \(p\)
that is orthogonal to \(p\) and has vertical component \(-p_y\).
Thus both \(A(p)\) and \(A(q)\) lie on the axis. Explicitly, with
\(U\coloneqq (X+tY)/(1+t^2)\) and \(V\coloneqq (Xt-Y)/(1+t^2)\), so that
\(p=U(1,t)\) and \(q=V(t,-1)\),
\[
\begin{gathered}
 J(p)=tU\,(t,-1),\qquad J(q)=\tfrac Vt\,(1,t),\\
 A(p)=(X+tY,0),\qquad A(q)=(X-Y/t,0),\qquad T_m=\bigl(t^2U+\tfrac Vt,\,-Y\bigr).
\end{gathered}
\]
The vertical component of \(T_m\) is
\(-tU+V=-Y(1+t^2)/(1+t^2)=-Y\).
At the horizontal root \(Y=0\) the projection formulas give \(J(p)=q\)
and \(J(q)=p\), so only \(p,q\) are added.

We give the exclusions used at each step. For a previously chosen
\(u\) and a moving vertex \(v\), a pair can generate a nondegenerate
source only on one of the three loci
\begin{equation}\label{eq:forbidden-loci}
       v\cdot u=0,\qquad v\cdot u=|u|^2,\qquad |v|^2=v\cdot u .
\end{equation}
The first two are lines perpendicular to \(u\), and the third is the circle
with diameter \([0,u]\). Each of \(p,q\) traverses the circle with diameter
\([0,m]\), so the only persistent relation of this kind for an old point
is the prescribed one with \(u=m\). Each of \(J(p),J(q)\) traces an
unbounded curve contained in no line. For example, writing
\(U\coloneqq (X+tY)/(1+t^2)\), one has \(J(p)=(t^2U,-tU)\). As \(t\to\infty\)
the first coordinate is unbounded and the second has a finite limit and
is nonconstant. The calculation for \(J(q)\) uses \(t\to0\).
The remaining three new points traverse horizontal lines, and \(u_x\ne0\)
rules out either line in \eqref{eq:forbidden-loci}. Consequently each
unwanted old--new condition eliminates finitely many values of \(t\).
Indeed, each new point is a rational function of \(t\), and each locus in
\eqref{eq:forbidden-loci} is the zero set of a polynomial of degree at
most two in \(v\). Along a new point the condition becomes a rational
function of \(t\), which either has finitely many zeros or vanishes
identically, in which case the whole curve lies in the locus. The
curve traced by \(p\) or \(q\) is the circle with diameter \([0,m]\) (up to
finitely many points); it lies in no line, and in the circle with
diameter \([0,u]\) only if \(u=m\). The curves traced by \(J(p),J(q)\) are
unbounded, hence lie in no circle, and lie in no line, as shown above.
The points \(A(p),A(q),T_m\) move nonconstantly on horizontal lines; a
horizontal line lies in no circle, and in a line \(v\cdot u=\mathrm{const}\)
only if \(u_x=0\), which is excluded.

For pairs within a single split the equations, after clearing
denominators and scaling by \(X\), are polynomials in \(Y/X,t\).
The following one rational evaluation tests all of them simultaneously.
For \(m=(1,2)\), \(t=3\), multiply the ordered list
\((m,p,q,J(p),J(q),A(p),A(q),T_m)\) by 30 to obtain
\[
 (30,60),\ (21,63),\ (9,-3),\ (189,-63),\
 (1,3),\ (210,0),\ (10,0),\ (190,-60).
\]
Numbering from zero, the pairs with vanishing dot product are
\((1,2),(1,3),(2,4),(3,4)\). The ordered pairs \((i,j)\), \(i\ne j\),
with \(v_i\cdot(v_j-v_i)=0\) are
\((1,0),(2,0),(1,5),(3,5),(2,6),(4,6),(3,7),(4,7)\).
All horizontal components in this test have different absolute values.
For each unwanted relation, write its cleared polynomial in powers of
\(t\), with coefficients polynomial in the parent slope \(Y/X\).
The test just given makes at least one coefficient a nonzero polynomial.
Its finitely many roots contain all the slopes for which the relation
could hold for every \(t\). Taking the union over the finite list of
unwanted relations still excludes only finitely many slopes. At the birth of every primary
vertex avoid these exceptional slopes, as well as vanishing coordinates
and coincident or opposite horizontal components. The horizontal
components of the new points are nonconstant rational functions of \(t\),
so this also avoids coincidences with older vertices. The two child
slopes are \(t\) and \(-1/t\), which lets us impose in advance the
slope exclusions needed when each child later becomes a parent.
At the first split
the children move on the root circle, and the same finite exclusions
on their slopes and coordinates apply.

To ensure the arithmetical condition on axis vertices, note that the
two new horizontal outputs have abscissae
\[
              U_1\coloneqq X+tY,\qquad U_2\coloneqq X-Y/t.
\]
For a prime \(\pi\) and a nonzero rational number \(w\), write
\(\operatorname{ord}_\pi(w)\) for the exponent of \(\pi\) in its
numerator minus that in its denominator. Rational squares have even
valuation at every prime. For each old axial abscissa \(C\ne0\),
choose a fresh odd prime \(\pi\) for which \(Y,C\) are units (valuation
zero) and all fixed coordinates are integral. Require
\(t\equiv(C-X)/Y+\pi\pmod{\pi^2}\), with rational denominators read
modulo \(\pi^2\). Then \(\operatorname{ord}_\pi(U_1-C)=1\), while
\(U_1\) and \(C\) are units. Both \(C(U_1-C)\) and \(U_1(C-U_1)\)
therefore have odd valuation. At another fresh prime use the same
argument for \(U_2=X-Y/t\), with \(1/t\) as parameter. Here \(C\ne X\)
by the distinct horizontal components already imposed. We may therefore
choose this prime so that the prescribed residue of \(1/t\) is a unit
and invert it to obtain a congruence on \(t\). For \(U_1,U_2\) themselves
take a fresh prime congruent to \(1\bmod 4\) avoiding \(Y(X^2+Y^2)\).
A root of \(t^2+1\) modulo this prime is simple, so all but one of
its lifts modulo the prime squared give valuation one. Both \(U_i\) are units, and
\(U_2-U_1=-Y(t^2+1)/t\) has valuation one. These finitely many conditions
are compatible with prescribing the direction of \(p\) in any real
open interval. By the Chinese remainder theorem they prescribe a
residue \(t_0\) modulo a product \(M_0\) of prime squares. For an
integer \(b\equiv1\bmod M_0\), all rationals \(a/b\) with
\(a\equiv t_0\bmod M_0\) satisfy the congruences. Their spacing on
the real line is \(M_0/b\), so any prescribed open interval contains
such a parameter when \(b\) is sufficiently large.

Choose the angle between the parent and the line of \(p\) in a fixed
small interval about \(\pi/4\). This gives a common lower bound
\(p_0>0\) for both relative squared lengths by the cosine and sine
formulas. Removal of the finitely many forbidden directions leaves an
open subinterval. Continuing to depth \(d\) proves the construction over
the rationals. Multiplication of every velocity by a common denominator
preserves all incidence conditions and multiplies each expression
\(P_i(P_j-P_i)\) by a square.
\end{proof}

\subsection{An orbit of the finite system}

Fix the old coefficient \(g\) of \eqref{eq:bank}. Write the new wave
at \(k_v\), after removal of its free oscillation, as \(r a_v\),
where \(r>0\) is the common amplitude scale and \(a_v\) is the complex
model amplitude. For the quadratic terms induced by the old wave,
the two factors \(r\) on the right and the one on the left leave the
time scale \(r^{-1}\). With a dot denoting differentiation with respect
to \(rt\), the contribution to \(i\dot a_h\) of
\(p+q=h\) is \(2\bar g\,a_pa_q\). The contributions to \(i\dot a_p\)
and \(i\dot a_q\) are \(2g\,a_h\bar a_q\) and \(2g\,a_h\bar a_p\).
The factor two counts the two orders of the positive inputs at the
output \(h\), and is the coefficient of the mixed density term at the
other outputs. Writing \(a_v=i(g/|g|)X_v\), the real amplitudes \(X_v\)
satisfy in the model time \(z\coloneqq 2|g|rt\)
\begin{equation}\label{eq:single-triangle}
  X_h'=X_pX_q,\qquad X_p'=-X_qX_h,\qquad X_q'=-X_pX_h .
\end{equation}
At a vertex in several triangles, all of the displayed contributions are
added. From here to the end of this section a prime denotes
differentiation with respect to \(z\).

A primary vertex at height \(j\) above the leaves has the same
incidences as any other primary vertex at this height. We give these
vertices a common amplitude \(x_j\), and denote the root by \(H\coloneqq x_d\).
For \(0\le j\le d-2\), let \(w_j,h_j,t_j\) denote the common amplitudes
of the vertices \(J(v), A(v), T_m\), respectively, for the split in
which \(v\) is a child of height \(j\). Thus \(t_j\) is one real
amplitude per sibling pair. Equation \eqref{eq:single-triangle}
preserves this symmetry. In these variables it reads
\begin{align}
 x_j'&=x_{j-1}^2-x_jx_{j+1}-w_jh_j\quad(0\le j<d),
       & H'&=x_{d-1}^2, \label{eq:tree-primary}\\
 h_j'&=x_jw_j,\qquad t_j'=w_j^2,
       & w_j'&=-x_jh_j-w_jt_j\quad(0\le j\le d-2).
       \label{eq:tree-aux}
\end{align}
Here \(x_{-1}\coloneqq 0\), and the term \(w_jh_j\) is omitted at \(j=d-1\).
An empty sum is zero, so the formulas also cover \(d=1\).

For a single split (\(d=1\)) the equations can be solved explicitly.
If \(x_0(0)=a>0\) and \(H(0)=0\), then \(x_0^2+H^2=a^2\), so
\(H'=a^2-H^2\). Hence
\[
          H(z)=a\tanh(az),\qquad x_0(z)=\frac{a}{\cosh(az)} .
\]
The two leaves decay and the root increases to \(a\). At larger depth
the auxiliary triangles also exchange amplitude with the primary
vertices. The next lemma chooses their initial amplitudes so that the
same final root value is attained in the limit.

\begin{lemma}\label{lem:tree-orbit}
For every \(a>0\) and \(b_0,\ldots,b_{d-2}>0\), the system
\eqref{eq:tree-primary}--\eqref{eq:tree-aux} has a solution starting
with \(x_0>0\) and \(x_1=\cdots=x_d=0\). For \(z\ge0\),
\(0\le H\le a\) and \(-b_j\le t_j\le0\). The solution converges
exponentially to an equilibrium at which \(H=a\), \(h_j=\pm b_j\),
and every other amplitude is zero. Its first two derivatives decay
exponentially as well.
\end{lemma}

\begin{proof}
Direct differentiation gives the conserved quantities
\begin{equation}\label{eq:tree-integrals}
 w_j^2+h_j^2+t_j^2=b_j^2,\qquad
 \sum_{j=0}^d x_j^2+\sum_{j=0}^{d-2}h_j^2=a^2+\sum_{j=0}^{d-2}b_j^2 .
\end{equation}
To see the cancellation across levels, sum only up to height \(j<d\).
Equations \eqref{eq:tree-primary}--\eqref{eq:tree-aux} give
\[
 \frac{d}{dz}\left(\sum_{k=0}^{j}x_k^2+
              \sum_{k=0}^{\min(j,d-2)}h_k^2\right)
       =-2x_j^2x_{j+1}.
\]
At each level \(2h_kh_k'\) cancels the term \(-2x_kw_kh_k\).
The terms linking consecutive primary levels telescope, leaving only
the displayed boundary term. At \(j=d-1\) this is \(-(H^2)'\).
The first invariant follows by differentiating the three squares
on one auxiliary sphere. These identities are for the reduced
amplitudes, one per level. The multiplicities of the individual
Fourier modes enter later in the Sobolev norm.

For each initial point on the product of the indicated spheres,
set all primary amplitudes above height zero equal to zero, and set
\(x_0(0)=(a^2+\sum_j(b_j^2-h_j(0)^2))^{1/2}\).
This gives a continuous choice of initial conditions with \(x_0(0)\ge a\).
The trajectories are global by \eqref{eq:tree-integrals}. Let
\(\mathcal P\coloneqq \prod_{j=0}^{d-2}\{(w,h,t)\in\mathbb R^3:
w^2+h^2+t^2=b_j^2\}\). For \(Z>0\) let \(F_Z:\mathcal P\to\mathcal P\)
map these initial sphere coordinates to their values at time \(Z\).
The maps \(F_{uZ}\), \(0\le u\le1\), give a homotopy to the identity.
Thus \(F_Z\) has degree one and is onto. Indeed a map of a closed
oriented manifold missing a point has degree zero at that point.
Take a preimage of the point with coordinates \((0,b_j,0)\) for
every \(j\). At time \(Z\) the second invariant reduces to
\(H(Z)^2+\sum_{k<d}x_k(Z)^2=a^2\), where the sum starts at zero.
Since \(H(0)=0\) and \(H'\ge0\), we have \(0\le H\le a\)
on \([0,Z]\). Likewise \(t_j'\ge0\) and \(t_j(Z)=0\), while the
sphere bounds \(t_j\) below by \(-b_j\). For a sequence \(Z\to\infty\)
these preimages have a convergent subsequence in \(\mathcal P\).
Fix any finite observation time. Continuous dependence passes all the
inequalities at that time to the trajectory from the limiting initial
point. The same initial point therefore works on the whole half-line.
This passage to the limit cannot be omitted. For finite \(Z\) the
selected trajectory passes through the points \((0,b_j,0)\) at time
\(Z\), but in general it leaves them afterwards: there
\(w_j'=-x_jb_j\), and as soon as \(w_j\ne0\) the identity
\(t_j'=w_j^2\) makes \(t_j\) positive. For finite \(Z\) the bounds above
therefore hold only on \([0,Z]\), and the limit is what gives them on
the whole half-line.
For \(d=1\) the product of spheres consists of a single point and the
same argument uses the primary invariant alone.

The omega limit set consists of limits of states at times tending to
infinity. It is invariant under this finite-dimensional flow. The
bounded monotone coordinates \(t_j,H\) each have a limit, so these
coordinates are constant on the limit set. Their derivatives
\(w_j^2,x_{d-1}^2\) therefore vanish there. Invariance and
\eqref{eq:tree-primary} then give \(x_{d-2}=\cdots=x_0=0\).
The two identities in \eqref{eq:tree-integrals} give
\(H^2=a^2+\sum_jt_j^2\), so \(H=a\), all \(t_j=0\), and
\(h_j=\pm b_j\). The omega limit set is connected, being the
intersection of the nested compact connected closures of trajectory
tails. It therefore consists of a single point. If the trajectory left
some fixed neighborhood of this point at arbitrarily late times,
compactness would supply a second limit point. This proves convergence.

For the decay rate, the vector of the lower \(x_j\)'s and the \(w_j\)'s
solves a linear equation \(Y'=(A+E(z))Y\) along this trajectory, with
\(E(z)\to0\). The matrix \(A\) consists of the scalar \(-a\) on
\(x_{d-1}\) and the symmetric blocks
\(\left(\begin{smallmatrix}0&-h_j(\infty)\\-h_j(\infty)&0\end{smallmatrix}\right)\)
on \((x_j,w_j)\). For example, write \(x_{j-1}^2\) as the unknown
\(x_{j-1}\) times a coefficient tending to zero. Fix \(\lambda>0\)
less than the smallest absolute value of the eigenvalues of \(A\).
We use the supremum norm weighted by \(e^{\lambda(z-Z)}\) on a late
half-line \(z\ge Z\). Write \(P_-,P_+\) for the orthogonal projections
onto the negative and positive eigenspaces of \(A\), and write
\(\lambda_0\coloneqq \min_{\mu\in\operatorname{spec}(A)}|\mu|\). The equation
on \(z\ge Z\) takes the form
\[
 \begin{aligned}
 Y(z)={}&e^{(z-Z)A}P_-Y(Z)
     +\int_Z^z e^{(z-v)A}P_-E(v)Y(v)\,dv\\
 &-\int_z^\infty e^{(z-v)A}P_+E(v)Y(v)\,dv .
 \end{aligned}
\]
In obtaining the last integral, the boundary term at infinity vanishes because \(Y\) is bounded.
In the weighted norm the two integrals have combined operator norm at most
\(2(\lambda_0-\lambda)^{-1}\sup_{v\ge Z}\|E(v)\|\). They also contract
in the unweighted norm when \(Z\) is sufficiently large. Uniqueness
in the unweighted space identifies the bounded solution with the
weighted one. Integrating the equations for \(H,h_j,t_j\) from infinity
gives decay of the remaining differences. The polynomial equations and their first derivatives
give the assertion about derivatives.
\end{proof}

The Sobolev cost of this orbit at its initial time is particularly
simple. At each primary split, write \(\zeta,1-\zeta\) for the two
ratios of squared lengths. For \(s>1\) define
\[
 c_s\coloneqq \max_{p_0\le\zeta\le1-p_0}\big(\zeta^s+(1-\zeta)^s\big)<1 .
\]
If the tree has root \((P,0)\), iteration of the split gives
\begin{equation}\label{eq:leaf-cost}
  \sum_{\text{primary leaves }v}|v|^{2s}|x_0(0)|^2
       \le |P|^{2s}\Big(a^2+\sum_j b_j^2\Big)c_s^d .
\end{equation}
All higher primary levels are initially zero. Each auxiliary amplitude
at height \(j\) is bounded by \(b_j\) for all time. Since the geometry
is fixed before choosing the positive \(b_j\), the total weighted cost of
all the auxiliary vertices can be made arbitrarily small.

The identity \(H'=x_{d-1}^2\) gives a quantitative estimate on the time
spent below a prescribed fraction of the terminal amplitude. The estimate
is uniform over the trajectories selected by the degree argument.

\begin{proposition}\label{prop:clock}
Take \(a=1\), \(0<\theta<1\), and put
\[
 C_\theta\coloneqq \left(\frac{100}{1-\theta^2}\right)^4,\qquad
 Z\coloneqq C_\theta d^{3/2}.
\]
If \(\sum_{j=0}^{d-2}b_j^{3/2}\le Z^{-1/2}\), then
\(H(z)\ge\theta\) for every \(z\ge Z\). In physical time the root
has thus reached this fraction of its terminal amplitude by
\(C_\theta d^{3/2}/(2|g|r)\).
\end{proposition}
\begin{proof}
Write \(E_j\coloneqq Z^{-1}\int_0^Z x_j^2\,dz\), including \(x_d=H\).
From \eqref{eq:tree-integrals} and monotonicity,
\begin{equation}\label{eq:root-average}
  1-H(Z)^2\le\sum_{j<d}E_j,\qquad E_{d-1}\le Z^{-1},\qquad E_d\le1.
\end{equation}
Indeed at each time \(1-H^2=\sum_{j<d}x_j^2-\sum_j(w_j^2+t_j^2)\),
and \(H'=x_{d-1}^2\). For \(1\le j<d\), integrate
\eqref{eq:tree-primary}, using \(x_j(0)=0\),
\(\int_0^Z w_j^2\le b_j\), and \(|h_j|\le b_j\), to get
\begin{equation}\label{eq:mean-recursion}
 E_{j-1}\le (E_jE_{j+1})^{1/2}+e_j,\qquad
 e_j\coloneqq \frac{|x_j(Z)|}{Z}+\frac{b_j^{3/2}}{\sqrt Z},
\end{equation}
where \(b_{d-1}\coloneqq 0\). For example, the mixed term in that integration
is bounded by
\[
 Z^{-1}\left|\int_0^Z w_j h_j\,dz\right|
 \le b_j Z^{-1/2}\left(t_j(Z)-t_j(0)\right)^{1/2}
 \le b_j^{3/2} Z^{-1/2}.
\]
At height \(d-1\) this term is absent. As \(Z\ge1\), our choice of radii implies
\(\sum_jb_j^2\le1\). Hence
\[
      e_j\le(\sqrt2+1)/Z,\qquad
      \sum_{j=1}^{d-1}e_j\le(\sqrt{2d}+1)/Z .
\]
The second estimate uses Cauchy--Schwarz at the single endpoint \(Z\).

To exploit the geometric mean in \eqref{eq:mean-recursion}, order the
levels downward, putting \(m_k\coloneqq E_{d-1-k}\) for \(-1\le k<d\).
The starting bounds are \(m_{-1}\le1\) and \(m_0\le Z^{-1}\).
If the two preceding averages are bounded by powers of \(Z^{-1}\),
their geometric mean has the average of those exponents. This suggests
\(\alpha_{-1}\coloneqq 0,\ \alpha_0\coloneqq 1,\ 
\alpha_k\coloneqq (\alpha_{k-1}+\alpha_{k-2})/2\), giving
\(\alpha_k=2/3+(-1/2)^k/3\). For \(k\ge1\) they belong to
\([1/2,3/4]\). Multiplying \eqref{eq:mean-recursion} by
\(Z^{\alpha_k}\) shows that \(m_kZ^{\alpha_k}\) exceeds the larger
of the preceding two scaled averages by at most
\((\sqrt2+1)Z^{-1/4}\). Let
\(n\coloneqq 1+\lceil\log_2(2+\log Z)\rceil\), with \(\log\) the natural
logarithm. Since \(Z\ge100^4\), we have \(n+1\le4Z^{1/4}\) and
\(nZ^{-1/4}<5\) (use \(n\le 3+\log Z\)). Induction gives
\(m_k\le15Z^{-\alpha_k}\) up to \(k=\min(n,d-1)\).
These levels contribute at most \(60Z^{-1/4}\) to \eqref{eq:root-average}.

\smallskip
\noindent\emph{Details of the induction.} The invariant
\eqref{eq:tree-integrals} gives
\(|x_j(Z)|\le(1+\sum_jb_j^2)^{1/2}\le\sqrt2\), and
\(b_j^{3/2}\le Z^{-1/2}\); this is the bound
\(e_j\le(\sqrt2+1)/Z\). Also
\(\sum_{j\ge1}|x_j(Z)|\le\sqrt d\,(\sum_jx_j(Z)^2)^{1/2}\le\sqrt{2d}\) and
\(\sum_jb_j^{3/2}/\sqrt Z\le Z^{-1}\). With \(j=d-k\),
\eqref{eq:mean-recursion} reads
\(m_k\le(m_{k-1}m_{k-2})^{1/2}+e_{d-k}\) for \(1\le k\le d-1\). Put
\(s_k\coloneqq m_kZ^{\alpha_k}\). Since
\(\alpha_k=(\alpha_{k-1}+\alpha_{k-2})/2\) and \(\alpha_k\le3/4\),
\[
 s_k\le(s_{k-1}s_{k-2})^{1/2}+(\sqrt2+1)Z^{\alpha_k-1}
 \le\max(s_{k-1},s_{k-2})+(\sqrt2+1)Z^{-1/4},
\]
and \(s_{-1},s_0\le1\). Hence
\(s_k\le1+k(\sqrt2+1)Z^{-1/4}<15\) for \(k\le n\), because
\(n\le2+\log_2(2+\log Z)\le3+\log Z\) (as \(2+L\le2^{1+L}\) for \(L\ge0\))
and \((3+\log Z)Z^{-1/4}<1\) for \(Z\ge100^4\); the same numbers give
\(n+1\le4Z^{1/4}\). As \(\alpha_k\ge1/2\) for \(k\ge1\) and
\(m_0\le Z^{-1}\), the levels \(0\le k\le\min(n,d-1)\) contribute at most
\(15(n+1)Z^{-1/2}\le60Z^{-1/4}\). Next, \(|\alpha_k-2/3|=2^{-k}/3\) and
\(2^{n-1}\ge2+\log Z\), so that \(Z^{-\alpha_k}\le e^{1/3}Z^{-2/3}\) for
\(k\in\{n-1,n\}\), and \(15e^{1/3}<21\). On the remaining levels
\(m_k\le\max(m_{k-1},m_{k-2})+e_{d-k}\) gives by induction
\(m_k\le21Z^{-2/3}+\sum_je_j\). Finally \(Z=C_\theta d^{3/2}\) gives
\(21dZ^{-2/3}=21C_\theta^{-2/3}\) and
\(d(\sqrt{2d}+1)/Z=(\sqrt2+d^{-1/2})C_\theta^{-1}\le3C_\theta^{-1}\), while
\(60C_\theta^{-1/4}=0.6(1-\theta^2)\),
\(21C_\theta^{-2/3}\le21\cdot10^{-16/3}(1-\theta^2)\) and
\(3C_\theta^{-1}\le3\cdot10^{-8}(1-\theta^2)\). The sum is less than
\(1-\theta^2\).

If there are further levels, the last two exponents already satisfy
\(|\alpha_k-2/3|\log Z\le1/3\), for \(k=n-1,n\). Thus both averages
are at most \(21Z^{-2/3}\). For the remaining levels use
\(\sqrt{uv}\le\max(u,v)\) in \eqref{eq:mean-recursion}. Each is bounded
by \(21Z^{-2/3}+\sum_j e_j\). Substitution into \eqref{eq:root-average}
yields
\[
  1-H(Z)^2\le 60Z^{-1/4}+21dZ^{-2/3}
                  +\frac{d(\sqrt{2d}+1)}{Z}
   \le 60C_\theta^{-1/4}+21C_\theta^{-2/3}+3C_\theta^{-1}
   <1-\theta^2 .
\]
Monotonicity of \(H\) proves the proposition.
\end{proof}

There is also a lower bound within the finite system. The moment
\[
  \sum_{j=0}^d j x_j^2+\sum_{j=0}^{d-2}j h_j^2
\]
has derivative \(2\sum_{j<d}x_j^2x_{j+1}\le2B^3\), where
\(B^2\coloneqq 1+\sum_jb_j^2\). Indeed, by
\eqref{eq:tree-primary}--\eqref{eq:tree-aux} the terms
\(\mp2jx_jw_jh_j\) cancel, and
\(2\sum_jj(x_jx_{j-1}^2-x_j^2x_{j+1})
=2\sum_{j<d}((j+1)-j)x_j^2x_{j+1}\); the bound uses
\(\sum_jx_j^2\le B^2\) and \(|x_{j+1}|\le B\). Its initial value is at most
\(d\sum_jb_j^2\). Thus, if \(\sum_jb_j^2\le\theta^2/2\), the model
time of the first crossing \(H=\theta\) is at least
\(d\theta^2/(4B^3)\). Both estimates concern a single choice of tree,
radii, and initial angles.

\section{Insertion of a tree into the two-dimensional flow}\label{sec:insertion}

Let \(V\) be any smooth axial solution with finite Birkhoff support whose
action vector lies in the ball of Lemma~\ref{lem:locking}. Assume its
occupied frequencies in the mass gauge are rational and its space-time
coefficient \eqref{eq:bank} is nonzero. Proposition~\ref{prop:bank}
provides an initial solution with these properties. Throughout this
section axial solutions are written in their mass gauge. For mass \(M\) the equation has the additional constant
potential \(-2M\). Write \(g\) for the \((1,1)\) space-time Fourier
coefficient of \(V\), and \(F\) for its finite set of occupied Birkhoff
indices. Fix the tree of Proposition~\ref{prop:tree} and an orbit given
by Lemma~\ref{lem:tree-orbit}. Let \(\mathcal V\subset\mathbb Z^2\) be
its full set of velocities, including the auxiliary vertices. For
\(v=(v_x,v_y)\in\mathcal V\), write \(x_v(z)\) for the real amplitude
assigned to that vertex in \eqref{eq:tree-primary}--\eqref{eq:tree-aux}.
We use three time variables. The equation is solved in physical time \(t\),
the error estimates use slow time \(\tau\coloneqq rt\), and the tree uses time
\(z\coloneqq 2|g|\tau\), in which its interaction coefficient equals one. The
positive amplitude scale \(r\) is specified in the next proposition.
In complex amplitudes these functions all carry the same constant
unit factor \(\gamma\coloneqq i g/|g|\). Enumerate the horizontal vertices as
\((P_i,0)\), and let \(c_i\coloneqq \lim_{z\to\infty}x_{(P_i,0)}(z)\).
At the main root \(P_i=P\), this limit is \(a>0\). At each other
horizontal vertex it is one of the signed auxiliary radii. Both the tree and its
positive auxiliary radii are fixed throughout the proof below.

\begin{proposition}\label{prop:insertion}
Let \(1<s_0<s<2\) and \(a_0>0\). Set \(r\coloneqq a_0R^{-s}\) and
\(k_v\coloneqq (1,0)+Rv\). For all sufficiently large integers \(R\) there are
a periodic axial solution \(V_R\) and a solution \(W_R\) of the full
two-dimensional equation in the mass gauge of \(V_R\), with the following
properties.
\begin{enumerate}
\item On \(F\), \(V_R\) and \(V\) have the same mass-gauged frequencies and
the same Birkhoff angles at \(t=0\). The newly occupied Birkhoff indices
of \(V_R\) are \(K_i\coloneqq 1+RP_i\). At \(K_i\) the modulus of the Birkhoff
coordinate is \(r|c_i+e_i|\) for some real \(e_i\) with \(|e_i|\le R^{-10}\).
\item \(W_R(t)-V_R(t)\) tends to zero exponentially in \(H^s\) as
\(t\to+\infty\). In the respective gauges of \(W_R\) and \(V\),
\[
       \sup_{t\ge0}\|W_R(t)-V(t)\|_{H^{s_0}}=o(1).
\]
\item There are real numbers \(\kappa_v=|k_v|^2+O(R^{-2}+r^2)\) such that
\begin{equation}\label{eq:insertion-profile}
 \sup_{t\ge0}
 \left\|W_R(t)-V(t)
 -r\gamma\sum_{v\in\mathcal V} x_v(2|g|rt)
                   e^{i k_v\cdot(x,y)-i\kappa_v t}\right\|_{H^s}
       \longrightarrow0 .
\end{equation}
In particular
\[
 \|W_R(0)-V(0)\|_{H^s}\le a_0\Big(\sum_v |v|^{2s}|x_v(0)|^2\Big)^{1/2}+o(1).
\]
At the main root
\(K\coloneqq 1+RP\) one has, writing \(H\coloneqq x_{(P,0)}\) and using
\(\widehat W_R(K,0,t)\) for the coefficient at \((K,0)\) at time \(t\),
\[
 \sup_{t\ge0}\left|\,|\widehat W_R(K,0,t)|-r H(2|g|rt)\,\right|=o(r).
\]
\end{enumerate}
All limits in this statement are taken after the finite tree and the
old orbit have been fixed.
\end{proposition}

The comparison in \eqref{eq:insertion-profile} involves three exact
solutions. The old periodic field \(V\) supplies the coefficient \(g\).
The enlarged periodic field \(V_R\) prescribes the final axial amplitudes
and phases. It is \(W_R\) that undergoes the two-dimensional transfer
and converges to \(V_R\). In particular the reference field in the
finite-plane expansion \eqref{eq:insertion-profile} is \(V\).

There are two different roles for the linear frequencies in this proof.
An error in the frequency of a persistent axial mode accumulates for
all future time. We therefore evolve these modes at their exact Birkhoff
frequencies. Off the axis, the limiting field vanishes, but its linearized
equation has both growing and decaying directions. We identify these
directions explicitly and use them to prescribe the boundary conditions
for the error. Unless indicated otherwise, constants in this section
may depend on the fixed old orbit, tree, and \(a_0,s_0,s\).
They are independent of \(R\) and of the terminal time.

We prove Proposition~\ref{prop:insertion} by first prescribing the
axial endpoint and then solving for the remaining Fourier coefficients.
Section~\ref{sec:ins-norms} gives the equations for those coefficients.
Section~\ref{sec:ins-propagator} determines which transverse directions
can be integrated from the left and which must be integrated from the
right. In Section~\ref{sec:ins-terminal} we substitute the tree orbit
into the full equations and estimate the errors in the transverse
field, the new axial field, and the old actions and angles.
Section~\ref{sec:ins-half} solves for them first on
\([T,S]\) in slow time, sends \(S\) to infinity with \(R\) fixed, and
extends the resulting solution back to time zero.

\subsection{The norm and the axial equations}\label{sec:ins-norms}

Use Lemma~\ref{lem:locking} to adjust the old actions by \(O(r^2)\)
and to choose the new radii in the intervals of part (1).
Only finitely many new actions are involved, and the old ones stay positive.
Let \(I^\infty\) denote the resulting action vector, let
\(c_R\coloneqq -2\sum_j I_j^\infty\), and take the new initial phases so that the
complex terminal amplitudes are \(r\gamma(c_i+e_i)\).
All renormalized occupied frequencies are rational, so the resulting axial
solution \(V_R\) is periodic in this gauge. In particular its constant
potential differs from that of \(V\) by \(O(r^2)\).

For one insertion it suffices to work in the exact invariant subspace
with vertical Fourier indices in \(R\mathbb Z\). These are the functions
fixed by translation of \(y\) by \(2\pi/R\). Write \(p\coloneqq (s+2)/2\) and set
\begin{equation}\label{eq:parameters-ins}
 \alpha\coloneqq \tfrac12\min(2-s,p+s-2),\qquad L\coloneqq R^\alpha,\qquad
 \rho\coloneqq R^{-3\alpha/4}.
\end{equation}
Here \(s<p<2\), \(L\to\infty\), and \(L/R\to0\). For
\(n,q\in\mathbb Z\) the integer \(q\) labels the vertical frequency
\(qR\). We weight the Fourier coefficient at \((n,qR)\) by
\[
 w(n,q)\coloneqq (1+|n|/R+|q|)^p
              \bigl(1+L^{-1}\operatorname{dist}(n,R\mathbb Z)\bigr).
\]
Thus \(\operatorname{dist}(n,R\mathbb Z)\coloneqq \min_{j\in\mathbb Z}|n-Rj|\).
The first factor measures the frequency in units of \(R\). The second
measures its horizontal displacement from a multiple of \(R\).
For \(f_{n,q}\coloneqq \widehat f(n,qR)\) put
\[
  \|f\|_Z^2\coloneqq \sum_{n,q}w(n,q)^2|f_{n,q}|^2,\qquad
  \|f\|_B\coloneqq \sum_{n,q}w(n,q)|f_{n,q}|.
\]
For sequences on the horizontal axis only the terms with \(q=0\) occur.
Reflection, conjugation and Fourier coordinate projections are bounded
in both norms. Since \(\langle(n,qR)\rangle^b\le C R^b w(n,q)\) for
\(0\le b\le p\), we have \(\|f\|_{H^b}\le C R^b\|f\|_Z\).

\begin{lemma}\label{lem:product-insertion}
Uniformly in \(R\),
\[
 \|fh\|_Z\le C\|f\|_B\|h\|_Z,\qquad
 \|fh\|_Z\le C\sqrt L\,\|f\|_Z\|h\|_Z .
\]
A plane wave at \(k_v\) times an axial function in a bounded smooth
family has bounded \(B\) norm.
\end{lemma}
\begin{proof}
In the coordinates \((n,q)\), the scaled length \(|n|/R+|q|\)
and the distance \(\operatorname{dist}(n,R\mathbb Z)\) are both
subadditive. Since \(1+a+b\le(1+a)(1+b)\) for \(a,b\ge0\),
the weight satisfies \(w(k)\le w(l)w(k-l)\). The first estimate
now follows from \(\ell^1*\ell^2\to\ell^2\).
For the second, write \(k=(n,q)\) for the two indices of the weight.
Cauchy--Schwarz in each convolution output gives
\[
 \|fh\|_Z^2\le
 \left(\sup_k\sum_l\frac{w(k)^2}{w(l)^2 w(k-l)^2}\right)
 \sum_{k,l}w(l)^2|f_l|^2 w(k-l)^2|h_{k-l}|^2 .
\]
The second factor is \(\|f\|_Z^2\|h\|_Z^2\).
To bound the first, write \(n=Rj+d\), \(-R/2<d\le R/2\).
The first factor in \(w\) is comparable to \(a(j,q)\coloneqq (1+|j|+|q|)^p\).
If \(k=l+(k-l)\), the two input \((j,q)\) indices add to the output
up to a carry of \(0,1\), or \(-1\) in \(j\). For each fixed input
remainder, the carry is fixed by the output. Consequently the squared
ratio of the three \(a\) weights is bounded by a constant times the sum
of the two reciprocal squared input weights, by
\((A+B)^p\le 2^{p-1}(A^p+B^p)\) for \(A,B\ge0\). Summing over the first
input \((j,q)\) costs a constant since \(\sum_{j,q}a(j,q)^{-2}<\infty\).
For the second factor put \(b(d)=1+|d|/L\) on the centered
representatives of \(\mathbb Z/R\mathbb Z\). For an input remainder
\(d_1\), the other remainder \(d_2\) is the centered representative
of \(d-d_1\). The triangle inequality on this cyclic group gives
\[
 \frac{b(d)^2}{b(d_1)^2b(d_2)^2}
 \le 2\bigl(b(d_1)^{-2}+b(d_2)^{-2}\bigr).
\]
The sum of this right-hand side over \(d_1\) is at most
\(4\sum_{|h|\le R/2}(1+|h|/L)^{-2}\le CL\). The preceding
sum over \((j,q)\) was uniform in \(d_1\), so the entire squared
convolution bound costs \(CL\). Taking its square root gives the
stated factor \(\sqrt L\).
Finally the velocity components are integers, so shifting an axial envelope
by \(k_v\) changes its remainder by one and its large-scale weight by at
most a fixed polynomial in the envelope index.
\end{proof}

Let \(P_0\) denote averaging in \(y\), equivalently the Fourier
projection onto \(q=0\), and let \(P_\perp\coloneqq 1-P_0\).
For a solution \(u\) in gauge \(c_R\), write \(U\coloneqq P_0 u,\ w_\perp\coloneqq P_\perp u\), and
\(z_j\coloneqq (\Phi(U))_j\), where \(\Phi\) is the axial Birkhoff map and
\(\Psi=\Phi^{-1}\). The index on \(z_j\) distinguishes the Birkhoff
variables from the model time \(z=2|g|rt\). We keep \(c_R\)
fixed when projecting the equation. The transverse component can
exchange mass with \(U\), so the axial mass alone need not determine
the prescribed gauge at intermediate times. In the expansion of
\(|U+w_\perp|^2(U+w_\perp)\), the two terms linear in \(w_\perp\) are
\(2|U|^2 w_\perp+U^2\overline{w_\perp}\). Both have zero \(y\)-average.
Projecting the remaining terms and differentiating the real map \(\Phi\)
on the axial component gives the exact equations
\begin{align}
 \dot z_j&=-i(\omega_j(I)+c_R)z_j+
                   D\Phi_j(U)[-i f],\label{eq:axis-insertion}\\
 f&\coloneqq 2U P_0|w_\perp|^2+\overline U P_0(w_\perp^2)
                          +P_0(|w_\perp|^2w_\perp),\notag\\
 i\dot w_\perp&=(-\Delta+c_R+2|U|^2)w_\perp+U^2\overline{w_\perp}
                           +N(U,w_\perp),\label{eq:trans-insertion}\\
 N(U,w)&\coloneqq P_\perp(2U|w|^2+\overline U w^2+|w|^2w).\notag
\end{align}
Dots denote physical time derivatives, \(I_j=|z_j|^2\), and \(\omega_j\)
are the ungauged frequencies. To justify \eqref{eq:axis-insertion} for
\(H^p\) solutions, first prescribe a continuous axial \(H^p\) forcing.
Approximate the forcing and initial data by smooth functions and use the
integral equation in Birkhoff variables after removing the fixed diagonal
\(j^2\). Indeed, if \(I=|z|^2\), then
\[
 \omega_j(I)+c_R-j^2=2\|z\|_{\ell^2}^2+c_R+\chi_j(I),\qquad
 \big\||z|^2-|z'|^2\big\|_{\ell^2}
 \le(\|z\|_{\ell^2}+\|z'\|_{\ell^2})\|z-z'\|_{\ell^2}.
\]
The first formula keeps track of the mass of the instantaneous
axial field. The second, together with the local analyticity of
\(I\mapsto\chi(I)\) into \(\ell^2\subset\ell^\infty\)
\cite{Molnar} (Appendix~\ref{app:molnar-four} at \(p=4\)), makes the
remaining diagonal locally Lipschitz as
a multiplier on \(h^p\). The chart and its differential are locally
Lipschitz at this regularity \cite{GKbook,KST}. The resulting
integral equation therefore has a unique local solution and passes
to the limit under the smooth approximations.
The same reasoning allows reconstruction by \(\Psi\).

On \(F\) write \(z_j=\sqrt{I_j}\exp(-i\theta_j)\), and put
\begin{equation}\label{eq:old-mod-insertion}
 J\coloneqq (I_F-I_F^\infty)/r,\qquad \Theta\coloneqq \theta_F-\theta_F^\infty(t),
 \qquad V_*\coloneqq \Psi(z|_F).
\end{equation}
Here \(I_F=(I_j)_{j\in F}\), and similarly for the other restrictions.
The terminal angles are the continuous real lifts
\(\theta_j^\infty(t)\coloneqq \theta_j^\infty(0)+(\omega_j(I^\infty)+c_R)t\)
on \(F\). We use the lift of the actual angles near these nominal angles,
so that \(\Theta\) measures their displacement. In \(z|_F\) we set all
coordinates outside \(F\) to zero. The profile \(V_*\) uses the actual old
coordinates and is smooth on a fixed compact parameter neighborhood.

Two chart estimates will be used repeatedly. Write \(\mathcal A\) for
either chart and \(\mathcal F_{\mathcal A}\coloneqq D\mathcal A(0)\) for its linear
Fourier part. We identify spatial functions with their Fourier coefficients in norms
on either side of the chart. For a point \(x_*\) on the old smooth family
and an axial perturbation \(y_*\) of \(Z\) norm at most \(Cr\), we have
\begin{align}
 \|D\mathcal A(x_*+y_*)-D\mathcal A(x_*)\|_{Z\to Z}
      &\le C r(1+R^{2-p}),\label{eq:chart-tube}\\
 \|(D\mathcal A(x_*)-\mathcal F_{\mathcal A})h\|_Z
      &\le C(R^{-1}+L^{-1})\label{eq:chart-head}
\end{align}
in the second line for \(h=e^{iR v_x x}h_0(x)\) with fixed \(v_x\in\mathbb Z\setminus\{0\}\)
and \(h_0\) in a bounded smooth family. The same bound applies to finite
sums with summable smooth-envelope bounds, by linearity over the reals.

For clarity, first drop the residue factor. On the axis the resulting
norm \(Z_0\) is equivalent to
\(\|\cdot\|_{L^2}+R^{-p}\|\cdot\|_{H^p}\) (and to the corresponding
sum with \(\ell^2\) and \(h^p\) norms for coordinate sequences). Thus
\(\|h\|_{H^1}\lesssim R\|h\|_{Z_0}\), and \(Rr=o(1)\).
The full chart is analytic in \(L^2\), while its difference from Fourier
is analytic \(H^1\to H^2\) in both directions \cite{KST}.
Let \(\|h\|_Z\le1\) and write \(u_t=x_*+t y_*\), \(0\le t\le1\).
The base points stay in a fixed neighborhood of the compact old family
in \(H^1\), since \(\|y_*\|_{H^1}\le CRr=o(1)\). The fundamental
theorem of calculus gives
\[
 (D\mathcal A(x_*+y_*)-D\mathcal A(x_*))h
   =\int_0^1D^2\mathcal A(u_t)[y_*,h]\,dt.
\]
The \(L^2\) bound on the integrand is \(Cr\). The \(H^p\) bound
uses \(H^2\subset H^p\) and the second differential of the
one-smoothing part, since the Fourier part has zero second differential.
This gives \(C\|y_*\|_{H^1}\|h\|_{H^1}\le CrR^2\).
Multiplying by \(R^{-p}\) gives the second term in the claimed
\(Z_0\) estimate. To restore the residue factor, note that
\(\|h\|_Z\asymp\|h\|_{Z_0}+\|\delta h\|_{Z_0}\), where
\(\delta\coloneqq (R/L)(T_{2\pi/R}-1)\) and \(T_\ell\) is spatial translation,
or its diagonal action on Birkhoff coordinates. Indeed the modulus of
the Fourier multiplier of \(\delta\) is
\(2(R/L)|\sin(\pi n/R)|\), comparable to
\(\operatorname{dist}(n,R\mathbb Z)/L\). The charts commute with translation.
To see the corresponding difference of the integrand, abbreviate
\(T=T_{2\pi/R}\), \(u=u_t\). Equivariance and bilinearity give
\[
\begin{aligned}
 \delta\bigl(D^2\mathcal A(u)[y_*,h]\bigr)
 ={}&D^2\mathcal A(Tu)[\delta y_*,Th]
      +D^2\mathcal A(Tu)[y_*,\delta h]\\
 &+\int_0^1D^3\mathcal A\bigl(u+v(Tu-u)\bigr)
                 [\delta u,y_*,h]\,dv .
\end{aligned}
\]
In the first two terms \(\|\delta y_*\|_{Z_0}\le Cr\) and
\(\|\delta h\|_{Z_0}\le C\) give the same bounds as before.
In the last term \(\|\delta u\|_{H^1}\le C(L^{-1}+Rr)\)
because \(x_*\) is uniformly smooth, and
\(\|\delta u\|_{L^2}\le C(L^{-1}+r)\). Both are bounded
independently of \(R\). Applying the \(L^2\) and the one-smoothing
bounds for the third differential completes \eqref{eq:chart-tube}. At the smooth base \(x_*\),
\(D\mathcal A(x_*)\) is bounded on \(L^2\) and on \(H^p\), hence also
on \(Z_0\). Applying the same translation identity to the first
differential adds a term involving \(D^2\mathcal A\) and
\(\delta x_*=O(L^{-1})\) in \(H^1\). Its scaled \(H^p\) norm is
\(O(R^{1-p}/L)\) for an input of unit \(Z\) norm. This also gives
\(\|D\mathcal A(x_*+y_*)\|_{Z\to Z}\le C\), since \(rR^{2-p}=o(1)\).

To prove \eqref{eq:chart-head}, at a smooth base the differential minus
its Fourier part maps \(H^1\to H^2\). Its real Hilbert adjoint gains
the same derivative. Indeed if \(S\) is the full chart differential and
\(\mathcal J_{\rm in},\mathcal J_{\rm out}\) are the constant symplectic
matrices in the two sets of coordinates, then
\(S^*\mathcal J_{\rm out}S=\mathcal J_{\rm in}\). Applying this formula
also to \(S_0=\mathcal F_{\mathcal A}\) gives explicitly
\[
 (S-S_0)^*=\mathcal J_{\rm in}(S^{-1}-S_0^{-1})
                         \mathcal J_{\rm out}^{-1}.
\]
Here \(S^{-1}\) is the differential of the inverse chart at the image
of the same smooth point. Each \(\mathcal J\) and its inverse preserve Sobolev norms.
Duality first gives \(H^{-2}\to H^{-1}\) for the chart differential minus
Fourier. Interpolation with \(H^1\to H^2\) gives both \(H^{-1}\to L^2\)
and \(H^{p-1}\to H^p\). A shifted smooth envelope \(h\) as in
\eqref{eq:chart-head} has \(\|h\|_{H^{-1}}=O(R^{-1})\) and
\(\|h\|_{H^{p-1}}=O(R^{p-1})\), by rapid Fourier decay of the envelope.
In the two terms of \(Z_0\) these bounds both give \(O(R^{-1})\).
If \(A(x_*)\coloneqq D\mathcal A(x_*)-\mathcal F_{\mathcal A}\), equivariance gives
\[
 \delta(A(x_*)h)=A(T_{2\pi/R}x_*)\delta h+
   (R/L)\bigl(D\mathcal A(T_{2\pi/R}x_*)-D\mathcal A(x_*)\bigr)h .
\]
Since \(v_x\) is an integer, translation by \(2\pi/R\) preserves the
carrier and \(\delta h=e^{iR v_x x}\delta h_0\). The mean value theorem
on the smooth envelope gives \(\|\delta h_0\|_{H^p}=O(L^{-1})\).
Modulation by this carrier is uniformly bounded in \(Z_0\), with a
constant depending on \(v_x\). Thus the first term costs \(O(L^{-1})\).
In the second term \(\|\delta x_*\|_{H^1}=O(L^{-1})\). The \(L^2\)
and scaled \(H^p\) estimates cost \(O(L^{-1})\) and
\(O(R^{1-p}L^{-1})\). This proves the claimed bound.

\subsection{The transverse propagator}\label{sec:ins-propagator}

We first retain from the linear part of \eqref{eq:trans-insertion} the
old coefficients at \(V_*\), together with all terms linear in the
terminal axial planes \(rY_\infty\). Here
\[
 Y_\infty(t,x)\coloneqq \gamma\sum_i(c_i+e_i)
               e^{iK_i x-i(K_i^2+\chi_{K_i}(I^\infty))t},
 \qquad \chi_j(I)=\omega_j(I)-j^2-2\sum_k I_k .
\]
By Lemma~\ref{lem:empty-frequency} and the \(O(r^2)\) change of the actions,
\(\max_i|\chi_{K_i}(I^\infty)|=O(R^{-2}+r^2)\).
The scalar \(2\sum_{j\in F} I_j+c_R\) goes into the remainder because it
equals \(2r\sum_{j\in F}J_j-2\sum_{j\notin F} I_j^\infty\) and is
\(O(r|J|+r^2)\). For each \(m\coloneqq qR>0\), use the two horizontal sequences
\[
 \bigl(\widehat{w_\perp}(n,m),\ \overline{\widehat{w_\perp}(n,-m)}\bigr)_{n\in\mathbb Z}.
\]
Thus the label in the second copy is the index \emph{before} conjugation.
These two copies, taken over \(q=1,2,\ldots\), represent all transverse
Fourier coefficients. Multiplication by an axial potential acts within
a copy. The term \(U^2\overline{w_\perp}\) acts between copies and its
upper-right entry at \((n,n')\) is \(\widehat{U^2}(n+n')\).
In particular the blocks for different \(m>0\) evolve independently
under the axial linearization.

We make the following three changes to the transverse linear equation.
First remove the off-diagonal term formed solely from \(V_*\). Next
diagonalize the scalar old potential up to a prescribed order at large
horizontal indices. Finally use the space-time harmonics of the old
orbit to separate the terms linear in \(rY_\infty\). All three steps
are performed in a single spatial norm.

We use the following version of the Schur test. Given a positive Hilbert
weight \(v_\iota\), let the size of a matrix be the larger of the two
suprema of its absolute row and column sums, after multiplying entry
\((\iota,\jmath)\) by
\(\max(v_\iota/v_\jmath,v_\jmath/v_\iota)\).
These sizes are submultiplicative, are unchanged by adjoints, and bound
the operator norms on weighted \(\ell^1\) and \(\ell^2\).
It is convenient to replace the weight in Lemma~\ref{lem:product-insertion}
by the equivalent weight
\begin{equation}\label{eq:flat-weight-ins}
 \widetilde w(n,q)\coloneqq 
 (1+|q|+(|n|/R-C_*)_+)^p
 \bigl(1+(\operatorname{dist}(n,R\mathbb Z)/L-2)_+\bigr),
\end{equation}
where \(a_+\coloneqq \max(a,0)\). We take \(C_*\coloneqq 2+\max_{v\in\mathcal V}|v_x|\),
so all axial indices \(1+RP_i\) lie in the flat part of the first
factor for large \(R\). Under a shift or reflection
\(n\mapsto\pm n+RP_i+b\), its ratios are bounded by
\(C(1+|b|)^{p+1}\), and by a constant independent of \(R\) when
\(|b|\le C'\log R\). The same holds with \(P_i=0\). This bound covers
both factors of the product weight in all Schur estimates below.

\emph{The old axial potential.} We begin with the order-one terms from \(V_*\). In the doubled
generator the two free blocks are \(A\coloneqq D_x^2+m^2\) and \(-A\),
\(D_x\coloneqq -i\partial_x\). Write \(C_{nn'}\coloneqq \widehat{V_*^2}(n+n')\),
\(D_{\rm free}\coloneqq \operatorname{diag}(A,-A)\), and put
\[
 X_{nn'}\coloneqq -\frac{C_{nn'}}{n^2+(n')^2+2m^2},\qquad
 J_{\rm off}\coloneqq \begin{pmatrix}0&X\\X^*&0\end{pmatrix}.
\]
A star denotes the ordinary Hilbert space adjoint on the Fourier
sequences. Since \(AX+XA=-C\), the off-diagonal entries of
\([D_{\rm free},J_{\rm off}]\) are \(-C,C^*\), where
\([D,J]=DJ-JD\). The denominator is at least \(2R^2\), and the smooth
old coefficients decay rapidly in \(n+n'\). It follows from the Schur
test above that \(J_{\rm off}\) and its derivatives in the finitely many
old angles and actions have size \(O(R^{-2})\), uniformly in \(q\).
Use \(e^{J_{\rm off}}\) as the map from new to physical variables.
In \(e^{-J_{\rm off}}D_{\rm free}e^{J_{\rm off}}\) the first
commutator cancels the old off-diagonal pair \(C,-C^*\).
Each further commutator has an additional factor \(J_{\rm off}\).
The time derivative costs \(O(R^{-2})\) at bounded parameter speeds.
The resulting scalar potential, after removal of the scalar specified
above, is \(B_0\coloneqq 2(|V_*|^2-\sum_{j\in F}I_j)\) in the upper copy and
its signed conjugate in the lower copy. It is real and has zero spatial
mean, by the mass identity for the axial chart.

For this scalar potential we construct an orthonormal basis depending
on the old angles and actions. Holding the actions fixed gives the nominal angular derivative
\(\partial_*\coloneqq \sum_{j\in F}(\omega_j(I^\infty)+c_R)\partial_{\theta_j}\). Fix an integer \(k_*\)
with \(k_*\alpha\ge s+2\alpha\). For \(|n|\ge L\) seek a column
\(e^{inx}f_n\) satisfying, to order \(O(|n|^{-k_*})\),
\[
 (2nD_x+D_x^2+B_0-i\partial_*)f_n=d_n f_n.
\]
To construct the column, write \(\eta\coloneqq (2n)^{-1}\),
\(\mathcal H\coloneqq D_x^2+B_0-i\partial_*\), and use the finite expansions
\(f_n\coloneqq \sum_{j=0}^{k_*}\eta^j F_j\),
\(d_n^{\rm app}\coloneqq \sum_{j=0}^{k_*-1}\eta^j\lambda_j\), with \(F_0=1\).
Let \(\langle h\rangle_x\) denote the spatial mean. If the coefficients
through \(F_j\) have been found, set
\[
 S_j\coloneqq \mathcal H F_j-\sum_{a=0}^{j-1}\lambda_a F_{j-a},
 \qquad \lambda_j\coloneqq \langle S_j\rangle_x,\qquad
 F_{j+1}\coloneqq -D_x^{-1}(S_j-\lambda_j).
\]
The sum is empty for \(j=0\). We take the inverse on mean-zero
functions, with zero-mean output. Thus the coefficient of \(\eta^j\)
in \((\eta^{-1}D_x+\mathcal H)f_n-d_n^{\rm app}f_n\) vanishes for
\(0\le j<k_*\). In particular \(\lambda_0=0\),
\(F_1=-D_x^{-1}B_0\), and
\(\lambda_1=-\langle B_0D_x^{-1}B_0\rangle_x=0\). For the last
identity pair the Fourier labels \(l,-l\), whose denominators have
opposite sign. Only finitely many smooth coefficients are needed, so
the same bounds hold after any fixed number of old-parameter derivatives.
Normalize each column in \(L^2\).
Its component of the defect perpendicular to the column is still
\(O(|n|^{-k_*})\). The parallel component of the normalized
\((D_x^2+B_0-i\partial_*)\)-column is real since the derivative of its
squared norm vanishes. Taking this real expectation minus \(n^2\) as
\(d_n\) gives the asserted residual and \(d_n=O(|n|^{-2})\). Indeed
the normalizing factor is \(1+O(n^{-2})\) since the first correction
to \(f_n\) had zero mean. The residual Fourier kernels are bounded by
\(C_j|n|^{-k_*}(1+|n'-n|)^{-j}\) for every fixed \(j\). The same
argument applies in the conjugate copy.

The approximate columns need not be mutually orthogonal. We correct
their Gram matrix, including at \(|n|=L\). In the next paragraph \(P\)
is the projection onto the horizontal indices \(|n|\ge L\), \(Q\coloneqq 1-P\),
and \(T:P\ell^2\to\ell^2\) maps the \(n\)-th unit sequence to the
normalized column just constructed. Set \(D_h\coloneqq (D_x^2+\operatorname{diag}
d_n)|_{P\ell^2}\) and \(H_*\coloneqq D_x^2+B_0\). We have \(H_*T-i\partial_*T-TD_h=E_*\) with
Schur size \(O(L^{-k_*})\). Multiplication by \(T^*\) and subtraction
of the adjoint identity give
\[
 [D_h,T^*T]-i\partial_*(T^*T)=T^*E_*-E_*^*T.
\]
Put \(G\coloneqq T^*T\) on \(P\ell^2\), and let
\(\mathfrak D B\coloneqq [D_h,B]-i\partial_*B\). This is a derivation, so
\(\mathfrak D(G-1)^j\) is the sum of \(j\) terms with one factor
\(\mathfrak DG\). Since \(G-1=O(L^{-1})\) in Schur size, the binomial
series gives \(\mathfrak D(G^{-1/2})=O(L^{-k_*})\).
Here the Leibniz rule is used on products for which the combined
commutator is Schur bounded. To check it for powers of \(G-1\), insert
finite horizontal cutoffs between factors. The cutoffs commute with
\(D_h\) and are independent of the old angles. The entrywise product
rule passes to the limit in distributions in those angles by absolute
convergence of the matrix products and their uniform Schur bounds.
The same argument applies to the absolutely convergent series below. Writing
\(G^{-1/2}=\sum_{j\ge0}\binom{-1/2}{j}(G-1)^j\) and
\(\mathfrak D(G-1)^j=\sum_{i=0}^{j-1}(G-1)^i(\mathfrak DG)(G-1)^{j-1-i}\), one gets
in Schur size
\[
 \|\mathfrak D(G^{-1/2})\|\le\sum_{j\ge1}\Bigl|\binom{-1/2}{j}\Bigr|\,j\,
 \|G-1\|^{j-1}\|\mathfrak DG\|
 =\tfrac12\bigl(1-\|G-1\|\bigr)^{-3/2}\|\mathfrak DG\|,
\]
and \(\mathfrak DG=T^*E_*-E_*^*T=O(L^{-k_*})\) by the displayed identity.
Moreover \(U_h^*U_h=G^{-1/2}GG^{-1/2}=1\), and the residual of \(U_h\) is
\[
\begin{aligned}
 H_*U_h-i\partial_*U_h-U_hD_h
 &=(H_*T-i\partial_*T-TD_h)G^{-1/2}
   +T\bigl([D_h,G^{-1/2}]-i\partial_*G^{-1/2}\bigr)\\
 &=E_*G^{-1/2}+T\mathfrak D(G^{-1/2}).
\end{aligned}
\]
The high-column map \(U_h\coloneqq TG^{-1/2}\) is an isometry in unweighted
\(\ell^2\), and its residual equals
\(E_*G^{-1/2}+T\mathfrak D(G^{-1/2})\). It is still \(O(L^{-k_*})\).
To extend \(U_h\), put \(\widetilde Q=1-U_hU_h^*\), the projection
onto the orthogonal complement of its range. Since \(U_h\) differs
from the inclusion of the high subspace by \(O(L^{-1})\), we have
\(\widetilde Q-Q=O(L^{-1})\). On the low subspace \(Q\ell^2\)
the operator \(Q\widetilde Q Q\) is therefore invertible. Put
\[
 U_l=\widetilde Q Q\,(Q\widetilde Q Q)^{-1/2}.
\]
The identity \(\widetilde Q^2=\widetilde Q\) gives \(U_l^*U_l=1\)
on \(Q\ell^2\), and \(U_l^*U_h=0\). The combined isometry
\((U_h,U_l)\) differs from the identity by \(O(L^{-1})\), so it
is onto for large \(L\). In particular, applying \(U_l^*\) to
the residual equation for \(U_h\) shows that the block from high
inputs to low outputs is \(O(L^{-k_*})\). The reverse block is
its adjoint, since the generator in this unitary frame is Hermitian.
Commuting \(D_x^2\) through a high column of \(T-P\) costs \(O(1)\)
in Schur size, since the factor \(2n(n'-n)+(n'-n)^2\) multiplies an
\(O(1/n)\) smooth envelope. Product rules and the same binomial series
bound the commutator with the full extension by \(O(1)\).
In the new basis the scalar operator, less \(D_x^2\) and including
the nominal connection, is Hermitian. Its high columns, and hence high
rows, are the real diagonal \(d_n\) plus \(O(L^{-k_*})\). The remaining
bounded Hermitian matrix acts solely on \(|n|<L\).

Denote the full scalar column map in the two copies by
\(\mathcal U\), with the conjugate construction in the second copy.
Then \(\mathcal S\coloneqq e^{J_{\rm off}}\mathcal U\) maps the resulting
coordinates back to the original field. In an equation
\(i\dot w=\mathcal L w\), the substitution \(w=\mathcal S a\)
gives the generator
\(\mathcal S^{-1}\mathcal L\mathcal S-i\mathcal S^{-1}\dot{\mathcal S}\).
Both \(\mathcal S-1\) and its old-parameter derivatives have Schur size
\(O(L^{-1})\). Along nonnominal old parameters the difference of the actual
and nominal connections has size
\begin{equation}\label{eq:oldconnection}
 C L^{-1}\bigl(|\dot I_F|+|\dot\theta_F-\dot\theta_F^\infty|\bigr).
\end{equation}
In particular variations of the order-one low Hermitian matrix may be
retained without any loss on the slow time scale.

It remains to decide which terms of order \(r\) are resonant.
Use the Fourier series \(\sum_{b,\nu}v_{b,\nu}e^{ibx-i\nu t}\) of the
old axial profile with its retuned old actions and nominal angles.
Here \(b\in\mathbb Z\) is a spatial frequency and
\(\nu\in Q_F^{-1}\mathbb Z\) a time frequency, with \(Q_F\) a common
denominator of the old occupied frequencies in the mass gauge. Thus
\(2\pi Q_F\) is an old period. The coefficient \(v_{b,\nu}\) is the
sum of all angular harmonics with that space-time frequency. A displacement \(O(r^2)\) of the actions
changes the summed coefficients by \(O(r^2)\) in every fixed Fourier
moment used here. Indeed the finite-angle parametrization is analytic
into each fixed Sobolev space on a complex angular strip of positive
width. For a multi-index \(l=(l_j)_{j\in F}\in\mathbb Z^F\), the
angular harmonic \(e^{-i l\cdot\theta}\) has the single spatial
frequency \(b=\sum_{j\in F}j l_j\), by translation covariance, and
time frequency \(\nu=\sum_{j\in F}l_j(\omega_j(I^\infty)+c_R)\).
Its coefficient decays exponentially in \(|l|_1\coloneqq \sum_j|l_j|\).
Truncation at \(|l|_1\le C_A\log R\)
leaves an \(O(R^{-A})\) tail in Schur and multiplier norms, with
\(|b|+|\nu|\le C'_A\log R\) in the truncation. Choose \(A\) larger
than every fixed power needed below. Changing the old phases and
actions to their actual values in a term already of order \(r\) costs
\(O(r(|\Theta|+r|J|))\).

\emph{Resonances of the new axial planes.} Pairing the terminal plane at \(K_i\) with \(v_{b,\nu}\) gives an
opposite-copy entry between \(n\) and \(n'\coloneqq K_i+b-n\). Here and below a bare divisor uses the free root frequency
\(K_i^2\), with its shift \(\chi_{K_i}(I^\infty)\) treated separately.
For this entry the divisor is
\[
 D\coloneqq n^2+(n')^2+2(qR)^2-K_i^2-\nu .
\]
If \(D\ne0\), then \(|D|\ge Q_F^{-1}\). If an endpoint has
\(|n|<L\), write \(n'=K_i+(b-n)\). Then
\[
 D=2q^2R^2+2K_i(b-n)+(b-n)^2+n^2-\nu .
\]
Since \(|b|+|\nu|=O(\log R)\) and \(|K_i|=O(R)\), all terms after
\(2q^2R^2\) sum to \(o(R^2)\), uniformly in \(q\). Hence
\(|D|\gtrsim R^2\) on every edge touching this low subspace. The
same calculation with \(n,n'\) interchanged covers \(|n'|<L\). On a zero-divisor edge \(q\) is bounded and both
\(|n|\) and \(|n'|\) are between fixed positive multiples of \(R\).
Indeed, substituting \(n'=K_i+b-n\) gives
\[
  2\big(n-(K_i+b)/2\big)^2+2(qR)^2
        =K_i^2+\nu-(K_i+b)^2/2 .
\]
The right-hand side is \(O(R^2)\), so \(n/R\) and \(q\) are bounded.
Along a subsequence \(q\) is a fixed nonzero integer and \(n/R\)
has a limit \(X\). Dividing by \(R^2\) gives \(q^2=X(P_i-X)\),
so \(X\) and \(P_i-X\) are both nonzero. In particular
\(P_i=X+q^2/X\). Since the
\(P_i\) are distinct, a row cannot serve two roots for arbitrarily large
\(R\). For the same root, changing \(b\) by a nonzero integer \(h\)
changes \((n')^2\) by \(h(2n'+h)\). Its modulus is at least
\(cR\) since \(|n'|\gtrsim R\) and \(|h|=O(\log R)\). Two retained
temporal labels can differ by only \(O(\log R)\). Once \(b\) is fixed,
\(D=0\) determines \(\nu\), and all harmonics with these same labels
have already been added into \(v_{b,\nu}\). The argument applies at
either endpoint. Thus the exact resonances form disjoint pairs.

The same-copy terms and their adjoints take \(n\) to \(n+K_i-b\).
Their bare divisor is the difference
\[
 (n+K_i-b)^2-n^2-(K_i^2-\nu)=2(K_i-b)(n-b)+\nu-b^2 .
\] It has modulus
\(\gtrsim R\) unless \(n=b\). Keep every entry \(b\leftrightarrow K_i\)
and its adjoint as a Hermitian coupling, with their actual time
oscillations. The vertex \(b\) cannot meet an opposite-copy resonance.
Intersection at \(K_i\) would require \(P_i(P_j-P_i)=q^2\) in the limit,
excluded by Proposition~\ref{prop:tree}. This separates all opposite-copy
pairs from the low Hermitian matrix and the retained same-copy terms.

For the discarded entries solve the homological equation using the bare
divisors above. More explicitly, an entry \(B_{uv}e^{-i\zeta t}\) in the
generator with free energies \(E_u,E_v\) is cancelled to first order by
the entry \(-B_{uv}e^{-i\zeta t}/(E_u-E_v-\zeta_{\rm bare})\) in the
exponent of the change of variables. The omitted root shift
\(\zeta-\zeta_{\rm bare}\) is \(O(R^{-2}+r^2)\). Each Fourier harmonic
defines a partial shift or reflection with at most one entry per row
and column, so the summed generator has Schur size \(O(r)\). Its terms
touching the low subspace have size \(O(r/R)\) for same-copy entries
and \(O(r/R^2)\) for opposite-copy entries. Its commutator with the
bounded low matrix therefore costs \(O(r/R)\), and that with the high
real diagonal costs \(O(r L^{-2})\). The terminal planes oscillate at frequencies of order \(R^2\).
To keep those frequencies out of the error bound, combine the free
commutator with the time derivative before estimating either one.
Write the change of variables as \(\mathcal M\coloneqq e^{\mathcal B(t)}\),
where \(\mathcal B\) is the exponent just constructed. Differentiating
the conjugation with respect to \(\theta\) gives the exact identity
\[
 \mathcal M^{-1}D_{\rm free}\mathcal M
   -i\mathcal M^{-1}\partial_t\mathcal M-D_{\rm free}
 =\int_0^1e^{-\theta\mathcal B}
   \bigl([D_{\rm free},\mathcal B]-i\partial_t\mathcal B\bigr)
   e^{\theta\mathcal B}\,d\theta .
\]
On a Fourier entry the factor in parentheses multiplies the exponent's
coefficient by its actual divisor. It equals the negative of the
discarded entry up to the root shift estimated above, and has total
Schur size \(O(r)\). Since \(\mathcal B\) itself has size \(O(r)\),
conjugation inside the integral changes this factor by \(O(r^2)\).
This gives \(\mathcal M=1+O(r)\) and the required higher-order bound. Conjugation of the order-\(r\) terms by the old frame costs
\(O(r/L)\).

\emph{The mixed boundary conditions.} These constructions give a bounded inverse for the transverse
inhomogeneous problem with mixed initial and terminal data. On the
complement of the opposite-copy pairs, the free diagonal, the old low
matrix, the real \(d_n(t)\) and the retained same-copy terms generate
an isometry in \eqref{eq:flat-weight-ins}. Indeed each coupling within a
given \(q\) and copy touches only indices where the weight equals
\((1+|q|)^p\). This holds on \(|n|<L\), at \(n=b\) with
\(|b|=O(\log R)\), and at \(n=K_i\) with
\(\operatorname{dist}(K_i,R\mathbb Z)\le1\). Every other entry is real diagonal. This holds for the
old matrix evaluated at the actual, time-dependent parameters. On
each opposite-copy pair remove the two free phases. Let
\(\delta_i\coloneqq \chi_{K_i}(I^\infty)\). Since the bare divisor vanishes, the
upper-right entry in the equation for \(i\partial_t\) is now
\(r\eta e^{-i\delta_i t}\), where \(\eta\coloneqq 2\gamma(c_i+e_i)v_{b,\nu}\) is the coupling coefficient of this pair.
The lower-left entry is its negative conjugate. The coordinate change
\(\operatorname{diag}(e^{-i\delta_i t/2},e^{i\delta_i t/2})\) removes
these phases and adds a diagonal of size \(|\delta_i|/2\) in physical
time. Thus in slow time \(\tau=rt\) the principal matrix for
\(\partial_\tau\) is
\[
 \begin{pmatrix}0&-i\eta\\i\bar\eta&0\end{pmatrix}.
\]
It is Hermitian and its square equals \(|\eta|^2\) times the identity,
so its eigenvalues are \(\pm|\eta|\). To compare the two endpoint
weights, use \(n'=RP_i+1+b-n\). The distances of \(n\) and \(n'\)
from \(R\mathbb Z\) differ by at most \(|1+b|=O(\log R)=o(L)\).
Their scaled horizontal frequencies stay bounded and their vertical
indices have the same absolute value. Thus the weights in
\eqref{eq:flat-weight-ins} have bounded ratio in both directions,
and the ordinary orthogonal spectral projections are bounded in \(Z\).
Discarding the half-angle diagonal correction and the \(d_n\)'s on
these pairs costs \(O((R^{-2}+r^2)/r)\) in slow time.

We collect the losses when the old parameters follow
their nominal motion. After division by \(r\), the errors from
removing the old off-diagonal pair and constructing the scalar
columns are bounded by \(C(R^{-2}+L^{-k_*})/r\).
The error from conjugating the terms of order \(r\) by the old frame
costs \(C/L\). The mixed conjugation contributes at most
\(C(R^{-1}+L^{-2}+r+R^{-2}+R^{-A})\) in the same units.
Together with the discarded diagonal on the resonant pairs these
are \(O(L^{-1})\), since \(R^{-2}/r=O(L^{-2})\) and
\(L^{-k_*}/r=O(L^{-2})\) by \eqref{eq:parameters-ins}.
When the old parameters move away from the nominal motion, one
must also keep \eqref{eq:oldconnection} and the changes of the
coefficients caused by the old action and angle displacements.

The set of exact pairs can change with \(R\). To choose one time
weight for all large \(R\), take the coefficients \(v^{\,0}_{b,\nu}\)
of the original periodic field \(V\), before retuning. For every
\(a>0\) the set of numbers
\[
       2|c_i v^{\,0}_{b,\nu}|,\qquad
       i\text{ an axial vertex},\quad
       (b,\nu)\in\mathbb Z\times Q_F^{-1}\mathbb Z ,
\]
contains only finitely many values above \(a\). This follows from
absolute summability and the finiteness of the tree. Choose \(\sigma>0\)
below the model decay rate in \(\tau\) and at positive distance from
this set and zero. Each rate on a retained pair differs from its
listed value by \(O(r^2+R^{-10}+R^{-A})\), by the retuning and
angular tail estimates. Thus the distance of all signed rates from
\(-\sigma\) has a positive lower bound independent of \(R\). For any \(T<S\) measure
inhomogeneities in the norm
\[
 \|f\|_{\sigma,Z}\coloneqq \sup_{T\le\tau\le S}e^{\sigma(\tau-T)}\|f(\tau)\|_Z.
\]
Write \(\mathcal U_0(\tau,u)\) for the propagator of this principal
system, including the changes of phase on the pairs, from slow time
\(u\) to \(\tau\). Denote by \(\Pi_f(u)\) the projection at time \(u\)
onto the decaying pair directions whose exponents are less than
\(-\sigma\), and put \(\Pi_c(u)\coloneqq 1-\Pi_f(u)\). These projections commute
with propagation in the sense that
\(\mathcal U_0(\tau,u)\Pi_f(u)=\Pi_f(\tau)\mathcal U_0(\tau,u)\).
For forcing \(f\) in the equation for \(\partial_\tau\), define
\[
 (\mathcal G_\perp f)(\tau)\coloneqq 
 \int_T^\tau \mathcal U_0(\tau,u)\Pi_f(u)f(u)\,du
 -\int_\tau^S\mathcal U_0(\tau,u)\Pi_c(u)f(u)\,du .
\]
This prescribes zero fast stable component at \(T\) and zero complementary
component at \(S\). A scalar exponent \(\mu\) contributes the time-weighted
factor \(e^{(\mu+\sigma)(\tau-u)}\). The first integral has \(\mu+\sigma<0\)
and \(\tau-u\ge0\), while the second has \(\mu+\sigma>0\) and
\(\tau-u\le0\). The isometric part obeys the second estimate with
\(\mu=0\). Since all the exponents stay a fixed distance from
\(-\sigma\), both weighted kernels are bounded in operator norm by
\(C e^{-\delta_{\rm gap}|\tau-u|}\) for some \(\delta_{\rm gap}>0\)
independent of \(R,T,S\). Here \(\delta_{\rm gap}\) can be any positive lower bound for
\(\sigma\) and the distances of the signed pair rates from \(-\sigma\).
Taking the Hilbert direct sum of the blocks before time integration gives
\(\|\mathcal G_\perp f\|_{\sigma,Z}\le 2C\delta_{\rm gap}^{-1}
\|f\|_{\sigma,Z}\), uniformly in the length of the interval.
The placement of the two traces also gives uniqueness. On a fast
direction the homogeneous solution is determined by its value at
\(T\), and on each remaining direction it is determined by its
value at \(S\). On the half-line the backward integral converges
with the same kernel bound. In that case membership in the
exponentially weighted space excludes every homogeneous
complementary direction, including the isometric part. On a fixed
bounded slow interval the same principal propagator is bounded in
both time directions. Large unitary rotations
used to obtain the pairwise estimate are needed only to define these
propagators and projections.

\subsection{Resonant forcing and the nonlinear terminal problem}\label{sec:ins-terminal}

\emph{Carrier phases and the prescribed forcing.}
We work first on a finite slow-time interval \([T,S]\), with bounds
independent of \(S\). There are three parts to the construction.
We choose phases so that each selected interaction has the time phase
of its output. We then solve the free inhomogeneous equation for
the remaining prescribed oscillations. Finally we estimate the errors
in the transverse field, the new axial coordinates, and the old actions
and angles. Only the last step uses a boundary value problem.

We first choose the time phase of each carrier. A transverse wave
has a decaying amplitude, but a new axial wave survives at positive
infinity. We therefore use its exact terminal frequency, including
the shift due to all occupied axial modes. Put
\(\kappa_v\coloneqq |k_v|^2+\varepsilon_v\), with
\(\varepsilon_{(P_i,0)}\coloneqq \chi_{K_i}(I^\infty)\) at every axial vertex.
For each selected triangle \(v+v'=u\) impose
\(\varepsilon_v+\varepsilon_{v'}=\varepsilon_u\). These conditions
have a common solution. Starting from the root, split the assigned
shift of each primary parent \(m=p+q\) between its children in any
fixed way, for instance \(\varepsilon_p=0\) and
\(\varepsilon_q=\varepsilon_m\). When \(m\) is off the axis its
split contains three additional triangles. Their equations read
\[
 \varepsilon_{J(p)}=\varepsilon_{A(p)}-\varepsilon_p,\qquad
 \varepsilon_{J(q)}=\varepsilon_{A(q)}-\varepsilon_q,\qquad
 \varepsilon_{T_m}=\varepsilon_{J(p)}+\varepsilon_{J(q)} .
\]
Here \(A(p),A(q)\) are on the axis and their shifts have already been
specified. The auxiliary vertices are distinct at different splits
by Proposition~\ref{prop:tree}. This recursion gives
\(\varepsilon_v=O(R^{-2}+r^2)\) on the fixed tree.

Write \(\mathbf x=(x,y)\), \(\tau=rt\), and separate the known waves
into their transverse and axial parts. In the following formulas we identify the axial polynomial \(Y\) with
its sequence of Fourier coefficients.
\begin{align*}
 W(t,\mathbf x)
   &\coloneqq \gamma\sum_{v_y\ne0}x_v(2|g|\tau)
                    e^{ik_v\cdot\mathbf x-i\kappa_v t},\\
 Y(t,x)
   &\coloneqq \gamma\sum_i\bigl(x_{(P_i,0)}(2|g|\tau)+e_i\bigr)
                         e^{iK_ix-i\kappa_{(P_i,0)}t}.
\end{align*}
The first sum is over the transverse vertices of \(\mathcal V\).
Here and below \(Y_\infty(t)\) is obtained from \(Y(t)\) by replacing
each \(x_{(P_i,0)}\) by its limit \(c_i\) and keeping its phase
\(e^{-i\kappa_{(P_i,0)}t}\). Thus the axial comparison at infinity
is to a time-dependent solution. The retuning errors \(e_i\) put the limiting axial amplitudes at their
prescribed values in \(Y_\infty\). By Lemma~\ref{lem:tree-orbit}, for
some \(\lambda>\sigma\) and every \(\epsilon>0\) one can choose \(T\)
so that the amplitudes of \(W,Y-Y_\infty\) and their first two
derivatives in \(\tau\) are bounded by
\(\epsilon e^{-\lambda(\tau-T)}\) for \(\tau\ge T\).
Such an amplitude derivative keeps the rapid carrier phase fixed.

In expanding \(|V_*+r(Y+W)|^2(V_*+r(Y+W))\), the coefficient of \(r^2\)
is \(2V_*|Y+W|^2+\overline{V_*}(Y+W)^2\). Terms involving only \(Y\)
are axial and are already included in the Birkhoff frequencies.
Terms involving just one transverse factor also have zero axial
projection. Thus the order-\(r^2\) transverse forcing and the
axial forcing after subtracting the mean transverse density are given by
\begin{align}
 Q_w&\coloneqq P_\perp\{2V_*(|Y+W|^2-|Y|^2)
                     +\overline{V_*}((Y+W)^2-Y^2)\},\notag\\
 Q_z&\coloneqq P_0\{2V_*(|W|^2-S_W)+\overline{V_*}W^2\},
 \qquad S_W\coloneqq \sum_{v:\,v_y\ne0}|x_v(2|g|\tau)|^2 .\label{eq:forcingpol-ins}
\end{align}
In fact \(S_W\) is the spatial mean of \(|W|^2\), by orthogonality of its
distinct Fourier modes. To handle this mean exactly, take \(2r^2S_W U\)
out of the full axial force \(f\) in \eqref{eq:axis-insertion}.
Gauge equivariance, differentiated along the curve
\(a\mapsto e^{-ia}U\) at zero, gives
\[
 D\Phi(U)[-2ir^2S_W U]=-2ir^2S_W z,\qquad
 2\operatorname{Re}\bigl(\bar z_j(-2ir^2S_W z_j)\bigr)=0 .
\]
Thus the mean changes each angular speed by \(2r^2 S_W\).
Its contribution on a new axial mode of size \(O(r)\) is
\(O(r^3S_W)\). The remaining known force has leading term \(r^2Q_z\).
Its outputs have nonzero horizontal velocity. Indeed a zero difference in
horizontal velocity would require the same vertex, whose diagonal
density has just been removed, and a zero sum would require opposite
horizontal velocities. Both statements use the distinct absolute
horizontal components in Proposition~\ref{prop:tree}. This observation
permits use of \eqref{eq:chart-head} before projecting to the old
coordinates. In extracting the prescribed harmonics of \(Q_w,Q_z\)
we use the nominal old angles. Replacing the actual \(V_*\) by that
profile costs \(C(|\Theta|+r|J|+r^2)\) times a decaying amplitude.
The low linear operator continues to use the actual old coordinates.

We determine the time frequency \(\zeta\) and spatial frequency \(k\)
of each remaining product. Put \(e=(1,0)\) and write \(v_x\) for the
horizontal component of \(v\). For two positive carriers and a conjugated
old harmonic \((b,\nu)\), one has \(k=(2-b)e+R(v+v')\) and
\(\zeta=\kappa_v+\kappa_{v'}-\nu\). Direct expansion gives
\[
 \zeta-|k|^2=-2R^2v\cdot v'+2R(b-1)(v_x+v'_x)
       +2-\nu-(2-b)^2+\varepsilon_v+\varepsilon_{v'} .
\]
For a carrier times a conjugated carrier and a positive old harmonic,
\(k=be+R(v-v')\), \(\zeta=\kappa_v-\kappa_{v'}+\nu\), and
\[
 \zeta-|k|^2=2R^2(v-v')\cdot v'+2R(1-b)(v_x-v'_x)
             +\nu-b^2+\varepsilon_v-\varepsilon_{v'} .
\]
When the coefficient of \(R^2\) is nonzero it is an integer of modulus
at least one. If it vanishes, Proposition~\ref{prop:tree} puts the
nonconstant output velocity in the tree, with nonzero horizontal component.

When \(b\ne1\) the linear term then dominates the logarithmically
bounded remainder. When \(b=1\), the only desired temporal harmonic is
\(\nu=1\). All others have separation at least \(1/Q_F\) before detuning.
Consider the surviving harmonic \((b,\nu)=(1,1)\) on one
triangle \(p+q=h\). The assigned frequency corrections give
\(\kappa_p+\kappa_q-\kappa_h=1\), so the product at each output has
exactly that carrier's phase. After removing \(r^2\) and the carrier,
the two orderings at the parent and the density term at a child have
coefficients
\[
 \begin{aligned}
 2\bar g\,\gamma^2 x_p x_q
       &=2i|g|\gamma\,x_p x_q,\\
 2g\,|\gamma|^2 x_h x_q
       &=-2i|g|\gamma\,x_h x_q ,
 \end{aligned}
 \qquad \gamma=i g/|g| .
\]
On the left of the evolution equation the amplitude derivative is
\(2i|g|\gamma\,x_v'\), with the prime taken in \(z=2|g|rt\).
The first line gives the parent equation and the second gives the
child equation in \eqref{eq:single-triangle}. Their contributions add
at shared vertices. The change in \(g\) from old-action retuning is
\(O(r^2)\), and the omitted angular tail has the bound above.
Every occurrence of a new radius error \(e_i\) in these products is
paired with a decaying transverse amplitude.

There is also a linear term at each surviving axial index \(K_i\).
Its carrier has the exact frequency
\(\kappa_{(P_i,0)}=\omega_{K_i}(I^\infty)+c_R\).
When the actual diagonal operator acts on that carrier, the remaining
frequency difference is \(\omega_{K_i}(I)-\omega_{K_i}(I^\infty)\).
It tends to zero with \(I-I^\infty\), at the rate estimated below.
Using the free frequency \(K_i^2\) here would leave the constant
terminal shift multiplying a nondecaying wave.

For each remaining retained harmonic \(a(\tau)e^{ik\cdot(x,y)-i\zeta t}\)
in the nominal forcing put \(E_{k,\zeta}\coloneqq e^{ik\cdot\mathbf x-i\zeta t}\)
and use
\[
 (i\partial_t+\Delta)
       \frac{a(\tau)E_{k,\zeta}}{\zeta-|k|^2}
 =a(\tau)E_{k,\zeta}
       +\frac{ir\,(\partial_\tau a)(\tau)E_{k,\zeta}}{\zeta-|k|^2}.
\]
For example, an order-\(r^2\) source of this form in the physical
equation is canceled by adding \(r^2aE_{k,\zeta}/(\zeta-|k|^2)\)
to the prescribed field. Its uncanceled derivative has size
\(r^3|\partial_\tau a|/|\zeta-|k|^2|\). The field correction itself
has size \(r\) in the original variables per unit of normalized
error, and \(\partial_t=r\partial_\tau\) on its slow amplitude.
The normalized error equation therefore divides the physical residual
by \(r^2\), leaving one factor of \(r\) in this derivative term. Write \(P_w,P_z\) for the resulting sums, where on the axis
we only take indices outside \(F\). For large \(R\) all retained axial
outputs here have \(|n|\gtrsim R\). The absolute divisors are at least
\(1/(2Q_F)\). A transverse output with zero horizontal velocity must have
a nonzero \(R^2\) coefficient, so its particular carries \(O(R^{-2})\).
This pays for applying the entire bounded old low matrix. Every other
particular lies at \(|n|\gtrsim R\), where the real diagonal shifts are
\(O(R^{-2})\). On the axis the same bound holds for the empty-site
renormalized frequency of the smooth old state. Its further change at
the actual actions is \(O(r|J|+r^2)\) in the gauged frequency. The scalar
mass discrepancy is included in this last estimate. Thus the bounded
divisors are used in the dressed transverse coordinates and in the
true axial coordinates. The \(B\) norms of the particular sums and of the
sums with one slow amplitude derivative are bounded by
\(C\epsilon e^{-\lambda(\tau-T)}\), by absolute old-harmonic summation.
This construction leaves errors of size \(O(r)\) after division by \(r^2\)
when a slow derivative falls on the particular amplitudes.

\emph{Equations for the four errors.}
We now seek the actual field as the sum of these prescribed terms
and an error. Write \(F^c\coloneqq \mathbb Z\setminus F\). Let \(h\) be the
transverse error after both coordinate changes, and \(b=(b_j)_{j\in F^c}\)
the axial error in Birkhoff coordinates. Each error is divided by \(r\).
With \(B_w,B_z\) denoting the prescribed profiles, set
\begin{equation}\label{eq:reconstruct-ins}
 w_\perp=\mathcal S(B_w+r\mathcal M h),\qquad
 z|_{F^c}=B_z+r b,\qquad
 B_w\coloneqq rW+r^2P_w,\quad B_z\coloneqq rY+r^2P_z .
\end{equation}
In the first identity the left side denotes the doubled transverse
field described above, and \(W,P_w\) represent their doubled versions.
Every transverse index and every axial index in \(F^c\) is retained. In particular \(\mathcal M\)
only acts on the unknown \(h\). The terms of the mixed operator on the
known fields were already included in \(Q_w\).

For reference we specify the exact equations that the terminal integral
map solves. Reconstruct the old coordinates from \eqref{eq:old-mod-insertion}
as \(z_j=\sqrt{I_j^\infty+rJ_j}\,\exp[-i(\theta_j^\infty+\Theta_j)]\),
and put \(U=\Psi(z)\). Let
\(\Lambda\coloneqq \operatorname{diag}(\omega_j(I)+c_R)\),
\(D_f\coloneqq D\Phi(U)[-if]\), and let \(\mathcal L(U)\) be the doubled linear
generator in \eqref{eq:trans-insertion}. Denote the doubled nonlinear
term by \(\mathcal N\coloneqq (N,-\overline N)\), with the appropriate vertical
components taken in each copy. Set
\(\mathcal K\coloneqq \mathcal S^{-1}(\mathcal L(U)\mathcal S-i\dot{\mathcal S})\).
Write \(\Lambda_{F^c}\) for the restriction of \(\Lambda\) to \(F^c\). All four equations below have the same slow time on the left.
The factors of \(r\) can be seen, for example,
from \(\partial_t(B_z+rb)=\dot B_z+r^2\partial_\tau b\).
For \(j\in F\), \(\dot I_j=2\operatorname{Re}(\bar z_j\dot z_j)\)
and \(\dot\theta_j=-\operatorname{Im}(\dot z_j/z_j)\).
Substituting the reconstruction into the equations therefore gives
\begin{align}
 i\partial_\tau h
   &=r^{-1}\mathcal M^{-1}(\mathcal K\mathcal M-i\dot{\mathcal M})h
      +r^{-2}\mathcal M^{-1}
             (\mathcal K B_w-i\dot B_w+\mathcal S^{-1}\mathcal N),
                  \label{eq:exact-h-ins}\\
 \partial_\tau b
   &=-ir^{-1}\Lambda_{F^c}b
       +r^{-2}(-i\Lambda_{F^c}B_z-\dot B_z+(D_f)|_{F^c}),\notag\\
 \partial_\tau J_j&=2r^{-2}\operatorname{Re}(\bar z_j(D_f)_j),
                  \qquad j\in F,\notag\\
 \partial_\tau\Theta_j
   &=r^{-1}\bigl(\omega_j(I)-\omega_j(I^\infty)
                         -\operatorname{Im}((D_f)_j/z_j)\bigr),
                         \qquad j\in F .\notag
\end{align}
On a trial path evaluate \(\dot{\mathcal S}\) by the chain rule using
the instantaneous physical velocities
\(2\operatorname{Re}(\bar z_j(D_f)_j)\) and
\(\omega_j(I)+c_R-\operatorname{Im}((D_f)_j/z_j)\).
This prescription only uses values of the trial path. At a fixed point
these are its actual velocities, so that \eqref{eq:exact-h-ins} is
equivalent to \eqref{eq:axis-insertion}--\eqref{eq:trans-insertion}.
The dots on \(\mathcal M,B_w,B_z\) are ordinary physical-time derivatives
of prescribed functions.

For a normed space \(D\) put
\(\|a\|_{\sigma,D}\coloneqq \sup_{T\le\tau\le S}e^{\sigma(\tau-T)}\|a(\tau)\|_D\).
Set
\[
 H_1\coloneqq \|h\|_{\sigma,Z},\quad A_1\coloneqq \|b\|_{\sigma,Z},\quad
 J_1\coloneqq \|J\|_{\sigma,\mathbb R^{F}},\quad
 L_1\coloneqq \|\Theta\|_{\sigma,\mathbb R^{F}},\quad X\coloneqq H_1+A_1+J_1+L_1 .
\]
Here \(\mathbb R^F\) carries its Euclidean norm. For a vanishing trial
error the residual is bounded in these norms by a constant times
\begin{equation}\label{eq:defect-scale-ins}
 e_R\coloneqq L^{-1}+R^{-1}+R^{-2}/r+r(1+R^{2-p})+L^{-k_*}/r
       =O(R^{-\alpha}).
\end{equation}
Indeed \(R^{-2}/r=a_0^{-1}R^{-(2-s)}\),
\(rR^{2-p}=a_0R^{-(p+s-2)}\), and \(L^{-k_*}/r\le a_0^{-1}R^{-2\alpha}\).
The scale \(L^{-1}=R^{-\alpha}\) bounds their sum along with the
remaining terms.

On \([T,S]\) we define a map from trial errors to new errors. Evaluate
the coefficients and remainders in \eqref{eq:exact-h-ins} at the trial
path. Solve for the new \(h\) using the transverse Green operator,
keeping the Hermitian low matrix at its trial values. For the new
\(b\), retain its exact real diagonal and prescribe \(b(S)=0\).
Thus a source \(f_j(\tau)\) in its equation contributes, at
\(j\in F^c\), the integral
\[
 -\int_\tau^S
 \exp\!\left(\frac{i}{r}\int_\tau^u
                  [\omega_j(I(v))+c_R]\,dv\right) f_j(u)\,du .
\]
Every diagonal exponential has modulus one. Taking the sequence
norm before integrating and then the time weight gives the factor
\[
 \int_\tau^S
 e^{\sigma(\tau-T)} e^{-\sigma(u-T)}\,du
 =\frac{1-e^{-\sigma(S-\tau)}}{\sigma}\le\sigma^{-1}.
\]
This bound is uniform in \(S\), even if an axial frequency vanishes.
The old actions and angles also use zero terminal errors and the
same integral bound. For the transverse error the fast stable
components instead have zero trace at \(T\), as prescribed by
\(\mathcal G_\perp\). These conditions determine which components
are integrated in each time direction. Write
\(H_1',A_1',J_1',L_1'\) for the norms of the four outputs of this
map. For \(X\le1,\ \sqrt L X\le1\), they satisfy
\begin{align}
 H_1'&\le C(e_R+\epsilon X+\sqrt L X^2+o(1)H_1),\notag\\
 A_1'&\le C(e_R+J_1+\epsilon X+\sqrt L X^2+o(1)A_1),\notag\\
 J_1'&\le C(e_R+\epsilon H_1+\epsilon^2 A_1+\sqrt L X^2),\label{eq:four-rows-ins}\\
 L_1'&\le C(e_R+J_1+r(A_1+H_1)+\sqrt L X^2).\notag
\end{align}
The constants here depend on the fixed background, tree and weight
\(\sigma\). They are uniform in \(R\) sufficiently large, in \(S>T\),
and in the starting time for the small tail. The \(o(1)\) coefficients
tend to zero as \(R\to\infty\) for these fixed choices.

For the proof put \(E\coloneqq e^{-\sigma(\tau-T)}\). First isolate the
error of the axial chart on the prescribed waves. Define
\[
 \mathfrak e_0\coloneqq r^{-1}\big[\Psi(z|_F+rY)-V_*\big]-Y,\qquad
 \mathfrak e_1\coloneqq r^{-1}\big[U-\Psi(z|_F+rY)\big],
\]
using Fourier series on the field side. In the transverse variables set
\(W_1\coloneqq (\mathcal S-1)W+r\mathcal S P_w\) and \(h_1\coloneqq \mathcal S\mathcal M h\).
Then
\[
 U=V_*+rY+r \mathfrak e_0+r \mathfrak e_1,\qquad w_\perp/r=W+W_1+h_1,
\]
where
\[
 \|\mathfrak e_0\|_Z\le Ce_R,\quad \|\mathfrak e_1\|_Z\le C(A_1+r\epsilon)E,\qquad
 \|W_1\|_B\le C(L^{-1}+r)\epsilon E,\quad
 \|h_1\|_Z\le C H_1 E .
\]
Indeed \(\mathfrak e_0\) is the Taylor error in replacing
\(\Psi(z|_F+rY)\) by \(V_*+rY\), divided by \(r\). Equations
\eqref{eq:chart-tube}--\eqref{eq:chart-head} at the actual old state
give its bound. For \(\mathfrak e_1\) use the bounded full chart differential
on the \(O(r)\) tube. For \(W_1,h_1\) use the two Schur bounds in
\eqref{eq:reconstruct-ins}. All prescribed Fourier polynomials have
uniform multiplier norms at their indicated amplitude.

\emph{Transverse field.} Put \(v\coloneqq U-V_*\).
The additional linear terms on a physical field \(w\) are exactly
\[
 2(\overline{V_*}v+V_*\bar v)w+2V_*v\bar w+
                            2|v|^2 w+v^2\bar w .
\]
Their order-\(r^2\) contributions on \(W\) have already been
included in \(Q_w\). Their order-\(r\) coefficients on \(h\) from
\(Y_\infty\) are precisely the mixed terms of the preceding
subsection. Replacing \(Y_\infty\) by \(Y\) adds \(O(\epsilon)H_1\)
in slow units. The factor \(\mathfrak e_0\) against the known multiplier \(W\)
costs \(O(e_R)\). Against a general error it costs
\(O(\sqrt L\,e_R)H_1=o(1)H_1\).
The factor \(\mathfrak e_1\) against \(W\) costs
\(O(\epsilon(A_1+r\epsilon))\), while against \(h_1\) it costs
\(O(\sqrt L(A_1+r\epsilon)H_1)\).
For example, two transverse unknowns and a smooth old factor satisfy
\[
 r^{-2}\|P_\perp(V_*|r h_1|^2)\|_Z
          \le C\sqrt L\, H_1^2 E^2 ,
 \qquad
 \|W h_1\|_Z\le C\epsilon H_1 E^2 .
\]
In the second inequality the known \(W\) is a \(B\) multiplier. The
different estimates for these two products explain the single loss
\(\sqrt L\) at quadratic order. Terms with two factors of \(v\)
have an additional \(r\).
In \(N(U,w)\) the two known transverse factors paired with \(V_*\)
are included in \(Q_w\). A single unknown pairs with a decaying
known factor, and two unknowns cost \(O(\sqrt L X^2)\).
The cubic all-error term is \(O(r L X^3)\). These bounds follow by
applying Lemma~\ref{lem:product-insertion} once for each additional
arbitrary factor, using \(B\) for the known ones.

The old-frame residual on a known transverse mode costs
\(O((R^{-2}+L^{-k_*})/r)\). Its horizontal index is comparable to \(R\),
so its own scalar diagonal is only \(O(R^{-2})\).
The scalar \(2\sum_F I_j+c_R\) on the known modes costs
\(C\epsilon(J_1+r)\), and on \(h\) costs \(C(J_1+r)H_1\).
The extra connection \eqref{eq:oldconnection}, divided by \(r\), is small. The true velocities at the trial point give
\(|\dot I_F|+|\dot\theta_F-\dot\theta_F^\infty|\le C(rJ_1 E+r^2)\)
on our tube. To see the \(r^2\) in this estimate, use
\(\|w_\perp\|_Z\le Cr(\epsilon+H_1)E\), expand products with
known multipliers first, and use \(\sqrt L X\le1\) on the remaining
pair of arbitrary factors. Terms due to freezing the mixed coefficients
at the nominal phases carry \(|\Theta|+r|J|\). The discarded nonresonant
homogeneous terms have slow operator norm
\(O((R^{-2}+L^{-k_*})/r+L^{-1}+r)\) in addition to these modulation
errors. These observations give the first line of \eqref{eq:four-rows-ins}.

\emph{Axial modes outside \(F\).} The \(\ell^1\) difference of the
actual and terminal actions is bounded by
\[
 C\bigl(rJ_1+r^2(\epsilon+A_1+r\epsilon)\bigr)E.
\]
Indeed on \(F\) we have \(I_j-I_j^\infty=rJ_j\). On \(F^c\) put
\(z^\infty\coloneqq rY_\infty\), with its true time phases, and use
\[
 \sum_{j\in F^c}\big||z_j|^2-|z_j^\infty|^2\big|
 \le \|z|_{F^c}-z^\infty\|_{\ell^2}
          \big(\|z|_{F^c}\|_{\ell^2}+\|z^\infty\|_{\ell^2}\big).
\]
The last factor is \(O(r)\), while the first contains \(r(Y-Y_\infty)\),
\(r^2P_z\), and \(rb\). To see the effect on one persistent
coordinate, write \(\Delta\omega_i=\omega_{K_i}(I)-\omega_{K_i}(I^\infty)\).
The mass formula and the local \(\ell^2\)-Lipschitz bound for \(\chi\)
give \(|\Delta\omega_i|\le C\|I-I^\infty\|_{\ell^1}\), uniformly in
\(K_i\). In the equation for \(b_{K_i}\) this mismatch contributes
\[
 \frac{|\Delta\omega_i|\,|rY_{K_i}|}{r^2}
 \le C\bigl[J_1+r(\epsilon+A_1+r\epsilon)\bigr]E ,
\]
since \(|Y_{K_i}|\) is bounded on the fixed tree. This explains why
the coefficient of \(J_1\) in the second row need not be small even
on the far tail. The part involving \(A_1\) carries an
additional \(r\). The diagonal on the error \(b\) itself is retained.
After removing the exact density rotation, \(f\) contains two transverse
factors or three. One unknown transverse factor times a known one costs
\(C\epsilon H_1\), two arbitrary ones \(C\sqrt L X^2\).
On the leading known terms the price for replacing \(D\Phi(U)\)
by \(D\Phi(V_*)\) is \(O(r(1+R^{2-p}))\), and that for replacing the
latter by Fourier is \(O(R^{-1}+L^{-1})\) on the high horizontal
outputs. All these bounds apply to the strong norm \(Z\).

\emph{Old actions and angles.} For each \(j\in F\) the density
rotation drops out of the action equation by
\(\operatorname{Re}(-i|z_j|^2)=0\).
The other known order-\(r^2\) terms have high horizontal index after
the old angular truncation, so their Fourier projection to \(F\)
vanishes. The same chart estimate gives \(O(e_R)\).
This works at the actual old angles, hence produces no order-one
coefficient of \(L_1\). An axial error in the old-coordinate forcing
has to act on a transverse pair. The latter either consists of two
known decaying factors, or has an additional unknown. Differentiating
the chart with respect to an axial change of size \(r A_1 E\) costs at
most \(C r(1+R^{2-p})A_1E\) by the proof of
\eqref{eq:chart-tube}. These estimates give the third line.
For the angular equation division is by \(r\), and the old \(z_j\)
are bounded away from zero for this fixed insertion. The frequency
shift just obtained gives \(C[J_1+r(\epsilon+A_1)]E\).
The remaining force gains an additional factor \(r\) compared with
the action equation. Since \(r\lesssim e_R\), this proves the last
line. Finally, the bounds for the particulars, their differentiated
amplitudes, and the omitted angular tails give contributions \(O(e_R)E\)
in each line. Every pure error evaluated in \eqref{eq:exact-h-ins}
has a decaying factor because the axial state at infinity is exactly \(V_R\).
In particular, with \(h=b=J=\Theta=0\), replacing \(Y\) by
\(Y_\infty\) and setting \(W=P_w=P_z=0\) reconstructs precisely
the axial Birkhoff solution \(V_R\). Substitution into all four
equations then gives zero residual. This verifies cancellation at the
terminal state itself, including on every persistent axial mode.

\smallskip
\noindent\emph{The terms of the four rows in detail.} We collect the
contributions to \eqref{eq:four-rows-ins} with their sizes. They are
pointwise bounds in the slow-time units of \eqref{eq:exact-h-ins}, and
each carries the factor \(E\). The transverse Green operator
\(\mathcal G_\perp\) and the backward integrals then multiply them by
\(2C\delta_{\rm gap}^{-1}\) and \(\sigma^{-1}\), respectively. Recall
\(\sqrt L\,e_R=O(R^{-\alpha/2})\) and \(\sqrt L\,r=a_0R^{\alpha/2-s}\).
\begin{itemize}
\item \emph{First row.} Pure sources: the old-frame residual on the
known modes, \(O((R^{-2}+L^{-k_*})/r)\); the carrier corrections
\(\varepsilon_v\), \(O((R^{-2}+r^2)/r)\); \(\mathfrak e_0\) against
\(W\), \(O(e_R)\); the terms with \(W_1\), \(O((L^{-1}+r)\epsilon)\); a
slow derivative of a particular and the products of three known
factors, \(O(r)\); the angular tail, \(O(R^{-A})\). Their sum is
\(O(e_R)\). Linear terms with a decaying known factor: \(Y-Y_\infty\)
in the mixed coefficients, \(O(\epsilon H_1)\); \(\mathfrak e_1\)
against \(W\), \(O(\epsilon(A_1+r\epsilon))\); the nominal old angles
and actions in the prescribed harmonics,
\(O(\epsilon(L_1+rJ_1+r))\); the scalar \(2\sum_FI_j+c_R\) on the known
modes, \(O(\epsilon(J_1+r))\). Linear terms in \(h\) with a small
coefficient: \(\mathfrak e_0\) against \(h_1\),
\(O(\sqrt L\,e_R)H_1\); the part \(r\epsilon\) of \(\mathfrak e_1\)
against \(h_1\), \(O(\sqrt L\,r\epsilon)H_1\); the parts of size \(r\)
of the scalar and of the connection \eqref{eq:oldconnection},
\(O(r)H_1\); and the discarded nonresonant homogeneous terms,
including the half-angle diagonals on the pairs, \(O(e_R)H_1\). This is
the coefficient \(o(1)=O(R^{-\alpha/2})\) of \(H_1\). Products of two
errors: \(J_1H_1\) from the scalar and the connection, \(A_1H_1\) from
\(\mathfrak e_1\) against \(h_1\), \(L_1H_1\) from the frozen phases,
and two transverse unknowns, each \(O(\sqrt L\,X^2)\). The cubic error
term is \(O(rLX^3)\le O(\sqrt L\,X^2)\), since \(r\sqrt L\,X\le1\).
\item \emph{Second row.} The mismatch \(\Delta\omega_i\) on the
persistent roots gives \(C[J_1+r(\epsilon+A_1+r\epsilon)]\), with a
coefficient of order one in front of \(J_1\). The force \(f\) gives
\(O(e_R)\) from two known factors, after the finite equation and the
particulars, \(O(\epsilon H_1)\) from one known and one unknown
transverse factor, and \(O(\sqrt L\,X^2)\) from two unknown ones.
The nominal angles in the prescribed harmonics give
\(O(\epsilon(L_1+rJ_1+r))\). The dependence of \(U\) and of
\(D\Phi(U)\) on \(b\) gives \(O(r)A_1\) and \(O(r(1+R^{2-p}))A_1\);
this is the coefficient \(o(1)\) of \(A_1\).
\item \emph{Third row.} The density rotation drops out. The known
two-wave terms project to \(F\) only through the chart corrections at
the actual old state, \(O(e_R)\), and therefore produce no term linear
in \(J_1\) or \(L_1\). A transverse unknown with a known factor gives
\(O(\epsilon H_1)\). An axial unknown enters through \(U\) and
\(D\Phi(U)\), multiplied by a transverse pair, giving
\(O(\epsilon^2A_1)\) with two known factors and a quadratic term
otherwise.
\item \emph{Fourth row.} The frequency shift gives
\(C[J_1+r(\epsilon+A_1+r\epsilon)]\), again with a coefficient of order
one in front of \(J_1\). The density rotation gives
\(2rS_W=O(r\epsilon^2)\). The remaining force is \(r\) times terms of
the type estimated in the third row, and \(r\lesssim e_R\).
\end{itemize}
Consequently the linear part of the map is triangular up to small
terms. The only linear terms with coefficients of order one lead from
\(J_1\) to \(A_1'\) and \(L_1'\), and the row for \(J_1'\) has no such
input. In particular \(L_1\) enters the other rows only with a factor
\(\epsilon\), through the prescribed harmonics, or quadratically.

\subsection{Passage to the half-line}\label{sec:ins-half}

\emph{Closing the four bounds.}
The only linear terms without a small coefficient on the right of
\eqref{eq:four-rows-ins} lead from \(J_1\) to \(A_1'\) and \(L_1'\).
Accordingly give the latter two components more room. Choose constants
\(M\ge1\) and \(M'\ge M\) with \(M'\ge 8CM\), increasing \(C\)
in \eqref{eq:four-rows-ins} to at least one. On the set
\[
 H_1,J_1\le M\rho,\qquad A_1,L_1\le M'\rho
\]
we have \(X\le 4M'\rho\). The terms \(CJ_1\) use at most
\(M'\rho/8\) of the corresponding bounds. Choose
\(4C\epsilon M'\le M/8\), with \(\epsilon\le1\), and then choose \(T\)
to ensure the required tail bound for the model. Keep these choices fixed.
Since
\[
 \frac{e_R}{\rho}=O(R^{-\alpha/4}),\qquad
          \sqrt L\,\rho=R^{-\alpha/4},
\]
taking \(R\) sufficiently large makes
\(C(e_R+\sqrt L X^2)=o(\rho)\), \(rX=o(\rho)\), and \(\sqrt L X\le1\).
The \(o(1)\) coefficients in \eqref{eq:four-rows-ins} are absorbed
at this last step. Thus the terminal map takes the displayed set into
itself for every \(S>T\).

\emph{Existence on a finite interval.}
To apply the compact fixed point theorem, fix \(R\) and \(S\). For
this fixed \(R\), the space \(Z\) is equivalent to \(H^p\) on the
chosen lattice, with the corresponding sequence interpretation on
the axis. Use continuous paths in \(H^s\), and ordinary continuous
paths for the two finite-dimensional variables, as the ambient space.
Intersect the four bounds above with a common Lipschitz bound in
\(H^{p-2}\) for the two field variables and in Euclidean norm for the
old variables. This intersection contains zero and is convex.
The image of the terminal map satisfies such a bound with a constant
depending on \(R\) and the tree. Indeed an output derivative in slow
time is its principal generator times the output plus its source.
The free generators divided by \(r\) map \(H^p\) to \(H^{p-2}\) at
fixed \(R\), and the remaining bounded operators and sources have
uniform bounds on each compact interval. In particular the local
bounds can be taken independently of the right endpoint \(S\).

For completeness, spatial tails beyond Fourier radius \(N\) of an
\(H^p\)-bounded path have \(H^s\) norm at most
\(\langle N\rangle^{s-p}\) times the \(H^p\) bound. Time increments
have \(H^s\) norm \(O(|\tau-\tau'|^{(p-s)/2})\) by interpolating
the uniform \(H^p\) bound and Lipschitz \(H^{p-2}\) bound. Thus
Arzel\`a--Ascoli gives relative compactness in the path topology.
Lower semicontinuity, applied first to finite Fourier sums, shows that
the four bounds and the Lipschitz constraint pass to the limit.
We have obtained a compact convex domain for the map.

We check continuity on this domain at fixed \(R,S\). Convergence in
\(C_\tau H^s\) suffices for all the pointwise products since \(s>1\).
The chart and its differential are continuous at this regularity.
The principal low matrix and \(\mathcal S\) depend smoothly on the
finitely many old parameters. The differentiated frame uses the
instantaneous velocities defined after \eqref{eq:exact-h-ins}.
There is one unbounded operator to check in the field equation.
After separating \(D_{\rm free}\), it enters the remainder through
commutators with the changes of coordinates. The preceding Schur
bounds for these commutators hold also with unit weight. With the
weight \(w\) they give bounded operators on \(H^p\) at fixed \(R\).
Interpolation gives boundedness on \(H^s\) and, by the same argument
applied to operator differences, continuity there. For \(\mathcal M\)
the computed commutator includes its time derivative, which is bounded
for fixed \(R\) by the finite angular truncation. Thus convergence of
trial paths in \(C([T,S],H^s)\) suffices for the bounded remainder.
The free generator is kept in the propagator.
The differences of axial diagonal generators are bounded operators
tending to zero since the actions converge in \(\ell^1\) and the
map \(\chi\) is locally Lipschitz on \(\ell^2\).
For two principal evolutions with the same free generator, Duhamel's
formula bounds the difference of propagators on the finite interval
by the integral of the operator-norm difference of their bounded
parts times the two evolution bounds. The boundary projections on
the opposite-copy pairs are fixed, while the trial-dependent low
matrix acts only on their complement. This proves continuity of the
terminal integral map.

\emph{The limit as the terminal time tends to infinity.}
Fix the same \(R\), take \(S_n=T+n\), and choose a fixed point on
each \([T,S_n]\). On \([T,T+m]\) the bounds just proved hold for
every \(n\ge m\) with constants independent of \(n\). Successive
extraction for integer \(m\) gives a single subsequence for which
the field errors converge in \(C([T,T+m],H^s)\) and the old
parameters converge uniformly, for every \(m\). The pointwise
\(Z\) bounds pass to the limit by lower semicontinuity at fixed
\(R\). In particular, for the limiting errors,
\[
 \|h(\tau)\|_Z+\|b(\tau)\|_Z+|J(\tau)|+|\Theta(\tau)|
 \le 2(M+M')\rho\, e^{-\sigma(\tau-T)},\qquad \tau\ge T .
\]
To identify the equation, reconstruct the physical field from each
fixed point by \eqref{eq:reconstruct-ins}. With \(R\) fixed, this
reconstruction is continuous in the stated topology, so the fields
converge in \(C_t H^s\) on every compact physical time interval in
\([T/r,\infty)\). Each satisfies Duhamel's formula there for the
same constant potential \(c_R\). The cubic term converges in
\(C_t H^s\) since \(s>1\). Passing to the limit in that formula
identifies a solution of the full equation throughout the half-line.

The fast stable trace at \(T\) passes to the limit because its
projection is bounded. At the right end the exponential estimate
above supplies the required decay for every component. The known
transverse profile and its particulars also decay, and the known
axial profile tends to the prescribed rotating coordinates of
\(V_R\). Reconstructing at fixed \(R\) therefore gives exponential
convergence to \(V_R\) in \(H^s\), with the factor
\(e^{-\sigma(rt-T)}\) in physical time. Both solutions conserve mass.
Taking their \(L^2\) distance to zero as \(t\to\infty\) gives
\(M(W_R)=M(V_R)\), so the constant \(c_R=-2M(V_R)\) used in the
construction is also the mass gauge of \(W_R\).

\emph{Backward continuation and reconstruction.}
At \(t=T/r\), lower semicontinuity gives \(H^p\) data with
normalized errors of size \(O(\rho)\). Solving the Cauchy problem
from these data gives a path continuous in \(H^p\), which by uniqueness
in \(H^s\) agrees with the compactness limit on their common interval.
Thus at fixed \(R\) its \(Z\) norm is continuous, as needed in the
following continuation argument. Extend the solution backwards by
the defocusing Cauchy flow \cite{HerrKwak}. On the fixed
slow interval \([0,T]\) the model
coefficients and their derivatives are bounded by a finite constant.
The same calculation along this actual trajectory now has normalized
source bounded by \(C_T(e_R+X+\sqrt L X^2)\), where \(X\) here denotes
the sum of the instantaneous error sizes. The principal transverse
propagator is bounded by \(C_T\) in both directions, and the axial
diagonal is isometric. On any continuation interval on which
\(\sqrt L X\le1\), the nonlinear term obeys \(\sqrt L X^2\le X\).
Duhamel's formula, integrated backwards, gives
\[
 X(\tau)\le C_T(\rho+T e_R)+C_T\int_\tau^T X(u)\,du .
\]
Gronwall yields \(X(\tau)\le C'_T\rho\). For large \(R\) this gives
\(\sqrt L X<1/2\) on the whole continuation interval, so by continuity
the interval reaches zero. Here persistence and uniqueness in \(H^p\)
allow continuation of the same field at the auxiliary order.

To read off the physical field, \(\mathcal S-1=O(L^{-1})\) on each
known transverse plane, the particulars cost \(O(r^2)\), and each free
high error is \(O_T(r\rho)\) in \(Z\), uniformly in positive time.
The axial chart estimates replace the added finite axial planes by
their physical Fourier series with error \(O(r e_R+r\rho)\) in \(Z\).
At old indices the action differs from the original one by
\(O(r^2+r|J|)\) and the angle by \(\Theta\), uniformly in time because
the old terminal angular speeds are locked. On this finite-dimensional
torus the radii stay positive and the chart is smooth into every
fixed Sobolev space. Consequently
\(\sup_{t\ge0}\|V_*(t)-V(t)\|_{H^s}=O_T(\rho+r^2)\). The \(e_i\) contribute \(O(rR^{-10})\) in \(Z\). Let \(f_R\)
denote the difference between \(W_R-V\) and the finite sum with
its prefactor in \eqref{eq:insertion-profile}. Combining the
old-coordinate estimate in physical norm with the conversion
\(\|h\|_{H^s}\le CR^s\|h\|_Z\) for the high errors gives
\[
 \sup_{t\ge0}\|f_R(t)\|_{H^s}
 \le C_T(\rho+r^2)
   +C_T R^s r\bigl(L^{-1}+r+e_R+\rho+R^{-10}\bigr)=o(1).
\]
Here \(R^sr=a_0\) and each term in parentheses tends to zero.
The finite sum itself has \(H^{s_0}\) norm
\(O(rR^{s_0})=O(R^{s_0-s})\), so the same decomposition also
gives the uniform lower-order comparison. At the main
root the old smooth field has Fourier coefficient \(o(r)\), uniformly
over its old angle torus. For \(f_R\) the single-coefficient estimate gives
\[
 |\widehat f_R((K,0),t)|\le \langle K\rangle^{-s}\|f_R(t)\|_{H^s}
       =o(r)
\]
uniformly in \(t\), since \(K=1+RP\) and \(r=a_0R^{-s}\) with \(P,a_0\)
fixed. Orthogonality of the finitely many distinct
planes at time zero gives the asserted input cost. This completes the
proof of Proposition~\ref{prop:insertion}.

\section{Persistence along one trajectory}\label{sec:limit}

We now return to all vertical frequencies. To use Proposition~\ref{prop:insertion}
a second time, we must preserve the part of the previously constructed
orbit that led to its periodic asymptote. We first describe the
linearized equation at this asymptote.

\subsection{The previous orbit}

Let \(V\) be smooth, axial and periodic of period \(T_0>0\) in the
equation with real constant potential \(c\). Linearizing the cubic
product at \(V\) gives the equation for a perturbation \(h\),
\begin{equation}\label{eq:previous-linear}
 i\partial_t h=(-\Delta+c+2|V|^2)h+V^2\bar h .
\end{equation}
Let \(L(t,v)\) take the value of \(h\) at time \(v\) to its value at
time \(t\). It is real linear because the equation involves \(\bar h\).
The operator norms below are taken on the indicated Sobolev spaces.
The one-period map \(P\coloneqq L(T_0,0)\) is its monodromy. If \(\Pi\) is a
spectral projection of \(P\), its transport to time \(t\) is
\(\Pi(t)\coloneqq L(t,0)\Pi L(0,t)\). Since \(\Pi\) commutes with \(P\),
this transported projection is \(T_0\)-periodic.
\begin{lemma}\label{lem:floquet}
On each \(H^b(\T^2)\), \(b\ge0\), the spectrum of the complexification
of \(P\) off the unit circle consists of isolated eigenvalues of finite
algebraic multiplicity. These eigenvalues and their generalized
eigenspaces are independent of \(b\), and the eigenspaces consist of
smooth functions. Finite sums of the associated spectral projections,
transported by \(L(t,0)\), are bounded from \(H^{b_0}\) to \(H^b\)
uniformly in \(t\), for \(0\le b_0\le b\).
\end{lemma}
\begin{proof}
Double the equation using \(h,\bar h\) as independent fields, both
expanded in the usual spatial Fourier basis. The generator for
\(i\partial_t\) is the sum of the free diagonal
\(D\coloneqq \operatorname{diag}(1-\Delta,-(1-\Delta))\), a diagonal matrix
of real smooth potentials, and the off-diagonal entries \(V^2,-\bar V^2\).
In these doubled coordinates the spectra and projections below are
taken over the complex numbers. Set
\(J\coloneqq \left(\begin{smallmatrix}0&X\\X^*&0\end{smallmatrix}\right)\).
For \(k,l\in\Z^2\) the matrix of \(X\) in the Fourier basis is
\[
 X_{kl}\coloneqq -\frac{\widehat {V^2}(k-l)}
                         {\langle k\rangle^2+\langle l\rangle^2}.
\]
Multiplying \(X_{kl}\) by \(\langle k\rangle^2+\langle l\rangle^2\)
gives the negative of the upper off-diagonal entry. The lower entry
cancels by taking adjoints with the opposite sign. For the smoothing
bound, \(\langle k\rangle^2/(\langle k\rangle^2+\langle l\rangle^2)\le1\)
and \(\langle k\rangle^b/\langle l\rangle^b\lesssim_b\langle k-l\rangle^b\).
After inserting the Sobolev weights, every entry is bounded in modulus
by a fixed summable function of \(k-l\), since \(V\) is smooth and
periodic. The same calculation with \(\partial_t(V^2)\) proves
\(J,\partial_tJ:H^b\to H^{b+2}\) uniformly over a period. In the coordinates
given by \(e^J\), the conjugated generator is therefore the scalar
diagonal Schr\"odinger operator plus a remainder gaining two derivatives.
Indeed the first commutator with the free diagonal is bounded at every
Sobolev order and every further commutator has a factor \(J\).
The differentiated coordinate change has the same gain. Duhamel's
formula on a period and \(e^{\pm J}-1:H^b\to H^{b+2}\) now give
\(P=P_D+K\), where \(P_D\) is the scalar diagonal monodromy and
\(K:H^b\to H^{b+2}\).

Here and below a scalar evolution refers to either of the two diagonal equations.
Both scalar evolutions are unitary on \(L^2\). Commuting a spatial
derivative of integer order \(m\) through the potential produces
only derivatives of order at most \(m-1\) on the solution. Integration
and induction give polynomial bounds in \(|t-v|\) on \(H^m\), in
either time direction. Interpolation gives the same conclusion for
real \(b\ge0\). At integer multiples of \(T_0\) these are polynomial
bounds on both \(P_D^n\) and \(P_D^{-n}\), \(n\ge1\). The resolvent of \(P_D\) thus exists off the unit
circle, by the geometric series for \(P_D\) or \(P_D^{-1}\).
Compactness of \(K\) and the analytic Fredholm theorem give the
spectral assertion, applied separately inside and outside the unit
circle. Indeed, for \(|\lambda|\ne1\),
\[
 \lambda-P=\bigl[1-K(\lambda-P_D)^{-1}\bigr](\lambda-P_D).
\]
The bracket is the identity plus a compact operator. The product is
invertible at \(\lambda=0\), since the linear evolution is reversible,
and for \(|\lambda|>\|P\|\) by a geometric series. These give one
invertible point in each of the two connected regions. For an eigenvector the
identity \(h=(\lambda-P_D)^{-1}Kh\) gains two derivatives. Induction
gives the same conclusion along a Jordan chain. The resolvents on
different Sobolev spaces agree on their common domains wherever
both exist, so the corresponding finite-rank projections agree.
A bounded finite-rank projection at the lower order with smooth
range is bounded into \(H^b\). Transport and periodicity complete the
proof. Projections defined by conjugation-invariant spectral sets
preserve the physical real subspace.
\end{proof}

Here is the nonlinear consequence we need. Its weak-norm assumption
controls the finitely many rapidly expanding old directions. The
strong-norm assumption is made only at the initial time.
\begin{proposition}\label{prop:history}
Let \(1<s_0<s\). Suppose \(U,V\) are solutions in the same constant
gauge \(c\), \(V\) is as above, and
\(\|U(t)-V(t)\|_s\le C e^{-\gamma t}\) for \(t\ge0\), with \(\gamma>0\).
Let \(W_n\) be solutions in gauges \(c_n\), and suppose
\[
 \varepsilon_n\coloneqq |c_n-c|+\sup_{t\ge0}\|W_n(t)-V(t)\|_{s_0}\to0,\qquad
 d_n\coloneqq \|W_n(0)-V(0)\|_s\to0 .
\]
There are solutions \(U_n=W_n+b_n\) in the gauges \(c_n\) such that
\[
 \sup_{t\ge0}\|b_n(t)-(U(t)-V(t))\|_s\to0,\qquad
 \|b_n(t)\|_s\le C' e^{-\nu t}
\]
for some \(\nu>0\) and all large \(n\). In particular
\(U_n(0)\to U(0)\) in \(H^s\), and \(U_n-U\to0\) uniformly in
\(H^{s_0}\) on the positive half-line.
\end{proposition}
\begin{proof}
For an eigenvalue \(\lambda\) call \(T_0^{-1}\log|\lambda|\)
its growth rate. Choose \(0<\nu<\beta<\gamma\) with \([-\beta,-\nu]\)
in a single gap of the rates, using Lemma~\ref{lem:floquet}.
Let \(P_f(t)\) be the transported projection onto the generalized eigenspaces
of rates below \(-\beta\) (the fast stable directions), and \(P_c(t)\coloneqq 1-P_f(t)\).
For some \(k>0\), in each of the two Sobolev spaces we have
\[
 \begin{aligned}
 \|L(t,v)P_f(v)\|&\le C e^{-(\beta+k)(t-v)}&& (v\le t),\\
 \|L(t,v)P_c(v)\|&\le C e^{(\nu-k)(v-t)}&& (t\le v).
 \end{aligned}
\]
For integer periods this is the spectral radius formula with a
strict gap, applied to the monodromy or its inverse on the indicated
subspace. The propagators over the remaining fraction of a period
are uniformly bounded.

Fix \(0<\eta<\beta-\nu\). We claim that
\begin{equation}\label{eq:delayed}
 \sup_{t\ge0}e^{-\eta t}\|W_n(t)-V(t)\|_s
                    \le C_\eta(\varepsilon_n+d_n).
\end{equation}
Writing \(e\coloneqq W_n-V\), the remainders about the old linearization are
\(2V|e|^2+\bar V e^2+|e|^2 e+(c_n-c)(V+e)\), up to a unit factor.
Their \(H^s\) norm is at most
\(C\varepsilon_n(1+\|e\|_s)\). Here and below we use
\(\|ab\|_s\le C(\|a\|_{s_0}\|b\|_s+\|a\|_s\|b\|_{s_0})\),
which follows from weighted Fourier convolution and
\(\ell^2(\langle k\rangle^{2s_0})\subset\ell^1\).
Project away the finitely many old rates above a threshold in
\((0,\eta/3)\) outside the rate set. The removed projection of \(e\)
is bounded by \(C_\eta\varepsilon_n\) in \(H^s\) at every time by
Lemma~\ref{lem:floquet}. The complementary propagator grows at most
as \(C_\eta e^{\eta(t-v)/2}\) forwards. To see explicitly where
smallness enters, let \(m(t)\coloneqq \sup_{0\le v\le t}e^{-\eta v}\|e(v)\|_s\).
Apply Duhamel's formula to the complementary component, multiply by
\(e^{-\eta t}\), and use
\(\int_0^t e^{-\eta(t-v)/2}\,dv\le2/\eta\).
Adding back the removed component yields
\(m(t)\le C_\eta(d_n+\varepsilon_n)+C_\eta\varepsilon_n m(t)\).
Absorb the last term for large \(n\) to obtain \eqref{eq:delayed}.

Write \(d\coloneqq U-V\). Starting at a late time \(A\), prescribe
\(P_f(A)b(A)=P_f(A)d(A)\) and set
\[
 (\mathcal Gf)(t)\coloneqq \int_A^t L(t,v)P_f(v)f(v)\,dv
                  -\int_t^\infty L(t,v)P_c(v)f(v)\,dv .
\]
For \(\xi=\beta,\nu\) denote by
\(\|h\|_{b,\xi}\coloneqq \sup_{t\ge A}e^{\xi(t-A)}\|h(t)\|_b\).
For \(\nu\le\xi\le\beta\), multiplying the forward kernel by
\(e^{\xi(t-v)}\) leaves decay at rate \(\beta+k-\xi\).
On the backward kernel the corresponding rate is \(\xi-\nu+k\).
Both are at least \(k\), so integration bounds \(\mathcal G\) in the
two norms \(\|\cdot\|_{s_0,\beta}\) and \(\|\cdot\|_{s,\nu}\). We seek a fixed point of
\[
 b\longmapsto L(t,A)P_f(A)d(A)+\mathcal GF_n(b),
\]
where
\[
 iF_n(b)=(c_n-c)b+2(|W_n|^2-|V|^2)b+(W_n^2-V^2)\bar b
                  +2W_n|b|^2+\bar W_n b^2+|b|^2 b .
\]
By \eqref{eq:delayed}, for fixed \(A\) the bound
\(
 D_*\coloneqq \sup_{t\ge A}e^{-\eta(t-A)}(1+\|V(t)\|_s+\|W_n(t)\|_s)
 \le 2+2\sup_t\|V(t)\|_s
\)
holds for all sufficiently large \(n\). On a ball of radius \(a\le1\) in
both norms, write \(h\coloneqq b-b'\). The product inequality gives
\[
 \begin{aligned}
 \|F_n(b)-F_n(b')\|_{s_0,\beta}
 &\le C(\varepsilon_n+a)\|h\|_{s_0,\beta},\\
 \|F_n(b)-F_n(b')\|_{s,\nu}
 &\le C(\varepsilon_n+a)\|h\|_{s,\nu}
                   +C(D_*+1)(1+a)^2\|h\|_{s_0,\beta}.
 \end{aligned}
\]
For instance, put \(a_n\coloneqq |W_n|^2-|V|^2\).
Then \(\|a_n(t)\|_{s_0}\lesssim\varepsilon_n\) and
\(\|a_n(t)\|_s\lesssim (D_*+1)e^{\eta(t-A)}\). The product estimate gives
\[
 e^{\nu(t-A)}\|a_n h\|_s
 \lesssim \varepsilon_n\|h\|_{s,\nu}
 +(D_*+1)e^{-(\beta-\nu-\eta)(t-A)}\|h\|_{s_0,\beta}.
\]
The strict inequality \(\eta<\beta-\nu\) bounds the second term at
every time. In the quadratic terms of \(F_n\), expanding the difference
between \(b\) and \(b'\) leaves one factor of \(b\) or \(b'\), of size
at most \(a\). The cubic term leaves two such factors. In the combined norm
\(\|h\|_{s,\nu}+\Lambda\|h\|_{s_0,\beta}\), first choose \(\Lambda\)
large compared with the second constant and the Green bound, then take
a small ball, \(A\) large so the free datum lies in half the ball,
and \(n\) sufficiently large. This gives a contraction. For the limiting coefficients \(V,c\),
the fixed point is \(d=U-V\), since the backward boundary term
vanishes by the decay \(e^{-\gamma t}\). To compare the two maps at
this fixed \(d\), keep \(A\) fixed. By \eqref{eq:delayed},
\[
 \sup_{t\ge A}e^{-\eta(t-A)}\|W_n(t)-V(t)\|_s
          \le C_\eta e^{\eta A}(\varepsilon_n+d_n)\longrightarrow0 .
\]
Factor \(W_n^2-V^2=(W_n-V)(W_n+V)\) and use the same expansion
with a conjugate for \(|W_n|^2-|V|^2\). By the product inequality,
their \(H^{s_0}\) norms tend uniformly to zero and their \(H^s\) norms
are \(o(1)e^{\eta(t-A)}\), uniformly for \(t\ge A\).
Every coefficient difference in \(F_n(d)-F_\infty(d)\), where
\(F_\infty\) uses \(V,c\), multiplies at least one \(d\) or \(\bar d\).
Its decay absorbs the factor \(e^{\eta(t-A)}\) since
\(\nu+\eta<\gamma\), and also allows the weak weight \(\beta<\gamma\).
In the combined norm \(\|\cdot\|_*=\|\cdot\|_{s,\nu}
+\Lambda\|\cdot\|_{s_0,\beta}\) the two fixed point maps have
the same free term \(L(t,A)P_f(A)d(A)\). If \(\vartheta<1\) is
their common contraction bound, subtraction gives
\[
 (1-\vartheta)\|b_n-d\|_*
 \le \|\mathcal G(F_n(d)-F_\infty(d))\|_*
 \longrightarrow 0 .
\]
In particular \(b_n(A)\to U(A)-V(A)\) in \(H^s\). On \([0,A]\)
continuous dependence from time zero gives \(W_n\to V\) in
\(C([0,A],H^s)\), since \(d_n+|c_n-c|\to0\). Thus
\(W_n(A)+b_n(A)\to U(A)\) in \(H^s\). Continuous dependence
backwards from \(A\) gives \(W_n+b_n\to U\) on that same interval.
Together with the weighted estimate on \([A,\infty)\) this proves
the proposition. If \(W_n\) tends to a periodic solution \(V_n\), then
\(U_n\) tends to the very same \(V_n\). Equality of their masses follows
on taking the limit of the difference, so their mass gauges also agree.
\end{proof}

\subsection{Choosing the data and the times}

Two indices enter the construction. The stage number \(j\) counts
transfers on the final orbit. Within a single stage we vary the tree
depth \(d\) to reduce the initial cost and choose a frequency spacing
\(R=R(d)\) to approximate that particular tree. Only one pair
\((d,R(d))\) is retained at stage \(j\).

Fix \(1<s<2\), put \(s_0\coloneqq (1+s)/2\), and fix \(G,\ell\) as in
Theorem~\ref{thm:clock}. Begin with \(U_0=V_0\) as in
Proposition~\ref{prop:bank}. The initial field \(V_0(0)\) may be taken
arbitrarily small in \(H^s\) because its Birkhoff indices are fixed, the five
actions tend to zero, and the inverse chart is continuous in \(H^s\).
Write \(g_0\ne0\) for its coefficient in \eqref{eq:bank}.
At any finite stage we have a true solution \(U\) exponentially
asymptotic to a smooth periodic axial solution \(V\), both in the
mass gauge of \(V\), whose coefficient \(g\) in \eqref{eq:bank}
satisfies \(|g|\ge |g_0|/2\).

We explain first how to make one choice with arbitrarily small error.
For an increasing sequence of depths \(d\), apply
Proposition~\ref{prop:tree} and Lemma~\ref{lem:tree-orbit} with main
terminal amplitude \(a=1\). Once the geometry is fixed, take the radii
\(b_j>0\) so small that \(\sum_j b_j^2\le1\), that the hypothesis of
Proposition~\ref{prop:clock} holds for \(\theta=3/4\), and that the
sum of \(|v/P|^{2s} b_j^2\) over all auxiliary vertices of the respective
levels is at most \(1/d\). Here \(P\) is the positive root abscissa.
Use Proposition~\ref{prop:insertion} with \(a_0=P^{-s}\), so
\(r=(RP)^{-s}\), and call its finite plane-wave sum
\(C_{d,R}\), including its prefactor \(r\). Orthonormality gives
\begin{equation}\label{eq:carrier-norm}
 \|C_{d,R}(t)\|_s^2
  =\sum_{v\in\mathcal V} r^2\langle e+Rv\rangle^{2s}
                   |x_v(2|g|rt)|^2,\qquad e=(1,0).
\end{equation}
For each fixed tree the coefficients \(r^2\langle e+Rv\rangle^{2s}\)
tend to \(|v/P|^{2s}\). The summed main weight \(\sum_v|v/P|^{2s}\)
over any single primary generation is at most one, since each split of a
weight \(w\) replaces it by \(w(\zeta^s+(1-\zeta)^s)\le w\).
The primary amplitudes are constant within a generation and their
squares sum over the generations to at most two by \eqref{eq:tree-integrals}.
Thus \(\sup_{t\ge0}\|C_{d,R}(t)\|_s\le 3\) for all sufficiently large
\(R\) at each \(d\). Its initial squared norm is at most
\(2c_s^d+1/d+o_{R\to\infty}(1)\), by \eqref{eq:leaf-cost}.

We now apply Proposition~\ref{prop:history} with the old pair
\(U,V\) fixed and with a sequence indexed by the depth \(d\).
Write \(c\) for the gauge of \(V\) and \(W_{d,R}\) for the solution
of Proposition~\ref{prop:insertion} at this depth. For each \(d\)
choose \(R=R(d)\) large enough that
\[
\begin{aligned}
 |c_R-c|+\sup_{t\ge0}\|W_{d,R}(t)-V(t)\|_{s_0}
        &\le 2^{-d},\\
 \sup_{t\ge0}\|W_{d,R}(t)-V(t)-C_{d,R}(t)\|_s
        &\le 2^{-d},\\
 \|C_{d,R}(0)\|_s &\le (2c_s^d+1/d)^{1/2}+2^{-d}.
\end{aligned}
\]
The first two choices use the fixed-tree limits, and the third uses
the leaf estimate above. At time zero the last two inequalities give
\(\|W_{d,R}(0)-V(0)\|_s\le (2c_s^d+1/d)^{1/2}+2^{1-d}\to0\).
Thus the global weak comparison and the strong initial comparison
hold on the same sequence. The strong norm of the new contribution
can still be of order one at late times.

Let \(b_{d,R}\) be the correction given by
Proposition~\ref{prop:history} and put \(\widetilde U\coloneqq W_{d,R}+b_{d,R}\).
Subtracting the old solution gives
\[
 \widetilde U-U-C_{d,R}
   =(W_{d,R}-V-C_{d,R})+\bigl[b_{d,R}-(U-V)\bigr].
\]
The first term is the error for a single tree, made small by choosing
\(R\) after \(d\). The second tends uniformly to zero by
Proposition~\ref{prop:history}, since the strong initial cost tends
to zero with depth and the weak error is arbitrarily small for each
fixed tree. We conclude that
\begin{equation}\label{eq:alltimecomparison}
 \sup_{t\ge0}\|\widetilde U(t)-U(t)-C_{d,R}(t)\|_s\longrightarrow0
\end{equation}
along the choices. The terms in this sum are compared in their specified
gauges. Multiplication by any of these scalar phases preserves each
Fourier modulus. We also have
\(\sup_{t\ge0}\|\widetilde U(t)-U(t)\|_{s_0}\to0\),
and \(\widetilde U\) is exponentially asymptotic to \(V_R\). The
initial distances in \(H^s\) tend to zero. All these conclusions hold
with any additional lower requirements on \(R(d)\), since for each
fixed tree the limits in Proposition~\ref{prop:insertion} hold as
\(R\to\infty\).

The forward orbit of \(U\) is relatively compact in \(H^s\) because its
restriction to a bounded interval is compact and its late tail
approaches the compact range of \(V\). Let \(P_{\ge N}\) retain the Fourier coefficients at \(|k|\ge N\).
Then \(\sup_{t\ge0}\|P_{\ge N}U(t)\|_s\to0\) as \(N\to\infty\).
Indeed, cover the compact orbit closure by finitely many \(H^s\)
balls of radius \(\varepsilon\) and take \(N\) large enough that
the tail of each center has norm below \(\varepsilon\). The projection
has norm at most one, giving \(2\varepsilon\) for every point of the
orbit. Each frequency in \(C_{d,R}\) has length at least \(R-1\).
Thus its cross term with \(U\) is bounded by
\[
 |\langle U(t),C_{d,R}(t)\rangle_{H^s}|
 \le \|P_{\ge R-1}U(t)\|_s\|C_{d,R}(t)\|_s
 \le 3\|P_{\ge R-1}U(t)\|_s .
\]
Here \(\langle\,\cdot\,,\,\cdot\,\rangle_{H^s}\) is the inner product
associated with our Fourier norm. The other error in expanding the
squared norm is bounded by the left side of \eqref{eq:alltimecomparison}
times \(2\sup_t\|U(t)\|_s+6+o(1)\). It follows that
\begin{equation}\label{eq:orth-switch}
 \sup_{t\ge0}\left|
   \|\widetilde U(t)\|_s^2-\|U(t)\|_s^2-\|C_{d,R}(t)\|_s^2
                 \right|\longrightarrow0 .
\end{equation}
This estimate includes the period during which the transfer is taking
place.

Set \(K\coloneqq 1+RP\), \(D\coloneqq K^s\), and consider the window
\([D^{3/4},D^{5/4}]\). With the finite tree fixed, \(rt\) tends uniformly
to zero on \(0\le t\le D^{3/4}\) as \(R\to\infty\). The early value of
\eqref{eq:carrier-norm} is therefore bounded by
\(2c_s^d+1/d+o_{R\to\infty}(1)\). On \(t\ge D^{5/4}\) the slow time
tends uniformly to infinity. The terminal amplitudes of
Lemma~\ref{lem:tree-orbit} then imply
\[
  \sup_{t\ge D^{5/4}}\big|\|C_{d,R}(t)\|_s^2-1\big|
                         \le 1/d+o_{R\to\infty}(1).
\]
Only the main root and the finitely many axial auxiliary amplitudes
survive. A bound on the time of capture itself comes from the finite system.
Put
\begin{equation}\label{eq:physical-clock}
 t_*\coloneqq \frac{C_{3/4}d^{3/2}(K-1)^s}{2|g|}.
\end{equation}
By this time the model root has reached \(3/4\) and stays there or above,
by Proposition~\ref{prop:clock}. The old field \(U\) has uniformly
small normalized Fourier coefficient at \((K,0)\) by compactness.
The root coefficient of \(C_{d,R}\) has normalized modulus
\(r\langle K\rangle^s H(2|g|rt)\), where \(r\langle K\rangle^s\to1\).
Together with \eqref{eq:alltimecomparison} this gives
\(\langle K\rangle^s|\widehat{\widetilde U}((K,0),t)|\ge 2/3\)
for every \(t\ge t_*\), if the choice is sufficiently far along the
sequence.

We give the simultaneous choices explicitly. At stage \(j\ge1\) use
the preceding procedure with \(U=U_{j-1}\), \(V=V_{j-1}\). Denote
the new solutions by \(U_j,V_j\), the root index by \(K_j\), the time by
\(D_j\coloneqq K_j^s\), and the plane sum by \(C_j\). Fix any \(\delta>0\).
Use the convergence in depth and, at each depth, the freedom to enlarge
\(R\), to meet the following conditions.
\begin{enumerate}
\item The initial distance \(\|U_j(0)-U_{j-1}(0)\|_s\) is less than
\(\delta 2^{-j}\). The error in \eqref{eq:orth-switch} is at most
\(2^{-j}\). The early squared norm of \(C_j\) and the late error from
one in its squared norm are at most \(2^{-j}\).
\item Require \(G(D_j)\ge j\) at every stage. For \(j>1\), also
require \(K_j>K_{j-1}\) and \(D_j^{3/4}>D_{j-1}^{5/4}+1\). Take \(\ell(K_j)\ge C_{3/4}d^{3/2}/(2|g|)\),
which makes \eqref{eq:physical-clock} at most \(K_j^s\ell(K_j)\).
\item At every earlier recorded frequency \((K_i,0)\), \(i<j\),
the change in normalized Fourier modulus, uniformly for all
\(t\ge0\), is at most \(2^{-j}/6\). Indeed the modulus difference
is at most \(\langle K_i\rangle^{s-s_0}\sup_t\|U_j(t)-U_{j-1}(t)\|_{s_0}\)
and only finitely many indices have been specified. At the new frequency
impose the \(2/3\) bound above.
\item Keep the successive changes of the periodic targets in uniform
\(H^1\) below \(|g_0|2^{-j-2}\), and the actions in the ball of
Lemma~\ref{lem:locking}. At each fixed tree this follows by taking
\(R\) large, since the total \(h^1\) norm of the finitely many new coordinates
is \(O_d(Rr)\to0\), the old action change is \(O_d(r^2)\), and their
angle speeds are unchanged. The common-period integral \eqref{eq:bank}
then keeps \(|g|\ge |g_0|/2\) throughout the construction.
\end{enumerate}
There is no uniformity in the finite-tree constants required by these
choices. More formally one can first impose all the lower requirements
on \(R(d)\) involving the given geometry, including the conditions on the
clock. The input errors for Proposition~\ref{prop:history} still tend to
zero as \(d\to\infty\). Taking one sufficiently large depth in this
sequence imposes the remaining strong and weak error tolerances.

The initial fields converge in \(H^s\) to a datum \(u_0\) with
\(\|u_0-U_0(0)\|_s\le\delta\). Global well-posedness and persistence
of regularity for the defocusing flow \cite{BourgainRestrict,HerrKwak}
give convergence of the physical trajectories in \(H^s\) on every
compact time interval. Equivalently, one can use the \(H^s\) algebra
estimate for local existence, uniqueness and continuous dependence at
this regularity and the global defocusing theory to continue the limit.
The masses converge, so the convergences on bounded intervals also hold
in the respective mass gauges. For fixed \(i\) and \(t\ge t_i\), where
\(t_i\) is \eqref{eq:physical-clock} for its tree, condition (3) and convergence
at the fixed frequency \(K_i\) give
\[
 \langle K_i\rangle^s|\widehat u((K_i,0),t)|
 \ge \frac23-\sum_{j=i+1}^\infty \frac{2^{-j}}6
 \ge \frac12 .
\]
In particular, for any \(m\ge1\) the squared norm is at least \(m/4\)
at every time \(t\ge\max_{1\le i\le m} t_i\), by summing over those
\(m\) distinct frequencies. This already proves the divergence in
Theorem~\ref{thm:main}. The more precise time dependence follows by
also accounting for what happens during a transfer.

Finally sum \eqref{eq:orth-switch} from stage one to stage \(n\).
The sum of the absolute errors is at most one. At any fixed \(t\)
the windows are disjoint, so at most one term in
\(\sum_{j\le n}\|C_j(t)\|_s^2\) is in a transition window. That term
is between zero and nine. At stage \(j\) compare the squared norm to the indicator
\(\mathbf 1_{\{D_j\le t\}}\), which is one if \(D_j\le t\) and zero
otherwise. Before the window both are within \(2^{-j}\) of zero, and
after the window their difference from each other is at most \(2^{-j}\).
Inside the one possible active window the difference is at most nine.
Summing these errors gives ten as a uniform bound for the difference
between the sum of the squared carrier norms and the counting function.
Adding the error of at most one from \eqref{eq:orth-switch} gives
\[
 \left|\|U_n(t)\|_s^2-\#\{j\le n:D_j\le t\}\right|
         \le 11+\sup_{v\ge0}\|U_0(v)\|_s^2 .
\]
The right-hand side is finite by periodicity and independent of \(n,t\).
For each \(t\), the left-hand norm converges by continuous dependence,
and the count stabilizes since \(D_j\to\infty\). This proves
\eqref{eq:count} uniformly in \(t\). On
\([D_j,D_{j+1})\) it gives \(\|u(t)\|_s\le (j+O(1))^{1/2}\)
and \(G(t)\ge j\). Hence \(\|u(t)\|_s/G(t)\to0\).
The counting function tends to infinity, and the initial norm can be
made arbitrarily small by the choices of \(U_0\) and \(\delta\).
This proves Theorems~\ref{thm:main} and \ref{thm:clock} for \(1<s<2\).

\subsection{The limiting axial solution and the height of a transfer}

The periodic axial solutions constructed above have a limit in every
lower Sobolev norm if the weak tolerances are chosen summably. One
can also measure the speed at which the two-dimensional solution
approaches this limit. Let \(\alpha\) be as in \eqref{eq:parameters-ins}
and put \(k_s\coloneqq 3\alpha/4\), so the error scale in that section is
\(\rho=R^{-k_s}\). We use almost periodicity in the following sense. An \(H^\sigma\)-valued function is almost periodic if it is a uniform
limit on \(\R\) of finite sums \(\sum_m f_m e^{i\omega_m t}\), with
\(f_m\in H^\sigma\) and \(\omega_m\in\R\).

\begin{corollary}\label{cor:asymptote}
For \(1<s<2\) the data in Theorem~\ref{thm:clock} can be chosen so that there exists
an axial solution \(v\), almost periodic in time in \(H^\sigma\) for
each \(1<\sigma<s\), satisfying
\[
 \lim_{t\to+\infty}(1+t)^b\|u(t)-v(t)\|_{H^\sigma}=0
 \qquad\left(0\le b<\frac{\min(s-\sigma,k_s)}s\right).
\]
At every \(t\), \(\langle K_j\rangle^s|\widehat v((K_j,0),t)|\ge1/2\)
for every \(j\). In particular \(v(t)\in\bigcap_{1<\sigma<s}H^\sigma
\setminus H^s\).
\end{corollary}
\begin{proof}
Fix a finite tree and one pair \((\sigma,b)\) as in the statement.
In the notation of Section~\ref{sec:insertion}, on the terminal interval
\(rt\ge T\) the difference of the old Birkhoff variables from their
terminal values at the same time is \(O(\rho e^{-\sigma_*(rt-T)})\)
in any fixed sequence norm on \(F\). Here \(\sigma_*>0\) denotes the
time weight used for this tree. On \(F^c\) the difference has \(Z\)
norm \(O(r e^{-\sigma_*(rt-T)})\), by subtracting \(rY_\infty\) in
\eqref{eq:reconstruct-ins}. The transverse part obeys the same bound.
On \(0\le rt\le T\) the backward estimate gives the respective bounds
\(O(\rho)\) and \(O(r)\). In \(h^\sigma\) the new axial coordinates
are \(O(R^\sigma r)\), tending to zero. The compact set of old angle
values thus has a neighborhood where the differential of the inverse
chart \(h^\sigma\to H^\sigma\) is uniformly bounded. Applying it to
the difference of the actual and terminal coordinates at the same time
gives, for \(R\) large,
\[
 \|W_R(t)-V_R(t)\|_{H^\sigma}
 \le C(\rho+R^{\sigma-s})\exp[-\sigma_*(rt-T)_+],
 \qquad t\ge0 .
\]
The constants here may depend on the tree. Writing \(\tau=rt\) and
using \(0<r\le1\), we have
\((1+t)^b\le r^{-b}(1+\tau)^b\). The supremum of
\((1+\tau)^b e^{-\sigma_*(\tau-T)_+}\) over \(\tau\ge0\) is finite.
Therefore the weighted error is at most
\(C(R^{sb-k_s}+R^{sb+\sigma-s})\). Both exponents are negative
precisely under the two strict upper bounds on \(b\) in the statement.

Use a countable dense collection of pairs \((\sigma,b)\) strictly
inside the allowed range, with \(\sigma,b\) rational. In the
construction of stage \(j\), impose estimates for its first \(j\)
pairs, taking \(s_0<s\) above their spatial exponents if necessary.
For each pair the preceding bound can be made arbitrarily small.
On the history correction, the convergence of \(b_n-(U-V)\) in the
exponential strong norm in the proof of Proposition~\ref{prop:history}
and on the fixed initial interval gives convergence in the polynomially
weighted \(H^\sigma\) norm. It follows that the successive differences of
\(U_j-V_j\) are summably small in that norm, each field in its own
mass gauge. Indeed
\((U_j-V_j)-(U-V)=(W_R-V_R)+(b_n-(U-V))\).
For each fixed pair, the weighted differences are bounded, for example,
by \(2^{-j}\) at all stages after it enters the list. For any \(n\)
beyond that stage, the weighted distance from the limit to \(U_n-V_n\)
is at most \(\sum_{j>n}2^{-j}\). Let \(t\to\infty\) first, using the
exponential decay of this fixed \(U_n-V_n\), and then let \(n\to\infty\).
This proves the asserted vanishing of the weighted limit.

We also require summable uniform \(H^\sigma\) differences of the targets
\(V_j\) on the whole real line. The locked old angles and the small
\(h^\sigma\) norm of the added coordinates give this by
Lemma~\ref{lem:locking} and the bounded chart differential.
Each \(V_j\) is approximable uniformly in \(H^\sigma\) on its period
by its time Fourier means (convolution with the Fej\'er kernels).
The uniform limit of these periodic targets is therefore almost periodic
in the sense just defined. Local well-posedness identifies it in its limiting gauge
as a solution independent of \(y\). At each stage \(U_j,V_j\) have
the same mass. Their limits consequently have the same mass, and
removing the common gauge preserves the norm of the difference.
It preserves almost periodicity too, since multiplication by
\(e^{ict}\) preserves uniform approximation by trigonometric
polynomials in time. The strict inequalities on \(\sigma,b\) allow
each real pair in the statement to be dominated by one of the rational
pairs. To see why the bound for the axial limit holds even at early times,
fix \(i\), \(j\ge i\), and \(t\in\R\). Denote the period of \(V_j\)
in its mass gauge by \(T_j^{\rm ax}\). For all sufficiently large
integers \(m\), the time \(t+mT_j^{\rm ax}\) is past the capture time
of \(K_i\). At these times \(U_j\) satisfies the lower bound \(1/2\)
at that index and \(U_j-V_j\to0\) in \(H^s\). Periodicity gives the
same lower bound for \(V_j(t)\). Let \(j\to\infty\), using continuity
of the single Fourier coefficient in \(H^\sigma\). The limit \(v(t)\)
has normalized modulus at least \(1/2\) at every \(K_i\), and hence the
sum defining \(\|v(t)\|_{H^s}^2\) diverges.
\end{proof}

There is a second way to normalize the trees. A small increase in the
initial cost at each leaf allows the Sobolev contribution of a single
root to grow exponentially with its depth.

\begin{corollary}\label{cor:height}
Fix \(1<s<2\) and \(0<\kappa<-\tfrac12\log c_s\), with \(c_s\)
as in the leaf estimate \eqref{eq:leaf-cost}. There is a single \(H^s\) datum,
a strictly increasing sequence of depths \(d_j\), and distinct
frequencies \((K_j,0)\), with the following property. Put
\(B_j\coloneqq e^{\kappa d_j}\). At that frequency the normalized Fourier modulus
\(\langle K_j\rangle^s|\widehat u((K_j,0),t)|\)
of the solution is at least \(B_j/2\) at every time
\[
 t\ge \frac{C_{3/4}d_j^{3/2}(K_j-1)^s}{2|g_j| B_j},
\]
where \(g_j\) is the coefficient \eqref{eq:bank} of the preceding
periodic orbit and \(\inf_j |g_j|>0\).
\end{corollary}
\begin{proof}
Keep \(a=1\) but replace \(a_0\) by \(B(d)/P^s\) at depth \(d\),
where \(B(d)\coloneqq e^{\kappa d}\). The primary initial cost is now
\(O(\exp[(\kappa+\tfrac12\log c_s)d])=o(1)\).
For each fixed tree, choose the positive auxiliary radii so that their
total normalized initial cost is at most \(2^{-d}\) and
Proposition~\ref{prop:clock} applies. Then take \(R=R(d)\) large enough
that the weak distance (including the gauge difference) and the absolute
\(H^s\) error in \eqref{eq:insertion-profile} are each at most \(2^{-d}\).
Thus the strong initial distance tends to zero with \(d\), as required by
Proposition~\ref{prop:history}, even though the late root height grows.
The same subtraction as in \eqref{eq:alltimecomparison} gives a uniform
absolute \(o(1)\) error after the history correction.

At the new root, \(r=B(d)/(K-1)^s\), so \(r\langle K\rangle^s\ge B(d)\).
Once \(2|g|rt\ge C_{3/4}d^{3/2}\), the model coefficient therefore has
normalized modulus at least \(3B(d)/4\). The old field has uniformly
vanishing normalized coefficients as \(K\to\infty\) by precompactness
in \(H^s\). Subtracting this coefficient and the absolute comparison
error leaves at least \(2B(d)/3\) on the new orbit. At each later stage
\(j\), impose, for every retained root \(i<j\),
\[
 \langle K_i\rangle^{s-s_0}\sup_{t\ge0}\|U_j(t)-U_{j-1}(t)\|_{s_0}
       \le \tfrac16 B_i 2^{-(j-i)} .
\]
There are finitely many such requirements at each stage. Their sum in
\(j>i\) costs at most \(B_i/6\) at root \(i\), giving the asserted
\(B_i/2\) bound in the limit. Convergence of the initial data and the
lower bound on \(|g_j|\) follow by the same summable choices as above.
\end{proof}

\section{The endpoint \texorpdfstring{\(H^2\)}{H2}}\label{sec:h2-extension}

We prove the qualitative assertion at \(s=2\). In this section \(L>0\)
denotes the height of one transfer. The depth, the axial background and
the strength of its selected Fourier coefficients will be fixed before
the dilation tends to infinity. The time of detection may depend on all
these choices.

\subsection{The second-order normal frequency}
Use the individual mass gauge of Section~\ref{sec:insertion}.
With \(A,B,\mathcal{Q}\) as in Lemma~\ref{lem:empty-frequency}, put
\[
 \Xi(v)=A/v_x^2-B/|v|^2,\qquad d(v,w)=\Xi(v)+\Xi(w)-\Xi(v+w).
\]
The Liouville primitive in that lemma is unchanged by the common
mass gauge, since the two canonical forms have the same pairing with
its phase generator. In particular \(A\) is the time mean of the
spatial variance of \(|V|^2\), also for resonant old frequencies.
Indeed \(M=\langle|U|^2\rangle\) is conserved, so the spatial variance
equals \(\langle|U|^4\rangle-M^2=\mathcal Q-M^2+2\mathcal LG\) by the
identity in that proof. Along any trajectory \(\mathcal LG\) is the
time derivative of the bounded function \(G\), and its mean over
\([0,T_0]\) is \((G(T_0)-G(0))/T_0\to0\) as \(T_0\to\infty\).
On a small two-site torus with positive comparable radii its first
nonzero term is positive. In detail, the phase symmetry leaves only
terms of odd degree in the Taylor expansion of \(\Psi\) at zero. On
two sites \(n_1\ne n_2\) with actions \(I_1,I_2\) this gives
\[
\begin{gathered}
 U=\sqrt{I_1}\,e^{i(n_1x+\varphi_1)}+\sqrt{I_2}\,e^{i(n_2x+\varphi_2)}
     +O(|I|^{3/2}),\\
 |U|^2-M=2\sqrt{I_1I_2}\cos\bigl((n_1-n_2)x+\varphi_1-\varphi_2\bigr)
     +O(|I|^2),
\end{gathered}
\]
because \(M=I_1+I_2\) exactly. The spatial mean of the square of the
last line is \(2I_1I_2+O(|I|^3)\) at every time. Therefore
\(A=2I_1I_2+O(|I|^3)\) and \(B=(A+M^2)/2\) are positive for
comparable small radii. Start there and retain \(A,B>0\) under
sufficiently small perturbations. In the mass gauge the lemma gives
high transverse diagonal \(n^2+m^2+A/n^2-B/(n^2+m^2)\), with
\(O_h(R^{-3})\) residual on \(|n|\ge hR,\, |m|\ge R\) and arbitrary
off-diagonal spatial decay at fixed smooth background.

At an empty axial index use the same doubled computation with \(m=0\).
To see also the asymptotics needed for a finite difference below,
continue the spatial recursion to order four, keeping the instantaneous
spatial mean as scalar \(a(t,n)=\sum_{j=2}^4 a_j(t)n^{-j}\).
At each order the upper equation sets its mean and the zero-mean
coefficient at the next order by \(2D_x\); the lower leading divisor
is \(2n^2\) and determines the next lower coefficient directly. Taking
one more upper and two more lower coefficients gives residual
\(O(|n|^{-5})\) in \(L^2\). Explicitly, in the notation of the proof
of Lemma~\ref{lem:empty-frequency} with \(E=n^2\), put
\(f=1+\sum_{j=1}^5n^{-j}f_j\) and \(g=\sum_{j=2}^6n^{-j}g_j\), with
\(f_0=1\), \(g_0=g_1=0\) and \(\langle f_j\rangle=0\) for \(j\ge1\).
The coefficients of \(n^{-j}\), \(0\le j\le4\), in the two residuals
\((2nD+K_+-a)f+Cg\) and \((2n^2+2nD+K_-+a)g+\bar Cf\) vanish
precisely when
\[
\begin{aligned}
 2Df_{j+1}+K_+f_j-\sum_{i=2}^{j}a_if_{j-i}+Cg_j&=0,\\
 2g_{j+2}+2Dg_{j+1}+K_-g_j+\sum_{i=2}^{j-2}a_ig_{j-i}+\bar Cf_j&=0 .
\end{aligned}
\]
The second line defines \(g_{j+2}\); for instance \(g_2=-\bar C/2\).
The spatial mean of the first line is \(\langle K_+f_j+Cg_j\rangle-a_j\)
for \(j\ge2\), which defines \(a_j\), and its mean-free part defines
\(f_{j+1}\). For \(j=0,1\) the mean vanishes as in the lemma, since
\(\langle K_+f_j\rangle=\langle V_0f_j\rangle\). Thus
\[
 a_2=\langle V_0f_2\rangle-\tfrac12\langle|C|^2\rangle,\qquad
 a_3=\langle V_0f_3+Cg_3\rangle,\qquad
 a_4=\langle V_0f_4+Cg_4\rangle,
\]
and the first of these is the scalar of the lemma at \(E=n^2\).
The terms left in the two residuals are of order \(n^{-5}\), with
coefficients involving finitely many derivatives of the smooth torus.
Apply \(D\Phi(V)\) to this complexified
normal column. With its scalar phase removed its empty \(n\)th
coefficient \(c_n(t)\) is \(1+O(1/|n|)\) uniformly: apply integer
one-smoothing \(H^1\to h^2\) to a wave of \(H^1\) size \(O(|n|)\)
and then take the \(n\)th coefficient \cite{KST}. In this coordinate
the exact linearized equation is a scalar rotation. Divide its
inhomogeneous equation by \(c_n\) and average in time. The logarithmic
derivative has zero average because a bounded continuous branch of
\(\log c_n\) exists. In formulas, as in the last display of the proof
of Lemma~\ref{lem:empty-frequency},
\[
 i\dot c_n=\bigl(\chi_n-a(t,n)\bigr)c_n+O(|n|^{-5}),\qquad
 \chi_n=\frac1{T_0}\int_0^{T_0}a(t,n)\,dt
        +\frac{i}{T_0}\bigl[\log c_n\bigr]_0^{T_0}+O(|n|^{-5}),
\]
and the bracket is bounded uniformly in \(T_0\) since
\(|c_n-1|\le1/2\) for large \(|n|\). The coefficients \(a_j\) are smooth functions on
one finite torus and thus have time means even at resonant frequencies.
Consequently at a fixed axial state
\[
 \chi_n=(A-B)n^{-2}+c_3n^{-3}+c_4n^{-4}+O_V(|n|^{-5})
 \qquad (|n|\longrightarrow\infty).
\]
Only the first coefficient requires a primitive uniform on the torus.

\subsection{A finite tree and distinct driving frequencies}
Use Proposition~\ref{prop:tree}. We can further arrange that the
numbers \(d(v,w)\) on the primary and auxiliary triangles are
nonzero and pairwise distinct, and differ from every mixed triangle
defect. To verify this, use the rational split parameter \(t\) in that
proposition. At a nonhorizontal parent \((X,Y)\) and \(t\to0\), the
primary defect has a pole of order two, the auxiliary defect at \(p\)
a pole of order four, and the one at \(q\) has order at most two.
Explicitly, put \(U=(X+tY)/(1+t^2)\) and \(W=(Xt-Y)/(1+t^2)\). Then
\[
\begin{gathered}
 p=U(1,t),\qquad q=W(t,-1),\qquad J(p)=(t^2U,-tU),\qquad J(q)=(W/t,W),\\
 A(p)=(X+tY,0),\qquad A(q)=(X-Y/t,0),\qquad
 T_m=\bigl((Yt^2+Xt-Y)/t,\,-Y\bigr),
\end{gathered}
\]
and, writing \(\Xi\) with the coefficients \(A,B\) of this section,
\[
\begin{aligned}
 \Xi(p)&=\frac{A}{U^2}-\frac{B}{U^2(1+t^2)},&
 \Xi(J(p))&=\frac{A}{t^4U^2}-\frac{B}{t^2U^2(1+t^2)},&
 \Xi(A(p))&=\frac{A-B}{(X+tY)^2},\\
 \Xi(q)&=\frac{A}{t^2W^2}-\frac{B}{W^2(1+t^2)},&
 \Xi(J(q))&=\frac{At^2}{W^2}-\frac{Bt^2}{W^2(1+t^2)},&
 \Xi(A(q))&=\frac{(A-B)t^2}{(Xt-Y)^2}.
\end{aligned}
\]
As \(t\to0\) one has \(U\to X\ne0\) and \(W\to-Y\ne0\). Hence
\(d(p,q)=A/(t^2Y^2)+O(|t|^{-1})\),
\(d(p,J(p))=A/(t^4X^2)+O(|t|^{-3})\) and
\(d(q,J(q))=A/(t^2Y^2)+O(|t|^{-1})\), because \(A>0\).
Interchanging the children by \(t\mapsto-1/t\) proves the remaining
separations and nonconstancies. Indeed this substitution exchanges
\(p\) and \(q\), so as \(t\to\infty\) the defect \(d(q,J(q))\) has a
pole of order four in \(t\), whereas \(d(p,q)\) and \(d(p,J(p))\) have
order at most two. The three intended defects are thus nonconstant
and pairwise distinct rational functions of \(t\). At a root of \(Yt^2+Xt-Y=0\),
the mixed output \(T_m\) has zero horizontal coordinate and nonzero
vertical coordinate, while the intended defects at this split are
finite. Thus its defect has a pole there. Here the roots are real,
\(t\ne0\), and \(U\ne0\ne W\): substituting \(t=-X/Y\) or \(t=Y/X\)
into \(Yt^2+Xt-Y\) gives \(-Y\) and \(Y^3/X^2\), respectively.
The term \(-\Xi(T_m)\) of the mixed defect contains
\(-A/(T_m)_x^2\), while \(\Xi(J(p))\), \(\Xi(J(q))\) and all terms of
the intended defects are finite at such a root. At the horizontal root the
single defect is nonconstant by the same projection formula: for
\(Y=0\) it is
\(d(p,q)=(1+t^2)^2(1+t^{-4})A/X^2-(1+t^2)(1+t^{-2})B/X^2-(A-B)/X^2\),
with a pole of order four at \(t=0\). At each
split avoid equality with every previously fixed relevant value, and
impose the within-split exclusions just proved. This removes finitely
many values and is compatible with the real balance and finite-place
conditions in Proposition~\ref{prop:tree}. All these inequalities
persist under small changes of \(A,B\).

First keep only the primary and auxiliary triangles, with equal absolute
couplings. After time and phase normalization the symmetric real
equations are \eqref{eq:tree-primary}--\eqref{eq:tree-aux} with all
mixed edges omitted. In detail the incoming primary term is \(x_{j-1}^2\),
and the two losses are \(x_jx_{j+1}\) and \(w_jh_j\) when that vertex
has a private pair; \(w_j'=-x_jh_j,\ h_j'=x_jw_j,\ H'=x_{d-1}^2\).
Prescribe \(w_j^2+h_j^2=b_j^2\) for \(j\le d-2\) and
\(H^2+\sum_{j<d}x_j^2+\sum h_j^2=1+\sum b_j^2\).
Both quantities are conserved. First,
\((w_j^2+h_j^2)'=2x_j(-w_jh_j+h_jw_j)=0\). Next, with \(x_d=H\),
\(x_{-1}=0\) and the term \(w_jh_j\) absent for \(j=d-1\), half the
derivative of the second sum is
\[
 \sum_{j<d}x_j\bigl(x_{j-1}^2-x_jx_{j+1}-w_jh_j\bigr)+x_dx_{d-1}^2
   +\sum_{j\le d-2}h_jx_jw_j
 =\sum_{j=1}^{d}x_jx_{j-1}^2-\sum_{j=0}^{d-1}x_j^2x_{j+1}=0 .
\]
Start with only \(x_0>0\) among the primaries, so that
\(x_0(0)^2=1+\sum_j(b_j^2-h_j(0)^2)\ge1\). At every fixed time
\(S\) the circle-angle map from launch to time \(S\) is homotopic
to the identity. Choose a preimage of \(h_j(S)=b_j\) for all \(j\).
Monotonicity gives \(0\le H\le1\) up to \(S\); compactness as
\(S\to\infty\) gives such an orbit on the half-line. On its invariant
omega-limit set \(H\) is constant, so \(H'=x_{d-1}^2\) and the
equation for \(x_{d-1}\) give \(x_{d-1}=x_{d-2}=0\).
Suppose \(x_{j+1}=x_j=0\) on a limit orbit, with
\(1\le j\le d-2\). The \(j\)th circle is stationary, and
\(x_{j-1}^2=w_jh_j\) is constant. If this constant were positive,
\(x_{j-1}\) would be a nonzero constant and the circle of radius
\(b_{j-1}>0\) would rotate at that angular speed. The equation
for \(x_{j-1}\) gives \(w_{j-1}h_{j-1}=x_{j-2}^2\ge0\), where
the incoming square is zero at the bottom. This contradicts the
change of sign of the product on a rotating circle.
All primaries there vanish. The invariants then force \(H=1,\ w_j=0,\
h_j=\pm b_j\). Indeed the second invariant reads
\(H^2=1+\sum_j(b_j^2-h_j^2)\ge1\) on the limit set, whereas
\(H\le1\) along the orbit; hence \(H=1\) and \(h_j^2=b_j^2\) for every
\(j\). The omega-limit is connected and hence is one point.
In local coordinates on these invariant circles and the last energy
surface the stable and unstable transverse rates are \(-1,\pm b_j\),
so convergence is exponential. To see the rates, linearize at the
limit point: \(x_{d-1}'=-x_{d-1}\) up to quadratic terms, and for
\(j\le d-2\) the pair \((x_j,w_j)\) satisfies
\(x_j'=-h_j(\infty)w_j\), \(w_j'=-h_j(\infty)x_j\) up to quadratic
terms, with eigenvalues \(\pm b_j\). The weighted argument at the
end of the proof of Lemma~\ref{lem:tree-orbit} then applies verbatim
to the vector of the \(x_j\), \(j<d\), and the \(w_j\). In the full complex graph the two
top donors and the private donor pairs are disjoint real-linear
hyperbolic pairs at this endpoint; all axes and the unused mixed
heads are centers. Phases of the endpoints may be prescribed: assign
sibling phases recursively at each primary triangle, and assign each
private partner phase after its endpoint phase is fixed.

For root \((P,0)\) put \(\eta=L/P^2\) and \(r=\eta/R^2\).
A plane wave of amplitude \(rx_v\) at \(k_v=(1,0)+Rv\) has \(H^2\) norm
\(r\langle k_v\rangle^2|x_v|=\eta|v|^2|x_v|(1+O(R^{-1}))\), which is
\(L(1+O(R^{-1}))\) at the root when \(H=1\).
The relative fourth-power weight of two children is
\(z^2+(1-z)^2\le c<1\) by balance: if the children of \(m\) have
squared lengths \(z|m|^2\) and \((1-z)|m|^2\), then
\(|p|^4+|q|^4=(z^2+(1-z)^2)|m|^4\) and \(p_0\le z\le1-p_0\).
Thus the physical \(H^2\) cost
at launch of the bottom layer is at most
\(L(1+\sum b_j^2)^{1/2}c^{d/2}+o_R(1)\), by the value of \(x_0(0)\)
recorded above. The sum of the fourth
weights at level \(k\) from the top is at most \(P^4c^k\).
The contribution of all primary levels at all times is therefore
\(O(L)\) in norm: the second invariant gives
\(x_j^2\le1+\sum b_i^2\) for every \(j\), so the squared norm is at
most \(\eta^2P^4\sum_kc^k(1+\sum b_i^2)(1+o_R(1))
\le L^2(1+\sum b_i^2)(1-c)^{-1}(1+o_R(1))\).
Choose \(d\) first and then all \(b_j>0\) small
enough to make the cost of every private partner and endpoint as
small as desired.

\subsection{Constructing the driving coefficients on one axial torus}
At \(s=2\) a shift of order \(R^{-2}\) persists on the transfer time
scale. The finite system in the next subsection needs a separate
driving harmonic for each primary or private triangle of the fixed
tree. We make the following simultaneous choice.

\smallskip
\noindent\emph{Choice of the axial torus.}
All previously locked frequencies keep their prescribed values.
At the frozen frequencies the group of harmonics selected for each
triangle has modulus \(\lambda(1+o(1))\) at the space-time key
\((1,1)\), as \(\lambda\downarrow0\) with the passive data fixed.
For each fixed sufficiently small \(\lambda\), take \(R\) large.
On the resulting terminal torus the weak rates for row \(j\)
supply its prescribed shift \(d_j/R^2+O(R^{-4})\), where \(d_j\)
is the defect of that triangle at the frozen base. The same torus
has the prescribed hot axial actions up to \(o(r^2)\).
The disjoint weak sites allow us to fit these shifts while
keeping the strong rates fixed.

\begin{proof}
\emph{Positive actions and the frozen frequency constraints.}
Write the old gauged frequencies as \(n^2+\tau c_n\), \(c_n\in\mathbb Q\),
with a common denominator making all old numerators \(1\bmod 5\).
Here \(0<\tau\ll1\) is Diophantine relative to \(1\); countably many
additional excluded values are specified below. Start for example with
two distinct occupied sites of slopes \(-4,-4\). Positivity follows
from \(D\chi(0)=-\mathrm{Id}\).

For each desired triangle take disjoint sets of new sites, with words
\[
 h=e_{1-T}+e_{1+3T}+e_{1+4T}-e_{1+T}-e_{1+5T},
 \qquad d_*=e_K-3e_{K+1}+3e_{K+2}-e_{K+3}.
\]
Their triples (charge, spatial momentum, free frequency) are
\((1,1,1)\) and \((0,0,0)\), respectively. Indeed for \(h\) the
charge is \(3-2=1\), and the signed sums of the sites and of their
squares are
\[
\begin{aligned}
 (1-T)+(1+3T)+(1+4T)-(1+T)-(1+5T)&=1+(-1+3+4-1-5)T=1,\\
 (1-T)^2+(1+3T)^2+(1+4T)^2-(1+T)^2-(1+5T)^2
   &=1+2(-1+3+4-1-5)T\\
   &\quad+(1+9+16-1-25)T^2=1 .
\end{aligned}
\]
The coefficients \((1,-3,3,-1)\) of \(d_*\) are those of a third
difference, which annihilates \(1\), \(n\) and \(n^2\). At the old torus each
of the five weak frequencies contributes \(O(T^{-2})\). On the
four moderate sites the coefficients of \(d_*\) take a third
difference. Applied to \(n^{-j}\) it is \(O(K^{-j-3})\), by
three applications of the mean value theorem. The remainder in
the empty-site expansion is already \(O(K^{-5})\). Thus
\(|(h+d_*)\cdot\chi|\le C(K^{-5}+T^{-2})\). Use comparable,
arbitrarily large \(K\)'s and choose \(T\ge K^3\), subject to
the congruences below, with mutually disjoint supports.
Let \(J\) bound these defects, of the displayed order. Set positive
actions \(L_0J\) on the four moderate sites in each row, initially
zero at the five weak sites. Lock the old rates, and vary one unit
moderate pivot per row to obtain \((h+d_*)\cdot\chi=0\) simultaneously.
Indeed \(D\chi+\mathrm{Id}\) is compact on \(\ell^2\) at the old torus
and \(D\chi\) is close to \(-\mathrm{Id}\) throughout a common small
ball \cite[Theorem~18.1 and Corollary~18.2]{Molnar} (see
Appendix~\ref{app:molnar-four} for \(p=4\)). Eliminating the
old coordinates therefore leaves the fixed high pivot matrix close to
its diagonal principal part. In fact, let \(D=D\chi\) at the old torus and
write \(D_{FE}\) for its block from the new sites \(E\) to the old sites \(F\),
with the other blocks defined likewise. The derivative on the new sites
after locking the old rates is
\[
 D_{EE}-D_{EF}D_{FF}^{-1}D_{FE}.
\]
The projections onto \(E\) converge strongly to zero as the sites recede.
Compactness of \(D+\mathrm{Id}\) gives
\(\|D_{EE}+\mathrm{Id}_E\|+\|D_{EF}\|+\|D_{FE}\|=o(1)\).
Since \(D_{FF}^{-1}\) is bounded, the same \(o(1)\) estimate holds for
the difference of the locked derivative from \(-\mathrm{Id}_E\).
The uniform initial increment cancels in
each principal moderate row: the principal part \(-\mathrm{Id}_E\)
changes \(d_*\cdot\chi\) by \(-L_0J(1-3+3-1)=0\). A change \(\delta\)
of the action at a pivot site where \(d_*\) has coefficient \(\pm1\)
(for instance \(K\)) changes the same row by \(\mp\delta(1+o(1))\)
and the other rows by \(o(|\delta|)\). The required
\(|\delta|\) is therefore at most \(J+L_0J/10\) up to a factor
\(1+o(1)\), which is less than the initial action \(L_0J\) once
\(L_0\ge2\). Choose \(L_0\) large, then the sites
large enough that the compressed error is less than \(L_0J/10\)
per row. The inverse function theorem, or contraction against the
unit pivots, keeps all moderate actions positive. Their \(H^2\) bill
is \(O(K^2\sqrt J)=o(1)\) per site. The independent nonpivot rates
may now be chosen in \(\tau\mathbb Q\) by density.

\smallskip
\noindent\emph{The pivot adjustment in detail.}
We record the contraction behind the last paragraph. Let
\(N_{\rm r}\) be the number of rows, one for each primary or private
triangle of the fixed tree; it is fixed before the sites are chosen.
In row \(j\) let \(M_j=\{K_j,\ldots,K_j+3\}\) be the moderate sites,
\(W_j\) the five weak sites, \(p_j=K_j\) the pivot, where
\(d_{*,j}\) has coefficient \(+1\), and \(w_j=h_j+d_{*,j}\). The
supports of the \(w_j\) are disjoint and
\(\|w_j\|_{\ell^2}^2=5+(1+9+9+1)=25\), so the row operator
\(Hx=(w_j\cdot x)_{j\le N_{\rm r}}\) has norm
\(\|H\|_{\ell^2\to\ell^2}=5\). Let \(\chi^{\rm lk}\) denote the map
from new actions \(x\in\ell^2(E)\), supported on the moderate sites,
to the new-site components of \(\chi\), after the old actions have
been solved from the locked old rates by the implicit function
theorem. It is analytic on a fixed small ball
\(\|x\|_{\ell^2}\le\varrho\), and its derivative at \(x=0\) is the
locked derivative displayed above, which we write as
\(-\mathrm{Id}_E+\mathcal K\) with \(\nu\coloneqq\|\mathcal K\|=o(1)\)
as the sites recede. Fix \(C_2\) with \(\|D^2\chi^{\rm lk}\|\le2C_2\)
on that ball. The Taylor remainder
\(\mathcal R(x)=\chi^{\rm lk}(x)-\chi^{\rm lk}(0)+x-\mathcal Kx\)
then satisfies
\[
 \|\mathcal R(x)\|\le C_2\|x\|^2,\qquad
 \|\mathcal R(x)-\mathcal R(x')\|\le2C_2\max(\|x\|,\|x'\|)\,\|x-x'\| .
\]
Put \(a=L_0J\), \(u=a\sum_j\mathbf 1_{M_j}\) and
\(P\delta=\sum_j\delta_je_{p_j}\) for \(\delta\in\mathbb R^{N_{\rm r}}\).
Since \(w_j\cdot\mathbf 1_{M_j}=1-3+3-1=0\) and
\(w_j\cdot e_{p_k}=\delta_{jk}\), we have \(Hu=0\), \(HP\delta=\delta\),
and
\[
 H\chi^{\rm lk}(u+P\delta)=\rho^0-\delta+H\mathcal K(u+P\delta)
   +H\mathcal R(u+P\delta),
\]
where \(\rho^0=H\chi^{\rm lk}(0)\) collects the initial defects of the
rows, so that \(\|\rho^0\|_\infty\le J\). The row equations
\((h_j+d_{*,j})\cdot\chi=0\) are therefore equivalent to the fixed point
problem
\[
 \delta=\mathcal T(\delta)\coloneqq\rho^0+H\mathcal K(u+P\delta)
   +H\mathcal R(u+P\delta).
\]
Let \(L_0\ge2\) and \(\|\delta\|_\infty\le b\coloneqq J+a/10\). Then
\(b\le a\), and \(\|u\|_{\ell^2}=2a\sqrt{N_{\rm r}}\),
\(\|P\delta\|_{\ell^2}\le\sqrt{N_{\rm r}}\|\delta\|_\infty\) give
\(\|u+P\delta\|_{\ell^2}\le3a\sqrt{N_{\rm r}}\). Using
\(\|\cdot\|_\infty\le\|\cdot\|_{\ell^2}\) on \(\mathbb R^{N_{\rm r}}\),
we obtain
\[
\begin{aligned}
 \|\mathcal T(\delta)\|_\infty
   &\le J+5\nu\cdot3a\sqrt{N_{\rm r}}+5C_2\cdot9a^2N_{\rm r}
    =J+\bigl(15\nu\sqrt{N_{\rm r}}+45C_2N_{\rm r}a\bigr)a,\\
 \|\mathcal T(\delta)-\mathcal T(\delta')\|_\infty
   &\le\bigl(5\nu\sqrt{N_{\rm r}}+30C_2N_{\rm r}a\bigr)
       \|\delta-\delta'\|_\infty .
\end{aligned}
\]
If
\[
 \nu\sqrt{N_{\rm r}}\le\frac1{300},\qquad C_2N_{\rm r}L_0J\le\frac1{900},
\]
then \(\mathcal T\) maps the box \(\|\delta\|_\infty\le b\) into itself
and contracts with constant at most \(1/60+1/30=1/20\). Both conditions,
and the ball condition \(3a\sqrt{N_{\rm r}}\le\varrho\), hold for all
sufficiently large sites, because \(N_{\rm r}\) and \(L_0\) are fixed while
\(\nu\to0\) and \(J\le C(K^{-5}+T^{-2})\to0\). The unique fixed point
has \(|\delta_j|\le J+L_0J/10\). Consequently every pivot action is at
least \(L_0J-J-L_0J/10\ge4J/5>0\), the other moderate actions equal
\(L_0J\), and the weak actions vanish. Since \(T\ge K^3\) gives
\(J\le2CK^{-5}\), every moderate site \(n\), \(|n|\le K+3\), costs
\(\langle n\rangle^2\sqrt{I_n}\le C'K^2\sqrt{L_0J}=O(K^{-1/2})\) in
\(H^2\).

For the density statement let the three nonpivot actions of each row
vary, \(I_n=L_0J+\vartheta_n\) with \(|\vartheta_n|\le a/100\). The
equation acquires the term \(-H\vartheta\), and
\(\|H\vartheta\|_\infty\le(3+3+1)a/100\). The same estimates on the box
\(\|\delta\|_\infty\le J+17a/100\) give a contraction with the same
constant, since \(\|u+\vartheta+P\delta\|_{\ell^2}\le3a\sqrt{N_{\rm r}}\)
there. Hence the solution \(\delta(\vartheta)\) is analytic, the pivot
actions stay above \((83L_0/100-1)J>0\). To bound the derivative of
\(\delta(\vartheta)\), use Euclidean norms in the differentiated fixed
point equation. At \(x=u+\vartheta+P\delta\), the bounds above give
\(\|H(\mathcal K+D\mathcal R(x))\|_{2\to2}
\le5(\nu+6C_2a\sqrt{N_{\rm r}})\le1/20\).
Since \(P\) is an isometry and \(\|H\|_{2\to2}=5\), inversion in the
pivot variables gives
\(\|D\delta(\vartheta)\|_{2\to2}\le(20/19)(5+1/20)<6\).
The nonpivot rates are \(n^2\) plus the nonpivot components of
\(\chi^{\rm lk}(u+\vartheta+P\delta(\vartheta))\). Since
\(-\mathrm{Id}_E\) has no entries between distinct sites, the pivot
readjustment enters their derivative only through \(\mathcal K\) and
\(D\mathcal R\); this derivative is therefore
\(-\mathrm{Id}+O(\nu+C_2a\sqrt{N_{\rm r}})\), and it differs from
\(-\mathrm{Id}\) by at most \(7(\nu+6C_2a\sqrt{N_{\rm r}})\le7/100\).
By the inverse function theorem in these \(3N_{\rm r}\) variables, the
nonpivot rates fill a neighborhood of their current values. We choose
them of the form \(n^2+\tau c_n\) with \(c_n\in\mathbb Q\), arbitrarily
close to the current values, which keeps all actions positive.

\emph{Nonvanishing of the driving coefficients.}
We first compute a Taylor coefficient at the vacuum and then follow
the corresponding angular coefficient to the occupied background.
For disjoint
positive and negative supports, with multisets
\(\{p_1,\ldots,p_{m+1}\},\{q_1,\ldots,q_m\}\), the coefficient of
minimal degree in the physical inverse chart is
\[
 (-1/2)^m\sum_{\rm distinct\ orders}
       \prod_{\alpha=1}^m
       \frac{1}{(p_\alpha-q_\alpha)(p_{\alpha+1}-q_\alpha)}.
 \tag{H1}\label{eq:h2path}
\]
This formula also applies at a bare kinetic resonance. Put \(J_0=(i/2)
\operatorname{diag}(1,-1)\), \(C_{pq}=1/(p-q)\), and let \(a_p,b_q\) be independent small constants. With \(t_1=x\),
the diagonal entries of \(A_*,B_*\) are respectively
\(A_p=a_p\exp(i\sum_\ell t_\ell p^\ell)\) and
\(B_q=b_q\exp(-i\sum_\ell t_\ell q^\ell)\).
The upper field is \(\mathbf1^\top
(1+A_*CB_*C^\top/2)^{-1}A_*\mathbf1\).
To verify the identity, introduce the rational matrix
\[
 D(z)=1+\sum_p\frac{r_p e_2^\top}{z-p}
           +\sum_q\frac{s_q e_1^\top}{z-q},\quad
 r_p=\frac{iA_p}{\sqrt2}\left(e_1+\sum_q\frac{s_q}{p-q}\right),\quad
 s_q=\frac{-iB_q}{\sqrt2}\left(e_2-\sum_p\frac{r_p}{p-q}\right).
\]
Multiplication by the inverse elementary nilpotent pole factor
removes each pole and shows \(\det D=1\). For
\(D\exp(\sum_\ell t_\ell z^\ell J_0)\), its \(\ell\)th logarithmic
time derivative has no pole, by the displayed phases, and equals the
polynomial part of \(z^\ell D J_0 D^{-1}\).
Writing \(W=D J_0 D^{-1}\), the identities \(W^2=-1/4\) and
\(W_x=[zJ_0+[D_1,J_0],W]\), where \(D_1\) is the coefficient of
\(z^{-1}\), recursively determine its diagonal and off-diagonal
Laurent coefficients. These are the homogeneous NLS hierarchy
matrices. The upper and lower off-diagonal entries of \([D_1,J_0]\)
are \(u/\sqrt2,v/\sqrt2\); the reversed second flow is
\(iu_t=-u_{xx}+u^2v\). The recursion and this choice of the second
time are derived in Section~\ref{rate:sec:axial-mark}. On simultaneous Birkhoff coordinates each
hierarchy Hamiltonian depends only on the actions
\cite[Theorem~1.1(3)]{GKsym}; apply its complexified Taylor identity
at zero on the finitely many input variables. Each monomial in any
action must have all signed power sums zero. The Vandermonde matrix
on the disjoint supports, including degree zero, excludes every
nonconstant monomial. With all actions zero the signed frequencies
in every hierarchy time are their vacuum values. The same Vandermonde
argument after adjoining an output site leaves only its prescribed
linear input in that Birkhoff coordinate. Expanding the upper matrix
inverse now gives \eqref{eq:h2path}. In more detail, a monomial
\(\prod_pa_p^{\alpha_p}\prod_qb_q^{\beta_q}\) has phase
\(\exp\bigl(i\sum_\ell t_\ell(\sum_p\alpha_pp^\ell-\sum_q\beta_qq^\ell)\bigr)\),
where \(\ell=0\) records the charge. An action is constant in every
time, so each of its monomials has
\(\sum_p\alpha_pp^\ell-\sum_q\beta_qq^\ell=0\) for all \(\ell\ge0\);
the Vandermonde determinant on the distinct points of the two
supports then gives \(\alpha=\beta=0\). With vanishing actions the
Birkhoff coordinate at a site \(n\) has phase
\(\exp(i\sum_\ell t_\ell n^\ell)\), and the same determinant with
\(n\) adjoined leaves only the monomial \(a_n\) when \(n\) is a
positive site and none otherwise; the conjugate coordinate at a
negative site \(q\) is likewise \(b_q\). The normalizations are
fixed by the linear part \(\mathbf1^\top A_*\mathbf1=\sum_pA_p\) of
the upper field. Finally the term of degree \(m\) in \(B_*\) of the
Neumann series is
\[
 \bigl(-\tfrac12\bigr)^m\mathbf1^\top(A_*CB_*C^\top)^mA_*\mathbf1
 =\bigl(-\tfrac12\bigr)^m\sum A_{p_1}\prod_{\alpha=1}^m
   \frac{B_{q_\alpha}A_{p_{\alpha+1}}}
        {(p_\alpha-q_\alpha)(p_{\alpha+1}-q_\alpha)},
\]
summed over all sequences \(p_1,\ldots,p_{m+1}\) and
\(q_1,\ldots,q_m\) in the two supports. Its spatial frequency is
\(\sum_\alpha p_\alpha-\sum_\alpha q_\alpha\). The sequences that are
orderings of two prescribed multisets carry the same monomial, and
their sum is \eqref{eq:h2path}.

Two special cases agree with the direct computations in the proof of
Proposition~\ref{prop:bank}. For \(m=1\) and positive sites \(p_1\ne p_2\),
formula \eqref{eq:h2path} gives
\(-1/((p_1-q_1)(p_2-q_1))=2/D(p_1,q_1,p_2)\), since
\(D(a,b,c)=-2(a-b)(c-b)\); this is the cubic coefficient \(1/D\) per
ordered triple, summed over the two orders of the positive inputs. For
the positive multiset \(\{0,4,5\}\) and the negative multiset \(\{2,6\}\)
it gives \(-1/9\), the quintic coefficient obtained there from the third
flow. This second agreement is a nontrivial check, since the computation
in that proof uses the third flow and not the dressing construction.

For the congruences \(K\equiv10,\ T\equiv31\pmod{1009}\), the
vacuum coefficient in \eqref{eq:h2path} for \(h+d_*\) is
\(54\bmod1009\). To compute it, use an initial
positive site \(i\), remaining positive multiplicities \(a\) and
negative multiplicities \(b\) in the recursion
\[
 W(i,0,0)=1,\quad W(i,a,b)=
 \sum_{b_j a_k>0}\frac{W(k,a-e_k,b-e_j)}
 {(p_i-q_j)(p_k-q_j)} .
\]
Use \(p=(1-T,1+3T,1+4T,K,K+2)\), \(q=(1+T,1+5T,K+1,K+3)\),
\(a=(1,1,1,1,3)\), \(b=(1,1,3,1)\). In the displayed order of
the positive sites the five values \(W(i,a-e_i,b)\) modulo \(1009\)
are \(678,979,232,347,211\), whose sum is \(429\bmod1009\).
Since \(2^{-6}=741\bmod1009\), the coefficient is
\(429\cdot741=54\bmod1009\). All opposite differences are units.
The same congruences allow
arbitrarily large groups.

The congruence gives a nonzero starting value for continuation.
On the occupied torus different old angular labels can contribute
to the same Fourier coefficient, so we must continue their sum
at the prescribed phase origins. Fix those origins (zero suffices).
For each finite recipe of old sites and rational slopes, new groups
and free moderate slopes, replace \(\tau\) by a scalar \(\zeta\).
At passive weak sites solve the old and free rate equations and the
row equations through zero. Divide the selected angular coefficient
by the required products of the new radii and sum the aliases
differing by old labels \(\beta\) with
\(\beta\cdot\mathbf1=\beta\cdot n=\beta\cdot n^2=
\beta\cdot c_{\rm old}=0\), at the prescribed origins. At zero every
nonzero such label has an extra old radius and disappears. The
vacuum value is \eqref{eq:h2path}. To justify continuation, in the
independent frequency coordinates replace each pivot coordinate by
\(\chi_{\rm piv}+h\cdot\chi_{\rm weak}/(d_*)_{\rm piv}\).
On the subspace of zero weak actions the derivative differs from
\(-\mathrm{Id}\) by at most \((1+\sqrt5)\|D\chi+\mathrm{Id}\|\),
by disjointness. The inverse exists on a common small complex
\(\ell^2\) ball: the prescribed vector along the segment is small
because its endpoint has small total action and all row operators
have uniformly bounded norm. The local inverse chart
\(\ell^4\to\mathcal F\ell^4\) is analytic near zero
\cite[Theorems~14.1 and 15.2]{Molnar}. Thus it can be evaluated
along any Cartesian square-root branch of this curve; its
\(\ell^4\) size squared is the \(\ell^2\) action size. The alias
projection is Haar averaging over the annihilator of the integer
subgroup. Factoring the required radii is analytic by the two-sign
Taylor series. Multiplying over the signs of the old square roots
gives an even holomorphic function of those radii, hence a
single-valued function of the actions, nonzero at zero. The isolated
zeros of these functions for all countably many recipes can be
excluded in the original choice of \(\tau\). At the selected point,
restore the moderate radii and denote each nonzero modulus, with the
five weak radii still divided out, by \(m_j\).

\smallskip
\noindent\emph{The continuation in detail.}
Fix a recipe \(\mathfrak r\): the old support \(F\) with its rational
slopes, the new rows with their sites and pivots \(p_j\), and rational
slopes \(c_n\) at the free moderate sites. Let \(S\) be the set of old
and moderate sites; the weak actions vanish throughout this step. For
complex actions \(x\in\mathbb C^S\) define the frequency coordinates
\[
 y_n(x)=\chi_n(x)\quad(n\in S\setminus\{p_j\}_j),\qquad
 y_{p_j}(x)=\chi_{p_j}(x)+\frac{h_j\cdot\chi_{W_j}(x)}{(d_{*,j})_{p_j}} .
\]
The required values are linear in the scalar \(\zeta\): \(y_n=\zeta c_n\)
at the old and free moderate sites and, by the row equation
\((h_j+d_{*,j})\cdot\chi=0\),
\(y_{p_j}=-\zeta\sum_{n\in M_j\setminus\{p_j\}}(d_{*,j})_nc_n/(d_{*,j})_{p_j}\).
Write this vector as \(\zeta v_{\mathfrak r}\). On \(\mathbb C^S\) one has
\(Dy=-\mathrm{Id}_S+\mathcal K'\) with
\[
 \mathcal K'x=P_S(D\chi+\mathrm{Id})x
   +\sum_j\frac{h_j\cdot\bigl[(D\chi+\mathrm{Id})x\bigr]_{W_j}}
              {(d_{*,j})_{p_j}}\,e_{p_j},
\]
because \(\mathrm{Id}\) has no entries from \(S\) to the weak sites. As
the \(W_j\) are disjoint and \(\|h_j\|_{\ell^2}=\sqrt5\), the second
operator has norm at most \(\sqrt5\|D\chi+\mathrm{Id}\|\), so that
\(\|\mathcal K'\|\le(1+\sqrt5)\|D\chi+\mathrm{Id}\|\). Choose a complex
\(\ell^2\) ball \(\|x\|\le\varrho_0\) on which
\(\|D\chi+\mathrm{Id}\|\le1/(2+2\sqrt5)\); it does not depend on the
recipe. For \(\|y_*\|<\varrho_0/2\) the map \(x\mapsto x+y(x)-y_*\) has
Lipschitz constant at most \(1/2\) on the closed ball and maps it into
itself, since \(y(0)=0\). Hence \(y\) has an analytic inverse on
\(\|y_*\|<\varrho_0/2\) with values in the ball, uniformly in the
recipe. Let \(\mathcal D(\mathfrak r)=\{\zeta\in\mathbb C:
|\zeta|\,\|v_{\mathfrak r}\|<\varrho_0/2\}\) and
\(x(\zeta)=y^{-1}(\zeta v_{\mathfrak r})\), an analytic curve with
\(x(0)=0\). If \(\mathfrak r\) is the recipe realized at \(\tau\) by a
torus with actions \(x_\tau\), then
\(\tau\|v_{\mathfrak r}\|=\|y(x_\tau)\|\le\tfrac32\|x_\tau\|\); as the
total action is small, \(\tau\in\mathcal D(\mathfrak r)\) and
\(x(\tau)=x_\tau\) by uniqueness in the ball.

To evaluate the field along the curve, lift every complex action to
Cartesian coordinates with real angles,
\((z_n,w_n)=(\rho_ne^{i\theta_n},x_n\rho_n^{-1}e^{-i\theta_n})\) with
\(|\rho_n|^2=|x_n|\) (and \(z_n=w_n=0\) if \(x_n=0\)), where \(w_n\)
replaces \(\bar z_n\) after complexification. Both entries have modulus
\(|x_n|^{1/2}\), so the lift has \(\ell^4\) norm
\((2\|x\|_{\ell^2}^2)^{1/4}\) and lies in the domain of the analytic
inverse chart once \(\varrho_0\) is small enough. At a new site with label \(\alpha_n\), the terms of the
Taylor series in \((z_n,w_n)\) with angular frequency \(\alpha_n\) are
\(z_n^{\alpha_n}\sum_bc_b(z_nw_n)^b\) if \(\alpha_n\ge0\), and the
analogous series with \(w_n^{-\alpha_n}\) otherwise. Dividing by
\(\rho_n^{\alpha_n}\), respectively by \((x_n/\rho_n)^{-\alpha_n}\), leaves a
power series in \(x_n\), independent of the lift; at a weak site
\(x_n=0\) and only its constant term remains. At the old sites take
\(\rho_n=x_n^{1/2}\) for a choice of the square roots, so that both
entries are \(\rho_ne^{\pm i\theta_n}\) and depend analytically on
\(\rho_n\). Let \(\Lambda\subset\mathbb Z^F\) be the group of old labels
\(\beta\) with \(\beta\cdot\mathbf1=\beta\cdot n=\beta\cdot n^2
=\beta\cdot c_{\rm old}=0\), and let \(\mathbb K\subset\mathbb T^F\) be
its annihilator, with Haar probability measure \(d\varphi\). By
Pontryagin duality \(\int_{\mathbb K}e^{i\gamma\cdot\varphi}\,d\varphi=1\)
if \(\gamma\in\Lambda\) and \(0\) otherwise. The selected label is zero
at the old sites. Hence, if \(g(\theta_F)\) is the normalized new-site
coefficient as a function of the old angles, the sum of its aliases at
the origin \(\theta_F=0\) is
\[
 \sum_{\beta\in\Lambda}\hat g(\beta)=\int_{\mathbb K}g(\varphi)\,d\varphi .
\]
The integrand is analytic in \((\rho_F,x_{S\setminus F})\), uniformly in
\(\varphi\), so the grouped coefficient
\(G_{\mathfrak r}(\rho_F,x_{S\setminus F})\) is analytic on a polydisc.
At \(\rho_F=0\), \(x_{S\setminus F}=0\), every \(\beta\ne0\) carries a
factor \(\rho_n\) with \(n\in F\) and drops out, so
\(G_{\mathfrak r}(0,0)\) is the vacuum coefficient \eqref{eq:h2path},
which is nonzero. Changing the sign of \(\rho_n\) replaces \(\theta_n\) by
\(\theta_n+\pi\), which in general changes the coefficient. The product
\(f_{\mathfrak r}=\prod_{\sigma\in\{\pm1\}^F}
G_{\mathfrak r}(\sigma\rho_F,x_{S\setminus F})\) is even in each
\(\rho_n\), and hence a power series in \(x_n=\rho_n^2\), \(n\in F\), and
in \(x_{S\setminus F}\). Thus \(\zeta\mapsto f_{\mathfrak r}(x(\zeta))\) is
holomorphic on \(\mathcal D(\mathfrak r)\) with value
\(\eqref{eq:h2path}^{2^{|F|}}\ne0\) at \(\zeta=0\), and its zeros form a
discrete, hence countable, subset of \(\mathcal D(\mathfrak r)\). There
are countably many recipes, so all these zeros are avoided by a
Diophantine \(\tau\). For the realized recipe, \(\tau\in\mathcal
D(\mathfrak r)\) and \(x(\tau)=x_\tau\); hence
\(f_{\mathfrak r}(x_\tau)\ne0\), and its factor with the positive radii,
which is the selected grouped coefficient with the radii divided out,
is nonzero.

\emph{Separating the rows and fitting the terminal rates.}
Choose the five weak radii in row \(j\) equal to
\((\lambda/m_j)^{1/5}\) and restore the old and moderate rate
conditions as above; the changes tend to zero at fixed sites when
\(\lambda\downarrow0\). The selected grouped coefficients now have
modulus \(\lambda(1+o(1))\). On the flat \((h_j+d_{*,j})\cdot c=0\)
approximate the resulting new rates by \(\tau\) times rationals.
Take the new common denominator to be \(Q_{\rm old}D_0\) with
\(D_0=1\bmod5\) divisible by distinct private primes \(\pi_j>30000\)
outside the old denominator. Give strong numerators residue \(1\bmod5\),
weak ones zero. Modulo \(\pi_j\) impose residues
\((1,13,169,2197,-2014)\) on its five weak numerators, and zero on all
others. Solve for a unit weak pivot on the flat. These conditions
are compatible by the Chinese remainder theorem and dense with
arbitrarily large denominators. All new rates are legitimate
coordinates at this positive point.

An integer-time charge-one label of total weak order at most five
has weak part exactly one \(h_j\). Indeed its scalar product with
the residues is an integer of modulus less than \(\pi_j\) divisible
by \(\pi_j\). Eliminate the last entry; all four base-\(13\) digits
have absolute value at most five and hence vanish. In detail, an
integer-time label \(l\) has \(l\cdot c=0\), because \(\tau\) is
irrational. Multiplying by the common denominator and reducing
modulo \(\pi_j\) leaves only the entries \(l_1,\ldots,l_5\) of \(l\)
on the five weak sites of row \(j\), taken in the order
\(1-T,1+3T,1+4T,1+T,1+5T\) of the displayed residues. Since
\(\sum_i|l_i|\le5\), the integer
\(l_1+13l_2+169l_3+2197l_4-2014l_5\) has modulus at most
\(5\cdot2197<\pi_j\) and is therefore zero. As \(-2014=1+13+169-2197\),
\[
 0=(l_1+l_5)+13(l_2+l_5)+169(l_3+l_5)+2197(l_4-l_5),
\]
where every bracket has modulus at most five. Reducing successively
modulo \(13\) shows that the four brackets vanish. Hence
\((l_1,\ldots,l_5)=-l_5(1,1,1,-1,-1)\), an integer multiple of the weak
part of \(h_j\), and the multiple is \(0\) or \(\pm1\) by the bound on
the order. The weak order five is then exhausted by one row. The charge
modulo five then specifies the positive sign and rules out an empty
weak part: the strong numerators are \(1\bmod5\) and the weak ones
\(0\bmod5\), so \(l\cdot c=0\) forces the strong charge to vanish
modulo five, while the total charge one equals the strong charge plus
the multiple just found. Only the multiple \(+1\) is possible.
For any finite angular cutoff avoid all other aliases
with nonzero moderate difference at key \((1,1)\): all moderate
entries are free on the flat after eliminating a weak pivot. The
remaining normalized strong tails have arbitrarily small absolute
sum by angular analyticity after factoring the five radii. Choose
their cutoff first and then the rational approximation. Over all
integer keys \((b,\nu)\) the absolute coefficient sum of these
order-five groups is \(O(\lambda)\), uniform for small \(\lambda\)
at fixed passive data; labels of weak order at least six have sum
\(O(\lambda^{6/5})\) by Cauchy estimates on a fixed weak polydisc.

\smallskip
\noindent\emph{Derivation of the two absolute sums.} Fix the passive
data (old and moderate sites, their actions and angles, and all
rates); as \(\lambda\downarrow0\) they move by \(o(1)\) and stay in a
fixed compact set, and the number \(f_p\) of passive occupied sites is
fixed. Complexify the Cartesian pair at each weak site
\(n\in W\coloneqq\bigcup_jW_j\) into independent variables
\((\zeta_n,\zeta_n^*)\), with \(\zeta_n^*=\overline{\zeta_n}\) on the
real slice. There are \(\rho_w>0\), \(a_0>0\) and \(C_\Psi\), independent
of \(\lambda\), such that the physical field \(\Psi(z)\in L^2(\T)\) is
analytic in the passive angles on \(|\operatorname{Im}\varphi|<a_0\)
and in \((\zeta,\zeta^*)\) on the polydisc
\(|\zeta_n|,|\zeta_n^*|<\rho_w\), with \(\|\Psi(z)\|_{L^2}\le C_\Psi\)
there. Indeed the number \(10N_{\rm r}\) of weak variables is fixed,
so for small \(\rho_w\) this polydisc lies in the complex \(L^2\)
chart ball. Expand \(\Psi(z)\) in the monomials
\(\zeta^\alpha\zeta^{*\beta}e^{-il'\cdot\varphi}\), \(l'\in\Z^{f_p}\).
By translation covariance each monomial carries a single spatial
exponential \(e^{ibx}\), with \(b\) determined by \(l'\) and
\(\alpha-\beta\), and a single time frequency. Hence each monomial
contributes to exactly one key \((b,\nu)\), and the absolute
coefficient sum over all keys is at most the sum of the absolute values
of all monomial coefficients. Cauchy's estimate in the weak variables
on the polydisc and in the passive angles on the strip of width
\(a_0\) bounds the coefficient \(c_{\alpha\beta l'}\) of
\(\zeta^\alpha\zeta^{*\beta}e^{-il'\cdot\varphi}e^{ibx}\) by
\[
 |c_{\alpha\beta l'}|\le C_\Psi\,\rho_w^{-|\alpha|-|\beta|}e^{-a_0|l'|_1/2}.
\]
On the torus every weak variable of row \(j\) has modulus
\(r_{w,j}=(\lambda/m_j)^{1/5}\le C_r\lambda^{1/5}\), because the moduli
\(m_j\) are bounded below at fixed passive data. A label of weak order
\(w\) arises only from monomials of weak degree
\(k=|\alpha|+|\beta|\ge w\). There are at most
\(\binom{k+10N_{\rm r}-1}{k}\le(k+1)^{10N_{\rm r}}\) monomials of
degree \(k\), and \(\sum_{l'}e^{-a_0|l'|_1/2}=\coth(a_0/4)^{f_p}\).
Hence, if \(x\coloneqq C_r\lambda^{1/5}/\rho_w\le1/2\),
\[
 \sum_{k\ge w}\ \sum_{|\alpha|+|\beta|=k}\ \sum_{l'}
   |c_{\alpha\beta l'}|\,(C_r\lambda^{1/5})^{k}
 \le C_\Psi\coth(a_0/4)^{f_p}\,x^w\sum_{i\ge0}(w+i+1)^{10N_{\rm r}}2^{-i}
 =C_w\,\lambda^{w/5},
\]
with \(C_w\) independent of \(\lambda\). The choice \(w=5\) bounds all
labels of weak order at least five, in particular the order-five
groups, by \(O(\lambda)\), and \(w=6\) gives \(O(\lambda^{6/5})\) for
weak order at least six. In particular the absolute sum of the order-five groups is at most \(C_w\lambda\) with \(w=5\); this is the bound \(C_0\lambda\) used in the proof of Lemma~\ref{lem:common-endpoint-green}.

At this point all the groups have been selected at the frozen
frequency vector \(\Omega^0\). To fit the prescribed defects, shift
only its weak rates to obtain \(\Omega^\infty\), keeping all
strong rates fixed, so that
\[
 h_j\cdot(\Omega^\infty-\Omega^0)=d_j/R^2+O(R^{-4}),
\]
where \(d_j\) is the intended defect evaluated at the frozen base.
Simultaneously impose the hot axial actions \(r^2|a_i|^2\) at
\(K_i=1+RP_i\) using the compressed frequency inverse. Here \(a_i\)
is the prescribed terminal amplitude at \(K_i\) in the normalized finite orbit above. Finally
approximate the weak rates and the hot axial rates, including the root
rate, while holding every strong rate at its frozen value. Choose the
accuracy so that the hot actions change by \(o(r^2)\), all occupied
actions stay positive, and the displayed shift retains its accuracy.
Enlarging the common denominator by a factor congruent to one modulo
five preserves the strong numerators' residues and allows every new
numerator to be chosen congruent to one modulo five for the next stage.

\smallskip
\noindent\emph{The last fit in detail.}
Let \(S\) be the strong sites, \(W=\bigcup_jW_j\) the weak sites,
\(\mathsf w_j\in W_j\) the weak pivot of row \(j\), so that
\((h_j)_{\mathsf w_j}=\pm1\), and \(\mathcal H\) the hot sites \(K_i\). We
solve the mixed problem in which the actions on \(\mathcal H\) and the
rates on \(S\cup W\) are prescribed,
\[
 I_{K_i}=r^2|a_i|^2,\qquad
 \Omega_n=\Omega^0_n\ \ (n\in S\cup W,\ n\ne\mathsf w_j),\qquad
 \Omega_{\mathsf w_j}=\Omega^0_{\mathsf w_j}+(h_j)_{\mathsf w_j}\,d_j/R^2 ,
\]
all other actions being zero. Then
\(h_j\cdot(\Omega-\Omega^0)=(h_j)_{\mathsf w_j}^2d_j/R^2=d_j/R^2\) and the
strong rates are unchanged. At fixed hot actions the derivative of the
rates on \(S\cup W\) with respect to the actions there is
\(-\mathrm{Id}+P(D\chi+\mathrm{Id})P\), where \(P\) projects to
\(S\cup W\), and the hot actions move these rates by at most
\(\|D\chi+\mathrm{Id}\|\,\|I_{\mathcal H}\|_{\ell^2}=O(r^2)=O(R^{-4})\).
The quantitative inverse function theorem used above therefore gives a
unique solution near the torus of the previous step, with
\[
 \|I_{S\cup W}-I^0_{S\cup W}\|_{\ell^2}
   \le2\bigl(R^{-2}\|(d_j)_j\|_{\ell^2}+CR^{-4}\bigr)=O(R^{-2}).
\]
The weak actions are \((\lambda/m_j)^{2/5}(1+o(1))\) by the choice of the
radii, so their relative change is \(O(\lambda^{-2/5}R^{-2})\), which
tends to zero as \(R\to\infty\) at fixed \(\lambda\); the strong actions
are bounded below at the fixed sites and change by \(O(R^{-2})\). The
normalized grouped coefficients are analytic in all actions (the hot
sites enter them with label zero), and the product of the five weak
radii changes by a factor \(1+O(\lambda^{-2/5}R^{-2})\). The selected
groups therefore keep modulus \(\lambda(1+o(1))\).

The congruences imposed at \(\Omega^0\) in the paragraph on separating
the rows are compatible with the flats. In the order
\(1-T,1+3T,1+4T,1+T,1+5T\) the weak part of \(h_j\) is
\((1,1,1,-1,-1)\), and
\[
 (1,1,1,-1,-1)\cdot(1,13,169,2197,-2014)=1+13+169-2197+2014=0,
\]
while \(d_{*,j}\) has zero residues modulo \(\pi_j\) and coefficient sum
\(1-3+3-1=0\) against the strong residues \(1\bmod5\); the weak residues
are \(0\bmod5\). For the final approximation write the frozen strong
rates as \(n^2+\tau c_n\) with common slope denominator \(Q\). On passing
to denominator \(QM\), \(M\equiv1\bmod5\), every strong numerator is
multiplied by \(M\) and retains its residue \(1\bmod5\). Only the weak
and hot slopes are rounded, to values \((5k+1)/(QM)\). The spacing of
the corresponding rates is \(5\tau/(QM)\), so their errors can be made
\(o(r^2)\). The frequency inverse gives action errors \(o(r^2)\), and
\(h_j\cdot(\Omega^\infty-\Omega^0)=d_j/R^2+o(R^{-4})\). In particular,
two frozen labels with the same space-time key and weak part still have
identical time frequencies, since their difference is supported on
the strong sites, whose rates were held fixed.
\end{proof}

\subsection{Finite dynamics and the full normal operator}
Use the exact Birkhoff frequencies on the terminal axes. On mobile
heads begin with \(|k_v|^2+\Xi(v)/R^2\); corrections \(O(R^{-3})\)
make each edge identity exact for its assigned row. Indeed, distribute
each parent's discrepancy between its children, starting at the root,
then solve one equation per private partner. The unused mixed-head
frequency needs no correction. The empty-axis expansion above gives
the stated accuracy on axes; the newly occupied very small actions
change the empty frequencies by \(O_V(r^2)\) after accounting for the
old action variation.

\smallskip
\noindent\emph{The carrier frequencies in detail.} Write \(\kappa_v\)
for the carrier frequency of a vertex \(v\) in the common mass gauge
and \(\kappa_v=|k_v|^2+\Xi(v)/R^2+\varepsilon_v\). In \(\Xi\) the
coefficients \(A,B\) are those of the frozen base. The terminal actions
differ from the frozen ones by \(O(R^{-2})\) by the last fit, and
\(A,B\) are smooth functions of the actions, so the terminal values of
\(\Xi(v)/R^2\) differ from these by \(O(R^{-4})\). (i) On an axial
vertex \(v=(P_i,0)\) the carrier is the exact Birkhoff frequency
\(K_i^2+\chi_{K_i}(I^\infty)\), with \(K_i=1+RP_i\) and
\(|k_v|^2=K_i^2\). The higher-order empty expansion gives
\(\chi_{K_i}=(A-B)K_i^{-2}+O(K_i^{-3})\) at the terminal torus with
its hot actions removed. The new actions, namely the hot actions \(r^2|a_i|^2\) and the induced change
\(O(r^2)\) of the other actions, change \(\chi_{K_i}\) by at most
\(\|D\chi\|_{\ell^2\to\ell^2}\|\Delta I\|_{\ell^2}=O(r^2)=O(R^{-4})\),
because \(\chi\) is Lipschitz from \(\ell^2\) to
\(\ell^2\subset\ell^\infty\) on the action ball. As
\(K_i^{-2}=(RP_i)^{-2}(1+O(R^{-1}))\) and
\(\Xi((P_i,0))=(A-B)/P_i^2\), the axial \(\varepsilon_v\) is fixed and
\(O(R^{-3})\). (ii) On a transverse head, \(n=1+Rv_x\), \(m=Rv_y\), the
true high diagonal is \(n^2+m^2+A/n^2-B/(n^2+m^2)+O(R^{-3})\), and
\(A/n^2=A/(Rv_x)^2+O(R^{-3})\),
\(B/(n^2+m^2)=B/(R^2|v|^2)+O(R^{-3})\). Hence it equals
\(|k_v|^2+\Xi(v)/R^2+O(R^{-3})\). (iii) The group of row \(j\) has
actual time frequency \(1+d_j/R^2+O(R^{-4})\): its word
\(h_j+d_{*,j}\) has free frequency \(1\) and zero \(\chi\)-sum at the
frozen point, and only the weak rates moved, by
\(h_j\cdot(\Omega^\infty-\Omega^0)=d_j/R^2+O(R^{-4})\). For an intended
edge \(e=(v,w;v+w)\) assigned to row \(j\) we have
\(|k_v|^2+|k_w|^2-|k_{v+w}|^2=1\) and
\(\Xi(v)+\Xi(w)-\Xi(v+w)=d_j\). The exact edge identity
\(\kappa_v+\kappa_w-\kappa_{v+w}=1+d_j/R^2+O(R^{-4})\) therefore reads
\[
 \varepsilon_v+\varepsilon_w-\varepsilon_{v+w}=\delta_e,\qquad
 \delta_e=O(R^{-4}),
\]
where \(\delta_e\) is the remainder of the actual group frequency. Solve
these equations from the root: \(\varepsilon_{\rm root}\) is given by
(i); for a primary parent \(m=p+q\) with \(\varepsilon_m\) known put
\(\varepsilon_p=0\) and \(\varepsilon_q=\varepsilon_m+\delta_e\); for a
private triangle \((v,J(v);A(v))\) the values \(\varepsilon_v\) and
\(\varepsilon_{A(v)}\) are known, and
\(\varepsilon_{J(v)}=\varepsilon_{A(v)}-\varepsilon_v+\delta_e\). Every
nonroot primary occurs as a child in exactly one primary triangle,
and each \(J(v)\) belongs to exactly one private triangle. Thus each
mobile correction is assigned only once. Along any root path at most \(d\) terms are added, so all
mobile corrections are \(O_d(R^{-3})\). By (ii) the mobile carriers
then differ from the true high diagonals by \(O(R^{-3})\). The mixed
head \(T_m\) lies in no intended triangle and keeps
\(\varepsilon_{T_m}=0\).

Use all angular harmonics with frozen key \((1,1)\) in the finite
equation, on all triangles including the mixed triangles.
For order five the extra actual phase is determined by \(h_j\) alone:
the strong rates were kept exact. In slow time \(\theta=r\lambda t\)
the intended coupling is constant with magnitude \(2+o(1)\), after
the assigned phase rotations. Every other order-five harmonic group
on that triangle, and every such group on a mixed triangle, oscillates
at frequency of modulus at least \(c_{\rm tree}/(\eta\lambda)\) by
the strict defect separations. Indeed, on a triangle \(e=(v,w;v+w)\)
with defect \(d_e\) the bare mismatch
\(|k_v|^2+|k_w|^2-|k_{v+w}|^2\) equals one, the corrected head
frequencies add \(d_e/R^2+O(R^{-3})\) to it, and the group of row
\(j'\) has actual frequency \(1+d_{j'}/R^2+O(R^{-4})\). Its physical
detuning on \(e\) is therefore \((d_e-d_{j'})/R^2+O(R^{-3})\), which
in slow time becomes
\[
 \frac{d_e-d_{j'}}{R^2r\lambda}+O\Big(\frac1{R^3r\lambda}\Big)
 =\frac{d_e-d_{j'}}{\eta\lambda}\bigl(1+O(R^{-1})\bigr),
\]
because \(r=\eta/R^2\). It vanishes only for the assigned row of an
intended triangle, where the \(O(R^{-3})\) corrections were chosen
to cancel it exactly. Divide each of these oscillatory vector
fields by its signed frequency to obtain a polynomial change of
coordinates \(O_{\rm tree}(\eta\lambda)\) in \(C^1\) on a bounded
neighborhood. Substitution leaves a \(C^1\) remainder of this same
order; the undivided higher-order bank has size \(O(\lambda^{1/5})\).
(In the slow time \(\theta=r\lambda t\) every physical coupling
coefficient is divided by \(\lambda\). The order-five groups therefore
have normalized absolute sum \(O(1)\), and the labels of weak order at
least six have normalized absolute sum \(O(\lambda^{6/5})/\lambda
=O(\lambda^{1/5})\), by the two absolute sums derived above. Each label
multiplies a quadratic monomial in the finitely many amplitudes, so the
same bounds hold for the corresponding vector fields in \(C^3\) on a
bounded neighborhood.) Every term vanishes on the axial plane. The detuning correction
and endpoint rounding cost \(o_R(1)\) at fixed \(\lambda\).
The resulting equation has a solution converging exponentially to
its exact rounded axial endpoint: take a time on the ideal landing
where the tail is small, use the linearized endpoint splitting with
a positive weight strictly below its hyperbolic rates, prescribe
the stable projection there, and integrate the other components
from \(+\infty\). Small \(C^1\) error and vanishing of the field at
the exact endpoint give a contraction in the weighted tube, including
the center directions. On the preceding fixed interval the inverse
propagates backwards by Gronwall. First small \(\lambda\), then large
\(R\), make the launch as close as prescribed to the ideal one;
undo the change, which fixes the endpoint. Here and in the subsequent
insertion the exact finite solution with all these terms is used.

\smallskip
\noindent\emph{The averaging and the landing in detail.}
Write the finite equation in slow time as
\[
 X'=F_0(X)+\sum_kG_k(X)e^{i\omega_k\theta}+B(X,\theta),
\]
where \(F_0\) is the ideal vector field with coupling of magnitude
\(2\) on every intended edge, in the phase convention above. The orbit
constructed with equal couplings is parametrized here in these slow-time
units. Denote by \(c_e^{\rm id}\) and \(c_e^{\rm fit}\) the ideal and
fitted coefficients on edge \(e\) after the assigned phase rotations,
and set \(\varepsilon_{\rm fit}=\max_e|c_e^{\rm fit}-c_e^{\rm id}|\). The
coefficient estimates above give
\(\lim_{\lambda\downarrow0}\limsup_{R\to\infty}\varepsilon_{\rm fit}=0\)
along the admissible choices. We include the difference between the
fitted and ideal constant edge fields in \(B\), together with the groups
of weak order at least six and the \(o_R(1)\) corrections. The \(G_k\)
are the oscillating order-five groups, with
\(|\omega_k|\ge c_{\rm tree}/(\eta\lambda)\). Let
\(\mathcal U\) be a fixed bounded neighborhood of the ideal orbit. All
these fields are polynomial vector fields of bounded degree, and
\(\sum_k\|G_k\|_{C^3(\mathcal U)}\le C\), because the absolute coefficient
sum of the order-five groups is \(O(\lambda)\) before division by
\(\lambda\). As the tree is finite and its edge fields are quadratic,
\(\|B\|_{C^1(\mathcal U)}\le C(\lambda^{1/5}+\varepsilon_{\rm fit})+o_R(1)\).
Each of these fields is a sum of triangle terms, each a product of the
amplitudes at two vertices of a triangle with at most one axial
vertex, so all of them vanish on the axial plane \(\mathcal A\). Put
\[
 \Phi(Y,\theta)=\sum_k\frac{G_k(Y)}{i\omega_k}e^{i\omega_k\theta},
 \qquad \|\Phi\|_{C^2(\mathcal U)}\le C\eta\lambda/c_{\rm tree},
\]
and \(X=Y+\Phi(Y,\theta)\). Since
\(\partial_\theta\Phi=\sum_kG_ke^{i\omega_k\theta}\),
\[
 (1+D_Y\Phi)Y'=F_0(Y+\Phi)
   +\sum_k\bigl(G_k(Y+\Phi)-G_k(Y)\bigr)e^{i\omega_k\theta}
   +B(Y+\Phi,\theta).
\]
Inverting \(1+D_Y\Phi\) by a Neumann series gives \(Y'=F_0(Y)+S(Y,\theta)\)
with \(\|S\|_{C^1(\mathcal U)}\le
C(\eta\lambda+\lambda^{1/5}+\varepsilon_{\rm fit})+o_R(1)\).
As \(\Phi\) and \(S\) vanish on \(\mathcal A\), the change is the identity
there. Since \(\mathcal A\) is a linear subspace and
\(\|D_Y\Phi\|\le C\eta\lambda\), also \(|\Phi(Y,\theta)|\le
C\eta\lambda\,\operatorname{dist}(Y,\mathcal A)\).

Let \(Y_{\rm id}\) be the ideal orbit, converging exponentially to
\(p_{\rm id}\in\mathcal A\), and let \(p_*\in\mathcal A\) be the exact
rounded endpoint, \(|p_*-p_{\rm id}|=o_R(1)\). Put \(\mathsf A=DF_0(p_*)\).
By the description of the endpoint in the preceding subsection, the
phase space splits \(\mathsf A\)-invariantly as \(E_s\oplus E_u\oplus E_c\),
with
\[
 \|e^{t\mathsf A}P_s\|+\|e^{-t\mathsf A}P_u\|\le Ce^{-\mu_*t}\ (t\ge0),
 \qquad \|e^{t\mathsf A}P_c\|\le C(1+|t|)^{m}\ (t\in\mathbb R),
\]
where \(\mu_*>0\) is the least hyperbolic rate and \(\mathsf A\) is
nilpotent on \(E_c\) (its only eigenvalue there is zero). These bounds
are uniform for \(p_*\) near \(p_{\rm id}\). Fix \(0<\sigma_1<\mu_*\). For
\(Y=p_*+Z\) we have \(Z'=\mathsf AZ+\mathcal N(Z,\theta)\) with
\(\mathcal N(Z,\theta)=F_0(p_*+Z)-\mathsf AZ+S(p_*+Z,\theta)\), and
\(F_0(p_*)=S(p_*,\theta)=0\) give
\[
 \mathcal N(0,\theta)=0,\qquad
 \|D_Z\mathcal N(Z,\theta)\|\le C|Z|+\|S\|_{C^1}\quad(|Z|\le\varepsilon_0).
\]
For \(\xi\in E_s\) consider
\[
 \mathcal L Z(\theta)=e^{(\theta-\theta_0)\mathsf A}\xi
   +\int_{\theta_0}^{\theta}e^{(\theta-u)\mathsf A}P_s\mathcal N(Z(u),u)\,du
   -\int_\theta^\infty e^{(\theta-u)\mathsf A}(P_u+P_c)\mathcal N(Z(u),u)\,du
\]
in the norm \(\|Z\|_{\sigma_1}=\sup_{\theta\ge\theta_0}
e^{\sigma_1(\theta-\theta_0)}|Z(\theta)|\). The kernels are bounded by
\(Ce^{-\mu_*(\theta-u)}\) for \(u\le\theta\), and by
\(C(1+u-\theta)^m\) for \(u\ge\theta\). Moreover
\(|\mathcal N(Z(u),u)|\le(C\|Z\|_{\sigma_1}+\|S\|_{C^1})\|Z\|_{\sigma_1}
e^{-\sigma_1(u-\theta_0)}\), and
\(\int_\theta^\infty(1+u-\theta)^me^{-\sigma_1(u-\theta_0)}du\le
C_{m,\sigma_1}e^{-\sigma_1(\theta-\theta_0)}\). Hence, for
\(\varepsilon\) and \(\|S\|_{C^1}\) small and \(|\xi|\le\varepsilon/(2C)\),
\(\mathcal L\) is a contraction of the ball \(\|Z\|_{\sigma_1}\le\varepsilon\).
Its fixed point solves the equation, has \(P_sZ(\theta_0)=\xi\), and
decays at the rate \(\sigma_1\); since the center components are
integrated from infinity, the limit is exactly \(p_*\). The fixed point
depends continuously on \((p_*,S,\xi)\). For \(S=0\), \(p_*=p_{\rm id}\)
and \(\xi_{\rm id}=P_s(Y_{\rm id}(\theta_0)-p_{\rm id})\) it equals
\(Y_{\rm id}-p_{\rm id}\) by uniqueness, once \(\theta_0\) is so large
that this tail lies in the ball; fix such a \(\theta_0\). Given any
tolerance, the solution with \(\xi=P_s(Y_{\rm id}(\theta_0)-p_*)\) is then
as close to \(Y_{\rm id}\) at \(\theta_0\) as required, provided
\(\|S\|_{C^1}\) and \(|p_*-p_{\rm id}|\) are small enough. On
\([0,\theta_0]\) the difference \(\Delta\) of the two solutions satisfies
\(|\Delta(\theta)|\le e^{L_F(\theta_0-\theta)}
\bigl(|\Delta(\theta_0)|+\theta_0\|S\|_{C^0}\bigr)\) by Gronwall, where
\(L_F\) is the Lipschitz constant of \(F_0\) on \(\mathcal U\). This
gives the stated accuracy by first taking \(\lambda\) small and then
taking \(R\) large, using the iterated bound on \(\varepsilon_{\rm fit}\).
Finally
\(X=Y+\Phi(Y,\theta)\) also converges to \(p_*\) at the rate \(\sigma_1\),
since \(|\Phi(Y,\theta)|\le C\eta\lambda|Y-p_*|\).

On the invariant vertical lattice \(m=Rq\) use the norm \(Z^p\)
with weight \((1+|n|/R+|q|)^p\) and its weighted Wiener counterpart;
\(p=2\) or \(1<p=p_0<2\). The sum of reciprocal squared weights
is \(O_p(R)\). Hence two generic factors cost at most \(C_p\sqrt R\)
in the scaled algebra estimate, while any prescribed plane with
smooth envelope acts with its Wiener multiplier bound.
In detail, write \(\omega(n,q)=(1+|n|/R+|q|)^p\). Grouping \(n\) in
blocks of length \(R\) gives
\(\sum_{n,q}\omega(n,q)^{-2}\le CR\sum_{j,q}(1+|j|+|q|)^{-2p}\le C_pR\),
as \(p>1\). Subadditivity of \(|n|/R+|q|\) gives
\(\omega(k)\le2^{p-1}(\omega(l)+\omega(k-l))\), and therefore
\[
 \|fh\|_{Z^p}\le2^{p-1}\bigl(\|f\|_{Z^p}\|h\|_{\ell^1}
        +\|f\|_{\ell^1}\|h\|_{Z^p}\bigr),\qquad
 \|h\|_{\ell^1}\le\Big(\sum\omega^{-2}\Big)^{1/2}\|h\|_{Z^p}
        \le C_p\sqrt R\,\|h\|_{Z^p},
\]
by Young's inequality and Cauchy--Schwarz. If one factor is a
prescribed plane with smooth envelope, its weighted \(\ell^1\) norm
replaces the factor \(\sqrt R\).

We describe the linear normal operator at the true terminal state.
First remove the old axial pairing \(V^2\) from every doubled \(m,-m\)
block. Its first off-diagonal Sylvester generator has kernel
\(-C_{n-n'}/(E_n+E_{n'})\) with \(E_n=n^2+m^2\), with conjugate
index reflected. Its weighted Schur norm, and that of its old
derivatives, is \(O(R^{-2})\). Exponentiate in the group preserving
\(\operatorname{diag}(1,-1)\). In the first diagonal commutator
replace the divisor by \(2E_{n'}\); the reciprocal error is
\(O(R^{-3})\) times a fixed polynomial in \(|n-n'|\), absorbed by
the smooth coefficient. Indeed \(E_n+E_{n'}\ge2m^2\ge2R^2\) and
\(|n'|/E_{n'}\le1/(2|m|)\le1/(2R)\), so that
\[
 \left|\frac1{E_n+E_{n'}}-\frac1{2E_{n'}}\right|
 =\frac{|n'-n|\,|n'+n|}{2E_{n'}(E_n+E_{n'})}
 \le\frac{|n-n'|\,(2|n'|+|n-n'|)}{4R^2E_{n'}}
 \le\frac{|n-n'|}{4R^3}+\frac{|n-n'|^2}{4R^4}.
\]
The correction on the positive branch is
\(-|C|^2/(2E(D_x))\) to this accuracy. The remaining off-diagonal
error including the time connection is \(O(R^{-2})\); a second
Sylvester step of size \(O(R^{-4})\) leaves error \(O(R^{-4})\).

We verify these three assertions. Keep the conventions of
Section~\ref{sec:ins-propagator}, namely
\(C_{nn'}=\widehat{V_*^2}(n+n')\), \(X_{nn'}=-C_{nn'}/(E_n+E_{n'})\),
\[
 J_{\rm off}=\begin{pmatrix}0&X\\X^*&0\end{pmatrix},\qquad
 C_{\rm off}=\begin{pmatrix}0&C\\-C^*&0\end{pmatrix},
\]
and write \(\mathcal B_{\rm d}\) for the block-diagonal scalar potential,
so that the doubled generator of one block is
\(\mathcal K_0=D_{\rm free}+\mathcal B_{\rm d}+C_{\rm off}\). Since
\(E_n\) is even in \(n\), the reflection of the conjugate index does
not affect the divisors. As \([D_{\rm free},J_{\rm off}]=-C_{\rm off}\),
the Lie series for the conjugation and for the time connection give
\[
\begin{aligned}
 e^{-J_{\rm off}}\mathcal K_0e^{J_{\rm off}}
   -ie^{-J_{\rm off}}\partial_te^{J_{\rm off}}
 ={}&D_{\rm free}+\mathcal B_{\rm d}+\tfrac12[C_{\rm off},J_{\rm off}]
   +\bigl([\mathcal B_{\rm d},J_{\rm off}]-i\dot J_{\rm off}\bigr)\\
 &+\tfrac12[[\mathcal B_{\rm d},J_{\rm off}],J_{\rm off}]
   +\tfrac13[[C_{\rm off},J_{\rm off}],J_{\rm off}]
   -\tfrac i2[\dot J_{\rm off},J_{\rm off}]+\cdots .
\end{aligned}
\]
Indeed the terms of first order are
\([D_{\rm free},J_{\rm off}]+C_{\rm off}=0\), and those of second order
containing \(C_{\rm off}\) are
\([C_{\rm off},J_{\rm off}]+\tfrac12[[D_{\rm free},J_{\rm off}],J_{\rm off}]
=\tfrac12[C_{\rm off},J_{\rm off}]\); at third order
\(\tfrac12-\tfrac16=\tfrac13\). The omitted terms contain at least three
factors \(J_{\rm off}\) or \(\dot J_{\rm off}\); the first of them are
\(\tfrac16[[[\mathcal B_{\rm d},J_{\rm off}],J_{\rm off}],J_{\rm off}]\),
\(\tfrac18[[[C_{\rm off},J_{\rm off}],J_{\rm off}],J_{\rm off}]\) and
\(-\tfrac i6[[\dot J_{\rm off},J_{\rm off}],J_{\rm off}]\). Each factor
\(J_{\rm off}\) or \(\dot J_{\rm off}\) costs \(O(R^{-2})\) in weighted
Schur size, so the second line is \(O(R^{-4})\) and the omitted terms
are \(O(R^{-6})\). The term
\(\tfrac12[C_{\rm off},J_{\rm off}]\) is block diagonal, and
\[
\begin{aligned}
 \tfrac12[C_{\rm off},J_{\rm off}]
 &=\tfrac12\operatorname{diag}\bigl(CX^*+XC^*,\,-C^*X-X^*C\bigr),\\
 \tfrac12(CX^*+XC^*)_{nn'}
 &=-\frac12\sum_lC_{nl}\overline{C_{n'l}}
   \Bigl(\frac1{E_{n'}+E_l}+\frac1{E_n+E_l}\Bigr).
\end{aligned}
\]
For integers \(a,a'\) put \(d=\bigl||a|-|a'|\bigr|\). The computation
displayed above gives
\(|(E_a+E_{a'})^{-1}-(2E_a)^{-1}|\le d/(4R^3)+d^2/(4R^4)\), and likewise
\(|(2E_a)^{-1}-(2E_{a'})^{-1}|\le d/(2R^3)+d^2/(2R^4)\). In the last
sum \(\bigl||n|-|l|\bigr|\le|n+l|\), \(\bigl||n'|-|l|\bigr|\le|n'+l|\)
and \(\bigl||n|-|n'|\bigr|\le|n+l|+|n'+l|\),
and the two coefficients decay rapidly in \(|n+l|\) and \(|n'+l|\).
Both reciprocals may therefore be replaced by \((2E_{n'})^{-1}\) with an
error \(O(R^{-3})\) in weighted Schur size. Since
\(\sum_lC_{nl}\overline{C_{n'l}}=\widehat{|V_*|^4}(n-n')\) by Parseval,
the upper block is multiplication by \(-|V_*|^4\) composed with
\((2E(D_x))^{-1}\), up to \(O(R^{-3})\), and the lower block is its
negative conjugate. This is the term \(-n^2|C|^2/(2E)\) of \(L_2\) in
the proof of Lemma~\ref{lem:empty-frequency}. The remaining term
\(Y_{\rm off}\coloneqq[\mathcal B_{\rm d},J_{\rm off}]-i\dot J_{\rm off}\)
is block off-diagonal of size \(O(R^{-2})\), together with its
derivatives in the old parameters. Let \(J_2\) be the block
off-diagonal matrix with \([D_{\rm free},J_2]=-Y_{\rm off}\); its
entries are those of \(Y_{\rm off}\) divided by \(\pm(E_n+E_{n'})\),
of modulus at least \(2R^2\), so \(J_2\) and its parameter
derivatives are \(O(R^{-4})\). Conjugation by \(e^{J_2}\) removes
\(Y_{\rm off}\) and adds only
\([\mathcal B_{\rm d}+\tfrac12[C_{\rm off},J_{\rm off}],J_2]-i\dot J_2\)
and terms of higher order, all of size \(O(R^{-4})\). After division
by the slow rate \(r\lambda=\eta\lambda R^{-2}\) the total
off-diagonal error is \(O(R^{-2}/(\eta\lambda))\).

For \(|n|\ge hR\) use the spatial columns above in their corresponding
Hermitian copy, with exact diagonal to error \(O(R^{-3})\). Projection
to that copy after the two changes preserves the residual.
Indeed, let \((e^{inx}f,e^{inx}g)\) be the doubled column of
Lemma~\ref{lem:empty-frequency} with \(|n|\ge hR\), so that
\(g=-\bar Cf/(2E)+O(R^{-3})\) in the notation of that proof. In the
reflected labelling of the second copy, \(X^*\) applied to the first
component is \(-\bar Cf/(2E)+O(R^{-3})\), by the reciprocal estimate
just used. Hence the second component \(g-X^*f+O(R^{-4})\) of the image
of the column under \(e^{-J_{\rm off}}\) is \(O(R^{-3})\), while its
first component differs from \(e^{inx}f\) by \(O(R^{-4})\); the map
\(e^{-J_2}\) changes both by \(O(R^{-4})\). The conjugated generator is
block diagonal up to \(O(R^{-4})\). The first component is therefore a
quasimode of the upper block with the same energy and with residual
\(O(R^{-3})\), and its envelope is still \(1+O(R^{-1})\).
For the
high column matrix \(T\), \(G=T^*T=1+O(R^{-1})\). The quasimode
identity and its adjoint bound the covariant commutator
\([D,G]-i\dot G\) with the high diagonal \(D\) by \(O(R^{-3})\). Leibniz and the power
series transfer this estimate to \(G^{-1/2}\). Complete \(TG^{-1/2}\)
unitarily by the polar part of the projection of the low subspace
onto its orthogonal complement. Spatial off-diagonal decay and the
\(O(1/|n|)\) high envelopes bound the commutator with \(D_x^2\) by
a bounded operator. The low generator is therefore free kinetic energy
plus a bounded Hermitian operator. It is retained at the actual trial
actions and angles in the nonlinear problem.

The terminal high axial waves also give terms of size \(r\), both
Hermitian and anomalous, in the normal operator. Their leading physical
shapes are the planes at \(K_i=1+RP_i\). This follows from
\eqref{eq:chart-tube}--\eqref{eq:chart-head} at \(p=2,p_0\):
the new axial input is \(O(r)\) in \(Z^2\), \(O(rR)\) in \(h^1\);
one-smoothing for both chart differentials and symplectic adjointness
give \(O(R^{-1})\) on a high normalized plane. Changing the input
in the differential changes its scaled operator norm by
\(O(r+rR^{2-p})\); this uses the unweighted analytic \(L^2\) bound
and the difference bound from \(h^1\) to \(H^2\), scaled on input
and output. Thus the error on a pure physical wave is
\(O(r[R^{-1}+r(1+R^{2-p})])=o(rR^{-1/2})\).
In detail, let \(z_{\rm hot}\) be the hot Birkhoff coordinates, of
size \(r|a_i|\) at the sites \(K_i\). The weight
\((1+|n|/R)^p\) is bounded at \(n=K_i\) and \(\bk{K_i}\le CR\), which
gives the two sizes \(O(r)\) in \(Z^p\) and \(O(rR)\) in \(h^1\). The
axial field is
\[
 \Psi(z_F+z_{\rm hot})=\Psi(z_F)+\mathcal F^{-1}z_{\rm hot}
 +\bigl(D\Psi(z_F)-\mathcal F^{-1}\bigr)z_{\rm hot}
 +\int_0^1\bigl(D\Psi(z_F+\vartheta z_{\rm hot})-D\Psi(z_F)\bigr)
 z_{\rm hot}\,d\vartheta .
\]
The second term is the sum of the planes at \(K_i\). By
\eqref{eq:chart-head} without the residue term the third term has
\(Z^p\) norm \(O(R^{-1})\|z_{\rm hot}\|_{Z^p}=O(rR^{-1})\), and by
\eqref{eq:chart-tube} the last term has \(Z^p\) norm
\(O(r(1+R^{2-p}))\|z_{\rm hot}\|_{Z^p}=O(r^2(1+R^{2-p}))\); this is the
displayed error. With \(r=\eta/R^2\) and \(p\ge p_0>1\),
\(r[R^{-1}+r(1+R^{2-p})]\le r(R^{-1}+\eta R^{-2}+\eta R^{-p})
=O(rR^{-1})\). The terms of size \(r\) in the normal operator are
\(2(\overline{V}U_{\rm hot}+V\overline{U_{\rm hot}})w\) and
\(2VU_{\rm hot}\bar w\), where
\(U_{\rm hot}=\Psi(z_F+z_{\rm hot})-\Psi(z_F)\) and \(V=\Psi(z_F)\). To
use the replacement of \(U_{\rm hot}\) by its planes in these
multipliers we need the error in the weighted Wiener norm
\(\sum_k(1+|k|/R)^p|\widehat f(k)|\), which acts as a multiplier on
\(Z^p\), because \((1+|k+k'|/R)^p\le(1+|k|/R)^p(1+|k'|/R)^p\) and
Young's inequality; it follows from coefficientwise bounds. By
translation covariance (a translation of \(x\) by \(a\) multiplies each
Birkhoff coordinate \(z_n\) by \(e^{ina}\)), the parts of
\((D\Psi(z_F)-\mathcal F^{-1})z_{\rm hot}\) that are linear in
\(z_{K_i}\) and in \(\bar z_{K_i}\), expanded in the harmonics
\(e^{ij\cdot\theta}\), \(j\in\Z^F\), of the old angles, have the
\(j\)-th harmonic at the single Fourier index \(k=K_i+j\cdot n_F\),
respectively \(k=-K_i+j\cdot n_F\). The map
\(D(\Psi-\mathcal F^{-1})(z_F)\) is bounded from \(h^1\) to \(H^2\)
\cite[Theorem~1.1 and Corollary~1.1]{KST} and analytic on a complex
strip of old angles, so Cauchy's estimate in the angles bounds the
coefficient at \(k=\pm K_i+j\cdot n_F\) by
\(Ce^{-c|j|_1}\bk{K_i}\bk k^{-2}|z_{K_i}|\), using
\(\|e_{K_i}\|_{h^1}=\bk{K_i}\) and \(\|f\|_{H^2}\ge\bk k^2|\widehat f(k)|\).
For \(|j\cdot n_F|\le|K_i|/2\) one has \(\bk k\ge\bk{K_i}/2\), and the
coefficient is at most
\(4Ce^{-c|j|_1}\bk{K_i}^{-1}|z_{K_i}|\le Cre^{-c|j|_1}R^{-1}\), since
\(|z_{K_i}|\le Cr\) and \(\bk{K_i}\gtrsim R\); for the other \(j\) one has
\(|j|_1\ge|K_i|/(2\max_{n\in F}|n|)\ge c'R\), and the coefficient is at
most \(Cre^{-c|j|_1}\bk{K_i}\le Cre^{-c|j|_1/2}e^{-cc'R/2}R\). In both
cases \(|k|\le|K_i|+|j|_1\max_{n\in F}|n|\le C(1+|j|_1)R\), so the weight
is at most \(C(1+|j|_1)^p\), and summing over \(j\) and the finitely many
\(i\) shows that the third term has weighted Wiener norm \(O(rR^{-1})\).
The last term is the part of degree at least two in the hot variables
\(z_{K_i},\bar z_{K_i}\) of the expansion of \(\Psi(z_F+z_{\rm hot})\),
which converges for \(|z_{K_i}|\le\rho_0\) with a fixed \(\rho_0>0\),
since \(\Psi\) is real analytic from \(h^0\) to \(L^2\) and the old torus
is compact. By translation covariance the coefficient of a monomial of
degree \(m\) in the hot variables and of the harmonic \(j\) in the old
angles is a single Fourier mode, at an index \(k\) with
\(|k|\le m\max_i|K_i|+|j|_1\max_{n\in F}|n|\le C(m+|j|_1)R\), where the
weight is at most \(C(1+m)^p(1+|j|_1)^p\); Cauchy's estimate on the
polydisc of radius \(\rho_0\) and on the strip of old angles bounds its
modulus by \(C\rho_0^{-m}e^{-c|j|_1}\). There are at most
\((m+1)^{2n_{\rm h}}\) monomials of degree \(m\), \(n_{\rm h}\) being the
number of hot sites, and \(|z_{K_i}|\le Cr\), so the weighted Wiener norm
of the last term is at most
\[
 C\sum_{m\ge2}(m+1)^{2n_{\rm h}+p}\Bigl(\frac{Cr}{\rho_0}\Bigr)^m
 \sum_{j\in\Z^F}(1+|j|_1)^pe^{-c|j|_1}=O(r^2).
\]
Thus the multipliers change by
\(O(rR^{-1}+r^2)=O(rR^{-1})=o(rR^{-1/2})\).

At angular degree up to \(C_1\log R\) write \(b=l\cdot n\),
\(\nu=l\cdot n^2+\tau l\cdot c^0\); the discarded tails can cost
any fixed negative power by taking \(C_1\) large at fixed clock.
Every nonzero bare divisor \(\mathbb Z-\nu\) is bounded away from
zero by an inverse power of \(1+|l|\), by the fixed rational base
and the Diophantine choice of \(\tau\) (when the rational multiple
of \(\tau\) is zero the nonzero integer costs at least one).
An anomalous pair has \(n+n'=K_i+b\) and divisor
\[
 n^2+n'^2+2R^2q^2-K_i^2-\nu.
\]
If either leg is low-horizontal this is \(\gtrsim R^2\) in modulus:
choose \(h>0\) small, since the limiting polynomial at \(x=n/R\)
is \(2(x^2-P_i x+q^2)\), \(q\ne0\). Take \(h\) also below half
the least nonzero absolute horizontal velocity of the finitely many
heads and pure sources. The Hermitian divisor is
\(2(K_i-b)(n-b)+\nu-b^2\), of modulus \(\gtrsim R\) for \(n\ne b\).
To see the two formulas, note that the pairing term couples the
coefficient at \((n,Rq)\) to the conjugate coefficient at
\((n',-Rq)\) through the product of the terminal plane at \(K_i\)
with an old harmonic of spatial index \(b\) and time frequency
\(\nu\). The output index is \(n=K_i+b-n'\) and the mismatch is
\((n^2+R^2q^2)+(n'^2+R^2q^2)-K_i^2-\nu\). With \(x=n/R\) and
\(n'/R=P_i-x+O(R^{-1}\log R)\) its quotient by \(R^2\) tends to
\(x^2+(P_i-x)^2+2q^2-P_i^2=2(x^2-P_ix+q^2)\), which is at least
\(2(1-|P_i|h)\) when \(|x|\le h\) or \(|P_i-x|\le h\), because the
polynomial is symmetric under \(x\mapsto P_i-x\) and \(q^2\ge1\). The density
term maps the input \(n\) to the output \(n'=n+K_i-b\) at the same
\(q\), with mismatch
\[
 (n+K_i-b)^2-n^2-K_i^2+\nu=2(K_i-b)(n-b)+\nu-b^2 .
\]
Here \(|b|+|\nu|^{1/2}\le C\log R\) and \(|K_i-b|\ge cR\), which gives
the lower bound for \(n\ne b\).
Keep every star from \(b\) to \(K_i\), and its adjoint, including
off-key harmonics there. Divide the other nonresonant terms against
the free diagonal and the time derivative together. The residual
from the actual slow phase costs \(O(R^{-2}\log R)\) times the
divided kernel. The commutator with the bounded low body acts only
when a leg is low, where the preceding large divisors apply.
Thus the remaining Schur error is
\(O(R^{-3}+r^2(1+\log R)^C+(r/R)(1+\log R)^C)=o(r\lambda)\).
Each harmonic has one shift or reflection per row; the weight ratios
are bounded by fixed powers of its spatial displacement in the scaled
coordinates. Constants here may depend on the fixed clock.

The three terms of this bound arise as follows. Let \(\mathcal B\) be
the exponent of the change of variables that removes the nonresonant
terms of order \(r\), as in Section~\ref{sec:ins-propagator}; its
entries are the discarded entries divided by their bare divisors. An
opposite-copy entry with both legs high has divisor at least
\(c(1+|l|)^{-C}\ge c(1+\log R)^{-C}\) at the retained angular degrees,
an opposite-copy entry with a low leg has divisor \(\gtrsim R^2\), and
a same-copy entry with \(n\ne b\) has divisor \(\gtrsim R\). Summing the
partial shifts and reflections over the harmonics, \(\mathcal B\) has
Schur size \(O(r(1+\log R)^C)\), and its part touching the low
subspace has size \(O(r/R)\). By the identity for
\(\mathcal M=e^{\mathcal B}\) in Section~\ref{sec:ins-propagator}, the
new generator contains, besides the retained terms,
\begin{itemize}
\item the difference between the actual and the bare divisor times
\(\mathcal B\). The actual divisor contains the root shift
\(\chi_{K_i}(I^\infty)=O(R^{-2})\), the actual shift
\(l\cdot(\Omega^\infty-\Omega^0)=O(R^{-2}\log R)\) of the old harmonic
and the two high diagonals \(A/n^2-B/(n^2+m^2)=O(R^{-2})\). Since
\(R^{-2}=r/\eta\), this term is \(O(r^2(1+\log R)^C)\);
\item the commutators of \(\mathcal B\) with the retained stars and
pairs and all terms of second order in the generator of order \(r\),
of size \(O(r^2(1+\log R)^C)\);
\item the commutator of \(\mathcal B\) with the bounded low body, of
size \(O((r/R)(1+\log R)^C)\), since only the entries of \(\mathcal B\)
with a low leg contribute;
\item the residual \(O(R^{-3})\) of the high columns and the error
\(O(R^{-4})\) of the two Sylvester steps; conjugating the terms of
order \(r\) by the column frame, which differs from the identity by
\(O(R^{-1})\), costs \(O(r/R)\); and the angular tail costs
\(O(R^{-N_0})\).
\end{itemize}
After division by \(r\lambda=\eta\lambda R^{-2}\) the sum is
\(O\bigl((\eta\lambda R)^{-1}+\lambda^{-1}(r+R^{-1})(1+\log R)^C\bigr)\),
which tends to zero as \(R\to\infty\) at fixed \(\eta\) and \(\lambda\).

At an exact bare anomalous resonance, \(q\) is bounded and its
limiting \(x\) satisfies \(x(P_i-x)=q^2\). Both legs are high.
The row determines \(P_i\) in the limit; at that root changing the
integer \(b\) changes the partner's square by order \(R\), for the
logarithmic range above. The remaining pairs are disjoint, and do
not meet a Hermitian star, by the nonsquare condition of
Proposition~\ref{prop:tree}. Flatten the horizontal weight for
\(|n|\le C_{\rm tree} R\) as in \eqref{eq:flat-weight-ins}.
The low generator and all stars together are Hermitian on each
constant-\(q\) space in this equivalent norm.

In a surviving anomalous pair there may be several order-five
groups even at an integer key different from \((1,1)\). Its
limiting second-order defect is
\[
 E_*=\Xi((x,q))+\Xi((P_i-x,-q))-\Xi((P_i,0)).
\]
In the interaction frame removing the full nominal high diagonal,
retain the group with \(d_j\) nearest \(E_*\), if present. When there are at least two rows, their defects \(d_j\)
have distinct limits as \(\lambda\downarrow0\) at fixed passive data.
Thus their pairwise distances are bounded below by some \(\Delta>0\)
for all sufficiently small \(\lambda\). Every group except the chosen
nearest one then has \(|d_j-E_*|\ge\Delta/2\). The bounded integers
\(q\) and the roots of \(x(P_i-x)=q^2\) give only finitely many limiting
pairs. For each fixed \(\lambda\), taking \(R\ge R_0(\lambda)\) large
enough ensures that the absolute physical detuning of each discarded
group is at least \(\Delta/(4R^2)\). With a single row there are no
other groups. Divide the oscillations of the discarded groups in this
frame. The
normalized change is \(O(\eta\lambda)\), since the generator
amplitudes in this frame are bounded. In physical time the trial change in the second-order correction
\(A/n^2-B/(n^2+m^2)\) is \(O(R^{-2})(r|J_{\rm old}|+r^2)\),
where \(J_{\rm old}=(I_{\rm old}-I_{\rm old}^\infty)/r\).
In the fixed terminal mass gauge, the old-mass scalar on the positive
branch changes by \(2\sum_{j\in F}(I_j-I_j^\infty)
=2r\sum_{j\in F}(J_{\rm old})_j\), with the opposite sign on the
conjugate branch. This term is included in the old-mass estimate below.
After this averaging, the unitary half-phase rotation for the retained group
gives the two exponents
\(\pm\operatorname{Re}\sqrt{|g|^2-\delta^2/4}\), with its actual
normalized detuning \(\delta\) and normalized coupling \(g\). In particular no small residual
detuning has been suppressed. Indeed, after the half-phase rotation
the two legs \(u,v\) of the pair satisfy, in slow time,
\[
 iu'=\tfrac\delta2u+g\bar v,\qquad iv'=\tfrac\delta2v+g\bar u,
 \qquad\text{that is}\qquad
 \begin{pmatrix}u\\ \bar v\end{pmatrix}'
 =\begin{pmatrix}-i\delta/2&-ig\\ i\bar g&i\delta/2\end{pmatrix}
  \begin{pmatrix}u\\ \bar v\end{pmatrix}.
\]
The matrix has trace zero and determinant \(\delta^2/4-|g|^2\), so
its eigenvalues are \(\pm\sqrt{|g|^2-\delta^2/4}\). They are real
when \(|\delta|<2|g|\) and purely imaginary when \(|\delta|>2|g|\).

\begin{stagelemma}[A common Green right inverse]\label{lem:common-endpoint-green}
Fix a finite endpoint preparation as above, \(1<p_0<2\), and
\(0<a<b\) below the finite orbit's landing rate in slow time
\(\theta=r\lambda t\). For every sufficiently small admissible clock
strength \(\lambda\) and every sufficiently large admissible
\(R\ge R_0(\lambda)\), there is a closed
interval \(J\subset(a,b)\) of length bounded below in terms of the
fixed preparation and \([a,b]\). On the invariant vertical lattice
\(m\in R\mathbb Z\), the exact linearization at the terminal axial
torus has a common Green right inverse for \(p\in\{p_0,2\}\) and
\(\sigma\in J\). More precisely, on a translated late half-line
\([t_0,\infty)\) it maps each continuous forcing of finite norm to a mild solution
with a bound in the same norm
\[
 \|f\|_{p,\sigma,t_0}
 =\sup_{t\ge t_0}e^{\sigma r\lambda(t-t_0)}\|f(t)\|_{H^p}.
\]
The fast stable trace at \(t_0\) is zero in the reduced coordinates.
The inverses agree on common data. At the fixed clock and dilation
their physical-norm bounds are uniform in \(t_0\).
\end{stagelemma}

\begin{proof}
A retained decay rate in \([a,b]\) requires \(|g|\ge a\).
The preceding absolute sum over integer keys is bounded after division
by \(\lambda\), uniformly as \(\lambda\downarrow0\). At fixed pump
and integer key, bare resonance prescribes both \(n+n'\) and
\(n^2+n'^2+2R^2q^2\). The latter bounds \(|q|\) in terms of the
fixed tree since \(K_i=1+RP_i\) and \(|\nu|=O(\log R)\) at the
retained angular cutoff. For each \(q\) the two equations leave at
most two choices of \(n\). Thus the actual count is bounded by an
\(M\) depending only on this preparation and \(a\), for every
sufficiently large admissible \(R\). After adjoining \(a,b\) the
longest gap has length at least \((b-a)/(M+1)\). Take its middle
third for \(J\); its location may vary with the parameters.

Set \(k^2=|g|^2-\delta^2/4\). For bounded \(\delta\) the pair
propagator uses \(\cosh(k\Delta\theta)\) and
\(\sinh(k\Delta\theta)/k\), with value \(\Delta\theta\) for the
second factor at \(k=0\). When the cut separates the pair, its
spectral projections are bounded by the clearance. For large
\(|\delta|\) the oscillatory propagator is uniformly bounded. The
Hermitian complement is isometric in the flattened norm. Integrate
forward from zero on the fast stable subspace and backward from
infinity on its complement. The clearance makes both weighted
integrals bounded throughout \(J\); their integration directions
are the same at every weight. The transverse remainder
\(O(\lambda^{1/5}+\eta\lambda)+o_R(1)\) is absorbed by the same
Neumann series. The separate axial block in Birkhoff coordinates is
a rotation with an occupied action-angle shear. Backward integration
here is bounded since a linear polynomial times a decaying exponential
is integrable. These formulas commute with inclusion of the stronger
weighted space into the weaker one, first for the reduced inverse and
then term by term in the Neumann series. For data with crossed regularities and weights, the two solutions
agree by uniqueness at the lower order and smaller weight. The scaled norms are equivalent
to the physical norms at fixed \(R\), and the estimates depend only on
elapsed time. This proves the asserted compatibility and translation bound.

\smallskip
\noindent\emph{Details.} We expand the steps of this argument. Write
\(\sigma_-=a\) and \(\sigma_+=b\) for the endpoints in the statement,
because \(b\) also denotes the spatial index of a key. Time is
\(\theta=r\lambda t\), and all generators are divided by \(r\lambda\).

\emph{(i) Counting the rates.} A retained anomalous pair is labelled
by a pump \(K_i\), an integer key \((b,\nu)\), a vertical index
\(q\ne0\) and a horizontal index \(n\) with
\[
 n+n'=K_i+b,\qquad n^2+n'^2+2R^2q^2=K_i^2+\nu .
\]
Write \(a_{i,\mathrm{end}}\) for the complex amplitude of the pump
at the rounded endpoint in the normalized carrier coordinates. If \(v\)
is the coefficient of the retained order-five group at this key, the
normalized coupling is \(g=2a_{i,\mathrm{end}}v/\lambda\).
By the last estimate in the preceding subsection, the sum of \(|v|\)
over all order-five groups at all integer keys is at most
\(C_0\lambda\), for small \(\lambda\) at fixed passive data. Since
\(|a_{i,\mathrm{end}}|=|a_i|+o_R(1)\) for the prescribed amplitude \(a_i\),
we have \(\sum|g|\le C_1\coloneqq2C_0\sum_i(|a_i|+1)\), the sum being over all
pumps, integer keys and groups. In particular \(|g|\le C_1\) for every
pair, and at most \(C_1/\sigma_-\) of these triples have
\(|g|\ge\sigma_-\). For a given triple, the second equation and
\(|\nu|\le C\log R\) give \(2q^2R^2\le K_i^2+C\log R\), so
\(|q|\le|P_i|\) for large \(R\); for each \(q\), eliminating \(n'\)
leaves a quadratic equation for \(n\). Therefore at most
\[
 M\coloneqq4C_1\sigma_-^{-1}\max_i|P_i|
\]
retained pairs have \(|g|\ge\sigma_-\), for every large admissible
\(R\). The two exponents of a pair are \(\pm k_+\), where
\(k_+=\operatorname{Re}\sqrt{|g|^2-\delta^2/4}\le|g|\). Thus at most
\(M\) pairs have \(k_+\in[\sigma_-,\sigma_+]\). These rates and the
two endpoints divide \([\sigma_-,\sigma_+]\) into at most \(M+1\)
intervals. Let \((\alpha,\beta)\) be a longest one, and let \(J\) be
its closed middle third. No rate lies in \((\alpha,\beta)\), so for
every \(\sigma\in J\) and every retained pair
\begin{equation}\label{eq:h2-clearance}
 |\sigma-k_+|\ge\kappa_0\coloneqq\frac{\sigma_+-\sigma_-}{3(M+1)},
 \qquad \sigma\ge\sigma_- .
\end{equation}

\emph{(ii) One pair.} By the display preceding the lemma, the
generator of a pair in slow time is the constant matrix
\[
 \mathsf M=\begin{pmatrix}-i\delta/2&-ig\\ i\bar g&i\delta/2\end{pmatrix},
 \qquad \mathsf M^2=k^2,\qquad k^2=|g|^2-\delta^2/4 .
\]
Hence
\(e^{\Delta\mathsf M}=\cosh(k\Delta)+k^{-1}\sinh(k\Delta)\,\mathsf M\),
where the second coefficient is \(\Delta\) at \(k=0\). If
\(|\delta|\le4C_1\), then \(\|\mathsf M\|\le|\delta|/2+|g|\le3C_1\) and
\(|k^{-1}\sinh(k\Delta)|\le|\Delta|e^{k_+|\Delta|}\), so that
\(\|e^{\Delta\mathsf M}\|\le(1+3C_1|\Delta|)e^{k_+|\Delta|}\). If
\(|\delta|>4C_1\ge4|g|\), then \(k=i\kappa\) with
\(\kappa^2=\delta^2/4-|g|^2\ge3\delta^2/16\) and
\(\|\mathsf M\|\le3|\delta|/4\), so that
\(\|e^{\Delta\mathsf M}\|\le1+\|\mathsf M\|/\kappa\le1+\sqrt3\) and
\(k_+=0\). In both cases
\begin{equation}\label{eq:h2-pair-propagator}
 \|e^{\Delta\mathsf M}\|\le C(1+|\Delta|)e^{k_+|\Delta|},
 \qquad C=C(C_1).
\end{equation}
When \(k_+>\sigma\), the number \(k=k_+\ge\sigma_-\) is real and
\(|\delta|\le2|g|\). The spectral projections
\(\mathsf P_\pm=\tfrac12(1\pm\mathsf M/k)\) then satisfy
\(\|\mathsf P_\pm\|\le\tfrac12(1+2C_1/\sigma_-)\) and
\(e^{\Delta\mathsf M}\mathsf P_\pm=e^{\pm k\Delta}\mathsf P_\pm\).
For leg weights \(w_1,w_2\), conjugation by
\(\operatorname{diag}(w_1,w_2)\) multiplies the operator norm bound
by at most \(\max(w_1/w_2,w_2/w_1)\). This factor is uniformly bounded
for \(p\in\{p_0,2\}\), as in Section~\ref{sec:ins-propagator}, and is
absorbed in the weighted Green estimate below.

\emph{(iii) The principal Green operator.} Let \(\Pi_f\) be the sum
of the projections \(\mathsf P_-\) over the pairs with \(k_+>\sigma\),
transported to the reduced coordinates by the unitary half-phase
rotations, and let \(\Pi_c=1-\Pi_f\). By \eqref{eq:h2-clearance} the
same pairs occur for every \(\sigma\in J\). Let \(\mathcal U_0\) be
the propagator of the pairs and of the Hermitian complement, put
\(\theta_0=r\lambda t_0\) and, as in Section~\ref{sec:ins-propagator},
\[
 (\mathcal G_0f)(\theta)
 =\int_{\theta_0}^{\theta}\mathcal U_0(\theta,u)\Pi_f(u)f(u)\,du
  -\int_\theta^\infty\mathcal U_0(\theta,u)\Pi_c(u)f(u)\,du .
\]
Multiply by \(e^{\sigma(\theta-\theta_0)}\). On the range of \(\Pi_f\)
the weighted kernel is bounded by
\(Ce^{-(k_+-\sigma)(\theta-u)}\le Ce^{-\kappa_0(\theta-u)}\) for
\(u\le\theta\). On the range of \(\Pi_c\) and for \(u\ge\theta\) it is
bounded by \(Ce^{-(k_++\sigma)(u-\theta)}\) on the growing direction
of a pair with \(k_+>\sigma\); by
\(C(1+u-\theta)e^{-(\sigma-k_+)(u-\theta)}
\le C(1+u-\theta)e^{-\kappa_0(u-\theta)}\) on a pair with
\(k_+<\sigma\), using \eqref{eq:h2-pair-propagator}; and by
\(e^{-\sigma_-(u-\theta)}\) on the Hermitian complement, which is
isometric in the flattened norm. Taking the direct sum over the
blocks before integrating in time gives
\[
 \|\mathcal G_0f\|_{p,\sigma}\le C_J\|f\|_{p,\sigma},\qquad
 C_J=C\bigl(\kappa_0^{-1}+\kappa_0^{-2}+\sigma_-^{-1}\bigr),
\]
for \(p\in\{p_0,2\}\), every \(\sigma\in J\), every \(t_0\) and every
large admissible \(R\). Here \(\|\cdot\|_{p,\sigma}\) denotes the
weighted supremum of the flattened \(Z^p\) norm in the reduced
coordinates. A homogeneous solution of the principal system with
finite norm and zero fast stable trace at \(\theta_0\) vanishes: its
\(\Pi_f\) component is determined by its trace, and no nonzero
solution in the range of \(\Pi_c\) decays at the rate \(\sigma\).

\emph{(iv) The remainder.} In the reduced coordinates the exact
transverse generator is the principal one plus a bounded remainder
\(\mathcal R(\theta)\) with
\[
 \|\mathcal R\|_{Z^p\to Z^p}\le\epsilon_*
 \coloneqq C(\lambda^{1/5}+\eta\lambda)+o_R(1),\qquad p\in\{p_0,2\}.
\]
Its three parts are the harmonics of weak order at least six, of
total size \(O(r\lambda^{6/5})\) before normalization; the averaging
of the nonretained groups on the pairs, of normalized size
\(O(\eta\lambda)\); and the Schur error \(o(r\lambda)\) estimated
above. A mild solution of the exact equation with forcing \(f\),
finite norm and zero fast stable trace at \(\theta_0\) is the same as
a solution of \(h=\mathcal G_0(f+\mathcal Rh)\), by the uniqueness
statement in (iii). The rate count in (i) and the separation of the
discarded groups give bounds independent of small \(\lambda\) for
\(C_J\) and for the coefficient of \(\lambda^{1/5}+\eta\lambda\) in
\(\epsilon_*\). Choose \(\lambda\) and then \(R_0(\lambda)\) so
that \(C_J\epsilon_*\le1/2\). Then
\[
 \mathcal G_\perp=\sum_{j\ge0}(\mathcal G_0\mathcal R)^j\mathcal G_0,
 \qquad \|\mathcal G_\perp f\|_{p,\sigma}\le2C_J\|f\|_{p,\sigma}.
\]
This bound is uniform in \(R\ge R_0(\lambda)\). It is the bound used
for the transverse row in the next subsection.

\emph{(v) The axial block.} Let \(V'\) be the terminal axial solution,
\(F'\) its finite Birkhoff support, and
\(\omega_j'=\omega_j(I^\infty)+c_R\) its gauged frequencies. For
\(1\le p\le2\) the differential \(D\Phi(V'(t))\) is bounded from
\(H^p\) to \(h^p\) together with its inverse, uniformly in \(t\), by
one-smoothing and compactness of the torus. It conjugates the axial
linearization to the linearization of
\(\dot z_j=-i(\omega_j(I)+c_R)z_j\) at the torus. For \(j\notin F'\)
the latter is \(\dot y_j=-i\omega_j'y_j+f_j\), since the variation of
\(\omega_j\) multiplies \(z_j=0\). For \(j\in F'\) use the tangent
action-angle variables \((\delta I,\delta\theta)\), which are
equivalent to \(y_j\) because \(I_j^\infty>0\). The variation of an
empty action vanishes to first order, so
\[
 \partial_t\delta I=f_I,\qquad
 \partial_t\delta\theta=\Omega'\delta I+f_\theta,\qquad
 \Omega'=\bigl(\partial_{I_k}\omega_j(I^\infty)\bigr)_{j,k\in F'} .
\]
Integrate the three equations backwards from infinity:
\[
\begin{aligned}
 \delta I(t)&=-\int_t^\infty f_I(u)\,du,\qquad
 y_j(t)=-\int_t^\infty e^{i\omega_j'(u-t)}f_j(u)\,du,\\
 \delta\theta(t)&=-\int_t^\infty f_\theta(u)\,du
               +\int_t^\infty(u-t)\,\Omega'f_I(u)\,du .
\end{aligned}
\]
Since
\(\int_t^\infty e^{-\sigma r\lambda(u-t_0)}du
=(\sigma r\lambda)^{-1}e^{-\sigma r\lambda(t-t_0)}\) and
\(\int_t^\infty(u-t)e^{-\sigma r\lambda(u-t_0)}du
=(\sigma r\lambda)^{-2}e^{-\sigma r\lambda(t-t_0)}\), the axial
solution is bounded in \(\|\cdot\|_{p,\sigma,t_0}\) by
\(C_{V'}\bigl[(\sigma r\lambda)^{-1}+\|\Omega'\|(\sigma r\lambda)^{-2}\bigr]
\|f\|_{p,\sigma,t_0}\). This constant depends on \(\lambda\) and \(R\)
but not on \(t_0\), and the three formulas involve neither \(p\) nor
\(\sigma\).

\emph{(vi) Compatibility.} Let \(p_1\le p_2\) in \(\{p_0,2\}\) and
\(\sigma_1\le\sigma_2\) in \(J\), and let \(\iota\) be the inclusion
of the space with indices \((p_2,\sigma_2)\) into that with indices
\((p_1,\sigma_1)\). The operators \(\mathcal G_0\), \(\mathcal R\) and
the axial integrals are given by the same formulas on both spaces,
with the same projections by (iii). Hence \(\iota\) intertwines
\(\mathcal G_0\), each term of the Neumann series, and the axial
inverse; the series converges in both norms. If a forcing has finite
norm for two arbitrary pairs of indices, both of its images lie in the
space with the smaller regularity index and the smaller weight, and
both coincide there with the image under the inverse of that space.
Finally the changes of variables of this subsection are bounded with
bounded inverses on \(Z^p\) for both \(p\), uniformly in time, and
\(Z^p\) is equivalent to \(H^p\) on the lattice at fixed \(R\). All
kernel bounds depend only on \(\theta-u\). This gives the statement
in the physical norms, with bounds independent of \(t_0\).
\end{proof}

\subsection{Insertion and passage to a single initial state}
For completeness we record the endpoint error sizes in the original
terminal argument. Write the axial projection of the actual field as
\(U_0\), its transverse part as \(w\), and \(z=\Phi(U_0)\).
Equations \eqref{eq:axis-insertion}--\eqref{eq:old-mod-insertion}
are exact with common constant mass gauge of the new target.
Set \(J_{\rm old}=(I_{\rm old}-I_{\rm old}^\infty)/r\) and let
\(\Theta\) be its angle displacement. At these trial actions and
angles take the old field \(V\) and its column matrix. Let \(W,Y\)
denote respectively the normalized mobile transverse coordinates
and axial Birkhoff planes of the exact finite comparison. Add an
unrestricted transverse error \(r h\), an axial error \(r a\),
and the quadratic particular solutions in \eqref{eq:forcingpol-ins}
--\eqref{eq:exact-h-ins}. The pure quadratic sources, divided by
\(r^2\), are
\[
 P_\perp\{2V(|Y+W|^2-|Y|^2)+\bar V((Y+W)^2-Y^2)\},\qquad
 P_0\{2V(|W|^2-\textstyle\sum_v|W_v|^2)+\bar V W^2\}.
\]
The removed scalar density acts on the entire axial field as its
phase generator, so it does not change old actions. For a fusion the
first two kinetic tests are \(-2R^2p\cdot q,\,
2R(b-1)(p_x+q_x)\); for a reverse they are \(2R^2(p-q)\cdot q,\,
2R(1-b)(p_x-q_x)\), with the common sign convention for the two
legs. After these tests precisely the nondegenerate graph edges at
\((b,\nu)=(1,1)\) survive. In detail, a fusion of the mobile waves at
\(k_p,k_q\) against a conjugated old harmonic of spatial index \(b\)
and frequency \(\nu\) has output \(k=(2-b,0)+R(p+q)\), and a reverse
term with the wave at \(k_q\) conjugated has output
\(k=(b,0)+R(p-q)\). Since \(k_v=(1,0)+Rv\),
\[
\begin{aligned}
 |k_p|^2+|k_q|^2-|k|^2&=-2R^2p\cdot q+2R(b-1)(p_x+q_x)+2-(2-b)^2,\\
 |k_p|^2-|k_q|^2-|k|^2&=2R^2(p-q)\cdot q+2R(1-b)(p_x-q_x)-b^2 .
\end{aligned}
\]
The terms of order \(R^2\) vanish exactly on the forward and reverse
triangles. The horizontal components have distinct nonzero absolute
values, so \(p_x\pm q_x\ne0\) for distinct velocities and the terms of
order \(R\) vanish only for \(b=1\). The constants \(1\) and \(-1\)
that remain at \(b=1\) are matched by \(\nu=1\) alone, up to the
second-order shifts treated above. Their entire coefficient cancels with the
exact finite equation. All other sources are divided against the
bare nonzero divisor at angular degree \(\le C_1\log R\). A pure
low-horizontal transverse particular has divisor of order \(R^2\),
paying for the bounded low body; the other transverse particulars
use the scalar dressed columns. Nonzero axial velocities have exact
empty frequencies \(|k|^2+O(R^{-2})\). The differentiated envelopes
of the particulars cost only powers of \(\log R\): the fast
derivative of the carrier itself is taken together with its kinetic
divisor. The old actual shifts in that divisor are \(O(R^{-2}\log R)\).
The discarded tail can cost \(O(R^{-N_0})\) for any fixed \(N_0\).
In particular a bounded nonzero kinetic divisor contributes a
normalized error \(o(1)\), while an anomalous homogeneous off-key
resonance has been included separately in the normal operator.

Use on a late slow half-line an exponential weight in the stated
gap. The mobile waves and the differences from their axial endpoint
there are bounded in weighted norm by an arbitrarily small \(c\).
Put \(H_1,A_1,J_1,\Theta_1\) for the four weighted sup norms of
\(h,a,J_{\rm old},\Theta\), and \(X=H_1+A_1+J_1+\Theta_1\).
The Green images satisfy
\begin{align*}
 H_{\rm new}&\le C(e_R+cX+\sqrt R X^2+o_R(1)H_1),\\
 A_{\rm new}&\le C(e_R+J_1+c(H_1+A_1+\Theta_1)+\sqrt R X^2+o_R(1)A_1),\\
 J_{\rm new}&\le C(e_R+cH_1+c^2 A_1+\sqrt R X^2),\\
 \Theta_{\rm new}&\le C(e_R+J_1+r(A_1+H_1)+\sqrt R X^2),
 \qquad e_R=C R^{-1}(1+\log R)^C .
\end{align*}
Here the pure column error is \(O(R^{-3}/r)\) and the divided
sources have the error just estimated. A single prescribed factor
acts by its Wiener norm; two free factors use the \(\sqrt R\)
bound, and a third free factor has an additional \(r\).
In the first row the full trial-dependent Hermitian low body is
on the left. For each fixed \(q\), the flattened \(Z^p\) weight is constant
on this block, which is disjoint from the retained anomalous pairs.
Its propagator along the trial trajectory is therefore isometric, giving
the same backward kernel bound \(e^{-\sigma(u-\theta)}\) for \(u\ge\theta\)
as in Lemma~\ref{lem:common-endpoint-green}. Its order-one angle derivative is never applied
as an error to a free low field. Changing the old mass and
the nominal chart velocity costs \(O(r|J_{\rm old}|+r^2)\);
the derivative of the normal change costs this times \(O(R^{-1})\).
After division by the slow rate \(r\lambda\), the variation of the
old-mass scalar from its terminal value is \(2\lambda^{-1}\sum_{j\in F}(J_{\rm old})_j\) on
the positive transverse branch. Its product with the prescribed transverse
packet and its error has weighted norm at most \(C_{F,\lambda}J_1(c+H_1)\).
This gives \(O(cX+X^2)\) in the first row at the fixed preparation.
Every angle error against a prescribed transverse wave carries
its mobile factor. In the second row keep the whole real diagonal
at the trial actions. Their distance in \(\ell^1\) from the terminal
actions is \(O(r|J_{\rm old}|+r^2(c+A_1+o_R(1))e)\), where \(e\)
is the decaying exponential. Its multiplication by a persistent
root accounts for the unsuppressed \(J_1\). For the third row,
the equal-wave density changes no action. After truncation of the old angular harmonics, the remaining
leading axial terms quadratic in the prescribed transverse carriers have high
horizontal output. Their projection onto a fixed old coordinate gains
\(O(R^{-1})\) by the chart differential and its symplectic adjoint.
Three prescribed transverse waves can produce low output. An internal
triangle \(p+q=h\) gives \(k_p+k_q-k_h=(1,0)\). Each wave has
physical size \(O(rce)\) in the Wiener norm, so after division by
\(r^2\lambda\) their cubic contribution is \(O((r/\lambda)c^3e^3)\).
Since \(r=\eta R^{-2}\) with \(\eta,\lambda\) fixed, this is absorbed
by \(e_R e\). A generic transverse variation there needs a
mobile transverse factor; varying the axial field needs two such
factors. In the last row the occupied shear gives \(J_1\);
external forcing divided by \(r\lambda\), in place of \(r^2\lambda\),
gains \(r\). Pure errors vanish at the exact terminal axial solution.

\smallskip
\noindent\emph{Derivation of the four rows.} The proof follows that
of \eqref{eq:four-rows-ins} line by line; we record the changes and
then list the terms. Time is \(\theta=r\lambda t\), so
\eqref{eq:exact-h-ins} holds with \(\tau,r^{-1},r^{-2}\) replaced by
\(\theta,(r\lambda)^{-1},(r^2\lambda)^{-1}\). The spaces \(Z,B\) are
replaced by \(Z^p\) and its Wiener counterpart, \(\sqrt L\) by
\(C_p\sqrt R\), and \(\epsilon\) by \(c\). The weight has no residue
factor, so no term \(L^{-1}\) occurs in \eqref{eq:chart-head} or in the
frame estimates, and \eqref{eq:chart-tube}--\eqref{eq:chart-head} read
\(O(r(1+R^{2-p}))\) and \(O(R^{-1})\). Write
\(e=e^{-\sigma(\theta-\theta_0)}\). The constant \(C\) may depend on
the fixed preparation, on \(\lambda\) and on the weight \(\sigma\). The
transverse row uses the bound \(2C_J\) of step (iv) in the proof of
Lemma~\ref{lem:common-endpoint-green}, on \([\theta_0,S]\) with the
complementary trace at \(S\) and uniformly in \(S\); the other three
rows use backward integration from the terminal time and the factor
\(\sigma^{-1}\), as in Section~\ref{sec:ins-terminal}.

The essential difference from Section~\ref{sec:insertion} is that
\(R^{-2}/r=\eta^{-1}\) is not small. Each term bounded there by
\(R^{-2}/r\) has been removed before the estimate. First, the error
of the old frame is now \(O(R^{-4})\) off the diagonal and
\(O(R^{-3})\) on the high columns. Second, the scalar diagonal
\(\Xi(v)/R^2\) of a known transverse mode and the exact Birkhoff shift
of a known axial mode are part of its carrier frequency, and the edge
identities hold exactly for the assigned rows. Third, the detuning of
a resonant pair is kept in the principal operator. A divided particular has physical size
\(O(r^2(1+\log R)^C)\). The exponent \(\mathcal B\) has Schur norm
\(O(r(1+\log R)^C)\), so its action on a prescribed field of size
\(O(r)\) has size \(O(r^2(1+\log R)^C)\). In the homogeneous
generator, the divisor mismatch \(O(R^{-2}(1+\log R))\) multiplies
\(\mathcal B\). After division by \(r\lambda\) its contribution is
\(O(\lambda^{-1}R^{-2}(1+\log R)^{C+1})=o_R(1)\) at the fixed clock.

\emph{First row.} The forcing of the equation for \(h\) consists of
the following terms, each measured after division by \(r^2\lambda\)
and application of the Green operator.
\begin{itemize}
\item The residual of the prescribed fields. The high columns and the
two Sylvester steps give \(O((R^{-3}+R^{-4})/(r\lambda))=O(R^{-1})\) on
a known wave of amplitude \(r\). The chart error \(\mathfrak e_0\)
against the known multiplier \(W\) gives
\(O(R^{-1}+r(1+R^{2-p}))\). A slow derivative of a divided particular
gives \(O(r(1+\log R)^C)\), and the particulars with bounded nonzero
divisors have normalized size \(O(r(1+\log R)^C)\). The angular tail
is \(O(R^{-N_0})\). The products of three known factors carry an
additional \(r\). All these terms are \(O(e_R)e\).
\item Terms linear in the errors with a decaying prescribed factor:
replacing \(Y_\infty\) by \(Y\) in the mixed operator, \(O(c)H_1\); the
axial error \(\mathfrak e_1\) against \(W\), \(O(c(A_1+rc))\);
freezing the prescribed harmonics at the nominal old angles and
actions, \(O(c(\Theta_1+rJ_1+r))\); the variation of the old-mass
scalar on the prescribed packet, \(C_{F,\lambda}cJ_1\). Their sum is
\(O(cX+e_R)e\).
\item Terms with two free factors: the old-mass scalar on \(h\),
\(C_{F,\lambda}J_1H_1\); \(\mathfrak e_1\) against \(h_1\) and two
unknown transverse factors against a smooth old factor, each
\(O(\sqrt R\,X^2)\) by the scaled algebra estimate; the cubic term in
the errors, \(O(rRX^3)\le O(\sqrt R\,X^2)\) when \(\sqrt R\,X\le1\).
\item Terms linear in \(h\) with a coefficient tending to zero:
\(\mathfrak e_0\) against \(h_1\), \(O(\sqrt R\,e_R)H_1\); the
difference of the actual and nominal time connections of the frame,
which by \eqref{eq:oldconnection} with \(R^{-1}\) in place of
\(L^{-1}\), and after division by \(r\lambda\), is
\(O(R^{-1}(J_1+r))H_1\); and the variation \(O(R^{-2})(r|J_{\rm old}|+r^2)\) of the second-order
diagonal, which after the same division is \(O(R^{-2})H_1\).
\end{itemize}
The nonresonant homogeneous terms and the averaged groups are part of
the remainder \(\mathcal R\) and have been absorbed by the Neumann
series. This gives the first row.

\emph{Second row.} Keep the real diagonal at the trial actions and
integrate backwards. The \(\ell^1\) distance of the trial actions from
the terminal actions is bounded by \(C(rJ_1+r^2(c+A_1+o_R(1)))e\), by
the Cauchy--Schwarz estimate of Section~\ref{sec:ins-terminal}. On a
persistent root this changes the frequency by at most the same amount,
and the mismatch contributes
\[
 \frac{|\Delta\omega_i|\,|rY_{K_i}|}{r^2\lambda}
 \le C\lambda^{-1}\bigl[J_1+r(c+A_1+o_R(1))\bigr]e .
\]
This is the unsuppressed \(J_1\). The force \(f\), after removal of
the density rotation, contains two or three transverse factors: two
known factors are cancelled by the finite system at the key \((1,1)\)
and are divided otherwise, with error \(O(e_R)\); one known and one
unknown factor give \(O(cH_1)\); two unknown factors give
\(O(\sqrt R\,X^2)\). Replacing \(D\Phi(U)\) by \(D\Phi(V_*)\) and then
by its Fourier part on the known high outputs costs
\(O(r(1+R^{2-p})+R^{-1})\le O(e_R)\), and on the outputs containing
an error it costs \(o_R(1)A_1\) or a quadratic term. Evaluating the
prescribed harmonics at the nominal old angles costs
\(O(c(\Theta_1+rJ_1+r))\). This gives the second row.

\emph{Third row.} In
\(\partial_\theta J_j=2(r^2\lambda)^{-1}\operatorname{Re}(\bar z_j(D_f)_j)\),
\(j\in F\), the density rotation drops out. Every other product of two
known factors has horizontal output of modulus \(\gtrsim R\), so its
projection to \(F\) through \(D\Phi(U)\) is bounded by
\(O(R^{-1}+r(1+R^{2-p}))\le O(e_R)\), by
\eqref{eq:chart-tube}--\eqref{eq:chart-head} and symplectic
adjointness. A transverse error enters together with a known transverse
factor, which gives \(O(cH_1)\), or with a second error, which gives
\(O(\sqrt R\,X^2)\). An axial error enters only through the dependence
of \(U\) and of \(D\Phi(U)\) on it, multiplied by a transverse pair;
two known factors give \(O(c^2A_1)\), and the other cases are
quadratic. No term linear in \(\Theta\) or \(J_{\rm old}\) without a
small factor occurs, because the estimate is made at the actual old
angles and actions.

\emph{Fourth row.} In
\(\partial_\theta\Theta_j=(r\lambda)^{-1}\bigl(\omega_j(I)-\omega_j(I^\infty)
-\operatorname{Im}((D_f)_j/z_j)\bigr)\) the frequency difference is at
most \(C\|I-I^\infty\|_{\ell^1}\), which by the bound used for the
second row gives \(C\lambda^{-1}[J_1+r(c+A_1+o_R(1))]e\).
The mean density \(2r^2S_W U\) contributes
\(D\Phi(U)[-2ir^2S_W U]=-2ir^2S_W z\), hence its contribution to
\(\partial_\theta\Theta_j\) is \(2rS_W/\lambda\).
Since \(S_W\lesssim c^2e^2\), \(\lambda\) is fixed and \(r\lesssim e_R\),
this term is \(O_\lambda(rc^2e^2)=O(e_R e)\).
For the remaining force, the lower bounds on the old \(|z_j|\)
allow the product estimates from the third row to be used with
\((r\lambda)^{-1}\) in place of \((r^2\lambda)^{-1}\).
The resulting extra factor \(r\) gives the stated fourth-row bound.

Finally, with all four errors equal to zero, with \(Y\) replaced by
\(Y_\infty\) and with \(W\) and the particulars removed, the
reconstruction is the exact terminal axial solution and all four
residuals vanish. Hence every pure error carries the factor \(e\).

With \(\rho=R^{-3/5}\) take radii \(M_1\rho\) in the first and third
rows and \(M_2\rho\) in the other two. Choose \(M_1\), then \(M_2\)
relative to the shear, then \(c\) small, and finally \(R\) large.
Both \(e_R/\rho\) and \(\sqrt R\,\rho\) tend to zero. To check the
invariance explicitly, let \(C\ge1\) be the constant in the four
displayed rows and suppose \(H_1,J_1\le M_1\rho\) and
\(A_1,\Theta_1\le M_2\rho\), so that \(X\le2(M_1+M_2)\rho\). Then
\[
\begin{aligned}
 H_{\rm new}/\rho&\le C\bigl(2c(M_1+M_2)+o_R(1)\bigr),
 & J_{\rm new}/\rho&\le C\bigl(cM_1+c^2M_2+o_R(1)\bigr),\\
 A_{\rm new}/\rho&\le C\bigl(M_1+c(M_1+2M_2)+o_R(1)\bigr),
 & \Theta_{\rm new}/\rho&\le C\bigl(M_1+o_R(1)\bigr),
\end{aligned}
\]
where \(o_R(1)\) collects \(e_R/\rho\), the terms
\(\sqrt RX^2/\rho\le4(M_1+M_2)^2\sqrt R\,\rho\), the terms with a
factor \(r\), and the \(o_R(1)\) multiples of \(H_1,A_1\). With
\(M_1=1\), \(M_2=4C\) and \(c\le1/(16C(1+M_2))\) the four right
sides are at most \(M_1/2,\ M_1/2,\ M_2/2\) and \(M_2/2\) for large
\(R\). The same bounds applied to differences are not needed, since
the argument below uses compactness rather than contraction. On a finite
interval ending at \(S\), use the stable projection forward and the
remaining components backward. The low unitary is evaluated along
the trial, using the instantaneous velocity from its coordinate
equations in the time connection, as in Section~\ref{sec:insertion}.
At fixed \(R,S\) the image has bounded time derivative in \(Z^0\)
by its \(Z^2\) bound. Interpolation, compact embedding on the fixed
lattice, and the finite-time Duhamel formula give a compact convex
invariant set in the \(C_t Z^{p_0}\) topology and continuity of the
map there. Schauder's theorem and then local compactness as
\(S\to\infty\) give an exact solution. Strong \(H^2\) persistence
and uniqueness identify the physical limit. The same estimates with
bounded coefficients and finite Gronwall on the preceding fixed slow
interval give an \(O(\rho)\) normalized launch error. This constructs
an exact solution \(W_R\) exponentially approaching the new axial
target \(V'\) on its invariant lattice.

\smallskip
\noindent\emph{The compactness argument in detail.} Fix \(R\) and
\(S>\theta_0\), and set \(J_S=[\theta_0,S]\). In the trial set
\(\mathcal B_S\) the fields have regularity
\(C(J_S;Z^{p_0})\cap L^\infty(J_S;Z^2)\), with the corresponding axial
spaces for \(a\), and the finite-dimensional variables \(J_{\rm old},\Theta\)
are continuous. Impose the four weighted bounds above at each time.
These bounds are closed in the lower path topology: the squared
higher norm at each time is the supremum of its finite Fourier sums.
Let \(\mathcal M\) be the map given by the four Green images, with
the fast stable components prescribed at \(\theta_0\) and the remaining
components prescribed at \(S\). The finite-interval Duhamel formulas
define this map at order \(p_0\). The higher bounds also hold on
\(\mathcal B_S\): Fourier truncations of a trial converge uniformly on
\(J_S\) at order \(p_0\) and preserve its higher bounds, so continuity of
the formulas at order \(p_0\) and lower semicontinuity give the same
bounds on the image. Thus \(\mathcal M(\mathcal B_S)\subset\mathcal B_S\).
Every field image is bounded pointwise in \(Z^2\) and, by its equation,
has time derivative bounded in \(Z^0\), with bounds depending on \(R,S\). Interpolation gives, for \(\theta,\theta'\in[\theta_0,S]\),
\[
 \|f(\theta)-f(\theta')\|_{Z^{p_0}}
 \le\|f(\theta)-f(\theta')\|_{Z^0}^{1-p_0/2}
    \|f(\theta)-f(\theta')\|_{Z^2}^{p_0/2}
 \le C_{R,S}|\theta-\theta'|^{1-p_0/2},
\]
so the images are equicontinuous with values in \(Z^{p_0}\). On the
fixed lattice the ratio of the weights of \(Z^{p_0}\) and \(Z^2\) tends
to zero at infinity, so bounded sets of \(Z^2\) are precompact in
\(Z^{p_0}\); the same holds for the axial components in \(h^2\subset
h^{p_0}\), and the old components are finite dimensional. By
Arzel\`a--Ascoli the closed convex hull \(\mathcal K_S\) of
\(\mathcal M(\mathcal B_S)\) is compact in \(C([\theta_0,S];Z^{p_0})\).
It is contained in \(\mathcal B_S\), which is convex and closed in this
topology since the weighted bounds are lower semicontinuous. The map
\(\mathcal M\) is continuous on \(\mathcal K_S\) for this topology: the
products are continuous from \(Z^{p_0}\times Z^{p_0}\) to \(Z^{p_0}\) by
the scaled algebra estimate with \(p=p_0>1\), the coefficients depend
continuously on the finite-dimensional old variables, and the Green
operators on \([\theta_0,S]\) are bounded on the weighted
\(Z^{p_0}\) spaces. Schauder's theorem gives a fixed point
\(X_S\in\mathcal K_S\). The weighted bounds are uniform in \(S\). Hence
a diagonal sequence \(S_k\to\infty\) converges in
\(C([\theta_0,S'];Z^{p_0})\) for every \(S'\); the limit obeys the same
bounds, and it satisfies the Duhamel equations on the half-line because
the backward integrals converge by the exponential weights. After
reconstruction, the limit satisfies the physical Duhamel formula in
\(H^{p_0}\) and takes values in \(H^2\) at every time by the pointwise
bounds. Starting at any such time, the local \(H^2\) solution with that
datum agrees with the limit by uniqueness in \(C_tH^{p_0}\). This gives
strong \(H^2\) continuity on each compact time interval.

At fixed old target \(V\) and fixed height \(L\), take first the
depth to make the launch arbitrarily small, followed by the other
finite choices and finally \(R\). We obtain
\[
 W_R(0)\to V(0)\ \hbox{in }H^2,\qquad
 \sup_{t\ge0}\|W_R(t)-V(t)\|_{H^{p_0}}\to0,\qquad
 \sup_{R,t\ge0}\|W_R(t)\|_{H^2}<\infty
\]
along the resulting family in its respective mass gauges. Indeed
the old rates and origins agree exactly, the incoming clock cost is
arbitrarily small, \(rR^{p_0}\to0\), and the strong bound on all
primary levels above is independent of depth at fixed height. Apply
one-smoothing of \(\Psi-\mathcal F^{-1}\) from \(h^1\) to \(H^2\)
at the terminal torus. Its old image is compact in \(H^2\) and the
terminal high inputs have \(h^1\) size \(O(rR)\). The normalized
physical modulus at \(1+RP\) is thus \(L+o_R(1)\) uniformly in
time there. By exponential landing it persists on \(W_R\) with
arbitrarily small fixed slack at sufficiently late times.

\smallskip
\noindent\emph{The three limits in detail.} With \(W_R\) and \(V\) in
their respective mass gauges, decompose \(W_R(t)-V(t)\) into five parts
and estimate each uniformly in \(t\ge0\). Since \(\bk{(n,Rq)}\le CR(1+|n|/R+|q|)\), a field \(f\) on
the lattice satisfies \(\|f\|_{H^b}\le CR^b\|f\|_{Z^b}\le CR^b\|f\|_{Z^2}\)
for \(0\le b\le2\).
(a) \emph{The normalized errors} \(h,a\) have \(Z^2\) norms at most
\(M_2\rho\) on the late half-line and \(O(\rho)\) on the preceding
fixed slow interval, by the backward Gronwall estimate. In physical
units they carry the factor \(r=\eta R^{-2}\), so their \(H^2\) norm is
\(O(R^2r\rho)=O(\eta\rho)\) and their \(H^{p_0}\) norm is
\(O(R^{p_0}r\rho)\); both tend to zero since \(\rho=R^{-3/5}\).
(b) \emph{The particulars} are \(r^2\) times Wiener-bounded fields
with divisors at least \(c(1+\log R)^{-C}\), so their \(H^2\) norm is
\(O(R^2r^2(1+\log R)^C)=O(\eta r(1+\log R)^C)\to0\). The frame
corrections \((\mathcal S-1)\) on the known planes are \(O(R^{-1})\),
which costs \(O(R^2rR^{-1})=O(\eta R^{-1})\) in \(H^2\).
(c) \emph{The finite packet.} Write \(x_v(\theta)\) and
\(x_v^{\rm id}(\theta)\) for the fitted and ideal amplitudes in the
carrier coordinates, assigning zero ideal amplitude to unused heads.
By orthogonality the squared \(H^b\) norm of the fitted packet is
\(\sum_v r^2\bk{k_v}^{2b}|x_v(\theta)|^2\). At fixed tree and height put
\[
 e_d(\lambda,R)=\sup_{\theta\ge0}
 \left(\sum_v\eta^2|v|^4
       \bigl(|x_v(\theta)|-|x_v^{\rm id}(\theta)|\bigr)^2\right)^{1/2}.
\]
The weighted contraction on the late half-line gives continuous
dependence of the decaying part in its weighted norm, and the rounded
endpoint converges to the ideal endpoint. On the preceding fixed interval
Gronwall gives uniform convergence in the same iterated limit. Undoing
the finite change of coordinates costs \(O_{\rm tree}(\eta\lambda)\).
Since the tree is finite, these estimates give
\(\lim_{\lambda\downarrow0}\limsup_{R\to\infty}e_d(\lambda,R)=0\).
For \(b=2\), \(r\bk{k_v}^2=\eta|v|^2(1+O_{\rm tree}(R^{-1}))\).
The triangle inequality therefore bounds the launch norm by
\(L(1+\sum_i b_i^2)^{1/2}c^{d/2}+e_d(\lambda,R)+o_R(1)\), plus
the prescribed private cost. At every time the corresponding norm
bound is \(L(1+\sum_i b_i^2)^{1/2}(1-c)^{-1/2}
+e_d(\lambda,R)+o_R(1)\), plus the same uniform private cost.
For \(b=p_0\) the coefficient is
\(r\bk{k_v}^{p_0}=O(\eta R^{p_0-2})\), so at fixed depth the \(H^{p_0}\)
norm is \(O_d(R^{p_0-2})=O_d(rR^{p_0})\to0\).
(d) \emph{The hot axial shape.} The nonlinear part
\(\Psi-\mathcal F^{-1}\) is Lipschitz from a small \(h^1\) ball to
\(H^2\) \cite[Corollary~1.1]{KST}, and the hot coordinates have \(h^1\)
size \(O(rR)=O(\eta R^{-1})\). Hence the axial field differs from
\(\Psi(z_F)\) plus the hot planes by \(O(\eta R^{-1})\) in \(H^2\).
(e) \emph{The old field.} The old rates and angle origins are locked
exactly, so the reference angles of \(V\) and of the new old
coordinates coincide for every \(t\) and no linear drift occurs. The
old coordinates of \(W_R\) differ from those of \(V\) by the angle
displacement \(\Theta=O(\rho)\), the action displacement
\(rJ_{\rm old}=O(r\rho)\) and the preparation changes; the latter
consist of the locked old action change and the moderate and weak
sites, whose \(h^2\) cost was made arbitrarily small. The old torus is
a fixed smooth finite-dimensional family, on which \(\Psi\) is smooth
into \(H^2\), so this part is \(O(\rho+r\rho)\) plus the preparation
tolerance, uniformly in \(t\).
Given a launch tolerance, first choose the depth to bound the ideal
launch cost. At that depth choose the private costs and the remaining
finite preparation, then choose \(\lambda\) small and \(R\) large to
put the fitting error in (c) and the reconstruction errors within the
tolerance. Letting these tolerances tend to zero gives
\(W_R(0)\to V(0)\) in \(H^2\). For \(H^{p_0}\) all parts are small at
all times by (a), (b), (c) with \(b=p_0\), (d) and (e). Taking
\(\sum_i b_i^2\le1\) at each depth, the geometric bound in (c) and
these choices give a uniform \(H^2\) bound at fixed height. Finally, at the root
\(K=1+RP\) the packet contributes \(r\bk K^2|H(\theta)|
=L|H(\theta)|(1+O(R^{-1}))\), because \(\eta=L/P^2\). Every other part
changes the normalized coefficient \(\bk K^2|\widehat f(K,0)|\) by at
most its \(H^2\) norm, and the old field contributes
\(\bk K^2|\widehat V(K,0,t)|\to0\) uniformly in \(t\) by precompactness
of its range in \(H^2\). The fitted root modulus converges exponentially
to \(|H_{*,R}|=1+o_R(1)\) at the rounded endpoint. By the uniform
terminal-torus estimate above and exponential landing, for every
\(\delta>0\) the normalized root modulus of \(W_R(t)\) lies within
\(\delta\) of \(L\) for all sufficiently large \(R\) and all
\(t\ge T(R,\delta)\).

We finish by specifying how the historical differences in
Section~\ref{sec:limit} are reused here. Choose increasing heights
\(L_j=4j+4\) and take each new dilation divisible by every old one.
For a saved axial reference \(V\) use its normal inverse just proved
at two slow-time weights \(0<\widetilde\nu<\widetilde\beta\) in one
gap, with the same boundary conditions in \(H^{p_0},H^2\). With
\(r,\lambda\) fixed at this saved stage, set
\(\nu=r\lambda\widetilde\nu\) and \(\beta=r\lambda\widetilde\beta\)
for the corresponding physical-time rates. Suppose \(Y+d\) is an
earlier exact solution in the gauge of \(Y\), with
\(d=O_{H^2}(e^{-\beta t})\)
and \(Y\) in its saved uniform weak tube about \(V\). The new inputs
\(Y_n\) converge to \(Y\) strongly initially and uniformly weakly
on the half-line, have bounded strong norms at this fixed height,
and their gauges converge. Solve the equation for a replacement
correction \(d+h\) by writing
\[
 \begin{aligned} h&=\mathcal G_V\{F_n(d+h)-F(d)\},\\
 F_n(b)&=-i\big[(c_n-c_V)b+2(|Y_n|^2-|V|^2)b+
 (Y_n^2-V^2)\bar b+2Y_n|b|^2+\bar Y_n b^2+|b|^2b\big].\end{aligned}
\]
Here \(F\) uses \(Y\) and its gauge. To derive this equation, let
\(Y_n+b\) solve the cubic equation in the gauge \(c_n\) of \(Y_n\).
Subtracting the equation for \(Y_n\) and then the linearization
\(i\partial_tb=(-\Delta+c_V+2|V|^2)b+V^2\bar b\) at \(V\) gives
\[
 i\partial_tb-(-\Delta+c_V+2|V|^2)b-V^2\bar b
 =(c_n-c_V)b+2(|Y_n|^2-|V|^2)b+(Y_n^2-V^2)\bar b
   +2Y_n|b|^2+\bar Y_nb^2+|b|^2b ,
\]
that is \(\partial_tb=\mathcal L_Vb+F_n(b)\) with the linearized
generator \(\mathcal L_V\). The old correction satisfies
\(\partial_td=\mathcal L_Vd+F(d)\). Hence \(b=d+h\) solves the new
equation if and only if
\(\partial_th=\mathcal L_Vh+F_n(d+h)-F(d)\), and \(\mathcal G_V\)
inverts \(\partial_t-\mathcal L_V\) with the prescribed boundary
conditions. In the sum of weighted strong
norm and \(\Lambda\) times weighted weak norm the tame product
estimate bounds the linear coefficient by its small uniform weak
size plus \(C(1+\sup_n\|Y_n\|_{H^2})^2/\Lambda\). For instance, if
\(q_n=|Y_n|^2-|V|^2\) has uniform weak norm \(\epsilon\) and
\(\|h\|_*=\|h\|_{H^2}+\Lambda\|h\|_{H^{p_0}}\), then
\[
 \|q_nh\|_*\le C\bigl(\epsilon\|h\|_{H^2}
      +\|q_n\|_{H^2}\|h\|_{H^{p_0}}\bigr)
      +C\Lambda\epsilon\|h\|_{H^{p_0}}
 \le C\Bigl(\epsilon+\frac{\|q_n\|_{H^2}}{\Lambda}\Bigr)\|h\|_* ,
\]
and \(\|q_n\|_{H^2}\le C(1+\sup_n\|Y_n\|_{H^2}+\|V\|_{H^2})^2\).
Choose \(\Lambda\) for this fixed-height family and start the half-line
sufficiently late to make the shifted \(\beta\)-norm of \(d\)
arbitrarily small. This gives a contraction at \(\beta\).
At \(\nu\) its inhomogeneity tends to zero by finite-time strong
continuous dependence and the extra \(e^{-(\beta-\nu)t}\) tail.
Compatibility of the inverses gives convergence for the same
correction. Propagate backwards over the fixed initial interval.
This preserves the original rate \(\beta\) for future reuse.
We spell out the last steps. Work on \([A,\infty)\) with the weights
\(e^{\xi(t-A)}\), \(\xi\in\{\nu,\beta\}\), and the norms
\(\|h\|_{*,\xi}=\|h\|_{H^2,\xi}+\Lambda\|h\|_{H^{p_0},\xi}\). Write
\[
 F_n(d+h)-F(d)=\bigl[F_n(d+h)-F_n(d)\bigr]+\bigl[F_n(d)-F(d)\bigr].
\]
The first bracket is linear in \(h\) with the coefficients
\(c_n-c_V\), \(|Y_n|^2-|V|^2\), \(Y_n^2-V^2\) and terms containing \(d\), plus terms
quadratic and cubic in \(h\). By the displayed tame estimate the
first coefficients cost \(C(\epsilon+\sup_n\|Y_n\|_{H^2}^2/\Lambda)\) in either
norm, and every term containing \(d\) costs \(C\sup_{t\ge A}\|d(t)\|_{H^2}\le
Ce^{-\beta A}\). The second bracket equals
\((c_n-c_Y)d+2(|Y_n|^2-|Y|^2)d+(Y_n^2-Y^2)\bar d+2(Y_n-Y)|d|^2+(\bar Y_n-\bar Y)d^2\),
up to the factor \(-i\); every term contains \(d\). Hence, with
\(\Lambda\) chosen first and then \(A\) large, the map is a contraction of a
small ball in both norms, uniformly in \(n\), and its fixed point
\(h_n\) satisfies
\(\|h_n\|_{*,\xi}\le2\|\mathcal G_V\|\,\|F_n(d)-F(d)\|_{*,\xi}\). At \(\xi=\beta\)
the right side is \(O(e^{-\beta A})\); the solution therefore lies in the
\(\beta\)-weighted space and keeps the rate \(\beta\). At \(\xi=\nu\) the right
side tends to zero as \(n\to\infty\): on \([A+T,\infty)\) each term is at most
\(C(\sup_n\|Y_n\|_{H^2}+\|Y\|_{H^2})^2e^{\nu(t-A)}\|d(t)\|_{H^2}\le
Ce^{-\beta A}e^{-(\beta-\nu)T}\), small for large \(T\); on \([A,A+T]\) the
coefficient differences tend to zero in \(H^2\) by finite-time strong
continuous dependence, because \(Y_n(0)\to Y(0)\) in \(H^2\) and \(c_n\to c_Y\).
The fixed points at the two weights coincide, since the inverses are
compatible and the \(\beta\)-ball lies in the \(\nu\)-ball. Thus
\(d+h_n\to d\) in the \(\nu\)-weighted \(H^2\) norm on \([A,\infty)\), and
\((Y_n+d+h_n)(A)\to(Y+d)(A)\) in \(H^2\). Continuous dependence backwards on
\([0,A]\) gives convergence there, and in particular of the data at \(t=0\).

At each stage first lift the newest saved difference about its
axial reference, then lift the finitely many earlier differences in
reverse order. All inputs belong to the requisite older lattices.
Choose a sufficiently late member of the fixed-height family to make
the final initial \(H^2\) differences summable and all older weak
tube changes summably small. At each previously attained frequency
\(K_i\) require the remaining sum of weak changes times
\(\langle K_i\rangle^{2-p_0}\) to be less than \(i\). Each new
full solution then has normalized modulus \(>2i\) at its \(i\)th
frequency for all \(t\ge T_i\), and the limit retains modulus
\(\ge i\) there on that whole interval. Strong initial convergence
and finite-time well-posedness identify the limit with a single
global \(H^2\) solution. Undoing its mass gauge gives
\(\|u(t)\|_{H^2}\to\infty\). An arbitrarily small starting axial torus
and summably small initial increments give an arbitrarily small datum.

\paragraph{Unit heights and the counting assertion at the endpoint.}
Alternatively set \(L=1\) at every insertion. We give the
additional uniform-in-time comparison needed for the counting law.
Let \(\delta>0\) be prescribed after the previous orbit has been
fixed. On the entire half-line the raw change from the old axial
target is the sum \(P_R(t)\) of its finitely many pure new carrier
planes, with their actual comparison coefficients, plus an error of
\(H^2\) norm at most \(\delta\). To verify this, choose \(d\) and
the private costs first. The old occupied frequencies and origins
agree exactly; the finite old actions and the moderator and weak
inputs may change by arbitrarily little in \(h^2\).
The physical shape error for the high axial insertions tends to
zero by one-smoothing since their \(h^1\) cost is \(O(rR)\) at
this fixed preparation. The normalized transverse columns differ
from their individual planes by \(O(R^{-1})\) on the finite set;
the physical cost of the scaled remainder is \(O(rR^2\rho)=o(1)\).
The nonlinear old modulations, of size \(O(\rho)\) in the old
angles and \(O(r\rho)\) in their actions, obey the same conclusion
uniformly by compactness of this fixed old torus. The difference
of each old history correction from its replacement tends to zero
in uniform \(H^2\) by its convergence at the smaller positive
weight and finite-time strong convergence. Thus after the finitely
many lifts the full change is still \(P_R+O_{L^\infty_tH^2}(\delta)\),
shrinking the input tolerance if necessary.

Write \(u_{\rm old}\) for the old trajectory and let \(P_{\ge M}\)
be the Fourier projection to \(\langle k\rangle\ge M\). The range of
\(u_{\rm old}\) is precompact in \(H^2\), since it approaches a
compact torus. Consequently
\(\sup_{t\ge0}\|P_{\ge M}u_{\rm old}(t)\|_{H^2}\to0\)
by Fourier truncation on a finite net of the range. All modes
of \(P_R\) lie beyond a threshold \(M_R\to\infty\), so
\[
 \sup_{t\ge0}|\langle u_{\rm old}(t),P_R(t)\rangle_{H^2}|
 \le \sup_{t\ge0}\|P_{\ge M_R}u_{\rm old}(t)\|_{H^2}
       \sup_{t\ge0}\|P_R(t)\|_{H^2}\longrightarrow0 .
\]
The second factor is bounded at unit height. Expanding the new squared
norm now gives the old squared norm plus \(\|P_R(t)\|_{H^2}^2\)
up to any prescribed uniform error. Explicitly, if the new trajectory
is \(u_{\rm old}+P_R+e\) with \(\sup_t\|e(t)\|_{H^2}\le\delta\), then
\[
 \Bigl|\|u_{\rm old}+P_R+e\|_{H^2}^2-\|u_{\rm old}\|_{H^2}^2
       -\|P_R\|_{H^2}^2\Bigr|
 \le2\bigl|\langle u_{\rm old},P_R\rangle_{H^2}\bigr|
   +2\delta\bigl(\|u_{\rm old}\|_{H^2}+\|P_R\|_{H^2}\bigr)+\delta^2 .
\]
To obtain that error, first choose
\(\delta\) using the old strong bound and then enlarge \(R\). Take these errors summable
over the stages. Set \(D_j=K_j^2\). Since \(K_j=1+RP\) and
\(r=\eta/R^2\) at the fixed preparation,
\[
 r\lambda D_j^{3/4}\le\eta\lambda(2P)^{3/2}R^{-1/2},\qquad
 r\lambda D_j^{5/4}\ge\eta\lambda P^{5/2}R^{1/2}.
\]
On \(t\le D_j^{3/4}\) the slow time
\(r\lambda t\) tends uniformly to zero as \(R\to\infty\) at fixed
preparation. Include \(e_d(\lambda,R)\) from (c) in the launch
tolerance by choosing \(\lambda\) before \(R\). The bounded finite
vector field then keeps \(P_R\) within the prescribed cold cost.
On \(t\ge D_j^{5/4}\) the slow time tends uniformly to infinity.
Exponential landing, \(|H_{*,R}|=1+o_R(1)\), and the prescribed
cost of the private terminal axes make the squared norm of \(P_R\)
differ from one by the prescribed tolerance.
On the intervening interval its norm is bounded by an absolute
constant, using balance and the unit-height level bounds above.
Select \(R\) to make these intervening intervals disjoint. Summing
the squared increments gives \eqref{eq:count}, since at most one
interval is in transition. Taking \(D_j^{3/4}\) sufficiently late
that \(G(D_j^{3/4})\ge j^2\) gives the prescribed \(o(G(t))\)
bound. The required slow entrance time at fixed preparation is
bounded independently of sufficiently large \(R\); taking \(R\)
last thus also gives \(t_j\le K_j^2\ell(K_j)\), with the slack in
each older Fourier coefficient preserved by the weak tubes.

\section{Exponents above two}\label{sec:super}

Fix \(2<s<\infty\). All constants in a single insertion may depend on
\(s\), the finite axial background and the depth of the tree. In this
section only the final spatial dilation is sent to infinity at that
insertion. This order of choices is useful for distinguishing the
physical and angular small divisors. We use the notation \(A,B\) of
Lemma~\ref{lem:empty-frequency} and the tree with its private
horizontal endpoints and its mixed vertices from
Proposition~\ref{prop:tree}. An insertion of unit height will mean
that the final Fourier amplitude on its main axis has \(H^s\)
weight \(1+o(1)\).

\subsection{Preparing the finite axial solution}

At a smooth axial torus of fixed finite support,
for every fixed \(L\ge2\) the empty-coordinate frequencies admit
the expansion
\begin{equation}\label{eq:sup-empty}
 \chi_n=\sum_{j=2}^{L}c_j n^{-j}+O(|n|^{-L-1}).
\end{equation}
Continue the spatial recursions of Lemma~\ref{lem:empty-frequency}
at \(E=n^2\) through order \(L\), as in Section~\ref{sec:h2-extension}.
The spatial mean of the upper equation gives the instantaneous scalar
\(a(\phi,n)=\sum_{j=2}^L a_j(\phi)n^{-j}\); its mean-zero part
determines \(f_{j+1}\) by \(2D_x\), and the lower equation determines
\(g_{j+2}\) by its leading term \(2n^2\). Keeping $f$ through order $L+1$ and $g$ through order $L+2$ gives a residual
\(O(|n|^{-L-1})\) in \(L^2_x\), uniformly in the angles \(\phi\).
With the carrier phase factored out, the positive empty component
\(c_n(\phi)\) of the column under the complexified \(D\Phi\) satisfies
\[
 c_n=1+O(|n|^{-1}),\qquad
 i\mathcal L c_n=(\chi_n-a(\phi,n))c_n+O(|n|^{-L-1}),
\]
where \(\mathcal L\) differentiates along the angular flow. The first
estimate follows from one-smoothing as in Lemma~\ref{lem:empty-frequency}.
For large \(|n|\), \(\log c_n\) is single-valued on the full angle
torus. After division by \(c_n\), its logarithmic derivative has zero
Haar mean, since \(\mathcal L\) is a constant vector field in the angles.
Thus \(c_j=\int a_j(\phi)\,d\phi\) with normalized Haar measure, also
for resonant frequencies. The Liouville identity gives \(c_2=A-B\).
For fixed \(L\), the recursion uses only finitely many derivatives,
spatial means and mean-zero spatial primitives. These and the chart
bounds above are uniform on compact families with fixed finite support.

Use the five weak letters \(h\) of Section~\ref{sec:h2-extension}.
Choose an integer \(m>\max(2,2s-2)\), a prime \(\pi>20\,2^{m-1}\)
and the neutral moderator word with generating polynomial
\[
 d_*(X)=X^K(1-X^\pi)\prod_{\ell=0}^{m-2}(1-X^{2^\ell}).
\]
The subset sums of the powers \(2^\ell\), \(0\le\ell\le m-2\), give each integer
\(0\le i<2^{m-1}\) once, and the factor \(1-X^\pi\) gives its
oppositely signed partner at \(i+\pi\). Write \(e_a\) for the
resulting sign at \(K+a\). Explicitly, if \(\sigma(i)\) denotes the
number of nonzero binary digits of \(i\), then
\[
 X^{-K}d_*(X)=\sum_{i=0}^{2^{m-1}-1}(-1)^{\sigma(i)}\bigl(X^i-X^{i+\pi}\bigr),
\]
and the \(2^m\) exponents are distinct because \(\pi>2^{m-1}\). Thus
\(e_i=(-1)^{\sigma(i)}=-e_{i+\pi}\), every site carries a coefficient
\(\pm1\), and the word is neutral. The zero of order \(m\) at \(X=1\)
gives \(\sum_a e_a a^j=0\) for \(0\le j<m\), by applying
\((X\partial_X)^j\) to \(X^{-K}d_*(X)\) at \(1\).
For every \(2\le l\le m+1\), Taylor expansion of \((K+a)^{-l}\)
through degree \(m-1\) in \(a\) therefore gives
\(\sum_a e_a(K+a)^{-l}=O(K^{-l-m})\). Indeed
\[
 (K+a)^{-l}=\sum_{j=0}^{m-1}\binom{-l}{j}a^jK^{-l-j}+O_{l,m}(a^mK^{-l-m}),
\]
the polynomial part is annihilated by the vanishing moments, and
\(0\le a<\pi+2^{m-1}\). The shifts are fixed as
\(K\) increases. The remainder in \eqref{eq:sup-empty} with
\(L=m+1\) is \(O(K^{-m-2})\), and the five weak sites cost \(O(T^{-2})\).
Together these estimates give
\[
 |(h+d_*)\cdot\chi|\le C(K^{-m-2}+T^{-2}).
\]
Here \(T\) is the common scale of the five weak letters. At fixed depth
use disjoint groups and comparable \(K\)'s. Take \(T\gg K^{1+m/2}\).
A positive common action buffer on the moderator sites of size
\(C_1(K^{-m-2}+T^{-2})\) costs \(O(K^{s-1-m/2})=o(1)\) in
\(H^s\). Indeed an action \(I\) at one of the fixed number of sites
of size comparable to \(K\) costs \(O(K^s\sqrt I)\), here
\(T^{-2}\le K^{-m-2}\), and \(s-1-m/2<0\) because \(m>2s-2\).
As in the endpoint argument the common buffer is annihilated
by each signed moderator, and one unit moderator pivot per row
corrects the deficit. Compactness of \(D\chi+\mathrm{Id}\) bounds the
compressed high-coordinate Jacobian error by an arbitrarily small
number after the finitely many sites have been separated from the old
ones. The Schur complement that keeps the old frequencies fixed has
the same property.

We verify nonvanishing of the divided physical coefficient for this
choice of moderator. In the alternating-path formula of
Section~\ref{sec:h2-extension} the five-letter coefficient is
\(-1/(9T^4)\), by summing its twelve orders. Indeed \eqref{eq:h2path}
is invariant under a common translation of the sites and, for two
negative letters, homogeneous of degree \(-4\) under their dilation,
so the coefficient is \(T^{-4}\) times its value at the positive
sites \(-1,3,4\) and the negative sites \(1,5\). For the six orders
\((p_1,p_2,p_3;q_1,q_2)\) equal to
\[
 (-1,3,4;1,5),\ (-1,3,4;5,1),\ (-1,4,3;1,5),\ (-1,4,3;5,1),\
 (3,-1,4;1,5),\ (3,-1,4;5,1)
\]
the products \(\prod_{\alpha=1}^2
((p_\alpha-q_\alpha)(p_{\alpha+1}-q_\alpha))^{-1}\) are
\[
 -\tfrac18,\quad\tfrac1{72},\quad-\tfrac1{12},\quad\tfrac1{36},\quad
 -\tfrac1{24},\quad-\tfrac1{72},
\]
with sum \(-2/9\). The six reversed orders give the same values, and
\((-1/2)^2(-4/9)=-1/9\). Take \(T=1\bmod\pi\).
The five signed weak sites are distinct modulo \(\pi\).
For a new positive and negative pair \(a,b\), multiply the path
coefficient by \(a-b\) and let them coalesce at \(y\) distinct from
the remaining labels. Only paths traversing that edge remain. Deleting
the adjacent pair, its insertions into an internal gap \(v,w\) contribute
\[
 \frac{v-w}{(v-y)(y-w)}
       =\frac1{v-y}+\frac1{y-w},
\]
and endpoint insertions contribute the corresponding single terms.
Both ends of every remaining vertex are counted. Including the
common factor \(-1/2\), the residue is the old coefficient times
\(\sum_{v\ {\rm positive}}(y-v)^{-1}
          -\sum_{w\ {\rm negative}}(y-w)^{-1}\),
up to the sign from ordering the colliding pair. To verify the
count, write an old path as \(p_1q_1p_2\cdots q_\mu p_{\mu+1}\), with
weight the product of \(1/(p-q)\) over its \(2\mu\) edges. The new pair
is adjacent in the order \(b,a\) inside a gap \(p_\alpha,q_\alpha\) or
after \(p_{\mu+1}\), and in the order \(a,b\) inside a gap
\(q_\alpha,p_{\alpha+1}\) or before \(p_1\). After multiplication by
\(a-b\) and coalescence the old weight is multiplied, respectively, by
\[
 \frac{p_\alpha-q_\alpha}{(p_\alpha-y)(y-q_\alpha)},\qquad
 \frac1{p_{\mu+1}-y},\qquad
 \frac{p_{\alpha+1}-q_\alpha}{(y-q_\alpha)(p_{\alpha+1}-y)},\qquad
 \frac1{p_1-y}.
\]
By the displayed partial fractions the sum of all these factors is
\(2\sum_p(p-y)^{-1}+2\sum_q(y-q)^{-1}\): every vertex of the old path
receives one term from each of its two sides, the missing side of an
end vertex being supplied by the endpoint insertion. The new
coefficient carries one more factor \(-1/2\), which gives
\(\sum_p(y-p)^{-1}-\sum_q(y-q)^{-1}\) times the old coefficient
\(\mathcal C_{\rm old}\). In formulas, with \(\mathcal C_{\rm new}\) the coefficient
after adjoining the pair,
\[
 \lim_{a,b\to y}(a-b)\,\mathcal C_{\rm new}
 =\mathcal C_{\rm old}\,L(y),\qquad
 L(y)=\sum_{v\ {\rm positive}}\frac1{y-v}
      -\sum_{w\ {\rm negative}}\frac1{y-w}.
\]
An already
coalesced opposite pair contributes zero to this last expression.
The weak signed sum has five poles and numerator of degree four
with leading coefficient one. Choose \(K\bmod\pi\) so that all
\(2^{m-1}\) residues \(K+i\) avoid its poles and zeros. There are
fewer than \(\pi\) excluded residues. Multiply the full coefficient
by the product of the paired differences and reduce modulo \(\pi\).
All other differences are units; a path uses a given paired edge
at most once. Iterating the residue identity gives a nonzero result.
In detail, number the pairs \(a_i=K+i\), \(b_i=K+i+\pi\),
\(0\le i<2^{m-1}\), let \(\mathcal C_n\) be the coefficient
\eqref{eq:h2path} of the word consisting of \(h\) and the first \(n\)
pairs, and put \(G_n=\mathcal C_n\prod_{i<n}(a_i-b_i)\). A path uses the
edge between \(a_i\) and \(b_i\) at most once, since every site
occurs once in it, and all differences other than the \(a_i-b_i\)
are units modulo \(\pi\): the weak sites are congruent to
\(0,4,5,2,6\), the residues \(K+i\) are distinct from one another
because \(2^{m-1}<\pi\), and they avoid the weak residues by the choice of
\(K\). Hence every \(G_n\) is \(\pi\)-integral. Regard \(G_{n+1}\)
as a function of the position \(b_n\) of the last site, the other
sites being fixed. It is a rational function with \(\pi\)-integral
coefficients whose denominators are units at \(b_n\equiv a_n\); the
simple pole of \(\mathcal C_{n+1}\) at \(b_n=a_n\) is cancelled by the factor
\(a_n-b_n\). Its values at \(b_n=a_n+\pi\) and at \(b_n=a_n\) are
therefore congruent modulo \(\pi\), and the latter value is
\(G_nL_n(a_n)\) up to sign by the residue identity, where \(L_n\) is
the signed sum for \(h\) and the first \(n\) pairs. Each earlier pair
contributes
\(\pm(a_i-b_i)/((a_n-a_i)(a_n-b_i))\equiv0\) to \(L_n(a_n)\). Thus
\[
 G_{n+1}\equiv\pm G_nL_h(K+n)\pmod\pi,\qquad
 G_{2^{m-1}}\equiv\pm\mathcal C_0\prod_{i=0}^{2^{m-1}-1}L_h(K+i)\pmod\pi,
\]
where \(L_h\) is the signed sum over the five weak sites. Here
\(\mathcal C_0=-1/(9T^4)\equiv-1/9\) is a unit, and each \(L_h(K+i)\) is a finite
nonzero residue by the choice of \(K\).
This proves physical nonvanishing for the required disjoint supports.

For clarity, here the frozen renormalized angular frequencies are
\(\tau c_n\) with real algebraic \(c_n\), and we require the old
list to be rationally independent. Start with two nearby negative
independent algebraic slopes and a sufficiently small \(\tau>0\).
Keep every previously occupied slope fixed. On each new flat
\((h_e+d_e)\cdot c=0\) one may use a private unit pivot, independently
for each group \(e\). The passive calibration and the continuation
along the common scaling parameter in Section~\ref{sec:h2-extension}
apply without the congruence condition on the slopes. More explicitly,
while all weak actions are zero use a moderator pivot to impose each
row. Approximate the free new rates by \(\tau\) times algebraic slopes
within the open set where the newly occupied moderator actions are positive. The normalized coefficient with the
specified signed row and zero other angular labels has the computed
nonzero vacuum value on this analytic recipe. Division by its indicated
radii extends analytically across zero by expansion in the complex
Cartesian variables and angular averaging. The bounded change of
coordinates in the slope inverse replaces each pivot by its row;
the finite row lengths here depend only on \(s\), with disjoint
groups. The small action ball of the inverse and the complex
Fourier--Lebesgue chart used in the preceding section therefore
connect this recipe to zero. Take the original \(\tau\) outside the
exceptional sets of the countably many finite algebraic recipes,
including the earlier slopes. Activate the five weak radii in each
group by the fifth-root prescription at the actual divided
coefficient. Their cost at the fixed sites tends to zero with the
common target modulus \(\lambda\). Retune the pivots and eliminate one
weak pivot of coefficient \(\pm1\) in each row on the activated flat.
Choose the remaining new slopes successively inside the open
free-coordinate cube that preserves positivity and the relative
coefficient fit. At each step let \(\mathbb K\) be the number field
generated by the old slopes and the new free slopes already chosen.
For any real algebraic \(\alpha\notin\mathbb K\), the dense set
\(\alpha+\mathbb Q\) is disjoint from \(\mathbb K\). Choosing the next
slope in this set preserves rational independence of the free list.
The eliminated slopes are algebraic by the row equations. Subtracting
integer multiples of the rows \(h_e+d_e\) clears the private pivot
entries of any integral relation. The remaining relation vanishes by
independence, so these rows generate the full relation lattice.
After the physical tilt, use the same successive choice to round
the weak pivots and the new terminal root slopes off the flat,
keeping the other slopes fixed as specified below.

For every one of the countably many algebraic lists so obtained,
the nonzero quantities \(b+\tau l\cdot c\), \(b\in\mathbb Z\),
have a polynomial lower bound in \(1+|b|+|l|\) at this same
\(\tau\). For \(b=0\) use the algebraic norm in the finite number
field of the list. For \(b\ne0\), restrict \(\tau\) first to a
compact positive interval. An affine expression that can approach
zero there has \(|l\cdot c|\) bounded below, so the lengths of its
exceptional intervals are summable on taking a sufficiently large
fixed exponent. Use Borel--Cantelli and exclude exact zeros. A
frozen integer-time angular label is consequently precisely a sum
of the private rows. Its charge, space frequency and integer time
all equal the sum of its row coefficients. Indeed an integer-time
label \(l\) has \(\tau l\cdot c\in\mathbb Z\), hence \(l\cdot c=0\)
by the exclusion of exact zeros of the nonzero expressions
\(b+\tau l\cdot c\); the relation lattice of the frozen list is
generated by the rows, and each row \(h_e+d_e\) has the triple
\((1,1,1)+(0,0,0)\). For a charge-one term
the order-five weak characters in this lattice are the positive
singletons: a label \(\sum_en_e(h_e+d_e)\) has charge \(\sum_en_e\)
and exactly \(5\sum_e|n_e|\) weak letters, because the weak supports
are disjoint, so charge one and five weak letters force a single
\(n_e=1\). All others have more weak letters and their absolute
coefficient sum is \(O(\lambda^{6/5})\) by analyticity on a larger
weak polydisc at the fixed passive preparation. Each singleton has
the required modulus \(\lambda(1+o(1))\). This construction separates
old phase harmonics without assuming statistical independence.

\smallskip
\noindent\emph{The Diophantine bound in detail.} Fix one algebraic list
\(c=(c_1,\ldots,c_k)\) and a number field \(\mathbb F\supset\mathbb Q(c)\)
of degree \(\mathsf d\). Choose \(q\in\mathbb N\) such that the numbers
\(\alpha_i=qc_i\) are algebraic integers, and put
\(A_c=\max_{i,\iota}|\iota(\alpha_i)|\) over the embeddings
\(\iota:\mathbb F\to\mathbb C\). If \(l\in\mathbb Z^k\) and \(l\cdot c\ne0\),
then \(\beta=l\cdot\alpha\) is a nonzero algebraic integer, so
\(|N_{\mathbb F/\mathbb Q}(\beta)|=\prod_\iota|\iota(\beta)|\ge1\), while
\(|\iota(\beta)|\le A_c|l|_1\) for every embedding. Hence
\[
 |\tau\,l\cdot c|=\frac\tau q\,|\beta|\ge\frac\tau q\,(A_c|l|_1)^{1-\mathsf d},
\]
which is the case \(b=0\). For \(b\ne0\) restrict \(\tau\) to a compact
interval \([\tau_1,\tau_2]\subset(0,\infty)\) and fix \(\gamma>k+1\). For
\((b,l)\) with \(l\cdot c\ne0\) put
\(\varepsilon_{b,l}=(1+|b|+|l|)^{-\gamma}\le1/2\) and
\[
 E_{b,l}=\{\tau\in[\tau_1,\tau_2]:|b+\tau\,l\cdot c|<\varepsilon_{b,l}\}.
\]
This is an interval of length at most \(2\varepsilon_{b,l}/|l\cdot c|\).
If it is nonempty, then \(|b|<\tau_2|l\cdot c|+\tfrac12\), so that
\(|l\cdot c|\ge|b|/(2\tau_2)\ge1/(2\tau_2)\) and its length is at most
\(4\tau_2\varepsilon_{b,l}\). If \(l\cdot c=0\), then
\(|b+\tau\,l\cdot c|=|b|\ge1\). Since
\(\sum_{b,l}(1+|b|+|l|)^{-\gamma}<\infty\) over \(\mathbb Z^{1+k}\), the
Borel--Cantelli lemma shows that almost every \(\tau\in[\tau_1,\tau_2]\)
lies in only finitely many \(E_{b,l}\). Remove in addition the countably
many exact zeros \(\tau=-b/(l\cdot c)\). Every remaining \(\tau\) satisfies
\(|b+\tau\,l\cdot c|\ge\kappa_\tau(1+|b|+|l|)^{-\gamma}\) for all \(b\ne0\)
and all \(l\), where \(\kappa_\tau>0\) is the minimum of \(1\) and of the
finitely many exceptional ratios. The algebraic lists are countable, so
a common \(\tau\) exists outside a null set, and it may in addition be
taken outside the countably many exceptional sets of the continuation.

\subsection{Choice of directions and matching of the actual frequencies}

On a normalized real tree with main root \((1,0)\) put
\[
 \Xi(v)=A/v_x^2-B/|v|^2,\qquad
 D(p,q)=\Xi(p)+\Xi(q)-\Xi(p+q).
\]
At the frozen positive background the functions \(1,D_e\) of all
the real split parameters, for primary and private edges, are linearly
independent over the constant reals. At the last graft, with parent
\((X,Y)\), \(XY\ne0\), the private defect for \(p\) has a fourth
order pole at zero where the other private defect and primary have
poles of order at most two. At infinity interchange the children.
The remaining primary has a pole, because \(A>0\).
Remove the three coefficients in this order and proceed backwards.
In detail, suppose \(c_0+\sum_ec_eD_e=0\) identically. The parameter
\(t\) of the last graft enters only its own three defects, whose
expansions were displayed in Section~\ref{sec:h2-extension}:
\(D(p,J(p))\) is the only one with a pole of order four at \(t=0\),
and \(D(q,J(q))\) the only one with a pole of order four at
\(t=\infty\), so their coefficients vanish. The primary defect
\(D(p,q)=A/(t^2Y^2)+O(|t|^{-1})\) is then the only remaining term
with a pole at \(t=0\), and its coefficient vanishes as well. All
defects of the last graft have now been removed, and the parameter of
the preceding graft enters only its own three defects among the
remaining terms.
The initial horizontal split has nonconstant defect. Exclude also
the equality of any mixed defect with a desired singleton, using
the horizontal pole of the mixed output in the endpoint construction.
On a compact box of admissible split angles, multiplication by their
common rational-function denominator makes the map from coefficient
vectors to polynomial numerators injective, hence bounded below.
The polynomial sublevel bound on a box,
\(\lvert\{|p|\le h\}\rvert\le C(h/\|p\|_\infty)^\alpha\)
for polynomials of bounded degree, follows by iterating the
one-variable sublevel inequality (split according to the size of a
nonzero coefficient polynomial). Summing over integer coefficients
and applying Borel--Cantelli supplies real parameters with
\[
 |a\cdot D|\ge c(1+|a|)^{-C}\qquad(0\ne a\in\mathbb Z^{\#\{e\}}).
\]
All constants may depend on this finite preparation.

\smallskip
\noindent\emph{The sublevel estimate in detail.} Let
\(\mathsf t\in\mathsf Q\) denote the split parameters on the compact
box, \(\mathsf k\) their number, and let \(\Delta(\mathsf t)\) be a
polynomial common denominator of the defects with
\(|\Delta|\ge\delta_0>0\) on \(\mathsf Q\). Then the \(P_e=\Delta D_e\) are
polynomials of degree at most \(\mathsf m\). By the linear independence
just proved, \(a\mapsto a\cdot P\) is an injective linear map from
\(\mathbb R^{\#\{e\}}\) to the finite-dimensional space of polynomials of
degree at most \(\mathsf m\), hence
\(\|a\cdot P\|_{L^\infty(\mathsf Q)}\ge c_0|a|\). In one variable the
Remez inequality on an interval \(I\) states
\(\sup_I|p|\le(4|I|/|E|)^{\mathsf m}\sup_E|p|\) for every measurable
\(E\subset I\) of positive measure and every polynomial \(p\) of degree at
most \(\mathsf m\); thus
\(|\{t\in I:|p(t)|\le h\}|\le4|I|\,(h/\|p\|_{L^\infty(I)})^{1/\mathsf m}\).
In \(\mathsf k\) variables write
\(p=\sum_{i\le\mathsf m}p_i(\mathsf t')\,t_{\mathsf k}^i\). By the
equivalence of norms on polynomials of bounded degree there is an \(i\)
with \(\|p_i\|_{L^\infty}\ge c\|p\|_{L^\infty}\), and a polynomial in one
variable whose \(i\)th coefficient has modulus \(\varrho\) has supremum at
least \(c\varrho\) on the interval. Splitting the \(\mathsf t'\) according
to whether \(|p_i(\mathsf t')|\ge\varrho\|p\|_\infty\), and using induction
on the dimension for the other set, we get
\[
 |\{|p|\le h\}|\le C\bigl(h/(\varrho\|p\|_\infty)\bigr)^{1/\mathsf m}
   +C\varrho^{\alpha_{\mathsf k-1}} .
\]
The choice
\(\varrho=(h/\|p\|_\infty)^{1/(1+\mathsf m\alpha_{\mathsf k-1})}\) gives the
stated bound with \(\alpha_{\mathsf k}=\alpha_{\mathsf k-1}/
(1+\mathsf m\alpha_{\mathsf k-1})\) and \(\alpha_1=1/\mathsf m\). Now let the
exponent \(C\) of the display satisfy \(C\alpha_{\mathsf k}>\#\{e\}\), and
put \(B_a=\{\mathsf t\in\mathsf Q:|a\cdot D(\mathsf t)|\le(1+|a|)^{-C}\}\).
As \(|a\cdot P|\le\|\Delta\|_\infty|a\cdot D|\),
\[
 |B_a|\le C'\Bigl(\frac{\|\Delta\|_\infty(1+|a|)^{-C}}{c_0|a|}\Bigr)^{\alpha_{\mathsf k}}
 \le C''(1+|a|)^{-C\alpha_{\mathsf k}},
\]
which is summable over \(a\in\mathbb Z^{\#\{e\}}\setminus\{0\}\). By
Borel--Cantelli almost every \(\mathsf t\) lies in only finitely many
\(B_a\). The set where \(a\cdot D=0\) for some \(a\ne0\) is a countable
union of proper algebraic subsets, hence null. For the remaining
\(\mathsf t\) the displayed bound holds with \(c\) the minimum of \(1\) and
of finitely many positive numbers. These parameters have full measure,
so they meet the nonempty open set defined by the strict exclusions.

\begin{stagelemma}[Uniform rational approximation of the tree]\label{lem:rational-tree-uniformity}
Fix \(s>2\), the finite depth, the frozen preparation and a real tree as
above with \(|a\cdot D|\ge c(1+|a|)^{-\kappa}\) for every nonzero
integer coefficient vector \(a\). Given \(C_1>0\) and \(b_0>\kappa+1\),
for every sufficiently large multiple \(R\) of the old vertical spacing
there is a rational tree satisfying the congruence exclusions of
Proposition~\ref{prop:tree}. Its split parameters differ from the
chosen parameters by \(O((\log R)^{-b_0})\); its common spatial
denominator \(P\) is at most \(C(1+\log R)^{C_0}\). With \(N=RP\)
and rational-tree defects \(D(N)\), one has
\[
 |a\cdot D(N)|\ge \tfrac c2(1+|a|)^{-\kappa}
 \qquad (0<|a|\le C_1\log N).
\]
The nonzero horizontal components and the strict geometric exclusions
on the prescribed finite list retain positive margins in normalized
coordinates. The constants and the threshold for \(R\) may depend
on the fixed stage, \(C_1\) and \(b_0\), and are uniform in \(R\).
\end{stagelemma}

\begin{proof}
Use a single rational tree from Proposition~\ref{prop:tree} to witness
the finitely many odd valuations. Write its split parameters as integer
pairs. At each witness prime the valuation of each tested nonzero
rational function is locally constant at this tree. Hence fixing the
integer pairs modulo one product \(M\) of sufficiently high prime powers
preserves all the required valuations. For each desired real parameter
\(t\), choose a positive denominator \(b\asymp(\log R)^{b_0}\) in its
prescribed class and a numerator in its prescribed class nearest to
\(bt\). The error is at most \(M/b\). At fixed depth the vertex
coordinates are fixed rational functions of these parameters. Clearing
their denominators costs at most a fixed power of \(1+\log R\).
On a compact neighborhood of the real tree all relevant real denominators
stay away from zero, so the defect functions are Lipschitz there and
the finite strict exclusions persist. In particular
\(\|D(N)-D\|\le C(\log R)^{-b_0}\). Since
\(\log N=\log R+O(\log\log R)\), the error in \(a\cdot D\) is at most
\(C|a|(\log R)^{-b_0}\le(c/2)(1+|a|)^{-\kappa}\) on the stated
range for sufficiently large \(R\). Indeed \(|a|\le2C_1\log R\) there,
so
\[
 C|a|(1+|a|)^{\kappa}(\log R)^{-b_0}
 \le C(1+2C_1\log R)^{\kappa+1}(\log R)^{-b_0}\longrightarrow0,
\]
because \(b_0>\kappa+1\).
For completeness, let \(F=G/H\) be one of the tested rational functions of
the integer pairs \((a_i,b_i)\), with \(G,H\) integer polynomials, and let
\((a^0,b^0)\) be the integer pairs of the witness tree, at which
\(G,H\ne0\). For a witness prime \(\pi\) put
\(E_\pi=1+\max\bigl(v_\pi(G(a^0,b^0)),v_\pi(H(a^0,b^0))\bigr)\). If
\((a,b)\equiv(a^0,b^0)\pmod{\pi^{E_\pi}}\) componentwise, then
\(G(a,b)\equiv G(a^0,b^0)\) and \(H(a,b)\equiv H(a^0,b^0)\) modulo \(\pi^{E_\pi}\),
hence \(v_\pi(F(a,b))=v_\pi(F(a^0,b^0))\). Taking \(M\) as the product of the
\(\pi^{E_\pi}\) over the finitely many witness primes and tested functions,
together with the modulus of the congruences of
Proposition~\ref{prop:tree}, all required valuations and congruences are
preserved on the prescribed classes modulo \(M\). In each class the
fractions \(a/b\) with \(b\asymp(\log R)^{b_0}\) fixed have spacing \(M/b\),
which gives the stated approximation.
\end{proof}

Put \(r=N^{-s}\); the physical root is at \(K_*=1+N\), and the slow
time is \(\theta=r\lambda t\). Below \(\mathcal{L}\) denotes a fixed
power of \(1+\log N\), allowed to increase finitely many times.

Here are the changes to the normal reduction. Choose an integer
\(J>s+4\). For all \(m'\in R\mathbb Z\), \(m'\ne0\), and
\(|n|\ge cN/P^2\), expand the doubled columns to spatial order
greater than \(J+3\), before making any temporal division.
With \(f=e^w,\ \langle w\rangle_x=0,\ z=g/f\), the upper equation
determines the next mean-zero \(Dw\) by \(2n\), the next scalar by
its spatial mean; the lower equation determines the next coefficient
of \(z\) by \(2E\), \(E=n^2+(m')^2\). Equivalently use the
undivided envelope recursion of Lemma~\ref{lem:empty-frequency}
and normalize the spatial mean of \(f\) to one. The coefficients are
smooth old fields times bounded polynomials in \(n^2/E\). Every
lower coefficient has a factor \(n^2/E\), so the residual estimate
remains uniform even when \(|m'|/|n|\) is unbounded. Increasing the
finite order absorbs all factors \(P=O(\mathcal{L})\). Before normalizing a positive column \((f,g)\), put
\(S=\|f\|_{L^2_x}^2-\|g\|_{L^2_x}^2>0\). Its equation with instantaneous
scalar \(d\) gives \(\partial_t S=-2\operatorname{Im}(d)S+O(N^{-J-2}\mathcal{L})\).
Dividing the column by \(\sqrt S\) gives unit Krein square and changes
the scalar to \(d+i(\partial_t S)/(2S)\), real to the residual accuracy.
To verify these two assertions, write the column equations, up to
the residuals, as
\[
 i\partial_tf=(\mathcal A-d)f+Cg,\qquad
 -i\partial_tg=(\mathcal B+d)g+\bar Cf,
\]
in the mass gauge of the isolated torus. Here \(V=e^{2iMt}U\),
\(V_0=2(|V|^2-M)\), \(C=V^2\), and the selfadjoint spatial operators are
\[
 \mathcal A=2nD+D^2+V_0,\qquad \mathcal B=2E+2nD+D^2+V_0 .
\]
Indeed, at fixed actions the gauge changes the common carrier from
\(E+2M\) to \(E\). The upper envelope is unchanged and the lower
envelope of Lemma~\ref{lem:empty-frequency} becomes
\(g=e^{-4iMt}g_{\rm ung}\). Thus
\(i\partial_t g=e^{-4iMt}(i\partial_t+4M)g_{\rm ung}\),
so the two recursions give the same instantaneous scalar \(d\) at
each point of the fixed-action torus. Taking real parts of the
\(L^2_x\) pairings with \(f\) and \(g\),
\[
 \tfrac12\partial_t\|f\|^2=-\operatorname{Im}(d)\|f\|^2
        +\operatorname{Im}\langle f,Cg\rangle,\qquad
 \tfrac12\partial_t\|g\|^2=-\operatorname{Im}(d)\|g\|^2
        +\operatorname{Im}\langle f,Cg\rangle,
\]
and the difference is the stated equation for \(S\). For
\((\tilde f,\tilde g)=S^{-1/2}(f,g)\) the additional term
\(-(\partial_tS/2S)(\tilde f,\tilde g)\) in \(\partial_t\) is the
same as replacing \(d\) by \(\tilde d=d+i(\partial_tS)/(2S)\) in both
equations, and
\(\operatorname{Im}\tilde d=\operatorname{Im}d+\partial_tS/(2S)
=O(N^{-J-2}\mathcal L)\). Its averaged leading term is \(A/n^2-B/E\), and its
nonconstant part at that order is a bounded time derivative by
Lemma~\ref{lem:empty-frequency}.

Fit the carriers before taking the temporal primitives. At each high
off-axis head use its bare energy plus the angle mean of its full
instantaneous spatial scalar. At an axial root take its exact frequency
on the terminal torus, whose new actions are prescribed as \(r^2\)
times the desired squared radii. The off-axis action derivative of the
averaged shift is \(O(N^{-2}\mathcal{L})\), since this mean has no angular
divisor. The high-root frequency derivative with respect to the old
actions, after mass subtraction, is \(o(1)\) by compactness of
\(D\chi+\mathrm{Id}\). These bounds are uniform on the shrinking
neighborhood of the frozen point. The compressed frequency inverse
with all previously occupied rates, all moderator rates, and the
nonpivot weak rates held at their frozen values allows a change of
one private weak rate per row. These are independent coordinates:
prescribe the new terminal root actions first, then solve the
compressed positive-action chart for the entire remaining rate
vector. Write \(y_e\) for the signed rate displacement of row \(e\)
from the frozen value and \(I(y)\) for the actions given by this
inverse. For a head \(v\), let \(\Delta_v(I)\) be the correction
to its bare energy, evaluated as above. Then the equations for all
edges are the single fixed point equation
\[
 y_e=\Delta_v(I(y))+\Delta_w(I(y))-\Delta_{v+w}(I(y)),
 \qquad e=(v,w;v+w).
\]
At \(y=0\) the right side is \(O(N^{-2}\mathcal L)\).
Each head correction has derivative \(o(1)\) in the varying actions
by the estimates above. The compressed inverse and the number of
edges are bounded at this fixed preparation, so the derivative of
the full right side with respect to \(y\) is \(o(1)\).
It is therefore a contraction on a ball of radius
\(C N^{-2}\mathcal L\) for a sufficiently large \(C\).
At its fixed point each row changes by
\(D_e(N)/N^2+O(N^{-3}\mathcal{L})\), as verified below.
For the next stage keep the previously occupied, moderator and
nonpivot weak rates fixed in the final rounding as well. Those fixed
rates are already rationally independent together: dropping one
private weak pivot per row removes all relations. In the full
compressed rate chart round the weak pivots and the new terminal
root rates to independent algebraic slopes times \(\tau\), to
precision smaller than both \(N^{-J-3}\) and \(o(r^2)\).
This preserves the positive root radii up to relative \(o(1)\). The finite tree
equations can be made exact again by changing the auxiliary and
mobile carrier speeds by \(O(N^{-J-2})\), starting from the specified
axial speeds and splitting errors between children.
To see the size of the row displacement, note that at every action
vector \(I\) of the shrinking neighborhood the high-column expansion
gives, for an off-axis head \(k_v=(1,0)+Nv\),
\[
 \Delta_v(I)=\frac{A(I)}{(1+Nv_x)^2}-\frac{B(I)}{|k_v|^2}
   +O(N^{-3}\mathcal L)=\frac{\Xi_I(v)}{N^2}+O(N^{-3}\mathcal L),
\]
where \(\Xi_I\) is computed with the coefficients \(A(I),B(I)\) and the
horizontal components have positive margins. For an axial root
\(K=1+Np\), let \(I_{\rm lo}\) be \(I\) with the actions at all the
current hot roots set equal to zero. Their prescribed amplitudes are
bounded at fixed tree, so \(\|I-I_{\rm lo}\|_{\ell^1}=O(r^2)\).
The local \(\ell^2\) Lipschitz bound for \(\chi\) gives
\(\lvert\chi_K(I)-\chi_K(I_{\rm lo})\rvert=O(r^2)\).
The supports of the low backgrounds lie in a fixed finite set and
\(K\) is empty there. Applying \eqref{eq:sup-empty} on this family gives
\[
 \Delta_v(I)=\chi_K(I)
 =\frac{A(I_{\rm lo})-B(I_{\rm lo})}{K^2}
   +O(|K|^{-3}+r^2).
\]
The definitions in Lemma~\ref{lem:empty-frequency} give
\(d\mathcal Q=2\sum_j I_j\,d\omega_j
=4M\,dM+2\sum_j I_j\,d\chi_j\).
Thus \(A=\mathcal Q-M^2\) and \(B=\mathcal Q/2\) change by
\(O(r^2)\) when the hot actions are restored, uniformly in their indices.
Since \(r^2=N^{-2s}=o(N^{-3})\) and \(|p|\) is bounded away from zero at this fixed tree, this yields
\(\Delta_v(I)=\Xi_I((p,0))/N^2+O(N^{-3}\mathcal L)\). The fit
uses \(\chi_K(I)\) at the actual actions. Along the fixed low support
the coefficients \(A,B\) are smooth, and \(\|I(y)-I(0)\|=O(|y|)=O(N^{-2}\mathcal L)\) at the fixed point. Hence
\[
 y_e=\frac{\Xi(v)+\Xi(w)-\Xi(v+w)}{N^2}+O(N^{-3}\mathcal L)
   =\frac{D_e(N)}{N^2}+O(N^{-3}\mathcal L),\qquad e=(v,w;v+w),
\]
with \(\Xi\) evaluated at the frozen actions. Only the expansion at each
fixed action vector is used, not its derivative in the actions.

Let \(F\) be the support of the axial background used for the spatial
columns, including the moderators and the activated weak sites. The
hot roots enter separately through their carriers. Write
\(\Omega^\infty\) here for the terminal angular frequencies restricted
to \(F\), and keep these divisors fixed along each trial with moving
old actions and phases. For every nonzero \(l\in\mathbb Z^F\) of
charge zero and spatial momentum zero with \(|l|\le C_1\log N\), we have
\[
                      |l\cdot\Omega^\infty|\ge N^{-2}\mathcal{L}^{-1}.
\]
If its frozen divisor is nonzero use the polynomial bound and its
\(O(N^{-2}\mathcal{L})\) displacement. Otherwise it is a nontrivial
neutral combination of the rows, and use the bound for \(a\cdot D(N)\);
the private pivots bound \(|a|\) by \(|l|\). Higher terms
\(O(N^{-3}\mathcal{L})\) preserve the separation.
Let \(\Omega^0\) denote the frozen frequency vector on the same support
\(F\). There \(|\Omega^\infty-\Omega^0|=O(N^{-2}\mathcal L)\).
If \(l\cdot\Omega^0\ne0\), the polynomial lower bound
at the frozen list gives
\[
 |l\cdot\Omega^\infty|\ge c(1+|l|)^{-C}-C|l|N^{-2}\mathcal L
 \ge\tfrac c2(1+C_1\log N)^{-C}\ge N^{-2}\mathcal L^{-1}.
\]
If \(l\cdot\Omega^0=0\), then \(l=\sum_ea_e(h_e+d_e)\) with integers
\(a_e\), not all zero, and \(|a_e|\le|l|\) by reading \(l\) at the
private pivot of row \(e\). Each row was displaced by
\(D_e(N)/N^2+O(N^{-3}\mathcal L)\). Hence
\[
 l\cdot\Omega^\infty=\frac{a\cdot D(N)}{N^2}+O\bigl(|a|N^{-3}\mathcal L\bigr),
 \qquad
 |l\cdot\Omega^\infty|\ge\frac{c}{2N^2}(1+C\log N)^{-\kappa}
       -\frac{C\mathcal L\log N}{N^3}\ge N^{-2}\mathcal L^{-1}
\]
by Lemma~\ref{lem:rational-tree-uniformity}, applied with a fixed
multiple of \(C_1\) in its range \(|a|\le C\log N\), and after
enlarging the power \(\mathcal L\).
For orders \(k\ge3\), divide each nonconstant spatial mean at order
\(n^{-k}\) by \(il\cdot\Omega^\infty\) through this cutoff.
The real phase so obtained and all needed old-parameter derivatives
are \(O(N^2|n|^{-3}\mathcal{L})=O(N^{-1}\mathcal{L})\); the last
equality uses \(|n|\ge cN/P^2\) and \(P=O(\mathcal L)\), so that
\(N^2|n|^{-3}\le c^{-3}P^6N^{-1}\). The tail after
the cutoff contributes \(O(N^{-J-2})\) by angular analyticity, on
increasing \(C_1\). At order two use its bounded primitive instead.
The last conjugation is a scalar phase on the individual input
columns. It commutes with the bare diagonal; its moving connection
costs \(O(N^{-1}\mathcal{L})(r|J_o|+r^2)\), where \(rJ_o\) is the old
action displacement. Indeed the old isolated rates and the full
terminal old rates differ by \(O(r|J_o|+r^2)\) in this tube. The
time derivative in the spatial recursion uses the actual isolated
old flow at those actions. Its difference from the trial velocity
multiplies old-parameter derivatives of the near-identity column.
All denominators in this final scalar operation stay fixed.
In particular, a temporal divisor of size \(N^{-2}\mathcal{L}^{-1}\)
is never fed back into a spatial recursion.

On each doubled pair of nonzero vertical copies, first remove the
old off-diagonal pairing to arbitrary finite order as in the proof
of Proposition~\ref{prop:insertion}. In detail every off-diagonal
kernel is divided by \(n^2+(n')^2+2(m')^2\ge2R^2\) with the
opposite sign on the other copy. Use its Krein-preserving
exponential. The commutator with the free generator is exactly
the undivided kernel; successive commutators with the bounded old
block and the connection lose at most a fixed old norm and gain
\(R^{-2}\). All estimates hold with summable off-diagonal moments.
More precisely, repeat the two steps displayed in
Section~\ref{sec:h2-extension}. If the block off-diagonal part of the
generator after \(k\) steps is \(Y_k\), of Schur size \(O(R^{-2k})\)
together with its old-parameter derivatives, let \(J_{k+1}\) be the
block off-diagonal matrix with \([D_{\rm free},J_{k+1}]=-Y_k\). Its
entries are those of \(Y_k\) divided by
\(\pm(n^2+(n')^2+2(m')^2)\), so \(J_{k+1}=O(R^{-2k-2})\). Conjugation by
\(e^{J_{k+1}}\) removes \(Y_k\); the new off-diagonal part consists of
the commutator of \(J_{k+1}\) with the bounded block-diagonal part, of
\(-i\dot J_{k+1}\), and of terms of higher order, so that
\(Y_{k+1}=O(R^{-2k-2})\). The block-diagonal part changes by
\(O(R^{-2k}\cdot R^{-2k-2})\) at each later step and is determined to
any prescribed order by finitely many steps. Since \(R=N/P\) and
\(P=O(\mathcal L)\), taking \(2k>J+2\) makes the off-diagonal error
\(O(N^{-J-2}\mathcal L)\), which after division by
\(r\lambda=N^{-s}\lambda\) is \(O(N^{s-J-2}\mathcal L)=O(N^{-6}\mathcal L)\).
Transform the direct high columns above by this frame and project
onto their Hermitian copies. Their residuals keep the stated order.
For their matrix \(T\), the covariant commutator of \(T^*T\) with
the real high diagonal has the same order by the equation and its
adjoint. The binomial series thus gives the same order for the
commutator of \((T^*T)^{-1/2}\). Complete the high isometry
\(T(T^*T)^{-1/2}\) on its low complement exactly as in
Proposition~\ref{prop:insertion}. The kinetic commutator of an
\(O(1/|n|)\) spatial envelope is bounded with Schur moments, and
the scalar temporal phase is diagonal on the input. This leaves
a full Hermitian low block at the actual trial old angles and
actions. It is kept in the principal evolution. Freezing that
low block at the terminal phases would introduce an uncontrolled
\(1/(r\lambda)\) angle derivative.

\subsection{Finite dynamics and the normal time weights}

Use the complete finite equation on the tree and its mixed
hypotenuses. The assigned singleton gives a constant coefficient
on each intended triad. Singletons assigned to any other defect
there, or to a mixed edge, rotate at slow speeds bounded below
in magnitude by \(c N^{s-2}\) at this fixed preparation. Indeed the
physical detuning of the singleton of row \(e'\) on a triad with
defect \(D_e(N)\) is \((D_e(N)-D_{e'}(N))/N^2+O(N^{-3}\mathcal L)\),
and in the slow time \(\theta=r\lambda t\) with \(r=N^{-s}\) it
becomes \((D_e(N)-D_{e'}(N))N^{s-2}/\lambda\) up to a relative error
\(O(N^{-1}\mathcal L)\); the defects are separated by
Lemma~\ref{lem:rational-tree-uniformity} with \(a\) a difference of
two unit vectors, and from the mixed defects by the strict
exclusions.
Their bounded oscillatory primitives remove them by a near-identity
finite change; its remainders are \(o_R(1)\) in \(C^1\) on
bounded finite sets, by integration by parts with the original
bounded vector field. In particular only its first derivative along
an actual trajectory is used. Higher weak terms have normalized
size \(O(\lambda^{1/5})\). Both finite fields vanish on the
axial terminal plane in these rotating coordinates. The common
mixed inverse on the exponentially landing finite orbit from the
endpoint therefore yields the exact finite orbit with the rounded
axial radii and an arbitrarily small launch error: first fix the
finite tolerance, take \(\lambda\) small, then \(R\) large.
Equivalently translate the reference by the axial rounding; the
forcing at the new constant endpoint is zero and its difference
along the exponentially landing reference decays at the same weight.
The early interval has fixed slow length.
The argument is the averaging and landing argument written out at the
end of the finite dynamics in Section~\ref{sec:h2-extension}, with two
changes. First, the oscillating singletons have slow frequencies
\(|\omega_k|\ge cN^{s-2}/\lambda\), and there are finitely many of them
at the fixed preparation; the relative error \(O(N^{-1}\mathcal L)\) of
these frequencies does not affect the lower bound. The change of
variables \(\Phi(Y,\theta)=\sum_kG_k(Y)e^{i\omega_k\theta}/(i\omega_k)\) therefore
satisfies \(\sup_\theta\|\Phi(\cdot,\theta)\|_{C^2_Y}=O(\lambda N^{2-s})\). As in
Section~\ref{sec:h2-extension}, set
\(\varepsilon_{\rm fit}=\max_e|c_e^{\rm fit}-c_e^{\rm id}|\), where
the fitted and ideal coefficients on each intended edge \(e\) are
measured after the phase rotations in slow time.
The singleton modulus \(\lambda(1+o(1))\) and the terminal rounding give
\(\lim_{\lambda\downarrow0}\limsup_{R\to\infty}\varepsilon_{\rm fit}=0\).
Each constant edge discrepancy multiplies a quadratic field vanishing
on the axial plane. Including these discrepancies gives
\[
 \|S\|_{C^1}\le
 C\bigl(\lambda N^{2-s}+\lambda^{1/5}+\varepsilon_{\rm fit}\bigr)+o_R(1).
\]
For a prescribed landing tolerance choose \(\lambda\) small and then
\(R\) large in this bound.
Second, the endpoint is the exact rounded axial point \(p_*\), at
distance \(o_R(1)\) from the endpoint \(p_{\rm id}\) of the reference
\(Y_{\rm id}\). Either use the continuity of the weighted fixed point in
\(p_*\) proved there, or translate the reference: \(Y_{\rm id}+p_*-p_{\rm id}\)
tends to \(p_*\) and solves the finite equation up to the forcing
\(F_0(Y_{\rm id}+p_*-p_{\rm id})-F_0(Y_{\rm id})\). Since \(F_0\) vanishes
on the linear axial plane, its derivative in axial directions vanishes
there, and this forcing is
\(O(|p_*-p_{\rm id}|\,\operatorname{dist}(Y_{\rm id},\mathcal A))
=O(|p_*-p_{\rm id}|\,|Y_{\rm id}-p_{\rm id}|)\); it decays at the landing
rate. In both forms the weighted contraction gives the exact finite
orbit landing at \(p_*\), and Gronwall on the early interval of fixed
slow length controls the launch.

There is one extra issue for the full terminal linear equation.
For each old harmonic \((b,\nu)\) in the angular cutoff and each axial
root \(K_i\), first retain the same-copy entry \(b\to K_i\) and its
Hermitian adjoint with their actual oscillatory coefficients, including
all off-key and noninteger-time harmonics. Indeed the bare divisor
\(2(K_i-b)(n-b)+\nu-b^2\) reduces to \(\nu-b^2\) at \(n=b\), while
for \(n\ne b\) its modulus is \(\gtrsim R\). Divide only the remaining
nonresonant entries at their frozen divisors. Every charged integer-time old harmonic has
space and
bare time \((1,1)\), since it is an integral combination of the
private rows with charge one. Thus an anomalous bare pair from an
axial root \(K_i=1+RP_i\) has
\[
 (n-1)+(n'-1)=RP_i,\qquad (n-1)(n'-1)=R^2q^2,\qquad q\ne0 .
\]
Indeed the output index of the pairing term is \(n=K_i+1-n'\), which
is the first equation, and the bare divisor
\(n^2+n'^2+2R^2q^2-K_i^2-1\) of Section~\ref{sec:h2-extension}
vanishes precisely when
\[
 2R^2q^2=(n+n'-1)^2+1-n^2-n'^2=2(n-1)(n'-1).
\]
If present its legs are bounded below in absolute value by
\(cN/P^2\): with \(N=RP\), \(x=(n-1)/N\) and \(y=(n'-1)/N\) one has
\(x+y=P_i/P\) and \(xy=q^2/P^2\ge P^{-2}\), so \(x,y\) have the same
sign, \(|x|,|y|\le|P_i|/P\le C\) and \(|x|=xy/|y|\ge c/P^2\).
Its root is unique from the sum and product. An
intersection with a Hermitian star would give
\(P_l(P_i-P_l)=q^2\), excluded by the local nonsquare conditions:
if a leg is the upper end \(K_l=1+RP_l\) of a star, then
\(n-1=RP_l\) and \(n'-1=R(P_i-P_l)\), and the product equation
becomes the displayed relation.
Entries touching a low endpoint outside the retained Hermitian
stars have anomalous bare divisor \(\gtrsim R^2\) or same-copy bare
divisor \(\gtrsim R\): for example \(2(qR)^2\) dominates the
\(O(|K_i|(|n|+|b|))\) term at \(|n|<cN/P^2\), with \(c\) small
and \(|b|\lesssim\log N\). Where both endpoints are high, the
bare noninteger separation is \(\mathcal{L}^{-1}\), and the scalar
shifts are \(O(N^{-2}\mathcal{L})\). An entrywise primitive therefore
costs \(O(r\mathcal{L})\) there and \(O(r/R)\) when it meets the
nondiagonal low block; its commutator with that block retains
the latter gain. Sum the shifts and their adjoints with their
summable angle coefficients in the Schur norm. Combine the free
kinetic commutator and the time derivative before estimating the
remainder. This gives \(o(1)\) after division by the slow scale.
In detail, let \(\mathcal B\) be the exponent formed by these
primitives. By the three divisor bounds it has Schur size
\(O(r\mathcal L)\), and its entries with a low endpoint have size
\(O(r/R)\). As in the identity for \(\mathcal M=e^{\mathcal B}\) in
Section~\ref{sec:ins-propagator}, the term
\([D_{\rm free},\mathcal B]-i\partial_t\mathcal B\) cancels the
discarded entries up to the difference between the actual and the
bare divisor times \(\mathcal B\). That difference consists of the
scalar shifts of the two legs and of the root, and of the displacement
of the speed of the old harmonic, all \(O(N^{-2}\mathcal L)\); this
leaves \(O(rN^{-2}\mathcal L^2)\). The commutators of \(\mathcal B\)
with the retained stars and pairs, and all terms of second order,
are \(O(r^2\mathcal L^2)\). The commutator with the nondiagonal low
block is \(O(r\mathcal L/R)\), and conjugation of the terms of order
\(r\) by the column frame costs \(O(rN^{-1}\mathcal L)\). After
division by \(r\lambda\) the total is
\(O\bigl(\lambda^{-1}(N^{-2}+r+R^{-1}+N^{-1})\mathcal L^2\bigr)\),
which tends to zero at fixed \(\lambda\).

On each remaining pair use the full nominal high real diagonals.
At most one singleton detuning can be near zero, by separation
of distinct \(D_e(N)\); primitive division of the others has normalized
cost \(O(N^{2-s}\mathcal{L})\). Retain this one harmonic with its
actual detuning, and absorb the higher weak tails. To see that
one can choose a positive mixed time weight uniformly in \(R\)
and small \(\lambda\) at this fixed tree, put
\(u=(q/P)^2,\ p_i=P_i/P,\ x=(n-1)/N,\ y=(n'-1)/N\).
Then \(x+y=p_i,\ xy=u\). The two squared lengths are
\(x^2+u=xp_i\) and \(y^2+u=yp_i\). Thus the sum of their reciprocals
is \(1/u\), whereas \(x^{-2}+y^{-2}=(p_i^2-2u)/u^2\).
Subtracting the parent's correction gives the leading mismatch in
\(N^{-2}\) units,
\[
 F_i(u)=A p_i^2/u^2-(2A+B)/u-(A-B)/p_i^2 .
\]
In detail, the two legs are \((x,q/P)\) and \((y,-q/P)\) in
normalized coordinates and the root is \((p_i,0)\), so the mismatch
is \(\Xi((x,q/P))+\Xi((y,-q/P))-\Xi((p_i,0))\) with
\[
 \frac1{x^2}+\frac1{y^2}=\frac{(x+y)^2-2xy}{(xy)^2}=\frac{p_i^2-2u}{u^2},
 \qquad
 \frac1{x^2+u}+\frac1{y^2+u}=\frac{x+y}{xyp_i}=\frac1u ,
\]
and \(A(p_i^2-2u)/u^2-B/u-(A-B)/p_i^2\) is the displayed function.
An exponent of the two-by-two
block with nonzero real part requires
\(\lvert F_i(u)-D_e(N)\rvert\le C(N^{2-s}+N^{-1}\mathcal{L})\).
Indeed the normalized detuning is the physical one divided by
\(r\lambda=N^{-s}\lambda\), that is
\(\delta=\lambda^{-1}N^{s-2}(F_i(u)-D_e(N)+O(N^{-1}\mathcal L))\), and a
nonzero real part of \(\sqrt{|g|^2-\delta^2/4}\) requires
\(|\delta|<2|g|\) with the bounded normalized coupling \(g\).
On the permitted compact range \(u\le p_i^2/4\) the small values
of \(u\) are excluded by \(Ap_i^2>0\), uniformly for these
rational trees near their fixed real limit. The sublevel set in
\(1/u\) lies in at most two intervals of length
\(C(N^{2-s}+N^{-1}\mathcal{L})^{1/2}\), since this is a quadratic
with leading coefficient bounded below. Explicitly, if
\(Q(v)=\alpha(v-v_0)^2+\mu\) with \(\alpha\ge\alpha_0>0\), then
\(|Q(v)|\le\varepsilon\) confines \((v-v_0)^2\) to an interval
\([a,b]\) with \(0\le a\le b\) and \(b-a\le2\varepsilon/\alpha_0\), and
\(\sqrt b-\sqrt a\le\sqrt{b-a}\); the set is therefore contained in
at most two intervals of length at most
\(2\sqrt{2\varepsilon/\alpha_0}\). Here \(v=1/u\) and
\(\alpha=Ap_i^2\). The mesh of distinct
positive \(|q|\) there is \(\gtrsim1/P\) in \(1/u\). Indeed
\(1/u=P^2/q^2\), the permitted range \(u_0\le u\le p_i^2/4\) confines
\(|q|\) to \([P\sqrt{u_0},P|p_i|/2]\), and for consecutive values
\[
 \frac{P^2}{q^2}-\frac{P^2}{(|q|+1)^2}
 =\frac{P^2(2|q|+1)}{q^2(|q|+1)^2}\ge\frac{P^2}{2|q|^3}\ge\frac cP .
\]
For fixed \(s>2\), Lemma~\ref{lem:rational-tree-uniformity} and the polylogarithmic bound on \(\mathcal L\)
give \(P(N^{2-s}+N^{-1}\mathcal L)^{1/2}=o(1)\). Each interval therefore
contains at most one value of \(|q|\) for large \(R\), and \(x+y=p_i\),
\(xy=u\) then determine the unordered pair of horizontal coordinates.
Thus at most a fixed number of such rates occurs per root and
singleton. The count is independent of \(\lambda\) small at the
fixed passive preparation: normalized singleton strengths and their
number are bounded, and the geometry can be taken in fixed compact
intervals separated from zero. Partition a compact positive interval
strictly below the finite orbit's exponential rate into more than
twice this number of cells, summing over all roots and singletons.
For every sufficiently large \(R\) a cell has a smaller interior
interval with uniform clearance. Use three strictly ordered weights
in that interval for the exact landing and the two subsequent
history estimates. A subsequence with the same cell suffices.

The pair itself has constant coefficients after a unitary half-phase
rotation, with real growth rates
\(\pm\operatorname{Re}\sqrt{|g|^2-\delta^2/4}\) for actual
coupling \(g\) and detuning \(\delta\). Its two-sided Green bound
on the cleared weights is uniform: bounded \(g,\delta\) use the
matrix exponential, including the continuous Jordan limit;
large \(|\delta|\) gives uniformly conditioned imaginary
eigenvalues. The disjoint Hermitian complement propagates
isometrically in the flattened weights used in
Proposition~\ref{prop:insertion}. The axial root abscissae are
\(O(N)\); all affected axial--low stars are in a common flat
horizontal band of this size. The low old block preserves its
own vertical index. The discarded homogeneous operators cost
\(O(\lambda^{1/5})+o_R(1)\) at the fixed preparation and are
absorbed on this same inverse, choosing \(\lambda\) and \(R\)
in that order. This argument requires no lower bound on the
discriminant at a double root.
The estimates are those of steps (ii)--(vi) in the proof of
Lemma~\ref{lem:common-endpoint-green}, which use only the following
facts, all established above for the present scaling. The normalized
couplings \(g\) of the retained pairs are bounded by a constant
\(C_1\) independent of \(R\) and of small \(\lambda\). Every weight
\(\sigma\) in the chosen interior interval satisfies
\(|\sigma-k_+|\ge\kappa_0>0\) for every rate
\(k_+=\operatorname{Re}\sqrt{|g|^2-\delta^2/4}\), with \(\kappa_0\)
independent of \(R\). The complement of the pairs is isometric in the
flattened norm, and the two legs of a pair have weights of bounded
ratio. Let \(C_W\) bound \(\max(w_1/w_2,w_2/w_1)\) over these pairs
and both regularity indices. Absorb this factor into the constant in
\eqref{eq:h2-pair-propagator}. The weighted projections of a pair with \(k_+>\sigma\) are bounded by
\(\tfrac{C_W}{2}(1+2C_1/\sigma)\), and the principal Green operator is
bounded by \(C(\kappa_0^{-1}+\kappa_0^{-2}+\sigma^{-1})\) on the
weighted space, for both regularity indices and for all weights in
the interval. The remainder
\(O(\lambda^{1/5})+o_R(1)\) is then removed by the Neumann series,
the axial block is integrated backwards from infinity with the
action-angle shear, and the inverses at the different weights and
regularities agree on common data, exactly as in that proof.

\subsection{The full equation and the fixed datum}

Put \(p=\lceil s\rceil\), and use the scaled norm \(Z^p\) with
weight \((1+|n|/N+|q|/P)^p\), \(m'=Rq\), and its weak companion
\(Z^2\). Prescribed heads have bounded weighted Wiener norms at
fixed tree and normalized height. Two unrestricted factors cost
\(C\sqrt{NP}\) by summing the reciprocal squared weight over
\((n,q)\). Indeed, grouping \(n\) in blocks of length \(N\) and \(q\)
in blocks of length \(P\),
\[
 \sum_{n,q}(1+|n|/N+|q|/P)^{-2p}
 \le CNP\sum_{j,k}(1+|j|+|k|)^{-2p}\le C_pNP,
\]
and the argument of Section~\ref{sec:h2-extension} for the scaled
algebra applies with \(\sqrt{NP}\) in place of \(\sqrt R\).
A perturbation of size \(O(r)\) in \(Z^p\) has physical
\(H^{p-1}\) size \(O(rN^{p-1})=o(1)\): since \(Rq=Nq/P\), one has
\(\langle(n,Rq)\rangle\le CN(1+|n|/N+|q|/P)\), hence
\(\|f\|_{H^b}\le CN^b\|f\|_{Z^p}\) for \(0\le b\le p\), and
\(rN^{p-1}=N^{p-1-s}\) with \(p=\lceil s\rceil<s+1\).
The change of either
Birkhoff differential from its smooth old value is \(O(r)\) on
\(L^2\), and its nonlinear part changes by \(O(rN^{p-1})\)
from \(H^{p-1}\) to \(H^p\), by integer one-smoothing. Applied
to a scaled unit vector its \(Z^p\) size is consequently
\(O(r+N^{-p}(rN^{p-1})N^{p-1})=O(N^{-1})\). Here the three factors
in the second term are the weight \(N^{-p}\) of the \(H^p\) part of
the scaled norm on the axis, the change of the nonlinear part, and
the \(H^{p-1}\) norm of a unit vector of \(Z^p\); their product is
\(rN^{p-2}=N^{p-2-s}\le N^{-1}\), again because \(p<s+1\).
The old differential itself on a high single carrier equals its
Fourier differential up to \(O(N^{-1}\mathcal{L})\) in the scaled
norm; use one-smoothing and the negative-order dual estimate as in
Proposition~\ref{prop:insertion}. Multiplication by this high-carrier error
has operator norm \(O(\sqrt{NP}\,N^{-1}\mathcal L)=o(1)\)
on \(Z^p\). The physical terminal Fourier shape follows separately from
\(h^{p-1}\to H^p\) local Lipschitz continuity of the nonlinear part.

We detail which forcing terms enter the terminal ball. In the exact
old-action, old-angle, complementary axial and transverse equations
use the same split of the known finite packet and the unknown as
in Proposition~\ref{prop:insertion}. In
\(u=U+w,\ U=P_{m'=0}u\), the forcing of the axial equation is
\(2U P_0|w|^2+\bar U P_0w^2+P_0(|w|^2w)\), with the common mass
gauge, and the transverse quadratic forcing is
\(P_\perp(2U|w|^2+\bar U w^2)\).
The equal-head density in the first axial term is a common phase
on the whole axial field and has zero action derivative. For other
two-head terms the nonzero integer macro dot product gives a divisor
at least \(R^2-O(N\log N)\). At a right angle the tree contains
the output, with nonzero limiting horizontal velocity. The linear
term in \(N\) forces the spatial key \(b=1\); the remaining
integer-time key is \(\nu=1\). Every such harmonic is already
included in the exact finite equation. Other low harmonics with
\(b=1\) have the polynomial frozen divisor above. Subtract each
nonresonant quadratic particular once, with its fast phase intact.
The derivative of a true finite amplitude in physical time costs
\(O(r\lambda)\); the change of a cutoff old harmonic's speed costs
\(O(N^{-2}\mathcal{L})\). The \(O(1)\) lower-order part of the low
Hermitian operator can meet a low particular only after the
\(R^2\) kinetic division. The all-order column residual after
division by \(r\lambda\) is \(O(N^{-1}\mathcal{L})\), also for its
moving connection on known heads up to a linear factor in the old
scaled action. The \(O(r^2)\) shifts at the exact terminal axial
actions are included in its true frequencies. This gives the total
normalized pure error \(e_N=O(N^{-1}\mathcal{L})\).

For precision let \(H,A_1,J_o,\Theta\) measure the decaying
weighted normalized transverse error, complementary axial field,
old action and old phase, respectively, and \(X\) their sum.
In the terminal operator leave the whole Hermitian low matrix
at the trial state on the left. There is no small parameter in
its dependence on the old angles before division by slow time.
The scaled fields are normalized by \(r\), the old action by \(r\),
and the old angle without an amplitude factor. The four bounds for
the Green image are
\[
\begin{aligned}
 H_{\rm out}&\le C(e_N+cX+\sqrt{NP}\,X^2+o(1)H),\\
 (A_1)_{\rm out}&\le C(e_N+J_o+c(H+A_1+\Theta)
                                +\sqrt{NP}\,X^2+o(1)A_1),\\
 (J_o)_{\rm out}&\le C(e_N+cH+c^2(A_1+\Theta)+\sqrt{NP}\,X^2),\\
 \Theta_{\rm out}&\le C(e_N+J_o+r(A_1+H)+\sqrt{NP}\,X^2).
\end{aligned}
\]
Here \(c\) bounds the nonterminal heads on a fixed sufficiently
late slow-time interval. The term \(c^2\Theta\) in the third row
records the change of the known two-head source under a displacement
of the old angles, as explained below. It can be dispensed with by
evaluating the high-frequency gain at the actual old angles, as in
the proof of Proposition~\ref{prop:insertion} and at \(s=2\); we keep
it because it has the same form as the term \(c^2A_1\) and does not
affect the choice of the radii. The transverse one-error terms carry
a mobile factor apart from the homogeneous operator. A variation
of the old action in a high scalar has its \(O(N^{-2}\mathcal{L})\)
derivative. Variation of the frequency multiplying a persistent
axial root gives the displayed unsuppressed \(J_o\). In the old action equation, the common density is the only leading
quadratic term in the prescribed transverse carriers with zero horizontal
macro frequency. After truncation of the old angular harmonics, the other
such terms have high horizontal frequency, gaining \(O(N^{-1}\mathcal L)\)
on old coordinates by smoothing and duality. A change in the old angle
of this quadratic source retains two mobile factors (or a quadratic
error). The cubic term in three prescribed transverse waves can have low
output. After division by \(r^2\lambda\) it is \(O((r/\lambda)c^3)\)
in the exponentially weighted norm, since all three waves decay.
With \(r=N^{-s}\) and \(\lambda\) fixed, \(r/\lambda=o(N^{-1})\),
which fits the \(O(N^{-1}\mathcal L)\) pure error bound.
The ordinary axial shear gives \(J_o\) in the angular equation.
Each two-error product obeys the displayed lattice bound, and a
third arbitrary field costs an extra \(Cr\sqrt{NP}X\).
The constants on the unsuppressed linear terms are uniform for
large \(R\): the old smooth degrees of freedom, normalized finite
Wiener norms and terminal spectral gap are fixed, and each
homological or high-cut factor \(\mathcal{L}\) on the remaining
linear terms accompanies a fixed negative power of \(N\).

\smallskip
\noindent\emph{Derivation of the four rows.} The argument is that of
\eqref{eq:four-rows-ins} and of its endpoint version in
Section~\ref{sec:h2-extension}; we list the terms in the present
scaling. Time is \(\theta=r\lambda t\) with \(r=N^{-s}\), so
\eqref{eq:exact-h-ins} holds with \(r^{-1},r^{-2}\) replaced by
\((r\lambda)^{-1},(r^2\lambda)^{-1}\). The fields are measured in
\(Z^p\) with \(p=\lceil s\rceil\), prescribed factors in its Wiener
counterpart, and two free factors cost \(C\sqrt{NP}\). Write \(e\)
for the decaying exponential weight. Now
\(N^{-2}/(r\lambda)=N^{s-2}/\lambda\) is large. Every correction of
size \(N^{-2}\) to a frequency has therefore been matched in advance
to accuracy \(O(N^{-J-2}\mathcal L)\) with \(J>s+4\): the high columns
with their instantaneous scalars and scalar phases, the carriers of
the heads with the angle means of these scalars, the exact axial
frequencies on the terminal torus, and the rows through the fixed
point equation for \(y\). After division by \(r\lambda\) an error of
this size is \(O(N^{s-J-2}\mathcal L)=O(N^{-6}\mathcal L)\). The chart
estimates are the three stated at the beginning of this subsection:
the change of a Birkhoff differential on the tube is \(O(N^{-1})\) on
\(Z^p\); the old differential differs from its Fourier part by
\(O(N^{-1}\mathcal L)\) on a high carrier; and, by symplectic
adjointness, the same gain holds for the projection of a high carrier
to the old coordinates.

\emph{First row.} The forcing of the transverse equation, after
division by \(r^2\lambda\), consists of the following terms.
\begin{itemize}
\item The residual of the prescribed fields: the column residual and
its moving connection on the known heads,
\(O(N^{-1}\mathcal L)(1+J_o)\); the chart errors on the known heads,
\(O(N^{-1}\mathcal L)\); for each divided particular the derivative
of its finite amplitude and the change of speed of its cutoff
harmonic, of relative size \(O((r\lambda+N^{-2}\mathcal L)\mathcal L)\);
the angular tails, \(O(N^{-J-2})\); and the products of three known
factors, \(O(r)\). Their sum is \(O(e_N)e+o(1)J_oe\), and the last
term is at most \(cJ_oe\) for large \(R\).
\item One error with a mobile known factor: \(Y-Y_\infty\) in the
mixed operator, \(O(cH)\); the axial error against a head,
\(O(cA_1)\); the displacement of the old angles and actions in the
prescribed harmonics, \(O(c(\Theta+rJ_o))\); the variation
\(2\lambda^{-1}\sum_{j\in F}(J_o)_j\) of the old-mass scalar on the
packet, \(O_\lambda(cJ_o)\).
\item Two errors, \(O(\sqrt{NP}\,X^2)\), including the old-mass scalar
on \(h\); three errors, \(O(rNP\,X^3)\le O(\sqrt{NP}\,X^2)\) when
\(\sqrt{NP}\,X\le1\).
\item Terms linear in \(h\) with a coefficient tending to zero: the
chart error on the prescribed axial planes against \(h\),
\(O(\sqrt{NP}\,N^{-1}\mathcal L)H\); the variation of a high scalar
with the old actions, \(O(N^{-2}\mathcal L)(J_o+r)H\); and the moving
connection of the scalar phases and of the column frame,
\(O(N^{-1}\mathcal L)(J_o+r)H\).
\end{itemize}
The homogeneous operators removed from the principal part cost
\(O(\lambda^{1/5})+o_R(1)\) and have been absorbed in the Green
inverse by the Neumann series.

\emph{Second row.} Keep the real diagonal at the trial actions and
integrate backwards from the terminal time. The \(\ell^1\) distance of
the trial actions from the terminal ones is at most
\(C(rJ_o+r^2(c+A_1+o(1)))e\), so the frequency of a persistent root
changes by at most this amount and its mismatch contributes
\(C\lambda^{-1}[J_o+r(c+A_1+o(1))]e\); this is the unsuppressed
\(J_o\). In the force, two known heads are cancelled by the finite
equation at the key \((1,1)\) and divided otherwise; one known and
one unknown transverse factor give \(O(cH)\); two unknown factors
give \(O(\sqrt{NP}\,X^2)\). The replacement of the Birkhoff
differential by its old value and then by its Fourier part costs
\(O(N^{-1}\mathcal L)\) on the known high outputs, and
\(O(c\,\Theta)\) arises from the old angles in the prescribed
harmonics.

\emph{Third row.} In
\(\partial_\theta(J_o)_j=2(r^2\lambda)^{-1}\operatorname{Re}(\bar z_j(D_f)_j)\)
the equal-head density drops out. Every other product of two known
heads has horizontal frequency of modulus \(\gtrsim N/\mathcal L\), so
its projection to the old coordinates gains \(O(N^{-1}\mathcal L)\)
and contributes \(O(e_N)\). A transverse error enters with a known
head, \(O(cH)\), or with a second error, \(O(\sqrt{NP}\,X^2)\). The
axial error and the old angle displacement enter through the
coefficients of a product of two transverse factors; with two known
heads this gives \(O(c^2(A_1+\Theta))\), and otherwise a quadratic
term.

\emph{Fourth row.} In
\(\partial_\theta\Theta_j=(r\lambda)^{-1}\bigl(\omega_j(I)-\omega_j(I^\infty)
-\operatorname{Im}((D_f)_j/z_j)\bigr)\) the frequency difference is
bounded by the \(\ell^1\) distance of the actions, which gives
\(C\lambda^{-1}[J_o+r(c+A_1+o(1))]e\). The force term is that of the
third row with one factor \(r\) more, because the division is by
\(r\lambda\) and the old \(|z_j|\) are bounded below.

With all four errors zero, \(Y\) replaced by \(Y_\infty\) and the
heads and particulars removed, the reconstruction is the exact
terminal axial solution, and all four residuals vanish. Every pure
error therefore carries the factor \(e\).

Fix \(1/2<a<1,\ \rho=N^{-a}\). Then \(e_N=o(\rho)\) and
\(\sqrt{NP}\rho=o(1)\). Give \(H,J_o\) radius \(M\rho\) and
\(A_1,\Theta\) radius \(M'\rho\), choose \(M'\) large relative
to \(M\), then choose \(c\) sufficiently small, and finally \(R\)
large. The displayed bounds preserve this ball. Explicitly,
\(e_N/\rho=O(N^{a-1}\mathcal L)\) and
\(\sqrt{NP}\,\rho=N^{1/2-a}\sqrt P\) tend to zero because
\(1/2<a<1\) and \(P\) is polylogarithmic in \(N\). On the ball,
\(X\le2(M+M')\rho\), and with the constant \(C\ge1\) of the four rows
\[
\begin{aligned}
 H_{\rm out}/\rho&\le C\bigl(2c(M+M')+o(1)\bigr),
 & (J_o)_{\rm out}/\rho&\le C\bigl(cM+2c^2M'+o(1)\bigr),\\
 (A_1)_{\rm out}/\rho&\le C\bigl(M+c(M+2M')+o(1)\bigr),
 & \Theta_{\rm out}/\rho&\le C\bigl(M+o(1)\bigr),
\end{aligned}
\]
where \(o(1)\) collects \(e_N/\rho\), the quadratic terms
\(\sqrt{NP}\,X^2/\rho\le4(M+M')^2\sqrt{NP}\,\rho\), the terms with a
factor \(r\) and the \(o(1)\) multiples of \(H,A_1\). For \(M=1\),
\(M'=4C\) and \(c\le1/(16C(1+M'))\) the four right sides are at most
\(M/2,\ M/2,\ M'/2,\ M'/2\) for large \(R\). The finite-strip
construction uses old velocities evaluated from the exact vector
field at the trial state; its chain rule agrees with the actual
one at a fixed point. On a finite strip use the closed trials
in \(C_tZ^2\cap L^\infty_t Z^p\) with the four pointwise weighted
bounds. At fixed \(R\) the time derivative of a field image is bounded
in \(Z^{p-2}\), and the images are relatively compact in \(C_t Z^2\).
Continuity
of the trial-dependent low propagator here only needs a finite
interval at fixed \(R\) and follows by Duhamel. Schauder and then
a diagonal finite-strip limit give the decaying solution on the
late half-line. On the remaining fixed slow-time interval retain
the actual Hermitian propagator; the bounded finite coefficients
and the same quadratic estimate give \(O_{\rm stage}(\rho)\)
backward by Gronwall and continuity. The reconstructed solution
satisfies the physical Duhamel formula in \(H^2\) and takes values
in \(H^p\) at every time. The local \(H^p\) solution starting from
each such datum agrees with it by \(H^2\) uniqueness, giving strong
\(H^p\) continuity.
The compactness argument is the one written out after the invariance
check in Section~\ref{sec:h2-extension}, with \(Z^p\), \(Z^2\) and
\(Z^{p-2}\) in place of \(Z^2\), \(Z^{p_0}\) and \(Z^0\). The time
derivative of an image is bounded in \(Z^{p-2}\); for \(p=3\),
interpolation between \(Z^{1}\) and \(Z^3\) gives equicontinuity in
\(Z^2\), and for \(p\ge4\) this is immediate. Bounded sets of \(Z^p\) are
precompact in \(Z^2\) on the fixed lattice, since the ratio of the
weights tends to zero, and the products are continuous on \(Z^2\) by the
scaled algebra estimate with exponent \(2>1\) and the factor
\(\sqrt{NP}\).
At fixed \(R\) the homogeneous inverse at the asymptotic axial
orbit is also bounded in ordinary physical \(H^2,H^s\) at the
same two history weights by the pair and Schur bounds just proved.
This is the inverse at the final smooth torus used for the earlier
history replacements, where constants may depend on that prefix.
In detail, the principal Green operator is a time integral
\(\int\mathcal K(\tau,u)f(u)\,du\) whose kernel is bounded on \(Z^2\)
and on \(Z^p\), \(p=\lceil s\rceil\), by
\(C(1+|\tau-u|)e^{-\kappa_1|\tau-u|}\) with a fixed \(\kappa_1>0\)
after the time weight is included, by steps (ii)--(iii) of the proof
of Lemma~\ref{lem:common-endpoint-green} in the present scaling; the
remainder \(\mathcal R(\tau)\) is bounded on both spaces with norm at
most \(\epsilon_*\), and the axial block is given by the explicit
backward integrals of step~(v). With
\(\omega(n,q)=1+|n|/N+|q|/P\), these spaces are the weighted \(\ell^2\)
spaces \(Z^b=\ell^2(\omega^b)\), and complex interpolation gives
\([\ell^2(\omega^2),\ell^2(\omega^p)]_\vartheta=\ell^2(\omega^s)\) for
\(s=(1-\vartheta)2+\vartheta p\). Hence for each fixed \(\tau,u\)
the kernel and the remainder are bounded on \(Z^s\) by the larger of
their two norms, and integrating in time gives the same weighted bound
on \(Z^s\); the Neumann series converges there for the same reason.
The changes of variables of the reduction are bounded with bounded
inverses on \(Z^2\) and \(Z^p\) and are interpolated in the same way.
Finally, at fixed \(R\) the weight \(\omega(n,q)\) is comparable to
\(\bk{(n,Rq)}\) in both directions, with constants depending on
\(R,P\), so \(Z^b\) and \(H^b\) are equivalent on the lattice for
\(0\le b\le p\). This gives the boundedness in \(H^2\) and \(H^s\) at
both history weights, with constants depending on \(R\).

For a balanced split \(v+w=h\) write \(z=|v|^2/|h|^2\). Orthogonality
gives \(|w|^2/|h|^2=1-z\), so the ratio of the sum of the two
homogeneous squared \(H^s\) weights to the parent weight is
\(z^s+(1-z)^s\le c_s<1\), uniformly on the allowed interval of \(z\).
Summing by levels and using the level invariant gives a bound
independent of the depth for the ideal primary contribution and
\(O(c_s^{d/2})\) for its norm at launch. At fixed depth the weights
\(r\bk{k_v}^s\) on the finite set of heads are uniformly bounded
along the rational approximations to the chosen tree. Thus the
difference between the fitted and ideal packet norms is bounded
by a fixed-tree constant times their uniform amplitude error in the
carrier coordinates, both at launch and on the half-line. Given a
tolerance, choose the depth for the ideal launch bound, then the
auxiliary amplitudes. At that depth choose the axial preparation and
finite comparison accuracy to put the weighted fitting error within
the tolerance, and finally take \(R\) sufficiently large. In \(H^s\) the
physical reconstruction of the normalized error \(O(\rho)\)
costs \(O(rN^s\rho)=O(\rho)\). In \(H^2\) the prescribed
unit-height packet costs \(O_{\rm stage}(rN^2)=O_{\rm stage}(N^{2-s})\)
uniformly in time. Taking tolerances tending to zero thus gives
raw insertions that converge initially in \(H^s\) to the old target,
stay bounded in \(H^s\) at fixed new height, and converge uniformly
on the half-line in \(H^2\). Up to a uniform \(o(1)\) in \(H^s\) their
increment is the finite plane sum, by the preceding reconstruction
and exact locking of all earlier occupied frequencies and origins.
Apply the lifting argument of Proposition~\ref{prop:history} in the
two-norm form given at the end of Section~\ref{sec:h2-extension},
with \(H^2\) and \(H^s\) in place of \(H^{p_0}\) and \(H^2\), saving
compatible weights as in the endpoint argument. The targets of the
present section are quasi-periodic rather than periodic, so the
Floquet splitting of Lemma~\ref{lem:floquet} is not available and
the proposition is not applied literally. As at \(s=2\), that form of
the argument uses only two ingredients: a pair of compatible right
inverses of the old linearization at two weights \(0<\nu<\beta\),
with the same boundary conditions in the two norms, and new inputs
that are bounded in the strong norm, converge strongly at the initial
time, and converge uniformly on the half-line in the weak norm.
Here the inverses are the homogeneous inverse at the final smooth
torus constructed in the preceding paragraph, which plays the role of
Lemma~\ref{lem:common-endpoint-green}, and the inputs have the three
properties just listed. With \(\mathcal G_V\) denoting this inverse,
the replacement correction \(d+h\) solves the same equation
\(h=\mathcal G_V\{F_n(d+h)-F(d)\}\); the contraction at the weight
\(\beta\) in the sum of the weighted \(H^s\) norm and \(\Lambda\)
times the weighted \(H^2\) norm, and the convergence at the weight
\(\nu\), follow as there. The product bound is tame with \(H^2\hookrightarrow
L^\infty\); each saved correction decays at its original larger
weight and converges strongly at the smaller weight. Divisibility
of the vertical spacings puts every new input on the old lattice.
Each old trajectory has precompact \(H^s\) range, so its cross
term with a finite new sum supported beyond \(R/2\) tends uniformly
to zero. These preparations give summable strong errors initially,
summable weak errors and summable deviations of the squared
increment from its finite plane sum as in the counting argument.
Writing \(D_j=K_j^s\), the new slow time tends uniformly to zero
for \(t\le D_j^{3/4}\) and to infinity for \(t\ge D_j^{5/4}\)
at fixed preparation: since \(K_j=1+N\) and \(r=N^{-s}\),
\[
 r\lambda D_j^{3/4}\le2^{3s/4}\lambda N^{-s/4},\qquad
 r\lambda D_j^{5/4}\ge\lambda N^{s/4}.
\]
The middle windows can be taken disjoint,
and \(\ell(K_j)\) can absorb its fixed entrance time divided by
\(\lambda\). This proves the counting and slow-growth conclusions
for this \(s\) after passage to one strongly convergent sequence of
initial states. For each \(p'<s\) uniform lower-norm boundedness
follows by imposing summably small increments in a countable
sequence of exponents increasing to \(s\), using interpolation
with the uniform \(H^2\) smallness at each fixed height.

\subsection{The limiting axial solution and data near a plane wave}

\begin{corollary}\label{cor:super-asymptote}
For every finite \(s\ge2\), the solution in Theorem~\ref{thm:clock}
can be chosen with an axial solution \(v\) which is almost periodic in
\(H^\sigma\) for every \(1<\sigma<s\), such that
\[
 (1+t)^b\|u(t)-v(t)\|_{H^\sigma}\longrightarrow0
 \quad\text{for }\quad
 0\le b<\frac{\min\{s-\sigma,1\}}{s}.
\]
For every \(j\) and every real \(t\),
\(\langle K_j\rangle^s|\widehat v((K_j,0),t)|\ge1/2\).
In particular \(v(t)\notin H^s\) at every time.
\end{corollary}
\begin{proof}
Use unit height and put \(N=K-1\), so \(r=N^{-s}\).
For any finite list of pairs \((\sigma,b)\) as in the assertion choose
\(a\in(1/2,1)\) with \(a>sb\) for every pair on the list.
At the endpoint the ball estimate also admits \(\rho=N^{-a}\):
at fixed tree \(N=PR\), and both \(e_R/\rho\) and
\(\sqrt R\,\rho\) tend to zero. In Section~\ref{sec:super}
the same choice absorbs the indicated powers of \(\log N\).
Take three ordered weights in one gap, with the largest one used for
the new solution and the two smaller ones reserved for its history.
For the estimate as \(N\) increases one can use the same smaller
positive decay rate on a subsequence. Denote the actual raw solution by
\(W\) and its exact rounded axial target by \(V'\). For each fixed
preparation there are constants \(C,c,T>0\) such that
\[
 \|W(t)-V'(t)\|_{H^\sigma}\le C(1+\log N)^C
       (N^{-a}+N^{\sigma-s})\exp[-c(r\lambda t-T)_+] .
\]
Indeed the change in the finitely many old angles is bounded by
the radius of the correction, and that in the old actions by \(r\)
times this radius. All other coordinates of the difference have
scaled norm at most \(C r\) on the early interval, with exponential
decay thereafter. This uses the exact finite trajectory and its own
endpoint. Its fixed accuracy relative to the ideal trajectory creates
no stationary source. Converting a scaled field to \(H^\sigma\)
costs at most \(C N^\sigma\). The old axial coordinates range over
a compact smooth torus. For the inverse chart its nonlinear part is
locally Lipschitz from \(h^{m-1}\) to \(H^m\) with
\(m=\lceil\sigma\rceil\ge2\), and the added high coordinates are
small in \(h^{m-1}\). These bounds prove the displayed inequality
in the mass gauge of the new solution.

Multiplication by \((1+t)^b\) costs at most
\(C_{b,c,T}(r\lambda)^{-b}\), with \(\lambda\) fixed before the dilation.
Indeed, writing \(\theta=r\lambda t\) and using \(r\lambda\le1\),
\[
 (1+t)^b e^{-c(\theta-T)_+}
 \le(r\lambda)^{-b}(1+\theta)^be^{-c(\theta-T)_+}
 \le(r\lambda)^{-b}\sup_{\theta\ge0}(1+\theta)^be^{-c(\theta-T)_+}.
\]
Since \(r=N^{-s}\), the displayed bound multiplied by \((1+t)^b\) is
at most
\(C\lambda^{-b}(1+\log N)^C(N^{sb-a}+N^{sb+\sigma-s})\).
Both resulting powers \(sb-a\) and \(sb+\sigma-s\) are negative:
the first by the choice of \(a\), which is possible because \(sb<1\),
and the second because \(sb<s-\sigma\).
Every previously saved correction converges in its smaller exponential
weight on a late half-line and strongly on the preceding compact
interval. Its change can therefore also be made arbitrarily small in
the polynomially weighted norm. Impose summable errors on successively
larger finite lists of rational pairs as in Corollary~\ref{cor:asymptote}.
The same choice gives summable uniform changes of the axial targets in
every \(H^\sigma\) on the real line. Indeed the old angular velocities
and origins agree in their respective mass gauges, the new clock cost
is arbitrarily small, and the new axial roots cost
\(O_{\rm stage}(N^{\sigma-s})\) there.
Each target is a continuous function on a finite angle torus evaluated
along a linear flow. Approximation by trigonometric polynomials on
that torus gives almost periodicity of the uniform limit.
The identification of the limit and removal of the common mass gauge
are as in Corollary~\ref{cor:asymptote}. For a fixed finite target,
there are arbitrarily large translations returning all its angles
arbitrarily close to their initial values. They carry the persistent
lower bound for each previously selected Fourier coefficient back to
any prescribed real time. Passing to the limit in that coefficient
gives the last assertion.
\end{proof}

\begin{corollary}\label{cor:near-plane}
Fix \(2\le s<\infty\). Every \(H^s(\T^2)\) neighborhood of a
plane wave \(a e^{ik\cdot x}\), \(a\in\C\), \(k\in\Z^2\), contains a datum
whose solution satisfies \(\|u(t)\|_{H^s}\to\infty\).
\end{corollary}
\begin{proof}
For a small nonzero amplitude start with a single axial Birkhoff
coordinate at index one. Translation and constant-phase equivariance
force its physical inverse to be a single Fourier wave, and preservation
of mass gives its modulus. Open a second sufficiently small positive
action, with phase origins chosen in the prescribed neighborhood.
For a two-site torus of this kind \(A>0\). Indeed the Liouville identity
gives \(A\) as the angle mean of the spatial variance of the density.
If it vanished, every orbit on this torus would have spatially constant
positive density at all times. The mass continuity equation would
give a spatially constant phase derivative and hence a single plane
wave, which has only one occupied Birkhoff coordinate. In formulas,
\(\partial_t|U|^2=-2\partial_x\operatorname{Im}(\bar UU_x)\) for the
axial equation, and for \(U=\sqrt M\,e^{i\varphi}\) with constant
\(M>0\) this reads \(0=-2M\varphi_{xx}\). Hence \(\varphi_x\) is a
constant, an integer \(k\) by periodicity, and
\(U=\sqrt M\,e^{i(kx+\varphi_0(t))}\).
Since \(D\chi\) is close to minus the identity, the frequency image of the positive
two-site actions in this neighborhood contains an open box. At \(s=2\) choose rational slopes
with numerators congruent to one modulo five and one common denominator;
for \(s>2\) choose two rationally independent algebraic slopes.
The allowed lists are dense at either exponent. Choose them for one
small positive value of the common factor and then vary that factor
in an open interval within the box to avoid the countable exceptional
sets of the construction. Every subsequent initial increment can be
made summably small in this neighborhood. A constant phase and a
Galilean transformation give the same conclusion near a small constant.

For general \(a\ne0\) choose an integer \(L\) with \(a/L\) in this
small range. If \(v\) solves \eqref{eq:nls}, so does
\(L v(L^2t,Lx)\), and its datum is arbitrarily close to \(a\) when
\(v(0)\) is sufficiently close to \(a/L\). Its \(H^s\) norm is
at least \(L\|v(L^2t)\|_{H^s}\) by the Fourier formula: the
coefficient of \(Lv(L^2t,Lx)\) at \(Lk\) is \(L\widehat v(k,L^2t)\),
its other coefficients vanish, and \(\langle Lk\rangle\ge\langle k\rangle\).
The three terms \(i\partial_t\), \(-\Delta\) and the cubic term each
acquire the factor \(L^3\) under this substitution. Galilean
modulation by the fixed \(k\) preserves divergence and continuity
in \(H^s\). The case \(a=0\) is Theorem~\ref{thm:main}.
\end{proof}

\section{Two-sided growth for \texorpdfstring{$s>2$}{s above two}}\label{sec:two-sided}
Fix $s>2$. At stage $j$ we use three solutions of the equation in
the autonomous mass gauge. The axial solution $V_{j-1}$ is the
starting background. The two-dimensional solution $W_j$ stays close
to $V_{j-1}$ in $H^2$ and approaches a new axial solution $V_j$ in
$H^s$. Comparison with
$V_{j-1}$ controls the change at previously occupied frequencies. Comparison with $V_j$ detects the new Fourier coefficient at
late times. The proposition states the estimates needed on this one
solution. Write
\[
 \|g\|_{q,\beta}=\sup_{t\ge0} e^{\beta t}\|g(t)\|_{H^q}.
\]
For an inhomogeneous Green operator we use the induced norm from
forcing to solution. Constants may depend on $s$.

\begin{proposition}\label{rate:prop:joint-insert}
Fix $s>2$. There are scales $B_j=B_*Q^j$,
$d_j=\lfloor\alpha\sqrt{B_j/\log B_j}\rfloor$, axial solutions
$V_j$ and complete two-dimensional solutions $W_j$, for $j\ge1$,
starting from a finite axial $V_0$. There are fixed constants
$b,C,C_0>0$ and frequencies $k_j=(1+N_j,0)$,
$B_j+O(1)\le\log N_j\le2B_j$, with the following estimates.
\begin{enumerate}
\item The two-dimensional orbit lands at $V_j$ with
$\|W_j-V_j\|_{q,\beta_j}\le e^{C_0 B_j}$ for $q=2,s$ and
$\beta_j\ge e^{-C_0B_j}$. On the same invariant lattice and time
half-line the linearized physical equation at $W_j$ has a right
inverse with zero prescribed stable trace, of norm at most
$e^{C_0B_j}$ in both of these spaces, with compatible inverses.
\item For $T_j=e^{(C_0+1)B_j}$ the mode $k_j$ satisfies
$\langle k_j\rangle^s |\widehat W_j(k_j,t)|
\ge0.9e^{b d_j}$ at every $t\ge T_j$. At
$S_j=e^{L_0 B_{j-1}}$, the entire $H^s$ Fourier tail at
$\langle k\rangle\ge K_j=e^{L_1B_{j-1}}$ costs at most
$Ce^{-d_j/C}$. One has $T_{j-1}\le S_j<T_j$ and
$\sup_{t\ge0}\|W_j(t)\|_{H^s}\le Ce^{C d_j}$.
\item $\sup_{t\ge0}\|W_j(t)-V_{j-1}(t)\|_{H^2}
\le e^{-D_hB_{j-1}}$. The constants $L_0>C_0+2$,
$L_1>2L_0+2$, and the required historical precision $D_h$ may
be selected successively after a coarse $C_0$ and before
$Q,\alpha,B_*$. In particular $D_h$ can pay every price of the
earlier stages in the weak history comparison.
\end{enumerate}
\end{proposition}
Sections~\ref{rate:sec:geom}--\ref{rate:sec:parameters} construct
the trajectories in the proposition. To pass to one initial datum,
we will use the linearized inverse at $W_j$ to restore the earlier
two-dimensional transfers. The next two propositions explain how to
control the resulting solution at times between transfers.

\subsection{High-frequency energy of an exact trajectory}
The full flow used here is global in $H^s$ for $s>1$.
Local existence, uniqueness and continuous dependence follow by applying
the contraction mapping theorem to Duhamel's formula in the algebra
$H^s(\T^2)$. The conserved mass and positive energy bound $\|u\|_{H^1}$
by a constant $M\ge1$. Splitting the Fourier sum at radius $K\ge2$ gives
$\|\widehat u\|_1\le C_s(M\sqrt{\log(2+K)}+
K^{1-s}\|u\|_{H^s})$ by Cauchy--Schwarz on both pieces. Taking
$K=2+(1+\|u\|_{H^s}/M)^{1/(s-1)}$, the high-frequency term is at most $M$. In the interaction
representation the cubic convolution has $h^s$ norm at most
$C_s\|\widehat u\|_1^2\|u\|_{H^s}$. Squaring the preceding
$\ell^1$ bound therefore gives, for $Y=2+\|u\|_{H^s}$,
$|Y'|\le C_{s,M}Y\log Y$ almost everywhere. In detail, the two pieces
of the Fourier sum are bounded by
\[
\begin{aligned}
 \sum_{\bk k\le K}|\widehat u(k)|
 &\le\Big(\sum_{\bk k\le K}\bk k^{-2}\Big)^{1/2}\|u\|_{H^1}
 \le CM\sqrt{\log(2+K)},\\
 \sum_{\bk k>K}|\widehat u(k)|
 &\le\Big(\sum_{\bk k>K}\bk k^{-2s}\Big)^{1/2}\|u\|_{H^s}
 \le C_sK^{1-s}\|u\|_{H^s},
\end{aligned}
\]
and with the stated $K$ one has $K^{1-s}\|u\|_{H^s}\le M$ and
$\log(2+K)\le C_s\log Y+C_{s,M}$. The convolution bound follows from
$\bk k^s\le3^s(\bk{k_1}^s+\bk{k_2}^s+\bk{k_3}^s)$ for
$k=k_1-k_2+k_3$ and Young's inequality. Thus
$|Y'|\le C_s\|\widehat u\|_1^2\|u\|_{H^s}\le C_{s,M}Y\log Y$, and
$\log Y$ grows at most exponentially on any bounded time interval.
This prevents a finite
lifespan endpoint in either time direction. Mass, energy and the
inequality follow first for smooth solutions and then by approximation
on each local existence interval. The scalar mass gauge preserves the
norms and has the same continuous dependence. This argument also supplies the Cauchy theory used in the preceding sections.

Fix a smooth nondecreasing cutoff $0\le\chi\le1$, identically zero on $[0,1]$ and identically one on $[2,\infty)$. Put
\[
 f_N(k)=\chi(\bk k/N)\bk k^{2s},\qquad
 E_N(t)=\sum_{k\in\Z^2} f_N(k)|\widehat u(k,t)|^2.
\]
\begin{proposition}\label{rate:p:flux}
Let $s\ge\sigma>1$, $\rho=1-1/\sigma$, and let $u\in C(I;H^s)$ solve \eqref{eq:nls} on a compact interval $I$. If $\sup_I\norm{u}_{H^\sigma}\le M$, then for sufficiently large $N$,
\[
 |E_N(t)-E_N(t')|\le C_{s,\sigma,M}
 (N^{-2}+|t-t'|N^{-\rho})\sup_{I_{t,t'}}\norm u_{H^s}^2,\qquad t,t'\in I .
\]
Here $I_{t,t'}$ is the closed interval between the indicated endpoints.
\end{proposition}
\begin{proof}
Write $a_k=e^{it|k|^2}\widehat u(k,t)$. Symmetrization writes $E'_N$ as a fixed real multiple of
\[
 \operatorname{Im}\sum_{k_1-k_2+k_3-k_4=0}
 [f_N(k_1)-f_N(k_2)+f_N(k_3)-f_N(k_4)]e^{-it\Phi}
 a_{k_1}\overline a_{k_2}a_{k_3}\overline a_{k_4}
\]
with $\Phi=|k_1|^2-|k_2|^2+|k_3|^2-|k_4|^2$. Indeed the equation
gives, writing $k_4$ for the output index,
\[
\begin{aligned}
 \dot a_{k_4}&=-i\sum_{k_1-k_2+k_3=k_4}
 e^{-it\Phi}a_{k_1}\overline a_{k_2}a_{k_3},\\
 E_N'&=2\operatorname{Re}\sum_{k_4}f_N(k_4)\overline a_{k_4}\dot a_{k_4}
 =2\operatorname{Im}\sum_{k_1-k_2+k_3-k_4=0}f_N(k_4)e^{-it\Phi}
 a_{k_1}\overline a_{k_2}a_{k_3}\overline a_{k_4}.
\end{aligned}
\]
Let $S_i$ denote the imaginary part of the same sum with $f_N(k_i)$
in place of $f_N(k_4)$. The summand is symmetric under
$k_1\leftrightarrow k_3$ and under $k_2\leftrightarrow k_4$, so
$S_1=S_3$ and $S_2=S_4$. Complex conjugation together with the
relabelling $(k_1,k_2,k_3,k_4)\mapsto(k_2,k_1,k_4,k_3)$ changes the
sign of $\Phi$ and gives $S_3=-S_4$. Hence
$E_N'=2S_4=-\tfrac12(S_1-S_2+S_3-S_4)$. In a nonzero term let
$K$ be the largest of the four brackets $\bk{k_i}$. Then $K\ge N$,
and momentum conservation makes the second-largest bracket at least
$cK$. The multiplier in square brackets thus costs at most the
product of the two largest brackets raised to $s$. With these two
factors in $\ell^2$ and the others in $\ell^1$, the quartic sum is
bounded by $C_s\|a\|_{h^s}^2\|a\|_1^2$. We gain a power of $N$
in one of the following three ways.

Put $L=N^{1/\sigma}$. If a third bracket is at least $L$, one of
the $\ell^1$ factors is bounded by
\[
 \sum_{\bk k\ge L}|a_k|
 \le \Big(\sum_{\bk k\ge L}\bk k^{-2\sigma}\Big)^{1/2}
       \|a\|_{h^\sigma}
 \le C_\sigma M L^{1-\sigma}=C_\sigma M N^{-\rho}.
\]
If both smaller brackets are below $L$, their multipliers vanish.
Suppose first that the two large entries have opposite signs in
$\Phi$. Their momenta differ by at most $2L$. Since
$|\nabla f_N(x)|\le C_s\bk x^{2s-1}$, the mean value theorem bounds
the remaining difference of multipliers by $C_sK^{2s-1}L$.
Relative to the preceding quartic bound this gains
$L/K\le N^{-\rho}$.

It remains to consider two large entries of the same sign and two
small entries below $L$. Here $|\Phi|\ge cK^2$ for $N$ sufficiently
large, since $L=o(N)$. Integrate $e^{-it\Phi}$ by parts in time
on $I_{t,t'}$. The factor $1/\Phi$ gives $O(N^{-2})$ at each
endpoint. In the integral the derivative falls on one of the four
amplitudes. The cubic equation gives
\[
 \norm{\dot a}_{\ell^1}\le C_\sigma M^3,\qquad
 \norm{\dot a}_{h^s}\le C_{s,\sigma}M^2\norm u_{H^s}.
\]
Use the first bound when the differentiated factor is one of the two
$\ell^1$ factors and the second for either weighted $\ell^2$ factor.
The integral from this case costs
$C_{s,\sigma,M}|t-t'|N^{-2}\sup_{I_{t,t'}}\|u\|_{H^s}^2$.
Combining it with the endpoint terms and the first two cases gives
the assertion. Approximation by smooth solutions on this compact
time interval justifies the integration by parts.
\end{proof}

\subsection{Passage to one initial datum}
For the remaining results $s>2$. Write $B_j=B_*Q^j$ with $Q\ge2$, $B_*>e$, and $d_j=\lfloor\alpha\sqrt{B_j/\log B_j}\rfloor$ for fixed $\alpha>0$. Choosing scalar phases at time zero, where appropriate, is equivalent to selecting representatives of the gauge-invariant estimates.

\begin{proposition}\label{rate:p:back}
Let $(U_j)$ be global $H^s$ solutions with $M=\sup_{j\ge0,\,t\ge0}\|U_j(t)\|_{H^2}<\infty$. Assume that $U_0$ is bounded in $H^s$ for $t\ge0$. Let $2K_{j-1}<K_j$ increase and assume
\begin{align*}
 &\sup_{t\ge0}\inf_\theta\norm{U_j(t)-e^{i\theta}U_{j-1}(t)}_{H^2}\le\delta_j,\\
 &E_{K_j}(U_j(S_j))^{1/2}\le c_j,\qquad
 \sup_{t\ge S_j}\norm{U_j(t)}_{H^s}\le C e^{A d_j}.
\end{align*}
Here $\ell_j=(2K_j)^{s-2}\delta_j$, $c_j$ and $\sqrt{\omega_j}\,e^{A'd_{j-1}}$ are summable for some $A'>A$, with $\omega_j=C_{s,2,M}(K_j^{-2}+S_j K_j^{-1/2})$. Let $T_j$ increase with $\log T_j\asymp B_j$ and $T_{j-1}\le S_j$. Assume also that for distinct frequencies $k_j$, some $a>0$ and every $t\ge T_j$,
\[
 \bk{k_j}^s |\widehat U_j(k_j,t)|\ge e^{a d_j},
 \qquad
 \bk{k_j}^{s-2}\sum_{i>j}\delta_i=o(e^{a d_j}).
\]
Then, after alignment at time zero, $U_j(0)$ converge strongly in $H^s$ to one datum, and its trajectory satisfies the two-sided estimate \eqref{eq:two-sided} for every sufficiently large time.
\end{proposition}
\begin{proof}
Set $M_j=\sup_{t\ge0}\|U_j(t)\|_{H^s}$ and
$M_j^e=\sup_{0\le t\le S_j}\|U_j(t)\|_{H^s}$.
The early supremum is finite by continuity and the late supremum
by hypothesis. For $0\le t\le S_j$, the multiplier
$(1-\chi(\bk k/K_j))^{1/2}$ is supported in $\bk k<2K_j$.
After optimizing over the scalar phase, the $H^s$ norm of this
part of $U_j(t)$ is at most $\|U_{j-1}(t)\|_{H^s}+\ell_j$.
Indeed, if $\Pi$ denotes this Fourier multiplier, then for every
real $\theta$
\[
 \|\Pi U_j\|_{H^s}\le\|\Pi e^{i\theta}U_{j-1}\|_{H^s}
   +(2K_j)^{s-2}\|U_j-e^{i\theta}U_{j-1}\|_{H^2}
 \le\|U_{j-1}\|_{H^s}+(2K_j)^{s-2}\|U_j-e^{i\theta}U_{j-1}\|_{H^2},
\]
because $\bk k^{s}\le(2K_j)^{s-2}\bk k^2$ on the support of $\Pi$.
Proposition~\ref{rate:p:flux} applied between $t$ and $S_j$ bounds
the complementary squared norm by $c_j^2+\omega_j(M_j^e)^2$; here
$\sigma=2$, $\rho=1/2$, $|t-S_j|\le S_j$ and
$E_{K_j}(U_j(S_j))\le c_j^2$.
These two squared norms add to $\|U_j(t)\|_{H^s}^2$. Hence
\begin{equation}\label{rate:eq:early-bound}
 \|U_j(t)\|_{H^s}^2\le
 (\|U_{j-1}(t)\|_{H^s}+\ell_j)^2+c_j^2+\omega_j(M_j^e)^2 .
\end{equation}
Take the early supremum and absorb the last term for large $j$.
Together with the late bound this gives
\[
 M_j\le\max\big\{C e^{A d_j},
       (1-\omega_j)^{-1/2}(M_{j-1}+\ell_j+c_j)\big\}.
\]
In detail, the early supremum of \eqref{rate:eq:early-bound} reads
$(1-\omega_j)(M_j^e)^2\le(M_{j-1}+\ell_j)^2+c_j^2
\le(M_{j-1}+\ell_j+c_j)^2$, and $M_j$ is the larger of $M_j^e$ and
the late supremum.
The product of the factors $(1-\omega_j)^{-1/2}$ converges after a
finite prefix. Since $\ell_j+c_j$ is summable and $d_j$ is eventually
nondecreasing, iteration gives $M_j\le C'e^{A'd_j}$. Explicitly,
if $\Pi_*$ bounds all the partial products of the factors
$(1-\omega_l)^{-1/2}$, $l>j_0$, then induction on $j$ gives
\[
 M_j\le\Pi_*\Bigl(M_{j_0}+\sum_{j_0<i\le j}(\ell_i+c_i)
          +\max_{j_0<i\le j}Ce^{Ad_i}\Bigr)\qquad(j>j_0).
\]
Applying the early bound once more gives
$M_j^e\le C''(1+e^{A'd_{j-1}})$. In particular
\[
 \gamma_j:=c_j+\sqrt{\omega_j}\,M_j^e
 \qquad\text{satisfies}\qquad \sum_j\gamma_j<\infty .
\]

At time zero, the $H^s$ tail of $U_j$ on $\bk k\ge 2K_j$ has norm
at most $\gamma_j$, since $\chi=1$ there. The corresponding tail
of $U_{j-1}$ has norm at most $\gamma_{j-1}$ by the same argument
at its own, smaller cutoff. Choose the initial phases recursively
to realize the $H^2$ comparison at time zero. On
$\bk k<2K_j$ the aligned difference then has $H^s$ norm at most
$\ell_j$. The full distance between successive aligned data is
therefore at most $\ell_j+\gamma_j+\gamma_{j-1}$ for large $j$.
It follows that the data converge strongly in $H^s$ to one datum.
Phase covariance and finite-time continuous dependence give convergence
of the aligned solutions to its solution $u$ at each fixed time.
The proof applies in the mass gauge as well, since all norm and
Fourier-modulus estimates are invariant under scalar phases.

Fix any $t$ with $T_j\le t<T_{j+1}$. For any phase $\theta$ the
difference of two Fourier moduli is bounded by the modulus of the
corresponding coefficient difference after multiplication by
$e^{i\theta}$. Summing the $H^2$ comparisons over future stages and
passing to the limit at this $t$ gives
\[
 \bk{k_j}^s|\widehat u(k_j,t)|
 \ge e^{a d_j}-\bk{k_j}^{s-2}\sum_{i>j}\delta_i
 \ge \tfrac12 e^{a d_j}
\]
for sufficiently large $j$. For the upper bound we can allow the
whole contribution $M_{j+1}$ of stage $j+1$. Every subsequent
index $i\ge j+2$ satisfies $t<T_{j+1}\le T_{i-1}\le S_i$.
Taking square roots in \eqref{rate:eq:early-bound} at this same time,
that is
$\|U_i(t)\|_{H^s}\le\|U_{i-1}(t)\|_{H^s}+\ell_i+c_i+\sqrt{\omega_i}\,M_i^e
=\|U_{i-1}(t)\|_{H^s}+\ell_i+\gamma_i$ for $t\le S_i$,
and summing from $i=j+2$ yields
\[
 \|u(t)\|_{H^s}
 \le M_{j+1}+\sum_{i\ge j+2}(\ell_i+\gamma_i)
 \le C e^{A'd_{j+1}}+O(1).
\]
Throughout the interval $\log t\asymp B_j$ and
$\log\log t=\log B_j+O(1)$. Since $B_{j+1}=QB_j$ with $Q$ fixed,
both $d_j$ and $d_{j+1}$ are comparable to
$\sqrt{\log t/\log\log t}$. The two estimates give the conclusion.
\end{proof}

\subsection{A family of right triangles with shared children}\label{rate:sec:geom}
Identify $\mathbb R^2$ with $\mathbb C$, with horizontal and vertical
coordinates $x,y$. In this section the absolute square is the sum
of the squares of the two coordinates. Its complexification is a
rational function, and does not involve complex conjugation of the
parameters. Put
\[
 P=(1-im)^{-1},\quad Q=1-P=-imP,\quad u=(1+it)^{-1},\quad
 g_{nk}=u^n(it)^k,\quad 0\le k\le n<d .
\]
Use $1$ as top and $Pg_{nk},Qg_{nk}$ as the primary vertices. A
primary $V$ at an internal level $0\le n\le d-2$ has children
$uV,(1-u)V$. Every primary at $1\le n\le d-1$ also has
a private companion and an axial vertex
\[
 J(V)=A(V)-V,\qquad A(V)=|V|^2/V_x .
\]
The designated triples are the top and its children, all primary
parents and their children, and the triples $(V,J(V),A(V))$. Every
triple has nonzero perpendicular inputs summing to its output. At
real parameter values outside finitely many algebraic sets, all
vertices are distinct. Distinct primary directions at different
levels differ by a nonconstant integer multiple of $\arctan t$
modulo $\pi/2$. Within a level, length ratios between parallel
primaries are signed powers of $t$ or such powers times $m$. The
companion of a vector of angle $\theta$ is $-i\tan(\theta)V$.
Indeed $Q/P=-im$ and $(1-u)/u=it$ are imaginary for real parameters,
so that the two inputs of every primary triple are perpendicular and
\[
 ug_{nk}=g_{n+1,k},\qquad(1-u)g_{nk}=g_{n+1,k+1}.
\]
Thus $g_{n+1,k+1}$ is a common child of $g_{nk}$ and $g_{n,k+1}$, and
level $n$ of each top branch has $n+1$ primaries. For
$V=(x,y)=|V|(\cos\theta,\sin\theta)$ one has $A(V)=(x^2+y^2)/x$ and
\[
 J(V)=\Bigl(\frac{y^2}{x},-y\Bigr)
     =|V|\tan\theta\,(\sin\theta,-\cos\theta)=-i\tan(\theta)V,
\]
which is perpendicular to $V$ and has $V+J(V)=A(V)$ on the axis.
Its moving pole cannot identically equal the scale of a primary.
If two such companions at the same level, or two axes, have the
same pole, their scale ratio is one of the preceding primary scale
ratios. Axes and nonaxes are distinct by their vertical components.
The two top inputs have no private triples.

\smallskip
\noindent\emph{Distinctness of the vertices in detail.} It suffices to
show that any two marked vertices are distinct rational functions of
\((m,t)\): two distinct rational functions agree only on a proper
algebraic subset, and there are finitely many pairs. We use the
notation \(\tau_n=\tan(n\arctan t)\) and
\(T_n=\tan(\arg(Pu^n))=(m-\tau_n)/(1+m\tau_n)\) of the reverse test
below; for \(n\ge1\) and generic \(t\), \(T_n\) has its only zero at
the real point \(m=\tau_n\ne0\) and its only pole at
\(m=-1/\tau_n\), and these points are distinct for distinct levels.
(a) The axes and the top have vanishing vertical component. A primary
\(Pg_{nk}\) has direction \(\arctan m-n\arctan t+k\pi/2\) modulo
\(\pi\), and \(Qg_{nk}\) this minus \(\pi/2\); as a function of \(m\)
neither is identically a multiple of \(\pi\). A
companion \(J(V)=(V_y^2/V_x,-V_y)\) has vertical component
\(-V_y\not\equiv0\). So axes and the top differ from all other vertices.
(b) For two primaries of one branch,
\(Pg_{nk}/Pg_{n'k'}=u^{n-n'}(it)^{k-k'}\) has modulus
\((1+t^2)^{(n'-n)/2}t^{k-k'}\), which is identically one only if
\(n=n'\) and \(k=k'\). For two branches the ratio contains the factor
\(-im\), so its modulus depends on \(m\).
(c) A companion of a primary \(V\) of level \(n\ge1\) is
\(J(V)=-i\tan(\theta_V)V\) with \(\tan\theta_V\in\{T_n,-T_n^{-1}\}\).
Its modulus \(|\tan\theta_V||V|\) has a pole at the real point
\(m=-1/\tau_n\), respectively \(m=\tau_n\), where every primary is finite
and nonzero, since primary moduli are
\((1+m^2)^{-1/2}(1+t^2)^{-n/2}t^k\), possibly times \(|m|\). Hence a
companion differs from every primary.
(d) Companions at different levels have poles at different real points.
At one level, if \(V'=cV\) with real \(c\), then \(J(V')=cJ(V)\) and
\(c\not\equiv1\) by (b) unless \(V'=V\); if \(V'\perp V\), then \(J(V)\)
is parallel to \(iV\) and \(J(V')\) to \(iV'\parallel V\), so the two
nonzero companions are perpendicular.
(e) An axis \(A(V)=|V|^2/V_x\) has a real pole exactly where \(V\) is
vertical, and \(A(\lambda V)=\lambda A(V)\) for real \(\lambda\). If
\(A(V)\equiv A(V')\), the poles coincide, so \(V,V'\) are vertical at
the same \(m\) for all \(t\); directions at different levels differ by a
nonzero multiple of \(\arctan t\) modulo \(\pi/2\), so \(V,V'\) are
parallel at one level, \(V'=\lambda V\), and \(\lambda\equiv1\), whence
\(V=V'\) by (b). The top has no pole and differs from every axis.

Write $D(v,w)$ for the pair of rational functions obtained by
taking the signed sums (inputs positive, output negative) of
$(x^{-2},|v|^{-2})$ on $v,w,v+w$. For positive constants $A_0,B_0$
the intended scalar defect is $A_0D_1-B_0D_2$. In the next lemma
polynomial nonidentity is asserted before the slopes are specialized.

\begin{lemma}\label{rate:l:geom}
At fixed depth $d$, the designated component pairs
$D_e\in\mathbb R(m,t)^2$ are linearly independent over $\mathbb R$.
Spans and identities of pairs are taken before specialization. For
a fixed finite list of tests, \emph{generically} means at admissible
real $(m,t)$ outside the poles and the zero sets of the nonzero
numerators in those tests.
A nontrivial reverse rectangle test with an output
and one input marked yields a designated triple whenever it holds
identically. Every extra perpendicular pair of marked inputs has a
component pair different from every constant linear combination
of the designated pairs. For one marked output $h$ and one designated
pair $D_e$, a perpendicular split of $h$ with this same component pair
(the inputs may be unmarked) is generically impossible unless $h$ is the output of $e$. For its own output the
inputs are determined up to exchange and, for a horizontal output,
reflection across the horizontal line.
\end{lemma}
\begin{proof}
For the private row at $v=(x,y)$, $S=x^2+y^2$, the two defects are
\[
 x^{-2}+x^2/y^4-x^2/S^2,\qquad y^{-2}-x^2/S^2.
\]
Indeed $J(v)=(y^2/x,-y)$ has $|J(v)|^2=y^2S/x^2$ and $A(v)=(S/x,0)$,
so the first component is $x^{-2}+(y^2/x)^{-2}-(S/x)^{-2}$ and the
second is $S^{-1}+x^2/(y^2S)-x^2/S^2=y^{-2}-x^2/S^2$.
At the real moving pole $y=0$ the first has order four. Fix a level
and direction, factoring out its common defect. The coefficient
from the two top subtrees is $R(t)+m^{-2}S_1(t)$, with distinct
Laurent monomials in $t$ within each polynomial. At generic fixed
$t$ all other private rows, as well as all primary rows, have pole
order at most two there. Thus both the value and the first
$m$-derivative of this coefficient vanish at a nonzero finite
$m_0(t)$, implying $R=S_1=0$. Varying $t$ gives zero for every
coefficient. In more detail, both components of a defect pair are
homogeneous of degree $-2$ under real dilations of the triple. At a
fixed level $n$, the primaries $Pg_{nk}$ with a fixed parity of $k$
are the real multiples $\pm t^{k-k_0}$ of one of them, and the
primaries $Qg_{nk}=-imPg_{nk}$ with the opposite parity are its
multiples $\pm mt^{k-k_0}$. Their private pairs are therefore
$t^{-2(k-k_0)}$, respectively $m^{-2}t^{-2(k-k_0)}$, times one private
pair, which is the meaning of the coefficient $R(t)+m^{-2}S_1(t)$ in
a linear relation. The common direction is horizontal at
$m_0(t)=\tan(n\arctan t-k_0\pi/2)$, where $y$ has a simple zero as a
function of $m$. The pole of order four of $x^2/y^4$ can be reduced
to the order two of the remaining rows only if the coefficient
vanishes to second order at $m_0(t)$, and
$\partial_m(R+m^{-2}S_1)=-2m^{-3}S_1$. After eliminating all private rows, the vertical
orientations at the deepest remaining primary layer are double
poles of its first component, with all other directions regular.
The same argument imposes vanishing of every individual formal
vertex coefficient. The children of parent coefficients $c_j$ have
coefficients $c_j+c_{j-1}$, with $c_{-1}=0$. Thus every parent
coefficient vanishes. Descending leaves only the nonzero top row.

In detail, suppose that the private rows and the rows of all parents
at levels greater than $n$ have been removed, and let $c^P_j,c^Q_j$,
$0\le j\le n$, be the coefficients of the rows with parents
$Pg_{nj},Qg_{nj}$; put $c_{-1}=c_{n+1}=0$. A primary of level $n+1$
occurs in the remaining rows only as a child: $g_{n+1,k}$ is the child
$ug_{nk}$ of $g_{nk}$ and the child $(1-u)g_{n,k-1}$ of $g_{n,k-1}$.
Let $G'=Pu^{n+1}$, and at fixed generic $t$ let $m_1(t)$ be a value of
$m$ at which $G'$ is vertical; it is finite and nonzero. The primaries
of level $n+1$ parallel to $G'$ are $Pg_{n+1,k}=\pm t^kG'$ with $k$
even and $Qg_{n+1,k}=\pm mt^kG'$ with $k$ odd. Every other vertex of a
remaining row is perpendicular to $G'$ or has a direction differing
from that of $G'$ by a nonzero integer multiple of $\arctan t$ modulo
$\pi/2$, so its abscissa does not vanish at $m_1(t)$; the top has
abscissa one. The first
component of the relation is therefore
\[
\begin{gathered}
 \bigl(\tilde R(t)+m^{-2}\tilde S(t)\bigr)(G'_x)^{-2}
 +(\text{a function regular at }m_1(t)),\\
 \tilde R=\sum_{k\ {\rm even}}\frac{c^P_k+c^P_{k-1}}{t^{2k}},\qquad
 \tilde S=\sum_{k\ {\rm odd}}\frac{c^Q_k+c^Q_{k-1}}{t^{2k}} .
\end{gathered}
\]
Since $G'_x$ has a simple zero at $m_1(t)$ and no other term has a
pole there, both the value and the first $m$-derivative of
$\tilde R+m^{-2}\tilde S$ vanish at $m_1(t)$. As before this gives
$\tilde S=0$ and then $\tilde R=0$, identically in $t$. The vertical
orientation of $iG'$ gives the same conclusion with the two parities
exchanged. Hence $c_k+c_{k-1}=0$ for $0\le k\le n+1$ in each branch,
and successively $c_0=c_1=\dots=c_n=0$. In particular the first
components of the designated pairs are already linearly independent.

The reverse test for a nonaxis output $h$ and nonaxis input $v$
is $\operatorname{Re}(h/v)=1$, because
$v\cdot(h-v)=\operatorname{Re}(\bar vh)-|v|^2$. At different levels the angle
difference depends only on $t$. A companion contributes a factor
$T_n=\tan(\arg(Pu^n))$ or $T_n^{-1}$, with a distinct real
moving pole or zero at each strictly positive level. It cannot
cancel with a different level or a power of $m$. Pure primaries
in different top branches also leave a power of $m$. Within a
branch the remaining ratio is $(it)^a(1+it)^b$. If $b<0$, a
nonzero constant real part would require $a\le0$ at zero and
$a\ge-b$ at infinity. If $b\ge0$, expanding the binomial, the
constant real part one occurs only for $(b,a)=(0,0),(1,0),(1,-1)$:
for $b=2$ the sole odd term has coefficient two, and for $b\ge3$
each parity has at least two nonzero terms. These give coincidence
or the two assigned children. Explicitly, for $h=g_{nk}$ and
$v=g_{n'k'}$ in one branch, $h/v=(it)^a(1+it)^b$ with $a=k-k'$ and
$b=n'-n$, and
\[
 \operatorname{Re}\bigl[(it)^a(1+it)^b\bigr]
 =\sum_{\substack{0\le j\le b\\ j\equiv a\ (2)}}
   \binom bj(-1)^{(a+j)/2}t^{a+j}\qquad(b\ge0).
\]
The exponents are distinct, so this equals one identically only if
the sum has the single term $j=-a$ with $\binom b{-a}=1$. The three
solutions correspond to $v=h$, $v=uh$ and $v=(1-u)h$. Within a single level the directions
are parallel or perpendicular, giving no additional reverse.
An axial output must be $A(v)=A(J(v))$. An axial input would have
to equal the horizontal coordinate of the nonaxis output; the
latter has a pole at $m=i$ or $-i$ from its primary part, whereas
the axes have only real poles. For $P,Q$ the axial output is $1$.

\smallskip
\noindent\emph{The reverse test in detail.} Identities are identities
of rational functions of \((m,t)\). To locate poles we fix a generic
real \(t>0\) and regard each coordinate as a rational function of a
complex variable \(m\); the complexification replaces \(\overline z\)
by the expression \(z^*\) obtained by changing \(i\) into \(-i\) in the
explicit formula, so that for instance \(P^*=(1+im)^{-1}\) and
\(|V|^2=VV^*\). Since \(\arg P=\arctan m\), \(\arg u=-\arctan t\) and
\(\arg(it)=\pi/2\),
\[
 \arg(Pg_{nk})\equiv\arctan m-n\arctan t+k\tfrac\pi2,\qquad
 \arg(Qg_{nk})\equiv\arg(Pg_{nk})-\tfrac\pi2\pmod\pi .
\]
Put \(\tau_n=\tan(n\arctan t)\), finite and nonzero at generic \(t\) for
\(n\ge1\), and
\[
 T_n\coloneqq\tan\bigl(\arg(Pu^n)\bigr)=\frac{m-\tau_n}{1+m\tau_n},
 \qquad T_0=m .
\]
A marked primary \(V\) of level \(n\) has \(\tan(\arg V)\in\{T_n,-T_n^{-1}\}\),
so its companion is \(J(V)=-iT_nV\) or \(J(V)=iT_n^{-1}V\). As a
M\"obius function of \(m\), \(T_n\) has its only zero at \(m=\tau_n\) and its
only pole at \(m=-1/\tau_n\). For \(n\ne n'\) (both \(\ge1\)) the four real
points \(\tau_n,-1/\tau_n,\tau_{n'},-1/\tau_{n'}\) are distinct at generic
\(t\), since an identity \(\tau_n\equiv\tau_{n'}\) or \(\tau_n\tau_{n'}\equiv-1\)
would force \(n\arctan t\equiv n'\arctan t\) modulo \(\pi/2\) for all \(t\).
Moreover \(T_n(\pm i)=\pm i\), because
\(T_n-i=(1-i\tau_n)(m-i)/(1+m\tau_n)\), and similarly at \(-i\).

Let \(h,v\) be marked nonaxis vertices, \(h=c_hh_0\), \(v=c_vv_0\), with
primaries \(h_0,v_0\) and \(c_h,c_v\in\{1,-iT_n,iT_n^{-1}\}\), \(n\) being the
level of the respective primary. Then
\[
 h_0/v_0=u^{n_h-n_v}(it)^{k_h-k_v}(-im)^\gamma,\qquad \gamma\in\{0,\pm1\},
\]
and \(\gamma\ne0\) exactly when \(h_0,v_0\) lie in different top branches.
Hence at fixed \(t\) one has \(h/v=\kappa\,\rho(m)\), where \(\kappa=\kappa(t)\)
is a complex constant and \(\rho(m)=m^\gamma\prod T_n^{\pm1}\) is a real
rational function, the product running over the companion factors
present. For real \(m\), \(\operatorname{Re}(h/v)=\operatorname{Re}(\kappa)\rho(m)\),
so an identity \(\operatorname{Re}(h/v)\equiv1\) forces \(\rho\) to be a
nonzero constant. By the locations of the zeros and poles of \(m\) and of
the \(T_n\), this happens only if (i) neither \(h\) nor \(v\) is a companion
and \(\gamma=0\), or (ii) both are companions at the same level whose
factors cancel, and \(\gamma=0\). In case (ii) the factors cancel only if
\(c_h=c_v\) (the other possibility gives \(\rho=T_n^{\pm2}\)); then \(h_0,v_0\)
are parallel, \(h/v=h_0/v_0\) is real, and \(\operatorname{Re}(h/v)=1\)
means \(h=v\). Case (i) is the binomial computation above.

For an axial output \(h\) (an axis or the top) and a nonaxis input \(v\),
\(v\cdot(h-v)=hv_x-|v|^2\), so the test holds exactly when
\(h=|v|^2/v_x=A(v)\). For a primary \(v\) of level \(\ge1\) this is its own
axis; \(A(P)=A(Q)=1\), because \(|P|^2=P_x=(1+m^2)^{-1}\) and
\(|Q|^2=Q_x=m^2(1+m^2)^{-1}\); and for a companion, with \(w=(x,y)\) and
\(S=x^2+y^2\), \(A(J(w))=(y^2S/x^2)/(y^2/x)=S/x=A(w)\). In each case the
triple is designated, provided distinct marked horizontal vertices
are distinct functions. This holds: \(A(\lambda V)=\lambda A(V)\) for real
\(\lambda\), \(A(V)\) has a real pole exactly where \(V\) is vertical, and the
top has none, so two axes agree identically only if their primaries
are parallel at one level with equal scale, hence equal. For an axial
input \(v\) and a nonaxis output \(h\) the test reads \(v=h_x\). Write a
marked primary as \(V=Pw\), where \(w\) is \(g_{nk}\) or \(-img_{nk}\). Then
\[
 V_x=\frac{w(1+im)+w^*(1-im)}{2(1+m^2)},\qquad
 V_y=\frac{w(1+im)-w^*(1-im)}{2i(1+m^2)},
\]
\[
 A(V)=\frac{|V|^2}{V_x}=\frac{2ww^*}{w(1+im)+w^*(1-im)} .
\]
At \(m=i\) the numerators of \(V_x\) and \(V_y\) equal \(2w^*\) and \(-2w^*\),
which are nonzero, so \(V_x\) and \(J(V)_x=V_y^2/V_x\) have simple poles
there (the latter behaves like \(-w^*/(1+m^2)\)). The denominator of \(A(V)\)
is, up to a factor \(m\) when \(w=-img_{nk}\), a polynomial in \(m\) with real
coefficients of degree one, which vanishes only at a real \(m\); the
top is constant. Hence \(v\equiv h_x\) is impossible for marked \(v,h\).

Two marked nonaxis inputs can be generically perpendicular only at
the same level. In the basis $G=Pu^n,iG$ their signed relative
scale is
\[
 L=\pm t^{l-k}m^\gamma T_n^\delta,
 \qquad\gamma=0,1,-1,\quad\delta=0,1,2 .
\]
Indeed a companion on the first basis direction has a factor
$T_n^{-1}$ and a companion on the second direction a factor $T_n$.
If $L$ depends on $m$, any preimage $L=\pm i$ away from $m=\pm i$
is a new pole of the radial component of the output. All individual
lengths and abscissae are regular and nonzero at that preimage:
their exceptional points are real or $\pm i$. The output is a
nonzero isotropic vector, whose horizontal coordinate cannot vanish.
For $l\ne k$ such a preimage exists since $T_n$ fixes $\pm i$ and
$t^{l-k}$ is nonconstant. For $l=k$ and degree at least two, the
possible preimages at $\pm i$ are simple at generic $t$; the
logarithmic derivative there is proportional to
$\gamma+\delta T_n'(i)$, which cannot vanish identically at $n\ge1$
unless $L$ is constant. Every nonconstant factor has its zeros and
poles on the real projective line. At $l=k$ the degree-one possibilities
are a primary with its own companion or $Pg_{nk},Qg_{nk}$:
the four vectors in unit-$G$ coordinates are two primaries on
opposite axes and their two perpendicular companions. Opposite
branch with a single companion is parallel at equal $k$, and the
two compensating companions involve $T_n^2$ and a power of $m$.
For $Pg_{nk}+Qg_{nk}=g_{nk}$ the output is independent of $m$.
In any representation of this defect by designated rows, the fourth-order
pole comparison at the moving horizontal orientations first eliminates
the private rows, since the first component on the left has poles of order at
most two there. Comparing the principal parts at the double poles then
fixes the primary incidences of the two inputs at all levels $n\ge1$.
Project a node from each top branch to its leaf
polynomial $X^{n-k}Y^k(X+Y)^{d-1-n}$. Every internal primary incidence
maps to zero, and the top incidence maps to a scalar multiple of
$(X+Y)^{d-1}$ in each branch. The singled-out node has a different
polynomial for $n\ge1$.

\smallskip
\noindent\emph{Perpendicular pairs in detail.} In the notation of the
reverse test, a marked nonaxis vertex of level \(n\) has direction
\(\arctan m-n\arctan t+c\pi/2\) modulo \(\pi\), with an integer \(c\). Two
such vertices of levels \(n,n'\) are perpendicular for all \(t\) only if
\((n'-n)\arctan t+(c-c')\pi/2\equiv\pi/2\pmod\pi\) for all \(t\); letting
\(t\downarrow0\) shows that \(c-c'\) is odd, and then
\((n'-n)\arctan t\in\pi\mathbb Z\) for small \(t\) gives \(n=n'\). An axial
vertex is perpendicular to a nonaxis vertex only where the latter is
vertical, at isolated values of \(m\). At level \(n\) put \(G=Pu^n\). Then
\(Pg_{nk}=(it)^kG\) is \(\pm t^kG\) for even \(k\) and \(\pm t^k\,iG\) for odd
\(k\), while \(Qg_{nk}=-im(it)^kG\) is \(\pm mt^k\,iG\) for even \(k\) and
\(\pm mt^kG\) for odd \(k\). For companions,
\(J(\sigma G)=-\sigma T_n\,iG\) and \(J(\sigma\,iG)=-\sigma T_n^{-1}G\). An
identically perpendicular pair is therefore \(v=\sigma_1G\), \(w=\sigma_2\,iG\)
with real scales, and \(L=\sigma_2/\sigma_1=\pm t^{l-k}m^\gamma T_n^\delta\),
where \(k,l\) are the indices of the underlying primaries,
\(\gamma\in\{0,\pm1\}\) records their branches, and \(\delta\in\{0,1,2\}\) is
the number of companions in the pair. With \(h=v+w=\sigma_1G(1+iL)\) and
\(|G|^2=GG^*=(1+m^2)^{-1}(1+t^2)^{-n}\),
\[
 |v|^{-2}+|w|^{-2}-|h|^{-2}
 =\frac1{\sigma_1^2|G|^2}\Bigl(1+L^{-2}-\frac1{1+L^2}\Bigr).
\]
Every second component of a designated pair is regular at each
\(m\notin\mathbb R\cup\{\pm i\}\): \(|Pg_{nk}|^{-2}=(1+m^2)(1+t^2)^nt^{-2k}\) is a
polynomial in \(m\); the factors \(m^{-2}\) and \(T_n^{\mp2}\) of the
\(Q\)-vertices and companions have real exceptional points; and
\(|A(V)|^{-2}=\cos^2(\arg V)\,|V|^{-2}=|V|^{-2}(1+\tan^2(\arg V))^{-1}\) is
singular only where \(T_n=\pm i\), i.e.\ at \(m=\pm i\). The same holds for
\(\sigma_1^{-2}|G|^{-2}\) and for \(1+L^{-2}\) at points where \(L=\pm i\). Hence,
if \(L(m_*)=\pm i\) for some \(m_*\notin\mathbb R\cup\{\pm i\}\), the second
component of \(D(v,w)\) has a pole at \(m_*\) and \(D(v,w)\) is not a
combination of designated pairs. Suppose \(L\) depends on \(m\). It is a
real rational function of \(m\), so all solutions of \(L=\pm i\) are
nonreal. If \(l\ne k\), then \(|L(\pm i)|=t^{l-k}\ne1\) at generic \(t\), since
\(|T_n(\pm i)|=1\), so no solution lies at \(\pm i\). If \(l=k\) and
\(L=\pm m^\gamma T_n^\delta\) has degree at least two, then
\(L(\pm i)=\pm(\pm i)^{\gamma+\delta}\); for even \(\gamma+\delta\) this is real,
and for odd \(\gamma+\delta\) one has \(L(-i)=-L(i)\), so at most one of
\(\pm i\) solves \(L=i\), and it does so simply. Indeed
\((\log L)'=\gamma/m+\delta T_n'/T_n\) and
\(T_n'(m)=(1+\tau_n^2)/(1+m\tau_n)^2\), so at \(m=i\) the logarithmic
derivative is \(-i(\gamma+\delta T_n'(i))\) with
\(T_n'(i)=(1-i\tau_n)/(1+i\tau_n)\) of modulus one and not real for
\(\tau_n\notin\{0,\infty\}\); it vanishes only if \(\gamma=\delta=0\). Thus
\(L=i\) has at least \(\deg L-1\ge1\) solutions away from \(\pm i\). The
remaining cases with \(L\) depending on \(m\) have \(l=k\) and degree one,
i.e.\ \((\gamma,\delta)\in\{(\pm1,0),(0,1)\}\), because
\(\deg(m^{\pm1}T_n)=\deg(T_n^2)=2\) and \(\deg(m^{\pm1}T_n^2)=3\). For
\((0,1)\) the pair consists of a primary and the companion of a primary
of the same branch, level and index, i.e.\ of the same primary: a
designated private pair. For \((\pm1,0)\) it is \(\{Pg_{nk},Qg_{nk}\}\), with
output \(g_{nk}\); for \(n=0\) this is the top pair, and for \(n\ge1\) it is
excluded by the leaf polynomials below. If \(L\) does not depend on \(m\),
then \(\gamma=\delta=0\), and the case is treated after the leaf
polynomials.
In detail, suppose that
$D(Pg_{nk},Qg_{nk})=\sum_ec_eD_e$ with $n\ge1$. The first component of
the left side is
$(Pg_{nk})_x^{-2}+(Qg_{nk})_x^{-2}-(g_{nk})_x^{-2}$, and its last term
does not depend on $m$. It has no pole of order four at a horizontal
orientation of a level, so the private coefficients vanish as in the
first part of the proof. For a primary $Pg_{n'k'}$ of level
$n'\ge1$ let
\[
 \iota(n',k')=c^P_{n'-1,k'}+c^P_{n'-1,k'-1}-c^P_{n'k'}
\]
be the coefficient with which $(Pg_{n'k'})_x^{-2}$ occurs in
$\sum_ec_eD_e$; here $c^P_{n'j}$ is the coefficient of the row with
parent $Pg_{n'j}$, and coefficients with indices outside
$0\le j\le n'\le d-2$ are zero. Comparing the double and the simple
poles at the two vertical orientations of each level $n'\ge1$, as in
the first part, gives $\iota(n,k)=1$ and $\iota(n',k')=0$ for every
other node of the branch at a level $n'\ge1$. Put
$\phi(n',k')=X^{n'-k'}Y^{k'}(X+Y)^{d-1-n'}$. Then
$\phi(n'+1,j)+\phi(n'+1,j+1)-\phi(n',j)=0$ for every parent. Summing
this identity with the weights $c^P_{n'j}$ over the rows of the branch
and collecting the terms at each node gives
\[
\begin{gathered}
 \sum_{n'\ge1}\sum_{k'}\iota(n',k')\,\phi(n',k')
   -c^P_{00}\,\phi(0,0)=0,\\
 \text{that is}\qquad
 X^{n-k}Y^k(X+Y)^{d-1-n}=c^P_{00}(X+Y)^{d-1}.
\end{gathered}
\]
This is false for $n\ge1$, because $X^{n-k}Y^k$ is not a multiple of
$(X+Y)^n$. The top row does not enter, since its inputs have level
zero.
In the $m$-independent case the inputs are
pure same-branch primaries, with an odd difference $j=|l-k|$.
The case $j=1$ is designated. For $j\ge3$ the direction of the
output involves $\arg(1\pm it^j)$ besides the level direction.
Its expansion at zero excludes any equality modulo $\pi/2$ to
an integer multiple of $\arctan t$. This gives a new moving
horizontal reciprocal-square pole.
Explicitly, let \(v=Pg_{nk}\) and \(w=Pg_{n,k+j}\) with odd \(j\ge3\); the other
branch is identical after the factor \(-im\). Then
\(h=v+w=Pg_{nk}(1+(it)^j)\) with \((it)^j=\pm it^j\), so
\(\arg h\equiv\arg v\pm\arctan(t^j)\pmod\pi\). As a function of \(m\),
\(h_x=\operatorname{Re}(Pc)=(\operatorname{Re}c-m\operatorname{Im}c)/(1+m^2)\)
with \(c=g_{nk}(1+(it)^j)\) independent of \(m\), so \(h_x^{-2}\) has a double
pole at the real point \(m_h(t)\) where \(h\) is vertical. At \(m_h(t)\)
neither \(v\) nor \(w\) is vertical: their directions differ from that of
\(h\) by \(\mp\arctan(t^j)\) and \(\pi/2\mp\arctan(t^j)\), which are not
multiples of \(\pi/2\) for \(0<t<\infty\). Hence the first component of
\(D(v,w)\) has a pole at \(m_h(t)\). For finite real \(m\ne0\), every primary \(V\) has positive squared
length. Hence \(A(V)_x^{-2}=V_x^2/|V|^4\) is regular there, and the
top contributes the constant one. The remaining terms \(V_x^{-2}\)
and \(J(V)_x^{-2}=V_x^2/V_y^4\) can have poles only where \(V\) is
vertical or horizontal. Thus a designated first component can have
a pole at a finite real \(m\ne0\) only if
\(\arctan m\equiv n'\arctan t\pmod{\pi/2}\) for some \(0\le n'<d\).
At \(m=0\) the \(Q\)-branch degenerates, as seen from
\(A(Qu)_x^{-2}=(m-t)^2/m^2\).

The vertical orientation of \(h\) satisfies
\(\arctan m_h(t)\equiv n\arctan t\mp\arctan(t^j)\pmod{\pi/2}\).
An identity with one of the preceding directions would give
\(\pm\arctan(t^j)\equiv(n-n')\arctan t\pmod{\pi/2}\).
Near \(t=0\) both sides tend to zero. Comparing their first
derivatives gives \(n=n'\), since \(j\ge3\), and then the nonzero
term of order \(t^j\) gives a contradiction. Taking \(n'=0\)
shows in particular that the angle defining this moving pole differs
from \(0\pmod{\pi/2}\) for generic \(t\), so \(m_h(t)\) is finite and nonzero.
At this point every designated first component is regular,
whereas that of \(D(v,w)\) has a double pole.

For the last assertion, let $U,V$ denote a proposed pair, $h$ its
output, $z=h_y/h_x$, $g=(1+|h|^2 V)^{-1}$, and
$p=k_x l_x/h_x^2$, where $k+l=h$, $k\perp l$. The two equations
to be satisfied are
\begin{equation}\label{rate:eq:extgeo}
 (p-(1-z^2)g)^2=z^2g(1-4g),\qquad
 (1+h_x^2 U)p^2+2p-1=0.
\end{equation}
To derive the equations, write $a=|k|^2/|h|^2$. Orthogonality
gives $|l|^2=(1-a)|h|^2$ and hence
$1+|h|^2V=1/a+1/(1-a)=1/[a(1-a)]$.
Let $x=k_x/h_x$, so that $p=x(1-x)$. The decomposition of $k$
along $h$ and its perpendicular gives $x=a\pm z\sqrt g$. Thus
$p=(1-z^2)g\pm z\sqrt g\,(1-2a)$, with $(1-2a)^2=1-4g$.
The second equation follows from
$1+h_x^2 U=x^{-2}+(1-x)^{-2}=(1-2p)/p^2$.
For $z\ne0$ and $g>1/4$ the first equation has no real root.
When $2+h_x^2U>0$ the roots of the second are real
(its discriminant is $4(2+h_x^2 U)$, and the linear case is immediate).
There is then no common complex root. For $z=0$ the test is
\[
 |h|^4 V^2-|h|^2 U-2=0.
\]
One may use a real orientation to disprove the identically
vanishing complex resultant. For $0<g\le1/4$ all roots of the
first equation are real and arise from the two possible reflected
splits of the prescribed length fractions.

The case analysis below has the following logic. For a fixed marked
$h$ and a fixed designated $e\,$ with $(U,V)=D_e$, a perpendicular
split of $h$ with component pair $D_e$ exists for generic parameters
only if the resultant in $p$ of the two equations
\eqref{rate:eq:extgeo}, respectively the horizontal test when $z=0$,
vanishes identically as a rational function of $(m,t)$. It therefore
suffices to exhibit, in each case, real parameters at which there is
no common root; these are found near a degenerate orientation or
scale. Three kinds of contradiction occur.
\begin{itemize}
\item[(a)] $z$ is finite and nonzero, $g\notin(0,1/4]$ and
$2+h_x^2U>0$. Then $z^2g(1-4g)<0$, so the first equation has no real
root, while both roots of the second are real. Note that $g<0$, which
occurs when $1+|h|^2V<0$, is covered together with $g>1/4$.
\item[(b)] $z=0$ and $|h|^4V^2-|h|^2U-2\ne0$.
\item[(c)] $0<g\le1/4$, and the finitely many real splits allowed by
the first equation have a horizontal defect whose size or whose pole
structure, under rotation of the branch, differs from that of $U$.
\end{itemize}
In the order of the text: a target in the opposite top branch, and
the top row against a strict descendant, give (a) with $g\to1$ and
$h_x^2U\to0$, or (b) with limiting value $-2$ for an axial $h$; the
top row against a top child gives (a) with $g\to1/2$; a primary
target against a primary $h$ of its branch gives (a) near a vertical
child, or (c) by the signs of the three double poles; against
$h=J(w)$ it gives (a) with $g\to1$; against an axial $h$ it gives (b)
with $U\to+\infty$. A private target against a primary of another
level gives (c) by comparing $U/V$; against a primary of the same
level, against $v$ itself, against a companion and against $J(v)$ it
gives (a) with $g\to1$, $1/2$, $1$ and $1/2$ respectively; against an
axial $h$ it gives (b), by the uniqueness of the positive root.

For opposite top branches let the branch of $h$ shrink by varying
$m$, with the target bounded at generic $t$. This sends $h$ and
$h_x^2U$ to zero. Against the top row itself each strict
descendant shrinks as $t\to0$ for $k>0$, and as $t\to\infty$
for $k=0$, at generic fixed $m$; the limiting direction at
finite $d$ stays away from the coordinate axes. For a top child
against the top row, shrinking that child gives $|h|^2V\to1$
and $U>0$. These are exclusions by the preceding equations.

Consider a primary target in the same branch. Divide the lengths
by the common branch radius so the configuration rotates as $m$
varies. For a primary $h$, both $|h|$ and $V$ are now constant.
If $g>1/4$ approach a vertical target child to make $U>0$.
Otherwise the finitely many relative splitting angles are fixed.
The horizontal defect of each such split has three distinct double
poles with leading signs $+,+,-$ under rotation. Its unique
negative pole gives the output direction and, from its coefficient,
the output length. This requires $h$ to be the target parent up
to sign, which is absent for distinct marked primaries at generic
parameters. For $h=J(w)$ rotate $w$ to horizontal, so $h_x$
vanishes quadratically, $h\to0$, $V$ stays bounded and $U$ has
at most a double pole. For axial $h=A(w)$ try the two perpendicular
target-child vertical orientations. At one of them $|h|$ stays
finite and nonzero, $V$ stays bounded, and $U\to+\infty$,
violating the horizontal test. This also covers $h=1$: its
normalized value can diverge at at most one of the two orientations.

For a private target based at $v$, rotate the normalized branch
to make $v$ horizontal. Now $V\to\infty$ and $U/V\to\infty$.
A primary $h$ at a different level has limiting direction away
from both axes for generic $t$. Each split with the required $V$
has a short leg pointing asymptotically perpendicular to $h$;
its $U/V\to |h|^2/h_y^2$. At the same level, if $k'\ne k$,
shrink $|h|/|v|=|t|^{k'-k}$ while holding a generic direction
of $v$ fixed. Then $|h|^2V\to0$ and $h_x^2 U$ stays bounded
and nonnegative. For the same primary, rotate $v$ to vertical to
obtain $|h|^2V\to1$ and $U>0$. If $h=J(w)$ make $w$ horizontal;
unless $w$ is parallel to $v$, $V$ stays finite, $U>0$, and
$h\to0$. Parallel companions have the same power-of-$t$ length
ratio just used, or, at the identical node, $|J(v)|^2 V\to1$
when $v$ becomes horizontal. The horizontal test has just one
positive solution for $|h|^2$ when $V>0$, and it is already
$|A(v)|^2$. The axis magnitudes are generically distinct. Finally
if a marked nonhorizontal output is the correct one, the two $p$
values for a reflected split can both solve the second equation
only if the two candidates coincide or $g=(1+z^2)^{-2}$:
compare the sum and product of the roots using
\eqref{rate:eq:extgeo}. Indeed the first equation is
$p^2-2(1-z^2)gp+(1-z^2)^2g^2-z^2g(1-4g)=0$, and the roots of the
second have sum equal to twice their product. If the two equations
have the same two distinct roots, then
\[
 2(1-z^2)g=2\bigl[(1-z^2)^2g^2-z^2g(1-4g)\bigr],
 \qquad\text{that is}\qquad 1=g\bigl[(1-z^2)^2+4z^2\bigr]=g(1+z^2)^2 .
\]
The latter equality can be excluded by rotation.
For a horizontal output reflection is a real ambiguity, treated
by the lattice congruence below.

\smallskip
\noindent\emph{The limits in the case analysis.} The quantities
\(|h|^2V\), \(h_x^2U\) and \(z\) are invariant under a common real dilation
of the configuration. If the required split existed on an open set
of parameters, the resultant in \(p\) of \eqref{rate:eq:extgeo} (or the
horizontal test for \(z=0\)) would vanish identically in \((m,t)\).
Fix generic \(t\). In case (c), when the relative splitting angles
are fixed, the product of the differences between the candidate
horizontal defects and \(U\) would then vanish identically in \(m\).
Each factor is meromorphic under rotation, so at least one candidate
would have the same poles and leading coefficients as \(U\).
Put \(\varphi=\arctan m\). After division by the top branch radius
\(|P|\) or \(|Q|\), the primaries of that branch rotate rigidly with
\(\varphi\) at fixed \(t\). For a primary \(v\) of angle \(\theta\),
the lengths of \(J(v)\) and \(A(v)\) are
\(|v||\tan\theta|\) and \(|v|/|\cos\theta|\), respectively.
\begin{itemize}
\item \emph{Opposite branches.} Let \(h\) lie in the \(P\)-branch and let
\(m\to\infty\) at fixed generic \(t\). Then \(|P|=(1+m^2)^{-1/2}\to0\), while
\(Q=-im/(1-im)\to1\). At generic \(t\) every positive-level primary of
the limiting \(Q\)-branch has both components nonzero. Its companion and
axis have nonzero abscissae, as does the level-zero primary \(Q\to1\).
Thus \(U,V\) stay bounded. In the shrinking branch the positive-level
primaries have limiting directions off the coordinate axes, so their
companions and axes also shrink. Hence \(|h|^2V\to0\), \(h_x^2U\to0\),
\(g\to1\) and \(2+h_x^2U\to2\). For each nonaxial \(h\), both components
are nonzero for all sufficiently large finite \(m\), including \(h=P\),
for which \(z=m\). Case (a) applies at such \(m\). For an axial \(h\),
\(z=0\) and \(|h|^4V^2-|h|^2U-2\to-2\), giving (b).
If \(h\) lies in the \(Q\)-branch, take \(m\to0^+\). Then \(P\to1\),
\(|Q|\to0\), and the same bounds hold, with \(z=-1/m\) at \(h=Q\).
\item \emph{Top row against a strict descendant.} Here
\(U=P_x^{-2}+Q_x^{-2}-1\) and \(V=|P|^{-2}+|Q|^{-2}-1\) do not depend on \(t\).
At fixed generic \(m\) a vertex \(Pg_{nk}\) or \(Qg_{nk}\) with \(n\ge1\) has
length \(\asymp t^k(1+t^2)^{-n/2}\), which tends to zero as \(t\to0\) if
\(k\ge1\) and as \(t\to\infty\) if \(k=0\). Companions and axes carry the
additional factors \(|T_n|^{\pm1}\), which tend to \(|m|^{\pm1}\) as \(t\to0\),
or \(|\cos(\arg V)|^{-1}\), bounded at generic \(m\). The limiting
directions, \(\arctan m+k\pi/2\) and \(\arctan m-n\pi/2\) modulo \(\pi\), are
off the axes for generic \(m\). Hence \(g\to1\): case (a), or case (b) with
limit \(-2\).
\item \emph{Top row against a top child.} For \(h=P\),
\(|P|^2V=1+m^{-2}-(1+m^2)^{-1}\to1\) as \(m\to\infty\), so \(g\to1/2\); moreover
\(U=(1+m^2)^2(1+m^{-4})-1>0\). Thus (a) applies with \(z=m\) for
every sufficiently large finite \(m\). For \(h=Q\) one has
\(|Q|^2V=1+m^2-m^2(1+m^2)^{-1}\to1\) as \(m\to0^+\), so (a) again
applies with \(z=-1/m\) for each sufficiently small positive \(m\).
\item \emph{Primary target \(e\), with parent \(F\) and children \(c_1\perp c_2\),
against a primary \(h\) of its branch.} After normalization \(|h|\) and \(V\)
are constant and \(U=U(\varphi)\) has double poles exactly at the vertical
orientations of \(c_1,c_2,F\), with leading coefficients
\(|c_1|^{-2},|c_2|^{-2}>0\) and \(-|F|^{-2}<0\) (near a simple zero of
\(\cos(\varphi+\alpha)\), \(\cos^{-2}(\varphi+\alpha)\sim(\varphi-\varphi_0)^{-2}\)).
If \(g>1/4\), approach an orientation where \(c_1\) is vertical.
Its positive double pole gives \(U>0\) in a punctured neighborhood.
Throughout this rotation \(g\) is constant. Taking a nearby orientation
with \(h_xh_y\ne0\) gives \(2+h_x^2U>0\) and \(z\ne0\), so (a) applies
even when \(h\) is horizontal or vertical at the limiting orientation. If
\(0<g\le1/4\), the admissible splits are \(k=ah\pm\sqrt g\,|h|\,\hat h^\perp\),
\(a\in\{a_0,1-a_0\}\), \(a_0(1-a_0)=g\), and they rotate rigidly with \(h\).
The horizontal defect \(k_x^{-2}+l_x^{-2}-h_x^{-2}\) of such a split has
double poles at the vertical orientations of \(k,l,h\), with leading
coefficients \(|k|^{-2},|l|^{-2}>0\) and \(-|h|^{-2}<0\); these orientations
are distinct because \(0<a<1\). An identity with \(U(\varphi)\) forces the
negative poles to agree, so \(h\parallel F\) and \(|h|=|F|\), i.e.\
\(h=\pm F\), which does not occur for distinct marked primaries at
generic parameters.
\item \emph{Primary target against \(h=J(w)\).} At the orientation
\(\varphi_0\) where \(w=(x,y)\) is horizontal, \(h_x=y^2/x\) vanishes to second
order in \(\varphi-\varphi_0\) and \(|h|=|y|S^{1/2}/|x|\to0\), while \(V\) is
constant and \(U\) has at most a double pole there. Hence \(h_x^2U\to0\),
\(|h|^2V\to0\) and \(g\to1\), while \(z=-x/y\) is finite and nonzero for
\(\varphi\ne\varphi_0\): case (a).
\item \emph{Primary target against an axial \(h\).} At the orientations
where \(c_1\), respectively \(c_2\), is vertical, \(U\to+\infty\): the other
child is then horizontal, and \(F\) is not vertical, since its direction
differs from those of its children by \(\pm\arctan t\) modulo \(\pi/2\). The
normalized length of \(h=A(w)=|w|/|\cos(\arg w)|\) diverges only where
\(w\) is vertical, hence at most at one of these two orientations; at the
other one \(|h|\) is finite and nonzero, \(V\) is bounded, and
\(|h|^4V^2-|h|^2U-2\to-\infty\), giving case (b). For the top
\(h=1\), the normalized length is \(1/|P|=1/|\cos\varphi|\) in
the \(P\)-branch and \(1/|Q|=1/|\sin\varphi|\) in the \(Q\)-branch.
In either branch the normalized length stays finite and nonzero at
at least one of the two orientations, where the same test applies.
\item \emph{Private target at \(v=(x,y)\).} Then
\(U=x^{-2}+x^2y^{-4}-x^2S^{-2}\) and \(V=y^{-2}-x^2S^{-2}\), with \(S=x^2+y^2\);
note that \(U>0\), since \(x^{-2}-x^2S^{-2}=(S^2-x^4)/(x^2S^2)>0\). As \(v\)
becomes horizontal (\(y\to0\)), \(V\sim y^{-2}\), \(U\sim x^2y^{-4}\) and
\(U/V\to\infty\).
 \begin{itemize}
 \item \(h\) a primary of another level, off the axes at that
 orientation. For a split \(k+l=h\), \(k\perp l\), with \(V_{\rm split}=V\to\infty\)
 and \(|h|\) fixed, one leg, say \(k\), has \(|k|\to0\); then \(l\to h\), \(k\)
 becomes perpendicular to \(h\), \(k_x=\pm|k|h_y/|h|+o(|k|)\), and
 \(U_{\rm split}/V_{\rm split}\to|h|^2/h_y^2<\infty\), contradicting \(U/V\to\infty\);
 here \(g\to0\), so this is case (c).
 \item \(h\ne v\) a primary of the same level. The opposite-branch case
 was treated above. In the same branch, let \(k'\) and \(k\) be the
 respective indices of \(h\) and \(v\). Then \(k'\ne k\) and
 \(|h|/|v|=t^{k'-k}\). Keep the direction \(\theta\) of \(v\) at a
 generic fixed value and let this ratio tend to zero.
 Then \(|h|^2V=(|h|/|v|)^2(\sin^{-2}\theta-\cos^2\theta)\to0\) and \(h_x^2U\to0\):
 case (a).
 \item \(h=v\): let \(v\) become vertical (\(x\to0\)). Then \(|h|^2V\to1\),
 \(g\to1/2\) and \(U\to+\infty\): case (a).
 \item \(h=J(w)\) with \(w\) not parallel to \(v\): at the orientation where
 \(w\) is horizontal, \(V\) is finite, \(U>0\) has at most a double pole,
 \(h\to0\) and \(h_x^2U\to0\): case (a). If \(w=\lambda v\) with real
 \(\lambda\ne\pm1\), it suffices again to consider the same branch.
 Parallel primaries then have the same level at generic parameters.
 Their distinct indices give \(|h|/|J(v)|=|\lambda|=t^{k'-k}\),
 where \(k'\) is the index of \(w\) and \(k\) that of \(v\).
 At fixed generic \(\theta\), \(|J(v)|/|v|=|\tan\theta|\), so the
 same limit gives \(|h|^2V\to0\) and \(h_x^2U\to0\).
 If \(w=v\), then \(|J(v)|^2V=S/x^2-y^2/S\to1\) as \(y\to0\).
 Thus \(g>1/4\), \(U>0\), and \(z=-x/y\ne0\) for every sufficiently
 small nonzero \(y\), as required in (a).
 \item \(h\) axial: for \(V>0\) the horizontal test \(V^2X^2-UX-2=0\) in
 \(X=|h|^2\) has exactly one positive root, the product of the roots being
 \(-2/V^2\). The designated split of \(A(v)\) satisfies it, so
 \(X=|A(v)|^2\); distinct marked horizontal vertices have generically
 distinct lengths, so only \(h=A(v)\) remains: case (b).
 \end{itemize}
\end{itemize}
\end{proof}

\begin{lemma}\label{rate:l:denom}
Fix a word cutoff $L=B^{k_0}$ with $k_0$ fixed, and an integer
$m_*$ independent of $B$. The parameters can be chosen with
$1/2<t<3/4$, $m$ in a fixed compact positive interval, so that
all nonidentities of Lemma~\ref{rate:l:geom} needed for combinations
of at most $m_*$ rows, including repetitions and signs, hold.
All primary directions satisfy $|\sin\theta|,|\cos\theta|\ge c/d$.
There is a common integer denominator $P_{\mathrm{lat}}$ with
\[
 \log P_{\mathrm{lat}}\le C_{k_0,m_*}d^2\log(2+B)
 \quad(2\le d\le B).
\]
For every sum of at most $L$ signed marked macro vertices, with
repetitions allowed, a zero horizontal component implies a zero vertical component. On clearing
denominators, all nonaxial heads have nonzero residues on one
isotropic line in $\mathbb F_5^2$ and all axial heads are zero
there.
\end{lemma}
\begin{proof}
The bounded arity tests in Lemma~\ref{rate:l:geom} use at most
$d^{C_{m_*}}$ polynomials of degree at most $C_{m_*}d^{C}$.
For an external pair use the resultant in \eqref{rate:eq:extgeo}
with the horizontal test for $z=0$; for its own output test
the remaining reflection condition. Thus there is no enumeration
of external integer points before specializing. For each nonzero
component comparison choose just one nonzero polynomial, and
include the finitely many nonzero denominators and the determinants
distinguishing the $2(d-1)$ moving directions. Choose
$t=p/q$, $q$ a prime greater than a sufficiently large polynomial
in $d,L$, the number and the degree of these tests, yet bounded
by twice that polynomial. A nonzero bivariate polynomial of degree
$D$ specializes identically in $m$ at at most $D$ real values of
$t$. There are therefore enough even $p\equiv q\pmod5$ in the
stated interval, with $p>10L$. Write $S$ for the sum of $2,d,L$
and the test counts and degrees; $\log q=O(\log S)$.
For definiteness: the designated rows are the top row, the
\(\sum_{n=0}^{d-2}2(n+1)\) primary rows and the \(\sum_{n=1}^{d-1}2(n+1)\)
private rows, in all \(O(d^2)\), and the marked vertices are \(O(d^2)\) in
number. Each marked vertex is a rational function of \((m,t)\) whose
numerator and denominator have degree \(O(d)\): \(g_{nk}=(it)^k(1+it)^{-n}\),
\(P=(1-im)^{-1}\), and companions and axes arise from the complexified
coordinates through \(J(V)=(V_y^2/V_x,-V_y)\) and \(A(V)=(|V|^2/V_x,0)\). A
test involving at most \(m_*\) designated rows with integer coefficients
bounded by \(m_*\), and possibly one marked output of an external split
(through the resultant of \eqref{rate:eq:extgeo} in \(p\), or the horizontal
test), clears to a polynomial of degree at most \(C_{m_*}d\), and the
number of such tests is at most \((Cd^2m_*)^{Cm_*}\le d^{C_{m_*}}\). If a
nonzero polynomial \(F(m,t)=\sum_jf_j(t)m^j\) of degree at most \(\mathsf D\)
vanished identically in \(m\) at \(t=t_0\), then \(f_j(t_0)=0\) for a fixed
\(f_j\not\equiv0\) of degree at most \(\mathsf D\); this gives at most \(\mathsf D\)
exceptional values of \(t\) per test, together at most \(C_{m_*}d^{C_{m_*}}\).
The number of even \(p\equiv q\pmod5\) with \(q/2<p<3q/4\) is at least
\(q/40-2\). Taking for \(q\) a prime between \(\mathsf Q_0\) and \(2\mathsf Q_0\)
(Bertrand), where \(\mathsf Q_0\) is a sufficiently large fixed power of
\(S\), this number exceeds the number of exceptional values, and
\(q/2>10L\) gives \(p>10L\). Since \(q\) is an odd prime larger than \(p\),
\(\gcd(p,q)=1\), and \(\log q=O(\log S)\).

Write $a_n+i b_n=(q-ip)^n$. Its two integer components are
coprime: a common rational prime would divide this Gaussian
integer and its conjugate, whereas $q-ip$ and $q+ip$ are
coprime in $\mathbb Z[i]$ by opposite parity and
$\gcd(p,q)=1$. At a moving axis the denominator in $m=k/M$ is
of the form $D_i=A_i M+B_i k$, where $(A_i,B_i)$ is one of
the primitive signed permutations of $(a_n,b_n)$; denominators
proportional up to sign are counted once. To see the formulas,
$A(Pu^n)=q^n/(a_n-mb_n)$ and
$A(iPu^n)=-q^n/(b_n+ma_n)$, since $\Re Pu^n$ and $\Im Pu^n$
are the corresponding numerators divided by
$(1+m^2)(p^2+q^2)^n/q^n$. Explicitly, $t=p/q$ gives
$u=q(q-ip)/(p^2+q^2)$, and
\[
 Pu^n=\frac{q^n(1+im)(a_n+ib_n)}{(1+m^2)(p^2+q^2)^n},\qquad
 |Pu^n|^2=\frac{q^{2n}}{(1+m^2)(p^2+q^2)^n},
\]
so that $\Re Pu^n$ and $\Re(iPu^n)=-\Im Pu^n$ have the numerators
$q^n(a_n-mb_n)$ and $-q^n(b_n+ma_n)$, and $A(V)=|V|^2/V_x$.
With $m=k/M$ the two denominators are $(a_nM-b_nk)/M$ and
$(b_nM+a_nk)/M$. Magnitudes from the other branch
supply $m$, and a real dilation rescales $A$ by the same
real number. Impose $A_i B_i\ne0$ by the polynomial exclusions.
Let $R_0$ be the product of $10pq(p^2+q^2)$, all nonzero
$|A_i|,|B_i|$, and absolute values of all nonzero pair
determinants. Its logarithm is $O(d^3\log q)$.

Take $M$ prime with $\log M=C_0d\log S+O(1)$, larger than
each individual factor of $R_0$. Restrict $k/M$ to an interval
of width $\gtrsim1/d$ where all $2d$ direction projections
are bounded below by $c/d$ of the vector length. Such an interval
exists inside a fixed compact interval of rotations by deleting
angular neighborhoods of total length less than half that interval
and counting components. Choose a residue of $k$ modulo $25$
so that $v_5(1+(k/M)^2)=1$. There are $N_0\gtrsim M/d$ candidates.
In this interval $|D_i|\gtrsim M/d^C$ for every $i$.

For each fixed $i$, the expected logarithm of the part of $D_i$
supported on primes dividing $R_0$ is $O(\log S)$. For $\pi\mid
B_i$ the contribution is zero. Otherwise at each prime power
$\pi^h$ a single congruence is needed. Its frequency is
$\le\pi^{-h}+O(1/N_0)$, allowing a factor $25$ for $\pi=5$
and its first two powers without restriction. The summed endpoint
error is $O[(\log R_0)\log(\max_i|D_i|)/N_0]=o(1)$. Finally
$\sum_{\pi\mid R_0}\log\pi/(\pi-1)\lesssim1+\log\log(3R_0)$;
split at $\log^2(3R_0)$, use
$\sum_{\pi\le x}\log\pi=O(x)$ below and the divisor product above.
By Markov's inequality, increasing $C_0$, a positive fraction of
$k$ have summed bad logarithm over all $i$ below $\log M/3$.
Omit the boundedly many multiples of the prime $M$ inside the
rotation interval, so $\gcd(k,M)=1$.
At most polynomially many remaining $k$ violate a polynomial
test. Hence every selected $D_i$ has a factor $h_i>M^{1/2}$
coprime to $R_0$, all the other $D_j$, $M$ and $k$; use the
pairwise determinants and primitivity. It is also coprime to
$M^2+k^2$, because a prime common to both would divide
$A_i^2+B_i^2=(p^2+q^2)^n$.

\smallskip
\noindent\emph{The choice of $k$ in detail.} Let $\mathcal K_0$ be the
set of candidates, namely the integers $k$ with $k/M$ in the chosen
interval and $k$ in the chosen class modulo $25$, so that
$N_0=\#\mathcal K_0\gtrsim M/d$. Let $D_i^{\rm bad}$ be the largest
divisor of $D_i$ composed of primes dividing $R_0$, and
$\beta_i(k)=\log D_i^{\rm bad}$. Let $\pi\mid R_0$. If $\pi\mid B_i$, then
$\pi\nmid A_i$ by primitivity and $\pi\ne M$, so $\pi\nmid D_i$. If
$\pi\nmid B_i$, then $\pi^h\mid D_i$ means
$k\equiv-A_iMB_i^{-1}\pmod{\pi^h}$, a single residue class. The
candidates form an arithmetic progression of difference $25$ intersected
with an interval, so a fixed class modulo $\pi^h$ has relative frequency
at most $\pi^{-h}+O(1/N_0)$ for $\pi\ne5$, and at most
$25\cdot5^{-h}+O(1/N_0)$ for $\pi=5$ and $h\ge3$. Only powers with
$\pi^h\le|D_i|$ occur. Hence
\[
 \frac1{N_0}\sum_{k\in\mathcal K_0}\beta_i(k)
 \le C\sum_{\pi\mid R_0}\frac{\log\pi}{\pi-1}+C
   +\frac{C\,\omega(R_0)\log\max_i|D_i|}{N_0},
\]
where $\omega(R_0)\le\log R_0$ is the number of prime factors of $R_0$.
The last term is $o(1)$ because $N_0\gtrsim M/d$ while $\log R_0$ and
$\log|D_i|$ are polynomial in $d\log S$. In the first sum, the primes
$\pi\le\log^2(3R_0)$ contribute $O(\log\log(3R_0))$ by Mertens' bound
$\sum_{\pi\le x}\log\pi/\pi=O(\log x)$, and each of the at most
$\omega(R_0)$ larger prime factors contributes at most
$4\log\log(3R_0)/\log^2(3R_0)$, which gives $o(1)$ in total. As
$\log R_0=O(d^3\log q)$ and $\log q=O(\log S)$, the mean of $\beta_i$ is
$O(\log S)$. Summing over the $O(d)$ indices $i$, the mean of
$\beta=\sum_i\beta_i$ is at most $C'd\log S$, and Markov's inequality gives
\[
 \#\{k\in\mathcal K_0:\beta(k)>\tfrac13\log M\}
 \le\frac{3C'd\log S}{\log M}\,N_0\le\frac{N_0}4
\]
once $C_0\ge12C'$. The excluded multiples of $M$ and the zeros of the
polynomial tests number at most $\operatorname{poly}(S)<N_0/4$, since
$\log M=C_0d\log S+O(1)$. Every remaining $k$ satisfies $\gcd(k,M)=1$ and
all tests, and $h_i\coloneqq|D_i|/D_i^{\rm bad}\ge cMd^{-C}M^{-1/3}>M^{1/2}$.
A prime dividing $h_i$ and $D_j$, $j\ne i$, divides
$(A_iB_j-A_jB_i)M$ and $(A_iB_j-A_jB_i)k$, hence the nonzero determinant
$A_iB_j-A_jB_i$, which divides $R_0$; this is impossible. A prime dividing
$h_i$ and $M$ equals $M$ and divides $B_ik$, which is impossible because
$M>|B_i|$ and $M\nmid k$. A prime dividing $h_i$ and $k$ divides $A_iM$,
hence $A_i$ and $R_0$, which is again impossible.

Rewrite every companion as an axis minus a primary. The combined
axis part at $D_i$ of any word of original length $L$ is
$(M r_i+k s_i)/D_i$ with integer $|r_i|+|s_i|\le CL(2q)^d$:
the numerators before summing are signed $p^j q^{n-j}$.
The common primary denominators divide
$(M^2+k^2)(q^2+p^2)^d$. If the horizontal word vanishes,
clear denominators and reduce at $h_i$. This gives
$h_i\mid B_i r_i-A_i s_i$, hence equality to zero since
$h_i\gg L(2q)^{2d}$. Primitivity implies $r_i=c_i A_i$,
$s_i=c_i B_i$, with integer $c_i$. In detail, multiply the vanishing
horizontal word by $D_i$ and by the remaining denominators, which are
units modulo $h_i$. Every term other than the axis part at $D_i$
acquires the factor $D_i\equiv0$, so
$Mr_i+ks_i\equiv0\pmod{h_i}$. Since $A_iM\equiv-B_ik$ modulo $h_i$,
\[
 0\equiv A_i(Mr_i+ks_i)\equiv k(A_is_i-B_ir_i)\pmod{h_i},
\]
and $k$ is a unit. The integer $B_ir_i-A_is_i$ has modulus at most
$(|A_i|+|B_i|)CL(2q)^d\le CL(2q)^{2d}<M^{1/2}<h_i$, so it is zero, and
then $(Mr_i+ks_i)/D_i=c_i$. The actual total
$\sum|c_i|\lesssim d L$, since every physical real axis
in the chosen angular interval has modulus $O(d)$.

Subtract the resulting integer constant $C$. The remaining
primary part is $F(u)+P(G(u)-F(u))$ for integer polynomials
of degree $\le d$, given as sums of at most $2L$ signed
Bernstein monomials. The horizontal condition after multiplication
by $1+m^2$ gives a real quadratic in $m$ with constant term
$C+\Re G$, linear term $-\Im(G-F)$ and quadratic term
$C+\Re F$. Its coefficients become integers of size
$\le C d L(2q)^{2d}$ after multiplying by $(q^2+p^2)^d$.
The rational root $k/M$ forces all of these integers to vanish
successively since $M$ is greater than each. So $G(u)=F(u)$
and $\Re F(u)=-C$. In detail, $P=(1+im)/(1+m^2)$ gives
\[
 (1+m^2)\bigl(C+\Re[F+P(G-F)]\bigr)
 =(C+\Re G)-m\,\Im(G-F)+m^2(C+\Re F),
\]
and if integers $c_0,c_1,c_2$ of modulus less than the prime $M$
satisfy $c_0M^2+c_1kM+c_2k^2=0$ with $\gcd(k,M)=1$, then $M\mid c_2$,
so $c_2=0$; next $M\mid c_1$, so $c_1=0$; and finally $c_0=0$.
The three vanishing coefficients give $\Im(G-F)=0$,
$\Re(G-F)=(c_0-c_2)/(q^2+p^2)^d=0$ and $\Re F=-C$. Here $u=q/(q+ip)$ and
$1-u=ip/(q+ip)$, so a Bernstein monomial $u^{n-k}(1-u)^k$ equals
$q^{n-k}(ip)^k(q+ip)^{d-n}/(q+ip)^d$. Write $F(u)=a/(q+ip)^d$ with $a$ a
Gaussian integer. Then $i\Im F(u)(q+ip)^d=a+C(q+ip)^d$ is also a
Gaussian integer. The coprimality of the two components of the
power shows $\Im F(u)\in\mathbb Z$: the rational number $\Im F(u)$
times each of two coprime integers is an integer. Reduction modulo $(p)$
gives $u=1$ there and thus $p\mid\Im F(u)$; indeed every monomial with
$k\ge1$ is divisible by $p$ and the others are real integers
modulo $p$. The original
word bounds give $|\Im F(u)|\le 2L<p$, proving the claim: the
vertical component of the word is $\Im F(u)$, because the axis parts
are real and $G=F$.

Clearing both linear denominators at each level, as well as
the common primary denominators, costs $O(d\log M+d^2\log q)$
in their logarithm. At the prime five the numerator of $P$
lies on an isotropic line and has nonzero coordinates, while
the denominator has valuation one. $u,itu$ are units at both
Gaussian factors above five since $1+t^2\not\equiv0$ there;
the same holds for $-im$. Hence, after multiplying the whole
configuration by five in the localized rational lattice,
each primary belongs to the same nonzero line. Its squared
length vanishes modulo five and its horizontal component is
a unit, so $A$ vanishes modulo five and $J=-V$ modulo five.
Choosing a common integer denominator with exactly the needed
five-adic valuation proves the lattice statement.
In detail, with \(m=k/M\) one has \(P=(1+im)/(1+m^2)=M(M+ik)/(M^2+k^2)\)
and \(v_5(M^2+k^2)=1\). Hence \(5P\) is \(5\)-integral and congruent modulo
\(5\) to a \(5\)-adic unit times the vector \((M^2,Mk)\). As
\(M^2+k^2\equiv0\pmod5\) and \(5\nmid M\), this vector lies on one of the two
isotropic lines \(\{y\equiv\pm2x\}\) of \(\mathbb F_5^2\), and it has nonzero
coordinates. The two isotropic lines are the kernels of reduction
modulo the two Gaussian primes \(2\pm i\) above five, so multiplication by
a Gaussian number which is a unit at both of them preserves each line
and keeps nonzero vectors nonzero; on an isotropic line a nonzero
vector has both coordinates nonzero. The factors
\(u=q(q-ip)/(p^2+q^2)\), \(it=ip/q\) and \(-im=-ik/M\) are such units:
\(p\equiv q\not\equiv0\) gives \(p^2+q^2\equiv2q^2\not\equiv0\), and
\(k\equiv\pm2M\not\equiv0\). Thus \(W=5V\) lies on the chosen line, with
nonzero coordinates, for every primary \(V\). For such a \(5\)-integral
\(W\) the squared length \(W_x^2+W_y^2\) is divisible by five while \(W_x\)
is a unit, so \(A(W)=|W|^2/W_x\equiv0\) and \(J(W)=A(W)-W\equiv-W\pmod5\);
moreover \(5\cdot1\equiv0\). Multiplying by an integer prime to five
clears the remaining denominators without changing these residues.
\end{proof}

For later scalar selection, on a lattice with vertical period
spacing $R$ take physical marked indices $e_1+Rv$ where $v$
is an integer vector from the cleared configuration. An exterior
bare diameter-circle relation has
\[
 x^2-Xx+q(q-Y)=0,\qquad x=(n-1)/R,\quad q\in\mathbb Z,
 \quad h=(X,Y)\in\mathbb Z^2 .
\]
The monic equation forces the rational $x$ to be an integer. All
such pairs have length $O(dP_{\mathrm{lat}})$. Thus for every remaining
nonzero component comparison its rational height is at most
$(2dP_{\mathrm{lat}})^{C_{m_*}}$, including all arbitrary exterior
diameter-circle points. A reflection of a nonzero isotropic vector
at a horizontal output exchanges the two different lines mod five.
Only the chosen one belongs to the invariant translated coset
$k_x-1=l k_y\bmod5$, $l^2=-1\bmod5$. Choosing $5\nmid R$ preserves
this condition. The further choice of an axial parameter separating
the two-component defects in the single physical frequency is
described in the next section.

\subsection{A finite trajectory landing at a specified axial point}
All estimates in this section use the ordinary Euclidean norm on the
independent complex graph coordinates. The real sector introduced below
serves only to find a reference orbit. The boundary inverse is for
arbitrary complex perturbations satisfying the usual conjugate-copy
reality.

\begin{lemma}\label{rate:l:finite}
Number the two top children at level $1$ and the last rows at level
$d$, with $c_{nk}=\binom{n-1}{k}^{1/2}$. Let $D=C_*d\log(d+2)$,
$G_n=D^n$, $b=e^{-c_1D}$, $\kappa=(2+d+D)^{-10}$.
There are absolute $c_1>0$, $C_*<\infty$ for which the following
holds for sufficiently large $d$. Give the top edge coefficient one,
the edges below level $n\ge1$ coefficients
$G_nc_{nk}/(c_{n+1,k}c_{n+1,k+1})$, and the private triangles
coefficients $\kappa/(bc_{nk})$. The Hamiltonian of an edge has
the form $i g(\bar h v w-h\bar v\bar w)$, with $h$ its output;
choose the symplectic convention giving $\dot h=g v w$. There is an
orbit whose initial primary amplitudes vanish except on the last
rows and whose endpoint has only axial amplitudes: the top equals
one and each private axis equals $b c_{nk}$. It decays toward
that point at rate $\kappa/2$ from time $D$. The propagator on
any subinterval of $[0,D]$, in either direction, and the mixed
boundary inverse on $[0,\infty)$ with any exponential rate
$\kappa/8\le\sigma\le\kappa/3$, have norms at most $\exp(C dD)$.
The stable projection is prescribed at $D$ and all its real
quadratures are included.
\end{lemma}
\begin{proof}
On the real sector put $c_{nk}x_n,c_{nk}w_n,c_{nk}h_n$ on a
primary, companion and private axis, respectively. The top is $X$.
Both subtrees use the same scalar values. For the edge with parent
$(n,k)$ the product of the two children cancels the two denominator
factors in its coefficient, leaving $G_n c_{nk}x_{n+1}^2$.
When $(n,k)$ is a child, its two possible parents contribute together
\[
 -\frac{G_{n-1}}{c_{nk}}
 (c_{n-1,k-1}^2+c_{n-1,k}^2)x_{n-1}x_n
 =-G_{n-1}c_{nk}x_{n-1}x_n
 \qquad (n\ge2),
\]
by Pascal's identity, with nonexistent indices set to zero.
In detail, the edge Hamiltonian gives $\dot h=gvw$,
$\dot v=-gh\bar w$ and $\dot w=-gh\bar v$. The vertex $(n,k)$ is
the child $uV$ of $V=(n-1,k)$, with sibling $(n,k+1)$, and the child
$(1-u)V$ of $V=(n-1,k-1)$, with sibling $(n,k-1)$. On the real sector
the two contributions to its amplitude $c_{nk}x_n$ are
\[
 -\frac{G_{n-1}c_{n-1,k}}{c_{nk}c_{n,k+1}}\,
   c_{n-1,k}x_{n-1}\,c_{n,k+1}x_n
 \quad\text{and}\quad
 -\frac{G_{n-1}c_{n-1,k-1}}{c_{n,k-1}c_{nk}}\,
   c_{n-1,k-1}x_{n-1}\,c_{n,k-1}x_n ,
\]
and $c_{n-1,k-1}^2+c_{n-1,k}^2=\binom{n-2}{k-1}+\binom{n-2}{k}
=\binom{n-1}{k}=c_{nk}^2$. For a private triangle at $(n,k)$ the
coefficient $\kappa/(bc_{nk})$ multiplies two of the amplitudes
$c_{nk}x_n,c_{nk}w_n,c_{nk}h_n$, which gives the three private terms
below after division by $c_{nk}$.
At $n=1$ the only parent is the top. Division by $c_{nk}$ now
gives the same equation for every $k$ on a level. The sector is
therefore invariant and satisfies
\begin{align*}
 X'&=x_1^2,\qquad
 x_n'=G_n x_{n+1}^2-G_{n-1}x_{n-1}x_n+f_n,\\
 f_n&=-(\kappa/b) w_n h_n,\quad
 w_n'=-(\kappa/b)x_n h_n,\quad h_n'=(\kappa/b)x_n w_n .
\end{align*}
Here $x_0=X$, the last forward input is omitted and there
are no private variables at $n=1$. The identities
$w_n^2+h_n^2=b^2$ and
$X^2+\sum x_n^2+\sum_{n\ge2}h_n^2=1+(d-1)b^2$ hold when
imposed initially. Indeed
$(w_n^2+h_n^2)'=2(\kappa/b)x_n(-w_nh_n+h_nw_n)=0$, and, with
$x_0=X$ and $G_0=1$, half the derivative of the second sum is
\[
 x_0x_1^2+\sum_{n=1}^{d-1}G_nx_nx_{n+1}^2
   -\sum_{n=1}^{d}G_{n-1}x_{n-1}x_n^2
   +\sum_{n\ge2}\bigl(x_nf_n+(\kappa/b)x_nw_nh_n\bigr)=0,
\]
since the first three terms telescope and each bracket vanishes.
Start at $X=x_1=\cdots=x_{d-1}=0$ with
arbitrary circle phases, and choose the unique positive $x_d$
from this last constraint. The primary norm squared is
$X^2+\sum x_n^2=1+\sum w_n^2$, because
$(d-1)b^2-\sum h_n^2=\sum w_n^2$.

Set $T=D$, $a_0=1/10$, $Y_m^2=\sum_{n=m}^d x_n^2$. On
$[0,T]$ every $x_n\ge-2\kappa bT$: at a first exit, on a
negative excursion the incoming square is nonnegative and a
negative parent amplifies the negative part at rate at most
$2G_{d-1}\kappa bT$. Integration gives
$\kappa bT\exp(2G_{d-1}\kappa bT^2)<2\kappa bT$.
Here and below choose $c_1$ sufficiently small first and $C_*$
large relative to $1/c_1$, so $G_{d-1}b$ times every fixed
polynomial in $d,D$ tends to zero. Similarly integrating the
positive part gives $|x_j(t)|\le C G_jt+C\kappa bT$ for
$j<d$, and $X(t)\le Ct$. Thus the energy of ancestors up to
$m-2$, including $X$, on $t\le3/G_{m-1}$ is
$O(D^{-2}+d\kappa^2b^2T^2)$ for $m\ge2$, by a geometric sum.
In detail, for \(1\le j<d\) and while \(x_j>0\),
\(x_j'=G_jx_{j+1}^2-G_{j-1}x_{j-1}x_j+f_j\le2G_j+2\sqrt2\,G_{j-1}\kappa bT
+\kappa b/2\le3G_j\), using \(x_{j+1}^2\le2\) and \(0<x_j\le\sqrt2\)
from the mass identity, \(x_{j-1}\ge-2\kappa bT\), \(|f_j|\le\kappa b/2\)
and \(G_{j-1}\le G_j\). As \(x_j(0)=0\) for \(j<d\), this gives
\(x_j(t)\le3G_jt\), and with the lower bound
\(|x_j(t)|\le3G_jt+2\kappa bT\). Also \(0\le X(t)\le2t\), since
\(X'=x_1^2\le2\). For \(m\ge2\), \(t\le3/G_{m-1}\) and
\(1\le j\le m-2\) this yields \(|x_j|\le9D^{-(m-1-j)}+2\kappa bT\) and
\(X\le6D^{-1}\), because \(G_n=D^n\). Hence, for \(D\ge2\),
\[
 X^2+\sum_{j=1}^{m-2}x_j^2\le36D^{-2}+\sum_{i\ge1}162D^{-2i}
   +8d\kappa^2b^2T^2\le300D^{-2}+8d\kappa^2b^2T^2 .
\]
Telescoping at any $m$ gives
\[
 (Y_m^2)'=-2G_{m-1}x_{m-1}x_m^2+2\sum_{n=m}^d x_nf_n .
\]
Consequently $Y_m'\le 2G_{d-1}\kappa bT Y_m+\sqrt d\,\kappa b$.
Where $x_{m-1}\ge c>0$ and $Y_{m+1}\le E$ we have the stronger
$Y_m'\le-cG_{m-1}(Y_m-E)+\sqrt d\,\kappa b$; for $Y_m<E$
use nonnegative dissipation directly. These inequalities hold in
the upper Dini sense at zero. To verify them, note that
$|w_nh_n|\le b^2/2$ on the circle, so $|f_n|\le\kappa b/2$ and
$2\sum_{n\ge m}x_nf_n\le\sqrt d\,\kappa b\,Y_m$ by Cauchy--Schwarz.
The lower bound $x_{m-1}\ge-2\kappa bT$ gives
$-2G_{m-1}x_{m-1}x_m^2\le4G_{d-1}\kappa bT\,Y_m^2$, and division of
the telescoped identity by $2Y_m$ gives the first inequality. In the
second case $x_m^2=Y_m^2-Y_{m+1}^2\ge Y_m^2-E^2$, so that
\[
 2Y_mY_m'\le-2cG_{m-1}(Y_m-E)(Y_m+E)+\sqrt d\,\kappa b\,Y_m ,
\]
and $(Y_m+E)/Y_m\ge1$ when $Y_m\ge E$.

Descend from $m=d$ to $2$, with $Y_{m+1}\le E_{m+1}$ on
$[a_0/G_{m-1},T]$, starting at zero for $m=d$. On the shorter
interval up to $3/G_{m-1}$ the small ancestors and the mass
identity give $x_m^2\ge0.7$ while $x_{m-1}<1/2$ if
$E_{m+1}<1/20$. Hence $x'_{m-1}\ge G_{m-1}/2$ there until
the crossing; the same inequality prevents a downward crossing.
From the first crossing until $a_0/G_{m-2}$,
$x'_{m-1}\ge-2G_{m-2}x_{m-1}-\kappa b$ gives
$x_{m-1}\ge1/3$. Integrating the improved $Y_m$ inequality from
$3/G_{m-1}$ yields
\[
 Y_m(t)\le E_{m+1}+C\sqrt d\kappa b/G_{m-1}
        +2e^{-G_{m-1}(t-3/G_{m-1})/3}.
\]
At $a_0/G_{m-2}$ the last term is $O(e^{-c_0D})$ for an
absolute $c_0>0$. Using the general upper derivative for the
rest of the interval permits
\[
 E_m\le(1+O(G_{d-1}\kappa bT^2))E_{m+1}
       + C e^{-c_0D}+C\sqrt d\,\kappa bT.
\]
All $E_m<1/20$, as the product of the $d$ prefactors is bounded.
For the last step $Y_2$ is small from $a_0$; while $X<1/2$
one has $X'=x_1^2\ge0.7$. After time 3 the top remains above
$1/2$. The improved inequality on $[3,D]$ gives uniformly in
every initial circle phase
\[
 Y_1(D)\le C d e^{-c_0D}+C d^{3/2}\kappa bD=o(b)
 \quad\text{if }c_1<c_0/4.
\]
These same estimates give
\begin{equation}\label{rate:eq:l1G}
 \sum_{n<d}\int_0^D G_n(|x_n|+|x_{n+1}|)\,dt\le CdD .
\end{equation}
For $x_n$, before $3/G_{n-1}$ use the bounded primary amplitudes;
on the next plateau integrate the exponential at cost
$G_n/G_{n-1}=D$. The small constant terms during and after the
plateau cost $O(G_{d-1}D[d e^{-c_0D}+d^{3/2}\kappa bD])=o(1)$.
For $x_{n+1}$ use its own plateau and the smaller factor $G_n$.
At $n=1$ use the $[3,D]$ bound, and at $n=0$ use $G_0=1$.
In detail, the recursion for \(E_m\) recorded below, started at
\(E_{d+1}=0\), gives \(E_m\le\varepsilon_d\coloneqq
C(de^{-c_0D}+d^{3/2}\kappa bD)\) for every \(m\). Fix \(2\le n\le d-1\).
For \(t\le3/G_{n-1}\) use \(|x_n|\le\sqrt2\); for
\(3/G_{n-1}\le t\le a_0/G_{n-2}\) use \(|x_n|\le Y_n\) and the displayed
bound \(Y_n(t)\le E_{n+1}+C\sqrt d\,\kappa b/G_{n-1}
+2e^{-G_{n-1}(t-3/G_{n-1})/3}\); for \(t\ge a_0/G_{n-2}\) use
\(|x_n|\le Y_n\le E_n\). Since \(\int_{3/G_{n-1}}^\infty
2e^{-G_{n-1}(t-3/G_{n-1})/3}dt=6/G_{n-1}\),
\[
 \int_0^DG_n|x_n|\,dt\le3\sqrt2\,\frac{G_n}{G_{n-1}}
 +6\frac{G_n}{G_{n-1}}+G_nD\bigl(\varepsilon_d+C\sqrt d\,\kappa b\bigr)
 \le(3\sqrt2+6)D+o(1),
\]
because \(G_n/G_{n-1}=D\) and
\(G_{d-1}D(\varepsilon_d+\sqrt d\,\kappa b)\to0\) by the choice of
\(c_1\) and \(C_*\). The same three ranges for \(x_{n+1}\), shifted by
one level, give \(\int_0^DG_n|x_{n+1}|\,dt\le3\sqrt2+6+o(1)\). At
\(n=d-1\) the leaf amplitude \(x_d\) is treated in the same way with
\(m=d\). At \(n=1\), \(|x_1|\le\sqrt2\) on \([0,3]\) and the improved
inequality gives \(|x_1|\le Y_1\le2e^{-(t-3)/3}+\varepsilon_d\) on
\([3,D]\), so \(\int_0^DG_1(|x_1|+|x_2|)\,dt=O(D)\); at \(n=0\),
\(\int_0^D(|X|+|x_1|)\,dt\le2\sqrt2D\). Summing over the \(d\) levels
gives \eqref{rate:eq:l1G}.
We record the elementary steps behind the constants used above. The
mass identity gives $|X|,|x_n|\le(1+(d-1)b^2)^{1/2}\le\sqrt2$, and the
circle identity gives $|f_n|\le\kappa b/2$.
\emph{Lower bound.} Let $t_1\le T$ be the first time at which some
$x_n=-2\kappa bT$. On a maximal interval before $t_1$ on which $x_n<0$
the incoming square is nonnegative, so
$(-x_n)'\le G_{n-1}x_{n-1}x_n+\kappa b/2\le2G_{d-1}\kappa bT(-x_n)+\kappa b/2$,
using $x_{n-1}\ge-2\kappa bT$ when $x_{n-1}<0$ and dropping the term
otherwise. Gronwall from the start of the excursion gives
$-x_n\le\tfrac12\kappa bT\exp(2G_{d-1}\kappa bT^2)<\kappa bT$, a
contradiction; hence $t_1$ does not exist.
\emph{Plateau.} On $[a_0/G_{m-1},3/G_{m-1}]$ the ancestors carry energy
$O(D^{-2}+d\kappa^2b^2T^2)\le0.01$, and $Y_{m+1}^2\le E_{m+1}^2<0.0025$. As
long as $x_{m-1}<1/2$, the mass identity $X^2+\sum x_n^2\ge1$ gives
$x_m^2\ge1-0.01-0.25-0.0025>0.7$, and then
$x_{m-1}'\ge0.7G_{m-1}-\tfrac{\sqrt2}2G_{m-2}-\kappa b/2\ge G_{m-1}/2$,
because $G_{m-2}/G_{m-1}=D^{-1}$. Starting above $-2\kappa bT$ at time
$a_0/G_{m-1}$, the amplitude $x_{m-1}$ therefore reaches $1/2$ before
$(a_0+1+4\kappa bT)/G_{m-1}<3/G_{m-1}$, and it cannot cross $1/2$
downward before $3/G_{m-1}$. Afterwards
$x_{m-1}'\ge-\sqrt2\,G_{m-2}x_{m-1}-\kappa b/2$ while $x_{m-1}>0$, so up to
the time $a_0/G_{m-2}$ one has
$x_{m-1}\ge\tfrac12e^{-\sqrt2a_0}-\kappa ba_0/G_{m-2}>1/3$. With $c=1/3$
the improved inequality yields the displayed bound for $Y_m$, using
$Y_m\le\sqrt2$, and at $t=a_0/G_{m-2}$ its exponential term equals
$2\exp(-(a_0D-3)/3)\le Ce^{-D/30}$; thus $c_0=1/30$ is admissible.
\emph{Recursion.} On the rest of $[a_0/G_{m-2},T]$ the general inequality
gives $Y_m(t)\le(Y_m(t_*)+\sqrt d\,\kappa bT)\exp(2G_{d-1}\kappa bT^2)$ with
$t_*=a_0/G_{m-2}$, which is the stated recursion for $E_m$. Over the $d$
steps the product of the prefactors is at most
$\exp(CdG_{d-1}\kappa bT^2)\le2$ and the additive terms total
$O(de^{-c_0D}+d^{3/2}\kappa bT)$, so all $E_m<1/20$.

Near $X=1,h_n=b$ put $\xi_n=x_n/b$, $\zeta_n=w_n/b$.
The invariants give the positive branches
$h_n/b=\sqrt{1-\zeta_n^2}$ and
$X=[1+b^2(\sum\zeta_n^2-\sum\xi_n^2)]^{1/2}$.
The linear equations are $\dot\xi_1=-\xi_1$ and
$(\dot\xi_n,\dot\zeta_n)=-\kappa(\zeta_n,\xi_n)$ for $n\ge2$.
On a small Euclidean ball of radius $\delta$ the derivative of the
nonlinear remainder costs $C[bG_{d-1}\delta+(\kappa+b^2)\delta^2]$,
by the nearest-neighbor products and the two explicit square roots.
Integrate the stable variables forward from the section and the
unstable variables backward from infinity at weight $\kappa/2$;
the kernel has $L^1$ norm $O(1/\kappa)$ and both projections have
norm one. This gives a stable manifold which projects onto a
fixed small ball of all the $\xi$'s: the stable coordinates project
isomorphically onto them and the derivative of the graph over the
linear subspace is $O(\delta)$ by the same contraction.
Write the corresponding angles as $F(\xi)$, continuously lifted
about the positive-axis origins. If $\beta$ is the vector of initial
phases, the map taking $\beta$ to the actual phases at $D$ minus
$F(x(D;\beta)/b)$ is homotopic to the identity of the torus.
Indeed the first term is homotopic by the actual flow and the
second has a continuous real lift. This degree-one map is onto.
Choose a preimage of zero. This gives the claimed orbit and
exponential decay at the specified time.

In the full nonsymmetric graph every vertex has bounded degree and
each primary Hessian entry is bounded by $G_n|x_n|$ or
$G_n|x_{n+1}|$, by the binomial ratios. Indeed the entries of the
edge with parent $(n,k)$ are its coefficient times the amplitude of
the third vertex, and
\[
 \frac{c_{nk}^2}{c_{n+1,k}c_{n+1,k+1}}=\frac{\sqrt{(n-k)(k+1)}}{n},
 \qquad
 \frac{c_{nk}}{c_{n+1,k}}=\sqrt{\frac{n-k}{n}},
 \qquad
 \frac{c_{nk}}{c_{n+1,k+1}}=\sqrt{\frac{k+1}{n}},
\]
by $\binom nk=\frac n{n-k}\binom{n-1}k$ and
$\binom n{k+1}=\frac n{k+1}\binom{n-1}k$. All three ratios are at
most one for $0\le k\le n-1$. In a private row the
large $(\kappa/b)x_n$ multiplies the variations of partner and
axis into one another in a skew-Hermitian block, including both
complex quadratures; hold the primary fixed in the real edge
Hamiltonian to see this. All its other entries cost $C\kappa$.
The two-sided early bound follows from \eqref{rate:eq:l1G}. At the
endpoint the top daughters have hyperbolic rates $\pm1$, the
private primary/companion pairs $\pm\kappa$ in each quadrature,
and all axial coordinates are centers. The endpoint projections
have bounded norm. On the tail the derivative differs from the
constant endpoint generator by at most
\[
 C(G_{d-1}b+\kappa+b)(Y_1(D)/b)e^{-\kappa(t-D)/2}
   =o(\kappa)e^{-\kappa(t-D)/2}.
\]
For example the private difference between axis and partner
has its factor $c_{nk}$ canceled by the inverse factor in its
coefficient. The constant tail Green with stable data at $D$ and
its complement at the right endpoint has norm $C/\kappa$.
Neumann inversion and then backward propagation on $[0,D]$
give the result, uniformly if the right endpoint is finite.

We record these last steps in detail.
\emph{(1) The early propagator.} Split the variational generator as
\(\mathcal A(t)=\mathcal A_{\rm sk}(t)+\mathcal A_{\rm re}(t)\), where
\(\mathcal A_{\rm sk}\) consists of the blocks coupling the variations
of companion and axis in each private row. With \(\dot h=gvw\) and
\(\dot w=-gh\bar v\) such a block is
\(\delta\dot h=gv\,\delta w\), \(\delta\dot w=-g\bar v\,\delta h\), which
is skew-Hermitian; the blocks act on disjoint pairs, so the propagator
\(\mathcal U_{\rm sk}\) of \(\mathcal A_{\rm sk}\) is unitary. Writing the
full propagator as \(\mathcal U(t,s)=\mathcal U_{\rm sk}(t,s)\mathcal V(t,s)\),
\(\mathcal V\) solves
\(\partial_t\mathcal V=\mathcal U_{\rm sk}(s,t)\mathcal A_{\rm re}(t)
\mathcal U_{\rm sk}(t,s)\mathcal V\), so
\(\|\mathcal U(t,s)\|\le\exp|\int_s^t\|\mathcal A_{\rm re}\||\) in both
time directions. The entries of \(\mathcal A_{\rm re}\) are the primary
Hessian entries, at most \(G_n|x_n|\) or \(G_n|x_{n+1}|\), the top
entries, at most \(2\sqrt2\), and the remaining private entries, at most
\(C\kappa\). Every vertex has bounded degree, so the Schur test gives
\(\|\mathcal A_{\rm re}(t)\|\le C\sum_{n<d}G_n(|x_n|+|x_{n+1}|)+C\). By
\eqref{rate:eq:l1G}, \(\int_0^D\|\mathcal A_{\rm re}\|\le C'dD\), and
\(\|\mathcal U(t,s)\|\le e^{C'dD}\) for \(s,t\in[0,D]\).
\emph{(2) The tail perturbation.} For \(t\ge D\) set
\(\eta(t)=C(Y_1(D)/b)e^{-\kappa(t-D)/2}\), so that
\(\|(\xi,\zeta)\|_2\le\eta(t)\le1\). At \((n,k)\) the primary and
companion amplitudes are \(c_{nk}b\xi_n\) and \(c_{nk}b\zeta_n\).
The private axis differs from \(c_{nk}b\) by
\(O(c_{nk}b\eta(t)^2)\), and \(X-1=O(b^2\eta(t)^2)\), by the two
square-root formulas. On a primary edge the change in each Hessian
entry is therefore bounded by \(G_n b\eta(t)\): multiplying its
coefficient by the amplitude of the third vertex gives one of the
three binomial ratios computed above, all at most one. On a private
edge the coefficient \(\kappa/(bc_{nk})\) cancels the common factor
\(bc_{nk}\), giving at most \(C\kappa\eta(t)\); on the top edge the
bound is \(Cb\eta(t)\), since \(c_{1,0}=1\). The bounded vertex degree
and the Schur test now give \(C(G_{d-1}b+\kappa+b)\eta(t)\) in the
ordinary Euclidean operator norm. At \(t=D\) this bound is \(o(\kappa)\): with
\(Y_1(D)\le Cde^{-c_0D}+Cd^{3/2}\kappa bD\) and \(b=e^{-c_1D}\),
\(c_1<c_0/4\), the three products are
\(G_{d-1}Y_1(D)\), \(\kappa Y_1(D)/b\) and \(Y_1(D)\), each \(o(\kappa)\)
by the choice of \(c_1,C_*\) and \(\kappa D\to0\).
\emph{(3) The constant tail Green.} At the endpoint the top daughters
form a pair with rates \(\pm1\) and orthogonal eigenvectors, each private
primary/companion pair has rates \(\pm\kappa\) in the coordinates
\((\xi,\zeta)\) with orthogonal eigenvectors, in each quadrature, and the
top and the axes have zero generator, since their equations are products
of two vanishing amplitudes. All spectral projections therefore have
norm one. For \(\kappa/8\le\sigma\le\kappa/3\), in the norm
\(\sup_t e^{\sigma(t-D)}|\cdot|\), the forward kernel on the stable part
is at most \(e^{-(\kappa-\sigma)(t-s)}\) and the backward kernel on the
unstable and center parts at most \(e^{-\sigma(s-t)}\), so the Green
operator has norm at most \((\kappa-\sigma)^{-1}+\sigma^{-1}\le10/\kappa\).
\emph{(4) Assembly.} By (2) and (3) the tail perturbation contributes an
operator of norm \((10/\kappa)\,o(\kappa)<1/2\), and the Neumann series
gives the tail inverse with norm at most \(20/\kappa\). On \([0,D]\) the
solution is propagated backward from \(D\) with the bound of (1), and
the forcing on \([0,D]\) contributes at most \(De^{C'dD}\). The total
norm is at most \(e^{C'dD}(D+20/\kappa)\le e^{CdD}\), and the same
argument on \([0,L]\), with the complementary data at \(L\), gives the
bound uniformly in \(L\).
\end{proof}

\begin{lemma}\label{rate:l:round}
The finite mixed inverse and landing persist, with a correction
$e^{-C_2 dD}$ in ordinary Euclidean norm, under a perturbation
of the finite field of size $e^{-C_3dD}$ in $C^1$, and of the
prescribed axial endpoint of size $e^{-C_3dD}$. Here $C_3$ may
be chosen after any sufficiently large $C_2$. The actual field
in the carrier coordinates must vanish at its one specified
perturbed endpoint. Time-dependent perturbations are allowed.
\end{lemma}
\begin{proof}
The quadratic field $F_0$ has Hessian norm $L\le e^{CD}$ and
vanishes on the axial plane. If $z_*$ tends to $a$ as above,
then $\|z_*-a\|_\sigma\le e^{Cd}$ in the early-unweighted,
late-exponential norm at any of the indicated weights. Translate
$z_*$ by the endpoint displacement $\Delta$. The ideal residual
is $D^2F_0(z_*-a,\Delta)$. The remaining residual is estimated
by comparing $F(t,z_*+\Delta)-F_0(z_*+\Delta)$ to its value
zero at $a+\Delta$. Its weighted size is thus at most
$e^{Cd}(e^{-C_3dD}+L|\Delta|)$. Use the mixed Green of
Lemma~\ref{rate:l:finite} on a correction ball $e^{-C_2dD}$.
Its product with $O(Le^{-C_2dD}+L|\Delta|+e^{-C_3dD})$
is less than $1/2$ in the stated order of choices. This also
gives the perturbed variational inverse on the same boundary
traces by a Neumann series. The estimate compares the actual
field at the one prescribed endpoint; a different axial action
generally has a nonzero phase velocity in this frame.

In detail, write the perturbed field as $F(t,z)=F_0(z)+E(t,z)$ with
$\|E\|_{C^1}\le e^{-C_3dD}$ on the tube, let $\Delta$ be the axial
displacement of the endpoint, $|\Delta|\le e^{-C_3dD}$, and assume
$F(t,a+\Delta)=0$. Since $F_0$ is quadratic and vanishes on the axial
plane, $F_0(a+s\Delta)=0$ for all real $s$ gives $DF_0(a)\Delta=0$ and
$D^2F_0(\Delta,\Delta)=0$. Hence
\[
 F_0(z_*+\Delta)-F_0(z_*)=DF_0(z_*)\Delta=D^2F_0(z_*-a,\Delta),
\]
and $F_0(a+\Delta)=0$ forces $E(t,a+\Delta)=0$. Seek
$z=z_*+\Delta+\zeta$. Then
\[
 \zeta'-DF_0(z_*)\zeta
 =D^2F_0(z_*-a,\Delta)+D^2F_0(\Delta,\zeta)+\tfrac12D^2F_0(\zeta,\zeta)
  +E(t,z_*+\Delta+\zeta),
\]
with $|E(t,z_*+\Delta+\zeta)|\le e^{-C_3dD}(|z_*-a|+|\zeta|)$. In the
weighted norm $\|\cdot\|_\sigma$, whose weight is at least one, the right
side is at most
\[
 Le^{Cd}|\Delta|+L|\Delta|\,\|\zeta\|_\sigma+L\|\zeta\|_\sigma^2
 +e^{-C_3dD}\bigl(e^{Cd}+\|\zeta\|_\sigma\bigr),
\]
by $\|z_*-a\|_\sigma\le e^{Cd}$. Let $\mathcal G$ be the mixed inverse of
Lemma~\ref{rate:l:finite}, with $\|\mathcal G\|\le e^{CdD}$, and consider
$\zeta\mapsto\mathcal G(\text{right side})$ on the ball
$\|\zeta\|_\sigma\le e^{-C_2dD}$. Since $L\le e^{CD}$ and $d,D\ge1$, this
map sends the ball into itself and contracts with constant $1/2$ provided
\[
 2e^{3CdD}e^{-C_3dD}\le\tfrac12e^{-C_2dD},\qquad
 4e^{2CdD}e^{-C_2dD}\le\tfrac12 ;
\]
here the second condition bounds both the Lipschitz constant
\(e^{CdD}(L|\Delta|+2Le^{-C_2dD}+e^{-C_3dD})\) and the terms of the image
containing \(\zeta\). Both hold if $C_2\ge2C+3$ and then
$C_3\ge C_2+3C+2$. The fixed point is
the perturbed orbit: it lands at $a+\Delta$, has the prescribed stable
data at $D$, and lies within $e^{-C_2dD}$ of $z_*+\Delta$. The same
Neumann series inverts the variational equation along it.
\end{proof}

At the cold end the actual $H^s$ cost of the two ideal leaf layers
for Fourier amplitude $r=H N^{-s}$ is at most
$C_sH[(4/5)^s+(1/5)^s]^{(d-1)/2}$. Indeed the leaves of the branch
of $P$ are the vertices $Pg_{d-1,k}$, with squared amplitudes
$\binom{d-1}kx_d^2$ and
$|g_{d-1,k}|^2=\alpha^{d-1-k}(1-\alpha)^k$, where
$\alpha=|u|^2=(1+t^2)^{-1}$ and $1-\alpha=|1-u|^2$. Hence
\[
 \sum_{k=0}^{d-1}\binom{d-1}k|Pg_{d-1,k}|^{2s}
 =|P|^{2s}\bigl(\alpha^s+(1-\alpha)^s\bigr)^{d-1},
\]
and similarly for $Q$. For $1/2<t<3/4$ one has $16/25<\alpha<4/5$,
and $\alpha^s+(1-\alpha)^s$ increases on $[1/2,1]$. At all times the primary
cost is $C_sH$ by the same binomial expansion on each level and
the invariant sum of $x_n^2$. The actual physical frequencies are
$e_1+Nv$ in normalized graph units: for the primaries the relative
error in the homogeneous weight is $O_s(C^d/N)$. The two
subtree radii are bounded away from zero. All private frequencies
have normalized modulus $O(d)$ by the real angular cutoff, and their
entire strong cost is at most $C_s H d^{C_s}2^{d/2}e^{-c_1D}$.
A finite Euclidean correction $\eta$ has physical strong cost at
most $C_s H d^{C_s}\eta$, without multiplying this bill by a
previous-scale preparation price. With $H=e^{b_*d}$ and
$0<b_*<-\log[(4/5)^s+(1/5)^s]/2$, the strong initial price
is $O(e^{-c_sd})$ and the prescribed terminal top has height $H$.

\subsection{The marked axial Fourier coefficient}\label{rate:sec:axial-mark}
Use $x\in\mathbb R/(2\pi\mathbb Z)$, normalized Haar measure and
$i u_t=-u_{xx}+|u|^2u$. Write $\Psi$ for the inverse periodic
defocusing NLS Birkhoff map, normalized so its linear part is inverse
Fourier. Complexify the two conjugate Cartesian variables independently,
with actions $I_j=z_j^+z_j^-$. A positive input has physical spatial
factor $e^{ijx}$ and a negative input $e^{-ijx}$.

\subsubsection{Normalizations and auxiliary Hamiltonians}
We use the normalized one-dimensional chart of
Section~\ref{sec:tree}. The hierarchy Hamiltonians depend only
on the actions \cite[Theorem~1.1(3)]{GKsym}.
At every integer order $m\ge1$ each nonlinear chart remainder
gains one derivative \cite[Theorem~1.1 and Corollary~1.1]{KST}.
In the frequency map of \cite[Theorem~18.1 and Corollary~18.2]{Molnar}
(with Appendix~\ref{app:molnar-four} at \(p=4\)) we label modes by their physical Fourier indices, as in
\eqref{eq:renormalized-frequency}; in particular $D\chi(0)=-\mathrm{Id}$.
For the period-one spectral operator below put $a=2\pi/\sqrt2$.

We identify precisely the hierarchy normalization used below. Put
\[
 J=\tfrac i2\operatorname{diag}(1,-1),\qquad
 V=2^{-1/2}\begin{pmatrix}0&u\\ v&0\end{pmatrix},\qquad
 A=\zeta J+V ,
\]
where now $v$ denotes the independent physical conjugate field.
There is a unique traceless formal series $H=J+\sum_{j\ge1}H_j\zeta^{-j}$
with $H_x=[A,H]$, $H^2=-\mathrm{Id}/4$.
Indeed $[J,H_{j+1}]=(H_j)_x-[V,H_j]$ fixes the off-diagonal part,
and the square fixes the diagonal part. Existence follows recursively
by formally diagonalizing $\partial_x-A$: use at each power an
off-diagonal conjugating coefficient to remove that power's off-diagonal
part using $\operatorname{ad}J$, and conjugate $J$ back by the resulting
near-identity formal series. Each coefficient is a differential
polynomial. Give the fields and $\partial_x$ weight one. Then $H_j$
has weight $j$ and the vector field
\[
 V_{t_p}=[J,H_{p+1}]=(\partial_x-\operatorname{ad}A)(\zeta^pH)_+,
 \qquad p=0,1,\ldots ,
\]
has weight $p+1$. Here $(\cdot)_+$ retains the nonnegative powers
of $\zeta$, including the constant term. The off-diagonal monomials have odd field degree $g$
by charge symmetry, with $p+1-g$ spatial derivatives. The single upper
wave has rate $+i n^p$ and the lower wave has rate $-i n^p$ in these
times, as also follows from the dressing below. In particular
$u_{t_2}=-iu_{xx}+iu^2v$, the reverse of physical NLS time.
To check this normalization, note that $\operatorname{ad}J$
multiplies the upper and lower entries of an off-diagonal matrix by
$i$ and $-i$, and annihilates diagonal matrices. The recursion and
the coefficients of $\zeta^{-1},\zeta^{-2}$ in $H^2=-\mathrm{Id}/4$
give $H_1=V$ and
\[
 H_2=\frac1{\sqrt2}\begin{pmatrix}0&-iu_x\\ iv_x&0\end{pmatrix}
     +\frac{iuv}2\operatorname{diag}(1,-1),
\]
because $2JH_2^{\rm diag}+V^2=0$ and $V^2=\tfrac12uv\,\mathrm{Id}$.
The upper entry of $[J,H_3]=(H_2)_x-[V,H_2]$ is then
$2^{-1/2}(-iu_{xx})-2^{-1/2}u\,(-\tfrac{iuv}2-\tfrac{iuv}2)
=2^{-1/2}(-iu_{xx}+iu^2v)$, and the upper entry of $V_{t_2}$ is
$2^{-1/2}u_{t_2}$. With $v=\bar u$ this is the equation
$iu_t=-u_{xx}+|u|^2u$ in the time $t=-t_2$.

For an independent check that these polynomial generators preserve
the actions of the specific chart just cited, let $S$ be the formal
diagonalizer, $H=SJS^{-1}$ and
$A_d=S^{-1}AS-S^{-1}S_x$. For $Y=S^{-1}\delta S$,
$\delta A_d=S^{-1}\delta V S+[A_d,Y]-Y_x$. Hence
\[
 \delta\int\mathrm{tr}(J A_d)\,dx
     =\int\mathrm{tr}(H\,\delta V)\,dx .
\]
Taking the $\zeta^{-p-1}$ coefficient shows that $[J,H_{p+1}]$ is a
canonical Hamiltonian vector field, with a fixed constant set by the
choice of $u,v$ and symplectic form. Its Hamiltonian density may be
chosen by integrating this gradient radially in $V$; a degree-$g$
Hamiltonian term then has $p+2-g$ spatial derivatives. The polynomial
zero-curvature identity makes the derivative of the one-period monodromy
a commutator, since its time matrix is periodic. Thus the entire
Zakharov--Shabat discriminant and all its contour actions are invariant.
Explicitly, in the period-one system
$L=i\sigma_3\partial_y+\left(\begin{smallmatrix}0&a u\\a v&0\end{smallmatrix}\right)$,
with $\sigma_3=\operatorname{diag}(1,-1)$, set $x=2\pi y$,
$\lambda_{\mathrm{sp}}=-\pi\zeta$ and $f=\operatorname{diag}(i,1)g$.
Then $Lf=\lambda_{\mathrm{sp}} f$ is exactly $g_x=(\zeta J+V)g$.
Indeed, with $f=(ig_1,g_2)$, $\partial_y=2\pi\partial_x$ and
$a/(2\pi)=2^{-1/2}$, the two rows of $Lf=-\pi\zeta f$ read
\[
 -2\pi g_{1,x}+a\,ug_2=-i\pi\zeta g_1,\qquad
 -2\pi i\,g_{2,x}+ia\,vg_1=-\pi\zeta g_2,
\]
that is $g_{1,x}=\tfrac{i\zeta}2g_1+2^{-1/2}ug_2$ and
$g_{2,x}=-\tfrac{i\zeta}2g_2+2^{-1/2}vg_1$.
With the period-one fundamental matrix normalized to the identity at
\(y=0\), let \(\Delta(\lambda_{\mathrm{sp}})\) be the trace of its monodromy.
By \cite[(2.2)]{Molnar}, the period-one action at spectral index \(j\) is
\[
 I_j^{(1)}=\frac1\pi\oint_{\Gamma_j}
 \frac{\lambda_{\mathrm{sp}}\Delta'(\lambda_{\mathrm{sp}})}
      {\sqrt{\Delta(\lambda_{\mathrm{sp}})^2-4}}\,
 d\lambda_{\mathrm{sp}},
\]
where \(\Gamma_j\) encloses the \(j\)th gap counterclockwise and the
square root is the canonical branch. The corresponding physical action is
\(I_{-j}=2(2\pi)^{-2}I_j^{(1)}\) by the scaling in Section~\ref{sec:tree}. It follows that the Hamiltonian commutes with
every $I_j$ of this chart, so its restriction is invariant under each
Cartesian coordinate rotation and depends only on the actions.
All arguments can first be made on smooth real potentials, then as
analytic identities of finite-parameter complex germs. For an individual
$p$ one works initially at sufficiently high Sobolev order. Only finitely
many $p$'s will be used in any coefficient estimate.

\subsubsection{Square-free coefficients and a reciprocal-moment template}
For disjoint finite labels $P,Q$, $|P|=m+1,\ |Q|=m$, let $C(P,Q)$ be
the coefficient of $\prod_{P}z_p^+\prod_Q z_q^-$ in the physical upper
field, after factoring out its spatial exponential. Then
\begin{equation}\label{rate:eq:axial-path}
 C(P,Q)=(-1/2)^m\sum_{\substack{p_1,\ldots,p_{m+1}\\q_1,\ldots,q_m}}
   \prod_{j=1}^m\frac1{(p_j-q_j)(p_{j+1}-q_j)} ,
\end{equation}
where the sum runs over the permutations of each set. To verify the
formula even at kinetic resonance, flow the parameters independently as
$A_p=a_p\exp(i\sum t_jp^j)$ and $B_q=b_q\exp(-i\sum t_jq^j)$, with
$t_1=x$ and the zeroth character identically one. Let $A_P,B_Q$ be
their diagonal matrices and $C_{pq}=(p-q)^{-1}$. Take
\[
 D(\zeta)=1+\sum_p\frac{r_pe_2^t}{\zeta-p}
             +\sum_q\frac{s_qe_1^t}{\zeta-q},\qquad
 r_p=\frac{iA_p}{\sqrt2}\Big(e_1+\sum_q C_{pq}s_q\Big),\quad
 s_q=\frac{-iB_q}{\sqrt2}\Big(e_2-\sum_p C_{pq}r_p\Big).
\]
This finite linear system is invertible for small parameters.
At $p$ the matrix $D(1-N/(\zeta-p))$ is regular with
$N=(iA_p/\sqrt2)e_1e_2^t$, and analogously at $q$ with
$N=(-iB_q/\sqrt2)e_2e_1^t$. Consequently $\det D=1$.
The potential pole in
$D_{t_j}D^{-1}+\zeta^jDJD^{-1}$ is removable since
$N_{t_j}=p^j[J,N]$ at $p$, $NJN=0$ and $N^2=0$; at $q$ use $q^j$.
The resulting polynomial is exactly $(\zeta^jDJD^{-1})_+$ by expansion
at infinity. Its first spatial coefficient is $V=[D_1,J]$,
$D_1=\lim_{\zeta\to\infty}\zeta(D-1)$.
The upper physical field is therefore
\[
 u=\mathbf1^t(1+A_P C B_Q C^t/2)^{-1} A_P\mathbf1.
\]
In detail, $D_1=\sum_pr_pe_2^t+\sum_qs_qe_1^t$, and the $(1,2)$ entry
of $[D_1,J]$ is $-i(D_1)_{12}$. Since $V_{12}=2^{-1/2}u$, this gives
$u=-i\sqrt2\sum_p\rho_p$, where $\rho_p$ and $\sigma_q$ denote the
first components of $r_p$ and $s_q$. The first components of the
linear system are
\[
 \rho_p=\frac{iA_p}{\sqrt2}\Big(1+\sum_qC_{pq}\sigma_q\Big),\qquad
 \sigma_q=\frac{iB_q}{\sqrt2}\sum_{p'}C_{p'q}\rho_{p'} .
\]
Eliminating $\sigma$ gives
$(1+A_PCB_QC^t/2)\rho=(i/\sqrt2)A_P\mathbf1$ for the vector
$\rho=(\rho_p)$, and the displayed formula follows.
These are smooth periodic potentials for complex parameters near zero.
Their spatial fundamental solution away from the finitely many poles is
\[
  D(x,\zeta)e^{x\zeta J}D(0,\zeta)^{-1}.
\]
 Periodicity
of $D$ shows that its monodromy matrix is conjugate to
$e^{2\pi\zeta J}$. By analyticity the entire discriminant coincides
with that of the vacuum, including at every pole. The contour action
formula with its vacuum normalization therefore gives zero for each
action of the actual Birkhoff chart. Equivalently, invariance of an
action under every hierarchy time forces every moment of
$\sum\alpha_p\delta_p-\sum\beta_q\delta_q$ to vanish on each of its
parameter monomials. Vandermonde on the disjoint supports gives the
same zero-action conclusion.
On this complex zero-action face a particular positive coordinate
evolves at its linear rate $+i k^j$, and a negative one at $-i k^j$.
For completeness, on each finite Cartesian truncation rotational
invariance makes the Hamiltonian a power series in the pair products.
The quotient of its $z_k^-$ derivative by $z_k^+$ on the zero-action
face thus has its vacuum value. Truncating in an arbitrarily high
Sobolev order and passing to the smooth limit gives the same identity
for each coordinate. A second Vandermonde test, subtracting or adding
$\delta_k$, leaves only $a_k$ or $b_k$ with coefficient one from the
free normalization. Thus the dressing is the actual inverse chart on
this disjoint complex subspace, and its Neumann expansion gives
\eqref{rate:eq:axial-path}. An action-dependent rotation of coordinate
angles with unchanged vacuum differential has no effect on this
zero-action computation.

\smallskip
\noindent\emph{The dressing identity step by step.}
(a) \emph{Poles and determinant.} Near \(\zeta=p\) write
\(D(\zeta)=\widetilde D(\zeta)+r_pe_2^t/(\zeta-p)\) with \(\widetilde D\) regular at \(p\),
and let \(N=(iA_p/\sqrt2)e_1e_2^t\). Since \(e_2^te_1=0\), one has \(N^2=0\)
and \(r_pe_2^tN=0\), whence
\[
 D(\zeta)\Bigl(1-\frac N{\zeta-p}\Bigr)
 =\widetilde D(\zeta)-\frac{\widetilde D(\zeta)N-r_pe_2^t}{\zeta-p}.
\]
By the definitions of \(r_p\) and \(C_{pq}\),
\(\widetilde D(p)N=(iA_p/\sqrt2)\bigl(e_1+\sum_qC_{pq}s_q\bigr)e_2^t=r_pe_2^t\), so
the right side is regular at \(p\). The same computation with
\(N=(-iB_q/\sqrt2)e_2e_1^t\) and the definition of \(s_q\) gives regularity
at \(q\). Each factor \(1-N/(\zeta-p)\) has determinant one, so \(\det D\) is
entire, tends to one at infinity, and is identically one.
(b) \emph{Zero curvature.} Put \(\Phi(\zeta)=D(\zeta)\exp(\sum_jt_j\zeta^jJ)\),
with \(t_1=x\). Near \(p\) write \(D=R\,(1+N/(\zeta-p))\), with \(R\) regular and
invertible by (a) and \(N^2=0\). The polar part of
\(\Phi_{t_j}\Phi^{-1}=D_{t_j}D^{-1}+\zeta^jDJD^{-1}\) is
\[
 \frac{R\bigl(N_{t_j}-\zeta^j[J,N]\bigr)R^{-1}}{\zeta-p}
 -\frac{R\bigl(N_{t_j}N+\zeta^jNJN\bigr)R^{-1}}{(\zeta-p)^2} .
\]
Here \(N_{t_j}=ip^jN=p^j[J,N]\), because \(\partial_{t_j}A_p=ip^jA_p\) and
\([J,N]=iN\); moreover \(N_{t_j}N=0\) and \(NJN=0\), since \(J\) is diagonal.
The polar part is therefore \(R(p^j-\zeta^j)[J,N]R^{-1}/(\zeta-p)\), which
is regular, and likewise at every \(q\). Thus \(B_j\coloneqq\Phi_{t_j}\Phi^{-1}\)
is a polynomial in \(\zeta\), equal to the polynomial part of its
expansion at infinity. As \(D=1+D_1\zeta^{-1}+O(\zeta^{-2})\) and
\(D_{t_j}D^{-1}=O(\zeta^{-1})\), we get \(B_j=(\zeta^jH)_+\) with \(H=DJD^{-1}\),
and \(B_1=\zeta J+[D_1,J]\). Since \(H\) is traceless, \(H^2=-\mathrm{Id}/4\),
\(H=J+O(\zeta^{-1})\) and \(H_x=[B_1,H]\), it is the formal series of the
normalization paragraph with \(V=[D_1,J]\). The compatibility of
\(\Phi_x=B_1\Phi\) and \(\Phi_{t_j}=B_j\Phi\) is
\(V_{t_j}=(\partial_x-\operatorname{ad}B_1)(\zeta^jH)_+\), the \(j\)th
hierarchy flow.
(c) \emph{Actions.} Each action \(I_k\) of the Birkhoff chart is invariant
under every hierarchy flow, since the hierarchy Hamiltonians depend only
on the actions \cite[Theorem~1.1(3)]{GKsym}, and under the common phase
\cite[Theorem~1.2(i)]{GKsym}. Along the flows of (b) the parameters move
as \(a_p\mapsto a_pe^{i\sum_jt_jp^j}\), \(b_q\mapsto b_qe^{-i\sum_jt_jq^j}\), and
the common phase multiplies every \(a_p\) by \(e^{i\vartheta}\) and every \(b_q\)
by \(e^{-i\vartheta}\). Hence every Taylor monomial
\(\prod a_p^{\alpha_p}\prod b_q^{\beta_q}\) of the analytic function \(I_k(a,b)\)
satisfies \(\sum_p\alpha_pp^j-\sum_q\beta_qq^j=0\) for
\(0\le j<|P|+|Q|\). The measure \(\sum\alpha_p\delta_p-\sum\beta_q\delta_q\) on
the \(|P|+|Q|\) distinct points of \(P\cup Q\) thus has vanishing moments
of all orders below the number of points, so it vanishes by the
Vandermonde determinant. Only the constant monomial survives, and
\(I_k\equiv I_k(0,0)=0\).
(d) \emph{Coordinates.} On the zero-action face the \(j\)th flow multiplies
\(z_k^\pm\) by \(e^{\pm ik^jt_j}\), the vacuum rate being read off from the
linear part. Each Taylor monomial of \(z_k^+(a,b)\) therefore satisfies
\(\sum_p\alpha_pp^j-\sum_q\beta_qq^j=k^j\) for \(0\le j\le|P|+|Q|\), and the same
Vandermonde argument with the point \(k\) adjoined gives \(\alpha=e_k\),
\(\beta=0\) if \(k\in P\), and no monomial otherwise. Comparing linear parts,
\(z_k^+=a_k\) on \(P\) and \(z_k^+=0\) off \(P\); likewise \(z_q^-=b_q\) on \(Q\) and
\(z_q^-=0\) off \(Q\). The dressed potential is thus the image under the
inverse chart of the Cartesian point with these coordinates, and the
Neumann expansion of the displayed formula for \(u\) gives
\eqref{rate:eq:axial-path}.

\begin{lemma}\label{rate:lem:axial-template}
For each fixed integer $h\ge2$ there are distinct nonzero integers
$v_i$ with signs $e_i=\pm1$ such that
\begin{gather*}
 \sum e_i=1,\qquad \sum e_i v_i=\sum e_i v_i^2=0,\qquad
 \sum e_i/v_i\ne0,\\
 \sum e_i v_i^{-j}=0\quad(2\le j\le h),\qquad C(v,e)\ne0.
\end{gather*}
Five designated entries have the signs and offsets of
$(-1,3,4\,;\,1,5)$ up to a common dilation, positives listed first.
The other entries contain both signs.
\end{lemma}
\begin{proof}
For this five-site starting set $S$, direct substitution of the
twelve terms in \eqref{rate:eq:axial-path} gives $C_S=-1/9$ and
$\sum e_i/v_i=-97/60$. The twelve terms were listed in
Section~\ref{sec:super}; the other values are
\[
 3-2=1,\qquad -1+3+4-1-5=0,\qquad 1+9+16-1-25=0,\qquad
 -1+\tfrac13+\tfrac14-1-\tfrac15=-\tfrac{97}{60}.
\]
Suppose the result has been obtained through
$j=m-1$, $m\ge2$, using rationals. Write
$x_T=N_0+\sum_{j\in T}h_j$ for $T\subset\{1,\ldots,m\}$ with
$m!\prod h_j=1$. Adjoin to $S$ a copy $S/x_T$ with signs multiplied by
$w_T=-(-1)^{m-|T|}$ for each $T$. The $m$th divided difference gives
$\sum w_T x_T^j=0$ for $0\le j<m$ and $-1$ for $j=m$.
Indeed, with $\Delta_hf(x)=f(x+h)-f(x)$,
\[
 \sum_T(-1)^{m-|T|}f(x_T)=\Delta_{h_1}\cdots\Delta_{h_m}f(N_0),
\]
which vanishes for polynomials of degree less than $m$ and equals
$m!\prod_jh_j=1$ for $f(x)=x^m$.
The previous inverse moments are multiplied by $1+\sum w_T x_T^j$,
and in each copy the two positive moments vanish separately.
In detail the copy $S/x_T$ consists of the points $v_i/x_T$ with
signs $e_iw_T$, so that for every integer $j$
\[
 \sum_ie_iv_i^{-j}+\sum_Tw_T\sum_ie_i(v_i/x_T)^{-j}
 =\Bigl(1+\sum_Tw_Tx_T^{\,j}\Bigr)\sum_ie_iv_i^{-j}.
\]
The factor is $1$ for $0\le j<m$ and $0$ for $j=m$. Thus the charge,
the nonzero first inverse moment and the vanishing inverse moments of
orders $2\le j<m$ are preserved, the inverse moment of order $m$
becomes zero, and the two positive moments $\sum e_iv_i$ and
$\sum e_iv_i^2$ vanish on every copy by homogeneity.

It remains to preserve the nonzero path coefficient. Put
$h_1=t^{m-1},\ h_j=\Gamma_j/t$ for $j\ge2$, $N_0=\Gamma_0/t$, with
positive rationals $\Gamma_j$ satisfying the product constraint.
Avoid coincident subset sums among $\Gamma_2,\ldots,\Gamma_m$.
For $m>2$ each unwanted equality is a nonzero Laurent polynomial
after eliminating the last variable by the product (the resulting
monomials are distinct); the positive rational parameter set is
dense. For $m=2$ distinctness holds directly. Next vary $\Gamma_0$
to separate the numbers $v/(\Gamma_0+\sum_{j\in T}\Gamma_j)$ over
$v\in S$ and $T\subset\{2,\ldots,m\}$. For different $v$ only finitely
many values are excluded. For fixed $v$ use the distinct subset sums.
The new points are opposing pairs $(y_i+\delta_i,y_i)$, positive first,
with $|y_i|\lesssim t$, mutual center separations $\gtrsim t$ and
$0<|\delta_i|\lesssim t^{m+1}$. Indeed the product constraint reads
$m!\prod_{j\ge2}\Gamma_j=1$ and does not involve $t$. For
$T\subset\{2,\ldots,m\}$ put $\Gamma_T=\Gamma_0+\sum_{j\in T}\Gamma_j$,
so that $x_T=\Gamma_T/t$ and $x_{T\cup\{1\}}=x_T+t^{m-1}$, with
$w_{T\cup\{1\}}=-w_T$. The two copies of a point $v\in S$ are
\[
 \frac v{x_T}=\frac{tv}{\Gamma_T},\qquad
 \frac v{x_T+t^{m-1}}=\frac{tv}{\Gamma_T}
   -\frac{t^{m+1}v}{\Gamma_T(\Gamma_T+t^m)},
\]
with opposite signs. The centers $tv/\Gamma_T$ are distinct multiples
of $t$ by the choice of the $\Gamma_j$, and the offsets have size
$t^{m+1}$, which is $o(t^2)$ because $m\ge2$.

Here are the collision estimates for \eqref{rate:eq:axial-path}.
For one opposing pair,
$\lim_{\delta\to0}\delta\,C(S\cup(y+\delta,y))=C_S L_S(y)$,
where $L_S(y)=\sum_{v\in S}e_v/(y-v)$.
Only neighboring occurrences of the two added letters contribute.
If inserted between an old positive $a$ and old negative $b$ the
ratio of edge factors after removing their common pole is
$1/(a-y)+1/(y-b)$, in either order of this old edge. Each end insertion
contributes $1/(a-y)$ at its positive endpoint. Counting both
incidences of each old vertex and multiplying by the new factor
$-1/2$ gives exactly $L_S(y)$.
For $b$ simultaneous pairs put $\Delta=\min_{i\ne j}|y_i-y_j|$ capped
by one and a fixed small distance from $S$, and $\delta_*=\max_i|\delta_i|$.
The function $G=(\prod_i\delta_i)C$ is holomorphic on
$|\delta_i|\le c\Delta$. On the distinguished boundary it is
$O_S(\Delta^{1-b})$, since a path has at most $2b-1$ edges entirely
among the new vertices. Cauchy's estimate bounds its sum of terms
involving every variable by
$O_S(\Delta^{1-b}(\delta_*/\Delta)^b)=o(1)$ if
$\delta_*=o(\Delta^2)$. By inclusion-exclusion the remaining terms
are the signed sum of proper restrictions setting some $\delta_i=0$.
The single residue formula on each such restriction gives the analogous
function on the surviving pairs times, for each collapsed pair,
\[
 L_S(y_i)+\sum_{\text{surviving }j}
 \left[(y_i-y_j-\delta_j)^{-1}-(y_i-y_j)^{-1}\right]
 =L_S(y_i)+O_S(\delta_*/\Delta^2).
\]
Induction on $b$ therefore gives
$G=C_S\prod_i L_S(y_i)+o(1)$, with the case of one pair as above.

\smallskip
\noindent\emph{The collision estimate in detail.} Let
$|\delta_i|\le c\Delta$ with $c\le1/4$. An edge factor joining sites of
two different new pairs has modulus at most $2/\Delta$, an edge factor
joining a new site to a site of $S$ is $O_S(1)$, and only the edge
inside pair $i$ carries the pole $\delta_i^{-1}$. A path visits every
site once, so it uses each inner edge at most once and has at most
$2b-1$ edges between new sites. If it uses $k$ inner edges, its
contribution to $G$ on the distinguished boundary $|\delta_i|=c\Delta$
is therefore
$O_S\bigl((2/\Delta)^{2b-1-k}(c\Delta)^{b-k}\bigr)=O_S(\Delta^{1-b})$, the
factors $\delta_i$ of the unused inner edges being kept. Expand
$G=\sum_\alpha g_\alpha\delta^\alpha$. Cauchy's estimate
$|g_\alpha|\le\sup|G|\,(c\Delta)^{-|\alpha|}$ bounds the part $G_{\rm all}$
with every $\alpha_i\ge1$ by
$O_S\bigl(\Delta^{1-b}(\delta_*/\Delta)^b\bigr)
=O_S\bigl(\Delta(\delta_*/\Delta^2)^b\bigr)$ when $\delta_*\le c\Delta/2$.
The remaining monomials have some $\alpha_i=0$, and inclusion-exclusion
over the set $A$ of such indices gives
\[
 G-G_{\rm all}=\sum_{\emptyset\ne A\subset\{1,\ldots,b\}}(-1)^{|A|+1}\,
 G|_{\delta_A=0}.
\]
On $\delta_A=0$ the single residue formula, applied to the pairs of $A$
one at a time, gives
$G|_{\delta_A=0}=G_{A^c}\prod_{i\in A}\bigl(L_S(y_i)+O_S(\delta_*/\Delta^2)\bigr)$,
where $G_{A^c}$ is the corresponding function of the surviving pairs; an
already coalesced pair contributes zero to the signed sum. By induction
on $b$, $G_{A^c}=C_S\prod_{j\notin A}L_S(y_j)+o(1)$, and
$\sum_{\emptyset\ne A}(-1)^{|A|+1}=1$. Hence
$G=C_S\prod_iL_S(y_i)+o(1)$ as $\delta_*/\Delta^2\to0$. In the
application $\Delta\gtrsim t$, $\delta_*\lesssim t^{m+1}$ with $m\ge2$ and
$y_i=O(t)$, so that $\delta_*/\Delta^2=O(t^{m-1})\to0$ and
$L_S(y_i)\to L_S(0)\ne0$.
Since $L_S(0)=-\sum e_v/v\ne0$ and here $\delta_*=o(\Delta^2)$, a
sufficiently small rational $t>0$ works. Clear denominators at the end.
\end{proof}

\subsubsection{The coefficient on a nonzero occupied background}
Fix the template length $q$, and set $n_a=1+K V_a$ for all new sites
with distinct nonzero integers $V_a$, $\max|V_a|\le C d^C$. A particular
marked row has $V_a=l v_a$, $1\le l\le C d^C$. Write $F$ for the old sites, $n_F=(m)_{m\in F}$ and $f=|F|$.
The indices in the old support $F$ satisfy $|m|\le Z$, and their Cartesian amplitudes are in a
common small $h^1$ ball. For the physical upper field take its angular
coefficient with the row's signs on its marked sites, zero on all
others including the old sites, divided by the $q$ required Cartesian
radii. Its physical Fourier index is one by translation covariance.
The divided coefficient is extended analytically through vanishing
new radii. It is a single temporal/angular harmonic of that physical
Fourier coefficient. Other harmonics at the same spatial site can occur.

\begin{proposition}\label{rate:prop:axial-mark-true}
Suppose the new weak actions vanish and the other new actions have
sum $S\le K^{-(h+1)/2}$. Let $L_0$ be a fixed integer with
$(L_0+1)(h+1)/2>q+1$. There are fixed exponents, depending only
on $q,h,L_0$ and the fixed polynomial in the macro bound, for which
$E=(2+Z)^C e^{Cf}(2+d+\log K)^C$ has the following property.
For $K$ sufficiently large compared with these factors the divided
coefficient is
\[
 (Kl)^{-(q-1)}C(v,e)+
 O\big(E K^{-q}+ E S K^{-(q-1)}+(C S)^{L_0+1}\big).
\]
Since $l\le Cd^C$, we can increase $K$ so that
$l^{q-1}\big(E/K+ES+K^{q-1}(CS)^{L_0+1}\big)$, including the
implicit constant, is at most $|C(v,e)|/2$. The modulus is then
at least $\tfrac12(Kl)^{-(q-1)}|C(v,e)|$. If only five weak
radii are divided out and each remaining marked action is at least
$K^{-C_0}$, the resulting modulus is at least $K^{-C_1}$ for a fixed
$C_1$.
\end{proposition}
\begin{proof}
Set $M_0=q+2L_0$.
Expand at zero new Cartesian variables and keep up to degree $M_0$.
For a monomial $P$ write $l_a=\alpha_a-\beta_a$ and
$w=\sum|l_a|$. Let $u^\sigma_{P,j}$ be its coefficient at old angle
$j\in\mathbb Z^f$ in copy $\sigma=\pm1$. The spatial exponent is
$\nu_j=\sum l_a n_a+j\cdot n_F$ and the sign-copy spectral label
is $\sigma\nu_j$. Cauchy on a fixed old-angle strip and fixed Cartesian
discs in the at most $M_0$ tested variables gives
$|u^\sigma_{P,j}|\le C e^{-c|j|_1}$ by the analytic $L^2$ chart.
Shrinking the strip costs at most $C^f$ in the absolute angular sum.
Any prescribed spatial derivative adds a factor at most
$(C K d^C+Z|j|_1)$ per derivative. With $L=C'\log K$, the tail
$|j|_1>L$ can thus be made $O(K^{-A})$ for any preassigned finite
$A$, also after fixed derivatives, provided $K$ dominates $E$.
Constants $C$ in $E$ may increase a fixed number of times in the proof.
The identities involving high Hamiltonians are initially valid on a
sufficiently small higher-regularity ball for the finitely supported
arguments. They extend to the common low-norm polydiscs by angular
translation covariance and the convergent differentiated Fourier
series just described.

For a single new Cartesian argument the derivative of the nonlinear
part of $\Psi$ maps $h^1$ to $H^2$ uniformly on the old complex
$h^1$ ball. Its norm on the unweighted basis vector is $O(|n_a|)$.
For $|j|_1\le L$, $|\nu_j|\ge|n_a|/2$, so the error from its free
plane is at most $C|n_a|^{-1} e^{-c'|j|_1}$ by angular Cauchy;
the tail is absorbed above. In particular the positive
old-angle-zero entry on a positive test is $1+O(E/K)$.

For each $2\le p\le M_0+3$, write $\omega^{(p)}$ for the auxiliary
angular speeds, with free speed $n^p$. Differentiating its chart
identity at a single empty positive input and taking the upper,
old-angle-zero coefficient gives
$|\omega^{(p)}_{n_a}-n_a^p|\le E K^{p-2}$: the coefficient of this
frequency difference is $1+O(E/K)$, and each nonlinear physical
term has at most $p-2$ spatial derivatives. The mean of an old angular
derivative vanishes and new action differentiation requires a pair.
For the old speeds a bound $E$ suffices. Any fixed number of new action
derivatives of the tested old or new speeds costs $E K^p$. To see this,
use the action-only Hamiltonian and Cauchy's estimate in a bounded
number of tested action slots on fixed small discs. Realize a tested
empty pair as $(\rho,I_a/\rho)$, with a phase and its opposite,
using fixed small $\rho$. For a tested occupied index $b$, write
$\rho_b=\max(\sqrt{I_b},\rho)$ and use
$(\rho_b,(I_b+\eta_b)/\rho_b)$ on $|\eta_b|\le c\rho^2$.
Cauchy differentiation has a fixed cost also in this slot. The
squared Cartesian size of all tested inputs remains
in a common small $\ell^2$ ball. The untested old entries keep their
old angles. To bound the spatial integral defining the Hamiltonian,
take absolute angle sums on a smaller strip for all tested and old
angles. Each derivative then costs at most $K d^C+Z$ times an angular
moment, and a density term has at most $p$ derivatives in total.
With only old inputs use $Z$ instead of $K d^C$. There is no division
by a shrinking old action in this estimate.
For the first bound let \(e^+\) be the positive unit vector at the empty
site \(n_a\). The action \(I_{n_a}=z^+_{n_a}z^-_{n_a}\) has vanishing differential
at an empty site, so the linearization of the \(p\)th flow at the old
torus maps \(e^+\) to \(i\omega^{(p)}_{n_a}e^+\) in Birkhoff coordinates. Put
\(W(\phi)=D\Psi(z_F(\phi))e^+\), a function on the old angle torus.
Write \(X_{\mathrm B}^{(p)}(z)_m^\pm=\pm i\omega_m^{(p)}(I)z_m^\pm\)
for the vector field in Birkhoff coordinates. Differentiating the chart
identity \(D\Psi(z)X_{\mathrm B}^{(p)}(z)=X^{(p)}(\Psi(z))\) in
the direction \(e^+\) gives
\[
 i\omega^{(p)}_{n_a}W+\mathcal D_pW=DX^{(p)}\bigl(\Psi(z_F)\bigr)W,
 \qquad \mathcal D_p=\omega^{(p)}_F\cdot\partial_\phi .
\]
Take the upper spatial coefficient at \(n_a\) and old angle zero. The term
\(\mathcal D_pW\) contributes nothing, being an angular derivative. The
coefficient \(\widehat W(n_a,0)\) equals \(1+O(E/K)\) by the one-derivative
gain of the nonlinear part of \(D\Psi\), as shown above. The linear
part of \(DX^{(p)}\) gives \(in_a^p\widehat W(n_a,0)\). Each nonlinear term of
\(X^{(p)}\) has \(g\ge3\) factors and \(p+1-g\le p-2\) spatial derivatives; in
\(DX^{(p)}W\) one factor is \(W\) and the others are old, each derivative
costs at most \(CKd^C+Z|j|_1\), and the absolute angular sums cost \(E\).
Hence \(i(\omega^{(p)}_{n_a}-n_a^p)(1+O(E/K))=O(EK^{p-2})\). For the second
bound write \(\omega^{(p)}=\partial_IH^{(p)}\) with the action-only
Hamiltonian \(H^{(p)}\), a spatial integral of a density with at most \(p\)
derivatives. At the complex points described above all tested
Cartesian inputs are bounded in \(\ell^2\) by fixed constants, and the
integral is bounded by \((CKd^C+Z)^pE\) after taking absolute angular
sums. Cauchy's estimate in the boundedly many tested action slots, on
discs of radius \(c\rho^2\) with \(\rho\) fixed, gives the bound \(EK^p\) for
any fixed number of action derivatives.

For $w\ge2$ and a chosen output copy let $t=\sum l_aV_a$. The measure
$\sum l_a\delta_{V_a}-\sigma\delta_{\sigma t}$ has a nonzero
restriction to the nonzero integers: equality there would force $w=0$
or $1$ since the $V_a$ are distinct. If $b$ is its number of nonzero
support points, the moments of degrees $2,\ldots,b+1$ distinguish it,
by the Vandermonde matrix times the diagonal of the squared sites.
Some such moment has absolute value at least one, and $b\le M_0+1$.
Explicitly, if the restriction is $\sum_{i=1}^bc_i\delta_{x_i}$ with
distinct nonzero integers $x_i$ and nonzero integers $c_i$, then
\[
 \Bigl(\sum_ic_ix_i^{\,p}\Bigr)_{p=2}^{b+1}
 =\bigl(x_i^{\,p-2}\bigr)_{p,i}\,\operatorname{diag}(x_i^2)\,(c_i)_i,
\]
the first factor is an invertible Vandermonde matrix, and the left
side is a nonzero integer vector. Since $n_a=1+KV_a$ and
$\sigma\nu_j=\sigma(Kt+\sum l_a+j\cdot n_F)$, the coefficient of $K^p$
in the free part $\sum l_an_a^p-\sigma(\sigma\nu_j)^p$ of the divisor
below is $\sum l_aV_a^p-\sigma(\sigma t)^p$, which is this moment.
For this $p$, at $|j|_1\le L$ the signed divisor
\[
 j\cdot\omega_F^{(p)}+\sum l_a\omega_{n_a}^{(p)}
                  -\sigma(\sigma\nu_j)^p
\]
has magnitude at least $K^p/2$: its leading coefficient in $K^p$
is that integer and its remaining terms cost $E K^{p-1}$.
Choose $p$ separately in each output copy and divide its own row
before estimating the two coupled rows.

We claim $\sum_j|u^\sigma_{P,j}|\le E K^{-(w-1)_+}$.
For $w\le1$ use Cauchy. For larger $w$ induct on the total degree.
A derivative of a Birkhoff frequency removes a positive number of
balanced pairs from the input monomial, retaining its net label;
its $E K^p$ cost is paid by the divisor, with the old-angle derivative
cost at most $L$ inside the cutoff. For a nonlinear physical product of $g$ factors the number of
spatial derivatives is $p+1-g$, so division by the divisor of
size at least $K^p/2$ gains $K^{-(g-1)}$ up to the factor $E$.
Let $w_i$ be the weights of the factors and let $k$ be the number
of them with $w_i>0$. The triangle inequality for their signed
multiplicities gives $w\le\sum_i w_i$. Consequently the total
power gained satisfies
\[
 (g-1)+\sum_i(w_i-1)_+
 =\sum_i w_i-1+(g-k)\ \ge\ w-1 .
\]
For each $i$, $(w_i-1)_+=w_i-\mathbf{1}_{\{w_i>0\}}$ since $w_i$ is a
nonnegative integer. Summing gives $\sum_i(w_i-1)_+=\sum_iw_i-k$,
and $g\ge k$ gives the inequality.
The terms having the unknown of the same total degree have all their
other factors old. Their operator norm in absolute angular sums
after division is at most $E K^{-2}$, since $g\ge3$. This absorbs
the truncated unknown jointly in the copies. All high-angle
complements of these products were bounded beforehand by Cauchy;
the bounded degree induction closes.

For definiteness we write the identity used in this induction. Let
$X^{(p)}$ be the physical vector field of the $p$th hierarchy
Hamiltonian. Its linear part multiplies $e^{inx}$ by $in^p$, and by
homogeneity its nonlinear part is a sum of products of $g\ge3$ factors
$\partial_x^{k_1}u^{\pm},\ldots,\partial_x^{k_g}u^{\pm}$ with
$k_1+\cdots+k_g=p+1-g$. Differentiating the physical field
$u=\Psi(z)$ along the Birkhoff flow $\dot z=X_{\mathrm B}^{(p)}(z)$ gives
\[
 D\Psi(z)X_{\mathrm B}^{(p)}(z)=X^{(p)}(\Psi(z)).
\]
Take the coefficient of the monomial $P$ at old angle $j$ in copy
$\sigma$. On the left expand the frequencies in the new actions about
zero. The zeroth-order term gives
$i\bigl(\sum_al_a\omega^{(p)}_{n_a}+j\cdot\omega^{(p)}_F\bigr)u^\sigma_{P,j}$.
The higher terms give the coefficients $u^\sigma_{P',j}$ of the monomials
$P'$ from which $P$ arises by adjoining balanced pairs $z_a\bar z_a$,
multiplied by new-action derivatives of the frequencies; for the old
frequencies these derivatives carry the factor $j$, of size at most
$L$ inside the cutoff. On the right the linear part gives
$i\sigma(\sigma\nu_j)^pu^\sigma_{P,j}$. Hence
\[
 \Bigl(j\cdot\omega^{(p)}_F+\sum_al_a\omega^{(p)}_{n_a}-\sigma(\sigma\nu_j)^p\Bigr)
 u^\sigma_{P,j}
 =\sum_{P'}c_{P,P'}\,u^\sigma_{P',j}-i\bigl[N^{(p)}\bigr]^\sigma_{P,j},
\]
where the frequencies are evaluated at zero new actions,
$|c_{P,P'}|\le EK^pL$, the sum runs over boundedly many $P'$ of smaller
degree with the same net label, and $[N^{(p)}]^\sigma_{P,j}$ is the
coefficient of the nonlinear part. By the induction hypothesis the
first sum is at most $E^2LK^{p}K^{-(w-1)_+}$ in absolute angular sum. For the second sum, write $j=j_1+\cdots+j_g$ for the angular labels in a
product. Even when $|j|_1\le L$, some $|j_i|_1$ can exceed $L$. For each
such $i$, use the exponential tail of the differentiated $i$th factor
and the full absolute angular sums of the other factors. Since the
degree and derivative order are fixed, the strip estimate above makes
their total contribution $O(K^{-A})$ for any prescribed $A$ by increasing
$C'$ in $L=C'\log K$. In the remaining products each derivative on
factor $i$ costs at most $CKd^C+Z|j_i|_1\le EK$. Thus a product of
$g$ factors of smaller degrees and weights
$w_i$ costs at most $EK^{p+1-g}\prod_iEK^{-(w_i-1)_+}\le E^{g+1}K^pK^{-(w-1)}$
by the weight inequality above. Division by the divisor, of modulus at
least $K^p/2$, gives the claim for all terms without a factor of full
degree. In the remaining terms one factor has the full degree and the
others are old; they define a linear map of the unknown vector
$(u^\sigma_{P,j})_{j,\sigma}$ of norm at most $EK^{1-g}\le EK^{-2}$ in absolute
angular sum, which is absorbed by a Neumann series.

When $P$ consists of $b$ distinct signed variables at pairwise
distinct sites, a second induction compares its coefficient with the
vacuum coefficient concentrated at old angle zero. The difference in
absolute angular sum is at most $E K^{-b}$. The single-input case
was proved above. No action derivative can remove a balanced pair
in such a monomial. With at least one entirely old factor, the last
weight inequality gains a further power $K^{-1}$ because every
nonold child has weight equal to its degree. With no old factor,
replace the children one at a time by their vacuum values and use
the induction bound. The remaining all-vacuum product has only
$j=0$; the change between the true and vacuum divisors there costs
$E K^{p-1}$ times a vacuum coefficient $O(E K^{-(b-1)})$ before
division by $K^p$. The same-degree term with old factors may also be
estimated using the already proved first induction bound. This proves
the claimed comparison including subwords whose vacuum coefficient
vanishes by charge. Thus \eqref{rate:eq:axial-path} and its homogeneity
give $(Kl)^{-(q-1)}C(v,e)+O(E K^{-q})$ for the desired character
when all new actions vanish. Here the vacuum coefficient at the sites
$n_a=1+Klv_a$ is a sum of products of $q-1$ reciprocal differences
$n_a-n_b=Kl(v_a-v_b)$, so it equals $(Kl)^{-(q-1)}C(v,e)$, and the
comparison just proved is applied with $b=q$.

To restore the nonweak actions, the Taylor terms through $L_0$ added
pairs still have weight $q$. Summing their action products at each
degree $j$ costs at most $S^j$ with the ordinary Taylor coefficient
convention, since the multinomial expansion of $S^j$ dominates this
sum. For the remaining tail multiply all new nonweak actions by
$\xi$. At each signed marked site take its selected side to be
$\rho e^{\pm i\phi_a}$ and the other to be
$\xi I_a\rho^{-1} e^{\mp i\phi_a}$; at unmarked new sites use symmetric
square roots of $\xi I_a$. Projecting on the selected character and
dividing by $\rho^q$ gives precisely the desired action series and
removes all square-root choices. On $|\xi|\le c/S$ the total squared
size in the fixed $q$ tested slots and all other slots stays in the
common small $\ell^2$ ball, and normalized Haar integration bounds
the physical coefficient independently of the number of unmarked
sites. One-variable Cauchy now bounds the omitted tail by
$(CS)^{L_0+1}$. This proves the estimate and both lower bounds.
\end{proof}

For a row with $n_a=1+Klv_a$, the spatial and kinetic sums both equal
one. Indeed, by the three moment conditions of
Lemma~\ref{rate:lem:axial-template},
\[
 \sum_ae_an_a=\sum_ae_a+Kl\sum_ae_av_a=1,\qquad
 \sum_ae_an_a^2=\sum_ae_a+2Kl\sum_ae_av_a+(Kl)^2\sum_ae_av_a^2=1 .
\]
In applying Lemma~\ref{rate:lem:axial-template} to a passive deficit
one needs the true empty-frequency expansion to start at order $n^{-2}$
with the same inverse-power coefficients for positive and negative $n$.
Under that expansion through $h$, expanding $(1+Klv_a)^{-j}$ for
$j\ge2$ and using the reciprocal moments gives a deficit
$O(C_{\mathrm{old}}(Kl)^{-h-1})$. Explicitly, for $2\le j\le h$,
\[
 \sum_ae_an_a^{-j}
 =\sum_{i=0}^{h-j}\binom{-j}{i}(Kl)^{-j-i}\sum_ae_av_a^{-j-i}
   +O\bigl((Kl)^{-h-1}\bigr)=O\bigl((Kl)^{-h-1}\bigr),
\]
because $2\le j+i\le h$ in every term of the finite sum, and the
remainder of the expansion of $\chi_{n_a}$ is of the same order. An expansion merely allowing a nonzero
$n^{-1}$ term would not suffice, since $\sum e_a/v_a\ne0$.

\subsection{Scale, geometry and calibration}
The quantifiers in this section are for each fixed finite $s>2$ separately.
The Fourier lattice is measured on $\mathbb T^2$ of period $2\pi$ and all integrals use normalized measure. In the mass gauge the genuine axial frequencies are $n^2+\chi_n(I)$, where $I_n=|z_n|^2$ in a choice of Birkhoff normalization with free differential Fourier. The mass is $\sum I_n$. The map $\chi$ is analytic on a complex $\ell^2$ action ball and $D\chi(0)=-\mathrm{Id}$. This follows by rescaling the cubic coefficient in \cite[Theorem~18.1 and Corollary~18.2, $p=4$]{Molnar}, together with Appendix~\ref{app:molnar-four}. The Hessian of the action Hamiltonian is symmetric on positive finite-dimensional faces. We use translation and phase covariance, simultaneous Birkhoff coordinates for the one-dimensional NLS hierarchy, and the one-derivative gain of the nonlinear part of the one-dimensional coordinate map and its inverse on the integer Sobolev spaces of order at least one. These statements concern the axial equation.

Use $B=B_*Q^j$, $B_o=B/Q$, $d=\lfloor\alpha\sqrt{B/\log B}\rfloor$, and a vertical spacing $R$ with $\log R=B+O(1)$. It may be chosen divisible by all earlier spacings and prime to five since their product is at most $\exp(C B_o)$. A rational family of normalized heads has a common denominator $P$ and physical indices $(1,0)+R v$, where $v$ is the integer macro vector, and the top is $(1+N,0)$ with $N=RP$. Thus a normalized direction $x$ corresponds to integer macro $v=P x$. With $H=e^{b d}$, choose
\[
 r=HN^{-s},\qquad \lambda=e^{-\epsilon B},\qquad u=r\lambda,
 \qquad \tau=ut .
\]
Here $b>0$ is small relative to $s$ and the constant contraction of the chosen split. At fixed $s$ choose $a<\min(s-2,1/4)/8$. All fixed derivatives and high spatial truncation orders are specified before $\epsilon$. We distinguish factors $\lambda^{-C_s}$ from
\[
 E=\exp\{C_{s,\epsilon}(B_o+d^2\log(2+B))+o(B)\}.
\]
Constants in these two classes can increase through a fixed number of operations. Only the second class can depend on the cutoff in the number of weak angular letters. First choose $\epsilon$ small compared with $a$ and all the fixed powers in the first class. After the intermediate historical and preparation constants have been chosen, take $Q$ large and $\alpha$ small enough that each needed product of factors in the second class is $\le N^\gamma$, with $\gamma$ also much smaller than $\epsilon$. In particular $\alpha^2Q$ can be made as small as desired. No lower bound on the final constant in the growth exponent is prescribed.

\smallskip
\noindent\emph{A ledger of the stage prices.} The factors denoted
\(E\), \(E_o\), \(E_c\), \(\mathcal E\) in this section and in
Sections~\ref{rate:sec:transverse-principal}--\ref{rate:sec:ward-landing},
and the relation \(X\preceq Y\) of Section~\ref{rate:sec:exterior}, are
all bounded by means of the following elementary facts. In each item
the bound follows from the coarse caps registered at the end of the
previous stage (Section~\ref{rate:sec:parameters}) and from the
definitions of \(B,B_o,d,P\).
\begin{itemize}
\item[(Pr1)] \emph{Depth.} Since \(d\le\alpha\sqrt{B/\log B}\) and
\(\log(2+B)\le2\log B\) for large \(B\),
\(d^2\log(2+B)\le2\alpha^2B\). With \(D_g=C_*d\log(d+2)\) and
\(d\le B\) one has \(dD_g\le C_*d^2\log(2+B)\),
\(\log(1/\kappa)=10\log(2+d+D_g)\le C\log(2+B)\), and the binomial
weights \(c_{nk}\le2^{d/2}\). Hence \(e^{Cd}\), \(e^{CdD_g}\),
\(\kappa^{-C}\) and \(2^{Cd}\) are at most
\(\exp\{C(d^2\log(2+B)+\log(2+B))\}\).
\item[(Pr2)] \emph{The previous support.} The registered caps give
\(\log(2+\max|n|)\le2\log N_{j-1}+O(1)\) and
\(-\log I_{\min}\le(2s+4)\log N_{j-1}+O(1)\) on the support inherited
from stage \(j-1\), and \(\log N_{j-1}\le2B_{j-1}=2B_o\). Thus the old
sites have modulus at most \(e^{4B_o+O(1)}\) and the old positive
actions are at least \(e^{-(4s+8)B_o-O(1)}\). The initial torus has two sites, and stage \(i\) adds at most
\(C_s(1+d_i)^2\) sites. To sum these increments, put
\(v(x)=x/\log x\). For \(x\ge e^2\) one has
\(d\log v(x)/d\log x=1-1/\log x\ge1/2\), hence
\(v(B_i)\le Q^{-(j-i)/2}v(B_j)\) for \(i\le j\).
Take \(B_*\) large enough that \(\alpha\sqrt{v(B_*)}\ge2\).
Then \((1+d_i)^2\le4\alpha^2v(B_i)\), and summing the geometric
series for \(Q\ge2\) gives
\[
  2+C_s\sum_{i\le j}(1+d_i)^2
  \le C'_s(1+\alpha^2v(B_j))\le C''_s(1+d_j^2).
\]
The same estimate through stage \(j-1\) gives
\(f_{j-1}\le C_s(1+\alpha^2B_o/\log B_o)\le C_s B_o\);
at the first stage \(f_0=2\).
\item[(Pr3)] \emph{Sites added before the hot axes.} The spectator has
\(|\ell|=e^{C_1B_o+O(1)}\) and action at least
\(a_\ell=e^{-C_2B_o}\). The moderate and weak sites are
\(1+Kl_ev_i\) with \(K=e^{kB_o+O_s(\log B)}\), \(l_e\le Cd^2\) and
template entries \(|v_i|\le C_s\); there are at most \(C_sd^2\) of
them, of modulus at most \(e^{kB_o+C_s\log B}\). The moderate actions
are at least \(a_m\ge K^{-h-1}\ge e^{-(h+1)kB_o-C\log B}\), and the
weak actions are \(I_w=\lambda^{2/5}e^{O(\Gamma)}\) once activated.
\item[(Pr4)] \emph{Lattice.} \(\log P\le C_{k_0,m_*}d^2\log(2+B)\) by
Lemma~\ref{rate:l:denom}, hence \(H_{\rm geom}=C(1+d)P\) has
logarithm \(O(d^2\log(2+B))\), and
\(\log N=\log R+\log P=B+O(d^2\log(2+B))\).
\item[(Pr5)] \emph{Angular sums.} For \(\delta>0\) and \(k\ge0\),
\(\sum_{l\in\Z^f}e^{-\delta|l|_1}=\coth(\delta/2)^f\) and
\(\sum_{l\in\Z^f}|l|_1^ke^{-\delta|l|_1}
\le(2k/\delta)^k\coth(\delta/4)^f\), since
\(x^ke^{-\delta x/2}\le(2k/\delta)^k\). A spatial derivative of the
angular harmonic with label \(l\) of an old field multiplies it by
\(l\cdot n_F\), of modulus at most \(Z|l|_1\), where
\(Z=1+\max_{m\in F}|m|\).
\end{itemize}
Consequently, on the current support \(F\) (old sites, spectator,
moderate and weak sites), \(Z\le e^{CB_o+C\log B}\),
\(f\le CB_o+Cd^2\), and every nonweak positive action is at least
\(e^{-CB_o-C\log B}\), with \(C\) depending on \(s,k,C_1,C_2\) and the
fixed orders. Any product of a fixed number of the factors \(1+Z\),
\(e^f\), \(C^f\), inverse nonweak actions, \(e^{Cd}\),
\(e^{CdD_g}\), \(\kappa^{-1}\), \(P\), \(H_{\rm geom}\) and
\(2+B\) is therefore at most \(e^{C\Gamma}\), where
\(\Gamma=B_o+d^2\log(2+B)+\log(2+B)\); this is a factor of the class
\(E\), the term \(C\log(2+B)\) being part of the \(o(B)\) in its
exponent. Since \(B_o=B/Q\), (Pr1) and (Pr4) give
\(\log E\le C_{s,\epsilon}(Q^{-1}+2\alpha^2)B+o(B)\le\gamma\log N\)
once \(Q\) is large, \(\alpha\) small and \(B_*\) large. By (Pr3),
fixed inverse powers of the current positive weak actions cost at most
\(e^{C\Gamma}\lambda^{-C_s}\). The factor \(e^{C\Gamma}\) belongs to
the class \(E\), and \(\lambda^{-1}=e^{\epsilon B}\) is paid for by
the choice of \(\epsilon\). The change of time units contributes
\(u^{-1}=r^{-1}\lambda^{-1}\), where \(r^{-1}=H^{-1}N^s\).
The powers of \(r^{-1}\) are retained in the explicit terms
\(N^{-J_0}/u^2\), \(N^{A-A'}/u\) and in the conversion of norms in
the stage budget.

\subsubsection{A true axial frequency with an affordable interval of values}
On a positive old support $F$ let $M=\sum I_n$, $\mathcal H$ be the ungauged one-dimensional action Hamiltonian and $\omega_n=n^2+2M+\chi_n$. Set
\[
 \mathcal Q=2\left(\sum I_n\omega_n-\mathcal H\right),
 \qquad A=\mathcal Q-M^2,\qquad D=\mathcal Q/2.
\]
On a small initial two-site torus, $A>0$ and $2M^2-\mathcal Q>0$ with strict fixed margins. Indeed the first quartic normal form gives $\mathcal Q=M^2+2I_{n_1}I_{n_2}+O(M^3)$ when $n_1,n_2$ are the two comparably occupied sites. In detail, $\omega_n=n^2+2M-I_n+O(|I|^2)$ by $D\chi(0)=-\mathrm{Id}$, and integrating $d\mathcal H=\sum\omega_n\,dI_n$ from zero gives
\[
 \mathcal H=\sum n^2I_n+M^2-\tfrac12\sum I_n^2+O(M^3),\qquad
 \mathcal Q=2\Bigl(\sum I_n\omega_n-\mathcal H\Bigr)=2M^2-\sum I_n^2+O(M^3).
\]
On two sites $2M^2-I_{n_1}^2-I_{n_2}^2=M^2+2I_{n_1}I_{n_2}$, so that $A=2I_{n_1}I_{n_2}+O(M^3)$ and $2M^2-\mathcal Q=I_{n_1}^2+I_{n_2}^2+O(M^3)$. The identities $d\mathcal H=\sum\omega_n\,dI_n$ and
\[
 d\mathcal Q=4M\,dM+2\sum I_n\,d\chi_n
\]
show that sufficiently small cumulative changes in $\ell^1$ preserve these margins, independently of the new frequencies. The second identity follows from the first: $d\mathcal Q=2\sum I_n\,d\omega_n$ and $d\omega_n=2\,dM+d\chi_n$. Suppose $\max_F|n|\le e^{C B_o}$, every positive old $I_n\ge e^{-C B_o}$ and the old size in $h^1$ remains in the small ball. On an old action ball of radius $h_0\ge e^{-C'B_o}$ inside positivity and this $h^1$ chart, the high empty axial expansion gives $|\chi_l(I,0)|\le E_o/l^2$ with $E_o\le e^{C'' B_o}$. The order two gain here matters. In the doubled high quasimode the order-zero coefficient is the mean of $B_0=2(|V|^2-M)$, which vanishes. The order-one coefficient is a constant times the mean of $B_0\partial_x^{-1}B_0$, also zero. Projection of the approximate high column by the actual Birkhoff differential onto its positive free mode leaves a scalar $1+O(E_o/|l|)$ with a single-valued argument on the angle torus. Its directional derivative integrates to zero against Haar measure, giving the same two vanishing coefficients for the true empty frequency. The uniform cost of any fixed higher order is polynomial in the old maximum, inverse positive floors and $e^{|F|}$ by old-angle analyticity.

\smallskip
\noindent\emph{Derivation of the bound \(|\chi_l|\le E_o/l^2\).} Choose the initial torus and the summable \(h^1\) increments so that
each old torus lies in the concentric ball of half the chosen chart
radius. Take \(h_0=\tfrac12\min_{n\in F}I_n\), so
\(h_0\ge e^{-(4s+8)B_o-O(1)}\) by (Pr2). For real actions \(J\) on \(F\)
with \(\|J-I\|_2\le h_0\), we have \(I_n/2\le J_n\le 3I_n/2\).
Thus the Cartesian \(h^1\) norm increases by at most a factor
\(\sqrt{3/2}\), and the entire positive family lies in the chart
with the same \(Z\) and \(f\).
Lemma~\ref{rate:lem:trans-axial-expansion} with expansion order two
(its parameter, also called \(h_0\) there, equal to \(2\); this is not
the radius) gives \(\chi_l=c_2(I)l^{-2}+O(E_o'|l|^{-3})\) with
\(c_2=\mathcal Q/2-M^2\), and \(|c_2|\le C\) on the small ball because
\(\mathcal Q=2M^2-\sum I_n^2+O(M^3)\). No action derivative is taken,
so \(E_o'\) involves no positive floor: it bounds a fixed number of
coefficients of the spatial recursion and the one-smoothing projection.
The spatial recursion uses products, spatial and nominal time derivatives,
spatial means and the zero-mean primitive $D_x^{-1}$. At each fixed order
there are only finitely many such operations. The mean and $D_x^{-1}$
are bounded by one in the absolute Fourier sums, the latter having
multiplier $1/k$ for $k\ne0$ and zero at $k=0$. On the angular
harmonic with label \(l'\) a spatial derivative costs \(Z|l'|_1\) and a
nominal time derivative costs
\(|l'\cdot\Omega|\le C(1+Z^2)|l'|_1\), since the old angular speeds are
\(n^2+O(1)\). The absolute angular sums then cost
\((C_k)^f\) by (Pr5). Hence \(E_o'\le C(1+Z)^{C}C^f\le e^{CB_o}\) by (Pr2),
and \(|\chi_l|\le(C+E_o'|l|^{-1})l^{-2}\le E_o/l^2\) with
\(E_o\le e^{C''B_o}\).

For any real unit old direction apply the central difference formula at $h=\sqrt{E_o}/|l|\le h_0$, using the uniform bound for $D^2\chi$ on the common ball. This gives $\|P_F\nabla\chi_l\|_2\le C\sqrt{E_o}/|l|$. Indeed, for a unit vector $e$ supported on $F$,
\[
 |\partial_e\chi_l(I)|
 \le\frac{|\chi_l(I+he)-\chi_l(I-he)|}{2h}+\frac h2\sup|\partial_e^2\chi_l|
 \le\frac{E_o}{l^2h}+Ch=(1+C)\frac{\sqrt{E_o}}{|l|}.
\] Symmetry gives the corresponding bound on $P_F D\chi\,e_l$. Choose the empty $l\equiv1\pmod5$ with $|l|=e^{C_1 B_o+O(1)}$, increasing $C_1$ to pay $\sqrt{|F|}$. Lock $\chi_F$ by its own old actions while filling $I_l>0$; their Jacobian inverse has norm at most two on a sufficiently small common ball. The derivative of $M$ with respect to $I_l$ is between $3/4$ and $5/4$ on a short positive interval, since this holds at the empty endpoint and $D^2\chi$ and the compressed inverse are bounded. In detail, differentiating the lock $\chi_F(I_F(I_l),I_l)=\mathrm{const}$ gives $I_F'=-(D_{FF}\chi_F)^{-1}\partial_{I_l}\chi_F$, where $D_{FF}\chi_F$ is the compression of $D\chi$ to $F$ and $\|(D_{FF}\chi_F)^{-1}\|\le2$. The matrix $D\chi$ is symmetric, being the Hessian of the action Hamiltonian minus the constant matrix with all entries $2$, so $\partial_{I_l}\chi_F=P_F\nabla\chi_l$; its $\ell^2$ norm is at most $C\sqrt{E_o}/|l|$ at $I_l=0$ by the central difference bound, and it changes by at most $CI_l$ on the interval by the bound on $D^2\chi$. As $M=\sum_{n\in F}I_n+I_l$, Cauchy--Schwarz gives $|M'-1|\le|F|^{1/2}\|I_F'\|_2\le2C|F|^{1/2}(\sqrt{E_o}/|l|+I_l)$. This is at most $1/4$ if $|l|\ge16C(|F|E_o)^{1/2}$, which holds for $|l|=e^{C_1B_o+O(1)}$ with $C_1$ large since $|F|E_o\le e^{CB_o}$ by (Pr2), and if $I_l\le(16C|F|^{1/2})^{-1}$. Shrinking the interval pays any prescribed $e^{-C_h B_o}$ old norm tolerance, including the new cost $\bk l^s\sqrt{I_l}$, and leaves an interval of length $e^{-C_2B_o}$. On it
\[
 (A/D)'=4M(2M^2-\mathcal Q)M'/\mathcal Q^2+O(I_l).
\]
To verify the formula, write $A/D=2-2M^2/\mathcal Q$, with primes
denoting derivatives in $I_l$ along the lock. The old rates are
locked, so $d\mathcal Q=4M\,dM+2\sum I_n\,d\chi_n$ gives
$\mathcal Q'=4MM'+2I_l\chi_l'$, and
\[
 (A/D)'=-\frac{4MM'}{\mathcal Q}+\frac{2M^2\mathcal Q'}{\mathcal Q^2}
 =\frac{4MM'(2M^2-\mathcal Q)}{\mathcal Q^2}
   +\frac{4M^2I_l\chi_l'}{\mathcal Q^2}.
\]
The main term has a fixed sign and size comparable to $M'$ by the
two strict margins and $3/4\le M'\le5/4$.
Take the spectator action in $[a_\ell,2a_\ell]$ with
$a_\ell=\exp(-C_2 B_o)$, increasing $C_2$ to pay the prescribed old
tolerance. The image under $A/D$ has length comparable to $a_\ell$.
There are at most $e^{C_4B_o}$ distinct nonzero rational component
pairs $(x,y)$ of height at most $e^{C_3B_o}$. Arrange first
$\log P\le C_{s,\epsilon}d^2\log B\le C B_o$. The quotients $y/x$
with $x\ne0$ divide the image interval into at most $e^{C_4 B_o}+1$
pieces. Choose the middle half of a longest piece. Its preimage
has length at least $e^{-C_5B_o}$ and
$|A x-Dy|\ge e^{-C_5B_o}$ there, after increasing $C_5$: $D$ has
a fixed positive lower bound and nonzero $|x|\ge e^{-C_3B_o}$.
For $x=0$ use the corresponding bound on $|y|$.

\smallskip
\noindent\emph{The selection of the spectator interval in detail.}
(i) On the small ball the two strict margins give fixed numbers
\(0<M_-\le M\le M_+\) and \(\varpi>0\) with \(A\ge\varpi\),
\(2M^2-\mathcal Q\ge\varpi\), \(D=\mathcal Q/2\ge\varpi/2\) and
\(\mathcal Q\le2M^2\). Hence the main term
\(4M(2M^2-\mathcal Q)M'/\mathcal Q^2\) lies between
\(c_{\rm sp}=3M_-\varpi/(4M_+^4)\) and \(C_{\rm sp}=10M_+^3/\varpi^2\) (use \(M'\in[3/4,5/4]\),
\(\mathcal Q\le2M_+^2\) and \(\mathcal Q\ge2D\ge\varpi\)), and the remainder
\(4M^2I_l\chi_l'/\mathcal Q^2\) is \(O(I_l)\), because \(\chi_l'\) is
bounded. For \(I_l\le2a_\ell\) small, \((A/D)'\in[c_{\rm sp}/2,2C_{\rm sp}]\), so
\(A/D\) maps \([a_\ell,2a_\ell]\) increasingly onto an interval of
length between \(c_{\rm sp}a_\ell/2\) and \(2C_{\rm sp}a_\ell\).
(ii) The component pairs to be separated have height at most
\((2dP)^{C_{m_*}}\le e^{C_3B_o}\), by Lemma~\ref{rate:l:denom}, (Pr4) and
\(d^2\log(2+B)\le2d^2\log B\le2c_*B_o\) (the standing demand
\(d^2\log B\le c_*B_o\) of Section~\ref{rate:sec:parameters}). A pair of rationals of height at most
\(\mathfrak h\) is one of at most \((2\mathfrak h+1)^4\) pairs, so there
are at most \(e^{C_4B_o}\) of them. A nonzero rational of height at most
\(\mathfrak h\) has modulus at least \(1/\mathfrak h\).
(iii) The at most \(e^{C_4B_o}\) forbidden values \(y/x\), \(x\ne0\),
divide the image interval into at most \(e^{C_4B_o}+1\) pieces. A
longest piece has length
\(\ell_*\ge c_{\rm sp}a_\ell/(2(e^{C_4B_o}+1))\). On its middle half every value
\(v\) of \(A/D\) satisfies \(|v-y/x|\ge\ell_*/4\) for all forbidden
values, so
\[
 |Ax-Dy|=|x|\,D\,|A/D-y/x|\ge e^{-C_3B_o}\tfrac{\varpi}2\tfrac{\ell_*}4,
\]
and for \(x=0\ne y\) one has \(|Ax-Dy|=D|y|\ge\tfrac{\varpi}2e^{-C_3B_o}\).
The preimage of the middle half has length at least
\(\ell_*/(4C_{\rm sp})\). With \(a_\ell=e^{-C_2B_o}\) all these quantities are
at least \(e^{-C_5B_o}\) for a suitable \(C_5\) depending on
\(C_2,C_3,C_4\) and the margins.

There is an additional choice on this very interval. Assume the
previously accepted direct rates $\chi_F$ satisfy
\begin{equation}\label{rate:eq:registry-prefix}
 |n+l\cdot\chi_F|\ge\kappa_F(2+|n|+|l|_\infty)^{-2|F|-6},
 \quad(n,l)\in\mathbb Z^{1+|F|}\setminus\{0\},\quad
 \kappa_F\ge e^{-C_{\mathrm{reg}}B_o}.
\end{equation}
At the base choose the two positive actions in an interior cube on which
the quartic margins hold; the local direct frequency chart and the
slicing argument below give this bound with some positive constant.
Along the old lock, $\frac{d}{dI_\ell}\chi_\ell$ is the scalar Schur
complement of $D_F\chi_F$ in $D_{F\cup\{\ell\}}\chi_{F\cup\{\ell\}}$.
It differs from $-1$ by at most $\eta+\eta^2/(1-\eta)$ if
$\|D\chi+\mathrm{Id}\|\le\eta$. Indeed the Schur complement is
$D_{\ell\ell}-D_{\ell F}D_{FF}^{-1}D_{F\ell}$ with
$|D_{\ell\ell}+1|\le\eta$, $\|D_{\ell F}\|,\|D_{F\ell}\|\le\eta$ and
$\|D_{FF}^{-1}\|\le(1-\eta)^{-1}$, the last by a Neumann series for
$D_{FF}=-(\mathrm{Id}-(D_{FF}+\mathrm{Id}))$. Hence on the chosen subinterval
the direct spectator rate covers an interval of width
$\delta_\ell\ge e^{-C_6B_o}$. Choose it so that
\[
 |n+l\cdot\chi_F+b\chi_\ell|
 \ge c\,\delta_\ell\,2^{-C(|F|+1)}
 (2+|n|+|l|_\infty+|b|)^{-2(|F|+1)-6}
 \quad(b\ne0).
\]
Indeed each bad relative length is bounded by twice its threshold
divided by $|b|\delta_\ell$, and summing the integer tuples in dyadic
height shells costs less than half the interval for an absolute
choice of $c,C$. Explicitly, with $f_+=|F|+1$ and
$h=2+|n|+|l|_\infty+|b|$, there are at most $(2^{j+2})^{f_++1}$ tuples
$(n,l,b)$ with $2^j\le h<2^{j+1}$, and the sum of the relative lengths
is at most
\[
 \sum_{j\ge1}(2^{j+2})^{f_++1}\cdot2c\,2^{-Cf_+}2^{-j(2f_++6)}
 \le2c\,2^{(4-C)f_+}\sum_{j\ge1}2^{-5j}<\tfrac12
\]
for $C=4$ and $c=1$. Indeed
$(2^{j+2})^{f_++1}2^{-j(2f_++6)}=2^{2f_++2}2^{-j(f_++5)}
\le2^{f_++2}2^{-5j}\le2^{4f_+}2^{-5j}$ for $j,f_+\ge1$, and
$2\sum_{j\ge1}2^{-5j}=2/31$. Restricting beforehand to a middle subinterval
leaves the same spatial margins. For $b=0$ inherit
\eqref{rate:eq:registry-prefix}. Lock this spectator rate and its chosen
angle origin exactly through the rest of this stage and every later
stage. We write $F_+=F\cup\{\ell\}$. Its new prefactor and every
price of this preliminary spectator selection are $e^{-O(B_o)}$
and $e^{O(B_o)}$, respectively, with constants allowed to depend on
the requested old precision. Only now choose the moderator dilation
$K=e^{kB_o+O_s(\log B)}$. Thus a fixed spectator rate resonant with
the old rates is never silently included in the frozen relation module.

\subsubsection{Simultaneous positivity and affine frequencies}
Let $H_e$ here denote rows of bounded length with disjoint coordinates, each of charge one, horizontal momentum one and kinetic sum one, containing a positive and a negative unit moderator pivot and five private weak letters. Write $\mathsf H$ for the matrix with rows $H_e$. Keep all of $F_+$ locked. If $J$ denotes the prospective new
actions, the old action derivative is
$-(D_{F_+}\chi_{F_+})^{-1}D_J\chi_{F_+}$ of norm at most
$\eta/(1-\eta)$. The resulting derivative of the new direct rates is
their Schur complement, within $\eta+\eta^2/(1-\eta)$ of minus
identity in Euclidean operator norm. On each interior ball its local
inverse takes direct-rate increments to action increments with norm
at most two after shrinking the common vacuum neighborhood. This
last assertion holds on a ball of radii comparable to the available
Euclidean distance to the action boundary, by solving the near-identity
map by contraction. No product of one-coordinate inverse norms occurs.

Let the rows use $q$ fixed signs, and write $m$ for their number.
Initially every weak site is empty and all moderators get the same
positive floor $a_m$. For each row choose two moderator pivots with
opposite signs. Put $v_e^+$ on its positive and $v_e^-$ (the magnitude
of its negative part) on its negative pivot, obtaining $S(v)$.
Then $\mathsf H S(v)=v$, $\operatorname{Lip}_{\ell^2} S\le1$, and
$\|\mathsf H\|\le\sqrt q$. The map
$v\mapsto v+\mathsf H\chi_{\mathrm{new}}(I+S(v))$ on the old locked chart
contracts if $\sqrt q(\eta+\eta^2/(1-\eta))<1/2$. It can be solved
on the radius $2|\mathsf H\chi_{\mathrm{new}}(I)|_2$ before approaching any
previous action floor. Its nonnegative increments preserve all new
moderator floors. If each unprepared row deficit is $O(a_m)$,
adding the initial moderator floors changes its vector of deficits
in $\ell^2$ by at most $C_q\sqrt m\,a_m$. Thus all moderators have
total action $O_q(m a_m)$, and the locked old action variation in
$\ell^2$ is $O_q(\sqrt m\,a_m)$. A smooth chart about this interior
point uses one unit moderator pivot per row, with either sign of
increment. Its row Jacobian, after the signs are accounted for,
is within $\sqrt q(\eta+\eta^2/(1-\eta))$ of the identity.
In detail, write \(v=v^+-v^-\) with \(v^\pm\ge0\) componentwise, and let
\(S(v)\) place \(v_e^+\) on the positive and \(v_e^-\) on the negative pivot of
row \(e\). Row \(e\) has coefficients \(+1\) and \(-1\) at these pivots and the
supports are disjoint, so \(\mathsf HS(v)=v\). For real numbers
\((a^+-b^+)^2+(a^--b^-)^2\le(a-b)^2\), so \(S\) is \(1\)-Lipschitz in \(\ell^2\),
and \(\|\mathsf Hx\|_2^2=\sum_e(H_e\cdot x)^2\le q\|x\|_2^2\) by Cauchy--Schwarz on
each support. Put \(\Phi(v)=v+\mathsf H\chi_{\rm new}(I+S(v))\); its fixed points
are exactly the solutions of \(\mathsf H\chi_{\rm new}(I+S(v))=0\). Using
\(\mathsf H(S(v)-S(w))=v-w\),
\[
 \Phi(v)-\Phi(w)=\mathsf H\int_0^1\bigl(D\chi_{\rm new}+\mathrm{Id}\bigr)
   \bigl(I+S(w)+\theta(S(v)-S(w))\bigr)\,d\theta\;\bigl(S(v)-S(w)\bigr),
\]
where \(D\chi_{\rm new}\) is the locked derivative (the Schur complement),
within \(\eta+\eta^2/(1-\eta)\) of \(-\mathrm{Id}\). Hence
\(\operatorname{Lip}\Phi\le\sqrt q(\eta+\eta^2/(1-\eta))<1/2\), and the fixed point
satisfies \(\|v\|_2\le2\|\Phi(0)\|_2=2\|\mathsf H\chi_{\rm new}(I)\|_2\). To bound
\(\Phi(0)\), note that before the floors are added each row deficit is
\(O(a_m)\), so these deficits have \(\ell^2\) norm \(O(\sqrt m\,a_m)\). Adding the
floor \(a_m\) on the at most \(qm\) moderator sites changes \(\chi_{\rm new}\) by
\(-a_m\mathbf 1_{\rm mod}\) up to an \(\ell^2\) error
\((\eta+\eta^2/(1-\eta))\sqrt{qm}\,a_m\), and \(\|\mathsf H\mathbf 1_{\rm mod}\|_2\le q\sqrt m\).
Thus \(\|\Phi(0)\|_2\le C_q\sqrt m\,a_m\); the added moderator action is
\(\|S(v)\|_1\le\sqrt{2m}\,\|v\|_2\le C_qma_m\); the total moderator action is
\(O_q(ma_m)\); and the locked old actions move in \(\ell^2\) by at most
\(\eta(1-\eta)^{-1}\) times the \(\ell^2\) norm \(O_q(\sqrt m\,a_m)\) of the new
action increments.

For existence of such rows fix an integer $h>4s+5$. The reciprocal-moment construction supplies a finite signed template $v_i\in\Z\setminus\{0\}$, $e_i=\pm1$, all sites distinct, with
\[
 \sum e_i=1,\quad\sum e_i v_i=\sum e_i v_i^2=0,\qquad
 \sum e_i/v_i\ne0,\qquad\sum e_i v_i^{-m}=0\ (2\le m\le h).
\]
The actual minimal signed inverse-Birkhoff coefficient at vacuum is nonzero. Its order is $q-1$ in inverse site differences, where $q$ is the fixed number of letters in this template. Use distinct disjoint rows $1+K l_e v_i$ with $l_e$ positive integers bounded by a fixed multiple of the number of edges, chosen greedily; take $5\mid K$. The uniform common empty-frequency expansion through $h$, beginning at inverse order two as above, leaves passive deficits $O(e^{CB_o}K^{-h-1})$. Choose $a_m=C_o K^{-h-1}$ with $1\le C_o\le e^{C B_o}$ bounding these deficits on the $F_+$ chart. Increasing $k$ makes the total moderate action cost, including the polynomial factor $m$, at most $K^{-(h+1)/2}$ and preserves all $F_+$ floors. Their full $h^s$ cost is $O(K^{s-(h+1)/4}\operatorname{poly}(d))$, arbitrarily small on the scale of the previous tolerance. Both signs of unit moderator pivots lie among the added letters beyond the initial five.

\smallskip
\noindent\emph{The moderate costs in detail.} (i) \emph{Passive
deficits.} For a row with sites \(n_a=1+Kl_ev_a\) and signs \(e_a\),
Lemma~\ref{rate:lem:trans-axial-expansion} gives
\(\chi_n=\sum_{j=2}^hc_jn^{-j}+O(E_o|n|^{-h-1})\) with coefficients common
to both signs of \(n\) and \(|c_j|\le E_o\le e^{CB_o}\). The reciprocal
moments of the template vanish for \(2\le j\le h\), so, as displayed at
the end of Section~\ref{rate:sec:axial-mark},
\(\sum_ae_an_a^{-j}=O((Kl_e)^{-h-1})\). Since the row has a fixed number
\(q\) of letters, its passive deficit \(\sum_ae_a\chi_{n_a}\) is at most
\(e^{CB_o}K^{-h-1}\).
(ii) \emph{Total moderate action.} There is one row per designated
edge, so \(m\le Cd^2\). By the contraction of the preceding subsection,
with floors \(a_m=C_oK^{-h-1}\), the total moderate action is
\(O_q(ma_m)\le C_qd^2e^{CB_o}K^{-h-1}\). As
\(\log K=kB_o+O_s(\log B)\), one has
\(K^{(h+1)/2}\ge e^{(h+1)kB_o/2-C_s'\log B}\), and \(d^2\le B\); so this
is at most \(K^{-(h+1)/2}\) as soon as \((h+1)k/2>C+1\) and \(B\) is
large, because then
\(C_qd^2e^{CB_o}\le e^{(C+1)B_o}\le e^{(h+1)kB_o/2-C_s'\log B}\) for
large \(B\) (recall \(\log B=o(B_o)\) at fixed \(Q\)). The locked old
actions move by \(O_q(\sqrt m\,a_m)\le K^{-(h+1)/2}\) in \(\ell^2\), and
\(K^{-(h+1)/2}\le e^{-(h+1)kB_o/2+C_s'\log B}\). This is less than half
of every old floor, which is at least \(e^{-(4s+8)B_o-O(1)}\) by (Pr2),
and of the spectator floor \(a_\ell=e^{-C_2B_o}\), once
\((h+1)k/2>4s+9+C_2\) and \(B\) is large.
(iii) \emph{The \(h^s\) cost.} Each moderate site has
\(\bk n\le CKd^2\) by (Pr3) and action at most \(K^{-(h+1)/2}\), and
there are at most \(C_sd^2\) such sites. Hence
\[
 \sum_{\rm moderate}\bk n^s\sqrt{I_n}\le C_sd^2(CKd^2)^sK^{-(h+1)/4}
 =O\bigl(K^{s-(h+1)/4}\operatorname{poly}(d)\bigr),
\]
and \(s-(h+1)/4<-3/2\) because \(h>4s+5\). With \(K=e^{kB_o+O(\log B)}\)
this is below any prescribed \(e^{-D'B_o}\) after increasing \(k\).

Apply Proposition~\ref{rate:prop:axial-mark-true} with the positive old data just selected, including the spectator, and the total moderator action bounded above. Lemma~\ref{rate:lem:axial-template} supplies the nonzero vacuum coefficient, and Lemma~\ref{rate:lem:trans-axial-expansion} supplies the passive deficit. By increasing the fixed moderator placement constant $k$, the true physical coefficient, divided only by its five weak radii, has modulus between $K^{-C_s}$ and $K^{C_s}$. This estimate holds at the real passive torus with the old rates locked.
In detail, the lower bound is the last assertion of Proposition~\ref{rate:prop:axial-mark-true}: its hypotheses hold because the new weak actions vanish, the other new actions have sum at most $K^{-(h+1)/2}$ by the moderate costs above, and every remaining marked action is at least $a_m\ge K^{-h-1}$. For the upper bound, Cauchy's estimate in the five weak Cartesian variables, on a polydisc of fixed radius $\rho_w$ inside the analytic $L^2$ chart, bounds the coefficient divided by the five weak radii by $C\rho_w^{-5}$; the remaining marked radii are at most one, so $|F_e|\le C\le K^{C_s}$.

Let $\Gamma=B_o+d^2\log(2+B)+\log(2+B)$, and allow constants
multiplying $\Gamma$ to depend on the preceding precision and the fixed
orders. Let $D_g=C_*d\log(d+2)$ and fix a finite comparison precision
$\xi=\exp(-L dD_g)$, taking $L$ sufficiently large for the actual finite
comparison, its Hessian, all coordinate counts, and both quadratures.
For each row choose the \emph{product} of the five private radii to be
$\lambda |g_e|/|F_e|$, with any fixed coefficient convention of the
physical Hamiltonian included in $g_e$. Equal fifth roots suffice.
Here $F_e$ is the actual coefficient after division by just these five
radii at the passive torus, and $g_e\ne0$ is the required graph
coefficient. Thus
\[
 K^{-C_s}\le |F_e|\le K^{C_s},\qquad
 |\log |g_e||\le C d\log(d+2),\qquad
 I_{w}^{1/2}=\lambda^{1/5}e^{O(B_o+d\log(d+2)+\log B)}.
\]
The first bound is the conclusion of the preceding paragraph. Since
$\log K=kB_o+O_s(\log B)$ it gives in particular
$e^{-C\Gamma}\le|F_e|\le e^{C\Gamma}$, which is the form used in the
remaining estimates of this section. The second bound follows from
the coefficients in Lemma~\ref{rate:l:finite}. For the third, each of
the five radii equals $(\lambda|g_e|/|F_e|)^{1/5}$, and
$|\log|F_e||\le C_s\log K=O(B_o+\log B)$.
Turn on these weak actions keeping their values fixed, lock $F_+$,
and use the interior moderator pivot inverse to restore the flat
relations $\mathsf H\chi=0$. The aggregate action change in $\ell^1$ from
the passive base is at most $\lambda^{2/5}e^{C\Gamma}$, including
the dimension factors. This is less than half every moderator,
spectator and previous action floor after the scale choices. It does
not require a weak site to serve as the pivot of this activation.
In detail, each weak action equals
\((\lambda|g_e|/|F_e|)^{2/5}\le\lambda^{2/5}e^{C\Gamma}\) by the three
displayed bounds, and there are \(5m\le Cd^2\) weak sites, so the weak
actions have \(\ell^1\) norm at most \(Cd^2\lambda^{2/5}e^{C\Gamma}\). They
change the row deficits \(\mathsf H\chi\) by at most
\(\|\mathsf H\|\,\|D\chi\|\,\|I_W\|_2\le C\sqrt q\,\|I_W\|_1\) in
\(\ell^2\). The interior pivot inverse removes this change with pivot
increments of \(\ell^2\) norm at most twice as large, and the locked
\(F_+\) actions move by at most \(\eta(1-\eta)^{-1}\) times the new
increments. Summing, and converting \(\ell^2\) to \(\ell^1\) norms at
the cost of the square root of the number of affected sites, which is
at most \(e^{\Gamma}\) by the ledger, the aggregate change is at most
\(\lambda^{2/5}e^{C'\Gamma}\). On the other
hand, every moderator, spectator and previous floor is at least
\(e^{-C''B_o-C\log B}\) by (Pr2), (Pr3) and \(a_\ell=e^{-C_2B_o}\). As
\(\lambda^{2/5}=e^{-2\epsilon QB_o/5}\) and
\(\Gamma\le(1+2\alpha^2Q)B_o+\log(2+B)\), the required inequality holds
when \(\alpha^2Q\le1\) and \(Q>5(3C'+C''+1)/(2\epsilon)\), which is
compatible with the order in which \(\epsilon\), \(Q\) and \(\alpha\)
are chosen.

For later variations the following concrete checks suffice. Write
$F_{\min}=\min_e|F_e|$ and $I_{w,\min}$ for the least current
weak action. An action movement $\Delta I$ with
\begin{equation}\label{rate:eq:registry-cube-cost}
 \|\Delta I\|_\infty\le c_q\xi I_{w,\min},\qquad
 \|\Delta I\|_1\le c_q(\xi F_{\min})^2
\end{equation}
after activation changes every weak product by relative $O_q(\xi)$
and every divided coefficient by at most $O_q(\xi F_{\min})$. For
the latter extract its marked harmonic by integration over all
angles, replacing each selected weak Cartesian side by fixed small
radius $\rho$ and its opposite by $I_w/\rho$. The map on the common
complex $\ell^2$ ball is Lipschitz; all the ordinary square-root
coordinates move by at most $\|\Delta I\|_1^{1/2}$ in $\ell^2$.
This gives $C_q(\|\Delta I\|_1^{1/2}+\|\Delta I\|_1)$ uniformly in
the number of rows. For the activation itself use only this second
inequality against the prescribed fixed product. The resulting
relative error is $\lambda^{1/5}e^{C\Gamma}\ll\xi$.
In converting relative into absolute normalized finite-vector error
also multiply by $|g_e|$ and the necessary fixed powers of the number
and size of heads; the constant $L$ above covers these costs.

Keep a fixed fraction of the spatial separation margin throughout.
The exact differential of $\mathcal Q$ gives
$|\Delta M|+|\Delta\mathcal Q|\le C\|\Delta I\|_1$ on the small ball,
independently of the largest index. Including the maximum rational
component height, a movement $\|\Delta I\|_1\le e^{-C_7 B_o}$ suffices.
In detail, \(|\Delta M|\le\|\Delta I\|_1\), and
\(d\mathcal Q=4M\,dM+2\sum I_n\,d\chi_n\) with
\(|\sum I_n\,d\chi_n|\le\|I\|_2\|D\chi\|_{\ell^2\to\ell^2}\|dI\|_2
\le C\|dI\|_1\) give \(|\Delta\mathcal Q|\le C\|\Delta I\|_1\); hence
\(A=\mathcal Q-M^2\) and \(D=\mathcal Q/2\) move by at most
\(C\|\Delta I\|_1\). For a pair \((x,y)\) of height at most
\(e^{C_3B_o}\) one has \(|x|+|y|\le2e^{C_3B_o}\), so \(Ax-Dy\) moves by at
most \(2Ce^{C_3B_o}\|\Delta I\|_1\). This is at most half of the margin
\(e^{-C_5B_o}\) if \(\|\Delta I\|_1\le e^{-C_7B_o}\) with
\(C_7>C_3+C_5+1\) and \(B\) large.
The initial spectator already pays the old $h^s$ tolerance on $F$;
the moderator cost is $K^{s-(h+1)/4}\operatorname{poly}(d)$.
Weak-site, old-adjustment and subsequent rounding costs can be made
smaller than any prescribed $e^{-D_{\mathrm{hist}}B_o}$: their indices
are at most $e^{O(B_o)}$ and their action changes are bounded by the
above powers of $\lambda$, respectively by the later small powers of
$N^{-1}$. For an old radius change use
$\|(z'-z)_F\|_{h^s}\le \langle Z\rangle^s
 |F|^{1/4}\|\Delta I_F\|_2^{1/2}$ with identical time origins and
old direct rates. For these bounds one first chooses $k$ large for
the moderator bill and then takes $Q$ sufficiently large for the
remaining prices depending on $D_{\mathrm{hist}},k$. The true current
hot actions are kept as separate high carriers.

\subsubsection{The two frequency cubes and the terminal fit}
Put $f_1=|F_+|+\#\{\text{new clock sites}\}$. In the direct-rate flat
eliminate one unit new pivot per disjoint row. Independent increments
of size $\delta$ per coordinate give full Euclidean rate increments
at most $q\sqrt{f_1}\delta$, and hence total $\ell^1$ action movement
at most $C q f_1\delta$ on the locked chart. For example the cube of
halfwidth
\[
 \delta_0=\lambda^2\exp(-C_8\Gamma)
\]
is available by choosing $C_8$ after the fixed protection prices.
It preserves positivity, \eqref{rate:eq:registry-cube-cost},
the spatial margin, and the old tolerances. The exponent $2$ on
$\lambda$ is fixed independently of the weak word cutoff and
$D_{\mathrm{hist}}$. Changing $C_8$ affects the $\Gamma$ price only.
In detail, a row has \(q\) letters with integer coefficients of modulus
one and a unit pivot, so the row equation determines the pivot
increment as a signed sum of at most \(q-1\) independent increments,
of modulus at most \((q-1)\delta\). With \(m\le f_1\) rows the full
rate increment has Euclidean norm at most
\((f_1+m(q-1)^2)^{1/2}\delta\le q\sqrt{f_1}\,\delta\). The locked
direct-rate inverse has norm at most two, so the action increment has
\(\ell^2\) norm at most \(2q\sqrt{f_1}\delta\) on at most \(2f_1\)
coordinates (including the locked ones), and \(\ell^1\) norm at most
\(2\sqrt2\,qf_1\delta\). For \(\delta=\delta_0\) this is at most
\(Cqf_1\lambda^2e^{-C_8\Gamma}\), with \(f_1\le e^{\Gamma}\) by the ledger.
The requirements are: half of every floor, the least being a weak
floor \(\lambda^{2/5}e^{-C\Gamma}\); the two conditions of
\eqref{rate:eq:registry-cube-cost}, whose right sides are at least
\(\lambda^{2/5}e^{-C'\Gamma}\) because
\(\xi=e^{-LdD_g}\ge e^{-LC_*\Gamma}\) by (Pr1) and
\(F_{\min}\ge e^{-C\Gamma}\); the spatial margin \(e^{-C_7B_o}\); and
the old tolerance \(e^{-D_{\mathrm{hist}}B_o}\). All of them hold once
\(C_8\) exceeds a fixed combination of \(C,C',C_7,D_{\mathrm{hist}}\),
since \(\lambda^2\le\lambda^{2/5}\le1\).

Let $x$ denote the independent direct rates on this cube. Fix an
integer form $n+l\cdot\chi_{F_+}+b\cdot x$ with $b\ne0$ and choose
$i$ with $|b_i|\ge1$. For each fixed choice of the other coordinates,
the set where the form has modulus less than $w$ has length at most
$2w/|b_i|$ in an interval of length $2\delta_0$. Integration in the
other coordinates gives bad relative volume at most $w/\delta_0$.
At dyadic height $2^j$ the number of tuples is at most
$(C2^j)^{f_1+1}$. Summing with thresholds
\[
 c\,\delta_0\,2^{-C(f_1+1)}
 (2+|n|+|l|_\infty+|b|_\infty)^{-2f_1-6}
\]
leaves positive relative volume. With $b=0$ use the stored bound
on $F_+$, including the independently selected spectator. Pivot
elimination has integral coefficients bounded in terms of $q$,
so the only zero affine labels are integer combinations of the
declared rows (with zero constant when using the renormalized rates).
At the physical angular cutoff $|\ell|_1\le C_JB$, an integer kinetic
part of size more than $1+C|\ell|_1$ cannot be canceled by these
bounded renormalized frequencies; all remaining heights in the affine
test are $O_q(1+B)$. By (Pr2) the inherited sites and the $O_s(1+d^2)$ current sites
together give $f_1\le C_s(1+d^2)$. Thus the dimensional log loss
$O_q(f_1\log(2+B))$ is at most $C_s(1+d^2)\log(2+B)\le C_s\Gamma$.
Thus every relevant nonflat denominator has size at least
$w_0=\lambda^2\exp(-C_9\Gamma)$, after increasing $C_9$.
Only one power of the cube width is lost in this measure estimate.

Turn on the prescribed positive hot actions $J_\infty$ at the current
axial heads, keeping all the $F_+$ and clock direct rates fixed for
this initial change. The locked actions move by
$O(\|J_\infty\|_2)$, with $\|J_\infty\|_1\le r^2 e^{C d}$.
Then keep $J_\infty$ and all independent nonpivot rates fixed and
write $y_e=H_e\chi$ for the variable row displacements; one unit
pivot per row realizes these displacements. Denote by
$\mathcal D_e(y)$ the sum of the actual carrier frequency corrections
on the two inputs minus that on the output, in the rotating branch
convention. The kinetic key is one. Use the true one-dimensional
frequencies on hot axes, including their own actions. On transverse
heads use the full angle means of the required fixed-order instantaneous
spatial scalar. These means do not use the temporal divisors about to
be selected. At the flat center
$\|(\mathcal D_e(0))_e\|_2\le N^{-2}E$ after summing the heads.
Indeed (T1) gives the leading terms, all higher terms cost
$N^{-3}E$, and $\|J_\infty\|_1$ is smaller.
In detail, the kinetic key is exactly one on every designated edge, so
\(\mathcal D_e(0)\) is the signed sum of the three head corrections. A
marked head has horizontal index \(|n|\ge Ne^{-Cd}d^{-C}\), because the
normalized directions satisfy \(|\cos\theta|\ge c/d\) and the
normalized lengths are at least \(e^{-Cd}\). On a transverse head the
correction is \(d_2(n,m)+d_{\ge3}\) with
\(|d_2|\le C/n^2\) by (T1) and \(|d_{\ge3}|\le E_c|n|^{-3}\); on an
axial head with the hot actions set to zero, the correction is
\(\chi_K=(A-D)K^{-2}+O(E_o|K|^{-3})\) with \(A-D=\mathcal Q/2-M^2\), by
Lemma~\ref{rate:lem:trans-axial-expansion}. The hot actions change
these by at most \(C\|J_\infty\|_2\le Cr^2e^{Cd}\le N^{-3}\). Thus every
head correction is at most \(Ce^{Cd}d^{C}N^{-2}\), every
\(|\mathcal D_e(0)|\) at most three times this, and with at most
\(Cd^2\) edges \(\|(\mathcal D_e(0))_e\|_2\le Cd^{C+1}e^{Cd}N^{-2}\le
N^{-2}E\) by (Pr1).

On a real positive $F_+$-plus-clock neighborhood of Euclidean radius
$\lambda^{C_s}e^{-C\Gamma}$ with the hot actions temporarily zero,
a high empty axial site $n$ has
$|\chi_n|\le E'/n^2$, $E'\le\lambda^{-C_s}e^{C\Gamma}$.
Here $|n|\ge N e^{-C d}d^{-C}$ on every marked axis. Apply the
central-difference estimate at $\sqrt{E'}/|n|$ within this neighborhood
to bound every unit occupied-direction derivative by
$C\sqrt{E'}/|n|$. Restoring the fixed hot actions changes this bound
by at most $C\|J_\infty\|_2$ using the uniform second differential
of $\chi$. The transverse mean derivatives cost at most
$N^{-2}\lambda^{-C_s}E$ each. All these bounds are in a fixed number
of differentiations. The direct-rate inverse and the passage from
individual rows to the row-vector norm therefore give
$\|D_y(\mathcal D_e)_e\|<1/2$ on $\|y\|_2\le 2N^{-2}E$,
increasing the constant in $E$ before defining this ball.
Write $\mathcal D=(\mathcal D_e)_e$. Since
$\|\mathcal D(0)\|_2\le N^{-2}E$, the derivative bound implies
$\|\mathcal D(y)\|_2\le N^{-2}E+\|y\|_2/2$ on this ball.
It follows that the ball is invariant and that $y=\mathcal D(y)$
has a unique solution there. This single vector fixes the
mismatches of all designated edges. Throughout this adjustment
the hot actions $J_\infty$, the rates on $F_+$ and the independent
nonpivot rates retain their prescribed values.
The derivative bound for a hot axial frequency is the central-difference
estimate used for the spectator above. For a unit occupied direction
$e$ put $h=\sqrt{E'}/|n|$; it is smaller than the radius
$\lambda^{C_s}e^{-C\Gamma}$ of the neighborhood since
$|n|\ge Ne^{-Cd}d^{-C}$. Then
\[
 |\partial_e\chi_n|\le\frac{|\chi_n(I+he)-\chi_n(I-he)|}{2h}
   +\frac h2\sup|\partial_e^2\chi_n|\le\frac{E'}{n^2h}+Ch
   =(1+C)\frac{\sqrt{E'}}{|n|}.
\]
Since the direct-rate inverse has norm at most two and each edge involves
three heads, $\|D_y\mathcal D\|$ is at most $6\sqrt{\#\{e\}}$ times the
largest $\ell^2$ norm of a head gradient, that is
$O\bigl(\#\{e\}^{1/2}[\sqrt{E'}N^{-1}e^{Cd}d^C+N^{-2}\lambda^{-C_s}E
+\|J_\infty\|_2]\bigr)$, which is below $1/2$ for large $B_*$.

Denote by $\Omega_*$ the final, rounded vector of physical rates on the
clock support $F$, and by $\eta_{\mathrm{rnd}}$ a common bound for the
rounding error in the tested frequency combinations. For each relevant
nonflat label the actual perturbation of its divisor is at most $C_J B N^{-2}E$ (plus the tiny final rounding). This is
$o(w_0)$, including the additional real corrections of any fixed
number of moving carriers. On the frozen module, charge-zero
combinations of bounded weak length obey, in root-normalized
coordinates,
\[
 \sum_e b_e H_e\cdot\Omega_*
   =N^{-2}\sum_e b_e\mathfrak d_e
       +O_{s,\epsilon}\big(N^{-3}\lambda^{-C_s}E+
                             r^2 E+\eta_{\mathrm{rnd}}\big),
 \qquad \sum_e b_e=0 .
\]
The spatial choice guarantees
$|\sum b_e\mathfrak d_e|\ge e^{-C B_o}$ for each required nonzero
word, after allowing the protected movement of the spectator ratio.
Here $\mathfrak d_e$ is the rational geometric defect paired with
the actual $(A,D)$ from (T1). The total $\sum|b_e|$ is bounded by
the weak length since the row weak supports are disjoint. The same
estimate, with the tested extra head corrections, applies to required
wrong charged coincidences. In all cases the displayed remainder
is $o(N^{-2}e^{-C B_o})$ after the scale choices. This explains
separately why both the larger nonflat gap and the smaller geometric
gap survive. The final data are not required to lie on the original
flat.

\smallskip
\noindent\emph{The two gaps in detail.} (i) \emph{Nonflat labels.} The
tilt from the flat center to \(\Omega_*\) keeps the rates on \(F_+\)
and the nonpivot rates fixed and moves each pivot rate by
\(|y_e|\le\|y\|_2\le2N^{-2}E\). A label with \(|\ell|_1\le C_JB\)
therefore changes by at most \(2C_JBN^{-2}E\), plus the rounding
\(\eta_{\mathrm{rnd}}\). As \(w_0=\lambda^2e^{-C_9\Gamma}\) and
\(\lambda^2=e^{-2\epsilon B}\), while \(N\ge R=e^B\), \(\log N\le2B\) by
(Pr4) and \(E\le N^\gamma\), one has
\(N^{-2}EB\le e^{-2B+2\gamma B+o(B)}\); with
\(\Gamma\le(Q^{-1}+2\alpha^2)B+o(B)\) the ratio of the change to \(w_0\)
is at most \(e^{-2B+2\gamma B+2\epsilon B+C_9(Q^{-1}+2\alpha^2)B+o(B)}\to0\).
Likewise \(\eta_{\mathrm{rnd}}/w_0\le u/w_0
=HN^{-s}\lambda^{-1}e^{C_9\Gamma}\le
e^{-sB+\epsilon B+C_9(Q^{-1}+2\alpha^2)B+o(B)}\to0\).
(ii) \emph{Module combinations.} On the flat center every row has
\(H_e\cdot\Omega=1\), since its word has free frequency one and zero
\(\chi\)-sum there. At the fixed point the row displacement is
\(y_e=\mathcal D_e(y)\), the signed sum of the three head corrections,
which by (T1), Lemma~\ref{rate:lem:trans-axial-expansion} and the hot
action bound equals \(N^{-2}\mathfrak d_e+O(N^{-3}\lambda^{-C_s}E+r^2E)\),
where \(\mathfrak d_e\) is evaluated with the actual
\((\mathcal A,\mathcal D)\) in root-normalized directions. For
\(\sum_eb_e=0\) the constant \(1\) cancels, which gives the displayed
identity, the rounding contributing \(\eta_{\mathrm{rnd}}\).
(iii) \emph{Separation.} Writing \((U_e,V_e)\) for the designated
component pair of the edge \(e\), \(\mathfrak d_e=\mathcal AU_e-\mathcal DV_e\),
so \(\sum_eb_e\mathfrak d_e=\mathcal Ax-\mathcal Dy'\) with
\((x,y')=\sum_eb_e(U_e,V_e)\). For \(b\ne0\) with
\(\sum|b_e|\le m_*\), Lemma~\ref{rate:l:geom} at the parameters of
Lemma~\ref{rate:l:denom} gives \((x,y')\ne(0,0)\), and its height is
at most \(e^{C_3B_o}\). The spectator choice gives
\(|\mathcal Ax-\mathcal Dy'|\ge e^{-C_5B_o}\), and the protected
movement \(\|\Delta I\|_1\le e^{-C_7B_o}\) keeps at least half of it.
(iv) \emph{The remainder} is \(o(N^{-2}e^{-CB_o})\): indeed
\(N^{-1}\lambda^{-C_s}Ee^{CB_o}\le
e^{-B+C_s\epsilon B+2\gamma B+CB/Q+o(B)}\to0\) for \(\epsilon,\gamma\)
small and \(Q\) large, \(r^2EN^2e^{CB_o}=H^2N^{2-2s}Ee^{CB_o}\to0\), and
\(\eta_{\mathrm{rnd}}N^2e^{CB_o}\le uN^2e^{CB_o}\le N^{2-s+o(1)}e^{CB/Q}\to0\)
for \(Q\) large, because \(\eta_{\mathrm{rnd}}\le u=HN^{-s}\lambda\) and
\(s>2\).

Finally release all new clock and hot direct rates on a full cube,
still pinning $F_+$. Write $f_2$ for the dimension including every
accepted occupied site. The full inverse has Euclidean norm at most
two as before, so $\|\Delta I\|_1\le C f_2\delta_1$ for cube
halfwidth $\delta_1$. For example one may take
$\delta_1=N^{-A_1}$ with the fixed exponent $A_1=4s+30$.
Here are the requirements this choice pays. The preliminary minimum
hot floor is at least $N^{-2s}e^{-C d\log(d+2)}$ from the finite
endpoint (root amplitude $r$, private axis radii $r c_{n,k}b_p$,
with $b_p=e^{-c_1D_g}$ and $c_{n,k}\ge1$). Let $I_{\min}$ be the minimum
of all positive actions at the center of this cube, including the new hot axes.
Require
\[
 C f_2\delta_1\le\min\left\{
 \tfrac14 I_{\min},\ c\,r^2\xi^2,\ c\,\xi I_{w,\min},
 c\,(\xi F_{\min})^2,\ e^{-C_7B_o}\right\}.
\]
Then the normalized hot endpoint moves by at most $\xi$, by the
square-root inequality after dividing by $r$, and the weak amplitudes,
divided coefficients, old actions and all floors keep their earlier
prices. In detail, $|\sqrt x-\sqrt y|\le\sqrt{|x-y|}$ and
$\|\Delta I\|_1\le Cf_2\delta_1\le c\,r^2\xi^2$ give a change of at
most $\sqrt c\,r\xi$ in every hot radius; the second, third and
fourth entries of the minimum are the hypotheses of
\eqref{rate:eq:registry-cube-cost} and of the hot rounding, and the
last one is the movement allowed by the spatial margin. With
$\delta_1=N^{-A_1}$ the left side is at most $N^{-4s-29}$ for large
$B_*$, whereas each of the five quantities on the right is at least
$N^{-2s-4}$ once the fractions of $B$ spent on $B_o$, on
$d^2\log B$ and on $\epsilon B$ are small: the hot floor is
$N^{-2s}e^{-Cd\log(d+2)}$, $r^2\xi^2=H^2N^{-2s}e^{-2LdD_g}$, and the
remaining factors are $e^{-O(\Gamma)}$ times fixed powers of
$\lambda$ and $\xi$. In the full chart the high axial frequency is a direct
coordinate and the transverse mean has the preceding small derivative
in the clock variables. In the hot directions the pure old transverse
means change only through the compensating occupied actions on the
direct chart. Hence every fitted discrepancy moves by at most
$C q\sqrt{f_2}\delta_1$ times a bounded derivative factor, and all
old-clock nonflat labels by $C_J B\delta_1$. In particular one can
take $\eta_{\mathrm{rnd}}\le u\,2^{-d}\exp(-L dD_g)$ in the finite
incidence and divisor estimates. To check the displayed budgets,
take $Q$ sufficiently large that every fixed multiple of
$C(D_{\mathrm{hist}},k)B_o$ is smaller than the previously assigned
fraction of $\epsilon B$, and take $\alpha$ small enough for the
same to hold of the finitely many multiples of $d^2\log B$.
The prefactors $\operatorname{poly}(f_2,d,B)$ and $e^{C d}$ are
absorbed by the remaining slack, increasing $B_*$ last. This
choice of $A_1$ does not depend on $D_{\mathrm{hist}}$.

Choose the final point in this cube by the same affine slicing with
$\delta_1,f_2$ in place of $\delta_0,f_1$, using the stored $F_+$
bound for pure-prefix forms. This registers \eqref{rate:eq:registry-prefix}
on all occupied sites for the next stage. For example the new
prefactor is at least $N^{-(A_1+1)}$ once $B_*$ is large: the new
factor $2^{-C f_2}$ costs less than $N$, and the two inherited
prefactors exceed this bound after enlarging $Q$. Since
$\log N\le(1+\gamma)B+O(1)$ this fixes a common $C_{\mathrm{reg}}$ at
the outset, e.g. $2(A_1+2)$. Its use on the previous scale is not
an input to $A_1$. The actual final least action is at least
$N^{-(2s+4)}$ and the greatest site at most $N^2$. Incremental
changes on each existing coordinate may be prescribed below
$e^{-D_{\mathrm{hist}}B_o}$ with $D_{\mathrm{hist}}$ greater than these
coarse floor and historical exponents. Their later sum at a site
born at stage $i$ is $O(e^{-D_{\mathrm{hist}} B_i})$, smaller than
half its birth floor after an initial enlargement of $B_*$. Rates
and origins already accepted at previous stages remain exactly fixed.
At the final frequencies use one fresh weak angle per disjoint row
to set the constant argument of its actual nonzero driving harmonic.
These private choices commute, do not move any action or rate, and
are stored with all the other newly accepted angle origins.

\subsection{The axial potential on the transverse Fourier space}
\label{rate:sec:transverse-principal}
We prove the estimates for the part of the transverse Hessian generated
solely by the occupied axial field. They concern the operator on arbitrary
inputs, including its full low horizontal subspace. No Fourier support or
smallness condition on an input in that subspace enters this calculation.

Write $F$ for the occupied sites including the prepared clock, before
the current hot Cartesian axes are added, $f=|F|$,
$n_F=(m)_{m\in F}$, $Z=1+\max_{m\in F}|m|$, and
\[
 V(I,\phi)=\Psi((\sqrt{I_j}e^{-i\phi_j})_{j\in F}),\quad
 M=\sum I_j,\quad a=2(|V|^2-M),\quad c=V^2,\quad
 {\mathcal{D}}=\Omega(I,0)\cdot\partial_\phi .
\]
In this sign convention $\dot\phi=\Omega(I,0)$ on the isolated torus
in its own mass gauge. Use the stage scales $m=Rq$, $N=RP$,
$B_o=B/Q$, $\lambda=e^{-\epsilon B}$, and
$E_c=\exp(C_{s,\epsilon}(B_o+d^2\log(2+B))+o(B))$.
A fixed number of changes in the constants of $E_c$ is allowed.

\subsubsection{The Schur estimates and the role of the angles}
On real phases $a$ is real. For all complex angle estimates use the
analytic reflection of the conjugate field. The old Cartesian coordinates
lie in a smaller complex $\ell^2$ chart ball on
$|\operatorname{Im}\phi_j|\le a_0$ with an extra $\epsilon B/6$
in each weak angle. Indeed their new weak radii are
$\lambda^{1/5}e^{O(B_o+d\log d)}$ and their number is polynomial in
$d$. Old and moderator coordinates use the fixed strip and their
real total squared size in the small ball. The same conclusion holds
on a sufficiently small relative action neighborhood. Fix the desired
final number $k$ of angle and action derivatives and any finite
asymptotic orders before choosing $\epsilon$. A fixed reduction of
the common strip controls every such angular derivative, leaving,
for example, width $\epsilon B/7$ in the weak angles above a smaller
fixed width. Each genuine action derivative costs
$E_c\lambda^{-C_s}$ from Cartesian square roots and the common
action-frequency chart. Previous occupied floors and new moderator
floors cost $e^{O_s(B_o)}$. Current hot floors are absent from $F$.
In detail, at fixed complex angles in the strip the map
\(I_j\mapsto V(I,\phi)\) is analytic on the disc
\(|I_j-I_j^0|\le I_j^0/2\): the square root is analytic there, and the
Cartesian point stays in the chart ball since
\(|\sqrt{I_j}|\le(3/2)^{1/2}\sqrt{I_j^0}\). Cauchy's estimate on the
product of these discs bounds any derivative of order \(k\) in the
actions by \(k!\,(I_{\min}/2)^{-k}\) times the supremum over the
polydisc, and that supremum obeys the same strip and moment bounds,
hence costs \(E_c\). By the ledger the nonweak positive actions are at
least \(e^{-CB_o-C\log B}\) and the weak ones at least
\(\lambda^{2/5}e^{-C\Gamma}\), so \((I_{\min}/2)^{-k}\le
2^k\lambda^{-2k/5}e^{Ck\Gamma}\), which is of the form
\(E_c\lambda^{-C_s}\) for every fixed \(k\) with \(C_s\ge2k/5\). The
mass-gauged frequencies \(\Omega(I,0)=(n^2+\chi_n(I,0))_{n\in F}\) have bounded
action derivatives of every fixed positive order on the common \(\ell^2\) ball,
since \(\chi\) is analytic there, so the action derivatives of
\({\mathcal D}=\Omega(I,0)\cdot\partial_\phi\) cost only angular
derivatives, which are paid by the strip.

In reflected physical labels each angular coefficient of the old
matrix shifts the horizontal coordinate by a prescribed signed old
momentum $h=\pm\ell\cdot n_F$. Thus its absolute sum over angular
labels, after the supremum over input frequency, bounds both weighted
Schur sums. Use both directions of the weight ratio to include
adjoints. For example the packet weights with polynomial factor
$(1+|n|/N+|q|/P)^p$ and any of the capped exponential word factors
have ratio bounded by
$(1+|h|/N)^p\exp(|\ell|_1/B^{K_0})$, up to a fixed constant for the
choice of signed generators. Here $K_0>5$ is the fixed exponent in the word weight. Multiplicative constants for a fixed
change of copy can be included by using the word weight on the
reflected copy induced by actual conjugation. Fixed spatial moments
and nominal time derivatives cost powers of $Z$ and $Z^2$, respectively,
and powers of $|\ell|_1$. On a smaller strip the absolute sum costs
$C^f$, since $\sum_{\ell}e^{-\delta|\ell|_1}=C_\delta^f$.
These losses lie in $E_c$ by the stage bounds. They also control
derivatives on the enlarged weak strip standing behind the weak-label
tails. In the following, matrices bounded in this translation algebra
act both on ordinary $L^2$ and on the indicated weighted spaces.
The assertion does not use a word-length lower bound on low indices.

\subsubsection{A spatial column, its scalar, and the empty axial frequency}
With $D_x=-i\partial_x$, reflect the lower horizontal index so that the
two envelopes have the same spatial carrier. The mass-gauged old operator
on a nonzero vertical slice is
\[
 L_m=\begin{pmatrix}D_x^2+m^2+a&c\\
                -\bar c&-(D_x^2+m^2+a)\end{pmatrix},
 \qquad \Sigma=\operatorname{diag}(1,-1).
\]
It is selfadjoint for the form defined by $\Sigma$. This doubled
equation follows by differentiating the physical cubic and adjoining
the conjugate equation. In the original, unreflected second index the
pairing entry is $\widehat c(n+n')$.

Fix $H_{\rm geom}=C(1+d)P$, with an absolute $C\ge1$ chosen so that
$|X|+|Y|\le H_{\rm geom}$ for every marked head. Throughout this stage
use $H_x=c_0R/H_{\rm geom}$ with a fixed $0<c_0<1/100$.
On each nonzero vertical slice the low input projection is
$Q_l=\mathbf 1_{\{|n|<H_x\}}$. For a high input $|n|\ge H_x$, put
$E_n=n^2+m^2$, $\theta=n^2/E_n$, $\delta=(2n)^{-1}$.
A column $e^{inx}(f,g)$ with instantaneous scalar $E_n+e$ satisfies
the two equations, up to the remainder of the truncation,
\[
 (2nD_x+D_x^2+a-i{\mathcal{D}}-e)f+cg=0,\qquad
 (2E_n+2nD_x+D_x^2+a+i{\mathcal{D}}+e)g+\bar c f=0 .
\]
Expand $f,g,e$ in powers of $\delta$, with $f_0=1$, $g_0=g_1=0$,
and $\langle f_j\rangle_x=0$ for $j\ge1$. With missing coefficients
zero the recursion is
\begin{align*}
 e_j&=\langle a f_j+c g_j\rangle_x,\\
 D_x f_{j+1}&=-(D_x^2+a-i{\mathcal{D}})f_j-cg_j+
                         \sum_{l+k=j}e_l f_k,\\
 g_{j+2}&=-2\theta\Big(\bar c f_j+D_xg_{j+1}
       +(D_x^2+a+i{\mathcal{D}})g_j+\sum_{l+k=j}e_lg_k\Big).
\end{align*}
Indeed $2n=\delta^{-1}$ and $2E_n=(2\theta\delta^2)^{-1}$, so the
coefficients of $\delta^j$ in the two column equations are
\begin{align*}
 D_xf_{j+1}+(D_x^2+a-i{\mathcal{D}})f_j
      -\sum_{l+k=j}e_lf_k+cg_j&=0,\\
 \frac{g_{j+2}}{2\theta}+D_xg_{j+1}+(D_x^2+a+i{\mathcal{D}})g_j
      +\sum_{l+k=j}e_lg_k+\bar cf_j&=0 .
\end{align*}
The spatial mean of the first line is
$\langle af_j+cg_j\rangle_x-e_j$, because $f_0=1$ and
$\langle f_k\rangle_x=0$ for $k\ge1$; this is the formula for $e_j$,
and the mean-free part determines $f_{j+1}$.
The upper inverse is $D_x^{-1}$ on zero spatial mean and costs at
most one in its Fourier norm. The lower recurrence contains $\theta$
as a factor in every $g_j$. Even as $|m|\to\infty$, a nominal
$2E_n$ acting on such a term costs at most $2n^2$ times its bounded
coefficient with this factor removed. Retaining $f$ through order
$M_0+1$, $g$ through $M_0+2$ and $e$ through $M_0$ gives residual
$E_c |n|^{-M_0}$ with the translation moments and angular derivatives
above. Each coefficient uses at most $2(M_0+3)$ spatial and $M_0+3$
nominal time derivatives on any of its factors. Only a bounded number,
depending on $M_0$, of factors occurs. Supremum in the input column
followed by the absolute angle sum gives the same residual as an
operator estimate.
To see the residual, insert the truncations
\(f^{\rm tr}=\sum_{j\le M_0+1}\delta^jf_j\),
\(g^{\rm tr}=\sum_{j\le M_0+2}\delta^jg_j\) and
\(e^{\rm tr}=\sum_{j\le M_0}\delta^je_j\) into the two column equations,
written with \(2n=\delta^{-1}\) and \(2E_n=(2\theta\delta^2)^{-1}\). By
the recursion the coefficients of \(\delta^j\), \(j\le M_0\), vanish in
both equations; in the lower one this uses \(g_{j+2}\) and
\(D_xg_{j+1}\), which are retained up to \(j=M_0\). What remains are the
terms of order \(\delta^j\), \(M_0+1\le j\le2M_0+2\), namely
\((D_x^2+a\mp i{\mathcal D})\) applied to the last retained
coefficients, \(D_xg_{M_0+2}\), the products \(e_lf_k\), \(e_lg_k\)
with \(l+k\ge M_0+1\), and \(cg_j\), \(\bar cf_j\) with
\(j\ge M_0+1\). The recursion forms each coefficient from \(a,c,\bar c\)
by multiplication, spatial and nominal time differentiation, spatial averaging,
and application of \(D_x^{-1}\) to zero-mean terms. At fixed \(M_0\) the number
of operations is bounded. The last two operators have norm at most one in
the absolute Fourier sums, and every lower coefficient carries the factor
\(\theta\le1\). By the strip and moment
bounds each such product has translation norm at most \(E_c\), so the
residual is at most \(C_{M_0}E_c|\delta|^{M_0+1}\le E_c|n|^{-M_0}\),
uniformly in \(m\). The term \(2E_ng^{\rm tr}\) never appears in the
remainder: it is matched exactly through \(g_{M_0+2}\).

For the column $h=e^{inx}(f,g)$ put
$\mu=\langle h,\Sigma h\rangle=1+O(E_c/|n|)$ on real phases.
Taking the signed norm derivative in
$(L_m-i{\mathcal{D}})h=(E_n+e)h+R$ gives
${\mathcal{D}}\mu=-2\operatorname{Im}(e)\mu+O(E_c|n|^{-M_0})$.
The normalized column $\mu^{-1/2}h$ therefore has scalar $\Re e$
up to the same remainder. The analytic reflection used in $\mu$
is its holomorphic continuation on the narrower strip.
Here $e_0=e_1=0$, $f_1=-D_x^{-1}a$ is purely imaginary, and
$g_2=-2\theta\bar c$. Bilinear integration by parts gives
\[
 e_2=\langle a^2\rangle-2\theta\langle|c|^2\rangle+
                           \tfrac i2{\mathcal{D}}\langle f_1^2\rangle .
\]
Indeed $\langle a f_2\rangle=
-\langle f_1(D_x^2+a-i{\mathcal{D}})f_1\rangle$ and
$\langle a f_1^2\rangle=0$. Thus the real order-two term in physical
units is $\{(1-\theta/2)\langle|V|^4\rangle-M^2\}/n^2$.
In detail, $e_0=\langle a\rangle=0$ and
$e_1=\langle af_1\rangle=-\langle aD_x^{-1}a\rangle=0$. The pairing
$\langle\cdot\rangle$ is bilinear, so $D_x$ is antisymmetric for it.
With $a=-D_xf_1$ and $D_xf_2=-(D_x^2+a-i{\mathcal{D}})f_1$ this gives
$\langle af_2\rangle=\langle f_1D_xf_2\rangle$ and
\[
 -\langle f_1D_x^2f_1\rangle=\langle(D_xf_1)^2\rangle=\langle a^2\rangle,
 \qquad
 \langle af_1^2\rangle=-\tfrac13\langle D_x(f_1^3)\rangle=0,
 \qquad
 i\langle f_1{\mathcal{D}}f_1\rangle=\tfrac i2{\mathcal{D}}\langle f_1^2\rangle ,
\]
while $\langle cg_2\rangle=-2\theta\langle|c|^2\rangle$. The last
term of $e_2$ is imaginary on real phases, since $f_1^2$ is real.
Finally $\delta^2=1/(4n^2)$, $\langle a^2\rangle=4(\langle|V|^4\rangle-M^2)$
and $|c|^2=|V|^4$, so
$\operatorname{Re}(e_2)\delta^2
=\{\langle|V|^4\rangle-M^2-\tfrac\theta2\langle|V|^4\rangle\}/n^2$.

Let $\alpha_{\mathrm{phys}}=-\operatorname{Im}\langle\bar V\,dV\rangle$ and
$\alpha_{\mathrm{act}}=-\operatorname{Im}\sum\bar z_j\,dz_j$.
Their difference pulled back to Cartesian coordinates is closed by
canonicity. On the chart ball the radial formula gives its primitive
$G(z)=\int_0^1(\Psi^*\alpha_{\mathrm{phys}}-\alpha_{\mathrm{act}})_{tz}(z)\,dt$
with $G(0)=0$. This primitive is invariant under the common phase and
spatial translation, since both forms and radial contraction are.
Its fixed derivatives have the same chart prices already established.
In the ungauged physical flow the two evaluations are
$\langle |V_x|^2+|V|^4\rangle$ and $\sum I_j\omega_j$. Passing to
the mass gauge subtracts the same multiple of $M$ from both.
It follows with $\mathcal H=\langle |V_x|^2+|V|^4/2\rangle$ that
\[
 \mathcal Q=2\Big(\sum I_j\omega_j-\mathcal H\Big),\qquad
 \langle|V|^4\rangle=\mathcal Q+2{\mathcal{D}}G,\qquad
 d\mathcal Q=4M\,dM+2\sum I_j\,d\chi_j .
\]
For the last identity use $d\mathcal H=\sum\omega_j\,dI_j$ and
$\omega_j=j^2+2M+\chi_j$. Define $b_2=(2/n^2-1/E_n)G$.
Right multiplication of the column by $e^{-i b_2}$ subtracts
${\mathcal{D}}b_2$ from its scalar. The order-two constant and its
action derivative along any fixed parameter tangent are precisely
\[
 d_2(n,m)=\frac{\mathcal Q-M^2}{n^2}
              -\frac{\mathcal Q}{2E_n},\qquad
 d d_2=\frac{d\mathcal Q-2M\,dM}{n^2}
              -\frac{d\mathcal Q}{2E_n}.
 \tag{T1}
\]
There is no temporal divisor in these formulas. To verify (T1),
substitute $\langle|V|^4\rangle=\mathcal Q+2{\mathcal{D}}G$ into the
real order-two term and use $\theta/n^2=1/E_n$:
\[
 \frac{(1-\theta/2)\langle|V|^4\rangle-M^2}{n^2}
 =\frac{\mathcal Q-M^2}{n^2}-\frac{\mathcal Q}{2E_n}
   +{\mathcal{D}}\Big[\Big(\frac2{n^2}-\frac1{E_n}\Big)G\Big].
\]
The last term is ${\mathcal{D}}b_2$. The identity for the quartic
mean is the evaluation of $\Psi^*\alpha_{\mathrm{phys}}-\alpha_{\mathrm{act}}=dG$
on the flow: in the ungauged flow
$\alpha_{\mathrm{phys}}(\dot V)=\langle|V_x|^2+|V|^4\rangle$ and
$\alpha_{\mathrm{act}}(\dot z)=\sum I_j\omega_j$, so that
$\langle|V_x|^2+|V|^4\rangle-\sum I_j\omega_j={\mathcal{D}}G$, and
subtracting $\mathcal H=\langle|V_x|^2+\tfrac12|V|^4\rangle$ gives
$\tfrac12\langle|V|^4\rangle=\tfrac12\mathcal Q+{\mathcal{D}}G$.
For a single plane wave $V=\sqrt M\,e^{ikx}$ the orbit evolves by a
common phase, so $\mathcal DG=0$. The quartic identity then gives
$\mathcal Q=M^2$ and hence $d_2=-M^2/(2E_n)$, which is the
second-order term of the Bogoliubov dispersion and serves as a check
of the normalization.

\begin{lemma}\label{rate:lem:trans-axial-expansion}
On the isolated smooth occupied torus, the empty one-dimensional mass-gauged
frequencies have an expansion, through any fixed order $h_0$, of the
form $\chi_n=\sum_{j=2}^{h_0}c_j(I)n^{-j}+O(E_o|n|^{-h_0-1})$,
with common coefficients for positive and negative $n$.
For each fixed truncation order $h_0$ and fixed derivative orders,
$E_o$ is polynomial in $1+Z$ and $e^f$, with additional polynomial
dependence on the inverse occupied floors for action derivatives.
In particular $c_2=\mathcal Q/2-M^2$.
\end{lemma}
\begin{proof}
Use the direct envelope recurrence just proved at $m=0$ on this one
large input $n$, so $\theta=1$, and truncate two orders past $h_0$.
This step does not use any global off-copy diagonalization at $m=0$.
The scalar coefficients depend on $\theta$ and the old field but not
on the sign of $n$. The axial mass derivative in the full physical
linearization adds a rank-one term proportional, in the time generator,
to the common phase tangent of $V$. Its forward Birkhoff differential
is the corresponding phase tangent at the old Cartesian vector and
has zero component at an empty site.

Apply the actual forward Birkhoff differential at $V$ to the column
at $m=0$ and take its upper empty component $c_n(\phi)$. This is
$1+O(E_o/|n|)$: the input has $H^1$ norm $O(E_o |n|)$ and the
nonlinear part of that differential maps $H^1$ to $h^2$, while the
free upper carrier coefficient is $1+O(E_o/|n|)$. The projected
$L^2$ residual has size $O(E_o |n|^{-h_0-1})$ by the $L^2$ chart
differential. An empty coordinate has zero term from the variation of
its frequency multiplying its background value. Therefore the actual
linearized equation, with the carrier phase factored out, gives
$i{\mathcal{D}}c_n=(\chi_n-e)c_n+O(E_o |n|^{-h_0-1})$.
The periodic angle function $c_n$ has a single-valued logarithm for
large $|n|$. Divide by $c_n$ and take the full angular Haar mean.
The derivative of its logarithm integrates to zero for any real
frequency vector. This gives the stated expansion; (T1) at $m=0$
gives its leading coefficient. The remainder is uniform on the
smaller old real or complex chart with the corresponding derivatives.
\end{proof}

\subsubsection{A low complement with a bounded lower-order part}
Now $m=Rq\ne0$ throughout. First remove old anomalous entries in all
horizontal frequencies. For $D_0=\operatorname{diag}(D_x^2+m^2,
-D_x^2-m^2)$ and the present off-copy upper block $C_{nj}$ take
\[
 K_{+,-}(n,j)=-C_{nj}/(n^2+j^2+2m^2),\qquad K_{-,+}=K_{+,-}^* .
\]
The signed adjoint relation gives $K^*\Sigma+\Sigma K=0$ and
$[D_0,K]=-B_{\mathrm{off}}$. Indeed the $(+,-)$ entry of $[D_0,K]$
at $(n,j)$ is
$\{(n^2+m^2)+(j^2+m^2)\}K_{+,-}(n,j)=-C_{nj}$, and likewise for the
other block. Its norm gains $(2R^2)^{-1}$ even on two
low horizontal indices, because $n^2+j^2+2m^2\ge2R^2$. With $B$ the current bounded perturbation,
the exact conjugation identity reads
\begin{align*}
 &e^{-K}(D_0+B-i{\mathcal{D}})e^K-(D_0+B-i{\mathcal{D}})\\
 &\hspace{9mm}=\int_0^1e^{-\theta' K}
 \big([D_0,K]+[B,K]-i{\mathcal{D}}K\big)e^{\theta' K}\,d\theta'.
\end{align*}
The bounded kinetic commutator preserves its graph domain, first
by approximation on finite Fourier inputs. The new off-copy remainder
gains $E_cR^{-2}$ per step after a further fixed strip loss. At $j_0$
steps it is $E_cR^{-2j_0}$. The diagonal addition stays bounded,
and each diagonal copy is selfadjoint for ordinary Hilbert space on
real phases. The conjugate branch is obtained from the reality
involution, and the resulting changes preserve the actual Fourier coset.

Apply the inverse of this change to the direct high columns above and
project to the positive Hermitian copy. Write its generator as $A$,
the resulting map from the high input subspace as $T_h$, and its real
instantaneous high scalar as $\Lambda$. We can leave the exact phase
$b_2$ and all other temporal phases until after the following completion.
From the full column equation,
\[
 A T_h-i{\mathcal{D}}T_h-T_h\Lambda=R_h,\qquad
 \|R_h\|\le E_c(H_x^{-M_0}+R^{-2j_0}).
\]
The discarded branch feeds back into this row only through the
small off-copy remainder. The high columns are allowed to extend
into the low horizontal output. For $G_h=T_h^*T_h$,
\[
 [\Lambda,G_h]-i{\mathcal{D}}G_h=T_h^*R_h-R_h^*T_h .
\]
This derivation identity is valid with angle-dependent real $\Lambda$.
To obtain it, take the Hilbert adjoint of the column equation on
real phases, $T_h^*A+i{\mathcal{D}}T_h^*-\Lambda T_h^*=R_h^*$,
multiply the column equation on the left by $T_h^*$ and the adjoint
equation on the right by $T_h$, and subtract: the terms with $A$
cancel, and
$-iT_h^*{\mathcal{D}}T_h-i({\mathcal{D}}T_h^*)T_h=-i{\mathcal{D}}G_h$.
The matrix $T_h$ differs from the high inclusion by at most
$E_c(H_x^{-1}+R^{-2})$; the convergent inverse-square-root series,
with the product rule for the same derivation, proves that
$U_h=T_hG_h^{-1/2}$ is isometric with the same residual price.
Denote its new residual by $R'_h$. For the fixed low input projection
$Q_l$ put $\Pi=1-U_hU_h^*$ and define
\[
 U_l=\Pi Q_l\,(Q_l\Pi Q_l|_{\mathrm{ran}\,Q_l})^{-1/2}.
\]
Since $\|\Pi-Q_l\|<1$, $\Pi Q_l\Pi$ is invertible on
$\mathrm{ran}\,\Pi$, so the range of $U_l$ is all of that complement.
The completed matrix $U=(U_h,U_l)$ is unitary. Its high-to-low block
in $U^* A U-iU^*{\mathcal{D}}U$ is $U_l^* R'_h$, and its low-to-high
block is its Hilbert adjoint. No truncation of a high column at its
output is performed.

The compressed low block is free kinetic plus a bounded Hermitian
operator of norm at most $E_c$. In fact a spatial shift by $j$ in a
nonconstant upper envelope costs $E_c/|n|$ with enough shift moments
to pay $|(n+j)^2-n^2|\le2|nj|+j^2$. On individual entries of a
Sylvester generator the unsigned horizontal kinetic commutator is
bounded since $|n^2-(n')^2|\le n^2+(n')^2+2m^2$. Its products,
adjoints, the Gram series and the low completion inherit this bound;
both fixed projections commute with $D_x^2$. The vertical scalar
$m^2$ commutes within a copy. This proves the assertion also for
the indicated angle and action derivatives, with the latter's
separate fixed floor prices. No smallness is asserted for the
bounded part of this low generator in slow time.
Explicitly, let \(K_h=T_h-\iota\), where \(\iota\) is the inclusion of the
high inputs. The column of \(K_h\) at a high input \(n\) is \(e^{inx}\) times
the correction of the envelope, of size \(E_c/|n|\) with arbitrary decay
in the spatial shift; so the entry of \(K_h\) from \(n\) to \(n+j\) is at most
\(E_c|n|^{-1}(1+|j|)^{-M}\) for any fixed \(M\), with the prices of the column
construction. The kinetic commutator multiplies this entry by
\((n+j)^2-n^2=2nj+j^2\), so the corresponding entry of \([D_x^2,K_h]\) is at
most \(E_c(2|j|+j^2/|n|)(1+|j|)^{-M}\), whose Schur sums are \(O(E_c)\) for
\(M>3\). The same bound holds for adjoints, for products with matrices
bounded in the translation algebra, for \(G_h^{-1/2}-1\) by its convergent
series, and for \(U_l-Q_l\), because \(D_x^2\) commutes with \(Q_l\) and with
\(m^2\). Since \(U\) is unitary,
\(U^*(D_x^2+m^2)U-(D_x^2+m^2)=U^*[D_x^2,U-1]\) is bounded by \(E_c\) in the
translation algebra. Hence the compressed low block is the free
kinetic term plus a bounded Hermitian operator.

\subsubsection{Evaluation of the diagonal at the final frequencies}
Let $\Omega^0$ be the frozen flat clock, $\Omega_*$ the actual tilted
terminal rates on $F$, and $H_e$ its private row labels. In the retained
angular range $|\ell|_1\le C_J B$ nonflat integer combinations have
distance at least $\lambda^{C_s}/E_c$ from zero at the frozen vector.
The change on $F$ in passing to the physical rates is $O(N^{-2}E_c)$ and the final
rounding is smaller, hence these lower bounds, with a factor $1/2$,
persist on this range. A nonzero frozen-flat label that is neutral
and spatially constant is $\ell=\sum b_e H_e$, $\sum b_e=0$. For
bounded weak length its coefficients have bounded $\sum|b_e|$
by their disjoint private supports. In the convention where the
signed unperturbed time key is $1$, the fitted rates satisfy
\[
 H_e\cdot\Omega_*=1+N^{-2}\mathfrak d_e+
 O(E_c\lambda^{-C_s}N^{-3}+E_c r^2+\eta_{\mathrm{rnd}}) .
\]
Here $\mathfrak d_e$ is the signed sum of
$(\mathcal Q-M^2)/v_x^2-\mathcal Q/(2|v|^2)$ on the two inputs and
output of that triangle, in root-normalized unshifted directions $v$.
The anchor $e_1$, higher powers in the spatial expansion, and the
terminal-axis action corrections are included in the displayed error;
the actual finite fit uses the full means, not just (T1).
The verified geometry and spectator separation required here are
$|\sum b_e\mathfrak d_e|\ge E_c^{-1}$ on these nonzero tested words.
The tilt and rounding errors above are $o(N^{-2}E_c^{-1})$ after
their bounded summation, so $|\ell\cdot\Omega_*|\ge
\frac12N^{-2}E_c^{-1}$. No such lower bound for all weak lengths
is used.

The real higher-order scalar satisfies $|d_{\ge3}|\le E_c |n|^{-3}$.
Its angular harmonics are neutral and have zero horizontal shift by
covariance of the recurrence. Separate its full angle
mean. Divide its remaining coefficients in the indicated temporal
range and weak length $\le m_*$ by $i\ell\cdot\Omega_*$, always at
this one fixed vector, obtaining a real phase $b_3$. Omitted
labels have sufficiently small \emph{undivided} sum: by the two
strip widths the two tails are bounded by
\[ E_c\big(e^{-a C_J B}+e^{-\epsilon Bm_*/7}\big) \]
on a smaller fixed strip,
with $a>0$. Choose $C_J,m_*$ for the desired fixed target power,
also paying any specified fixed derivatives. Then
\[
 \|\partial_\phi^j b_3\|\le E_c N^2 H_x^{-3}\le E_c' N^{-1},
 \qquad
 \|\partial_I^k\partial_\phi^j b_3\|
 \le E_c' \lambda^{-C_s k}N^{-1}.
\]
For the first bound, the coefficients of $d_{\ge3}$ are at most
$E_cH_x^{-3}$ on high inputs and the divisors are at least
$\tfrac12N^{-2}E_c^{-1}$; with $H_x=c_0R/H_{\rm geom}$,
$H_{\rm geom}=C(1+d)P$ and $N=RP$,
\[
 N^2H_x^{-3}=c_0^{-3}C^3(1+d)^3P^6N^{-1},
\]
and the factor $(1+d)^3P^6$ belongs to $E_c'$.
In either translation norm the entire spatial frame, including
the two high scalar phases, differs from the identity by $E_c'/R$,
also on free low inputs. Indeed
$\|T_h-\mathrm{inclusion}\|\le E_c(H_x^{-1}+R^{-2})$ implies
the same estimate for $G_h-1$, $\Pi-Q_l$ and $U_l-Q_l$ by the
convergent square-root formulas. The two high scalar phases add
at most $E_c'N^{-1}$, and $R/H_x=O(dP)$ fits the stated price.
In particular every extra mixed cubic obtained by flipping the
copy of a transverse argument with this transformation costs at
most $E_c'/(R\lambda)$ in normalized field units.

Right-multiply the completed frame by $e^{-i(b_2+b_3)}$ on the
high inputs, with the identity on the low inputs. The low generator
is untouched. The high angle-dependent error under ${\mathcal{D}}$ is
the sum of the spectral/angle tail and
$-(\Omega(I,0)-\Omega_*)\cdot\partial_\phi b_3$.
The high scalar mean itself and its action derivatives have size
$E_c H_x^{-2}$ and $E_c\lambda^{-C_s k}H_x^{-2}$ after its free term
has been removed. In particular differentiating the mean requires
no temporal division. For truly moving coordinates one must add the
chain-rule connection from $\dot I$ and
$\dot\phi-\Omega(I,0)$ on the entire frame, including all high phases.
The denominators in this differentiation stay fixed.

For use in real Hamiltonian coordinates write the finished symplectic
frame as $T$, the constant complex structure as ${\mathcal{J}}$ with
${\mathcal{J}}^2=-1$, and the real untransformed Hessian as $B_{\mathrm{old}}$.
Define $C=T^*B_{\mathrm{old}}T+T^*{\mathcal{J}}{\mathcal{D}}T$ and
$G_a=T^*{\mathcal{J}}\partial_{\phi_a}T$. Differentiating along an angle,
\[
 \partial_{\phi_a} C=T^*(\partial_{\phi_a}B_{\mathrm{old}})T+
 {\mathcal{D}}G_a+({\mathcal{J}}C)^*G_a+G_a{\mathcal{J}}C .
 \tag{T2}
\]
Indeed $T^{-1}\partial_{\phi_a}T=-{\mathcal{J}}G_a$, $G_a$ is symmetric,
and $\partial_{\phi_a}$ commutes with ${\mathcal{D}}$ at fixed actions.
In detail, differentiating $T^*{\mathcal{J}}T={\mathcal{J}}$ along
$\phi_a$ and along the flow gives $G_a^*=G_a$ and the symmetry of
$T^*{\mathcal{J}}{\mathcal{D}}T$, hence of $C$; moreover
$T^{-1}=-{\mathcal{J}}T^*{\mathcal{J}}$. Thus
$\partial_{\phi_a}T=-T{\mathcal{J}}G_a$ and
$(\partial_{\phi_a}T)^*=G_a{\mathcal{J}}T^*$, and the four terms of
$\partial_{\phi_a}C$ in which $T$ or $T^*$ is differentiated are
\[
\begin{aligned}
 &G_a{\mathcal{J}}T^*B_{\mathrm{old}}T-T^*B_{\mathrm{old}}T{\mathcal{J}}G_a
  +G_a{\mathcal{J}}T^*{\mathcal{J}}{\mathcal{D}}T
  -T^*{\mathcal{J}}({\mathcal{D}}T){\mathcal{J}}G_a
  -T^*{\mathcal{J}}T{\mathcal{J}}\,{\mathcal{D}}G_a\\
 &\qquad=G_a{\mathcal{J}}C-C{\mathcal{J}}G_a+{\mathcal{D}}G_a ,
\end{aligned}
\]
because $T^*{\mathcal{J}}T{\mathcal{J}}=-1$. Finally
$({\mathcal{J}}C)^*=-C{\mathcal{J}}$ by the symmetry of $C$.
For an action derivative the extra commutator
$(\partial_{I_a}\Omega)\cdot\partial_\phi$ must also be included.
The angle derivative of the high block of $C$ has only the residual
and frequency-difference prices just derived; an angle derivative
of the compressed low bounded part can have size $E_c$.

To obtain a physical residual $N^{-J}$, one admissible choice is
$M_0=\lceil2J+12\rceil$, $j_0=\lceil J+6\rceil$, and angular tails
initially $N^{-2J-8}$. For small enough macro-height exponent,
$H_x,R\ge N^{3/4}$. A spatial derivative reservoir of
$2(M_0+j_0+6)$ and a nominal time derivative reservoir of
$M_0+j_0+6$ beyond the requested final derivatives suffice in these
finite operations. First choose $\epsilon$ for the fixed action
derivative losses, then $Q$ large and the depth coefficient small
to pay the finitely many $E_c$ factors. All kinetic estimates hold
uniformly in $|q|\ge1$. Lemma~\ref{rate:lem:trans-axial-expansion} used
only a large single horizontal carrier and the true axial Birkhoff
differential, independently of this vertical gap.

\subsection{Estimates in the two-dimensional variables}
Throughout a stage put $N=RP_{\mathrm{lat}}$, $r=HN^{-s}$,
$\lambda=e^{-\epsilon B}$ and $u=r\lambda$; the slow time is
$\tau=ut$. Use $\mathcal E$ for any one of the fixed products
of prices $\exp(C_{s,\epsilon}[B_o+d^2\log(2+B)]+o(B))$
involved below and $\mathcal P$ when also a fixed power
$\lambda^{-C_s}$ may be included. A coefficient multiplied by a
small power of $\lambda$ will first be absorbed using only
$\mathcal E$. In these notations the cutoff in the number of
weak letters can enter $C_{s,\epsilon}$, but does not enter
$C_s$ in an uncompensated power of $\lambda^{-1}$. The latter
powers arise from a fixed number of coordinate derivatives,
the change of time and the noninteger-frequency gap. All spatial
derivative and normal-form orders are fixed before $\epsilon$;
they can be taken larger than the specified cap powers as necessary.

\subsubsection{A Fourier algebra on the full lattice}
Write $k=(n,Rq)$ on the lattice with fixed vertical spacing.
The charge-one coset inside it is invariant for the cubic flow:
$k_1-k_2+k_3$ stays in the coset. Use the ambient lattice for
products of other charges. Let $\ell(k)$ be the smallest number
of signed generators summing to $k$, with value $+\infty$
outside their subgroup. Take as generators $e_1=(1,0)$, the axial sites $(m,0)$ for
$m\in F$ (including the current clock), and every physical marked head. On the signed second copy reflect
$k$ if using the same horizontal carrier convention. Put
\[
 b(k)=1+|n|/N+|q|/P_{\mathrm{lat}},\quad
 W_A(k)=\min(N^A,e^{\ell(k)/B^{K_0}}),\quad
 X_{p,A}=\ell^2(b^p W_A).
\]
Here $p>2\lceil s\rceil+10$, $A>2s+12$, $A'>A+3s+15$
are fixed, $K_0>5$ is fixed, and the geometric word test of
Lemma~\ref{rate:l:denom} is used through $L=B^{K_0+3}$.
These inequalities are the requirements of the present subsection.
The values fixed once and for all in
Section~\ref{rate:sec:parameters}, namely $A>4s+20$, $A'=A+3s+20$
and $J_0>4(A'+s+p+10)$, satisfy them as well as the requirements
$A>2s+12$, $A'=A+3s+20$, $J_0>2s+12$ of
Section~\ref{rate:sec:ward-landing}, and these values are understood
throughout.
Both $W$ and $b$ are submultiplicative; the polynomial
alternative $b(k+k')^p\le 2^p(b(k)^p+b(k')^p)$ will also be
used. Indeed $1+x+y\le(1+x)(1+y)$ and
$1+x+y\le(1+x)+(1+y)$ for $x,y\ge0$, the word length is subadditive,
and $\min(N^A,\xi\eta)\le\min(N^A,\xi)\min(N^A,\eta)$ for
$\xi,\eta\ge1$. The generator count satisfies $g\le C_s(1+d^2+|F_{\mathrm{old}}|)$,
where $F_{\mathrm{old}}$ is the inherited support and the template
has fixed length. By (Pr2) this gives $g\le C'_s(1+d^2)$.
A sum of length $\le C B^{K_0+1}$ is determined by the signed
multiplicity of each generator, so the number of such positions is
at most $(2C B^{K_0+1}+1)^g$. Its logarithm is
$O_{s,K_0}((1+d^2)\log(2+B))$, as required for the bound $\mathcal E$.

For either cap, if at least one input lies in these unsaturated
positions, put that input in $\ell^1$ after weighting using this
count, and use submultiplicativity for the output. If both are
saturated, the output word weight divided by the product of the two
input word weights is at most $N^{-A}$ in the convolution bound. On ordinary $\ell^2(b^p)$
the product costs $O_p(\sqrt{NP_{\mathrm{lat}}})$, by the polynomial
alternative and Cauchy--Schwarz with
$\sum b^{-2p}\le C_p NP_{\mathrm{lat}}$. Thus $X_{p,A}$ is a
convolution algebra with price $\mathcal E$. Its embedding into
ordinary $\ell^1(b^q)$, $q\le p-2$, also costs $\mathcal E$:
use the count on the unsaturated part, and pay the same square
root by $N^{-A}$ on the rest. A signed partial shift by one
angular label of the old torus and boundedly many heads has
word-weight ratio at most
$\exp[(|l|_1+C)/B^{K_0}]$ in either direction. Its polynomial
weight ratio costs a power of
$1+d+(|l|_1+1)Z/N$, where $Z=1+\max_{m\in F}|m|$.
These estimates hold also on the stronger cap. The short-word
faithfulness shows that every nonaxial vector with
$|n|<H_x$ is saturated at both caps, since an
old shift of $O(B^{K_0+3}Z)=o(R)$ cannot cancel a nonzero
integer horizontal macro displacement.
In formulas, let \(\Lambda\) be the set of unsaturated positions, where
\(e^{\ell(k)/B^{K_0}}<N^A\); by the count above \(\#\Lambda\le\mathcal E\). Split
\(f=f_\Lambda+f_{\Lambda^c}\), \(h=h_\Lambda+h_{\Lambda^c}\), and write \(\omega=b^pW_A\).
Submultiplicativity \(\omega(k_1+k_2)\le\omega(k_1)\omega(k_2)\) gives
\[
 \|f_\Lambda*h\|_{X_{p,A}}\le\|\omega f_\Lambda\|_{\ell^1}\|h\|_{X_{p,A}}
 \le(\#\Lambda)^{1/2}\|f\|_{X_{p,A}}\|h\|_{X_{p,A}},
\]
and likewise with the roles exchanged. If both inputs are saturated,
\(W_A(k_1+k_2)\le N^A=N^{-A}W_A(k_1)W_A(k_2)\), and the polynomial
alternative with Cauchy--Schwarz gives
\[
\begin{aligned}
 \|f_{\Lambda^c}*h_{\Lambda^c}\|_{X_{p,A}}
 &\le N^A2^p\bigl(\|b^pf_{\Lambda^c}\|_{\ell^2}\|h_{\Lambda^c}\|_{\ell^1}
     +\|f_{\Lambda^c}\|_{\ell^1}\|b^ph_{\Lambda^c}\|_{\ell^2}\bigr)\\
 &\le C_p\sqrt{NP_{\mathrm{lat}}}\,N^{-A}\|f\|_{X_{p,A}}\|h\|_{X_{p,A}},
\end{aligned}
\]
because \(\|g\|_{\ell^1}\le(\sum b^{-2p})^{1/2}\|b^pg\|_{\ell^2}\) and
\(\|b^pg_{\Lambda^c}\|_{\ell^2}=N^{-A}\|g_{\Lambda^c}\|_{X_{p,A}}\). As
\(\sqrt{NP_{\mathrm{lat}}}\,N^{-A}\le1\), the product constant is at most
\(C_p\mathcal E^{1/2}\). For the embedding,
\(\|b^qf_\Lambda\|_{\ell^1}\le(\#\Lambda)^{1/2}\|f\|_{X_{p,A}}\) and, for \(q\le p-2\),
\(\|b^qf_{\Lambda^c}\|_{\ell^1}\le N^{-A}(\sum b^{2(q-p)})^{1/2}\|f\|_{X_{p,A}}
\le C\sqrt{NP_{\mathrm{lat}}}\,N^{-A}\|f\|_{X_{p,A}}\).

Here is the version for analytic Cartesian axial coordinates.
The inverse Birkhoff map and its inverse-coordinate counterpart
are analytic on small complex $h^0$ balls with the Fourier
linear normalization fixed above. For a derivative of order $m$
divided by $m!$ about the old torus, polarization of the Cauchy
estimate bounds its value on unit coordinate vectors in $L^2$
by $C^m$. On an old angle strip, a Fourier coefficient in these
coordinates is supported at the single spatial frequency
specified by the signed inputs plus $l\cdot n_F$, by translation
equivariance. Cauchy in the angles supplies $e^{-a_0|l|_1}$
on a smaller strip. Taking absolute sums costs $C^{|F|}$;
fixed translation moments are paid by the preceding polynomial
weight. The same holds for the reverse signed kernel (or the
canonical adjoint) by exchanging an output and an input in the
momentum condition. Consequently the absolute kernel on $m$
copies of $X_{p,A}$ has norm at most $C^m\mathcal E^{m+1}$;
here the extra fixed power incorporates the linear derivative
and can be taken once outside the convergent Taylor sum.
In detail, write the \(m\)th Taylor term, divided by \(m!\), as
\(T_m[\zeta_1,\ldots,\zeta_m](k)=\sum_l\sum t_l(k;n_1,\ldots,n_m)
\prod_i\zeta_i(n_i)\,e^{-il\cdot\phi}\), the inner sum running over
\(k=\sum_i\sigma_in_i+l\cdot n_F\), with \(\sigma_i=\pm1\) the signs of
the copies. The Cauchy estimate on the chart ball of fixed radius,
polarization (whose constant \(m^m/m!\le e^m\)) and Cauchy in the
angles give \(|t_l|\le C^me^{-a_0|l|_1/2}\). For fixed \(l\) the sum over
the \(n_i\) is an \(m\)-fold convolution of the moduli
\(|\zeta_i|\), translated by \(l\cdot n_F\). Applying the algebra
estimate \(m-1\) times gives
\(\||\zeta_1|*\cdots*|\zeta_m|\|_{X_{p,A}}\le(C_p\mathcal E^{1/2})^{m-1}
\prod_i\|\zeta_i\|_{X_{p,A}}\), and the translation by \(l\cdot n_F\)
multiplies the weight by at most
\(C(1+(|l|_1+1)Z/N)^pe^{(|l|_1+C)/B^{K_0}}\), as shown above. Summing
over \(l\) with the factor \(e^{-a_0|l|_1/2}\) costs at most
\(C^f(1+Z/N)^p\le\mathcal E\) by (Pr5). Hence the absolute kernel maps
\(X_{p,A}^m\) to \(X_{p,A}\) with norm at most
\(C^m\mathcal E(C_p\mathcal E^{1/2})^{m-1}\le C^m\mathcal E^{m+1}\). With
inputs of size \(r\|Y\|\) the Taylor series converges absolutely as
soon as \(Cr\mathcal E\|Y\|\le1/2\).
Each successive input is accompanied by the physical factor $r$,
and $r\mathcal E^C\ll1$. For its derivative on one free argument
in ordinary $\ell^2(b^q)$ put all other physical arguments in
$\ell^1(b^q)$ using $0\le q\le p-2$. General symmetry or bounded
analyticity alone on an unweighted Hilbert space would not imply
this convolution estimate; the single signed output is essential.
For a derivative on occupied actions, differentiating the old
Cartesian radii may cost the inverse positive floors; for any
fixed derivative order this costs $\mathcal P$. Current high hot
axes have Cartesian, not polar, coordinates.

\subsubsection{The projected reference and its actual coefficients}
In the autonomous mass gauge subtract $M^2$ from the Hamiltonian.
Let $z_F$ have occupied polar coordinates $(I,\phi)$ and vanish
off $F$, let $a$ be supported on $F^c$, and write $w$ for the
physical nonaxial part. Put $z=z_F+a$, $V=\Psi(z)$, $V_F=\Psi(z_F)$,
$M_a=\|z\|_{\ell^2}^2=\|V\|_{L^2}^2$ and $M_w=\|w\|_{L^2}^2$.
The exact pure axial Hamiltonian is the one-dimensional action
Hamiltonian in the same gauge. Apart from the transverse kinetic
term, its mixed extension is
\[
 2\langle(|V|^2-M_a)|w|^2\rangle
 +\operatorname{Re}\langle V^2\bar w^2\rangle
 +2\operatorname{Re}\langle\bar V w|w|^2\rangle
 +\tfrac12\langle|w|^4\rangle-M_w^2,
\]
where the brackets denote spatial integral. Vertical momentum
eliminates every term linear in $w$. In detail, the quartic part of
the Hamiltonian is $\tfrac12\langle|V+w|^4\rangle-(M_a+M_w)^2$, and
pointwise, using $2(\operatorname{Re}\zeta)^2=|\zeta|^2+\operatorname{Re}\zeta^2$
with $\zeta=\bar Vw$,
\[
 \tfrac12|V+w|^4=\tfrac12|V|^4+2|V|^2\operatorname{Re}(\bar Vw)
   +2|V|^2|w|^2+\operatorname{Re}(\bar V^2w^2)
   +2\operatorname{Re}(\bar Vw)|w|^2+\tfrac12|w|^4 .
\]
The second term has zero integral because $V$ is independent of $y$
and $w$ has no zero vertical mode, and
$-2M_aM_w=-2M_a\langle|w|^2\rangle$ combines with the third.
In particular $D_a M_a=0$
at the empty Cartesian origin. At $a=0$ remove the pure quadratic
transverse part by the spatial construction. The leading remaining
mixed Hamiltonian, homogeneous of field degree three, is
\begin{align*}
 &4\langle\operatorname{Re}(\bar V_FD\Psi(z_F)a)|w|^2\rangle
 +2\operatorname{Re}\langle V_F D\Psi(z_F)a\,\bar w^2\rangle
 +2\operatorname{Re}\langle\bar V_F w|w|^2\rangle .
\end{align*}
There are at least two transverse factors in every term. The
other terms start at degree four. After substituting the spatial
frame its derivative, canonical signed adjoints and small axial
coordinate corrections obey the preceding kernel estimates.
In particular at the terminal configuration with $w=0$ the
mixed field vanishes identically. The rank-one mass factors
are bounded by Cauchy--Schwarz on their scalar input, putting
a positive spatial weight on their known output when the free
variable is in the scalar. The exact pure axial Cartesian field
is diagonal, with coefficient $n^2+\chi_n(I,|a|^2)$; its
nonlinear part has bounded first and second action derivatives
as a map on $\ell^2$. Multiplying a frequency derivative by
a known vector $x_n$ in $\ell^2(b^q)$ costs at most
$\|b^q x\|_\infty\|D\chi\|_{2\to2}\|\delta I\|_2$.
No weighted mapping theorem for $D\chi$ is needed.

For the finite reference keep the full projection of this leading
degree-three Hamiltonian on the marked finite planes, in the
Fourier part of the mixed coefficients with $I=I^*$, the rounded
terminal action vector restricted to $F$, and angles advancing at their accepted terminal rates. Keep the actual pure
axial action detuning on the marked axes. Retain the exact real
constant spatial means on high transverse carriers to the
prescribed order and the true Birkhoff frequencies on terminal
axes. Corrections from the high spatial frame and $D\Psi$ on
the actual marked inputs are $\mathcal E/R$ in the projected
coefficients and the corresponding two-known sources, as the
marked nonzero horizontal macro is an integer. The general
unmarked homogeneous axial approximation above $N^{1/4}$ costs
only $\mathcal E N^{-1/4}$ and is used solely where applicable.
The full one-dimensional projection in low empty modes is
treated separately in the exterior estimate. Every assigned
triad of distinct heads has coefficient $2\overline{v_e}$ on
$z_v z_w\bar z_h$, with the charge-one old Fourier harmonic
$v_e$ specified in the calibration. This factor follows by
permuting the two positive and the two negative entries of the
physical quartic with prefactor $1/2$. The same coefficient
results from each of the three displayed mixed terms in its
appropriate axial/nonaxial case.
In detail, the quartic is
\(\tfrac12\sum_{k_1-k_2+k_3-k_4=0}\widehat u(k_1)\overline{\widehat u(k_2)}
\widehat u(k_3)\overline{\widehat u(k_4)}\). The monomial
\(z_vz_w\bar z_h\overline{v_e}\), with \(v\ne w\), arises from the two
orders of \(k_v,k_w\) in the slots \(k_1,k_3\) and the two orders of
\(k_h\) and the index \(b\) of the old harmonic in the slots
\(k_2,k_4\); the four ordered choices each carry \(\tfrac12\). In the
three displayed mixed terms, with \(a'=D\Psi(z_F)a\) replaced by its
Fourier part for the leading coefficient: if \(v,w,h\) are transverse,
\(2\operatorname{Re}\langle\bar V_Fw|w|^2\rangle\) contains
\(\langle\bar V_Fw^2\bar w\rangle\), and \(w^2\bar w\) contains
\(2w_vw_w\bar w_h\); if the input \(v\) is axial,
\(4\langle\operatorname{Re}(\bar V_Fa')|w|^2\rangle\) contains
\(2\langle\bar V_Fa'w\bar w\rangle\), and \(w\bar w\) contains
\(w_w\bar w_h\) once; if the output \(h\) is axial,
\(2\operatorname{Re}\langle V_Fa'\bar w^2\rangle\) contains
\(\langle\bar V_F\bar a'w^2\rangle\), and \(w^2\) contains \(2w_vw_w\).
In each case the coefficient of \(z_vz_w\bar z_h\) is \(2\) times the
component \(\overline{v_e}\) of \(\bar V_F\).
Fit $2\overline{v_e}/\lambda$
to the imaginary coefficient of Lemma~\ref{rate:l:finite} using one
new angle per row. The target accuracies in the amplitude are
absolute relative to the prescribed row strengths, whose
logarithms cost $O(D+d\log d)$.

In the common mass gauge, at the rounded terminal actions let $\nu_j$ be the true axial Birkhoff
frequency on an axial head and the nominal carrier frequency including
the retained constant diagonal corrections on a transverse head.
For the edge $e=(v,w;h)$ put $\Lambda_e=H_e\cdot\Omega_*$ and
$d_e=\nu_v+\nu_w-\nu_h-\Lambda_e$. The old harmonic carries the
factor $e^{-i\Lambda_e t}$. Its conjugate in the Hamiltonian term
$2\overline{v_e}z_vz_w\bar z_h$ therefore contributes $e^{+i\Lambda_e t}$. After the true physical fit in the axial parameters
one can insist on
$d_0=\max_e|d_e|\le u2^{-d}e^{-C_4dD}$. Move only the
nonaxial rotating frequencies, writing $z_j=e^{-i(\nu_j+x_j)t}y_j$
with $x_j=0$ on every axial head. The phase of the displayed term is
$e^{-it d_e^{\rm rot}}$, where
\[
 d_e^{\rm rot}=(\nu_v+x_v)+(\nu_w+x_w)-(\nu_h+x_h)-\Lambda_e
              =d_e+x_v+x_w-x_h .
\]
Thus the choice $x_h=x_v+x_w+d_e$ gives $d_e^{\rm rot}=0$ on
each primary edge, and the conjugate monomial has the opposite, also
zero, mismatch. At the fixed top choose shifts on the two top children
whose sum is $-d_{\mathrm{top}}$ and whose magnitudes are at most
$d_0$. In each branch give all its leaves the same shift
$\xi$, then recurse upward. A head $m$ levels above the leaves
has shift $2^m\xi$ plus a binomial sum of mismatches with
coefficient sum at most $2^m-1$. The branch root condition has
a solution $|\xi|\le d_0$, and all intermediate shifts are
$\le3\cdot2^d d_0$. Indeed every primary parent is the output of
exactly one designated triple, so the relation $x_h=x_v+x_w+d_e$
defines its shift from those of its two children, and by induction
on $m$ the shift is $2^m\xi+\sigma_m$ with
$|\sigma_m|\le2|\sigma_{m-1}|+d_0\le(2^m-1)d_0$. At the root of a
branch, $m=d-1$ and the prescribed shift $x_*$ has $|x_*|\le d_0$, so
$\xi=2^{-(d-1)}(x_*-\sigma_{d-1})$ has modulus at most $d_0$, and
$|2^m\xi+\sigma_m|\le(2^{m+1}-1)d_0$ at every level. The degree-one private companion absorbs
the remaining mismatch in its own row. For $e=(v,J(v);A(v))$, put
$x_{J(v)}=-d_e-x_v$. Since $x_{A(v)}=0$, this gives
$d_e^{\rm rot}=d_e+x_v+(-d_e-x_v)-0=0$.
The diagonal perturbation in slow units is thus $O(e^{-C_4dD})$.
The designated monomials thus have no explicit time phase in these rotations.
The action-dependent detuning stays in the nonlinear vector field, and every
hot axis uses its actual terminal frequency; all old phases are
still advanced in the original physical time.

In slow time write the finite vector field as $F_0+S+G$ on the
Euclidean tube of radius one about the reference and on segments
to its prescribed endpoint, contained in a ball $e^{Cd}$. The
small part $S$ contains the fitting errors just described, the
true hot action detuning of $C^1$ norm $\mathcal E r/\lambda$,
and discarded angular/weak tails when already small. (For the
detuning: a hot axis rotates at its terminal frequency, and its true
frequency differs from it by at most \(C\|I-I^\infty\|_1\). In the tube
the hot actions differ from their terminal values by at most
\(r^2e^{Cd}\) times the normalized deviation, as
\(|z|^2-|z^\infty|^2=(|z|-|z^\infty|)(|z|+|z^\infty|)\). The old occupied actions are fixed at \(I^*\) in the
finite reference system. Multiplying by the
normalized hot amplitude and dividing by the slow rate \(u=r\lambda\)
gives a field of \(C^1\) norm \(Cr^2e^{Cd}/u\le\mathcal Er/\lambda\).) The other
time-dependent terms have a primitive $P$ at fixed state. For
a noninteger class the original physical old coefficient is
bounded by $\mathcal E$ and its inverse gap by $\mathcal P$;
for each nearly resonant wrong charged row its summed physical
old coefficient is bounded by $\lambda\mathcal E$ and inverse
gap by $N^2\mathcal E$. The very long words in the projected
system may be kept as a small $C^1$ part. The two gap and tail
estimates are justified by the axial choice and the external
classifications. Taking two finite state derivatives still costs
only $\mathcal E$. Hence in the slow equation
\[
 \|P\|_{C^2}\le\mathcal E r\lambda N^2+\mathcal P r,
 \qquad \|G\|_{C^2}\le\mathcal E/\lambda .
\]
The explicit substitution $z=y+P(\tau,y)$ gives the exact
field $(1+D P)^{-1}[F_0(y+P)+S(\tau,y+P)
+G(\tau,y+P)-G(\tau,y)]$. The error relative to $F_0$ is
bounded in $C^1$ by the small $S$ plus
$\mathcal E rN^2+\mathcal P r$; in particular there is no
division by $\lambda$ on the former price. Every $G$ and $P$
vanishes at the all-axial point, by vertical momentum. The
true pure axial detuning also vanishes there. Undo the finite
change after applying Lemma~\ref{rate:l:round} at $\sigma=\kappa/3$. Its landing then supplies the two strictly smaller weights $\kappa/4$ and $\kappa/8$ for the full equation. The conjugacy and
its differential have uniformly bounded norm. This gives the
true finite reference and its full complex Green in the original
rotations. In the latter its slow envelope derivative, after
subtracting the endpoint on the decaying tail, can cost
$\mathcal E/\lambda$. In differentiating each covariant
particular or homogeneous primitive below pay this once; never
differentiate a fast physical phase separately. The finite
reference and its derivative may be used on the entire interval
from its cold launch to its prescribed hot endpoint.

To see the two bounds, note that with amplitudes normalized by $r$ and
slow time $\tau=r\lambda t$, a physical cubic term with old coefficient
$c$ contributes $c/\lambda$ times a quadratic monomial to the slow field,
and if it oscillates with physical frequency $\omega$, its slow frequency
is $\omega/(r\lambda)$. Its primitive at fixed state is therefore
$(c/\lambda)(r\lambda/\omega)=cr/\omega$, and two state derivatives on the
tube cost $\mathcal E$. Summing the angular coefficients, a noninteger
class contributes at most $\mathcal E\mathcal Pr$, which is of the form
$\mathcal Pr$, and a nearly resonant wrong charged row at most
$\lambda\mathcal E\cdot N^2\mathcal E\cdot r$, which is of the form
$\mathcal Er\lambda N^2$; the undivided field has size at most
$c/\lambda\le\mathcal E/\lambda$. For the substitution, $\partial_\tau P=G$ at
fixed state gives
\[
 y'-F_0(y)=(1+DP)^{-1}\bigl[F_0(y+P)-F_0(y)-DP\,F_0(y)
   +S(\tau,y+P)+G(\tau,y+P)-G(\tau,y)\bigr],
\]
whose $C^1$ norm is at most $C\|S\|_{C^1}+C\mathcal E\|P\|_{C^2}(1+\|G\|_{C^2})$,
because $F_0$ and its first two derivatives are bounded by
$e^{CD}\le\mathcal E$ on the tube. Inserting the two bounds,
$\|P\|_{C^2}\|G\|_{C^2}\le\mathcal E^2rN^2+\mathcal P\mathcal Er/\lambda$; the
factor $\lambda$ in the first term of $\|P\|_{C^2}$ cancels the division by
$\lambda$ in $\|G\|_{C^2}$, and the second product is again of the form
$\mathcal Pr$, since $\mathcal P$ may contain fixed powers of
$\lambda^{-1}$. This gives the stated error $\mathcal ErN^2+\mathcal Pr$.
We record where the two gap bounds come from. A time-dependent term
of the projected finite field comes from a triad of marked heads and
an old angular harmonic of charge one, or from their conjugates.
For a retained label \(|\ell|_1\le C_JB\), the spatial frequency
\(b=\ell\cdot n_F\) satisfies \(|b|\le C_JBZ=o(R)\). Write the physical
heads as \(k_v=e_1+Rv\), \(k_w=e_1+Rw\), \(k_h=e_1+Rh\) with integer
vectors \(v,w,h\) from the cleared configuration. Momentum conservation
\(k_v+k_w-k_h=(b,0)\) gives \(b=1\) and \(h=v+w\), since
\(|b-1|<R\). The kinetic key is therefore
\[
 |k_v|^2+|k_w|^2-|k_h|^2=1-2R^2 v\cdot w .
\]
If \(v\cdot w\ne0\), its absolute value is at least \(2R^2-1\).
The angular frequency is \(O(B(Z^2+1))=o(R^2)\), and the marked
diagonal and rotation corrections are \(o(1)\). Thus the physical
mismatch has modulus at least \(R^2\) for large \(R\), and the
primitive contributes at most \(\mathcal E r/R^2\).
If \(v\cdot w=0\), the reverse rectangle test of
Lemmas~\ref{rate:l:geom} and~\ref{rate:l:denom} identifies the
marked triad with a designated edge \(e\). Its kinetic key is \(1\).
In the rotating frame its physical frequency is \(\ell\cdot\Omega_*\)
minus this key and the calibrated head corrections, which are
\(N^{-2}\mathfrak d_e\) up to the errors of the terminal fit.
Conjugation reverses all frequency signs. (i) If \(\ell\) is not in the frozen module, the
affine selection on the cube of halfwidth \(\delta_0\) gives a divisor of
modulus at least \(w_0=\lambda^2e^{-C_9\Gamma}\), and the tilt, the rounding
and the head corrections cost \(o(w_0)\); so the inverse gap is at most
\(\mathcal P\). The coefficient is an angular Fourier coefficient of the old
field, at most \(\mathcal E\). (ii) If \(\ell=\sum_{e'}b_{e'}H_{e'}\) is a charged module
label of bounded weak length (\(\sum b_{e'}=1\)) acting on a triad \(e\) and
different from the singleton \(H_e\), its frequency is
\(N^{-2}(\sum b_{e'}\mathfrak d_{e'}-\mathfrak d_e)\) up to the remainder of the terminal
fit. Since \(b-e_e\) is a nonzero neutral word, the spatial choice of the
spectator gives modulus at least \(\tfrac12N^{-2}e^{-CB_o}\), so the inverse
gap is at most \(N^2\mathcal E\). Its coefficient is at most \(\lambda\mathcal E\): the
fitted singletons have modulus \(\lambda|g_{e'}|e^{O(1)}\), and module labels
with at least fifteen weak letters are smaller. (iii) Labels whose weak
length exceeds the cutoff have total coefficient at most
\(\mathcal Ee^{-\epsilon Bm_*/6}\), which after division by \(\lambda\) is a
small \(C^1\) term and is put into \(S\). The singleton \(H_e\) on its own triad
is stationary in the rotating frame and belongs to \(F_0\).

\subsection{The exterior field and boundary interactions}\label{rate:sec:exterior}
Fix \(s>2\). We give the comparison in physical time before passing
to the slow time \(\tau=u t\), where \(u=r\lambda\),
\(r=e^{bd}N^{-s}\), \(N=RP\), and \(\lambda=e^{-\epsilon B}\).
Write \(X\preceq Y\) for a bound \(X\le E^C Y\), where
\[
 E=\exp\{C_{s,\epsilon}(B_o+d^2\log(2+B))+o(B)\}.
\]
Each use below costs a fixed power, independent of the depth. Write
\(L=\lambda^{-C_s}\) for a fixed power sufficient for the
noninteger and differentiated small-angle prices. Fixed powers of
\(L\) are allowed upon decreasing \(\epsilon\) first. The spatial
word cutoff and the weak angular cutoff enter the exponent of \(E\)
only. Let \(F\) denote the old occupied, spectator, moderator and
weak axial sites, excluding current hot axial heads. The field of this torus is \(V_F=\Psi(z_F)\). Its complex angular strip gives
absolute Fourier sums of size \(\preceq1\), including a fixed
number of angle derivatives. We cut at \(|j|_1\le C_J B\) in
every mixed principal matrix and primitive. For larger \(|j|_1\)
the ordinary strip bound, including
\(\sum_j e^{-c|j|_1}\le C^{|F|}\), pays any prescribed
\(N^{-J}\) by increasing \(C_J\). A translated word weight costs
at most \(e^{O(|j|_1/B^{K_0})}\) and is included in this estimate.
In a product of finitely many angular series one may first cut
each factor separately, then enlarge the cutoff for their total
label by this fixed number of factors. Indeed the convolution
with any discarded tail remains small on the narrower strip,
even if large opposing labels have small sum.
Only the compressed pure-old low body, at each fixed vertical
label, will be kept without an angular cut.

We use the following specific inputs from the geometric and
calibration arguments. The integer macros \((X,Y)\) of the
marked heads obey \(X\ne0,\ |X|+|Y|\le H_{\rm geom}\). The
common physical coset is
\[
 (k_x-1)\equiv l k_y\pmod5,\qquad l^2\equiv-1\pmod5,\qquad
 k_y\in R\Z,\quad (R,5)=1 .
\]
For every nonaxial marked macro \(Y\not\equiv0\pmod5\).
A signed macro word of length at most \(B^{K_0+3}\) with zero
horizontal coordinate has zero vertical coordinate. Old axial
translations of this length have size \(o(R)\).
On a retained charged old angular label, with physical space/time
keys \(b,\nu\), we have \(b\equiv1\pmod5\) and
\(|b|^2+|\nu|+1\le R^\eta\), with a small fixed \(\eta\).
Before the final \(O(N^{-2}E)\) tilt, a nonmodule gap of the form
integer plus old frequency has reciprocal \(\preceq L\).
A charge-one label in the frozen module has bare space and time
keys \(1,1\). Its five-weak-letter terms are the designated
singletons, of summed size \(\preceq\lambda\).
All nonsingleton charged module terms have at least fifteen weak
letters. After normalization and summation their contribution to a
homogeneous field row is \(\preceq N^{-h_w}\), where \(h_w>0\) is
the weak angular tail exponent. On a source with two
prescribed inputs one retains every weak length up to the larger
fixed forcing cutoff before using its tail estimate.

\smallskip
\noindent\emph{Where these inputs come from.} (1) \(X\ne0\) and the
coset statements are Lemma~\ref{rate:l:denom} (nonzero horizontal
components, the isotropic line modulo five) together with the choice
\(5\nmid R\); on the chosen line a nonzero vector has both coordinates
nonzero, which gives \(Y\not\equiv0\pmod5\), and
\(|X|+|Y|\le H_{\rm geom}\) is the definition of \(H_{\rm geom}\).
(2) Short-word faithfulness is the word statement of
Lemma~\ref{rate:l:denom} with \(L=B^{K_0+3}\). An old axial
translation by a word of this length has modulus at most
\(B^{K_0+3}Z\le e^{CB_o+C\log B}=o(R)\), since \(\log R=QB_o+O(1)\).
(3) Every occupied axial site is \(\equiv1\pmod5\): the initial sites
\(1,6\), the spectator \(l\equiv1\), the template sites
\(1+Kl_ev_i\) with \(5\mid K\), and the hot axes of earlier stages,
\(1+RP_i\) with \(5\mid P_i\) because cleared axial heads vanish modulo
five by Lemma~\ref{rate:l:denom}. A charge-one label \(\ell\) has spatial
key \(b=\sum\ell_nn\equiv\sum\ell_n=1\pmod5\). With
\(|\ell|_1\le C_JB\), \(|b|\le C_JBZ\) and
\(|\nu|\le C_JB\max_F|\omega_n|\le C_JB(Z^2+C)\), so
\(|b|^2+|\nu|+1\le e^{CB_o+C\log B}\le R^\eta\) for \(Q\) large.
(4) A nonmodule gap is at least \(w_0=\lambda^2e^{-C_9\Gamma}\) by the
affine selection on the cube of halfwidth \(\delta_0\), so its
reciprocal is \(\preceq\lambda^{-2}\le L\).
(5) The module is generated by the rows, each with charge, momentum and
kinetic sum one. A charge-one element \(\sum_eb_eH_e\) has
\(\sum_eb_e=1\), hence momentum and kinetic sum one. The weak supports
of the rows are disjoint, so it has exactly \(5\sum_e|b_e|\) weak
letters. Five letters force \(\sum|b_e|=1\), a single \(b_e=1\); and
since \(\sum|b_e|\equiv\sum b_e\pmod2\), every other charge-one element
has \(\sum|b_e|\ge3\), that is, at least fifteen weak letters. The
fitted singletons have modulus \(\lambda|g_e|e^{O(1)}\), whose sum over
the rows is \(\preceq\lambda\).
(6) A label with \(m_w\) weak letters has coefficient at most
\(\mathcal Ee^{-(\epsilon B/6)m_w}\) by the enlarged weak strip of width
\(\epsilon B/6\). For \(m_w\ge15\), after division by \(\lambda\) for the
field normalization, this is at most
\(\mathcal Ee^{-15\epsilon B/6+\epsilon B}=\mathcal Ee^{-3\epsilon B/2}\le
N^{-h_w}\) with \(h_w=\epsilon/2\), even after the sum over labels,
which costs another factor \(\mathcal E\): indeed
\(\mathcal E^2e^{-3\epsilon B/2}\le e^{4\gamma B-3\epsilon B/2}\le
e^{-\epsilon B}\le N^{-\epsilon/2}\), using \(\mathcal E\le N^\gamma\),
\(\log N\le2B\) and \(\gamma\le\epsilon/8\). (The complexified weak
radii on the enlarged strip are \(\lambda^{1/5}e^{\epsilon B/6+O(\Gamma)}
=e^{-\epsilon B/30+O(\Gamma)}\), still small, which is why the width is
\(\epsilon B/6<\epsilon B/5\).)

For a bare anomalous split of an arbitrary unmarked high leg from
a marked head, the two-component test separates every undesignated
singleton. For a bare split or fusion involving two marked nodes
the test instead separates every undesignated charged label through
the stated larger cutoff. These are the two different nonidentity
statements required of the geometry. The spectator lies outside
the forbidden ratios of the resulting nonzero rational component
pairs. Marked temporal carriers have already been calibrated to
within \(O(u\eta_0)\) of the genuine constant high diagonals,
where \(\eta_0\) is the prescribed finite shadow tolerance. On
the marked coordinates we may put their calibrated velocities
\(\kappa_k\) in the nominal diagonal \(D\) and put the discrepancy
in the finite equation. On the unmarked high coordinates \(D\)
uses the true high diagonal through the fixed order of the spatial
reduction. On the transverse low coordinates its comparison value
is free kinetic. The actual bounded low matrix is left on the
left-hand side. Doubled coordinates have diagonal signs
\(\varepsilon_k=\pm1\).

\subsubsection{The two low density stars}
Set \(h=(h_x,RY)\), \(h_x=1+RX\), and write the input
\((n,Rq)\). The same-copy density shift goes to
\((n+h_x-b,R(q+Y))\), and the opposite-copy pairing goes to
\((h_x+b-n,R(Y-q))\). After subtracting the time of the moving
head, their exact free mismatches are
\begin{align}
 D_F&=2(n-b)(RX+1-b)+2R^2qY+\nu-b^2, \label{rate:ex:df}\\
 D_A&=2(n-b)(n-1-RX)+2R^2q(q-Y)+b^2-\nu. \label{rate:ex:da}
\end{align}
In particular the old harmonic need not have time \(b^2\).
To derive the two formulas, write \(k=(n,Rq)\) for the input and
\(k'\) for the output. The density term multiplies the input by the
moving head and by a conjugated old harmonic of spatial index \(b\)
and time frequency \(\nu\), so that
\(D_F=|k'|^2-|k|^2-|h|^2+\nu\); the pairing term gives
\(D_A=|k'|^2+|k|^2-|h|^2-\nu\). Since
\(|h|^2=h_x^2+R^2Y^2\),
\[
\begin{aligned}
 (n+h_x-b)^2-n^2-h_x^2+2R^2qY
   &=2(n-b)(h_x-b)-b^2+2R^2qY,\\
 (h_x+b-n)^2+n^2-h_x^2+2R^2q(q-Y)
   &=2(n-b)(n-h_x)+b^2+2R^2q(q-Y),
\end{aligned}
\]
and \(h_x=1+RX\).
The true constant high corrections in the comparisons made here
are \(\preceq N^{-2}\), uniformly in \(q\), since the high
horizontal abscissae have magnitude at least \(H_x\) and
\(P,H_{\rm geom}\preceq1\).
Explicitly, at a high input \(|n|\ge H_x\) the constant correction is
\(d_2(n,m)+\bar d_{\ge3}\) with
\(|d_2|\le|\mathcal Q-M^2|n^{-2}+\mathcal Q(2E_n)^{-1}\le Cn^{-2}\) by
(T1), uniformly in \(m\), and \(|\bar d_{\ge3}|\le E_c|n|^{-3}\). Since
\(H_x=c_0R/H_{\rm geom}\) and \(R=N/P\),
\(n^{-2}\le H_{\rm geom}^2P^2c_0^{-2}N^{-2}\), and
\(H_{\rm geom}^2P^2\preceq1\) by (Pr4).

A transverse end is low precisely when \(q\ne0\) and \(|n|<H_x\).
On a low input to
\eqref{rate:ex:df}, a nonzero integer \(qY\) dominates the remaining
terms, giving a gap \(\gtrsim R^2\). Otherwise \(Y=0\).
For \(n\ne b\) the gap is \(\gtrsim R\), and for \(n=b\)
retain the entire finite-angular density star
\[
 (b,q)\ \longleftrightarrow\ (1+RX,q)
\]
in coordinates (physical abscissa, vertical integer).
For a low output put \(\vartheta=n'-(2-b)\). Direct substitution
in \eqref{rate:ex:df} gives
\[
 2R^2(q'Y-Y^2-X^2)+2(RX+1-b)\vartheta
                +2(1-b)^2+\nu-b^2.
\]
Indeed the output is \((n',Rq')=(n+h_x-b,R(q+Y))\), so
\(n-b=\vartheta+(1-b)-RX\) and \(q=q'-Y\), and
\[
 2\bigl(\vartheta+(1-b)-RX\bigr)\bigl(RX+(1-b)\bigr)
 =2(RX+1-b)\vartheta+2(1-b)^2-2R^2X^2 .
\]
If the first parenthesis is nonzero it dominates. If it vanishes,
\(Y\ne0\), \(q'=Y+X^2/Y\) is an integer, and
\(\vartheta\ne0\) gives a gap \(\gtrsim R\). Retain the other
finite-angular star and its reverse,
\[
 (1-RX,X^2/Y)\ \longleftrightarrow\ (2-b,Y+X^2/Y).
\]
All retained star entries have their actual time-dependent
coefficients and Hermitian adjoints. In \eqref{rate:ex:da} the two
legs can be interchanged. At a low leg only \(q=Y\) can remove
the \(R^2\) term. Then \(n=b\) is impossible, since
\((n-1,Rq)=(0,RY)\bmod5\) lies outside the indicated line.
Every remaining gap is \(\gtrsim R\). Thus every mixed primitive
touching the true transverse low matrix has a large divisor.
The common choice \(C_JB\) in the stars is essential to the
word estimate below.

\subsubsection{Bounded components at arbitrarily large vertical indices}
At a zero bare same-copy high gap with a module label, put
\(w=((n-1)/R,q)\). The retained high graph has edges
\(w\leftrightarrow w+v,\ w\cdot v=0\), \(v=(X,Y)\),
after the finite graph and its divided boundary incidences are
removed. The first coordinate of \(w\) can be nonintegral. In a
run of parallel directions every edge connects a common
perpendicular foot to a point on the line through that foot in
the given direction. A simple parallel run thus needs at most two
steps. Each switch between nonparallel directions occurs at one
of the at most four intersections of the two base/translated
line pairs for the specified labels. For \(g\) heads a simple
shortest high-only path has at most \(2(1+4g^2)\) steps, also
after removing edges.
For a marked direction \(v\) put \(L_0(v)=\{w:w\cdot v=0\}\) and
\(L_1(v)=\{w:w\cdot v=|v|^2\}=L_0(v)+v\). An edge in direction \(v\) joins a
point of \(L_0(v)\) to its translate by \(v\), which lies in \(L_1(v)\). If
\(v'=cv\) with real \(c\ne0\), then \(L_0(v')=L_0(v)\) and
\(L_1(v')=\{w\cdot v=c|v|^2\}\). A point \(w+v\in L_1(v)\) therefore lies on
\(L_0(v')\cup L_1(v')\) only if \(c=1\), and then the edge in direction \(v'\)
from \(w+v\) leads back to \(w\). Consequently two consecutive edges of a
simple path in parallel directions share their endpoint in \(L_0(v)\),
and a maximal run of parallel directions has the form
\(w+v\to w\to w+v'\), with at most two steps. At a vertex where a path
changes from \(v\) to a nonparallel direction \(v'\), the vertex lies in
\((L_0(v)\cup L_1(v))\cap(L_0(v')\cup L_1(v'))\), which consists of at most
four points. For \(g\) heads there are at most \(4g^2\) such switching
points, and a simple path visits each at most once, so it has at most
\(4g^2+1\) maximal parallel runs and at most \(2(1+4g^2)\) steps. A shortest
path is simple, and removing edges does not affect the bound.

Join all transverse low coordinates at each fixed \(q\) for the
purpose of forming components. The pure-old off-copy block has
been removed to the prescribed order using
\(n^2+(n')^2+2R^2q^2\ge2R^2\) and the paired symplectic
transformation. In the remaining same-copy compressed low
potential the kinetic term \(R^2q^2\) cancels from all
commutators at fixed \(q\). The bounded low part has operator
norm \(\preceq1\) by the old angular moment estimates. It can
depend on the actual old parameters throughout.

A reverse nonaxial low star has both ordinate magnitudes bounded
by \(2H_{\rm geom}^2\). Hence a low slab with \(|q|>4H_{\rm geom}^2\) can leave
through an axial star only. Its high end has macro abscissa
\(|x|\le H_{\rm geom}\) and the same \(q\). On any nonaxial bare density
edge there one would have
\[
 |qv_y|\le |xv_x|+|v|^2\le2H_{\rm geom}^2,\qquad |v_y|\ge1,
\]
since \(|x|,|v_x|\le H_{\rm geom}\) and
\(|v|^2\le(|X|+|Y|)^2\le H_{\rm geom}^2\).
Axis-only steps keep \(q\) and visit only \(x=0\) or
\(|x|\le H_{\rm geom}\) before another slab. The entire component therefore
has this same \(q\). No component containing a low slab with
\(|q|\le4H_{\rm geom}^2\) can reach a larger low ordinate. On that
component every high vertex has a high-only path of at most
\(2(1+4g^2)\) steps to an incident star. Both its macro
coordinates are bounded by \(C(1+g^2H_{\rm geom}+H_{\rm geom}^2)\), regardless of
the number of attached low slabs. This proves that the ratios of
\(b_k^p=(1+|k_x|/N+|q|/P)^p\) on any such whole component
are \(\preceq1\). On a component disjoint from the low slabs
the same bound follows from the length of a shortest path, even
when all its vertices have large abscissa or ordinate.

For \(W(k)=\min(N^A,e^{\operatorname{length}(k)/B^{K_0}})\),
every transverse low vertex has length \(>B^{K_0+3}\).
Indeed the low horizontal coordinate forces zero macro abscissa
by integrality and the \(o(R)\) old shifts, contradicting
nonzero vertical coordinate by short-word faithfulness. To reach
a low slab from any high vertex in its component, at most
\(O(g^2)\) high-only steps and one star step with at most
\(C_JB+O(1)\) letters suffice. All vertices of such components
are saturated. Elsewhere the ratio of word weights is bounded
by \(\exp(O(g^2)/B^{K_0})\le2\). In particular neither low star
contains a marked high endpoint. In the reverse star, matching
a marked abscissa \(-X\) would force its ordinate to be
\(-Y\) by faithfulness, while the actual ordinate is \(X^2/Y\).

Replace \(W b^p\) by its component infimum, preserving an
equivalent Hilbert norm up to \(\preceq1\). The retained
principal same-copy matrices, including the low body, are
Hermitian with their actual phases. Their propagation is unitary
in this flattened norm. The kinetic diagonal also commutes
with the weight. In the slow-time norm with weight one until
\(T\) and \(e^{\sigma(\tau-T)}\) afterwards, the backward
integral has norm \(\preceq T+1/\sigma\). The same proof works
when the word weight is omitted, at any fixed order \(p\).
No price here grows with an upper cutoff for \(q\).

\subsubsection{The actual near divisors}
In a nonmodule high/high gap use the affine inequality for
\(a+j\cdot\chi_F\) at the frozen flat. The possibly very large
integer \(j\cdot(n_F^2)\) is included in \(a\). A dangerous
value requires \(|a|\le C(1+|j|_1)\) since \(\chi_F\) is
bounded. The cutoff therefore gives reciprocal \(\preceq L\),
independent of the two free end momenta. The final tilt, summed
over the retained labels, and the high diagonal changes are
\(\preceq B N^{-2}\), which preserve this gap. On a module
label a nonzero bare kinetic integer has absolute value at
least one. At a zero bare anomalous gap equation
\eqref{rate:ex:da} becomes
\[
 x^2-Xx+q(q-Y)=0,\qquad x=(n-1)/R.
\]
Indeed at \((b,\nu)=(1,1)\) the formula reads
\(D_A=2(n-1)(n-1-RX)+2R^2q(q-Y)=2R^2\{x(x-X)+q(q-Y)\}\), that is
\(w\cdot(w-v)=0\) for \(w=(x,q)\) and \(v=(X,Y)\). Similarly
\(D_F=2R^2(xX+qY)=2R^2\,w\cdot v\) at the module label, which is the
edge relation used above.
The rational root \(x\) is integral by monicity. All these legs
lie on the circle of diameter \([0,(X,Y)]\), giving at most
\(C H_{\rm geom}^2\) possible integer points per head. On a zero
same-copy boundary incidence from a marked head the remaining
macro is likewise integral by translation of two given ones.
At each high nondegenerate leg its two reciprocal-square
components have height at most a fixed power of \(H_{\rm geom}\).
For the undesignated pair let \(U,V\) be the differences from the
prescribed clock components in root-normalized coordinates. With
\((\mathcal A,\mathcal D)=(\mathcal Q-M^2,\mathcal Q/2)\) evaluated
at the final rounded background on \(F\) as in (T1), the terminal fit
in Section~7.6.3 gives the signed gap
\[
 \Delta=N^{-2}(\mathcal A U-\mathcal D V)
       +O_{s,\epsilon}\bigl(N^{-3}\lambda^{-C_s}E+r^2E+\eta_{\mathrm{rnd}}\bigr).
\]
The separation \(|\mathcal A U-\mathcal D V|\ge e^{-CB_o}\) persists
at this background, and the remainder is \(o(N^{-2}e^{-CB_o})\) for
the same \(C\) by Section~7.6.3. Thus \(|\Delta|\ge\tfrac12N^{-2}e^{-CB_o}\)
at all sufficiently large scales, giving \(|\Delta|^{-1}\le 2N^2e^{CB_o}\preceq N^2\).
The true high axial frequencies on empty
coordinates use the axial expansion at the old torus, with
terminal action changes \(O(r^2 E^C)\) harmless at this
precision. No Diophantine bound over arbitrary points of an
unbounded annulus is used.

\subsubsection{An empty axial Birkhoff direction}
Let \(h=(h_x,RY)\), \(Y\ne0\), be a moving head and let
\(e_n^{\tau_n}\) be the unit direction at \(n\notin F\) in the signed
Cartesian copy \(\tau_n=\pm1\). Write
\(\mathcal R(z)=|\Psi(z)|^2-\|z\|_{\ell^2}^2\). At \(a=0\), hence
\(z=z_F\), the two differentiated transverse coefficients are
\[
 2D\mathcal R(z_F)[e_n^{\tau_n}],\qquad
 2V_F\,D\Psi(z_F)[e_n^{\tau_n}].
\]
Here the real differentials are extended complex linearly to the signed
copies. For a real variation \(v\) supported on \(F^c\),
\(D\mathcal R(z_F)[v]=2\operatorname{Re}(\overline{V_F}D\Psi(z_F)[v])\),
since the mass differential vanishes in that direction.
By spatial translations their coefficient for old angular label
\(j\) and input copy \(\tau_n\) is supported only at
\(\ell=\tau_n n+j\cdot n_F\). Its modulus is bounded by
\(Ce^{-c|j|_1}\), uniformly in \(n\), by the analytic
\(L^2\) chart differential and the uniform old \(H^1\) bound
on the strip. In the density differential all coefficients
with \(\ell=0\) vanish exactly, since
\(\int_{\mathbb T}\mathcal R(z)\,dx/(2\pi)=0\) for every nearby \(z\).
This uses the full Birkhoff differential. Its old labels need not
have charge one. For example on the complex zero-action face
with a positive old coordinate \(z_p\) and a negative test at
\(n\ne p\), the dressing formula gives the term
\(-z_p^3e^{i(3p-n)x}/(p-n)^2\) in
\(D_{z_n^-}(\Psi^2)\). The following small-index estimates
allow any such charge. Only the free Fourier part of the
large-index differential uses the charge-one gap formulas.

Take \(M=N^{1/4}\) and ensure
\((1+\max_{m\in F}|m|)^2 C_JB\le N^{1/4}\). For \(|n|\le M\) we have
\(|\ell|\le2M,\ |j\cdot\Omega|\le C M\) and the true empty
axial frequency satisfies \(|\omega_n|\le C(1+n^2)\).
For density, the two transverse indices differ by
\((\ell,0)\) and their true gap has magnitude \(\gtrsim R\):
the nonzero bare difference is
\(|\ell(2h_x+\ell)|\gtrsim R\) and all extra times are
\(O(N^{1/2})=o(R)\). For the pairing differential the outgoing
index is \((-h_x+\ell,-RY)\), and the sum of the two transverse
energies is at least \(2R^2\). Both estimates hold also on
the reverse rows, as the real Hessian in canonical axial
coordinates has the corresponding signed transpose. For each
fixed \(j,\tau_n,h\) the discarded/divided matrix is a partial
translation or reflection on \(n\), so its weighted Schur price
has no factor counting \(|n|\le M\).

When \(|n|>M\), the coefficient at
\(k=\tau_n n+j\cdot n_F\) of the difference between the true
chart differential and its free Fourier part has modulus at most
\[
 C e^{-c|j|_1}\frac{\bk n}{\bk k^2}\le
 C' e^{-c|j|_1}|n|^{-1}
 \quad (|j|_1\le C_JB).
\]
This is the one-derivative gain \(h^1\to H^2\) of the
differential. Use the basic \(L^2\) differential bound on the
discarded angular tail. After summation the true minus free
one-head matrix in this range costs \(\preceq r/M\), including
both signed transposes. On a prescribed high axial input
\(|n|\ge H_x\) the cost improves to \(\preceq r/R\).
In detail, the displayed inequality uses \(|k|\ge|n|/2\), which holds
because \(|j\cdot n_F|\le C_JBZ\le M/Z\le|n|/2\) for \(|n|>M\) and
\(|j|_1\le C_JB\), by the choice of \(M\). The one-head matrix is the
product of the moving head, of physical size at most \(r|g|\) with
\(|g|\preceq1\), and the two differentiated coefficients. For a fixed
old label \(j\), input copy \(\tau_n\) and head, the map from the input
index \(n\) to the output index is a partial translation or
reflection, with at most one entry in every row and column. Its
entries in the range \(|n|>M\) are at most
\(rC'e^{-c|j|_1}M^{-1}\), and the ratio of the input and output
weights is at most a fixed power of
\(1+H_{\rm geom}/P+|j|_1Z/N\) times \(e^{O((1+|j|_1)/B^{K_0})}\). By
the Schur test its weighted norm is at most the product of these two
bounds. Summing over \(|j|_1\le C_JB\) costs \(C^f\preceq1\) by (Pr5),
and the sum over the at most \(Cd^2\) heads and two copies costs
\(\preceq1\). The labels \(|j|_1>C_JB\) carry only the \(L^2\) bound
\(Ce^{-c|j|_1}\), and their sum is \(\preceq e^{-cC_JB/2}\le N^{-J}\)
for \(C_J\) large. Hence the true minus free one-head matrix is
\(\preceq r/M\) in this range. On a prescribed high axial input,
\(|n|^{-1}\le H_x^{-1}=H_{\rm geom}/(c_0R)\preceq R^{-1}\) replaces
\(M^{-1}\).
For \(M<|n|<H_x\) the free Fourier part is covered directly
by \eqref{rate:ex:df}--\eqref{rate:ex:da} with axial endpoint \(q=0\).
If the other head is transverse, the forward density divisor
has main term \(2(n-b)(RX+1-b)\), and the pairing has main
term \(2(n-b)(n-1-RX)\), both \(\gtrsim MR\).
At a low axial destination in the reverse density placement,
the integer in the first parentheses of the reverse formula is
\(-Y^2-X^2\ne0\). Above that slab the high classification
applies. Terms with only axial fields are already represented
by the exact axial action-frequency equation. Varying the
axial chart point from the prescribed old torus by the other
physical packet changes these differentiated coefficients by
\(\preceq r\), by its analytic translated Taylor expansion
and the known packet Wiener norm. This gives a physical
one-head matrix error \(\preceq r^2\).

\subsubsection{The signed boundary Hessian and the divided source}
For an undesignated cubic Hamiltonian monomial in the field
coordinates, write the complexified term and its conjugate as
\[
 c_j e^{-i\nu t}
 \prod_{m=1}^3 z_{a_m}^{(\delta_m)}+\text{conjugate},
 \qquad z^{(+)}=z,\quad z^{(-)}=\bar z,\quad \delta_m=\pm1 .
\]
The momentum shifts use these same signs. Differentiate in slot
three to form its field equation. For each of the other two signed
inputs write \(b_m+\xi_m\), where \(b_m\) is its prescribed value
and \(\xi_m\) its variation in physical units. At a marked slot
\(b_m=r g_m(\tau)e^{-i\kappa_m t}\), with conjugation when required.
After factoring out the common coefficient and the Hamiltonian sign,
the two inputs contribute
\[
 (b_1+\xi_1)(b_2+\xi_2)
   =\underbrace{b_1\xi_2+\xi_1 b_2}_{\text{one prescribed input}}
    +\underbrace{b_1b_2}_{\text{two prescribed inputs}}
    +\xi_1\xi_2 .
\]
The first term gives a linear operator on the variation. Differentiating
in the other slots and including the conjugate monomial supplies the
reverse entries of the real Hessian. The second term is a forcing in
slot three that is already present at zero variation. In physical
units the matrix coefficient carries one factor \(r\), and the
forcing carries two. In the signed rotating variables both have the
same scalar mismatch up to the sign of the orientation,
\[
                     \Phi=\nu+\sum_m\delta_m D_{a_m}.
\]
Here the signs of \(\nu\) and of all slots are understood in
the displayed Hamiltonian convention. Taking
\(D_{a_m}=\kappa_{a_m}\) on marked slots makes this statement
exact, including when the free variation in one of those slots
belongs to the finite block. Using the genuine high diagonal
there instead changes all comparisons by at most \(3u\eta_0\),
well below the lower high/high gap. A repeated index carries
its elementary multiplicity in the Hessian. Every derivative
of a mixed finite/exterior edge, including an exterior input
and a finite output, is placed in the field matrix. Removing
only an inhomogeneous particular would leave one or both of
these homogeneous rows.

In doubled variables put \(D_{ii}=\varepsilon_i D_i\). If a
matrix entry is \(B_{ij}=A_{ij}(\tau)e^{-i\zeta t}\), set
\[
 \Delta_{ij}=\varepsilon_iD_i-\varepsilon_jD_j-\zeta,
 \qquad K_{ij}=-B_{ij}/\Delta_{ij}.
\]
Reverse signed Hessian entries satisfy
\(B_{ji}=\varepsilon_i\varepsilon_j\overline{B_{ij}}\)
with frequency \(-\zeta\), so
\(K_{ji}=-\varepsilon_i\varepsilon_j\overline{K_{ij}}\).
Indeed \(\Delta_{ji}=\varepsilon_jD_j-\varepsilon_iD_i+\zeta=-\Delta_{ij}\)
is real, and
\(K_{ji}=-B_{ji}/\Delta_{ji}
=\varepsilon_i\varepsilon_j\overline{B_{ij}}/\Delta_{ij}\).
Use \(\exp K\) including the reality-related entries. Each
retained head and angle acts by a partial shift or reflection
in either direction, bounded by its coefficient times the
ratio of radial weights and \(e^{O((1+|j|_1)/B^{K_0})}\).
The forward and reverse radial ratios are both bounded by
a fixed power of \(1+H_{\rm geom}/P+|j|\max_{m\in F}|m|/N\). The divided
homogeneous operator has accordingly no cardinality cost
from free endpoints.

The exact cancellation uses the covariant combination
\[
 ([D,K]-iK_t)_{ij}=-B_{ij}
       +i u\,\Delta_{ij}^{-1}A'_{ij}(\tau)e^{-i\zeta t}.
\]
To verify it, note that
\([D,K]_{ij}=(\varepsilon_iD_i-\varepsilon_jD_j)K_{ij}\) and, since
\(\tau=ut\),
\[
 -i\partial_tK_{ij}
 =\frac{i}{\Delta_{ij}}\bigl(uA'_{ij}(\tau)-i\zeta A_{ij}(\tau)\bigr)e^{-i\zeta t}
 =\frac{iu}{\Delta_{ij}}A'_{ij}(\tau)e^{-i\zeta t}-\zeta K_{ij},
\]
so that the sum is \(\Delta_{ij}K_{ij}\) plus the term with
\(A'_{ij}\).
For \(B_{ij}\) here its frozen coefficient \(c_j\) and
temporal divisor are evaluated at fixed terminal parameters.
The evolving prescribed envelopes satisfy
\(|g|+1\preceq1,\ |g'|\preceq\lambda^{-1}\) even before
removing their fast finite harmonics. Variation of the actual
old parameters in the original coefficients goes into the
old-to-field derivative on a known profile or a small
coefficient times an additional free field. In particular no
derivative of a nearly zero \(\Delta\) in an old action is used.
For the free part of the conjugation,
\[
 e^{-K}D e^K-i e^{-K}\partial_t e^K-D
 =\int_0^1 e^{-aK}\big([D,K]-iK_t\big)e^{aK}\,da .
\]
Thus subsequent terms of this conjugation use an additional
factor \(\norm K=o(1)\). Against a strong physical one-head
matrix of size \(\preceq r\), a near primitive costs
\(r(r\lambda N^2)/u=rN^2\) in slow units. The retained
bounded low operator commutes exactly with the part of \(K\)
supported on two high endpoints. Its commutator with the
remaining part costs only \(\preceq r/R\) in physical units,
even if that low operator is time dependent.

For the independent two-known source into a high exterior
coordinate, subtract its signed particular divided by the same
true mismatch. Its physical numerator is
\(r^2 c_j g_a(\tau)g_b(\tau)e^{-i\zeta t}\) with conjugations
as appropriate. In normalization by \(r\) the particular is
of the same order as \(K\) below. Applying a strong one-head
matrix and dividing the physical forcing by \(ru\) costs at
most the corresponding last-column price. Differentiating
\(g_ag_b\) costs at most one \(\lambda^{-1}\) in slow time.
Every independent mixed boundary source contains a decaying
transverse input and hence vanishes at the axial limit.
Contributions with only axial inputs are included in the
integrable axial flow. The operator \(K\) can have a nonzero
terminal value and is used as a change of variables, not as
a decaying independent forcing. Actual endpoint and section
coordinates are transformed along with it.

The following table suppresses factors \(\preceq1\).
In the middle column the physical primitive of a two-known
source is the listed quantity multiplied by \(r\).
The last column includes a strong one-head commutator and one
slow envelope derivative, in the normalized evolution units
\(u\) for matrices and \(ru\) for independent sources.
\[
\begin{array}{l|c|c|c}
\text{class}&
  |B_{\mathrm{phys}}|& |K|\text{ or }|p_{\mathrm{phys}}|/r
       &\text{normalized error}\\ \hline
\text{near high/high}&r\lambda&r\lambda N^2&rN^2\\
\text{noninteger or fast high/high}&r&rL&rL/\lambda\\
\text{divided low incident}&r&r/R&r/(R\lambda)\\
\text{low body commutator}& & &1/(R\lambda)\\
\text{empty axial differential error }(|n|>M)&r/M& &1/(M\lambda)
\end{array}
\]
For the column \(|B_{\mathrm{phys}}|\) the corresponding physical
two-known source bounds are \(r|B_{\mathrm{phys}}|\). The entries are
obtained as follows. The middle column is \(|B_{\mathrm{phys}}|\)
divided by the lower bound for the divisor: \(N^{-2}\) for a near
high/high entry, \(L^{-1}\) for a noninteger or fast one, and \(R\)
for a divided entry meeting the low body. The last column is the
product of a strong one-head matrix of size \(r\) with \(|K|\),
divided by \(u=r\lambda\); the term \(u\Delta^{-1}A'\) in the
covariant identity obeys the same bound for each class. Indeed, for a
summand with envelope \(A(\tau)=rc_jg(\tau)\) and frozen coefficient
\(c_j\), the bound on \(g'\) gives \(|A'|\preceq r|c_j|/\lambda\),
also at zeros of \(g\). Thus
\[
 \frac{r\cdot r\lambda N^2}{r\lambda}=rN^2,\qquad
 \frac{r\cdot rL}{r\lambda}=\frac{rL}\lambda,\qquad
 \frac{r\cdot r/R}{r\lambda}=\frac r{R\lambda},
\]
while the bounded low operator against a divided low incident
primitive gives \((r/R)/(r\lambda)=1/(R\lambda)\), and the empty
axial differential error is \((r/M)/(r\lambda)=1/(M\lambda)\). Long angular
tails cost any prescribed \(N^{-J}\) before a fixed finite
number of normalizations. Long weak tails in a homogeneous
field row and the stricter independent-source tail are
estimated on their appropriate separate cutoffs. Since \(s>2\),
\(rN^2=e^{bd}N^{2-s}\) has a fixed negative power margin.
The remaining displayed prices are small after the stated
order of \(\epsilon,Q,d/\sqrt{B/\log B}\) and the initial
scale. These estimates concern the exterior principal inverse
and its interface with the finite graph. The growing branch
and its prescribed stable data still require the finite
dynamical inverse and the separate nonlinear feedback estimates.

\subsubsection{A separate estimate for the evaluated source}
The linear term and the forcing term in the preceding expansion
require different weak angular cutoffs. A homogeneous tail of norm \(EN^{-h_w}\) multiplies
the unknown variation. Composing it with the field inverse, of norm
at most \(E\), costs \(E^2N^{-h_w}=o(1)\), so it can be absorbed.
The forcing with two prescribed inputs remains when the variation
is zero. Its image under the same inverse must fit inside the
field correction radius \(\rho=N^{-a}\) of
Section~\ref{rate:sec:ward-landing}. The exponent \(h_w\) can be
smaller than \(a\). We therefore keep and test the forcing terms
through the longer weak angular cutoff estimated below.

Let \(m(\tau)\) be the actual finite trajectory in normalized field
units \(r\), and \(\Pi\) its full marked projection. In defining
\(m\), keep on \(\Pi\) every term for which the only available
normalized estimate is \(EN^{-h_w}\). Such terms qualify for the
small perturbation in the finite contraction and are included in the
exact reference. Only the stronger independent-source bound permits
omitting them when the residual is evaluated. Form the normalized
inhomogeneous residual by inserting \(m\) into the full
field \emph{before} changing the unknown error. In these old
frame and rotating coordinates the actual substitution is
\(m+p+e^K x\), with \(p\) the divided particular in normalized
field units. Thus the matrix conjugation acts on \(x\). Its
time derivative does not act on \(m\) a second time.

On the homogeneous field block the retained density part is
physically compressed to \((1-\Pi)B_{\mathrm{ext}}(1-\Pi)\) with
its complete Hermitian reverse. It annihilates \(m\). If
marking either endpoint would have created a density edge into
an unmarked high coordinate, this is precisely a two-known
output of the cubic Hamiltonian. Its independent source must
be treated using the stronger cutoff, and its two homogeneous
Hessian orientations as in the preceding subsection. A
connection to a principal low star cannot be such a marked
endpoint by the word argument. During the integration of the
independent source, applying this exterior matrix to a known
particular costs its physical \(O(rE^C)\) one-head norm;
its large old low potential is supported entirely on the low
subspace. The known low part of a transformed profile is
estimated in the stronger cap \(A'\): on the low slab the
\(X_{p,A}\) norm is at most \(N^{A-A'}\) times
the \(X_{p,A'}\) norm. This pays the physical \(1/u\)
on that explicit known piece. The actual low potential on
an arbitrary unknown is the principal evolution throughout.

There is a sharper axial dichotomy for this source than the
one for the homogeneous field. If two marked transverse modes
give an empty axial output \(n\), translation covariance of
the full \(D\Psi\) tensor yields
\[
  \pm n+j\cdot n_F=\pm h_{a,x}\pm h_{b,x}.
\]
At \(|j|_1\le C_JB\) the right side is \(R\) times an integer
plus \(0\) or \(\pm2\), because \(h_x=1+RX\) on every marked head.
Hence either
\(|n|\le 2+C_JB\max_{m\in F}|m|\ll N^{1/4}\), or \(|n|\ge R/2\):
the first case is that of a vanishing integer multiple of \(R\), and
otherwise \(|n|\ge R-2-C_JB\max_{m\in F}|m|\).
In the second case the one-smoothing comparison costs
\(\preceq R^{-1}\) also when the axial coordinate is the
\emph{output}, by the signed Hessian transpose above.
In the first case the density term requires identical
horizontal macros, and equal vertical labels by its projection
to the axis. Thus the prescribed transverse plane waves are
identical, their density has only spatial harmonic zero, and
\[
 \int_{\mathbb T}D\mathcal R(z_F)[e_n^{\tau_n}]\,\frac{dx}{2\pi}=0
\]
for every old angle and every signed empty input. If one
prescribed plane has been deformed by the old transverse
frame, the corresponding term with \(T-1\) instead costs
\(\preceq R^{-1}\). The pairing term in the low case can
have opposite macro heads but its mismatch contains the sum
of the transverse kinetic energies and is \(\gtrsim R^2\).
Thus no source is left with only the homogeneous \(N^{-1/4}\)
gain. If the output is transverse low, a two-known short
angular word cannot reach it at all, by the horizontal word
test. A high transverse output with a prescribed axial input
already gains \(\preceq R^{-1}\) on its chart-differential
remainder, since the prescribed axial abscissa has nonzero
integral macro.

For any fixed target \(J_0\), choose the source weak cutoff
large enough that the enlarged weak strip gives normalized
cost \(\preceq N^{-J_0}\) after division by \(\lambda\).
For instance the factor
\(\lambda^{-1}\exp[-\epsilon B m_w/6]\) from weak length
\(m_w\) pays this for \(m_w>6(2J_0+2)/\epsilon\) when
\(1\le\log N/B\le2\). Indeed \(\lambda^{-1}=e^{\epsilon B}\) and
\(N^{-J_0}\ge e^{-2J_0B}\), so it suffices that
\(\epsilon m_w/6\ge\epsilon+2J_0\), which follows from the stated
inequality because \(\epsilon\le2\). All stationary words below this
cutoff are explicitly tested as two-known words. In particular
discarding every nonsingleton on an independent source using
its homogeneous \(EN^{-h_w}\) estimate alone would not suffice.
A retained stationary charged coefficient has summed size
\(\preceq\lambda\): the exact five-letter labels are precisely
the fitted singletons by the disjoint-pivot module, and the
others have at least fifteen weak letters. Nonstationary
labels may have order-one coefficients but use the fast
denominator.

The finite-only block, both types of exterior output, repeated
indices, and the angular/weak tails now exhaust the degree-three
Hamiltonian evaluated on two known slots. Its noncancelled
particular or one-head residual costs
\(\preceq L^C(rN^2+R^{-1})\) in the field equation divided
by \(ru\). Physical field terms of degree at least three
in its vector field give \(O(r^3 E^C)\), hence
\(\preceq r/\lambda\) after the same normalization.
The prescribed spatial residual on a known field costs
\(\preceq N^{-J_0}/u\); one can use the larger
\(\preceq N^{-J_0}/u^2\) common allowance when the old
variables are included. At zero unknown old variables the
implicit old action displacement is \(O(r^2L^CE^C)\), and
the mismatch of old isolated and full terminal frequencies
is \(O(r^2E^C)\). Their field connections on a known wave
therefore cost at most \(\preceq L^C r/\lambda\), with
a further \(N^{-1}\) where the connection differentiates
a small spatial column or scalar primitive. Mixed terms
decay at the true terminal point by their transverse
factor, including the envelope derivative since the actual
bounded finite vector field vanishes there and has
\(C^1\) bound \(\preceq\lambda^{-1}\) on the joining
segment. Pure axial terms use the exact diagonal action flow
and its prescribed endpoint. In the common time weight the
field source is consequently bounded by
\[
 E^C L^C\big(rN^2+N^{-1}+r+N^{-J_0}/u^2\big),
\]
where \(R^{-1}=P/N\preceq N^{-1}\). Choosing \(J_0\)
above the fixed normalization prices gives the asserted power
for the independent source.
We list the arithmetic behind this bound; every source is divided
by \(ru=r^2\lambda\), and factors \(\preceq1\) are suppressed.
\begin{itemize}
\item A two-known source \(r^2c\,g_ag_b\,e^{-i\zeta t}\) with mismatch
\(\Delta\) is replaced by its particular of physical size
\(r^2|c|/|\Delta|\). By the covariant identity of the preceding
subsection there remain the term with the derivative of the envelope,
of physical size \(u\,r^2|c|\,|(g_ag_b)'|/|\Delta|\le
u\,r^2|c|/(\lambda|\Delta|)\), and the image of the particular under
the strong one-head matrices, of size \(r\cdot r^2|c|/|\Delta|\). For
a stationary charged word \(|c|\preceq\lambda\) and
\(|\Delta|^{-1}\preceq N^2\); both remainders are then
\(r^3\lambda N^2\), that is \(rN^2\) after normalization. For a
nonstationary word \(|c|\preceq1\) and \(|\Delta|^{-1}\preceq L\);
both are \(r^3L\), that is \(rL/\lambda\le L^Cr\).
\item In a two-known source whose output is a high axial coordinate,
or which has a prescribed axial input, the difference between the
true chart differential and its Fourier part has physical size
\(r^2/R\); so has every term containing the factor \(T-1\) of the
transverse frame on a marked head. These terms are left undivided,
and their normalized size is \(1/(R\lambda)\le L^CR^{-1}\).
\item The terms of degree at least three give
\(r^3/(r^2\lambda)=r/\lambda\).
\item The spatial residual \(N^{-J_0}\) of the high columns on a known
wave of amplitude \(r\) gives \(rN^{-J_0}/(r^2\lambda)=N^{-J_0}/u\).
\item At zero unknown old variables the old velocities differ from
their terminal values by \(O(r^2L^C)\). Multiplied by the parameter
derivative of the frame on a known wave of amplitude \(r\), this gives
\(r^3L^C/(r^2\lambda)=L^Cr/\lambda\), with a further factor \(N^{-1}\)
when a small column or scalar primitive is differentiated.
\end{itemize}
The sum is bounded by the displayed expression, since
\(N^{-J_0}/u\le N^{-J_0}/u^2\).
This estimate is for the inserted
field. The old action and angle equations require their own
differentiated Hamiltonian identities.

\subsection{The occupied actions and the full half-line problem}\label{rate:sec:ward-landing}
We work in the autonomous mass gauge. With actions \(I_k=|z_k|^2\)
and angles \(z_k=\sqrt{I_k}e^{-i\phi_k}\),
\(\dot\phi=\partial_I\mathcal K\) and
\(\dot I=-\partial_\phi\mathcal K\). On the physical transverse
space use the real scalar product
\(\langle f,g\rangle=2\operatorname{Re}\sum_k\bar f_k g_k\)
and complex structure \({\mathbb J}=-i\), transported together
to its real representation. Then \({\mathbb J}^*=-{\mathbb J}\),
\({\mathbb J}^2=-1\) and \(\dot w={\mathbb J}\nabla_w\mathcal K\)
with precisely this scalar product. Dots denote physical time
derivatives. Primes denote derivatives in \(\tau=ut\), where
\(u=r\lambda\).

\subsubsection{Field estimates used in this step}
Here are the precise estimates on the field problem used below. They distinguish the inverse on the finite calibrated graph from the inverse on the full equation at the solution we will select. Fix
\[
 a=\tfrac1{16}\min(s-2,\tfrac14),\qquad
 \rho=N^{-a},\qquad \rho_o=N^{-3a}.
\]
Choose \(A>2s+12\), \(A'=A+3s+20\), \(J_0>2s+12\) and \(p>2\lceil s\rceil+10\). These are the inequalities used in the present subsection; they hold for the values \(A>4s+20\), \(A'=A+3s+20\), \(J_0>4(A'+s+p+10)\) fixed in Section~\ref{rate:sec:parameters}, which are the ones understood throughout, and \(a\) is the same number as there. With the same \(K_0\) and the signed word length \(\ell\), use
\[
 X_c=\ell^2(b^pW_c),\qquad b(k)=1+(|k_x|+|k_y|)/N,\qquad
 W_c=\min(N^c,\exp(\ell(k)/B^{K_0})).
\]
Constants denoted by \(E\) have the scale dependence prescribed in the registration of the clocks. Write \(\mathcal P_0\) for a fixed product of these constants and fixed powers of \(\lambda^{-1}\), including the fixed occupied derivatives. The constants are chosen so that all the resulting products of type \(\mathcal P_0\) are at most \(N^\eta\), \(\eta<a/20\). We keep track of the smaller \(E\) prices at the field-only inversion: each fixed product there is \(o(N^{h_w})\) for the homogeneous weak tail exponent \(h_w\).

Use the same stable section at \(\tau=T\) at both time weights. For \(\beta\) equal to either weight put \(w_\beta(\tau)=\exp(\beta(\tau-T)_+)\) and \(\|f\|_\beta=\sup w_\beta\|f(\tau)\|\). The two weights are strictly smaller than the decay exponent of the actual finite reference. The following are the field estimates to be used, in physical field units \(r\):
\begin{enumerate}
\item The actual finite reference reaches its specified axial endpoint with \(\|Z_\perp\|_\beta+\|Z_{\mathrm{ax}}-A_\infty(t)\|_\beta\le E\), and \(\sup\|A_\infty(t)\|\le E\). Here \(A_\infty(t)=a_\infty(t)/r\) includes the fixed terminal rotations. Its full real finite linearization has a mixed inverse of norm \(E\) with arbitrary stable trace at \(T\). The estimate is uniform also on strips \([0,L]\), \(L\ge T\), with complementary data prescribed at the right end.
\item All known profiles, explicit two-known particulars, old frames, frozen mixed changes and their inverses needed in the field equations have their stated translated operator or vector bounds at both caps \(A,A'\). In particular the known transverse profiles and particulars have \(X_{A'}\) bounds \(E\) with the temporal decay above. The field coordinate change on an error, denoted below by \(U(t)\), costs \(E\) on \(X_A\) and on \(\ell^2(b^q)\), \(q=0,2,s\). The unweighted \(b^q\)-Wiener norm of any \(X_A\)-bounded reference costs at most \(E\) times its \(X_A\) norm. The analytic axial chart tensors have the corresponding translated convolution bounds. At either cap an elementary translation by an old angular label \(l\) and a fixed number of heads costs at most \(C(1+d+(|l|_1+1)Z/N)^p\exp((|l|_1+C)/B^{K_0})\). The resulting operator bounds pass to the exponentials and inverse square roots by their norm-convergent series, as verified below. The cardinality bound for short words changes only by replacing \(A\) by the fixed \(A'\).
\item In the pure old transverse frame \(w=T(I,\phi)y\), put \(D=\Omega(I)\cdot\partial_\phi\), \(C=T^*BT+T^*{\mathbb J}DT\), where \(B\) is the physical transverse Hessian at vanishing new variables and \(\Omega(I)=\partial_I\mathcal K_{\mathrm{ax}}(I,0)\). Let \(\Omega_*=\partial_I\mathcal K_{\mathrm{ax}}(I^*,|a_\infty|^2)\). If \(Q_L\) is the compressed low transverse projection, the high-high and high-low parts, including the transpose, obey
\begin{equation}\label{rate:eq:ward-col}
 \|(\partial_\phi C)_{HH,HL,LH}\|
 \le \mathcal P_0\big(N^{-J_0}+N^{-1}|\Omega(I)-\Omega_*|\big).
\end{equation}
The low part apart from the free kinetic is bounded and selfadjoint on each signed copy, with fixed derivatives bounded by \(\mathcal P_0\). The remainder in \eqref{rate:eq:ward-col} enjoys the same estimate after another angle derivative. On differentiating in \(I\) one must also differentiate \(\Omega(I)-\Omega_*\), allowing \(\mathcal P_0N^{-1}\|D_I\Omega\|\) separately. The high scalar mean is action-only, with first derivative bounded by \(\mathcal P_0\). The old spatial transformations and their fixed parameter derivatives are bounded by \(\mathcal P_0\). For a temporal primitive used in these assertions the divisor is evaluated at the \emph{fixed} terminal \(\Omega_*\).
\item The field-only principal inverse uses the finite graph and the exterior Hermitian components in the flattened norm. It has norm \(E\) for every path of admissible real trial \(I,\phi\) in its low coefficient, with the same component supports and the same two boundary conditions. This follows from the finite estimate and unitary backward integration in the exterior, costing \(E(T+\beta^{-1})\). After the physical divided rows and explicit particulars have been restored, the additional nominal homogeneous field operators have norm at most
\[
 E N^{-h_w}+\mathcal P_0\big(N^{-1/4}+rN^2+r+N^{-1}\big).
\]
The nominal normalized known field source away from the undifferentiated low principal coefficient costs
\[
 \mathcal P_0\big(rN^2+N^{-1}+r+N^{-J_0}/u^2\big).
\]
The \(rN^2\) term allows a physical strong row \(O(rE)\) times a divided near primitive \(O(r\lambda N^2 E)\), divided by \(u\). Near primitives have two high ends. A homogeneous primitive meeting the compressed low body has physical size \(Er/R\). Their adjoints obey the same estimates. The stronger cutoff on a two-known source is included in this last source budget.
\end{enumerate}
Items 1--4 follow, respectively, from the perturbed finite landing,
the translated Fourier algebra, the transverse spatial reduction,
and the exterior-component calculation. We now use them to solve
for the occupied actions jointly with the field and then to invert
the physical variational equation on the resulting solution.

\smallskip
\noindent\emph{Verification of items 1--4.} We indicate how each item
is obtained from the preceding sections.

\emph{Item 1.} By Lemma~\ref{rate:l:finite} the ideal field \(F_0\)
has a landing orbit \(z_*\), with \(\|z_*-a\|_\sigma\le e^{Cd}\) as
noted in the proof of Lemma~\ref{rate:l:round}, and a
mixed inverse of norm at most \(e^{CdD_g}\), \(D_g=C_*d\log(d+2)\),
for every weight \(\kappa/8\le\sigma\le\kappa/3\), uniformly on finite
strips; the section is at \(T=D_g\) in the slow time fixed by fitting
\(2\overline{v_e}/\lambda\) to its coefficients. By the last paragraph
of the subsection on the projected reference, the actual finite field
is conjugate, through \(z=y+P\) with
\(\|P\|_{C^2}\le\mathcal Er\lambda N^2+\mathcal Pr\), to \(F_0\) plus a
\(C^1\) error at most
\(\|S\|_{C^1}+\mathcal ErN^2+\mathcal Pr\), where \(S\) consists of
the fitting errors \(O(e^{-C_4dD_g})\), the hot action detuning
\(\mathcal Er/\lambda\) and the undivided angular tails. In the
marked projection the wrong charged harmonics within the cutoffs have
already been included in the primitive \(P\). The longer words retained
in \(S\) cost \(O(N^{-J_0})\) in normalized \(C^1\) by the angular
and weak-letter cutoffs. Even the weaker bound \(E^C N^{-h_w}\) for
a retained nonsingleton suffices here. In the exterior equation this
bound multiplies a field variation, as in item~4.

Choose \(C_4>C_3\). Once \(\epsilon\) and the cutoffs are fixed,
\(h_w>0\) is fixed. Take \(Q\) large and \(\alpha\) small, in the
order of Section~\ref{rate:sec:parameters}, so that
\[
 E^C\le N^{h_w/4},\qquad
 e^{C_3dD_g}\le N^{\min\{h_w/4,(s-2)/4,J_0/2\}} .
\]
The second inequality follows from \(dD_g\le C_*\alpha^2 B\) and
\(\log N\ge B-O(1)\). Since \(r=HN^{-s}\) with \(H=N^{o(1)}\), the
price bounds give \(\mathcal ErN^2+\mathcal Pr+\mathcal Er/\lambda
\le N^{-(s-2)/2}\). After multiplication by \(e^{C_3dD_g}\), this
bound, both tail estimates and the fitting error \(O(e^{-C_4dD_g})\)
all tend to zero. The perturbation therefore has size at most
\(e^{-C_3dD_g}\) for sufficiently large \(B_*\).
Lemma~\ref{rate:l:round}, which holds at every weight in
\([\kappa/8,\kappa/3]\), then gives the landing of the actual finite
field at the weight \(\kappa/3\) and its mixed inverses at the two
weights \(\beta\in\{\kappa/4,\kappa/8\}\), with the same bounds up to
a factor two. The conjugacy and its differential are bounded.
Finally \(e^{Cd}\) and \(e^{CdD_g}\) are at most \(E\), because
\(dD_g\le C_*d^2\log(2+B)\), and the endpoint amplitudes in units of
\(r\) are \(1\) on the top and \(b_pc_{nk}\le2^{d/2}\) on the private
axes.

\emph{Item 2.} The elementary kernels in the spatial changes and
the divided rows are sums \(\sum_{l,\mathsf h}c_{l,\mathsf h}\,
\tau_{l,\mathsf h}\) of partial shifts or reflections by an old
angular label \(l\) and a bounded collection \(\mathsf h\) of signed
heads. A shift by a word of length \(\ell_0\) changes the capped word
weight by a factor at most \(e^{\ell_0/B^{K_0}}\) for every cap,
because \(\min(N^c,\xi\eta)\le\min(N^c,\xi)\,\eta\) for
\(\xi,\eta\ge1\). With the polynomial ratio of the subsection on the
Fourier algebra,
\[
 \|\tau_{l,\mathsf h}\|_{X_c\to X_c}
 \le C\bigl(1+d+(|l|_1+1)Z/N\bigr)^{p}e^{(|l|_1+C)/B^{K_0}},
 \qquad c\in\{A,A'\},
\]
and the same bound holds on \(\ell^2(b^q)\), \(q=0,2,s\). Taking the
coefficient supremum over input indices and summing the angular
decay \(e^{-a_0|l|_1}\) costs at most \(C^{|F|}\le E\).
In each of these spaces the preceding coefficient and divisor bounds
give \(O(ER^{-2})\) for the off-copy Sylvester generators and
\(\preceq r\lambda N^2+r\lambda^{-C_s}+r/R\) for the divided mixed
generator, with the fixed exponent \(C_s\) from the noninteger divisor bound.
Both bounds are \(o(1)\). For any such generator \(K\), summing its
operator series gives \(\|e^{\pm K}\|\le e^{\|K\|}\) in the same norm.
The column and adjoint estimates of
Section~\ref{rate:sec:transverse-principal} likewise give
\(\|G-1\|=O(E(H_x^{-1}+R^{-2}))=o(1)\) on the high and low input
spaces for their respective Gram operators. The binomial series for
\(G^{\pm1/2}\) therefore converge in each operator norm. These bounds
pass to the finite compositions defining the changes and their inverses.
The known profiles are finite sums of marked plane waves, whose word
length is one, with bounded envelopes. Applying the divided rows to
pairs of these waves gives the stated vector bounds at the cap \(A'\). The Wiener bound is the embedding of \(X_A\) into
\(\ell^1(b^q)\), \(q\le p-2\), proved there: the positions of length at
most \(CB^{K_0+1}\) number at most \(E\), and on the remaining
positions the weight \(W_A\) equals \(N^A\), which pays the factor
\(\sqrt{NP_{\mathrm{lat}}}\) of the Cauchy--Schwarz inequality.
Replacing \(A\) by \(A'\) changes only this last power.

\emph{Item 3.} In the finished frame of
Section~\ref{rate:sec:transverse-principal} decompose
\(C=T^*BT+T^*{\mathbb J}DT\) along the high and low inputs. By the
construction there, its high-high block is the real diagonal
\(E_n+\bar d(n,m;I)\), where \(\bar d\) is the angle mean of the high
scalar and depends on the actions only, plus a remainder. The high-low
blocks are remainders. The low block consists of the free kinetic term,
a bounded Hermitian operator on each signed copy, and a residual
coupling between copies of size \(O(E_cR^{-2j_0})\). This coupling is
included in the field remainder. The high-block remainder has three
sources. The first is the residual \(R_h'\) of the completed columns,
at most \(E_c(H_x^{-M_0}+R^{-2j_0})\), which is at most
\(\mathcal P_0N^{-J_0}\) for the orders \(M_0,j_0\) chosen at the end
of that section, since \(H_x,R\ge N^{3/4}\). The second is the
undivided angular and weak tails, also at most
\(\mathcal P_0N^{-J_0}\). The third comes from the phase \(b_3\). If
\(\sum_{\ell\ne0}\widehat d_\ell e^{i\ell\cdot\phi}\) is the divided
part of the scalar, then
\(b_3=\sum_\ell\widehat d_\ell e^{i\ell\cdot\phi}/(i\ell\cdot\Omega_*)\)
with the fixed vector \(\Omega_*\), while \(D=\Omega(I)\cdot\partial_\phi\).
Hence
\[
 Db_3=\sum_{\ell\ne0}\widehat d_\ell e^{i\ell\cdot\phi}
      +(\Omega(I)-\Omega_*)\cdot\partial_\phi b_3 ,
 \qquad
 \|(\Omega(I)-\Omega_*)\cdot\partial_\phi b_3\|
 \le E_c'N^{-1}|\Omega(I)-\Omega_*| .
\]
The first term cancels the oscillating scalar and the second is the
stated mismatch; the phase \(b_2\) involves no divisor and leaves no
mismatch. An angle derivative annihilates the action-only diagonal
and, on the fixed strip, costs a fixed factor on each remainder.
This proves \eqref{rate:eq:ward-col} and its analogue with one more
angle derivative. Differentiating the mismatch in \(I\) gives the
additional term \((D_I\Omega)\cdot\partial_\phi b_3\), of size
\(\mathcal P_0N^{-1}\|D_I\Omega\|\). The bounds on the frame and on
the low block were proved in the same section.

\emph{Item 4.} After the divided rows have been removed the principal
operator is the direct sum of the finite graph and of the exterior
same-copy components. The first is inverted by item~1. On the second
the retained matrices are Hermitian in the flattened norm for every
real trial path, so the propagation is unitary there and the backward
integral has norm at most \(T+\beta^{-1}\), with
\(\beta^{-1}\le8\kappa^{-1}=8(2+d+D_g)^{10}\le E\). The additional
homogeneous operators are the entries of the last column of the table
in Section~\ref{rate:sec:exterior}: \(rN^2\) for the near primitives,
\(rL/\lambda\le\mathcal P_0r\) for the fast ones,
\(r/(R\lambda)+1/(R\lambda)\le\mathcal P_0N^{-1}\) for the primitives
meeting the low body, \(1/(M\lambda)=N^{-1/4}/\lambda\) for the empty
axial differential, together with the weak tail \(EN^{-h_w}\) and the
change \(r^2/u=r/\lambda\) of the one-head matrices with the chart
point. The source bound is the last display of that section.

\subsubsection{Recovery of the actions and differentiation of the weak harmonics}
The occupied list \(F\) contains the previous sites and the newly prepared spectator, moderator and weak sites. The \emph{current} high axial endpoints are in the Cartesian complement \(a_{\mathrm{ax}}\). In the frame preceding the frozen mixed change write
\[
 (a_{\mathrm{ax}}/r,Y)=Z_{\mathrm{kn}}(t)+U(t)e,\qquad Y=y/r .
\]
This includes the full correction and its frozen rotations. The particulars are included in \(Z_{\mathrm{kn}}\). In particular
\(\|Y\|_{\beta,X_A}+\|a_{\mathrm{ax}}/r-A_\infty(t)\|_{\beta,X_A}\le E\)
on \(\|e\|_{\beta,X_A}\le\rho\). The prescribed rotations, divisors, and arguments of the known fast phases are fixed at this stage.

Set \(G_i=T^*{\mathbb J}\partial_{\phi_i}T\), and prescribe normalized coordinates \(j,\theta\) by
\begin{equation}\label{rate:eq:ward-I}
 \phi=\phi_*(t)+\theta,\qquad
 I=I^*+u j+\tfrac{r^2}{2}
              \big(\langle Y,G_i(I,\phi)Y\rangle\big)_{i\in F},
 \qquad \phi_*(t)=\phi_*(0)+t\Omega_* .
\end{equation}
The scalar product is real. The \(I\)-Lipschitz constant of the last term is at most \(\mathcal P_0 r^2\|Y\|_2^2=o(1)\), with \(2\) here and in scalar products denoting ordinary \(L^2\). The finite Euclidean summation over \(F\) is included in \(E\). The right side stays inside a fixed relative fraction of each positive floor: the previous floors have inverse price \(\exp(O_s(B_o))\), and the inverse of the present minimum on \(F\) costs at most \(\lambda^{-C_s}E\). Thus \eqref{rate:eq:ward-I} has one solution in this neighborhood and
\begin{equation}\label{rate:eq:ward-DI}
 |I-I^*|\le u|j|+\mathcal P_0r^2\|Y\|_2^2,\quad
 \|I_j\|\le2u,\quad
 \|I_\theta\|\le \mathcal P_0r^2\|Y\|_2^2,\quad
 \|I_e\|\le \mathcal P_0r^2\|Y\|_2 .
\end{equation}
The last differential is bounded on ordinary \(L^2\) as well. Replacing \(Y\) by the finite graph in this formula would miss the actual low return.
In detail, write \(\Psi(I)=\tfrac{r^2}2\bigl(\langle Y,G_i(I,\phi)Y\rangle\bigr)_{i\in F}\),
so that \eqref{rate:eq:ward-I} reads \(I=I^*+uj+\Psi(I)\). On the stated
neighborhood \(\|G_i\|+\|\partial_IG_i\|+\|\partial_\phi G_i\|\le\mathcal P_0\) on
\(L^2\), so \(\operatorname{Lip}_I\Psi\le q\coloneqq\mathcal P_0r^2\|Y\|_2^2\le1/2\) and
\(|\Psi|\le\mathcal P_0r^2\|Y\|_2^2\). The contraction principle gives a unique
solution with \(|I-I^*|\le u|j|+\mathcal P_0r^2\|Y\|_2^2\), which is much smaller
than every positive floor since \(u\rho_o\) and \(r^2E\) are. Differentiating,
\(I_j=(1-\partial_I\Psi)^{-1}u\), \(I_\theta=(1-\partial_I\Psi)^{-1}\partial_\phi\Psi\) and
\(I_e=(1-\partial_I\Psi)^{-1}\partial_e\Psi\), with \(\|(1-\partial_I\Psi)^{-1}\|\le2\),
\(\|\partial_\phi\Psi\|\le\mathcal P_0r^2\|Y\|_2^2\) and
\(\partial_e\Psi[\delta e]=r^2\bigl(\langle Y,G_i\,\delta Y\rangle\bigr)_i\), where
\(\delta Y=U(t)\delta e\) has \(\|\delta Y\|_2\le E\|\delta e\|_2\). This gives
\eqref{rate:eq:ward-DI}, the last bound also on ordinary \(L^2\).

On a complex occupied-action ball of radius \(c\min_F I_i^*\), the signed square roots have the same bounds, up to fixed factors, as at \(I^*\). The enlarged strip \(|\operatorname{Im}\phi_i|<\epsilon B/6\) on the weak sites still belongs to the small Cartesian chart, by the initial \(O(\lambda^{1/5}E)\) radii and the fixed slack in the exponent. An angular coefficient at \(l\) therefore has the additional factor \(\exp[-(\epsilon B/6)m_w(l)]\), where \(m_w(l)\) counts its weak letters, besides ordinary exponential decay on a smaller angular strip. An angle derivative is paid by this ordinary strip. By Cauchy's formula, each fixed action derivative costs \(C(\min_F I_i^*)^{-1}\), independently of \(m_w(l)\) and of the cutoff used for the known forcing. The same Cauchy estimate applies along every complex \(\ell^2\) unit direction. An undivided non-singleton homogeneous weak tail has normalized norm \(EN^{-h_w}\) at the frozen parameters. Its action derivative may cost \(\mathcal P_0 N^{-h_w}\). We absorb only the undifferentiated tail at the field-only inversion. On parameter differences one may use the coarser normalized bound \(\mathcal P_0(|\theta|+|I-I^*|)\), even for the sum over the entire weak tail. Its contribution to a field derivative through \(I\) contains the \(I_e\) of \eqref{rate:eq:ward-DI}. In the old torque these same mixed Hamiltonian terms are taken before division and keep their physical factor \(r^3\).

\subsubsection{The exact torque and the angular shear}
In physical transverse variables the Hamiltonian takes the form
\[
 \mathcal K=\mathcal K_{\mathrm{ax}}(I,|a_{\mathrm{ax}}|^2)
       +\tfrac12\langle w,B(I,\phi)w\rangle
       +V(I,\phi,a_{\mathrm{ax}},w).
\]
Every term of \(V\) has at least two transverse factors and degree at least three in \((a_{\mathrm{ax}},w)\). The unbounded kinetic in \(B\) is independent of \(I,\phi\). Define
\[
 S=T^{-1}{\mathbb J}V_w
    -T^{-1}\big[(\partial_I T)\dot I+
             (\partial_\phi T)(\dot\phi-\Omega(I))\big]y .
\]
The velocities on the right are the original Hamiltonian vector field evaluated at the instantaneous state. Thus \(\dot y={\mathbb J}Cy+S\). Differentiating \(T^*{\mathbb J}T={\mathbb J}\) at fixed \(I\) gives
\[
 \partial_{\phi_i}C=T^*(\partial_{\phi_i}B)T+D G_i
                     +({\mathbb J}C)^*G_i+G_i{\mathbb J}C .
\]
Here \(T^{-1}\partial_{\phi_i}T=-{\mathbb J}G_i\), and \(D\) commutes
with \(\partial_{\phi_i}\) at fixed \(I\). To see which action to correct,
differentiate \(q_i=\langle y,G_iy\rangle/2\) along \(\dot y={\mathbb J}Cy+S\).
The terms from \({\mathbb J}Cy\) and \(DG_i\) sum to
\[
 \tfrac12\langle y,(\partial_{\phi_i}C-
                          T^*(\partial_{\phi_i}B)T)y\rangle .
\]
The original action satisfies
\(\dot I_i=-\langle y,T^*(\partial_{\phi_i}B)Ty\rangle/2-V_{\phi_i}\).
Subtracting \(\dot q_i\) therefore cancels the derivative of the physical
Hessian. The corrected action is \(I_i^\sharp=I_i-q_i\), and its equation
in the variables \eqref{rate:eq:ward-I} is
\begin{align}\label{rate:eq:ward-torque}
 u^2 j_i'={}&-\tfrac12\langle y,(\partial_{\phi_i}C)y\rangle
             -V_{\phi_i}-\langle S,G_i y\rangle \\
 &-\tfrac12\big\langle y,
       [(\partial_I G_i)\dot I+
       (\partial_\phi G_i)(\dot\phi-\Omega(I))]y\big\rangle .\nonumber
\end{align}
The derivatives of \(V\) here are taken at fixed physical \(w\). Formula \eqref{rate:eq:ward-torque} is a vector field on state space. No derivative of a trial path appears in it.

In detail, \(G_i\) is symmetric, so along \(\dot y={\mathbb J}Cy+S\)
\[
 \dot q_i=\langle\dot y,G_iy\rangle+\tfrac12\langle y,\dot G_iy\rangle
 =\tfrac12\big\langle y,\big[({\mathbb J}C)^*G_i+G_i{\mathbb J}C\big]y\big\rangle
   +\langle S,G_iy\rangle+\tfrac12\langle y,\dot G_iy\rangle ,
\]
and \(\dot G_i=DG_i+(\partial_IG_i)\dot I
+(\partial_\phi G_i)(\dot\phi-\Omega(I))\) by the definition of \(D\).
By (T2) the three operators \(({\mathbb J}C)^*G_i\), \(G_i{\mathbb J}C\)
and \(DG_i\) add up to
\(\partial_{\phi_i}C-T^*(\partial_{\phi_i}B)T\). On the other hand
\(\dot I_i=-\partial_{\phi_i}\mathcal K
=-\tfrac12\langle w,(\partial_{\phi_i}B)w\rangle-V_{\phi_i}\) with
\(w=Ty\). Subtracting, the two terms with \(\partial_{\phi_i}B\)
cancel and
\[
 \frac{d}{dt}(I_i-q_i)=-\tfrac12\langle y,(\partial_{\phi_i}C)y\rangle
   -V_{\phi_i}-\langle S,G_iy\rangle
   -\tfrac12\big\langle y,\big[(\partial_IG_i)\dot I
     +(\partial_\phi G_i)(\dot\phi-\Omega(I))\big]y\big\rangle .
\]
Finally \eqref{rate:eq:ward-I} reads \(I_i-q_i=I_i^*+uj_i\), because
\(y=rY\), and \(d/dt=u\,d/d\tau\); hence the left side is \(u^2j_i'\).

The bounded old potential and the degree structure, using the convergent chart expansion with all undifferentiated factors in their packet or Wiener norms, give
\[
 |\dot I|\le \mathcal P_0r^2\|Y\|_2^2,\qquad
 |\dot\phi-\Omega(I)|\le \mathcal P_0r^2,\qquad
 |V_\phi|\le \mathcal P_0r^3 w_\beta^{-1},\qquad
 \|V_w\|_2+\|S\|_2\le \mathcal P_0r^2 w_\beta^{-1}.
\]
To derive these bounds recall that \(w=rTY\) with \(T\) and \(T^{-1}\)
bounded by \(\mathcal P_0\), that \(\|Y(\tau)\|_2\) and the
\(b^q\)-Wiener norm of \(Y(\tau)\) are at most \(Ew_\beta^{-1}\) by
item~2, and that the axial Cartesian field has
\(\||a_{\mathrm{ax}}|^2\|_1\le r^2E\).
(i) Since \(\mathcal K_{\mathrm{ax}}\) does not depend on \(\phi\),
\(\dot I_i=-\tfrac12\langle w,(\partial_{\phi_i}B)w\rangle-V_{\phi_i}\).
The first term is at most
\(\mathcal P_0\|w\|_2^2\le\mathcal P_0r^2\|Y\|_2^2\), since
\(\|w\|_2\le r\|T\|\,\|Y\|_2\). Every term of
\(V\) and of \(V_\phi\) has at least two transverse factors and degree
at least three; bounding two transverse factors in \(L^2\) and the
others in the Wiener norm, each of physical size at most \(rE\), gives
\(|V_\phi|\le\mathcal P_0r^3\|Y\|_2^2\le\mathcal P_0r^3w_\beta^{-1}\).
Hence \(|\dot I|\le\mathcal P_0r^2\|Y\|_2^2\).
(ii) \(\dot\phi-\Omega(I)=[\partial_I\mathcal K_{\mathrm{ax}}(I,|a_{\mathrm{ax}}|^2)
-\partial_I\mathcal K_{\mathrm{ax}}(I,0)]+\tfrac12\langle w,B_Iw\rangle+V_I\).
The bracket equals
\(\bigl(\chi_n(I,|a_{\mathrm{ax}}|^2)-\chi_n(I,0)\bigr)_{n\in F}\).
Its \(\ell^2(F)\) norm is at most \(C\||a_{\mathrm{ax}}|^2\|_2\le Cr^2E\)
by the local Lipschitz bound for \(\chi\). The other two terms are at most \(\mathcal P_0r^2\|Y\|_2^2\) and
\(\mathcal P_0r^3\), the action derivatives costing \(\mathcal P_0\).
(iii) Differentiating \(V\) in \(w\) leaves at least one transverse
factor and degree at least two, so
\(\|V_w\|_2\le\mathcal P_0\,(rEw_\beta^{-1})(rE)\le\mathcal P_0r^2w_\beta^{-1}\).
In \(S\) the second term is at most
\(\mathcal P_0(|\dot I|+|\dot\phi-\Omega(I)|)\|y\|_2
\le\mathcal P_0r^2\cdot rEw_\beta^{-1}\).
(iv) The last line of \eqref{rate:eq:ward-torque} is at most
\(\mathcal P_0(|\dot I|+|\dot\phi-\Omega(I)|)\|y\|_2^2
\le\mathcal P_0r^2\cdot r^2E^2w_\beta^{-2}\).
(v) On the nonzero-vertical low sites the word length exceeds
\(B^{K_0+3}\), so \(W_A=N^A\) there and \(b\ge1\); hence
\(\|Q_LY\|_2\le N^{-A}\|Y\|_{X_A}\) at each time, and
\(\|Q_LY\|_{\beta,2}\le N^{-A}\|Y\|_{\beta,X_A}\le EN^{-A}\).
(vi) \(|\Omega(I)-\Omega(I^*)|\le C|I-I^*|\le C(u|j|+\mathcal P_0r^2\|Y\|_2^2)\)
by \eqref{rate:eq:ward-DI}, and
\(|\Omega(I^*,0)-\Omega_*|\le C\||a_\infty|^2\|_1\le Cr^2E\).
For example \(V_w\) has an undifferentiated transverse factor, and its remaining packet factor has physical size \(rE\) in Wiener norm. The last line of \eqref{rate:eq:ward-torque} has physical size at most \(\mathcal P_0r^4w_\beta^{-1}\). On the full trial ball, \(\|Q_LY\|_{\beta,2}\le E N^{-A}\), by the original cap on the actual pre-mixed field. From \eqref{rate:eq:ward-DI}, \(|\Omega(I)-\Omega_*|\le \mathcal P_0(u|j|+r^2)\). The nondecaying \(r^2\) in this last bound is multiplied by transverse fields in \eqref{rate:eq:ward-col}. Dividing \eqref{rate:eq:ward-torque} by \(u^2\) consequently bounds its normalized right hand side by
\begin{equation}\label{rate:eq:ward-jbound}
 \|R_j\|_\beta\le \mathcal P_0\lambda^{-2}
          (N^{-J_0}+N^{-1}u\|j\|_\beta+r+N^{-2A}).
\end{equation}
The four terms arise as follows, with \(y=rY\) and
\(u^2=r^2\lambda^2\). By \eqref{rate:eq:ward-col} the high-high and
mixed blocks of \(\tfrac12\langle y,(\partial_{\phi_i}C)y\rangle\) are
at most \(r^2\mathcal P_0(N^{-J_0}+N^{-1}(u|j|+r^2))\|Y\|_2^2\), and
its low block is at most
\(r^2\mathcal P_0\|Q_LY\|_2^2\le r^2\mathcal P_0N^{-2A}w_\beta^{-2}\).
The terms \(V_{\phi_i}\) and \(\langle S,G_iy\rangle\) are
\(O(\mathcal P_0r^3w_\beta^{-1})\), and the last line of
\eqref{rate:eq:ward-torque} is \(O(\mathcal P_0r^4w_\beta^{-1})\).
Division by \(u^2\) turns the factors \(r^2,r^3,r^4\) into
\(\lambda^{-2},\lambda^{-2}r,\lambda^{-2}r^2\).

In the angular equation retain the matrix \(D_I\Omega_*\) evaluated at \((I^*,|a_\infty|^2)\), differentiating the occupied actions at those fixed Cartesian endpoint actions. Then \(\theta'=(D_I\Omega_*)j+R_\theta\), where
\[
 uR_\theta=
 \partial_I\mathcal K_{\mathrm{ax}}(I,|a_{\mathrm{ax}}|^2)-\Omega_*
           -u(D_I\Omega_*)j
 +\tfrac12\langle w,B_I w\rangle+V_I .
\]
In \(\ell^2\) one has \(\||a_{\mathrm{ax}}|^2-|a_\infty|^2\|\le Er^2 w_\beta^{-1}\). Taylor expansion of the analytic mass-subtracted frequency, \eqref{rate:eq:ward-DI}, and the bounded derivative \(B_I\) yield
\begin{equation}\label{rate:eq:ward-thbound}
 \|R_\theta\|_\beta\le \mathcal P_0(r/\lambda+u\|j\|_\beta^2).
\end{equation}
Indeed \(I-I^*=uj+O(\mathcal P_0r^2\|Y\|_2^2)\) by
\eqref{rate:eq:ward-DI}. The first-order Taylor term of
\(\partial_I\mathcal K_{\mathrm{ax}}\) in \(I\) therefore cancels
\(u(D_I\Omega_*)j\) up to \(O(\mathcal P_0r^2w_\beta^{-1})\), the
second-order term is \(O(\mathcal P_0(u|j|+r^2w_\beta^{-1})^2)\), the
dependence on \(|a_{\mathrm{ax}}|^2\) contributes
\(O(Er^2w_\beta^{-1})\), and
\(\tfrac12\langle w,B_Iw\rangle+V_I=O(\mathcal P_0r^2w_\beta^{-1})\).
Division by \(u=r\lambda\) gives \(r^2/u=r/\lambda\) and
\((u|j|)^2/u=u|j|^2\).
An isolated \(\Omega(I^*,0)-\Omega_*\) in place of this subtraction would give a constant source at infinity.

\subsubsection{A compact self-map with unequal radii}
On the finite strip use the real closed ball
\[
 \|e\|_{\beta,X_A}\le\rho,\qquad
 \|j\|_\beta,\|\theta\|_\beta\le\rho_o .
\]
For each trial recover \eqref{rate:eq:ward-I} and evaluate the two preceding right hand sides. With \(\mathcal I_L f(\tau)=-\int_\tau^L f(s)\,ds\), put
\[
 j_{\mathrm{new}}=\mathcal I_L R_j,\qquad
 \theta_{\mathrm{new}}=\mathcal I_L[(D_I\Omega_*)j_{\mathrm{new}}+R_\theta].
\]
Using the new value in this linear shear performs two successive integrations. Since \(\|\mathcal I_L\|\le T+\beta^{-1}\le E\), \eqref{rate:eq:ward-jbound}--\eqref{rate:eq:ward-thbound} give \(o(\rho_o)\) in both coordinates, uniformly in \(L\).
In detail, for \(\tau<T\) one has
\(w_\beta(\tau)\int_\tau^\infty w_\beta(s)^{-1}ds\le T+\beta^{-1}\), and for
\(\tau\ge T\) the same integral is \(\beta^{-1}\); with \(T=D_g\) and
\(\beta^{-1}\le8/\kappa\) this is at most \(E\) by (Pr1). On the ball,
\[
 \|j_{\mathrm{new}}\|_\beta\le E\mathcal P_0\lambda^{-2}
 \bigl(N^{-J_0}+N^{-1}u\rho_o+r+N^{-2A}\bigr).
\]
The prefactor is at most \(N^{a/10}\) by the choice of the prices.
Inside the bracket \(N^{-J_0}\), \(N^{-1}u\rho_o\) and \(N^{-2A}\) are
at most \(N^{-1}\rho_o\), and \(r=HN^{-s}\le N^{-2}\). Thus
\(\|j_{\mathrm{new}}\|_\beta\le4N^{a/10}(N^{-2}+N^{-1}\rho_o)\), and, with
\(\rho_o=N^{-3a}\) and \(3a+a/10<1\), even
\(E^2\|j_{\mathrm{new}}\|_\beta\le4N^{2\gamma+a/10}(N^{-2}+N^{-1}\rho_o)
=o(\rho_o)\), since \(E\le N^\gamma\). The matrix \(D_I\Omega_*\) is the
compression to \(F\) of \(D\chi(I^*,|a_\infty|^2)\), so
\(\|D_I\Omega_*\|_{\ell^2(F)\to\ell^2(F)}\le C\le E\). It follows that
\[
 \|\theta_{\mathrm{new}}\|_\beta\le E\bigl(E\,\|j_{\mathrm{new}}\|_\beta
 +\mathcal P_0(r/\lambda+u\rho_o^2)\bigr)=o(\rho_o),
\]
because \(r/\lambda=HN^{-s}e^{\epsilon B}\le N^{-2+\epsilon+o(1)}\) and
\(u\rho_o^2\le\rho_o^2\).

Invert the field-only principal problem at the \emph{trial} occupied parameters, retaining the entire compressed bounded low operator on the left. The same-copy exterior stars may keep their actual old coefficients and both signed adjoints, so the exterior on each copy is Hermitian in the fixed flattened norm. A star kept at its nominal parameters instead has a coefficient difference bounded by \(r\mathcal P_0(|\theta|+|I-I^*|)\). The simultaneous field-only Green absorbs the nominal weak tail and the remaining small homogeneous rows by its \(E\) price. Denote this inverse by \(\mathcal G_{\mathrm{tr}}\).

The low factor in the construction norm requires particular care. On all nonzero-vertical low sites both \(W_A\) and \(W_{A'}\) are saturated, as their word lengths exceed \(B^{K_0+3}\). For each prescribed vector \(Z\) one therefore has
\[
 \|Q_L Z\|_{X_A}\le N^{A-A'}\|Z\|_{X_{A'}}.
\]
This pays \(1/u\) for the bounded low operator on known vectors, including its moving parameter derivatives. Its restriction preserves the low sector, where the output cap is constant. On a free field the commutator with a high-high localized mixed primitive vanishes at first order, and a primitive meeting the low body has norm \(Er/R\). Thus its normalized first commutator costs \(E/(R\lambda)\), also when the bounded low operator is evaluated at the trial. Further Lie commutators contain an extra small primitive. In contrast the low body on the arbitrary trial correction itself is never estimated by \(N^{-A}\) in \(X_A\).

We record the remaining parameter-to-field estimates. For the high action-only scalar and the Cartesian axial diagonal,
\(\delta I/u=O(|j|+\mathcal P_0r/\lambda)\) and
\(\delta |a_{\mathrm{ax}}|^2/u=O(E r w_\beta^{-1}/\lambda)\), where the latter difference is from the exact endpoint profile. The connection involving \(\dot I,\dot\phi-\Omega(I)\) costs \(\mathcal P_0r/\lambda\) on a transverse field. The mismatch from the frozen velocity in the high scalar costs \(\mathcal P_0N^{-1}(u|j|+r^2)/u\) there. On a persistent axial carrier take the full difference from its true terminal frequency. These terms have combined weighted cost \(\mathcal P_0(\rho_o+r/\lambda+N^{-J_0}/u)\). The change of the true one-head or two-known cubic coefficient on a known field is bounded in normalized units by \(\mathcal P_0(|\theta|+|I-I^*|)\), summing the full angle series. When it also acts on a correction its norm has an extra factor \(\rho\). Terms at least quadratic in the corrections cost \(\mathcal P_0\rho^2\), by the \(X_A\) product inequality and the physical \(r\) on further input factors. Every source has a decaying transverse factor, an axial action difference, or a decaying occupied variation. The residual vanishes at the exact axial endpoint.

\smallskip
\noindent\emph{These estimates in detail.} All terms below are physical
field terms divided by \(u=r\lambda\), the normalization of the field
equation in slow time. (i) By \eqref{rate:eq:ward-DI},
\(I-I^*=uj+O(\mathcal P_0r^2\|Y\|_2^2)\), so \(\delta I/u=j+O(\mathcal P_0r/\lambda)\);
and \(\||a_{\mathrm{ax}}|^2-|a_\infty|^2\|\le Er^2w_\beta^{-1}\) gives the
second bound. A high action-only scalar and the axial diagonal have
bounded action derivatives (the first by the price \(\mathcal P_0\) of
an action derivative, the second because \(\chi\) is Lipschitz on
\(\ell^2\)), so they change by \(\mathcal P_0\) times these quantities.
(ii) The connection
\(T^{-1}[(\partial_IT)\dot I+(\partial_\phi T)(\dot\phi-\Omega(I))]\) has
physical size \(\mathcal P_0(|\dot I|+|\dot\phi-\Omega(I)|)\le\mathcal P_0r^2\)
by the velocity bounds above, hence \(\mathcal P_0r/\lambda\) after
division by \(u\). (iii) The frozen-velocity mismatch
\((\Omega(I)-\Omega_*)\cdot\partial_\phi b_3\) is at most
\(E_c'N^{-1}|\Omega(I)-\Omega_*|\le\mathcal P_0N^{-1}(u|j|+r^2)\), that
is \(\mathcal P_0N^{-1}(|j|+r/\lambda)\) after division. (iv) On a
persistent axial carrier the frequency difference
\(\omega_K(I)-\omega_K(I^\infty)\) is at most \(C\|I-I^\infty\|_1\), and is
paid as in (i). With \(\|j\|_\beta\le\rho_o\), (i)--(iv) sum to
\(\mathcal P_0(\rho_o+r/\lambda)\), and the residual of the spatial
columns adds \(\mathcal P_0N^{-J_0}/u\). (v) The coefficients of the
one-head and two-known cubic terms are analytic in the old angles and
actions with bounded derivatives on the strip, so moving
\((I,\phi)\) by \((I-I^*,\theta)\) changes them by
\(\mathcal P_0(|\theta|+|I-I^*|)\) after summing the angle series; on a
correction of size \(\rho\) this acquires the factor \(\rho\). (vi) The
product inequality in \(X_A\) costs \(\mathcal E\) per product, and every
further input carries the physical factor \(r\), so terms at least
quadratic in the corrections are at most \(\mathcal P_0\rho^2\).

Apply \(\mathcal G_{\mathrm{tr}}\) to the remaining field source evaluated at the trial to define \(e_{\mathrm{new}}\). Combining these estimates and the nominal field budget gives
\[
 \|e_{\mathrm{new}}\|_{\beta,X_A}
 \le E \mathcal P_0\big(rN^2+N^{-1}+r/\lambda+N^{-J_0}/u^2
                 +N^{A-A'}/u+\rho_o+\rho^2\big)+o(\rho)
 \le \rho/2
\]
at sufficiently large scales. In particular
$rN^2=N^{-(s-2)+o(1)}$ and $s-2\ge16a$. Also
$N^{A-A'}/u=N^{-2s-20}\lambda^{-1}/H$ and
$N^{-J_0}/u^2=N^{-J_0+2s}\lambda^{-2}/H^2$.
The choice of $\epsilon$ in Section~\ref{rate:sec:parameters}
absorbs the fixed powers of $\lambda^{-1}$. Among the two terms
involving the radii, $\rho_o=\rho^3$ and $\rho^2$ each
leave at least one extra factor $\rho$ before division by the field radius.
These bounds are uniform in the right endpoint $L$.

Fix \(N,L\) and \(0<\delta<1\). Take \(e\in C([0,L];X_{p-\delta,A})\)
and \(j,\theta\in C([0,L];\mathbb R^F)\) with the three weighted bounds
above imposed at each time in their stated norms. With the physical
conjugation symmetry on the field, these trials form a closed convex
set in the underlying real space with the uniform path topology.
Indeed, finite Fourier truncation gives lower semicontinuity of the
pointwise \(X_{p,A}\) norm in \(X_{p-\delta,A}\). On this interval
Duhamel's formula in the free interaction frame gives continuous
dependence in $L^2$ of the low propagator on its bounded coefficient.
The other inverses are given by uniformly convergent Neumann
series. Interpolation with the uniform bound at order $p$ gives
continuity of the map at order $p-\delta$ as well, since the
weight $W_A$ is bounded above and below at this fixed scale.

Its image is relatively compact. At each time the spatial embedding
is compact since $b(k)\to\infty$. The equation bounds the time
derivatives of the image at order $p-2$ on the fixed interval.
Interpolating differences between two times gives equicontinuity at
order $p-\delta$, so Arzel\'a--Ascoli applies. Schauder's theorem
gives a solution on each interval $[0,L]$. As $L\to\infty$, the
equation gives the same local time-derivative bound on each
\([0,L']\) for all \(L\ge L'\), using the uniform top-order estimate.
Compactness on these intervals yields a diagonal subsequence.
In Duhamel's formula its limit satisfies the full
equation, and the stable trace passes to the limit at the fixed time
$T$. The exponential bound passes to each finite time, giving
the required axial limit as time tends to infinity. Local Sobolev
uniqueness identifies this solution with the full NLS trajectory
through its time-zero value.

\paragraph{Support of the undivided low operator.}
The projection \(Q_L\) here is taken in the completed pure old
coordinates \(y\). First form the entire transformed Hessian \(C\) in
those coordinates. Its low same-copy bounded part \(L_L=Q_L L_L Q_L\)
is Hermitian by compression there. The opposite-copy row and the
high-to-low remainder are included with their adjoints in the field
budget above. The fixed diagonal rotations commute with \(Q_L\) and
with the flattened weights, and the finite averaging on the marked
graph misses this low sector. Extend \(Q_L\) by zero on the axial field
coordinates and put \(Q_H=\mathrm{Id}-Q_L\) on the field space.
For a subsequent bounded mixed primitive
write \(K=Q_HKQ_H+K_L\), where all the rows \emph{and columns} in
\(K_L\) have total operator bound \(Er/R\). Irrespective of whether
\(\exp K\) is unitary, the exact identity
\[
 e^{-K}L_Le^K-L_L
   =\int_0^1e^{-\xi K}[L_L,K_L]e^{\xi K}\,d\xi
\]
costs \(\mathcal P_0/(R\lambda)\) after division by \(u\). No time derivative
of the trial \(L_L\) is involved. For the unbounded free operator one
instead uses the combined commutator and time derivative of the
primitive: for the equation \(i\partial_t v=D_0v+\cdots\) its exact
integrand is \(e^{-\xi K}([D_0,K]-i\partial_tK)e^{\xi K}\).
These identities first hold on smooth vectors with the prescribed
finite or absolutely summable translations, then by closure in the
equation. The free commutator alone is never priced as a bounded
potential.

A Hermitian star connecting low and high sites has no such support
cancellation with \(Q_HKQ_H\). Its physical size \(Er\), times the
near primitive \(Er\lambda N^2\), divided by \(u\), costs
\(E^2rN^2\). This is one place requiring \(s>2\). The fast primitive
costs at most \(r\mathcal P_0\) in this product. These terms are errors even
if the full conjugating transformation fails to be unitary. Finally
the connection \(T^{-1}(T_I\dot I+T_\phi(\dot\phi-\Omega(I)))\)
costs \(\mathcal P_0r^2\) in physical units even on a free low field. Its
actual parameter derivatives retain the small velocities or the
quadratic field dependence of their differences from the axial
frequencies. The old torque was evaluated before these mixed changes;
its two kinetic terms cancel in the exact formula
\eqref{rate:eq:ward-torque}. In the real slow-time equation for the
mixed field variation, the contribution from varying the low Hessian
coefficient \(L_L(I,\phi)\) is
\[
 \frac1u U(t)^{-1}{\mathbb J}(\delta L_L)Q_L
       \bigl(Z_{\mathrm{kn}}(t)+U(t)e\bigr).
\]
Here \(Q_L\) is extended by zero on the axial coordinates, and \(U(t)\)
is fixed under the variation. The term in which the field itself varies
belongs to the conjugated low operator treated above. The two-cap
estimate controls the known part \(Q_LZ_{\mathrm{kn}}\) in \(X_A\).
On the installed trajectory the full projected field is bounded by
\(EN^{-A}\) in the ordinary variational norms used below.

\subsubsection{The linear inverse at the installed state}
Linearize along the full solution just constructed, in the ordinary
scaled spaces \(\ell^2(b^q)\), \(q=0,2,s\), at each of the two time
weights. Write \(h\) for the field variation and
\(z=(\delta j,\delta\theta)\) for the occupied action and angle
variations. With forcing \(f,g\) and zero stable trace, the equations
have the form
\[
 h=\mathcal G_f(f+A_o z),\qquad
 z=\mathcal G_o(g+R_f h+R_o z).
\]
Here \(\mathcal G_f\) inverts the field equation with the occupied
variables held fixed, and \(\mathcal G_o\) integrates the old shear
backward. The operator \(A_o\) maps occupied variations into the
field equation. The operators \(R_f,R_o\) give the remaining terms
in the occupied equations.

Eliminating \(h\) leaves the operator
\(\mathcal G_o(R_o+R_f\mathcal G_f A_o)\) on \(z\). Indeed,
substituting the first equation into the second,
\[
 z=\mathcal G_o\bigl(g+R_f\mathcal G_ff\bigr)
   +\mathcal G_o\bigl(R_o+R_f\mathcal G_fA_o\bigr)z .
\]
Thus the
estimate we need is
\[
 \|\mathcal G_o(R_o+R_f\mathcal G_f A_o)\|<\tfrac12 .
\]
Thus even a large \(A_o\) is admissible if the feedback through
\(R_f\) is sufficiently small.

\emph{The field equation and its dependence on the occupied variables.}
When one input in a chart derivative or a cubic term varies, bound the other inputs in their actual \(b^q\)-Wiener norms. The installed field correction contributes at most \(E\rho\) to each such multiplier. Write \(A_{\mathrm{ax}}=a_{\mathrm{ax}}/r\) on this trajectory. For a frequency derivative known only on \(\ell^2\), use
\[
 \|b^q A_{\mathrm{ax}}\,D\chi[\delta I]\|_2
 \le \|b^q A_{\mathrm{ax}}\|_\infty\|D\chi\|_{2\to2}\|\delta I\|_2 .
\]
The field-field terms induced by the correction are bounded by
\(\mathcal P_0\rho\). Differentiating an action-only diagonal through
\(I_e\) contributes a further \(\mathcal P_0r/\lambda\).

In the low block the actual field before the mixed change has
\(\ell^2(b^q)\) norm at most \(EN^{-A}\), by its \(X_A\) bound.
The \(I_e\) contribution there is therefore at most
\(\mathcal P_0(r/\lambda)N^{-A}\). Differentiating in \(\theta\) or
\(j\) gives bounds \(\mathcal P_0N^{-A}/u\) and \(\mathcal P_0N^{-A}\),
respectively. In the high block the derivative of the frequency
mismatch \(\Omega(I)-\Omega_*\) costs \(\mathcal P_0N^{-1}\) on
the scaled \(j\) input and \(\mathcal P_0N^{-1}r/\lambda\) on a
field input. All other terms in \(A_o\) cost at most \(\mathcal P_0\).
The rotations and mixed changes use fixed frequencies, so their
state derivatives introduce no factor \(t\).

Absorb the nominal small tail first, and then the field perturbations
bounded by \(\mathcal P_0\rho\) and by the estimates for \(I_e\).
This gives \(\|\mathcal G_f\|\le E\), enlarging \(E\) within the
allowed bounds. The backward inverse for the old shear also has
norm at most \(E\). Collecting the occupied derivatives above gives
\(\|A_o\|\le\mathcal P_0(1+N^{-A}/u)\).

\emph{Feedback into the occupied equations.}
Consider the first term of \eqref{rate:eq:ward-torque}. A free low
variation pairs with an actual low factor at cost
\(\mathcal P_0\lambda^{-2}N^{-A}\). Indeed, at fixed \(I,\phi\) set
\(C_{i,L}=Q_L(\partial_{\phi_i}C)Q_L\). Since \(y=rY\), \(u=r\lambda\)
and \(C_{i,L}\) is selfadjoint, the field derivative of the normalized
low quadratic term in the direction \(h_Y=\delta Y\) is
\[
 -\lambda^{-2}\langle Q_L h_Y,C_{i,L}Q_LY\rangle .
\]
The factor \(N^{-A}\) comes from \(\|Q_LY\|_{\beta,2}\le EN^{-A}\)
on the installed trajectory. The variation satisfies
\(\|Q_L h_Y\|_2\le\|h_Y\|_2\) in the ordinary variational space.
The high part costs
\(\mathcal P_0\lambda^{-2}[N^{-J_0}+N^{-1}(u\rho_o+r^2)]\).
Differentiating its coefficient in a field direction inserts \(I_e\)
and costs a further \(\mathcal P_0\lambda^{-2}N^{-1}r^2\) in the high
mismatch. A free \(j\) derivative there costs
\(\mathcal P_0\lambda^{-2}N^{-1}u\). An old angle derivative
preserves the small high mismatch or the two low factors.

The mixed terms involving \(V_\phi\) or \(S\) retain physical size
\(\mathcal P_0r^3\) under normalized field differentiation. In the
last term of \eqref{rate:eq:ward-torque}, the velocities
\(\dot I,\dot\phi-\Omega(I)\) before the Ward correction and their
normalized field derivatives have size \(\mathcal P_0r^2\). In the
angle equation the field and \(\theta\) derivatives, as well as the
error after subtracting the retained \(j\)-shear from the \(j\)
derivative, cost at most \(\mathcal P_0(r/\lambda+u\rho_o)\).
Together these bounds give
\(\|R_f\|,\|R_o\|\le B_{\mathrm{ret}}\), where
\[
 B_{\mathrm{ret}}=\mathcal P_0\big[\lambda^{-2}(N^{-A}+N^{-J_0}+N^{-1}u+r)
                        +r/\lambda+u\big].
\]

\emph{Solving the coupled equations.}
The norm of the operator left after eliminating \(h\) is at most
\[
 \|\mathcal G_o\| B_{\mathrm{ret}}
        (1+\|\mathcal G_f\|\,\|A_o\|)<\tfrac12 .
\]
Indeed $N^{-A}/u=N^{s-A}\lambda^{-1}/H$, and each term in
$B_{\mathrm{ret}}$ retains a negative power of $N$ after the fixed
products of $E,\mathcal P_0$ and $\lambda^{-1}$ are accounted for
as in Section~\ref{rate:sec:parameters}. The Neumann series determines
$z$, and the first equation then determines $h$. A prescribed nonzero
stable trace adds its field solution to the first right-hand side
and gives the same inverse estimate. Both time weights and the three
spatial orders impose the same stable trace. Solutions constructed
on their common domain coincide by uniqueness at the weaker weight
and order. All coefficients in these equations are evaluated along
the full solution obtained above.

To return to physical variables, at each time let \(\mathcal F_t\) denote the instantaneous coordinate map used above. If \(\dot v=X(v)\) and \(\dot\xi=DX(v)\xi+f\), then \(k=D\mathcal F_t(v(t))\xi\) solves the transformed variational equation with forcing \(D\mathcal F_t(v(t))f\), with the additional constant factor for slow time. This follows by the chain rule on smooth inputs and Duhamel continuation at fixed scale. All coordinate differentials and their inverses have fixed polynomial prices in \(N\) on the three scaled spaces by the translated chart estimates, the occupied floors and \eqref{rate:eq:ward-DI}. Also \(\|f\|_{\ell^2(b^q)}\le C_q\|f\|_{H^q}\le C_q N^q\|f\|_{\ell^2(b^q)}\). Thus the physical bounds remain fixed powers of \(N\). If \(\Pi_{\mathrm{st}}\) is the complete finite stable trace in the transformed coordinates, the physical trace is
\[
 P\xi=\Pi_{\mathrm{st}}D\mathcal F_{T/u}(v(T/u))\,\xi(T/u).
\]
It includes both real stable quadratures of each complex pair. The unrestricted homogeneous unaugmented physical operator can have a nonzero decaying kernel, for example \(\xi=\partial_y v\) when the endpoint is axial and the solution is not. The invertibility just proved concerns precisely the displayed augmented trace.

\subsection{Physical comparison and iteration}
All comparisons before the very last limit use the autonomous
mass-gauged equation. The transformation back to the original
equation is a scalar phase depending only on the solution's mass.
For an exponent $q$ let $\|d\|_{q,\beta}$ denote
$\sup_{t\ge0}e^{\beta t}\|d(t)\|_{H^q}$ in physical time.

\subsubsection{The three estimates on each raw input}
Let $W_j$ be the raw field at a stage and $V_j$ its axial endpoint.
When comparing $V_{j-1}$ and $V_j$, all occupied preexisting
rates $n^2+\chi_n$ and angular origins are locked in this gauge.
Compare actual old angles to
$\phi_{j-1,0}+t\Omega_{j-1}$ with the common time orientation.
For each earlier site the reference angle is therefore the same on
both tori for every $t$. Along $W_j$ its displacement is precisely
$\Theta=\phi-\phi_*(t)$ from \eqref{rate:eq:ward-I}, whose bound also
holds on the initial interval. A frequency error, even a small
constant one, would instead give an angle error growing linearly
in $t$.

There are constants $C_0,C_1$ independent of the stage and the
requested prefix precision such that, for $q=2,s$,
\begin{equation}\label{rate:eq:raw-physical}
 \beta_j\ge e^{-C_0B_j},\qquad
 \|W_j-V_j\|_{q,\beta_j}+\|G_j\|_{q,\beta_j}\le e^{C_0B_j}.
\end{equation}
Here $G_j$ is the inhomogeneous physical variational right inverse
along $W_j$ with the same augmented stable data at the section in
both orders. To convert the section estimate to the physical
norm, $\|f\|_{b^q}\le C_s\|f\|_{H^q}\le C_sN^q\|f\|_{b^q}$.
The instantaneous coordinate differential is composed of the
axial Birkhoff coordinate differential and the polar map on the
occupied coordinates, the triangular action correction, the
spatial frame and the small mixed exponential, fixed unitary
coordinate rotations, and the scale multipliers $r^{\pm1},u^{\pm1}$.
Both directions in ordinary $\ell^2(b^q)$ cost a fixed power of $N$.
Indeed the known inputs in the forward and the inverse physical
chart expansions have scaled Wiener norms bounded by $r\mathcal E$
about their pure old values, the old positive floors cost at most
$N^{C_s}$, and no radial derivative of a hot Cartesian coordinate
is divided by its modulus. This conversion applies to the
inhomogeneous forcing after differentiating the exact coordinates
at $W_j$; the physical-time source introduces an additional
$1/u$. The same stable trace at a single physical time prescribes
the boundary condition. Converting the early/late weight costs
$\exp(\sigma D)$, and one may take
$\beta_j=u\kappa/8$. These losses in total have logarithm at
most $C_0B_j$, once prices $\mathcal E\le N^\gamma$ and
$\mathcal P\le N^{a/10}$ have been absorbed.
In detail, the factors are: \(N^q\le N^s\) from the section norms; a
fixed power \(N^{C_s}\) from the coordinate differentials, since the
smallest positive action at the final point satisfies \(I_{\min}^{-1}\le C N^{2s+4}\),
\(r^{-1}=N^s/H\le N^s\) and
\(u^{-1}=N^s\lambda^{-1}/H\le N^{s+1}\) from the scale multipliers and
the physical-time source; \(e^{\sigma D_g}\le e^{D_g}\le E\) from the
early/late weight, since \(e^{\beta_jt}=e^{\kappa\tau/8}\le
e^{\kappa T/8}w_{\kappa/8}(\tau)\) with \(\tau=ut\) and \(T=D_g\); and
the normalized inverse norms and field bounds, at most
\(\mathcal E\mathcal P\le N^{a}\). Moreover \(\beta_j=u\kappa/8\ge
N^{-s-1}\kappa/8\ge N^{-s-2}\). As \(\log N_j\le2B_j\), the logarithm of
the product is at most \(C\log N_j\le2CB_j\), with \(C\) depending only
on \(s\) and the fixed caps; any \(C_0\ge2C+2s+4\) serves. Exponential landing
uses the same coordinate conversions for the finitely many
parameters and the decaying field components. The stationary
Cartesian endpoints have their true terminal phases. In particular
one does not integrate a pointwise frequency approximation along
the half-line.

Here is the sharper budget for the unweighted raw strong profile.
Set $m=\lceil s\rceil-1\ge2$. The one-dimensional inverse chart
$\Psi=\mathcal F^{-1}+G$ has a uniformly Lipschitz remainder
from a small $h^m$ ball to $H^{m+1}$; take its center at zero
and shrink the ball using continuity of the derivative there.
This uses only an integer $m$ even for real $s$. Every retained
current carrier satisfies $|k|\le C dN$ and their normalized
Euclidean size is $e^{Cd}$. Thus their $h^m$ cost is at most
$H N^{m-s}e^{C_s d}$. Before the mixed coordinate change, the weighted field error of size
$N^{-a}$ has physical $H^q$ cost at most $C_sH N^{q-s-a}$,
$q\le s$. Each mixed change and particular on the known field
adds a normalized $\mathcal E(r\lambda N^2)+\mathcal P r$
or $\mathcal E/R$ in the scaled strong norm, and its action on
the unknown costs at most $\mathcal E$. Thus the physical
reconstruction error from them is at most $H N^{-a'}$ for
some fixed $a'>0$. The finite-only shadow and hot-action
rounding cost at most $C_sH d^{C_s}e^{-C_2dD}$ in physical
$H^s$, since they are supported on the marked planes.
In detail, a scaled field \(f\) normalized by \(r\) has physical
\(H^q\) norm at most \(rC_sN^q\|f\|_{b^q}=C_sHN^{q-s}\|f\|_{b^q}\), since
\(\bk k\le CN\,b(k)\). Every carrier has \(\bk k\le CdN\), which gives
the bound \(HN^{m-s}e^{C_sd}\) for their \(h^m\) cost, and the field
error of normalized size \(N^{-a}\) costs \(C_sHN^{q-s-a}\). The mixed
changes and particulars on the known field have normalized size
\(\mathcal E r\lambda N^2+\mathcal Pr+\mathcal E/R\le
N^{-\min(s-2,1)+2\gamma+a/10}\), because \(r\lambda N^2\le N^{2-s}\),
\(r\le N^{-s+o(1)}\) and \(R^{-1}=P/N\); their action on the unknown,
of size at most \(N^{-a}\), costs \(\mathcal EN^{-a}\). Hence the
reconstruction error is at most \(HN^{-a'}\) with
\(a'=\min(a,s-2,1)/2\), since \(a<\min(s-2,1/4)/8\) and
\(\gamma\) is much smaller. The shadow and the rounding have
normalized size at most \(e^{-C_2dD}\) on the at most \(Cd^2\) marked
planes, each of which has physical weight at most
\(r(CdN)^s=C^sd^sH\).

Let $F_o$ be the old support and $Z_o=\max_{n\in F_o}|n|$.
At an old site $n$, split the change in $\sqrt{I_n}e^{-i\phi_n}$
into a change of radius and a change of angle. The bounds
$|\sqrt{x}-\sqrt{y}|^2\le |x-y|$ and
$|e^{it}-1|\le |t|$ show that the $\ell^2$ norm of the radial
change is at most $\|I-I_o\|_{\ell^1(F_o)}^{1/2}$.
Cauchy--Schwarz on the finite set $F_o$ gives the following
estimate for the total difference on the old sites in $h^q$ for $0\le q\le s$
\[
 C\langle Z_o\rangle^q
 \bigl[|F_o|^{1/4}\|I-I_o\|_2^{1/2}
    +\|(\sqrt{I_{o,n}}|\Theta_n|)_{F_o}\|_2\bigr].
\]
In detail, at each old site
\[
\begin{aligned}
 \bigl|\sqrt{I_n}e^{-i\phi_n}-\sqrt{I_{o,n}}e^{-i\phi_{o,n}}\bigr|
 &\le\bigl|\sqrt{I_n}-\sqrt{I_{o,n}}\bigr|+\sqrt{I_{o,n}}\,|\Theta_n|,\\
 \sum_{n\in F_o}\bigl|\sqrt{I_n}-\sqrt{I_{o,n}}\bigr|^2
 &\le\|I-I_o\|_{\ell^1(F_o)}\le|F_o|^{1/2}\|I-I_o\|_2,
\end{aligned}
\]
and, since $q\ge0$, the weight of $h^q$ on $F_o$ is at most $\langle Z_o\rangle^q$.
The true $I$ here includes the restored quadratic action correction
of physical order $r^2\mathcal P$; the raw old difference from
the retuned endpoint has size at most $uN^{-3a}+r^2\mathcal P$.
All changes of the old actions during the choice of the endpoint,
and the $h^s$ costs of the new clocks and spectator, can be
made $\le e^{-D'B_o}$ for any specified $D'$ after allowing
the indicated $O_{D'}(B_o)$ expenses in calibration. For the
old-coordinate estimate use this adjustable precision before
the square root and the $Z_o^s$ factor. The uniform $h^m$
ball follows by summing these bills and
$H_iN_i^{m-s}e^{C_s d_i}$, since $m<s$; take the two-site
base small first. Consequently
\[
 \|G(z)-G(z_o)\|_{H^s}\le C_m\|z-z_o\|_{h^m}.
\]
In particular one can use the physical strong error for the
nonlinear part of the axial reconstruction without any uniform
chart assumption on unbounded $h^s$ balls.

Let $\mathcal C_j(t)$ denote the actual marked plane waves of
the ideal finite reference in their true carrier rotations.
In $H^s$ the raw error relative to $V_{j-1}+\mathcal C_j$ is
bounded uniformly by
\[
 C\{e^{-D' B_o}+H N^{-a'}+
 H N^{m-s}e^{C_s d}+H d^{C_s}e^{-C_2 dD}\}.
\]
For the weak comparison keep the finite perturbation on its
marked frequencies, where its $H^2$ cost carries an additional
factor $N^{2-s}$. Alternatively bound all the actual marked
coordinates together in $H^2$ by $HN^{2-s}e^{Cd}$.
The one-smoothing comparison in $H^2$ can likewise use the
small $h^2$ ball and this last $h^2$ bill. It follows that for
arbitrarily prescribed $D_h$ one has
\begin{equation}\label{rate:eq:raw-weak}
 \sup_{t\ge0}\|W_j(t)-V_{j-1}(t)\|_{H^2}
       \le e^{-D_h B_{j-1}} .
\end{equation}
The old angular terms and the current negative powers of $N$ satisfy
this bound once $Q$ is sufficiently large. To compare the two axial
solutions themselves, fix $t$ and take $t_n\to+\infty$ such that the
pair of phase vectors at $t_n$ tends to the pair at $t$. Such times
exist by recurrence of a linear flow on the product torus, including
when its frequencies satisfy rational relations. By continuity,
\[
 \|V_j(t)-V_{j-1}(t)\|_{H^2}
 \le \limsup_{n\to\infty}
 \big(\|V_j(t_n)-W_j(t_n)\|_{H^2}
      +\|W_j(t_n)-V_{j-1}(t_n)\|_{H^2}\big)
 \le e^{-D_hB_{j-1}} .
\]
Thus summing the increments gives a common $H^2$ bound for all
\(V_j,W_j\). In \(H^s\) the marked waves cost at most \(C H_j\) at all
times, with \(H_j=e^{bd_j}\). The preceding reconstruction estimate gives
\[
 \sup_{t\ge0}\|W_j(t)-V_{j-1}(t)\|_{H^s}\le C H_j+C e^{-c d_j}.
\]
Since \(\|W_j(t)-V_j(t)\|_{H^s}\to0\), the same recurrence argument at
order \(s\) gives the bound \(C H_j+C e^{-c d_j}\) for
\(\sup_{t\ge0}\|V_j(t)-V_{j-1}(t)\|_{H^s}\). Summing the axial increments
up to \(j-1\) and then adding \(W_j-V_{j-1}\) yields
\[
 \sup_{t\ge0}\|W_j(t)\|_{H^s}\le C e^{C d_j}.
\]
For a coarse old $H^{s+1}$ cap, the finite-support axial
parameterization on a fixed angular strip has a uniformly
bounded physical $L^2$ norm. Its angular Fourier coefficient
has spatial index $l\cdot n_F$ and size $C e^{-a_0|l|_1}$.
Summing with its weight gives a bound
$C_s^{|F|+1}(1+Z)^{s+1}\le e^{C_0 B_j}$ for $V_j$, after
increasing the coarse constant.

Let $T_j=\exp[(C_0+1)B_j]$. The physical root of $V_j$ has
Fourier modulus in the $H^s$ normalization at least $0.95H_j$:
the true axial coordinate has modulus $(1+o(1))r$, and its
one-smoothing error at that frequency in this normalization
is $O((1+N)^{s-m-1})=O(1)$. The latter includes integer $s$.
By \eqref{rate:eq:raw-physical} the error $W_j-V_j$ at all times
$t\ge T_j$ is at most $\exp(C_0B_j-e^{B_j})$ in $H^s$.
It leaves a persistent root of normalized height $0.9 H_j$.
Choose $L_0>C_0+2$, $L_1>\max(2L_0+2,C_0+1)$ and put
$S_j=e^{L_0 B_o}$, $K_j=e^{L_1B_o}$. For the cold estimate
at $S_j$ compare the actual finite profile with its launch in
their individual rotating coordinates. The slow derivative costs
at most $\mathcal E/\lambda$, so the physical $H^s$ price of
this comparison is at most $C_sH d^{C_s}\mathcal E r S_j$,
in addition to its cold launch. The free phases all have modulus
one. The old axial tail costs $e^{C_0B_o}/K_j$. These bounds,
the leaf estimate in Lemma~\ref{rate:l:finite}, and the raw strong
error above give
\begin{equation}\label{rate:eq:cold-physical}
 \|P_{\langle k\rangle\ge K_j}W_j(S_j)\|_{H^s}
      \le C e^{-c d_j},\qquad T_{j-1}\le S_j<T_j .
\end{equation}
Here $Q$ is taken sufficiently large and eventually fixed before
choosing the starting scale. In particular $B_o\gg d_j$ at
large stages with this fixed $Q$.

\smallskip
\noindent\emph{Derivation of \eqref{rate:eq:cold-physical}.} Write
\(W_j(S_j)=V_{j-1}(S_j)+\mathcal C_j(S_j)+\mathcal R_j(S_j)\), where
\(\mathcal R_j\) is the raw error relative to \(V_{j-1}+\mathcal C_j\)
estimated above.
(a) By the coarse old cap,
\(\|P_{\bk k\ge K_j}V_{j-1}(t)\|_{H^s}\le K_j^{-1}\|V_{j-1}(t)\|_{H^{s+1}}
\le e^{C_0B_o}/K_j=e^{-(L_1-C_0)B_o}\).
(b) Every marked frequency has \(|k|\ge R-1\), whereas
\(K_j=e^{L_1B_o}\le R/2\) once \(Q>L_1+1\), because
\(\log R=QB_o+O(1)\). Hence all of \(\mathcal C_j(S_j)\) lies in
\(\bk k\ge K_j\). In its rotating coordinates the normalized finite
profile \(Z(\tau)\) satisfies \(|Z(\tau)-Z(0)|\le\tau\sup|Z'|\le
\tau\mathcal E/\lambda\), and at \(t=S_j\) the slow time is
\(\tau=r\lambda S_j\), so the change is at most \(\mathcal ErS_j\). A
normalized vector \(\eta\) on the at most \(Cd^2\) marked planes has
physical \(H^s\) norm at most \(C_sHd^{C_s}|\eta|\), because
\(r\bk k^s\le r(CdN)^s=C^sd^sH\). The rotations have modulus one. Thus
\(\|\mathcal C_j(S_j)\|_{H^s}\le\|\mathcal C_j(0)\|_{H^s}
+C_sHd^{C_s}\mathcal ErS_j\). By the leaf estimate after
Lemma~\ref{rate:l:round} the launch is at most
\(C_sH[(4/5)^s+(1/5)^s]^{(d-1)/2}+C_sHd^{C_s}2^{d/2}e^{-c_1D_g}
\le Ce^{-c_sd}\) for \(H=e^{b_*d}\) with
\(b_*<-\tfrac12\log[(4/5)^s+(1/5)^s]\). The comparison term equals
\(C_sH^2d^{C_s}\mathcal EN^{-s}e^{L_0B_o}\), whose logarithm is at most
\[
 2b_*d+(\gamma-s)\log N+L_0B/Q+O(\log d)
 \le2b_*d+(\gamma-s+L_0/Q)B+O(\log d)\le-B/2
\]
for \(Q>L_0\), \(\gamma<1/4\) and large \(B\): here \(\mathcal E\le
N^\gamma\), \(\gamma-s<0\) and \(\log N\ge B\) give the first
inequality, and \(\gamma-s+L_0/Q<\tfrac14-2+1=-\tfrac34\) together with
\(d=o(B)\) gives the second. In particular the comparison term is at
most \(e^{-B/2}\le e^{-d}\).
(c) Each term of the raw error bound is at most \(Ce^{-d}\):
\(e^{-D'B_o}\le e^{-d}\) as \(B_o\gg d\); \(HN^{-a'}\) and
\(HN^{m-s}e^{C_sd}\) are at most \(N^{-\min(a',s-m)/2}\) since
\(d=o(B)\); and \(Hd^{C_s}e^{-C_2dD_g}\le e^{-d}\) for \(C_2>b_*+2\).
Adding (a)--(c) gives \eqref{rate:eq:cold-physical} with \(c=\min(c_s,1)\).
Finally \(T_{j-1}=e^{(C_0+1)B_o}<e^{L_0B_o}=S_j\) because \(L_0>C_0+2\),
and \(S_j<T_j=e^{(C_0+1)QB_o}\) because \(Q>L_0\).

\subsubsection{Composition in the smaller Sobolev space}
At stage $j$, $W_j$ is close to the preceding axial solution
$V_{j-1}$. We must also restore the decaying parts of the earlier
two-dimensional trajectories. Write $X_i^{(j)}$ for the result of
restoring the earlier stages down to index $i$, starting with
$X_j^{(j)}=W_j$ and $X_{j+1}^{(j)}=V_j$. For
$i=j-1,j-2,\ldots,1$ the three known solutions are
\[
 U=X_i^{(j-1)},\qquad Y=X_{i+1}^{(j-1)},\qquad Y'=X_{i+1}^{(j)}.
\]
The old difference $d=U-Y$ decays at rate $\beta_i$, and the new
change $v=Y'-Y$ is small uniformly in $H^2$. Seek
$X_i^{(j)}=Y'+d+e$. Let $F$ be the mass-gauged cubic vector
field, including its factor $-i$. Subtracting the equations for
$Y',U,Y$ cancels the free evolution of $Y'+U-Y$. Subtract and
add $DF(W_i)e$ in the remaining equation and apply the inverse
$G_i$ along $W_i$. With zero stable trace for the additional
correction $e$, the result is
\[
 e=G_i\{F(U+v+e)-F(Y+v)-F(U)+F(Y)-DF(W_i)e\}.
\]
The lattice at stage $j$ is contained in that at stage $i$, with
the same congruence condition. All these fields therefore lie in
the domain of $G_i$.

In detail, $X_i^{(j)}=U+v+e$ and the four fields
$X_i^{(j)},Y',U,Y$ solve the same equation
$\partial_tX=\mathcal L_0X+F(X)$ with the free generator
$\mathcal L_0$, so that $e=X_i^{(j)}-Y'-U+Y$ satisfies
\[
 \partial_te-\mathcal L_0e-DF(W_i)e
 =F(U+v+e)-F(Y+v)-F(U)+F(Y)-DF(W_i)e .
\]
Let $\mathcal T$ be the real symmetric trilinear polarization of
$F$, including its mass term. The independent source is
$6\mathcal T(Y,d,v)+3\mathcal T(d,d,v)+3\mathcal T(d,v,v)$.
The other three terms are
$3\mathcal T(U+v-W_i,U+v+W_i,e)+3\mathcal T(U+v,e,e)
+\mathcal T(e,e,e)$. Indeed $F(x)=\mathcal T(x,x,x)$ and $U=Y+d$
give
\[
\begin{aligned}
 &F(Y+d+v)-F(Y+v)-F(Y+d)+F(Y)\\
 &\qquad=6\mathcal T(Y,d,v)+3\mathcal T(d,d,v)+3\mathcal T(d,v,v),\\
 &F(U+v+e)-F(U+v)\\
 &\qquad=3\mathcal T(U+v,U+v,e)+3\mathcal T(U+v,e,e)+\mathcal T(e,e,e),
\end{aligned}
\]
because the terms of degree at most one in $(d,v)$ and the pure
cubes cancel in the first line, and
$3\mathcal T(U+v,U+v,e)-DF(W_i)e
=3\mathcal T(U+v-W_i,U+v+W_i,e)$ by symmetry. On a common $H^2$ ball this gives a
weighted weak source of size $C A_i\delta$, where
$A_i=\|d\|_{2,\beta_i}$, $\delta=\sup\|v\|_{H^2}$.
Write $k_i=\max(1,\|G_i\|_{2,\beta_i},\|G_i\|_{s,\beta_i})$.
The weak ball of radius $R_e=2C k_i A_i\delta$ is invariant
and contractive if $C k_i(\sup\|U-W_i\|_{H^2}+\delta+R_e)<1/2$.
At its first use $A_i\le e^{C_0B_i}$. We allow twice this bound in
later stages. The new difference from $Y'$ is $d+e$, so its
weighted change from $d$ costs at most $R_e$. The difference from
$U$ is $v+e$, whose uniform $H^2$ norm is at most $\delta+R_e$.
Both are bounded by $e^{C_wB_i}\delta$, with $C_w$ depending only
on the coarse bounds at stage $i$.
In detail, \(H^2(\T^2)\) is an algebra, so
\(\|\mathcal T(a,b,c)\|_{H^2}\le C\|a\|_{H^2}\|b\|_{H^2}\|c\|_{H^2}\).
Every term of the independent source contains one factor \(d\), whose
weighted norm is at most \(A_i\), and one factor \(v\), of norm at most
\(\delta\), the remaining factor lying in the common \(H^2\) ball; this
is the weighted source \(CA_i\delta\). The term linear in \(e\) has
coefficient of norm at most \(C(\sup\|U-W_i\|_{H^2}+\delta)\), and the
quadratic and cubic terms are at most \(CR_e\|e\|\) on the ball. Hence
the image of the ball has weighted norm at most
\(k_i\bigl(CA_i\delta+C(\sup\|U-W_i\|_{H^2}+\delta+R_e)R_e\bigr)\le R_e\)
under the stated condition, and the same bound with \(\|e-e'\|\) gives
the contraction. By \eqref{rate:eq:raw-physical},
\(k_i\le e^{C_0B_i}\), and \(A_i\le2e^{C_0B_i}\), so
\(R_e\le4Ce^{2C_0B_i}\delta\) and
\(\delta+R_e\le e^{C_wB_i}\delta\) with \(C_w=2C_0+1\) once \(B_i\) is
large. The condition itself holds because
\(\sup\|U-W_i\|_{H^2}\le\sum_{j>i}\varepsilon_j\), \(\delta\) and
\(R_e\) are much smaller than \(e^{-C_0B_i}/(4C)\), as shown next.
The first change to be propagated
is $W_j-V_{j-1}$, bounded by \eqref{rate:eq:raw-weak}. Descending
from $j-1$ to $i$ therefore bounds the change at each row by
\[
 e^{-D_h B_{j-1}}\prod_{k=i}^{j-1}e^{C_w B_k}
 \le \varepsilon_j,\qquad
 \varepsilon_j:=e^{-(D_h-2C_w)B_{j-1}},
\]
since $\sum_{k<j} B_k\le2B_{j-1}$; indeed $B_k=B_*Q^k$ with $Q\ge2$
gives $\sum_{k\le j-1}B_k\le B_{j-1}Q/(Q-1)$.
At fixed old $i$ the
accumulated change of either $U-W_i$ or the weighted norm of
$d$ is at most $\sum_{j>i}\varepsilon_j$. It is much smaller
than $e^{-(C_0+1)B_i}$ if $D_h$ is chosen large enough.
Induction maintains the ball, its contraction, the common weak
bound and the registered $A_i$. Choose $D_h-2C_w$ also larger
than $(s-2)\max(L_1,C_0)+2$.

For the strong order use the same weak Picard iterates. The tame
product estimate for $\mathcal T$ puts the $H^s$ norm on one
factor and $H^2$ on the others, by Fourier convolution and
the mass scalar product. Terms with an $e$ carry the exponential
weight on an $e$. If it is also the strong factor in the term
linear in $e$, its coefficient contains $U+v-W_i$, which is
weakly small; if there are two or three $e$'s the extra one
is weakly small. If the strong factor is different, the weak
weighted $e$ costs $R_e$. In the independent source put the
weight on a $d$. Thus the coefficient of the previous strong
iterate can also be made at most $1/2$, and such iteration gives
\[
 \|d+e\|_{s,\beta_i}\le e^{C_sB_i}
 \big(1+\|d\|_{s,\beta_i}+\|Y'\|_{L^\infty H^s}
 +\|Y\|_{L^\infty H^s}+\|U\|_{L^\infty H^s}\big).
\]
Pointwise weak compactness in $H^s$ and convergence in the lower
norm give this same bound for the actual weak fixed point.
Uniqueness in $H^2$ and local $H^s$ persistence identify it with
the true strong solution on each finite interval. The logarithm
$L_j'$ of two plus the maximum of all strong sup norms and
weighted pair prefactors through stage $j$ satisfies
$L_j'\le\max(L_{j-1}',C B_j)+Cj+C\sum_{i<j}B_i$.
To iterate this estimate, bound every new maximum through stage
$j$ by $C B_j$ or by its initial value, and add the accumulated
increments. The terms $Ck$ sum to $O(j^2)$ and the double sum is
at most $2C\sum_{k\le j}B_{k-1}=O(B_j)$. Thus $L_j'\le C'B_j$.
The use of a single strong factor in each product is what gives
an additive recurrence for these logarithms.

Write $U_j=X_1^{(j)}$. Telescoping the row pairs gives
\[
 \|U_j(t)-W_j(t)\|_{H^s}
 \le j e^{C'B_j}\exp[-e^{-C_0B_{j-1}}t].
\]
For $t\ge S_j$ the logarithm of the right-hand side is at most
\[
 \log j + C'Q B_{j-1}-\exp\big((L_0-C_0)B_{j-1}\big)
 \le -B_{j-1}
\]
after increasing $B_*$ with $Q$ fixed. This transfers the tail
bound at $S_j$, the Fourier lower bound for $t\ge T_j$, and the
late $H^s$ upper bound from $W_j$ to $U_j$. The bound
$\sup_{t\ge0}\|U_j(t)-U_{j-1}(t)\|_{H^2}\le\varepsilon_j$ above
controls the effect of subsequent stages.

\subsection{Order of parameters and construction on successive lattices}\label{rate:sec:parameters}
Fix a single real $s>2$ and the integers $m,p$ above. Use
$a=\min(s-2,1/4)/16$ for the normalized field error; for
comparisons where $s-m$ also occurs keep that strictly positive
number as a separate exponent. Choose, for example, $A>4s+20$,
$A'=A+3s+20$ and a spatial and angular remainder order
$J_0>4(A'+s+p+10)$. These values satisfy the weaker inequalities
$A>2s+12$, $A'>A+3s+15$ and $J_0>2s+12$ under which the same
constants were used in the subsection on the two-dimensional
variables and in Section~\ref{rate:sec:ward-landing}, and the present
$a$ satisfies the condition $a<\min(s-2,1/4)/8$ imposed in the
subsection on scale, geometry and calibration. In the spatial recurrence and the angular
tail cuts take their finite internal orders higher as indicated
there. Choose an integer $h>4s+5$ for the axial template, then
the template, its length $q$, the number of hierarchy derivatives
and the finite-order coefficient expansion. All these numbers
depend on $s$ alone. The relative range of the two splitting
angles and the constants in Lemma~\ref{rate:l:finite} can be fixed
throughout.

Initialize with a two-site torus on, for instance, $1,6$. The
actions range over a small open box of two positive comparable
values. The direct rate Jacobian is close to $-\mathrm{Id}$
on this box; shrinking it and the surrounding analytic ball
can guarantee the smallness required by the disjoint template
rows. Choose the two direct rates inside a small interior cube
by the affine slicing argument. All phases can be fixed. This
initializes the prefix affine bound, positive margins for
$\mathcal Q-M^2$ and $2M^2-\mathcal Q$, small integer-order
$h^m$ bounds and finite high norms. With $B_*$ large it obeys
the first coarse maximum-site, least-action and inverse-price
bounds. The sum of all subsequent small action differences
prescribed in the preceding sections keeps the two strict mass
margins by $d\mathcal Q=4M\,dM+2\sum I_n\,d\chi_n$.

Here are the choices for the prices involving the current
logarithmic scale. All fixed powers of $\lambda^{-1}$ on the
current scale arise from the slow normalization, the fixed
current positive-floor/weak-radius derivatives, and the reciprocal
of the finite-dimensional direct-rate cube at the frozen module.
The cube there has an explicit $\lambda^2$ factor with the
other prices involving $B_o+d^2\log B$. There is no
current-scale inverse of the much smaller final direct-rate
cube used in the nonlinear PDE estimate. That cube is saved
for the next stage only. Fix the positive fractions of $\log N$
available in the estimates involving inverse powers of $\lambda$,
including $a/100$ for the combined coordinate and derivative prices.
Let $C_s$ bound the total exponent of $\lambda^{-1}$ in each required
product, and choose $\epsilon>0$ so that $2C_s\epsilon$ is smaller
than each of these fixed fractions. Thereafter choose once the weak letter
cut $m_*$ so large that the absolute tail past that length,
using width $\epsilon B/6$ in the weak coordinates and
ordinary angle width on the rest, is $O(N^{-2J_0})$ after
the fixed normalizations. For an unmarked homogeneous
charge-one nonsingleton the width bound already pays
$15\epsilon B/6$ against a single divisor $\lambda^{-1}$
from the field normalization. Thus the fixed positive exponent
$h_w$ in its bound $\mathcal E N^{-h_w}$ may be much smaller
than $a$. That bound is used with a field-only inverse price
$\mathcal E$ before the upper triangular old-coordinate arrows.
The weak cutoff affects the polynomial counts/degrees of geometry
exclusions, the required spectator intervals, the fixed numerical
constants in the old-and-macro price and the lengths of bounded
sums, but never the exponent of an extra unabsorbed inverse weak
floor per weak letter.

For clarity about induction, register the following coarse caps
at stage $j$, now taking $I_{\min}$ over the final occupied actions:
\[
 \log(2+\max|n|)\le 2\log N_j+O(1),\quad
 -\log I_{\min}\le(2s+4)\log N_j+O(1),\quad
 f_j\le C_s\sum_{i\le j}(1+d_i)^2 .
\]
Register the final affine estimate with exponential prefactor
$e^{-C_{\mathrm{reg}} B_j}$ and polynomial exponent $2f_j+6$,
as in the direct-rate selection. If all current designated
products of $\mathcal E$ fit within $N^{\gamma}$, with
$\gamma\ll\min(a,h_w,\epsilon)$, then the conclusions
\eqref{rate:eq:raw-physical} and the old $H^{s+1}$ bound hold with
a coarse $C_0$ depending only on these fixed caps, $s$ and
the chosen orders. To see this before imposing any prefix
precision, all normalized inverse norms are at most
$N^{a}$ after the two successive absorptions; converting
physical time and coordinates costs fixed powers of
$r^{-1},u^{-1},N^p$, the old positive floors and the frame
bounds. Also $\log(1/\kappa)=O(\log d)$,
$\log(1/u)\le(s+1)\log N$, $\log N\le2B$ under these
coarse caps, and $f_j\le B_j$ can be imposed below. These
give a fixed numerical upper bound on $C_0$ independent of
the price of a future request for old accuracy. The raw times
and both physical Green bounds may use the same enlarged $C_0$.
The three caps propagate as follows. The sites added at stage \(j\) are
the spectator, of modulus \(e^{C_1B_o+O(1)}\), the moderate and weak
sites, of modulus at most \(e^{kB_o+C\log B}\le N_j\) for \(Q>k+1\),
and the hot axes \(1+RP_i\); a normalized axial vertex has modulus
\(|V|/|\cos\theta_V|\le Cd\), so \(|1+RP_i|\le CdN_j\). Hence
\(\log(2+\max|n|)\le\log N_j+\log(Cd)+O(1)\le2\log N_j+O(1)\). The new
hot actions at the center of the final cube are at least
\(r^2b_p^2=H^2N_j^{-2s}e^{-2c_1 D_g}\ge N_j^{-2s-1}\).
For \(Q\) large the weak and moderate actions there are bounded below,
up to fixed factors, by \(\lambda^{2/5}e^{-C\Gamma}\ge N_j^{-1}\)
and \(K^{-h-1}\ge N_j^{-1}\), respectively. By (Pr3) the spectator
was selected with \(I_\ell\ge a_\ell=e^{-C_2B_o}\). Subsequent locks
retain a fixed fraction of this floor, which exceeds \(N_j^{-1}\)
for \(Q>C_2\) and large \(B\). The inherited actions likewise retain
a fixed fraction of their (Pr2) lower bound. Moving within the final
cube costs at most \(Cf_2\delta_1\le N_j^{-4s-29}\), so
\(-\log I_{\min}\le(2s+4)\log N_j+O(1)\). Finally the stage adds
one spectator, at most \(C_sd_j^2\) moderate and weak sites and at most
\(Cd_j^2\) hot axes, so \(f_j\le f_{j-1}+C_s(1+d_j)^2\).

Choose $L_0,L_1,C_w,D_h$ in the order specified in the physical
history proof, and any larger prefix accuracy $D'$ needed
before reconstructing the old coordinates and clocks. Fix the
coarse accuracy with which the finite graph is perturbed,
including the designated phase/amplitude and terminal-action
precision $e^{-C_4dD}$. When selecting the spectator the
component heights are at most a fixed power of $dP_{\mathrm{lat}}$.
Demand $d^2\log B\le c_* B_o$ with $c_*>0$ small enough
for their logarithms to fit in the prescribed old-scale budget.
At this stage the initial spectator size, its geometric and
affine exclusion widths, and the passive moderator placement
constant $k$ in $K=e^{kB_o+O_s(\log B)}$ can be chosen in
that order. None depends on $Q,\alpha$. Each old floor,
maximum index, strip count, spectator interval and moderator
precision used here costs at most $e^{C(s,\epsilon,D_h)B_o}$
times a polynomial in $B$; for the angular count one also
allows $e^{C f\log B}$. The latter, including the old sites,
is dominated by $\exp(C_s d^2\log B)$ plus the harmless
initial two-site factor. The row arity is fixed throughout.

Finally take $Q\ge2$ sufficiently large for the $B_o/B=1/Q$
fractions of all the preceding fixed bills, including the old
Sobolev comparisons and $L_0 B_o$ at the cold time. Take
$\alpha>0$ so small that the $d^2\log B/B\le\alpha^2$
fractions and $\alpha^2 Q$ obey all the indicated budgets.
Increase $B_*$ last. This accommodates any fixed powers of
$d,B,e^{Cd}$ and $e^{CD}$ not already included. In particular
$e^{C dD}$ times the normalized finite perturbations tends to
zero. Uniform finite powers of high index polynomials times a
previous-site price are accounted for once at each construction;
there is no multiplication of one whole stage inverse at each
level inside the finite graph. The only multiplication of old
stage inverses across stages is the one estimated by
$\sum_{i<j}B_i\le2B_{j-1}$ in the history proof.
For orientation we record the arithmetic behind these budgets. Since
$d\le\alpha\sqrt{B/\log B}$, for large $B$
\[
\begin{aligned}
 d^2\log(2+B)&\le\alpha^2B\,\frac{\log(2+B)}{\log B},\\
 dD=C_*d^2\log(d+2)&\le C_*\alpha^2\,\frac{B}{\log B}\,\log(2+\alpha\sqrt B)
 \le C_*\alpha^2B,
\end{aligned}
\]
while $B_o=B/Q$ and $\log(1/\lambda)=\epsilon B$. Thus each of the four
logarithmic prices listed in the assembly below is at most
$\bigl(C_{s,\epsilon}\alpha^2+C_sC_*\alpha^2+C_{s,\epsilon,D_h}Q^{-1}+C_s\epsilon\bigr)B$
up to lower order terms, whereas the gains are fixed powers of $N^{-1}$,
with $\log N=B+O(d^2\log B)=B(1+O(\alpha^2))$. The choice of
$\epsilon$ has already reserved the required margins for the term
$C_s\epsilon B$. The field-only estimate for the homogeneous weak tail
contains only the first three prices, whose sum can also be required
to be below $h_w\log N/4$ once $h_w$ is fixed. Taking $Q$ large and
then $\alpha$ small with $\alpha^2Q\le c_*$ fits those three prices
within all the remaining margins. Increase $B_*$ last to absorb
the $o(B)$ terms.

At each stage a vertical spacing $R_j$ with
$\log R_j=B_j+O(1)$, prime to five and a multiple of
$R_{j-1}$, exists by rounding a suitable multiple of $R_{j-1}$
in $[e^{B_j},2e^{B_j}]$; initialize with $R_0=1$. Choose
the same mod-five isotropic line at every step. The signed
trilinear equation preserves these nested cosets, and the
preceding estimates compare the old and new autonomous flows
without external forcing.

\begin{proof}[Assembly of Proposition~\ref{rate:prop:joint-insert}]
At stage $j$ the two-tree geometry, the prescribed positive
axial radii, and the angles determining the edge signs are chosen
on one completed terminal torus. The bounds on the logarithmic
prices of this choice are
\[
\begin{array}{ll}
 \text{lattice denominator and counting the short words}
      & C_{s,\epsilon}d^2\log(2+B),\\
 \text{finite Green and required finite accuracy}
      & C_s dD,\quad D=C_*d\log(d+2),\\
 \text{selection on earlier axes, including the spectator}
      & C_{s,\epsilon,D_h}B_o,\\
 \text{inverse powers of the current weak floors or slow scales}
      & C_s\epsilon B\quad\text{for the powers of }\lambda^{-1}.
\end{array}
\]
The additional powers of $r^{-1}$ from the change of units occur
in the explicit terms $N^{-J_0}/u^2$, $N^{A-A'}/u$ and in the
conversion of norms. They are not charged as inverse powers of
$\lambda$. The three other small terms to check are
$rN^2=HN^{2-s}$, $N^{-1/4}$ and $H N^{m-s}$, with
$m=\lceil s\rceil-1<s$. The last occurs only in passing from
the nonlinear axial chart to physical Fourier coefficients. For
all required products the choice of the old and arithmetic fractions
can also reserve a fraction of $s-m$. For the homogeneous weak
angular tail the inverse used for its first absorption has only
the first three prices. This is why its small exponent need not
exceed $C_s\epsilon$.

Lemma~\ref{rate:l:finite} first gives a reference landing at the specified
axial radii. With the actual harmonics and the final rounded carrier
rates, the divided rows and Lemma~\ref{rate:l:round} give the finite
solution at weight $\kappa/3$. Section~\ref{rate:sec:exterior} divides
the independent two-known sources and both entries of their signed
boundary Hessians using the same physical mismatch. Its remaining
same-copy operator is Hermitian on the exterior components.
Section~\ref{rate:sec:ward-landing} then constructs $W_j$ at weight
$\kappa/4$ and its compatible physical variational inverses,
including the one at weight $u\kappa/8$. The original polar
action and the angle shear are both included there. This gives
the first assertion by \eqref{rate:eq:raw-physical}. The direct physical
root estimate and \eqref{rate:eq:cold-physical} give the hot and
cold assertions on this very $W_j$. The same trajectory satisfies
\eqref{rate:eq:raw-weak}, because the frequencies and angle origins of
every earlier occupied coordinate were fixed in the terminal choice.
The initial two-site torus and the final direct-rate cube verify
the initial and successive induction conditions, respectively.
\end{proof}

\begin{proof}[Proof of Theorem~\ref{thm:two-sided}]
The first restored orbit is \(U_1=W_1\), which is uniformly bounded
in \(H^s\) by Proposition~\ref{rate:prop:joint-insert}. The successive
comparisons in \(H^2\) give
\[
 \sup_{j\ge1,\,t\ge0}\|U_j(t)\|_{H^2}
 \le \sup_{t\ge0}\|W_1(t)\|_{H^2}+\sum_{i\ge2}\varepsilon_i<\infty .
\]
Apply Proposition~\ref{rate:p:back} to the solutions $U_j$ just
constructed. In its notation take $\delta_j=\varepsilon_j$,
$c_j=O(e^{-c d_j})$, $k_j=(1+N_j,0)$, and write
$D_*=D_h-2C_w$. We verify the summability and the persistence
conditions explicitly. First,
\[
 (2K_j)^{s-2}\delta_j
 =2^{s-2}\exp\big(((s-2)L_1-D_*)B_{j-1}\big).
\]
This is summable by the choice of $D_h$. Next,
$S_jK_j^{-1/2}=\exp(-(L_1/2-L_0)B_{j-1})$.
Since $L_1>2L_0+2$ and $d_{j-1}=o(B_{j-1})$, the series
$\sum_j \sqrt{\omega_j}\,e^{A'd_{j-1}}$ converges for every fixed
$A'$. Also $\sum_j c_j<\infty$, since $B_j$ grows geometrically.

The root of $U_j$ has normalized Fourier modulus at least
$0.9e^{bd_j}-e^{-B_{j-1}}$ for all $t\ge T_j$. As
$\log N_j\le 2B_j$, the possible loss at this frequency in all
later stages is at most
\[
 \bk{k_j}^{s-2}\sum_{i>j}\delta_i
 \le C_s e^{2(s-2)B_j}\sum_{i>j} e^{-D_* B_{i-1}}
 \le 2C_s e^{-(D_*-2(s-2))B_j}.
\]
Here $D_*>2(s-2)$ by the choice of $D_h$ and $C_0\ge3$.
The hypothesis of Proposition~\ref{rate:p:back} follows with any
$0<a<b$ after omitting finitely many stages. Keeping the actual
root height in its lower-bound argument gives, for the single
limiting solution and every $T_j\le t<T_{j+1}$,
\[
 bd_j-O(1)\le\log\|u(t)\|_{H^s}\le C'd_{j+1}+O(1).
\]
Throughout that interval,
$(C_0+1)B_j\le\log t<(C_0+1)Q B_j$.
Thus $\log\log t=\log B_j+O(1)$ with a bound independent of $j$.
The definition of $d_j$ now yields \eqref{eq:two-sided} on the
union of all sufficiently late intervals. Explicitly, put
$\mathsf L=\log t$ and $\varkappa=C_0+1$, so that
$\mathsf L/(\varkappa Q)<B_j\le\mathsf L/\varkappa$. Since $x/\log x$
increases for $x>e$ and
$\tfrac12\log\mathsf L\le\log(\mathsf L/\varkappa)\le\log\mathsf L$ for
$\mathsf L\ge\varkappa^2$,
\[
 \frac{\mathsf L}{\varkappa Q\log\mathsf L}\le\frac{B_j}{\log B_j}
 \le\frac{2\mathsf L}{\varkappa\log\mathsf L},\qquad
 \frac{B_{j+1}}{\log B_{j+1}}\le\frac{QB_j}{\log B_j}
 \le\frac{2Q\mathsf L}{\varkappa\log\mathsf L}.
\]
With $d_j=\lfloor\alpha\sqrt{B_j/\log B_j}\rfloor$ this gives
\[
 d_j\ge\frac{\alpha}{\sqrt{\varkappa Q}}
        \sqrt{\frac{\log t}{\log\log t}}-1,
 \qquad
 d_{j+1}\le\alpha\sqrt{\frac{2Q}{\varkappa}}
        \sqrt{\frac{\log t}{\log\log t}},
\]
and the two bounds for $\log\|u(t)\|_{H^s}$ become
\eqref{eq:two-sided} with any $c_s<b\alpha/\sqrt{\varkappa Q}$ and
$C_s>C'\alpha\sqrt{2Q/\varkappa}$, for $t\ge t_s$. Returning from the mass
gauge multiplies the solution by a scalar of modulus one.
\end{proof}

\appendix
\section{The local frequency estimate at \texorpdfstring{$p=4$}{p=4}}\label{app:molnar-four}

We use the period-one normalization of \cite{Molnar}. For finitely supported
nonnegative action vectors \(J\), put \(\omega_n^*=\omega_n-(2\pi n)^2-4\sum_j J_j\).
We show that \(R(J)=\omega^*(J)+2J\) extends holomorphically from these
data near zero to an open complex \(\ell^2\) ball, with values in
\(\ell^{3/2}\).

\emph{Spectral tails.}
Use the period-two Fourier basis for the periodic spectral problem and
write \(q_j=|\varphi_j^+|+|\varphi_j^-|\) for the two potential
coefficients in the signed indexing of \cite[(1.15)]{Molnar}.
Set \(A=\|q\|_4\), \(T_N=\|(q_j)_{|j|\ge N}\|_4\), and let
\(\lambda_n^\pm\) be the two periodic eigenvalues near \(n\pi\).
For \(A\) sufficiently small and \(N\ge1\), we claim
\begin{equation}\label{eq:app-periodic-tail}
 \sum_{|n|\ge N,\,\pm}|\lambda_n^\pm-n\pi|^4
 \lesssim \sum_{|n|\ge N}q_{2n}^4+N^{-1}A^8+T_N^8 .
\end{equation}
The first term retains the direct Fourier coefficients in (1.15),
which are missing from the quantitative bound in
\cite[Proposition~6.1]{Molnar}.

To see this, first take a trigonometric polynomial potential.
Project onto the free double eigenspace at \(n\pi\) with projection
\(P\), and put \(Q=1-P\). Write \(L_0\) for the free operator and \(V\)
for the potential matrix. For \(|\lambda-n\pi|\le\pi/4\), let \(S\)
be a diagonal square root of \((\lambda-QL_0Q)^{-1}\) on the range
of \(Q\). Put \(h_0=0\),
\(h_j=|j|^{-1}\) for \(j\ne0\). In the free basis the entries coupling
indices \(n+j,n+k\) are bounded by \(q_{2n+j+k}\). Hence
\[
 \|SQVQS\|_{\mathrm{HS}}^2
 \lesssim\sum_{j,k}h_jh_k q_{2n+j+k}^2
 \le \|h*h\|_2 A^2,\qquad \|h*h\|_2\le \|h\|_{4/3}^2 .
\]
The inverse of \(1-SQVQS\) is uniformly bounded. Each row of
\(PVQS\) and column of \(SQVP\) has squared norm at most
\(C L_n\), where \(L_n=\sum_j h_jq_{2n+j}^2\).
Thus the Schur correction \(PVQS(1-SQVQS)^{-1}SQVP\) to \(PVP\)
has norm at most \(C L_n\).
At either root the reduced equation gives
\(|\lambda_n^\pm-n\pi|\lesssim q_{2n}+L_n\).
The same estimates on the boundary keep both roots in the disc by
continuation from zero. Approximation gives the bound for complex
Fourier--Lebesgue potentials with coefficients in \(\ell^4\), by
continuity of the isolated spectral pair
\cite[Sections~3 and~6]{Molnar}. For \(|n|\ge N\), split the sum in
\(L_n\) at \(|j|=N\). Its large-\(j\) part has \(\ell^4_n\) norm
at most \(C N^{-1/4}A^2\) by Young's inequality. In the other part
\(|2n+j|\ge N\), giving norm at most \(C T_N^2\). This proves
\eqref{eq:app-periodic-tail}. In particular the gap lengths
\(\gamma_n=\lambda_n^+-\lambda_n^-\) are locally bounded in \(\ell^4\).
The bounded reflection extensions in \cite[Sections~4--5]{Molnar}
give the same qualitative localization for the Dirichlet and auxiliary spectra.

\emph{Normalization of the discriminant.}
For a period-one potential, let \(P_{\rm e}\) and \(P_{\rm o}\) be the
products over even and odd \(n\), respectively, of
\((\lambda_n^+-\lambda)(\lambda_n^--\lambda)/\pi_n^2\), taken in
symmetric order, with \(\pi_n=n\pi\) for \(n\ne0\) and \(\pi_0=1\).
With this normalization the corrected forms of \cite[(2.4)]{Molnar} are
\[
 \Delta-2=-P_{\rm e},\qquad \Delta+2=4P_{\rm o}.
\]
At the vacuum the products are \(4\sin^2(\lambda/2)\) and
\(\cos^2(\lambda/2)\), by Euler's product. We verify the constants
and the asymptotic estimate needed from \cite[Lemma~8.1]{Molnar}.

Put \(a_n^\pm=\lambda_n^\pm-n\pi\). Off discs of any fixed radius
\(0<\rho<\pi/4\) about \(\pi\mathbb Z\), the ratio of either product
to its vacuum value is a product of factors
\(1+a_n^\pm/(n\pi-\lambda)\). It tends to one as \(|\lambda|\to\infty\).
Indeed \(\|((n\pi-\lambda)^{-1})_n\|_{4/3}\le C_\rho\) there.
For any \(\varepsilon>0\), \eqref{eq:app-periodic-tail} makes the
\(\ell^4\) norm of the tail of \(a^\pm\) less than \(\varepsilon\)
throughout a sufficiently small neighborhood of each fixed potential.
Hölder's inequality controls the tail of the sum of absolute deviations
by \(C_\rho\varepsilon\), and the finite part tends to zero uniformly
on that neighborhood. The grouped products are holomorphic in the
potential since the two roots in each disc enter through the characteristic
polynomial of its Riesz spectral subspace, and the products converge
locally uniformly on compact \(\lambda\)-sets by the same tail estimate.

For smooth potentials the trace of the monodromy satisfies
\(\Delta=2\cos\lambda+o(e^{|\operatorname{Im}\lambda|})\)
by \cite[Theorem~2.2]{GKbook}. The zeros of \(\Delta-2\) and
\(\Delta+2\), with multiplicities, are the even and odd spectral
pairs, respectively. Thus \((\Delta-2)/P_{\rm e}\) and
\((\Delta+2)/P_{\rm o}\) extend to entire functions. On
\(|\lambda|=(j+1/2)\pi\), the vacuum products both have modulus
comparable to \(e^{|\operatorname{Im}\lambda|}\). The two quotients
therefore converge uniformly on these circles to \(-1\) and \(4\),
so the maximum principle gives the displayed identities.
Density extends \(2-P_{\rm e}=4P_{\rm o}-2\) to the small complex
Fourier--Lebesgue ball and defines \(\Delta\) there.

The even product now gives
\(\Delta-2\cos\lambda=o(e^{|\operatorname{Im}\lambda|})\)
off the radius-\(\rho\) discs. Cauchy's estimate on circles of radius
\((\pi/4-\rho)/2\) gives
\(\Delta'+2\sin\lambda=o(e^{|\operatorname{Im}\lambda|})\)
when \(\operatorname{dist}(\lambda,\pi\mathbb Z)\ge\pi/4\).
There \(|\sin\lambda|\gtrsim e^{|\operatorname{Im}\lambda|}\), giving
the relative estimate for \(\Delta'\) used to locate its zeros.
These arguments also give \(\Delta^2-4=-4P_{\rm e}P_{\rm o}\) with
the normalization used in the contour integrals below.

\emph{Cancellation in the frequency formula.}
Let \(\tau_n=(\lambda_n^++\lambda_n^-)/2\),
\(G_n=[\lambda_n^-,\lambda_n^+]\), and
\(w_n(\lambda)^2=(\tau_n-\lambda)^2-\gamma_n^2/4\), with
\(w_n(\lambda)\sim\tau_n-\lambda\) at infinity.
Write \(\Omega_{nk}=\Omega_{nk}^{(2)}\) for the second contour moment
of \cite[Section~20]{Molnar}. In its defining integral the
normalized differential is multiplied by \(F_k^2\), where \(F_k\)
is the local Abelian integral of \cite[Section~19]{Molnar} vanishing
at the endpoints of \(G_k\). For real finite-gap data the frequency formula is
\(\omega_n^*=-(2/\pi)\sum_k\Omega_{nk}\)
\cite[Lemma~20.2]{Molnar}. For an open \(n\)th gap the bracket identity
used in that proof is \(2\{\Delta(\lambda),\theta_n\}=\psi_n(\lambda)\)
\cite[Lemma~18.2]{GKbook}. If the \(n\)th gap is closed, open only
the \(n\)th Birkhoff coordinate and let it tend to zero with all other
coordinates fixed. The nonzero moments in the sum are then supported
on the original finite set of open gaps together with \(n\)
\cite[Lemma~20.1]{Molnar}. Continuity of the frequency on this finite
action face and of each contour moment gives the formula at the limit.

For each \(m\), the critical point \(\lambda_m^\bullet\) of the
discriminant and the zeros \(\sigma_m^n\) of the numerator \(\psi_n\)
of the normalized differential
(for \(m\ne n\)) differ from \(\tau_m\) by \(O(|\gamma_m|^2)\),
uniformly in the indices on a small complex ball
\cite[Lemma~10.10 and Theorem~12.1]{Molnar}.
Set \(U_n=\sum_{m\ne n}|\gamma_m|^2/|m-n|\).
On the disc at \(n\), each factor \((s_m-\lambda)/w_m(\lambda)\)
in the products (2.18) and (4.9) of \cite{Molnar}, with
\(m\ne n\) and \(s_m=\lambda_m^\bullet\) or \(\sigma_m^n\), equals
\(1+O(|\gamma_m|^2/|m-n|)\). Indeed, divide numerator and
denominator by \(\tau_m-\lambda\) and use
\(|\tau_m-\lambda|\gtrsim|m-n|\). The products are thus \(1+O(U_n)\),
since \(U\) is bounded. On either bank of the cut, integration of
(2.18) from an endpoint gives
\[
 F_n=iw_n+O\bigl(|\gamma_n|(U_n+|\gamma_n|)\bigr).
\]
Here \(w_n'=-(\tau_n-\lambda)/w_n\) and
\(\int_{G_n}|d\lambda/w_n|=\pi\) on one bank.
On the lower bank, oriented from \(\lambda_n^-\) to \(\lambda_n^+\),
\(\int_{G_n^-} w_n\,d\lambda=i\pi\gamma_n^2/8\).
The diagonal contour integral reduces to twice the lower-bank integral
of \(F_n^2\,i(1+O(U_n))\,d\lambda/w_n\).
Its leading term is therefore \(\pi\gamma_n^2/4\).
In the action formula \cite[(2.19)]{Molnar} the leading term is
\(\gamma_n^2/4\), since the arcsine mean of \(t^2\) is \(1/2\)
and \(2(\lambda_n^\bullet-\tau_n)/\gamma_n=O(|\gamma_n|)\).
Consequently
\[
 |\Omega_{nn}-\pi J_n|\lesssim
 |\gamma_n|^2(U_n+|\gamma_n|),\qquad
 |\Omega_{nk}|\lesssim |\gamma_k|^3/|n-k|\quad(k\ne n).
\]
For the latter bound, (4.9) on \(G_k\) has one factor
\((\sigma_k^n-\lambda)/w_k\), with numerator \(O(|\gamma_k|)\),
and an unpaired \(1/w_n=O(1/|n-k|)\). The remaining product is
bounded and \(F_k=O(|\gamma_k|)\). Closed gaps follow by continuity
of the actions and moments. Young's inequality gives
\(U=h*|\gamma|^2\in\ell^3\) since \(h\in\ell^{6/5}\).
The diagonal error thus belongs to \(\ell^{4/3}\), and the
off-diagonal sum is bounded in \(\ell^{3/2}\) by
\(C\|h\|_{12/11}\||\gamma|^3\|_{4/3}\).
In \(\omega_n^*+2J_n\) the term \(-(2/\pi)\pi J_n\) cancels.
The bounds are locally uniform on the complex potential ball. For
fixed \(n\) the moment sum is holomorphic there by the contour formula
\cite[Lemma~20.1]{Molnar} and locally uniform convergence of the sum
(the \(\ell^4\) tail of \((h_{n-k})_k\) tends to zero).

\emph{Differentiation at an old action.}
Use the local holomorphic inverse Birkhoff map on the complex
\(\ell^4\) Cartesian ball \cite[Theorems~14.1 and~15.2]{Molnar},
with coordinate pairs \(z_j^+,z_j^-\) normalized by \(J_j=z_j^+z_j^-\).
On each finite coordinate slice, the moment formula on real finite-gap
data is invariant under rotations of the individual angles.
Holomorphic continuation from the real slice therefore leaves only
Taylor monomials with equal powers of \(z_j^+\) and \(z_j^-\) at each
site. These give compatible holomorphic functions of the finite action
vectors. Every such complex \(J\) has a lift with
\(|z_j^+|=|z_j^-|=|J_j|^{1/2}\), so that
\(\|z^+\|_4^2=\|z^-\|_4^2=\|J\|_2\).
Thus the \(\ell^{3/2}\) bound holds on a common \(\ell^2\) action
ball independently of the slice. Cauchy estimates give a common
Lipschitz bound on each smaller ball in the \(\ell^{3/2}\) norm, first on finite
output projections and then by taking their supremum. By density the slice functions
extend uniquely there, and the one-variable Cauchy formula along lines
passes to the limit, giving a holomorphic map \(R\) into \(\ell^{3/2}\).

For every fixed action \(J_*\) in this ball,
\(DR(J_*):\ell^2\to\ell^{3/2}\) is bounded and hence compact by
Pitt's theorem, as in \cite[proof of Corollary~18.2]{Molnar}.
Let \((\mathcal S I)_n=I_{-n}\). In the normalization of
Section~\ref{sec:tree}, \(J_*=(2\pi)^2\mathcal S I_*/2\) and
\[
 D\chi(I_*)+\mathrm{Id}
   =\tfrac12\mathcal S^{-1}DR(J_*)\mathcal S
   \quad\text{on }\ell^2,
\]
where the range of \(DR\) is included in \(\ell^2\).
This proves compactness also at each fixed nonzero old action.

\end{document}